%% file: MeridonMaster.tex
\documentclass[11pt]{book}

\input{macros_Meridon.tex}

\title{Equidistribution and counting under equilibrium states in
  negative curvature and trees. Applications to
  non-Archimedean Diophantine approximation}
\author{Anne Broise-Alamichel \and Jouni Parkkonen \and Fr\'ed\'eric
  Paulin} 
\date{\today}

\makeindex
\makeglossaries

\begin{document}
\bibliographystyle{../alphanum}

\maketitle

\tableofcontents

\input{meridonI.tex}

\input{meridonII.tex}

\input{meridonIII.tex}

\input{meridonIV.tex}

\appendix
\input{appendixBuzzi.tex}

\newpage
\setlength\glsdescwidth{0.78\linewidth}%
\printglossary[type=symbols,style=superragged3col]

\newpage
\phantomsection
\input{MeridonMaster.ind}

\newpage
{\small \bibliography{../biblio} } 
 
 \bigskip
{\small\

\noindent \begin{tabular}{l}
Laboratoire de math\'ematique d'Orsay, UMR 8628 Universit\'e Paris-Sud et CNRS\\
Universit\'e Paris-Saclay, 91405 ORSAY Cedex, FRANCE\\
{\it e-mail: anne.broise@math.u-psud.fr}
\end{tabular}

\medskip

\noindent \begin{tabular}{l} 

Department of Mathematics and Statistics, P.O. Box 35\\ 
40014 University of Jyv\"askyl\"a, FINLAND.\\
{\it e-mail: jouni.t.parkkonen@jyu.fi}
\end{tabular}

\medskip
\noindent \begin{tabular}{l} Laboratoire de math\'ematique d'Orsay,
  UMR 8628 Universit\'e Paris-Sud et CNRS\\ Universit\'e Paris-Saclay,
  91405 ORSAY Cedex, FRANCE\\ {\it e-mail:
    frederic.paulin@math.u-psud.fr}
\end{tabular}

}

 \end{document}

%% file: macros_Meridon.tex
\usepackage[T1]{fontenc}
\usepackage[utf8]{inputenc}
\usepackage{amsmath,amssymb,amscd,epsfig,bbm}
\usepackage{stmaryrd}
\usepackage[colorlinks]{hyperref}
\usepackage{mathabx}
\usepackage{comment}
\usepackage{color}

\usepackage[acronym,toc]{glossaries}
\usepackage{glossary-longragged}
\usepackage{glossary-superragged}
\newglossary[slg]{symbols}{sym}{sbl}{List of Symbols}

\usepackage{titlesec}
\titleformat{\subsection}[runin]{\bfseries}{}{}{}[]

\usepackage[all,cmtip]{xy}

\usepackage[textsize=small,textwidth=2cm]{todonotes}
\usepackage{varwidth}

\usepackage{url}

\usepackage{enumitem}
\setlist[itemize]{labelindent=*, leftmargin=.5 truecm,nosep}
\setlist[enumerate]{labelindent=*,leftmargin=1 truecm}

\usepackage[nottoc]{tocbibind}

\usepackage{fancyhdr}
\usepackage[us,12hr]{datetime} 
\fancypagestyle{plain}{
\fancyhf{}
\rfoot{ {\ddmmyyyydate\today}}
\lfoot{\thepage}
}

\numberwithin{equation}{chapter}

\usepackage{amsthm}

\theoremstyle{plain}
\newtheorem{defi}{Definition}[chapter]
\newtheorem{prop}[defi]{Proposition}
\newtheorem{theo}[defi]{Theorem}
\newtheorem{conj}[defi]{Conjecture}
\newtheorem{lemm}[defi]{Lemma}
\newtheorem{coro}[defi]{Corollary}

\theoremstyle{definition}
\newtheorem{rema}[defi]{Remark}
\newtheorem{exem}[defi]{Example}
\newtheorem{exems}[defi]{Examples}

\newcommand{\bdefi}{\begin{defi}}
\newcommand{\edefi}{\end{defi}}
\newcommand{\bprop}{\begin{prop}}
\newcommand{\eprop}{\end{prop}}
\newcommand{\btheo}{\begin{theo}}
\newcommand{\etheo}{\end{theo}}
\newcommand{\blemm}{\begin{lemm}}
\newcommand{\brema}{\begin{rema}}
\newcommand{\erema}{\end{rema}}
\newcommand{\bexer}{\begin{exem}}
\newcommand{\eexer}{\end{exem}}
\newcommand{\bexem}{\begin{exem}}
\newcommand{\eexem}{\end{exem}}
\newcommand{\bexems}{\begin{exems}}
\newcommand{\eexems}{\end{exems}}
\newcommand{\bconj}{\begin{conj}}
\newcommand{\econj}{\end{conj}}
\newcommand{\elemm}{\end{lemm}}
\newcommand{\bcoro}{\begin{coro}}
\newcommand{\ecoro}{\end{coro}}
\newcommand{\dem}{\noindent{\bf Proof. }}
\newcommand{\rem}{\noindent{\bf Remark. }}

\usepackage{mathrsfs}
\renewcommand\mathcal{\mathscr}

\newcommand{\A}{{\mathcal A}}
\newcommand{\B}{{\mathcal B}}
\newcommand{\C}{{\mathcal C}}
\newcommand{\D}{{\mathcal D}}
\newcommand{\E}{{\mathcal E}}
\newcommand{\F}{{\mathcal F}}
\newcommand{\G}{{\mathcal G}}

\renewcommand{\H}{{\mathcal H}}
\newcommand{\I}{{\mathcal I}}

\renewcommand{\L}{{\mathcal L}}
\newcommand{\M}{{\mathcal M}}
\newcommand{\N}{{\mathcal N}}
\newcommand{\OOO}{{\mathcal O}}

\newcommand{\Q}{{\mathcal Q}}

\newcommand{\U}{{\mathcal U}}
\newcommand{\V}{{\mathcal V}}

\newcommand{\maths}[1]{{\mathbb #1}}  

\newcommand{\CC}{\maths{C}}
\newcommand{\DD}{\maths{D}}
\newcommand{\FF}{\maths{F}}
\newcommand{\GG}{\maths{G}}
\newcommand{\HH}{\maths{H}}
\newcommand{\KK}{\maths{K}}
\newcommand{\LL}{\maths{L}}
\newcommand{\NN}{\maths{N}}

\newcommand{\PP}{\maths{P}}
\newcommand{\QQ}{\maths{Q}}
\newcommand{\RR}{\maths{R}}
\newcommand{\SSS}{\maths{S}}
\newcommand{\TT}{\maths{T}}

\newcommand{\XX}{\maths{X}}
\newcommand{\YY}{\maths{Y}}
\newcommand{\ZZ}{\maths{Z}}

\newcommand{\KKK}{{\mathfrak K}}

\newcommand{\mmm}{{\mathfrak m}}

\newcommand{\ppp}{{\mathfrak p}}

\newcommand{\uA}{\underline{A}}

\newcommand{\bs}{\backslash}
\newcommand{\ga}{\gamma}
\newcommand{\Ga}{\Gamma}
\newcommand{\ov}[1]{{\overline #1}} 
\newcommand{\ra}{\rightarrow}
\newcommand{\wt}[1]{{\widetilde{#1}}}
\newcommand{\wh}[1]{{\widehat{#1}}}

\newcommand{\Aut}{\operatorname{Aut}}
\newcommand{\Ax}{\operatorname{Ax}}

\newcommand{\bigO}{\operatorname{O}}
\newcommand{\card}{{\operatorname{Card}}}
\newcommand{\CAT}{\operatorname{CAT}}
\newcommand{\codeg}{\operatorname{codeg}}

\newcommand{\cqfd}{\hfill$\Box$}
\newcommand{\cyl}[1]{{[#1]}}
\newcommand{\dbs}{\backslash\mbox{\!\!\!\;}\backslash}
\newcommand{\Deg}{\operatorname{deg}}

\newcommand{\Dirac}{{\Delta}}
\newcommand{\discr}{\operatorname{Disc}}
\newcommand{\dvol}{\;d\operatorname{vol}}

\newcommand{\gengeod}
{\operatorname{\widecheck{\G\,}\!\!}}
\newcommand{\Geven}{\G_{\rm even}}
\newcommand{\gengeodeven}{\gengeod_{\rm even}}
\newcommand{\gengeododd}{\gengeod_{\rm odd}}
\newcommand{\haar}{\operatorname{Haar}}

\newcommand{\Heis}{\operatorname{Heis}}

\newcommand{\height}{\operatorname{ht}}
\newcommand{\green}{\mathfrak G}
\newcommand{\id}{\operatorname{id}}
\renewcommand{\Im}{{\operatorname{Im}}}
\newcommand{\inj}{\operatorname{inj}}
\newcommand{\Isom}{\operatorname{Isom}}
\newcommand{\Laplacian}{{\mathbf \Delta}}
\newcommand{\Leb}{\operatorname{Leb}}
\renewcommand{\Leb}{\operatorname \LL}
\newcommand{\Lebmeas}{\mathcal L}
\newcommand\len{\lambda}
\newcommand{\lk}{\operatorname{lk}}
\newcommand{\loops}[1]{L_{#1}}

\newcommand{\n}{\operatorname{\tt n}}

\newcommand{\Par}{\operatorname{Par}}
\newcommand{\per}{\operatorname{Per}}
\newcommand{\periodic}{\mathbf{Per}}
\newcommand\Perp{\operatorname{Perp}}

\newcommand{\pressure}[1]{P_{#1}}

\renewcommand{\Re}{{\operatorname{Re}}}
\newcommand{\redn}{\operatorname{\tt N}}
\newcommand{\redtr}{\operatorname{\tt Tr}}
\newcommand{\reg}{\operatorname{{\tt R}\!}}

\newcommand{\smallo}{\operatorname{o}}
\newcommand{\stab}{\operatorname{Stab}}
\newcommand{\Tamagawa}{\rm T}
\newcommand{\tr}{\operatorname{\tt tr}}
\newcommand{\Tvol}{{\operatorname{Tvol}}}
\newcommand{\TVol}{{\operatorname{TVol}}}
\newcommand{\Veven}{V_{\rm even}}
\newcommand{\Vodd}{V_{\rm odd}}
\newcommand{\Vol}{\operatorname{Vol}}
\newcommand{\vol}{\operatorname{vol}}

\newcommand{\weakstar}{\overset{*}\rightharpoonup}

\newcommand{\hdr}{{\HH}^2_\RR}

\newcommand{\hnr}{{\HH}^n_\RR}

\newcommand{\hnc}{{\HH}^n_\CC}

\newcommand{\GL}{\operatorname{GL}}
\newcommand{\PGL}{\operatorname{PGL}}
\newcommand{\SL}{\operatorname{SL}}
\newcommand{\PSL}{\operatorname{PSL}}
\newcommand{\SO}{\operatorname{SO}}

\newcommand{\flow}[1]{{{\tt g}^{#1}}}  
\newcommand{\wss}{W^+} 
\newcommand{\wsu}{W^-} 
\newcommand{\wssu}{W^\pm} 
\newcommand{\wsus}{W^\mp}
\newcommand{\ws}{W^{0+}} 
\newcommand{\wu}{W^{0-}} 
\newcommand{\wosu}{W^{0\pm}} 
\newcommand{\wous}{W^{0\mp}} 
\newcommand{\muss}[1]{{\mu_{\wss(#1)}}}
\newcommand{\musu}[1]{{\mu_{\wsu(#1)}}}
\newcommand{\mussu}[1]{{\mu_{\wssu(#1)}}}

\newcommand{\normal}[1]{\partial^1_{+}{#1}}

\newcommand\normalout{\partial^1_{+}}
\newcommand\normalin{\partial^1_{-}}
\newcommand\normalpm{\partial^1_{\pm}}
\newcommand\normalmp{\partial^1_{\mp}}

\usepackage{mathtools}
\makeatletter
\DeclareRobustCommand\widecheck[1]{{\mathpalette\@widecheck{#1}}}
\def\@widecheck#1#2{%
    \setbox\z@\hbox{\m@th$#1#2$}%
    \setbox\tw@\hbox{\m@th$#1%
       \widehat{%
          \vrule\@width\z@\@height\ht\z@
          \vrule\@height\z@\@width\wd\z@}$}%
    \dp\tw@-\ht\z@
    \@tempdima\ht\z@ \advance\@tempdima2\ht\tw@ \divide\@tempdima\thr@@
    \setbox\tw@\hbox{%
       \raise\@tempdima\hbox{\scalebox{1}[-1]{\lower\@tempdima\box
\tw@}}}%
    {\ooalign{\box\tw@ \cr \box\z@}}}
\makeatother

\usepackage{scalerel,stackengine}
\stackMath
\newcommand\reallywidehat[1]{%
\savestack{\tmpbox}{\stretchto{%
  \scaleto{%
    \scalerel*[\widthof{\ensuremath{#1}}]{\kern-.6pt\bigwedge\kern-.6pt}%
    {\rule[-\textheight/2]{1ex}{\textheight}}
  }{\textheight}%
}{0.5ex}}%
\stackon[1pt]{#1}{\tmpbox}%
}

\newcommand\reallywidetilde[1]{%
\savestack{\tmpbox}{\stretchto{%
  \scaleto{%
    \scalerel*[\widthof{\ensuremath{#1}}]{\kern-.6pt\sim\kern-.6pt}%
    {\rule[-\textheight/2]{1ex}{\textheight}}
  }{\textheight}%
}{0.5ex}}%
\stackon[1pt]{#1}{\tmpbox}%
}

\newcommand{\project}[1]
{\medskip
\begin{center} {\color{red} \framebox{%
   \begin{varwidth}{15cm}
  {\bf Project:} {\color{black}#1}
 \end{varwidth}}  }
 \end{center}
\medskip }

\newcounter{fig}

\newglossaryentry{pointprojectifinfty}{type=symbols,
  name={\ensuremath{\infty}}, sort=00,
  description={standard point at infinity $[1:0]$ of a projective line}}

\newglossaryentry{weakstarconv}{type=symbols,
  name={\ensuremath{\weakstar}}, sort=000,
  description={weak-star convergence of measures on locally compact spaces}}

\newglossaryentry{characteristicfunction}{type=symbols,
  name={\ensuremath{\mathbbm{1}_A}}, sort=1,
  description={characteristic function of a subset $A$}}

\newglossaryentry{equivrelatequivfamil}{type=symbols,
  name={\ensuremath{\sim\;=\;\sim_\D}}, sort=2,
  description={equivalence relation on index set of an equivariant
    family $\D$}}

\newglossaryentry{normealphaholder}{type=symbols,
  name={\ensuremath{\|f\|_\alpha}}, sort=3,
  description={$\alpha$-H\"older norm of $f\in \C_{\rm c}^\alpha
    (Z)$}}

\newglossaryentry{normellsobolev}{type=symbols,
  name={\ensuremath{\|\psi\|_\ell}}, sort=4, description={Sobolev
    $W^{\ell,2}$-norm of $\psi\in \C_{\rm c}^\ell (N)$}}

\newglossaryentry{absolutevalue}{type=symbols,
  name={\ensuremath{|\cdot|_v}}, sort=5, description={(normalised)
    absolute value associated with a valuation $v$}}

\newglossaryentry{cylinderlength}{type=symbols,
  name={\ensuremath{|w|}}, sort=5, description={length of the 
cylinder $[w]$ associated with an admissible sequence $w$}}

\newglossaryentry{basepointBTt}{type=symbols, name={\ensuremath{*_v}},
  sort=6, description={base point $*_v=[\OOO_v\times\OOO_v]$ of
    the Bruhat-Tits tree $\XX_v$}}

\newcommand\souligne{}
\DeclareRobustCommand\souligne[1]{\underline{#1}}

\newglossaryentry{maxcompactopendiag}{type=symbols,
  name={\ensuremath{\souligne{A}(\OOO_v)}}, sort=aoo,
  description={maximal compact-open diagonal subgroup of 
$\PGL_2(K_v)$}}

\newglossaryentry{automorphismgroupmetrictree}{type=symbols,
  name={\ensuremath{\Aut(\XX,\lambda)}}, sort=aut,
  description={automorphism group (edge-preserving, without
    inversion) of a metric tree $(\XX,\lambda)$}}

\newglossaryentry{automorphismgroupsimptree}{type=symbols,
  name={\ensuremath{\Aut \XX}}, sort=aut,
  description={automorphism group (without inversion) of a
    simplicial tree $\XX$}}

\newglossaryentry{closedball}{type=symbols,
  name={\ensuremath{B(x,r)}}, sort=ball,
  description={closed ball of center $x$ and radius $r$ in a metric space}}

\newglossaryentry{ballhamenstadt}{type=symbols,
  name={\ensuremath{B^\pm(w,\eta')}}, sort=ballhamen,
  description={Hamenstädt's ball of radius $\eta'>0$ with center
    any geodesic line extension of $w\in \G_\pm X$}}

\newglossaryentry{dynamball}{type=symbols,
  name={\ensuremath{B(\ell;T,T',r)}}, sort=ballhb,
  description={dynamical ball in the space of geodesic lines $\G X$}}

\newglossaryentry{smoothprojcurve}{type=symbols,
  name={\ensuremath{\bf C}}, sort=c,
  description={geometrically connected
  smooth projective curve over $\FF_q$}}

\newglossaryentry{complementary}{type=symbols,
  name={\ensuremath{^c A}}, sort=ca,
  description={complementary set of a subset $A$}}

\newglossaryentry{conducpotential}{type=symbols,
  name={\ensuremath{\wt c_F}}, sort=caa,
  description={system of conductances on $\XX$ associated with 
a potential $\wt F$}}

\newglossaryentry{conducpotentialdown}{type=symbols,
  name={\ensuremath{c_F}}, sort=caaa,
  description={system of conductances on $\Ga\bs\XX$ associated with 
a potential $F$}}

\newglossaryentry{periodconduct}{type=symbols,
  name={\ensuremath{c(g)}}, sort=caaaa,
  description={period for a system of conductances $c$ 
of a closed orbit $g$ for the geodesic flow}}

\newglossaryentry{espacecontcompsupp}{type=symbols,
  name={\ensuremath{\C_{\rm c} (Z)}}, sort=caaaaaa,
  description={space of real-valued continuous maps with compact
    support on $Z$}}

\newglossaryentry{espacealphaholdeb}{type=symbols,
  name={\ensuremath{\C_{\rm b}^\alpha (Z)}}, sort=calphaa,
  description={space of bounded $\alpha$-H\"older-continuous 
    real-valued functions on $Z$}}

\newglossaryentry{espacealphaholdebk}{type=symbols,
  name={\ensuremath{\C^{k,\,\alpha}_{\rm b} (Z)}}, sort=calphab,
  description={space of real-valued functions on $Z$ with 
bounded $\alpha$-H\"older-continuous derivatives of order 
at most $k$ along the flow}}

\newglossaryentry{espacealphaholder}{type=symbols,
  name={\ensuremath{\C_{\rm c}^\alpha (Z)}}, sort=calphac,
  description={space of $\alpha$-H\"older-continuous real-valued
    functions with compact support on $Z$}}

\newglossaryentry{espacealphaholderkc}{type=symbols,
  name={\ensuremath{\C^{k,\,\alpha}_{\rm c} (Z)}}, sort=calphad,
  description={space of real-valued functions on $Z$ with 
bounded $\alpha$-H\"older-continuous derivatives of order 
at most $k$ along the flow and compact support}}

\newglossaryentry{espacellsobolev}{type=symbols,
  name={\ensuremath{\C_{\rm c}^\ell(T^1M)}}, sort=cell,
  description={space of real-valued $\C^\ell$-smooth functions
    with compact support on $T^1M$
    }}

\newglossaryentry{catmoinsun}{type=symbols,
  name={\ensuremath{\CAT(-1)}}, sort=cellcat, description={metric
  space satisfying the Alexandrov-Topogonov comparison property with
  the real hyperbolic space of constant curvature $-1$}}
  
\newglossaryentry{codegre}{type=symbols,
  name={\ensuremath{\codeg_\DD(x)}}, sort=codeg,
  description={codegree of a vertex $x$ with respect to a subtree
    $\DD$}}

\newglossaryentry{cocycle}{type=symbols,
  name={\ensuremath{C^\pm_{\Ga,F^\pm}}}, sort=cocycle,
  description={(normalised) Gibbs cocycle associated with the group
  $\Ga$ and the potential $\wt F^\pm$}}

\newglossaryentry{codegrecalD}{type=symbols,
  name={\ensuremath{\codeg_\D(x)}}, sort=codegD,
  description={codegree of a vertex $x$ with respect to a family
    of subtrees $\D$}}

\newglossaryentry{convexhull}{type=symbols,
  name={\ensuremath{\C\Lambda\Ga}}, sort=convexhull,
  description={convex hull in $X$ of the limit set $\Lambda\Ga$ of
    $\Ga$}}

\newglossaryentry{correlationcoeff}{type=symbols,
  name={\ensuremath{\operatorname{cov}_{m,\,n}(\phi,\psi)}},
  sort=conva, description={$n$-th correlation coefficient of two
  observables $\phi,\psi$ for the measure $m$ under a transformation}}

\newglossaryentry{correlationcoeffcont}{type=symbols,
  name={\ensuremath{\operatorname{cov}_{\mu,\,t}(\psi,\psi')}},
  sort=convab, description={correlation coefficient at time $t$ of
  $\psi,\psi'$ for the measure $\mu$ under a flow }}

\newglossaryentry{correlationcoeffalgeb}{type=symbols,
  name={\ensuremath{\operatorname{cov}_{\mu_v,\,g}}},
  sort=convac, description={correlation coefficient for $g\in G_v$ 
and a measure $\mu_v$ on $\Ga\bs G_v$}}

\newglossaryentry{conjugate}{type=symbols,
  name={\ensuremath{\overline{x}}}, sort=convjugate,
  description={conjugate of a quaternion $x$}}

\newglossaryentry{crossratio}{type=symbols,
  name={\ensuremath{[a,b,c,d]}}, sort=crossratio,
  description={crossratio of pairwise distinct points
    $a,b,c,d$ in $K_v$}}

\newglossaryentry{absolcrossratio}{type=symbols,
  name={\ensuremath{|a,b,c,d|_v}}, sort=crossratioa,
  description={absolute crossratio of pairwise distinct points
    $a,b,c,d$ for a valuation $v$}}

\newglossaryentry{crossratiotree}{type=symbols,
  name={\ensuremath{\ldbrack\xi_1,\xi_2,\xi_3,\xi_4\rdbrack}}, sort=crossratiot,
  description={(logarithmic) crossratio of pairwise distinct points
    $\xi_1,\xi_2,\xi_3,\xi_4$ in $\partial_\infty\XX$}}

\newglossaryentry{boundaryatinfinity}{type=symbols,
  name={\ensuremath{\partial_\infty X}}, sort=daaa,
  description={space at infinity of $X$}}

\newglossaryentry{boundaryatinfinitybis}{type=symbols,
  name={\ensuremath{\partial^2_\infty X}}, sort=daaaa,
  description={space of distinct ordered pairs of points at infinity of $X$}}

\newglossaryentry{boundarytowardse}{type=symbols,
  name={\ensuremath{\partial_e \XX}}, sort=daaaaa, description={set
    of points at infinity of the geodesic rays whose initial
    (oriented) edge is $e$}}

\newglossaryentry{boundaryVD}{type=symbols, name={\ensuremath{\partial
      V\DD}}, sort=daaaaaa, description={boundary
    of set of vertices of a simplicial subtree $\DD$}}

\newglossaryentry{boundaryD}{type=symbols, name={\ensuremath{\partial
      \DD}}, sort=daaaaaaa, description={maximal subgraph with set
    of vertices $\partial V\DD$}}

\newglossaryentry{outerunitnormalbundle}{type=symbols,
  name={\ensuremath{\normalout D}}, sort=dbp,
  description={outer unit normal bundle of a closed convex subset $D$}}

\newglossaryentry{innerunitnormalbundle}{type=symbols,
  name={\ensuremath{\normalin D}}, sort=dbm,
  description={inner unit normal bundle of a closed convex subset $D$}}

\newglossaryentry{degreevaluation}{type=symbols,
  name={\ensuremath{\Deg v}}, sort=degre,
  description={degree of valuation $v$, equal to $\dim_{\FF_q}k_v$}}

\newglossaryentry{criticalexponent}{type=symbols,
  name={\ensuremath{\delta=\delta_{\Ga,\,F}}}, sort=deltaa,
  description={critical exponent of $(\Ga,F)$}}

\newglossaryentry{laplacianconductanceup}{type=symbols,
  name={\ensuremath{\Laplacian_{\widetilde{c}^{\,\pm}}}}, sort=deltaba,
  description={Laplacian operator associated with a system of
    conductances $\widetilde{c}^\pm$ on a simplicial tree}}
    
    \newglossaryentry{laplacianconductancedown}{type=symbols,
  name={\ensuremath{\Laplacian_{c}}}, sort=deltabb,
  description={Laplacian operator associated with a system of
    conductances $c$ on a graph of groups}}

\newglossaryentry{diracmass}{type=symbols,
  name={\ensuremath{\Dirac_x}}, sort=deltac,
  description={unit Dirac mass at a point $x$}}

\newglossaryentry{discred}{type=symbols,
  name={\ensuremath{\discr_D}}, sort=disc,
  description={reduced discriminant of a quaternion algebra 
$D$ over $\QQ$}}

\newglossaryentry{distancelike}{type=symbols,
  name={\ensuremath{d_{D}}}, sort=distancea,
  description={distance-like map on $\partial_\infty X-
\partial_\infty D$ associated with a closed convex subset $D$}}

\newglossaryentry{distGchechX}{type=symbols,
  name={\ensuremath{d=d_{\gengeod X}}}, sort=distanceaa,
  description={distance on the space of generalised geodesic lines}}

\newglossaryentry{distgerm}{type=symbols,
  name={\ensuremath{d=d_{T^1X}}}, sort=distanceaaa,
  description={distance on the space of germs of geodesic lines}}

\newglossaryentry{Hamenstadtdistanceinfinity}{type=symbols, 
name={\ensuremath{d_\H}},
  sort=distanceb, description={Hamenst\"adt's
  distance at infinity associated with an horoball $\H$}}

\newglossaryentry{visualdist}{type=symbols, name={\ensuremath{d_x}},
  sort=distancec, description={visual distance on $\partial_\infty
    X$ seen from $x\in X$}}

\newglossaryentry{Hamenstadtdistancepm}{type=symbols, 
name={\ensuremath{d_{\wssu(w)}}},
  sort=distanced, description={Hamenst\"adt's
  distance on the strong stable/unstable ball of $w\in\G_\pm X$}}

\newglossaryentry{setofedges}{type=symbols, name={\ensuremath{E\XX}},
  sort=edge, description={set of edges of a graph $\XX$}}

\newglossaryentry{Eulerfunct}{type=symbols, name={\ensuremath{\varphi_{R_v}}},
  sort=euler, description={Euler function of function ring $R_v$}}

\newglossaryentry{negativepart}{type=symbols,
  name={\ensuremath{f^-}}, sort=fbb, description={negative part of a
  real-valued map $f$}}

\newglossaryentry{fibrationpm}{type=symbols,
  name={\ensuremath{f^\pm_D}}, sort=fib, description={fibration
    over $\normalpm D$ with fibers the stable/unstable leaves}}

\newglossaryentry{potentialup}{type=symbols,
  name={\ensuremath{\wt F}}, sort=fp, description={potential on $T^1X$}}

\newglossaryentry{potentialdown}{type=symbols,
  name={\ensuremath{F}}, sort=fpd, description={potential on $\Ga\bs T^1X$}}

\newglossaryentry{potentialconduc}{type=symbols,
  name={\ensuremath{\wt F_c}}, sort=fpda, description={potential on
  $T^1X$ associated with a system of conductances $\wt c$}}

\newglossaryentry{potentialconducdown}{type=symbols,
  name={\ensuremath{F_c}}, sort=fpdadown, description={potential on
  $\Ga\bs T^1X$ associated with a system of conductances $c$}}

\newglossaryentry{finitefield}{type=symbols,
  name={\ensuremath{\FF_q}}, sort=fpdfield, description={finite
    field of order a prime power $q$}}

\newglossaryentry{gammav}{type=symbols,
  name={\ensuremath{\Ga_v}}, sort=gammav, description={modular
group at a valuation $v$ of a function field}}

\newglossaryentry{genus}{type=symbols, name={\ensuremath{g}}, 
sort=gb, description={genus of the smooth projective curve $\bf C$}}

\newglossaryentry{spacegeodlin}{type=symbols, 
name={\ensuremath{\G X}}, sort=gengeo, description={space of 
    geodesic lines in $X$}}

\newglossaryentry{gengeod}{type=symbols, name={\ensuremath{\gengeod
      X}}, sort=gengeod, description={space of generalised
    geodesic lines in $X$}}

\newglossaryentry{gengeoddiscrete}{type=symbols,
  name={\ensuremath{\gengeod\XX}}, sort=gengeoddis,
  description={space of generalised discrete geodesic lines in a
    simplicial tree $\XX$}}

\newglossaryentry{evendiscgeodflowspace}{type=symbols,
  name={\ensuremath{\Geven\XX}}, sort=gengeoddiseven,
  description={space of discrete geodesic lines $\ell$ in $\XX$
  with $\ell(0)$ at even distance from $x_0$}}

\newglossaryentry{evendiscgenegeodflowspace}{type=symbols,
  name={\ensuremath{\gengeodeven\XX}}, sort=gengeoddiseven,
  description={space of generalised discrete geodesic lines 
  $\ell$ in $\XX$ with $d(\ell(0),x_0)$  even}}

\newglossaryentry{gengeodpm}{type=symbols, name={\ensuremath{\G_\pm
      X}}, sort=gengeodpm, description={space of generalised
    positive/negative geodesic rays in $X$}}

\newglossaryentry{gengeodpmz}{type=symbols,
  name={\ensuremath{\G_{\pm,\,0} X}}, sort=gengeodpmz,
  description={space of generalised geodesic lines in $X$
    isometric exactly on $\pm\mathopen{[}0,+\infty\mathclose{[}$}}

\newglossaryentry{geodesicflow}{type=symbols,
  name={\ensuremath{(\flow{t})_{t\in\RR}}}, sort=gt,
  description={(continuous time) geodesic flow on space of
    generalised geodesics $\gengeod X$}}

\newglossaryentry{geodesicflowdown}{type=symbols,
  name={\ensuremath{(\flow{t})_{t\in\RR}}}, sort=gta,
  description={(continuous time) geodesic flow on the quotient space of
    generalised geodesics $\Ga\bs\gengeod X$}}

\newglossaryentry{geodesicflowdis}{type=symbols,
  name={\ensuremath{(\flow{t})_{t\in\ZZ}}}, sort=gtdis,
  description={(discrete time) geodesic flow on space of
    generalised geodesics $\gengeod \XX$, as well as on 
$\Ga\bs\gengeod \XX$}}

\newglossaryentry{metricentropydisc}{type=symbols,
  name={\ensuremath{h_m(T)}}, sort=hab, description={metric
    entropy of a transformation $T$ with respect to a probability measure
    $m$}}

\newglossaryentry{metricentropy}{type=symbols,
  name={\ensuremath{h_m(\phi_1)}}, sort=hac, description={metric
    entropy of a flow $(\phi_t)_{t\in\RR}$ with respect to a probability measure
    $m$}}

\newglossaryentry{complexityquadirrat}{type=symbols,
  name={\ensuremath{h(\alpha)}}, sort=hacomplexityi, 
  description={complexity of a quadratic irrational $\alpha$ in $K_v$}}

\newglossaryentry{complexityquadirratqp}{type=symbols,
  name={\ensuremath{h(\alpha)}}, sort=hacomplexityiqp, 
  description={complexity of a quadratic irrational $\alpha$ in $\QQ_p$}}

\newglossaryentry{complexityloxo}{type=symbols,
  name={\ensuremath{h(\alpha)}}, sort=hacomplexityb, 
  description={complexity of a loxodromic fixed point $\alpha$}}

\newglossaryentry{relheightquadirrat}{type=symbols,
  name={\ensuremath{h_{\alpha}(\beta)}}, sort=hairratquad, 
  description={relative height of a quadratic irrational
$\beta$ with respect to $\alpha$}}

\newglossaryentry{relheightloxo}{type=symbols,
  name={\ensuremath{h_{\alpha}(\beta)}}, sort=hailoxo, 
  description={relative height of a loxodromic fixed point
$\beta$ with respect to $\alpha$}}

\newglossaryentry{haarKv}{type=symbols,
  name={\ensuremath{\haar_{K_v}}}, sort=har, 
  description={normalised Haar measure of $(K_v,+)$}}

\newglossaryentry{height}{type=symbols,
  name={\ensuremath{{\height_\infty}}}, sort=hat, 
  description={height of a horoball in the Bruhat-Tits tree $\XX_v$}}
  
\newglossaryentry{heisenberggroup}{type=symbols,
  name={\ensuremath{{\Heis_{2n-1}}}}, sort=hath, 
  description={Heisenberg group of dimension $2n-1$}}

\newglossaryentry{horoshrink}{type=symbols,
  name={\ensuremath{\H[t]}}, sort=hb,
  description={horoball contained in $\H$ whose boundary is at
    distance $t$ from the boundary of $\H$}}

\newglossaryentry{horoballstable}{type=symbols,
  name={\ensuremath{H\!B_{+}(w)}}, sort=horoballstable,
  description={stable horoball of $w\in\G_+ X$}}

\newglossaryentry{horoballunstable}{type=symbols,
  name={\ensuremath{H\!B_{-}(w)}}, sort=horoballunstable,
  description={unstable horoball of $w\in\G_- X$}}

\newglossaryentry{horospherestable}{type=symbols,
  name={\ensuremath{H_{+}(w)}}, sort=horospherestable,
  description={stable horosphere of $w\in\G_+ X$}}

\newglossaryentry{horosphereunstable}{type=symbols,
  name={\ensuremath{H_{-}(w)}}, sort=horosphereunstable,
  description={unstable horosphere of $w\in\G_- X$}}

\newglossaryentry{complexhyperbolicspace}{type=symbols,
  name={\ensuremath{{\hnc}}}, sort=hot, 
  description={complex hyperbolic space of dimension $n$}}

\newglossaryentry{antipodalmap}{type=symbols,
  name={\ensuremath{\iota}}, sort=iota, description={antipodal map
    $w\mapsto\{t\mapsto w(-t)\}$ on $\gengeod X$, as well as on $T^1X$}}

\newglossaryentry{antipodalmapdown}{type=symbols,
  name={\ensuremath{\iota}}, sort=iotad, description={antipodal map
    $\Ga w\mapsto\{t\mapsto \Ga w(-t)\}$  on
    $\Ga\bs\gengeod X$, as well as on $\Ga\bs T^1X$}}

\newglossaryentry{reciprocityindexloxo}{type=symbols,
  name={\ensuremath{\iota_\alpha}}, sort=iotare,
  description={reciprocity index of a loxodromic fixed point $\alpha$}}

\newglossaryentry{reciprocityindex}{type=symbols,
  name={\ensuremath{\iota_G(\beta)}}, sort=iotarec,
  description={$G$-reciprocity index of a quadratic irrational $\beta$}}

\newglossaryentry{isometrygroup}{type=symbols,
  name={\ensuremath{\Isom(X)}}, sort=isom, description={isometry
    group of $X$}}

\newglossaryentry{fracidealclassset}{type=symbols,
  name={\ensuremath{\I_v}}, sort=iv, description={set of 
classes of fractional ideals of $R_v$}}

\newglossaryentry{functionfield}{type=symbols, name={\ensuremath{K}},
  sort=kk, description={global function field over $\FF_q$}}

\newglossaryentry{completfunctionfield}{type=symbols,
  name={\ensuremath{K_v}}, sort=kkv, description={completion of
    function field $K$ for the valuation $v$}}

\newglossaryentry{residualfield}{type=symbols,
  name={\ensuremath{k_v}}, sort=kkvpetit, description={residual
    field of the valuation $v$ on $K$}}

\newglossaryentry{distancetranslat}{type=symbols,
  name={\ensuremath{\lambda(\ga)}}, sort=l, description={translation
    length of an isometry $\ga$ of $X$}}

\newglossaryentry{periodlength}{type=symbols,
  name={\ensuremath{\lambda(g)}}, sort=la, description={length of
    a periodic orbit $g$}}

\newglossaryentry{limitset}{type=symbols,
  name={\ensuremath{\Lambda\Ga}}, sort=lba, description={limit set
    of a discrete group of isometries $\Ga$ of $X$}}

\newglossaryentry{connicallimitset}{type=symbols,
  name={\ensuremath{\Lambda_{\rm c}\Ga}}, sort=lbaa, description={conical
    limit set of a discrete group of isometries $\Ga$ of $X$}}

\newglossaryentry{ldeuxgraphgroup}{type=symbols,
  name={\ensuremath{\Leb^2(\YY,G_*)}}, sort=lbaaa,
  description={Hilbert space of square integrable maps on $V\YY$
    for the measure $\vol_{(\YY,G_*)}$}}

\newglossaryentry{lebmeasperio}{type=symbols,
  name={\ensuremath{\Lebmeas_g}}, sort=lc,
  description={Lebesgue measure along a periodic orbit $g$}}

\newglossaryentry{link}{type=symbols, name={\ensuremath{\lk}},
  sort=lk, description={links of vertices in simplicial trees}}

\newglossaryentry{lengthspectrum}{type=symbols,
  name={\ensuremath{\loops{\Ga}}}, sort=lo, description={length
    spectrum of action of $\Ga$ on $X$}}

\newglossaryentry{ln}{type=symbols, name={\ensuremath{\ln}}, sort=log,
  description={natural logarithm (with $\ln(e)=1$)}}

\newglossaryentry{endpointmapline}{type=symbols,
  name={\ensuremath{\ell_\pm}}, sort=lpm,
  description={positive/negative endpoint of geodesic line
    $\ell$}}

\newglossaryentry{ellstar}{type=symbols,
  name={\ensuremath{\ell^*}}, sort=ls,
  description={standard basepoint in space of geodesic lines 
$\G\XX_v$}}

\newglossaryentry{multiplicitysimplicial}{type=symbols,
  name={\ensuremath{m_{\D}(x)}}, sort=msubd,
  description={multiplicity of a vertex $x$ with respect to an
    equivariant family $\D$ }}

\newglossaryentry{gibbsmeasure}{type=symbols,
  name={\ensuremath{\wt m_{F}}}, sort=msubf, description={Gibbs
    measure on the space of geodesic lines $\G X$}}

\newglossaryentry{gibbsmeasuredown}{type=symbols,
  name={\ensuremath{m_{F}}}, sort=msubfd, description={Gibbs
    measure on the quotient space of geodesic lines $\Ga\bs \G X$}}

\newglossaryentry{gibbsmeasuresimplicial}{type=symbols,
  name={\ensuremath{\wt m_{F}}}, sort=msubfsimp, description={Gibbs
    measure on the space of discrete geodesic lines $\G\XX$}}

\newglossaryentry{gibbsmeasuresimplicialdown}{type=symbols,
  name={\ensuremath{m_{F}}}, sort=msubfsimpd, description={Gibbs
    measure on the quotient space of discrete geodesic lines 
    $\Ga\bs\G\XX$}}

\newglossaryentry{gibbsmeasuredownrenorm}{type=symbols,
  name={\ensuremath{\overline{m_F}}}, sort=msubft,
  description={renormalised Gibbs measure $m_F/||m_F||$ on
    $\Ga\bs\G X$}}

\newglossaryentry{gibbsmeasuredownrenormc}{type=symbols,
  name={\ensuremath{\overline{m_c}}}, sort=msubftc,
  description={renormalised Gibbs measure $m_c/||m_c||$ on
    $\Ga\bs\G \XX$}}

\newglossaryentry{Pattersondensity}{type=symbols,
  name={\ensuremath{(\mu_x^\pm)_{x\in X}}}, sort=mu,
  description={(normalised) Patterson density for the pair
    $(\Ga,\wt F^\pm)$}}

\newglossaryentry{strongmeasure}{type=symbols,
  name={\ensuremath{\mussu{w}}}, sort=mua,
  description={skinning measures on the strong stable or strong
    unstable leaf $\wssu(w)$}}

\newglossaryentry{hausdorffmeasure}{type=symbols,
  name={\ensuremath{(\mu_x^{\rm Haus})_{x\in X}}}, sort=muhaus,
  description={Hausdorff measures of the visual distances $d_x$
  on $\Lambda\Ga$}}

\newglossaryentry{normideal}{type=symbols, name={\ensuremath{N(I)}},
  sort=ne, description={(absolute) norm of a nonzero
    ideal $I$ in a Dedekind ring}}

\newglossaryentry{neighbourhood}{type=symbols,
  name={\ensuremath{\N_\epsilon A}}, sort=neighbourho,
  description={closed $\epsilon$-neighbourhood of a subset $A$ of a
    metric space}}

\newglossaryentry{minusneighbourhood}{type=symbols,
  name={\ensuremath{\N_{-\epsilon} A}}, sort=neighbourhoo,
  description={set of points of $A$ at distance at least 
$\epsilon$ from the complement of $A$}}

\newglossaryentry{homeoinftystableaf}{type=symbols,
  name={\ensuremath{N^\pm_w}}, sort=neighbourhood,
  description={homeomorphism between stable/unstable leaves and
    inner/outer normal bundles of horoballs}}

\newglossaryentry{norm}{type=symbols,
  name={\ensuremath{\n(\beta)}}, sort=norm,
  description={relative norm of quadratic irrational $\beta$}}

\newglossaryentry{normquaternion}{type=symbols,
  name={\ensuremath{\redn(x)}}, sort=normq,
  description={reduced norm  of a quaternion $x$}}

\newglossaryentry{conditionalweak}{type=symbols,
  name={\ensuremath{\nu^\mp_w}}, sort=nu, description={conditional
    measure on the (weak) stable/unstable leaf $\wosu(w)$ of $w\in\G_\pm
    X$}}

\newglossaryentry{origin}{type=symbols, name={\ensuremath{o}},
  sort=o, description={initial vertex map $E\XX\ra V\XX$ in a
  graph $\XX$}}

\newglossaryentry{shadow}{type=symbols, name={\ensuremath{\OOO_xA}},
  sort=ombre, description={shadow of a subset $A$ of $X$ seen from
    $x\in X\cup\partial_\infty X$}}

\newglossaryentry{valuationring}{type=symbols,
  name={\ensuremath{\OOO_v}}, sort=oo, description={valuation ring
    of $v$ in $K_v$}}

\newglossaryentry{recurrentset}{type=symbols,
  name={\ensuremath{\Omega_{\rm c}}}, sort=ooh, description={two-sided recurrent
  set for the geodesic flow in $\Ga\bs\G X$}}

\newglossaryentry{footpointprojection}{type=symbols,
  name={\ensuremath{\pi}}, sort=pg, description={footpoint projection
    $w\mapsto w(0)$ of (generalised) geodesic lines, and of their
    germs at $0$}}

\newglossaryentry{naturalextension}{type=symbols,
  name={\ensuremath{\pi_+}}, sort=pga, description={natural extension 
from one-sided to two-sided shifts}}

\newglossaryentry{uniformizergene}{type=symbols,
  name={\ensuremath{\pi_v}}, sort=pgf, description={uniformiser 
of a field endowed with a valuation  $v$}}

\newglossaryentry{uniformizer}{type=symbols,
  name={\ensuremath{\pi_v}}, sort=pgg, description={uniformiser 
of a valuation  $v$ of a function field $K$ over $\FF_q$}}

\newglossaryentry{fundamentalgroupoid}{type=symbols,
  name={\ensuremath{\pi(\YY,G_*)}}, sort=pgy, description={fundamental 
groupoid of a graph of groups $(\YY,G_*)$}}

\newglossaryentry{totalmassPatterson}{type=symbols,
  name={\ensuremath{\wt \phi_{\mu^\pm}}}, sort=phi,
  description={total mass function of Patterson density
    $(\mu^\pm_x)_{x\in X}$}}

\newglossaryentry{closestpointmap}{type=symbols,
  name={\ensuremath{P_D}}, sort=phd,
  description={closest point map to a convex subset $D$}}

\newglossaryentry{closestpointmaptonormal}{type=symbols,
  name={\ensuremath{P^\pm_D}}, sort=phdpm, description={closest
    point map homeomorphism from $\partial_\infty X-\partial_\infty D$
    to outer/inner normal bundle of $D$}}

\newglossaryentry{PGLplus}{type=symbols,
  name={\ensuremath{\PGL_2(K_v)^+}}, sort=pplus, description={kernel 
of  morphism $[g]\mapsto v(\det g) \mod 2$ from 
$\PGL_2(K_v)$ to $\ZZ/2\ZZ$}}

\newglossaryentry{pressuredisc}{type=symbols,
  name={\ensuremath{\pressure{\phi}}}, sort=pressure,
  description={pressure of a potential $\phi$ under a transformation}}

\newglossaryentry{pressure}{type=symbols,
  name={\ensuremath{\pressure{\psi}}}, sort=pressurea,
  description={pressure of a potential $\psi$ under a flow}}

\newglossaryentry{metricpressuredisc}{type=symbols,
  name={\ensuremath{\pressure{\phi}(m)}}, sort=pressuremet,
  description={metric pressure for a potential $\phi$ of a
  probability measure $m$ invariant under a transformation}}

\newglossaryentry{metricpressure}{type=symbols,
  name={\ensuremath{\pressure{\psi}(m)}}, sort=pressuremeta,
  description={metric pressure for a potential $\psi$ of a 
flow-invariant probability measure $m$}}

\newglossaryentry{poincareseries}{type=symbols,
  name={\ensuremath{Q=Q_{\Ga,\,F,\,x,\,y}}}, sort=qg,
  description={Poincar\'e series of $(\Ga,F)$}}

\newglossaryentry{orderresidualfield}{type=symbols,
  name={\ensuremath{q_v}}, sort=qv, description={order of residual
    field $k_v$}}

\newglossaryentry{affinering}{type=symbols, name={\ensuremath{R_v}},
  sort=rv, description={affine algebra of the affine
    curve ${\bf C}-\{v\}$}}

\newglossaryentry{onesidedshift}{type=symbols,
  name={\ensuremath{\sigma_+}}, sort=sigma,
  description={one-sided shift in symbolic dynamics}}

\newglossaryentry{galconjloxo}{type=symbols,
  name={\ensuremath{\alpha^\sigma}}, sort=sigmaga,
  description={Other fixed point than $\alpha$ of a 
loxodromic element}}

\newglossaryentry{galconj}{type=symbols,
  name={\ensuremath{\beta^\sigma}}, sort=sigmakv,
  description={Galois conjugate of a quadratic
  irrational $\beta$ in $K_v$}}

\newglossaryentry{galconjqp}{type=symbols,
  name={\ensuremath{\alpha^\sigma}}, sort=sigmaqp,
  description={Galois conjugate of a quadratic
  irrational $\alpha$ in $\QQ_p$}}

\newglossaryentry{skinningpm}{type=symbols,
  name={\ensuremath{\wt\sigma^\pm_D}}, sort=skinningpm,
  description={skinning measure on outer/inner normal bundle of
    convex subset $D$}}

\newglossaryentry{skinningpmcalD}{type=symbols,
  name={\ensuremath{\wt\sigma^\pm_{\D}}}, sort=skinningpmcalD,
  description={outer/inner skinning measure on $\gengeod X$ of a
    family of closed convex subsets $\D$}}

\newglossaryentry{skinningpcalDdown}{type=symbols,
  name={\ensuremath{\sigma^+_{\D}}}, sort=skinningpmcalD,
  description={outer skinning measure on $\Ga\bs\gengeod X$ of a
    family of closed convex subsets $\D$}}

\newglossaryentry{skinningmcalDdown}{type=symbols,
  name={\ensuremath{\sigma^-_{\D}}}, sort=skinningpmcalD,
  description={inner skinning measure on $\Ga\bs\gengeod X$ of a
    family of closed convex subsets $\D$}}

\newglossaryentry{skinningpmOmega}{type=symbols,
  name={\ensuremath{\wt\sigma^\pm_{\Omega}}}, sort=skinningpmOmega,
  description={outer/inner skinning measure of a family
    $\Omega=(\Omega_i)_{i\in I}$ of subsets of $(\normalpm D_i)_{i\in
      I}$}}

\newglossaryentry{terminus}{type=symbols, name={\ensuremath{t}},
  sort=t, description={terminal vertex map $E\XX\ra V\XX$ in a
  graph $\XX$}}

\newglossaryentry{Tpi}{type=symbols,
  name={\ensuremath{T\pi}}, sort=ta,
  description={first edge map of a discrete time geodesic line}}

\newglossaryentry{unittangentspace}{type=symbols,
  name={\ensuremath{T^1X}}, sort=tangentspace,
  description={space of germs at $t=0$ of geodesic lines in $X$}}

\newglossaryentry{trace}{type=symbols,
  name={\ensuremath{\tr(\beta)}}, sort=trace,
  description={relative trace of quadratic irrational $\beta$}}

\newglossaryentry{tracequaternion}{type=symbols,
  name={\ensuremath{\redtr(x)}}, sort=traceq,
  description={reduced trace of a quaternion $x$}}

\newglossaryentry{volformgraphgroupedge}{type=symbols,
  name={\ensuremath{\Tvol_{(\YY,\,G_*)}}}, sort=Tvolform,
  description={volume form on the set of edges of a graph of
    finite groups $(\YY,\,G_*)$}}

\newglossaryentry{voltotgraphgroupedgemet}{type=symbols,
  name={\ensuremath{\Tvol_{(\YY,\,G_*,\lambda)}}}, sort=Tvolformet,
  description={volume form of the set of edges of a metric graph of
    finite groups $(\YY,\,G_*,\lambda)$}}

\newglossaryentry{voltotgraphgroupedge}{type=symbols,
  name={\ensuremath{\TVol(\YY,\,G_*)}}, sort=Tvoltot,
  description={total volume of the set of edges of a graph of
    finite groups $(\YY,\,G_*)$}}

\newglossaryentry{voltotmetgraphgroupedge}{type=symbols,
  name={\ensuremath{\TVol(\YY,\,G_*,\lambda)}}, sort=Tvoltotm,
  description={total volume of the set of edges of a metric graph of
    finite groups $(\YY,\,G_*,\lambda)$}}

\newglossaryentry{calUpmsubD}{type=symbols,
  name={\ensuremath{\U^\pm_D}}, sort=upm,
  description={domain of the fibration $f^\pm_D$}}

\newglossaryentry{dynamicalneighbourhood}{type=symbols,
  name={\ensuremath{V^\pm_{w,\,\eta,\,\eta'}}}, sort=va,
  description={dynamical neighbourhoods of a point $w\in\G_\pm  X$}}

\newglossaryentry{dynamicalneighbourhoodsubset}{type=symbols,
  name={\ensuremath{\V^\pm_{\eta,\,\eta'}(\Omega^\mp)}}, sort=vaa,
  description={dynamical neighbourhood of a subset $\Omega^\mp$ of
    $\G_\pm X$}}

\newglossaryentry{valuationatinfinity}{type=symbols,
  name={\ensuremath{v_\infty}}, sort=vell, description={valuation
    at infinity of $\FF_q(Y)$}}

\newglossaryentry{germgeodline}{type=symbols,
  name={\ensuremath{v_\ell}}, sort=vella,
  description={germ at $t=0$ of a geodesic line $\ell$}}

\newglossaryentry{volformgraphgroup}{type=symbols,
  name={\ensuremath{\vol_{(\YY,\,G_*)}}}, sort=volform,
  description={volume form on the set of vertices of a graph of
    finite groups $(\YY,\,G_*)$}}

\newglossaryentry{volumegraphgroup}{type=symbols, 
name={\ensuremath{\Vol(\YY,\,G_*)}}, sort=voltot, description={volume 
of a graph of finite groups $(\YY,\,G_*)$}}

\newglossaryentry{setofvertices}{type=symbols,
  name={\ensuremath{V\XX}}, sort=vx, description={set of vertices
  of a graph $\XX$}}
    
\newglossaryentry{setofevenvertices}{type=symbols,
  name={\ensuremath{\Veven\XX}}, sort=vxe, description={set of vertices
  of a pointed graph $(\XX,x_0)$ at even distances from $x_0$}}
  
\newglossaryentry{setofoddvertices}{type=symbols,
  name={\ensuremath{\Vodd\XX}}, sort=vxo, description={set of vertices
  of a pointed graph $(\XX,x_0)$ at odd distances from $x_0$}}

\newglossaryentry{endpointmap}{type=symbols,
  name={\ensuremath{w_\pm}}, sort=wpm,
  description={positive/negative endpoint of generalised
    geodesic line $w$}}

\newglossaryentry{Wss}{type=symbols, name={\ensuremath{W^+(w)}},
  sort=Wss, description={strong stable leaf of $w\in\G_+ X$}}

\newglossaryentry{Wsu}{type=symbols, name={\ensuremath{W^-(w)}},
  sort=Wsu, description={strong unstable leaf of $w\in\G_- X$}}

\newglossaryentry{Ws}{type=symbols, name={\ensuremath{W^{0+}(w)}},
  sort=Wssz, description={stable leaf of $w\in\G_+ X$}}

\newglossaryentry{Wu}{type=symbols, name={\ensuremath{W^{0-}(w)}},
  sort=Wsuz, description={unstable leaf of $w\in\G_- X$}}

\newglossaryentry{geomrealiz}{type=symbols,
  name={\ensuremath{|\XX|_\lambda}}, sort=X,
  description={geometric realisation of a metric tree
    $(\XX,\lambda)$}}

\newglossaryentry{BruhatTitstree}{type=symbols,
  name={\ensuremath{\XX_v}}, sort=Xv, description={Bruhat-Tits tree
    of $(\PGL_2,K_v)$}}

\newglossaryentry{zetafunctionfield}{type=symbols,
  name={\ensuremath{\zeta_K(s)}}, sort=zeta,
  description={Dedekind's zeta function of a function field $K$}}

%% file: meridonI.tex
\chapter{Introduction}
\label{sec:intro}

In this book, we study equidistribution and counting problems
concerning locally geodesic arcs in negatively curved spaces endowed
with potentials, and we deduce, from the special case of tree
quotients, various arithmetic applications to equidistribution and
counting problems in non-Archimedean local fields.

\medskip
For several decades, tools in ergodic theory and dynamical systems
have been used to obtain geometric equidistribution and counting
results on manifolds, both inspired by and with applications to
arithmetic and number theoretic problems, in particular in Diophantine
approximation. Especially pioneered by Margulis, this field has
produced a huge corpus of works, by Bowen, Cosentino, Clozel, Dani,
Einseidler, Eskin, Gorodnik, Ghosh, Guivarc'h, Kim, Kleinbock,
Kontorovich, Lindenstraus, Margulis, McMullen, Michel, Mohammadi,
Mozes, Nevo, Oh, Pollicott, Roblin, Shah, Sharp, Sullivan, Ullmo,
Weiss and the last two authors, just to mention a few contributors. We
refer for now to the surveys
\cite{Babillot02a,Oh10,ParPau16LMS,ParPauCIRM} and we will explain in
more details in this introduction the relation of our work with
previous works.

In this text, we consider geometric equidistribution and counting
problems weighted with a potential function in quotient spaces of
$\CAT(-1)$ spaces by discrete groups of isometries. The $\CAT(-1)$
spaces form a huge class of metric spaces that contains (but is not
restricted to) metric trees, hyperbolic buildings and simply connected
complete Riemannian manifolds with sectional curvature bounded above
by $-1$.  In Chapter \ref{sec:negcurv}, we review some basic
properties of these spaces and we refer to \cite{BriHae99} for more
details. Although some of the equidistribution and counting results
with potentials on negatively curved manifolds are known,\footnote{See
  for instance \cite{PauPolSha15}.} as well as some of such results on
$\CAT(-1)$ spaces without potential,\footnote{See for instance
  \cite{Roblin03}} bringing together these two aspects and producing
new results and applications is one of the goals of this book.

We extend the theory of Patterson-Sullivan, Bowen-Margulis and
skinning measures to $\CAT(-1)$ spaces with potentials, with a special
emphasis on trees endowed with a system of conductances. We prove that
under natural nondegeneracy, mixing and finiteness assumptions, the
pushforward under the geodesic flow of the skinning measure of
properly immersed locally convex closed subsets of locally $\CAT(-1)$
spaces equidistributes to the Gibbs measure, generalising the main
result of \cite{ParPau14ETDS}.

We also prove that the (appropriate generalisations of) the initial
and terminal tangent vectors of the common perpendiculars to any two
properly immersed locally convex closed subsets jointly equidistribute
to the skinning measures when the lengths of the common perpendiculars
tend to $+\infty$. This result is then used to prove asymptotic
results on weighted counting functions of common perpendiculars whose
lengths tend to $+\infty$. Common perpendiculars have been studied, in
various particular cases, sometimes not explicitly, by Basmajian,
Bridgeman, Bridgeman-Kahn, Eskin-McMullen, Herrmann, Huber,
Kontorovich-Oh, Margulis, Martin-McKee-Wambach, Meyerhoff, Mirzakhani,
Oh-Shah, Pollicott, Roblin, Shah, the last two authors and many
others.  See the comments after Theorem \ref{theo:mainintrocount}
below, and the survey \cite{ParPau16LMS} for references.

In Part \ref{sect:arithappli} of this book, we apply the geometric
results obtained for trees to deduce arithmetic applications in
non-Archimedean local fields. In particular, we prove equidistribution
and counting results for rationals and quadratic irrationals in any
completion of any function field over a finite field.

\bigskip Let us now describe more precisely the content of this book,
restricted to special cases for the sake of the exposition.

\section*{Geometric and dynamical tools}

Let $Y$ be a geodesically complete connected proper locally $\CAT(-1)$
space (or good orbispace), which is nonelementary, that is, whose
fundamental group is not virtually nilpotent. In this introduction, we
will mainly concentrate on the cases where $Y$ is either a metric
graph (or graph of finite groups in the sense of Bass and Serre, see
\cite{Serre83}) or a Riemannian manifold (or good orbifold) of
dimension at least $2$ with sectional curvature at most $-1$. Let $\G
Y$ be the space of locally geodesic lines of $Y$, on which the
geodesic flow $(\flow t)_{t\in\RR}$ acts by real translations on the
source. When $Y$ is a simplicial\footnote{that is, if its edges all
  have lengths $1$} graph (of finite groups), we consider the discrete
time geodesic flow $(\flow t)_{t\in\ZZ}$, see Section
\ref{subsec:trees}. If $Y$ is a Riemannian manifold, then $\G Y$ is
naturally identified with the unit tangent bundle $T^1Y$ by the map
that associates to a locally geodesic line its tangent vector at time
$0$.  In general, we define $T^1Y$ as the space of germs of locally
geodesic lines in $Y$, and $\G Y$ maps onto $T^1Y$ with possibly
uncountable fibers.

Let $F:T^1Y\to\RR$ be a continuous map, called a {\em potential},
which plays the same role in the construction of Gibbs
measures/equilibrium states as the energy function in Bowen's
treatment of the thermodynamic formalism of symbolic dynamical systems
in \cite[Sect.~1]{Bowen75}. We
define in Section \ref{subsec:criticexpo} the {\em critical exponent}
$\delta_F$ associated with $F$, which describes the logarithmic growth
of an orbit of the fundamental group on the universal cover of $Y$
weighted by the (lifted) potential $F$, and which coincides with the
classical critical exponent when $F=0$. When $Y$ is a metric graph, we
associate in Section \ref{subsec:cond} a potential $F_c$ to a {\it
  system of conductances}\index{system of conductances} $c$ (that is,
a map from the set of edges of $Y$ to $\RR$), in such a way that the
correspondence $c\mapsto F_c$ is bijective at the level of cohomology
classes, and we denote $\delta_{F_c}$ by $\delta_c$. 

\medskip
In this introduction, we assume that $F_{\,}$ is bounded and that
$\delta_F$ is finite and positive in order to simplify the statements.
  
\medskip
We say that the pair $(Y,F)$ satisfies the {\em
  \ref{eq:HC}-property}\index{property HC} if the integral of $F$ on
compact locally geodesic segments of $Y$ varies in a
H\"older-continuous way on its extremities (see Definition
\ref{defi:HCproperty}). The pairs which have the \ref{eq:HC}-property
include complete Riemannian manifolds with pinched sectional curvature
at most $-1$ and H\"older-continuous potentials, and metric graphs
with any potential. This \ref{eq:HC}-property is the new technical
idea compared to \cite{PauPolSha15} which allows the extensions to our
very general framework. See also \cite{ConLafTho16}, under the very
strong assumption that $Y$ is compact.

In Chapter \ref{sect:measures}, building on the works of
\cite{Roblin03}\footnote{itself building on the works of Patterson,
  Sullivan, Coornaert, Burger-Mozes, ...} when $F=0$ and of
\cite{PauPolSha15}\footnote{itself building on the works of Ledrappier
  \cite{Ledrappier95b}, Hamenstädt, Coudène, Mohsen} when $Y$ is a
Riemannian manifold, we generalise, to locally $\CAT(-1)$ spaces $Y$
endowed with a potential $F$ satisfying the \ref{eq:HC}-property, the
construction and basic properties of the {\em Patterson densities} at
infinity of the universal cover of $Y$ associated with $F$ and the
{\em Gibbs measure} $m_F$ on $\G Y$ associated with $F$.

Using the Patterson-Sullivan-Bowen-Margulis approach, the Patterson
densities are limits of renormalised measures on the orbit points of
the fundamental group on the universal cover of $Y$, weighted by the
potential, and the Gibbs measures on $\G Y$ are local products of
Patterson densities on the endpoints of the geodesic line, with the
Lebesgue measure on the time parameter, weighted by the Gibbs cocycle
defined by the potential.

Generalising a result of \cite{CooPap96}, we prove in Section
\ref{subsec:harmonicmeasure} that when $Y$ is a regular simplicial
graph and $c$ is an antireversible system of conductances, then the
Patterson measures, normalised to be probability measures, are
harmonic measures (or hitting measures) on the boundary at infinity of
the universal cover of $Y$ for a transient random walk on the
vertices, whose transition probabilities are constructed using the
total mass of the Patterson measures.

Gibbs measures were first introduced in statistical mechanics, and are
naturally associated via the thermodynamic formalism%
\index{thermodynamic formalism}\footnote{See for instance
  \cite{Ruelle04,Keller98,Zinsmeister96}.} with symbolic dynamics. We
prove in Section \ref{subsec:gibbsmeasure} that our Gibbs measures
satisfy a Gibbs property analogous to the one in symbolic dynamics. If
$F=0$, the Gibbs measure $m_F$ is the {\em Bowen-Margulis measure}
$m_{\rm BM}$. If $Y$ is a compact Riemannian manifold and $F$ is the
{\em strong unstable Jacobian} $v\mapsto-\,\frac{d}{dt}_{\mid t= 0}
\ln\operatorname{Jac} \big({\flow t}_{\mid \wsu(v)}\big)(v)$, then
$m_F$ is the Liouville measure and $\delta_F=0$ (see
\cite[Chap.~7]{PauPolSha15} for more general assumptions on $Y$).
Thus, one interesting aspect of Gibbs measure is that they form a
natural family of measures invariant under the geodesic flow that
interpolates between the Liouville measure and the Bowen-Margulis
measure (which in variable curvature are in general not in the same
measure class).  Another interesting point is that such measures are
plentiful: a recent result of Belarif \cite{Belarif16} proves that
when $Y$ is a geometrically finite Riemannian manifold with pinched
negative curvature and topologically mixing geodesic flow, the finite
and mixing Gibbs measures associated with bounded H\"older-continuous
potentials are, once normalised, dense (for the weak-star topology) in
the whole space of probability measures invariant under the geodesic
flow.

The Gibbs measures are remarkable measures for $\CAT(-1)$ spaces
endowed with potentials due to their unique ergodic-theoretic
properties. Let $(Z,(\phi_t)_{t\in\RR})$ be a topological space
endowed with a continuous one-parameter group of homeomorphisms and
let $\psi:Z\ra\RR$ be a bounded continuous map. Let $\M$ be the set of
Borel probability measures on $Z$ invariant under the flow
$(\phi_t)_{t\in\RR}$.  Let $h_m(\phi^1)$ be the (metric) entropy of
the geodesic flow with respect to $m\in\M$.  The {\em metric
  pressure}\index{metric!pressure}\index{pressure!metric} for $\psi$
of a measure $m\in\M$ and the {\em pressure}\index{pressure} of $\psi$
are respectively
$$
\gls{metricpressure}=h_m(\phi_1)+\int_{Z}\,\psi\;dm
\;\;\;\;\;{\rm and}\;\;\;\;\;
\gls{pressure}=\sup_{m\in\M}\;\pressure{\psi}(m)\;.
$$ 
An element $m\in\M$ is an {\em equilibrium state}%
\index{equilibrium state} for $\psi$ if the least upper bound defining
$\pressure{\psi}$ is attained on $m$.

Let $F^\sharp:\G Y\ra\RR$ be the composition of the canonical map $\G
Y\ra T^1Y$ with $F$, and note that $F^\sharp=F$ if $Y$ is a Riemannian
manifold. When $F=0$ and $Y$ is a Riemannian manifold, whose sectional
curvatures and their first derivatives are bounded, by
\cite[Theo.~2]{OtaPei04}, the pressure $\pressure{F}$ coincides with
the entropy of the geodesic flow, it is equal to the critical exponent
of the fundamental group of $Y$, and the Bowen-Margulis measure $m_F=
m_{\rm BM}$, normalised to be a probability measure, is the measure of
maximal entropy.  When $Y$ is a Riemannian manifold whose sectional
curvatures and their first derivatives are bounded and $F$ is
H\"older-continuous, by \cite[Theo.~6.1]{PauPolSha15}, we have
$\pressure{F}= \delta_F$. If furthermore the Gibbs measure $m_F$ is
finite and normalised to be a probability measure, then $m_F$ is an
equilibrium state for $F$.

In Section \ref{subsec:varprinc}, we prove an analog of these results
for the potential $F^\sharp$ when $Y$ is a metric graph of groups. The
case when $Y$ is a finite simplicial graph\footnote{that is, a finite
  graph of trivial groups with edge lengths $1$} is classical by the
work of Bowen \cite{Bowen75}, as it reduces to arguments of subshifts
of finite type (see for instance \cite{CooPap93}). When $Y$ is a
compact\footnote{a very strong assumption that we do not want to make
  in this text} locally $\CAT(-1)$-space,\footnote{not in the orbifold
  sense, hence this excludes for instance the case of graphs of groups
  with some nontrivial vertex stabiliser} a complete statement about
existence, uniqueness and Gibbs property of equilibrium states for any
H\"older-continuous potential is given in \cite{ConLafTho16}.

\btheo[The variational principle for metric graphs of
  groups] \label{theo:varprintreeintro} Assume that $Y$ is a metric
graph of finite groups, with a positive lower bound and finite upper
bound on the lengths of edges. If the critical exponent $\delta_F$ is
finite, if the Gibbs measure $m_F$ is finite, then
$\pressure{F^\sharp} = \delta_F$ and the Gibbs measure normalised to
be a probability measure is the unique equilibrium state for
$F^\sharp$.  
\etheo

The main tool is a natural coding of the discrete time geodesic flow
by a topological Markov shift (see Section \ref{subsec:TMS}). This
coding is delicate when the vertex stabilisers are nontrivial, in
particular as it does not satisfy in general the Markovian property of
dependence only on the immediate past (see Section
\ref{subsec:codagesimplicial}). We then apply results of Buzzi and
Sarig in symbolic dynamics over a countable alphabet (see Appendix
\ref{appendixBuzzi} written by J.~Buzzi), and suspension techniques
introduced in Section \ref{subsec:codagemetric}. See also
\cite{Kempton11}.

\medskip
Let $Y$ be any geodesically complete connected proper locally
$\CAT(-1)$ space, and let $D$ be any connected proper nonempty {\em
  properly immersed}\index{properly immersed}\footnote{By definition,
  $D$ is the image in $Y$, by the universal covering map, of a proper
  nonempty closed convex subset of the universal cover of $Y$, whose
  family of images under the universal covering group is locally
  finite.} closed locally convex subset of $Y$. In Chapter
\ref{sect:skinning}, we generalise for nonconstant potentials on $Y$
the construction of the {\em skinning measures} $\sigma^+_{D}$ and
$\sigma^{-}_{D}$ on the outer and inner unit normal bundles of $D$ in
$Y$.  We refer to Section \ref{subsect:nbhd} for the appropriate
definition of the outer and inner unit normal bundles of $D$ when the
boundary of $D$ is not smooth.  We construct these measures
$\sigma^+_{D}$ and $\sigma^{-}_{D}$ as pushforwards of the Patterson
densities associated with the potential $F$ to the outer and inner
unit normal bundles of the lift of $D$ in the universal cover of
$Y$. This construction generalises the one in \cite{ParPau14ETDS} when
$F=0$, which itself generalises the one in \cite{OhSha12,OhSha13} when
$M$ has constant curvature and $D$ is a ball, a horoball or a totally
geodesic submanifold.

In Section \ref{subsec:skinningwithpot}, we prove the following result
on the equidistribution of equidistant hypersurfaces in $\CAT(-1)$
spaces.  This result is a generalisation of
\cite[Theo.~1]{ParPau14ETDS} (valid in Riemannian manifolds with zero
potential) which itself generalised the ones in
\cite{Margulis04a,EskMcMul93,ParPau12JMD} when $Y$ has
constant curvature, $F=0$ and $D$ is a ball, a horoball or a totally
geodesic submanifold. See also \cite{Roblin03} when $Y$ is a $\CAT(-1)$
space, $F=0$ and $D$ is a ball or a horoball.

\btheo \label{theo:introequidequid} Let $Y,D$ be as above, and let $F$
be a potential of $Y$ satisfying the \ref{eq:HC}-property.  Assume
that the Gibbs measure $m_{F}$ on $\G Y$ is finite and mixing for the
geodesic flow $(\flow t)_{t\in\RR}$, and that the skinning measure
$\sigma^+_{D}$ is finite and nonzero. Then, as $t$ tends to $+\infty$,
the pushforwards $(\flow t)_*\sigma^+_{D}$ of the skinning measure of
$D$ by the geodesic flow weak-star converge towards the Gibbs measure
$m_F$ (after normalisation as probability measures).  \etheo

We prove in Theorem \ref{theo:equidsimplicial} an analog of Theorem
\ref{theo:introequidequid} for the discrete time geodesic flow on
simplicial graphs and, more generally, simplicial graphs of groups.
As a special case, we recover known results on nonbacktracking simple
random walks on regular graphs.  The equidistribution of the
pushforward of the skinning measure of a subgraph is a weighted
version of the following classical result, see for instance
\cite{AloBenLubSod07}, which under further assumptions on the spectral
properties on the graph gives precise rates of convergence.

\bcoro \label{coro:introrandomwalk}
Let $\YY$ be a finite regular graph which is not bipartite.  Let $\YY'$
be a nonempty connected subgraph. Then the $n$-th vertex of the
nonbacktracking simple random walk on $\YY$ starting transversally to
$\YY'$ converges in distribution to the uniform distribution as $n\ra
+\infty$.  
\ecoro

See Chapter \ref{sec:skinningwithpot} for more details and for the
extensions to nonzero potential and to graphs of groups, as well as
Section \ref{subsec:rateequidtrees} for error terms.

\section*{The distribution of common perpendiculars}

Let $D^-$ and $D^+$ be connected proper nonempty properly immersed
locally convex closed subsets of $Y$. A {\em common perpendicular}
from $D^-$ to $D^+$ is a locally geodesic path in $Y$ starting
perpendicularly from $D^-$ and arriving perpendicularly to
$D^+$.\footnote{See Section \ref{subsec:creatcommonperp} for
  explanations when the boundary of $D^-$ or $D^+$ is not smooth.}  We
denote the length of a common perpendicular $\alpha$ from $D^-$ to
$D^+$ by $\len(\alpha)$, and its initial and terminal unit tangent
vectors by $v^-_\alpha$ and $v^+_\alpha$. In the general $\CAT(-1)$
case, $v^\pm_\alpha$ are two different parametrisations (by
$\mp[0,\lambda(\alpha)]$) of $\alpha$, considered as elements of the
space $\gengeod Y$ of generalised locally geodesic lines in $Y$, see
\cite{BarLuc12} or Section \ref{subsec:lines}. For all $t> 0$, we
denote by $\Perp (D^-, D^+, t)$ the set of common perpendiculars from
$D^-$ to $D^+$ with length at most $t$ (considered with
multiplicities), and we define the counting function with weights by
$$
\N_{D^-,\,D^+,\,F}(t) =
\sum_{\alpha\in\Perp(D^-,\,D^+,\,t)}\;e^{\int_{\alpha} F}\;,
$$ 
where $\int_{\alpha} F=\int_0^{\len(\alpha)} F(\flow t v^-_\alpha)\;
dt$.  We refer to Section \ref{subsec:multandcount} for the definition
of the multiplicities in the manifold case, which are equal to $1$ if
$D^-$ and $D^+$ are embedded and disjoint. Higher multiplicities for
common perpendiculars $\alpha$ can occur for instance when $D^-$ is a
nonsimple closed geodesic and the initial point of $\alpha$ is a
multiple point of $D^-$.

Let $\Perp(D^-,D^+)$ be the set of all common perpendiculars from
$D^-$ to $D^+$ (considered with multiplicities). The family
$(\len(\alpha))_{\alpha\in\Perp(D^-,\,D^+)}$ is called the {\em marked
  ortho\-length spectrum}\index{ortholength spectrum!marked} from
$D^-$ to $D^+$. The set of lengths (with multiplicities) of elements
of $\Perp(D^-,D^+)$ is called the {\em ortholength
  spectrum}\index{ortholength spectrum} of $D^-,D^+$. This second set
has been introduced by Basmajian \cite{Basmajian93} (under the name
``full orthogonal spectrum'') when $M$ has constant curvature, and
$D^-$ and $D^+$ are disjoint or equal embedded totally geodesic
hypersurfaces or embedded horospherical cusp neighbourhoods or
embedded balls. We refer to the paper \cite{BriKah10} which proves
that the ortholength spectrum with $D^-=D^+=\partial M$ determines the
volume of a compact hyperbolic manifold $M$ with totally geodesic
boundary (see also \cite{Calegari11} and \cite{MasMcS13}).

We prove in Chapter \ref{sec:equidcountdownstairs} that the critical
exponent $\delta_F$ of $F$ is the exponential growth rate of
$\N_{D^-,\,D^+,\,F} (t)$, and we give an asymptotic formula of the
form $\N_{D^-,\,D^+,\,F}(t) \sim c\,e^{\delta_Ft}$ as $t\ra+\infty$,
with error term estimates in appropriate situations. The constants $c$
that will appear in such asymptotic formulas will be explicit, in
terms of the measures naturally associated with the (normalised)
potential $F\,$: the Gibbs measure $m_F$ and the skinning measures of
$D^-$ and $D^+$.

When $F=0$ and $Y$ is a Riemannian manifold with pinched sectional
curvature and finite and mixing Bowen-Margulis measure, the
asymptotics of the counting function $\N_{D^-,\,D^+,\,0}(t)$ are
described in \cite[Theo.~1]{ParPau16ETDS}.  The only restriction on the
two convex sets $D^\pm$ is that their skinning measures are finite.
Here, we generalise that result by allowing for nonzero potential and
more general $\CAT(-1)$ spaces than just manifolds.

The counting function $\N_{D^-,\,D^+,\,0}(t)$ has been studied in
negatively curved manifolds since the 1950's and in a number of more
recent works, sometimes in a different guise. A number of special
cases (all with $F=0$ and covered by the results of
\cite{ParPau16ETDS}) were known:
\begin{itemize}
\item $D^-$ and $D^+$ are reduced to points, by for instance \cite{Huber59},
  \cite{Margulis69} and \cite{Roblin03},
\item $D^-$ and $D^+$ are horoballs, by \cite{BelHerPau01},
  \cite{HerPau04}, \cite{Cosentino99} and \cite{Roblin03} without an
  explicit form of the constant in the asymptotic expression,
\item $D^-$ is a point and $D^+$ is a totally geodesic submanifold, by
  \cite{Herrmann62}, \cite{EskMcMul93} and \cite{OhSha16} in
  constant curvature,
\item $D^-$ is a point and $D^+$ is a horoball, by
  \cite{Kontorovich09} and \cite{KonOh11} in constant curvature, and
  \cite{Kim15} in rank one symmetric spaces,
\item $D^-$ is a horoball and $D^+$ is a totally geodesic submanifold,
  by \cite{OhSha12} and \cite{ParPau12JMD} in constant curvature, and
\item $D^-$ and $D^+$ are (properly immersed) locally geodesic lines
  in constant curvature and dimension $3$, by \cite{Pollicott11}.
\end{itemize}
We refer to the survey \cite{ParPau16LMS} for more details on the
manifold case.
  
When $X$ is a compact metric or simplicial graph and $D^\pm$ are
points, the asymptotics of $\N_{D^-,\,D^+,\,0}(t)$ as $t\ra+\infty$ is
treated in \cite{Guillope94}, as well as \cite{Roblin03}.  Under the
same setting, see also the work of Kiro-Smilansky-Smilansky announced
in \cite{KirSmiSmi16} for a counting result of paths (not assumed to
be locally geodesic) in finite metric graphs with rationally
independent edge lengths and vanishing potential.

\medskip The proofs of the asymptotic results on the counting function
$\N_{D^-,\,D^+,\,F}$ are based on the following simultaneous
equidistribution result that shows that the initial and terminal
tangent vectors of the common perpendiculars equidistribute to the
skinning measures of $D^-$ and $D^+$.  We denote the unit Dirac mass
at a point $z$ by $\Dirac_z$ and the total mass of any measure $m$ by
$\|m\|$.

\btheo \label{theo:mainintroequidis} Assume that $Y$ is a nonelementary
Riemannian manifold with pinched sectional curvature at most $-1$ or a
metric graph. Let $F:T^1Y\ra \RR$ be a potential, with finite and
positive critical exponent $\delta_F$, which is bounded and
H\"older-continuous when $Y$ is a manifold. Let $D^\pm$ be as
above. Assume that the Gibbs measure $m_F$ is finite and mixing for
the geodesic flow.  For the weak-star convergence of measures on
$\gengeod Y\times \gengeod Y$, we have
$$
\lim_{t\ra+\infty}\; \delta_F\;\|m_F\|\;e^{-\delta_F t}
\sum_{\alpha\in\Perp(D^-,\,D^+,\,t)} \;e^{\int_{\alpha} F}\; 
\Dirac_{v^-_\alpha} \otimes\Dirac_{v^+_\alpha}\;=\;
\sigma^+_{D^-}\otimes \sigma^-_{D^+}\,.
$$
\etheo

There is a similar statement for nonbipartite simplicial graphs and
for more general graphs of groups on which the discrete time geodesic
flow is mixing for the Gibbs measure, see the end of Chapter
\ref{sec:equidarcs} and Section
\ref{subsect:equicountmetricgraphgroup}. Again, the results can then be
interpreted in terms of nonbacktracking random walks.

In Section \ref{subsec:downstairs}, we deduce our counting results for
common perpendiculars between the subsets $D^-$ and $D^+$ from the above
simultaneous equidistribution theorem.

\btheo\label{theo:mainintrocount}
\smallskip\noindent (1) Let $Y,F,D^\pm$ be as in Theorem
\ref{theo:mainintroequidis}.  Assume that the Gibbs measure $m_F$ is
finite and mixing for the continuous time geodesic flow and that the
skinning measures $\sigma^+_{D^-}$ and $\sigma^{-}_{D^+}$ are finite
and nonzero. Then, as $s\to+\infty$,
$$
\N_{D^-,\,D^+,\,F}(s)\sim
\frac{\|\sigma^+_{D^-}\|\,\|\sigma^{-}_{D^+}\|}
{\|m_F\|}\,\frac{e^{\delta_{F} \,s}}{\delta_{F}}\,.
$$

\noindent (2) If $Y$ is a finite nonbipartite simplicial graph, then
$$
\N_{D^-,\,D^+,F}(n)\sim 
\frac{e^{\delta_{F}}\;\|\sigma^+_{D^-}\|\;\|\sigma^-_{D^+}\|}
{(e^{\delta_{F}}-1)\;\|m_{F}\|}\;\;e^{\delta_{F}\, n}\,.
$$
\etheo

The above Assertion (1) is valid when $Y$ is a good orbifold instead
of a manifold or a metric graph of finite groups instead of a metric
graph (for the appropriate notion of multiplicities), and when $D^-$
and $D^+$ are replaced by locally finite families. See Section
\ref{subsect:equicountmetricgraphgroup} for generalisations of
Assertion (2) to (possibly infinite) simplicial graphs of finite
groups and Sections \ref{subsect:erroterms} and
\ref{subsect:errormetricgraphgroup} for error terms.

We avoid any compactness assumption on $Y$, we only assume that the
Gibbs measure $m_F$ of $F$ is finite and that it is mixing for the
geodesic flow.  By Babillot's theorem \cite{Babillot02b}, if the
length spectrum of $Y$ is not contained in a discrete subgroup of
$\RR$, then $m_F$ is mixing if finite. If $Y$ is a Riemannian
manifold, this condition is satisfied for instance if the limit set of
a fundamental group of $Y$ is not totally disconnected, see for
instance \cite{Dalbo99,Dalbo00}. When $Y$ is a metric graph,
Babillot's mixing condition is in particular satisfied if the lengths
of the edges of $Y$ are rationally independent.

As in \cite{ParPau16ETDS}, we have very weak finiteness and curvature
assumptions on the space and the convex subsets we consider.
Furthermore, the space $Y$ is no longer required to be a manifold and
we extend the theory to nonconstant weights using equilibrium
states. Such a weighted counting has only been used in the
orbit-counting problem in manifolds with pinched negative curvature in
\cite{PauPolSha15}. The approach is based on ideas from Margulis's
thesis to use the mixing of the geodesic flow. Our skinning measures
are much more general than the Patterson measures appearing in earlier
works. As in \cite{ParPau16ETDS}, we push simultaneously the unit
normal vectors to the two convex sets $D^-$ and $D^+$ in opposite
directions.

\medskip
Classically, an important characterisation of the Bowen-Margulis
measure on closed negatively curved Riemannian manifolds ($F=0$) is
that it coincides with the weak-star limit of properly normalised sums
of Lebesgue measures supported on periodic orbits. The result was
extended to $\CAT(-1)$ spaces with zero potential in \cite{Roblin03}
and to Gibbs measures in the manifold case in
\cite[Theo.~9.11]{PauPolSha15}. As a corollary of the simultaneous
equidistribution result Theorem \ref{theo:mainintroequidis}, we obtain
a weighted version for simplicial and metric graphs of groups. The
following is a simplified version of such a result for Gibbs measures
of metric graphs.

Let $\periodic'(t)$ be the set of prime periodic orbits of the
geodesic flow on $Y$. Let $\lambda(g)$ denote the length of a closed
orbit $g$. Let $\Lebmeas_g$ be the Lebesgue measure along $g$ and let
$\Lebmeas_g(F)$ be the period of $g$ for the potential $F$.

\btheo\label{theo:introequidclosedorbitscont} Assume that $Y$ is a
finite metric graph, that the critical exponent $\delta_F$ is positive
and that the Gibbs measure $m_{F}$ is mixing for the continuous time
geodesic flow.  As $t\to+\infty$, the measures
$$
\delta_{F}\,e^{\delta_{F}t}\sum_{g\in\periodic'(t)}e^{\Lebmeas_g(F)}\,\Lebmeas_g
$$
and 
$$
\delta_{F}\,t\,e^{\delta_{F}t}\sum_{g\in\periodic'(t)}
e^{\Lebmeas_g(F)}\,\frac{\Lebmeas_g}{\lambda(g)}
$$ 
converge to $\frac{m_{F}}{\|m_{F}\|}$ for the weak-star convergence
of measures.  
\etheo

See Section \ref{sec:equidclosedorbits} for the proof of the full
result and for a similar statement for (possibly infinite) simplicial
graphs of finite groups. As a corollary, we obtain counting results of
simple loops in metric and simplicial graphs, generalising results of
\cite{ParPol90}, \cite{Guillope94}.

\bcoro
\smallskip\noindent Assume that $Y$ is a finite metric graph whose
vertices have degrees at least $3$, such that the critical exponent
$\delta_F$ is positive.
\begin{enumerate}
\item
If the Gibbs measure is mixing for the continuous time geodesic flow,
then
$$
\sum_{g\in\periodic'(t)}e^{\Lebmeas_g(F)}\sim \frac{e^{\delta_F\, t}}{\delta_F\, t}
$$
as $t\to+\infty$.
\item
If $Y$ is simplicial and if the Gibbs measure is mixing for the
discrete time geodesic flow, then
$$ 
\sum_{g\in\periodic'(t)}e^{\Lebmeas_g(F)}\sim
  \frac{e^{\delta_F}}{e^{\delta_F}-1}\,\frac{e^{\delta_F\, t}}{t}
$$
as $t\to+\infty$.
\end{enumerate}
\ecoro

\medskip

In the cases when error bounds are known for the mixing property of
the continuous time or discrete time geodesic flow on $\G Y$, we obtain
corresponding error terms in the equidistribution result of Theorem
\ref{theo:introequidequid} generalising \cite[Theo.~20]{ParPau14ETDS}
(where $F=0$) and in the approximation of the counting function
$\N_{D^-,\,D^+,\,0}$ by the expression introduced in Theorem
\ref{theo:mainintrocount}.  In the manifold case, see \cite{KleMar96},
\cite{Clozel03}, \cite{Dolgopyat98}, \cite{Stoyanov11},
\cite{Liverani04}, \cite{GiuLivPol13}, and Section
\ref{subsect:erroterms} for definitions and precise references.  Here
is an example of such a result in the manifold case.

\btheo\label{theo:4} 
Assume that $Y$ is a compact Riemannian manifold and $m_F$ is
exponentially mixing under the geodes\-ic flow for the H\"older
regularity, or that $Y$ is a locally symmetric space, the boundary of
$D^\pm$ is smooth, $m_F$ is finite, smooth, and exponentially mixing
under the geodesic flow for the Sobolev regularity. Assume that the
strong stable/unstable ball masses by the conditionals of $m_F$ are
H\"older-continuous in their radius.
\begin{enumerate}
\item 
As $t$ tends to $+\infty$, the pushforwards $(\flow
t)_*\sigma^+_{D^-}$ of the skinning measure of $D^-$ by the geodesic
flow equidistribute towards the Gibbs measure $m_F$ with exponential
speed.
\item 
  If the skinning measures $\sigma^+_{D^-}$ and $\sigma^-_{D^+}$
  are finite and nonzero, there exists $\kappa>0$ such that, as
  $t\ra+\infty$,
$$
\N_{D^-,\,D^+,\, F}(t)=
\frac{\|\sigma^+_{D^-}\|\;\|\sigma^-_{D^+}\|}{\delta_F\;\|m_F\|}\;
e^{\delta_F t}\big(1+\operatorname{O}(e^{-\kappa t})\big)\;.
$$
\end{enumerate}
\etheo 

See Section \ref{subsect:erroterms} for a discussion of the
assumptions and the dependence of $\operatorname{O}(\cdot)$ on the
data. Similar (sometimes more precise) error estimates were known
earlier for the counting function in special cases of $D^\pm$ in
constant curvature geometrically finite manifolds (often in small
dimension) through the work of Huber, Selberg, Patterson, Lax and
Phillips \cite{LaxPhi82}, Cosentino \cite{Cosentino99}, Kontorovich
and Oh \cite{KonOh11}, Lee and Oh \cite{LeeOh13}.

When $Y$ is a finite volume hyperbolic manifold and the potential $F$
is constant $0$, the Gibbs measure is proportional to the Liouville
measure and the skinning measures of totally geodesic submanifolds,
balls and horoballs are proportional to the induced Riemannian
measures of the unit normal bundles of their boundaries. In this
situation, there are very explicit forms of the counting results in
finite-volume hyperbolic manifolds, see \cite[Cor.21]{ParPau16ETDS},
\cite{ParPau16LMS}.  These results are extended to complex hyperbolic
space in \cite{ParPau16MA}.

As an example of this result, if $D^-$ and $D^+$ are closed geodesics
of $Y$ of lengths $\ell_-$ and $\ell_+$, respectively, then the number
$\N(s)=\N_{D^-,\,D^+,\,0}(s)$ of common perpendiculars (counted with
multiplicity) from $D^-$ to $D^+$ of length at most $s$ satisfies, as
$s\ra+\infty$,
\begin{equation}\label{eqcaseclosedgeod}
\N(s)\sim
\frac{\pi^{\frac{n}{2}-1}\Ga(\frac{n-1}{2})^2}
{2^{n-2}(n-1)\Ga(\frac{n}{2})}
\;\frac{\ell_-\ell_+}{\Vol(Y)}\;
e^{(n-1)s}\;.
\end{equation}

\section*{Counting in weighted graphs of groups}

From now on in this introduction, we only consider metric or
simplicial graphs or graphs of groups.

Let $\YY$ be a connected finite graph with set of vertices $V\YY$ and
set of edges $E\YY$ (see \cite{Serre83} for the conventions). We
assume that the degree of the graph at each vertex is at least
$3$. Let $\len:E\YY\ra\mathopen{]}0,+\infty\mathclose{[}$ with 
$\len(\ov{e})=\len(e)$ for every $e\in E\YY$ be an {\em edge length
      map}, let $Y=|\YY|_\lambda$ be the geometric realisation of
$\YY$ where the geometric realisation of every edge $e\in E\YY$ has
length $\len(e)$, and let $c:E\YY\ra \RR$ be a map, called a
(logarithmic) {\em system of conductances} in the analogy between
graphs and electrical networks, see for instance \cite{Zemanian91}.

Let $\YY^\pm$ be proper nonempty subgraphs of $\YY$.  For every $t\geq
0$, we denote by $\Perp(\YY^-,\YY^+,t)$ the set of edge paths
$\alpha=(e_1,\dots, e_k)$ in $\YY$ without backtracking, of length
$\len(\alpha)= \sum_{i=1}^k \len(e_i)$ at most $t$, of conductance
$c(\alpha)= \sum_{i=1}^k c(e_i)$, starting from a vertex of $\YY^-$
but not by an edge of $\YY^-$, ending at a vertex of $\YY^+$ but not
by an edge of $\YY^+$. Let
$$
\N_{\YY^-,\YY^+}(t)=\sum_{\alpha\in\Perp(\YY^-,\YY^+,\,t)} e^{c(\alpha)}
$$ 
be the number of paths without backtracking from $\YY^-$ to $\YY^+$
of length at most $t$, counted with weights defined by the system of
conductances.

Recall that a real number $x$ is {\em Diophantine}\index{Diophantine}
if it is badly approximable by rational numbers, that is, if there
exist $\alpha,\beta>0$ such that $|x-\frac{p}{q}|\geq \alpha\,
q^{-\beta}$ for all $p,q\in\ZZ$ with $q> 0$.  We obtain the
following result, which is a very simplified version of our results for
the sake of this introduction.

\btheo \label{theo:countintrograph} (1) If $Y$ has two cycles whose
ratio of lengths is Diophantine, then there exists $C>0$ such that for
every $k\in\NN-\{0\}$, as $t\ra+\infty$,
$$
\N_{\YY^-,\YY^+}(t)= C\;e^{\delta_c \,t}\big(1+\bigO(t^{-k})\big)\;.
$$

\smallskip\noindent (2) If $\len\equiv 1$, then there exist
$C',\kappa>0$ such that, as $n\in\NN$ tends to $+\infty$, 
$$
\N_{\YY^-,\YY^+}(n)=C'\;e^{\delta_c  \,n}\big(1 +\bigO(e^{-\kappa \,n})\big)\;.
$$ 
\etheo

Note that the Diophantine assumption on $Y$ in Theorem
\ref{theo:countintrograph} (1) is standard in the theory of quantum
graphs (see for instance \cite{BerKuc13}).

The constants $C= C_{\YY^\pm,\,c,\,\len} >0$ and $C'=C'_{\YY^\pm,\,c}
> 0$ in the above asymptotic formulas are explicit. When $c\equiv 0$
and $\len\equiv 1$, the constants can often be determined concretely,
as indicated in the two examples below.\footnote{See Section
  \ref{subsect:equicountmetricgraphgroup} for more examples.} Among
the ingredients in these computations are the explicit expressions of
the total mass of many Bowen-Margulis measures and skinning measures
obtained in Chapter \ref{sect:explicitcomputdiscret}.

See Sections \ref{subsect:equicountmetricgraphgroup},
\ref{subsec:graphgroupsequidcount} and
\ref{subsect:errormetricgraphgroup} for generalisations of Theorem
\ref{theo:countintrograph} when the graphs $\YY^\pm$ are not embedded
in $\YY$, and for versions in (possibly infinite) metric graphs of
finite groups. In particular, Assertion (2) remains valid if $Y$ is
the quotient of a uniform simplicial tree by a geometrically finite
lattice in the sense of \cite{Paulin04b}, such as an arithmetic
lattice in $\PGL_2$ over a non-Archimedian local field, see
\cite{Lubotzky91}. Recall that a locally finite metric tree $X$ is
{\em uniform} if it admits a discrete and cocompact group of
isometries, and that a {\it lattice} $\Ga$ of $X$ is a lattice in the
locally compact group of isometries of $X$ preserving without edge
inversions the simplicial structure. We refer for instance to
\cite{BasKul90, BasLub01} for uncountably many examples of tree
lattices.

\bexem (1) When $\YY$ is a $(q+1)$-regular finite graph with constant
edge length map $\len\equiv 1$ and vanishing system of conductances
$c\equiv 0$, then $\delta_c =\ln q$, and if furthermore $\YY^+$ and
$\YY^-$ are vertices, then (see Equation
\eqref{eq:numberofwalksingraph})
$$
C'=\frac{q+1}{(q-1)\,\card (V\YY)}\;.
$$
\smallskip\noindent (2) When $\YY$ is biregular of degrees $p+1$ and
$q+1$ with $p,q\geq 2$, when $\len\equiv 1$ and $c\equiv 0$, then
$\delta_c =\ln\sqrt{pq}\,$, and if furthermore the subgraphs
$\YY^\pm$ are simple cycles of lengths $L^\pm$, then (see Equation
\eqref{eq:bicyclette}) the number of common perpendiculars of even length 
at most $2N$ from $\YY^-$ to $\YY^+$ as $N\to+\infty$ is
asymptotic to
$$
\frac{(p+q)\;L^-\;L^+}
{2\,(pq-1)\;\card(E\YY)}\;(pq)^{N+1}
$$
\eexem

The main tool in order to obtain the error terms in Theorem
\ref{theo:countintrograph} and its more general versions is to study
the error terms in the mixing property of the geodesic flow.  Using
the already mentioned coding (given in Section
\ref{subsec:codagesimplicial}) of the discrete time geodesic flow by a
two-sided topological Markov shift, classical reduction to one-sided
topological Markov shift, and results of Young \cite{Young98} on the
decay of correlations for Young towers with exponentially small tails,
we in particular obtain the following simple criteria for the
exponential decay of correlation of the discrete time geodesic flow,
where we only assume $\YY$ to be locally finite (and maybe not
finite). See Theorem \ref{theo:critexpdecaysimpl} for the complete
result.

\btheo \label{theo:critexpdecaysimplintro} 
Assume that the Gibbs measure $m_F$ is finite and mixing for the
discrete time geodesic flow on $\YY$. Assume moreover that there exist
a finite subset $E$ of $V\YY$ and $C',\kappa' >0$ such that for all
$n\in\NN$, we have
$$ 
m_F\big(\{\ell\in \G\YY\;:\;\ell(0)\in E \;{\rm and}\;
\forall\;k\in\{1,\dots,n\},\;\ell(k)\notin E \}\big)\leq
C'\;e^{-\kappa' n}\;.
$$ 
Then the discrete time geodesic flow has exponential decay of
H\"older correlations for $m_F$.  
\etheo

The assumption of having exponentially small mass of geodesic lines
which have a big return time to a given finite subset of $V\YY$ is in
particular satisfied (see Section \ref{subsec:mixingratesimpgraphs})
if $\YY$ is the quotient of a uniform simplicial tree by a
geometrically finite lattice,\footnote{See for instance
  \cite{Paulin04b}.}  such as an arithmetic lattice in $\PGL_2$ over a
non-Archimedian local field, see \cite{Lubotzky91}, but also by many
other examples of $\YY$. This statement corrects the mistake in
\cite{Kwon15}, as indicated in its erratum.

These results allow to prove in Section
\ref{subsec:mixingratemetgraphs}, under Diophantine assumptions, the
rapid mixing property for the continuous time geodesic flow, that
leads to Assertion (1) of Theorem \ref{theo:countintrograph}, see
Section \ref{subsect:errormetricgraphgroup}. The proof uses suspension
techniques due to Dolgopyat \cite{Dolgopyat98ETDS} when $Y$ is a
compact metric tree, and to Melbourne \cite{Melbourne07} otherwise.

\bigskip
As a corollary of the general version of the counting result Theorem
\ref{theo:mainintrocount}, we have the following asymptotic for the
orbital counting function in conjugacy classes for groups acting on
trees.  Given $x_0\in X$ and a nontrivial conjugacy class $\KKK$ in a
discrete group $\Ga$ of isometries of $X$, we consider the counting
function
$$
N_{\KKK,\,x_0}(t)=\card\{\ga\in\KKK\;:\; d(x_0,\ga x_0)\leq t\}\;,
$$ 
introduced by Huber \cite{Huber56} when $X$ is replaced by the real
hyperbolic plane and $\Ga$ is a lattice. We refer to \cite{ParPau15MZ}
for many results on the asymptotic growth of such orbital counting
functions in conjugacy classes, when $X$ is replaced by a finitely
generated group with a word metric, or a complete simply connected
pinched negatively curved Riemannian manifold. See also
\cite{ChaPet15arxiv,ArzCasTao15,Pollicott15}.

\btheo \label{theo:genekenisonsharpintro} Let $X$ be a uniform metric
tree with vertices of degree $\geq 3$, let $\delta$ be the Hausdorff
dimension of its space of ends, let $\Ga$ be a discrete group of
isometries of $X$, let $x_0$ be a vertex of $X$ with trivial
stabiliser in $\Ga$, and let $\KKK$ be a loxodromic conjugacy class in
$\Ga$.

\smallskip\noindent (1) If the metric graph $\Ga\bs X$ is compact and
has two cycles whose ratio of lengths is Diophantine, then there
exists $C>0$ such that for every $k\in\NN-\{0\}$, as $t\ra+\infty$,
$$
\N_{\KKK,\,x_0}(t)= C\;e^{\frac{\delta}{2} \,t}\big(1+\bigO(t^{-k})\big)\;.
$$

\smallskip\noindent(2) If $X$ is simplicial and $\Ga$ is a
geometrically finite lattice of $X$, then there exist $C',\kappa>0$
such that, as $n\in\NN$ tends to $+\infty$,
$$
\N_{\KKK,\,x_0}(n)= C'\;e^{\delta\;\lfloor
  (n-\len(\ga))/2\rfloor}\big(1 +\bigO(e^{-\kappa \,n})\big)\;.
$$
\etheo

We refer to Theorem \ref{orbcountconjuggene} for a more general
version, including a version with a system of conductances in the
counting function, and when $\KKK$ is elliptic. When $\Ga\bs X$ is
compact and $\Ga$ is torsion free,\footnote{In particular, $\Ga$ then
  has the very restricted structure of a free group.} Assertion (1) of
this result is due to Kenison and Sharp \cite{KenSha15}, who proved it
using transfer operator techniques for suspensions of subshifts of
finite type. Up to strengthening the Diophantine assumption, using
work of Melbourne \cite{Melbourne07} on the decay of correlations of
suspensions of Young towers, we are able to extend Assertion (1) to
all geometrically finite lattices $\Ga$ of $X$ in Section 
\ref{sect:orbcountconjugacyclass}.

The constants $C=C_{\KKK, x_0}$ and $C'=C'_{\KKK, x_0}$ are explicit.
For instance in Assertion (2), if $X$ is the geometric realisation of
a regular simplicial tree $\XX$ of degree $q+1$, if $x_0$ is a vertex
of $\XX$, if $\KKK$ is the conjugacy class of $\ga_0$ with translation
length $\len(\ga_0)$ on $X$, if
$$
\Vol(\Ga\dbs \XX)=\sum_{[x]\in\Ga\bs V\XX}\frac{1}{|\Ga_x|}
$$ 
is the volume\footnote{See for instance \cite{BasKul90, BasLub01}.}
of the quotient graph of groups $\Ga\dbs \XX$ , then
$$
C'=\frac{\len(\ga_0)}
{[Z_\Ga(\ga_0):\ga_0^\ZZ]\;\Vol(\Ga\dbs \XX)}\;,
$$ 
where $Z_\Ga(\ga_0)$ is the centraliser of $\ga_0$ in $\Ga$.  When
furthermore $\Ga$ is torsion free, $\ga_0$ is not a proper power and
$\Ga\bs \XX$ is finite, as $\delta=\ln q$, we get that there exists
$\kappa>0$ such that
$$
\N_{\KKK,\,x_0}(n)= \frac{\lambda(\ga_0)}{\card(\Ga\bs \XX)}\;q^{\lfloor
  (n-\len(\ga_0))/2\rfloor} +\bigO(q^{(1-\kappa') \,n/2})\;
$$ 
as $n\in\NN$ tends to $+\infty$, thus recovering the result of
\cite{Douma11} who used the spectral theory of the discrete Laplacian.

\section*{Selected arithmetic applications}

We end this introduction by giving a sample of our arithmetic
applications (see Part \ref{sect:arithappli} of this book) of the
ergodic and dynamical results on the discrete time geodesic flow on
simplicial trees described in Part \ref{part:equid} of this book, as
summarized above. Our equidistribution and counting results of common
perpendiculars between subtrees indeed produce equidistribution and
counting results of rationals and quadratic irrationals in
non-Archimedean local fields. We refer to \cite{BroParPau16CRAS} for
an announcement of the results of Part \ref{sect:arithappli}, with a
presentation different from the one in this introduction.

\medskip To motivate what follows, consider $R=\ZZ$ the ring of integers,
$K=\QQ$ its field of fractions, $\wh K=\RR$ the completion of $\QQ$
for the usual Archimedean absolute value $|\cdot|$, and $\haar_{\wh
  K}$ the Lebesgue measure of $\RR$ (which is the Haar measure of the
additive group $\RR$ normalised so that $\haar_{\wh K}([0,1])=1$). 

The following equidistribution result of rationals, due to Neville
\cite{Neville49}, is a quantitative statement on the density of $K$ in
$\wh K$: For the weak-star convergence of measures on $\wh K$, as
$s\ra +\infty$, we have
$$
\lim_{s\ra+\infty}\;\frac{\pi}{6}\;s^{-2}
\sum_{p,q\in R\;:\; p R+q R =R,\;|q|\leq s}\Delta_{\frac{p}{q}}=
\haar_{\wh K}\;.
$$ 
Furthermore, there exists $\ell\in\NN$ such that for every smooth
function $\psi:\wh K\ra\CC$ with compact support, there is an error
term in the above equidistribution claim evaluated on $\psi$, of the
form $\bigO(s(\ln s)\|\psi\|_{\ell})$ where $\|\psi\|_{\ell}$ is the
Sobolev norm of $\psi$. The following counting result due to Mertens
on the asymptotic behaviour of the average of Euler's totient function
$\varphi:k\mapsto \card (R/kR)^\times$, follows from the above
equidistribution one:
$$
\sum_{k=1}^n\;\varphi(k)= \frac{3}{\pi}\;n^{2}+\bigO(n\ln n)\;.
$$ 
See \cite{ParPau14AFST} for an approach using methods similar to the
ones in this text, and for instance \cite[Theo.~330]{HarWri08} for a
more traditional proof, as well as \cite{Walfisz63} for a better error
term.

\medskip Let us now switch to a non-Archimedean setting, restricting
to positive characteristic in this introduction. See Part
\ref{sect:arithappli} for analogous applications in 
characteristic zero. 

Let $\FF_q$ be a finite field of order $q$. Let $R=\FF_q[Y]$ be the
ring of polynomials in one variable $Y$ with coefficients in $\FF_q$.
Let $K=\FF_q(Y)$ be the field of rational fractions in $Y$ with
coefficients in $\FF_q$, which is the field of fractions of $R$. Let
$\wh K=\FF_q((Y^{-1}))$ be the field of formal Laurent series in the
variable $Y^{-1}$ with coefficients in $\FF_q$, which is the
completion of $K$ for the (ultrametric) absolute value $|\frac{P}{Q}|
=q^{\deg P-\deg Q}$. Let $\OOO=\FF_q[[Y^{-1}]]$ be the ring of formal
power series in $Y^{-1}$ with coefficients in $\FF_q$, which is the
ball of centre $0$ and radius $1$ in $\wh K$ for this absolute
value. 

Note that $\wh K$ is locally compact, and we endow the additive
group $\wh K$ with the Haar measure $\haar_{\wh K}$ normalised so that
$\haar_{\wh K} (\OOO)=1$. The following results extend (with
appropriate constants) when $K$ is replaced by any function field of a
nonsingular projective curve over $\FF_q$ and $\wh K$ any completion
of $K$, see Part \ref{sect:arithappli}.

The following equidistribution result\footnote{See Theorem
  \ref{theo:Mertensfunctionfield} in Chapter \ref{sec:rattionalequid}
for a more general version.} of elements of $K$ in $\wh K$ gives an
analog of Neville's equidistribution result for function fields.  Note
that when $G=\GL_2(R)$, we have $(P,Q)\in G(1,0)$ if and only if
$\langle P,Q\rangle = R$. We denote by $H_x$ the stabiliser of any
element $x$ of any set endowed with any action of any group $H$.

\btheo\label{theo:mertensfunctionfieldintro} Let $G$ be any finite
index subgroup of $\GL_2(R)$. For the weak-star convergence of
measures on $\wh K$, we have
$$
\lim_{t\ra+\infty}\;\frac{(q+1)\;[\GL_2(R):G]}
{(q-1)\;q^2\;[\GL_2(R)_{(1,0)}:G_{(1,0)}]}\;q^{-2\,t}
\sum_{(P,Q)\in G(1,0),\;\deg Q\leq t}\Delta_{\frac{P}{Q}}\; 
=\;\haar_{\wh K}\;.
$$
\etheo

We emphazise the fact that we are not assuming $G$ to be a congruence
subgroup of $\GL_2(R)$.  This is made possible by our geometric and
ergodic methods.

The following variation of this result is more interesting when the
class number of the function field $K$ is larger than $1$ (see
Corollary \ref{coro:appliEulerFunctField} in Chapter
\ref{sec:rattionalequid}).

\btheo\label{theo:fractidealintro} Let $\mmm$ be a nonzero
fractional ideal of $R$ with norm $\redn(\mmm)$. For the weak-star
convergence of measures on $\wh K$, we have
$$
\lim_{t\ra+\infty}\;\frac{q+1}{(q-1)\;q^2}\;s^{-2}
\sum_{\substack{(x,y)\in\mmm\times\mmm\;\\
\;\redn(\mmm)^{-1}\redn(y)\leq s,
\ R x + R y= \mmm}} \Delta_{\frac xy}\;=\;\haar_{\wh K}\,.
$$
\etheo

If $\alpha\in\wh K$ is quadratic irrational over $K$,\footnote{that
  is, $\alpha$ does not belong to $K$ and satisfies a quadratic
  equation with coefficients in $K$} let $\alpha^\sigma$ be the Galois
conjugate of $\alpha$,\footnote{that is, the other root in $\wh K$ of
  the irreducible quadratic polynomial over $K$ defining $\alpha$} let
$\tr(\alpha) =\alpha+\alpha^\sigma$ and $\n(\alpha)
=\alpha\alpha^\sigma$, and let
$$
h(\alpha)=\frac{1}{|\alpha-\alpha^\sigma|}\;.
$$ 
This is an appropriate complexity for quadratic irrationals in a
given orbit by homographies under $\PGL_2(R)$. See Section
\ref{subsec:quadirratposchar} and for instance \cite[\S 6]{HerPau10}
for motivations and results. Note that although there are only
finitely many orbits by homographies of $\PGL_2(R)$ on $K$ (and
exactly one in the particular case of this introduction), there are
infinitely many orbits of $\PGL_2(R)$ in the set of quadratic
irrationals in $\wh K$ over $K$. The following result gives in
particular that any orbit of quadratic irrationals under $\PGL_2(R)$
equidistributes in $\wh K$, when the complexity tends to infinity. See
Theorem \ref{theo:caracposequidquadirr} in Section
\ref{subsec:quadirratposchar} for a more general version. We denote by
$\cdot$ the action by homographies of $\GL_2(\wh K)$ on $\PP_1(\wh K)=
\wh K\cup \{\infty= [1:0]\}$.

\btheo\label{theo:appliquadirratintro} Let $G$ be a finite index
subgroup of $\GL_2(R)$. Let $\alpha_0\in\wh K$ be a quadratic
irrational over $K$.  For the weak-star convergence of measures on
$\wh K$, we have
$$
\lim_{s\ra+\infty}\;\frac{(\ln q)\;(q+1)\;m_0\;[\GL_2(R):G]}
{2\;q^2\;(q-1)^3\;\big|\ln |\operatorname{tr} g_0|\big|}\;s^{- 1}
\sum_{\alpha\in G\cdot\alpha_0,\; h(\alpha)\leq s}\Delta_{\alpha} 
\;=\;\haar_{\wh K}\;.
$$
where $g_0\in G$ fixes $\alpha_0$ with $|\operatorname{tr} g_0|>1$,
and $m_0$ is the index of $g_0^\ZZ$ in $G_{\alpha_0}$.  
\etheo

Another equidistribution result of an orbit of quadratic irrationals
under $\PGL_2(R)$ is obtained by taking another complexity,
constructed using crossratios with a fixed quadratic irrational. We
denote by $[a,b,c,d]=\frac{(c-a)(d-b)}{(c-b)(d-a)}$ the crossratio of
four pairwise distinct elements in $\wh K$. If $\alpha,\beta\in \wh K$
are two quadratic irrationals over $K$ such that $\alpha\notin
\{\beta, \beta^\sigma\}$,\footnote{See Section
  \ref{subsec:crossratioloxofix} when this condition is not
  satisfied.} let
$$
h_\beta(\alpha)=\max\{|[\alpha,\beta,\beta^\sigma,\alpha^\sigma]|,\;
|[\alpha^\sigma,\beta,\beta^\sigma,\alpha]|\}\;,
$$
which is also an appropriate complexity when $\alpha$ varies in a
given orbit of quadratic irrationals by homographies under $\PGL_2(R)$.
See Section \ref{subsec:crossratioloxofix} and for instance \cite[\S
4]{ParPau14AFST} for motivations and results in the Archimedean case.

\btheo\label{theo:applitwoquadirratdintro} Let $G$ be a finite index
subgroup of $\GL_2(R)$. Let $\alpha_0,\beta \in\wh K$ be two quadratic
irrationals over $K$.  For the weak-star convergence of measures on
$\wh K-\{\beta,\beta^\sigma\}$, we have, with $g_0$ and $m_0$ as in
the statement of Theorem \ref{theo:appliquadirratintro},
\begin{align*}
\lim_{s\ra+\infty}\;&\frac{(\ln q)\;(q+1)\; m_0\;[\GL_2(R):G]}
{2\;q^2\;(q-1)^3\;|\beta-\beta^\sigma|\;\big|
\ln |\operatorname{tr} g_0|\big|}
\;s^{-1} \sum_{\alpha\in G\cdot \alpha_0,\; h_{\beta}(\alpha)\leq s}
\Delta_{\alpha} \\&\;=\;
\frac{d\haar_{\wh K}(z)}{|z-\beta|\,|z-\beta^\sigma|}\;.
\end{align*}
\etheo

The fact that the measure towards which we have an equidistribution is
only absolutely continuous with respect to the Haar measure is
explained by the invariance of $\alpha\mapsto h_\beta(\alpha)$ under
the (infinite) stabiliser of $\beta$ in $\PGL_2(R)$.  See Theorem
\ref{theo:mainrelatheight} in Section \ref{subsec:crossratioloxofix}
for a more general version.

\medskip
The last statement of this introduction is an equidistribution result
for the integral representations of quadratic norm forms 
$$
(x,y)\mapsto \n(x-y\alpha)
$$ 
on $K\times K$, where $\alpha\in \wh K$ is a quadratic irrational
over $K$. See Theorem \ref{theo:mainnormform} in Section
\ref{subsec:integralrepresentation} for a more general version, and
for instance \cite[\S 5.3]{ParPau14AFST} for motivations and results
in the Archimedean case.

There is an extensive bibliography on the integral representation of
norm forms and more generally decomposable forms over function fields,
see for instance \cite{Schmidt78,Mason81,Gyory83,Mason86}.  These
references mostly consider higher degrees, with an algebraically
closed ground field of characteristic $0$, instead of $\FF_q$.

\btheo \label{theo:normformintro} Let $G$ be a finite index
subgroup of $\GL_2(R)$ and let $\beta\in\wh K$ be a quadratic
irrational over $K$. For the weak-star convergence of measures on
$\wh K-\{\beta,\beta^\sigma\}$, we have
\begin{align*}
\lim_{s\ra+\infty}\; &\frac{(q+1)\;[\GL_2(R_v):G]}
{q^2\;(q-1)^3\;[\GL_2(R_v)_{(1,0)}:G_{(1,0)}]}
\;s^{-1}\sum_{\substack{(x,y)\in G(1,0),\;\\|x^2-xy\tr (\beta) +y^2\n(\beta)|\leq s}}
\Delta_{\frac{x}{y}} \\&\;=\;
\frac{d\haar_{\wh K}(z)}{|z-\beta|\,|z-\beta^\sigma|}\;.
\end{align*}
\etheo

Furthermore, we have error estimates in the arithmetic applications:
There exists $\kappa>0$ such that for every locally constant function
with compact support $\psi:\wh K\ra\CC$ in Theorems
\ref{theo:mertensfunctionfieldintro}, \ref{theo:fractidealintro} and
\ref{theo:appliquadirratintro}, or $\psi:\wh K- \{\beta,
\beta^\sigma\}\ra\CC$ in Theorems \ref{theo:applitwoquadirratdintro}
and \ref{theo:normformintro}, there are error terms in the above
equidistribution claims evaluated on $\psi$, of the form
$\bigO(s^{-\kappa})$ where $s=q^t$ in Theorem
\ref{theo:mertensfunctionfieldintro}, with for each result an explicit
control on the test function $\psi$ involving only some norm of
$\psi$, see in particular Section \ref{subsec:locconst}.

The link between the geometry described in the first part of this
introduction and the above arithmetic statements is provided by the
Bruhat-Tits tree of $(\PGL_2, \wh K)$, see \cite{Serre83} and Section
\ref{subsec:BruhatTitstrees} for background. We refer to Part
\ref{sect:arithappli} for more general arithmetic applications.

\bigskip\noindent{\small {\em Acknowledgements: } This work was
  partially supported by the NSF grants no 093207800 and DMS-1440140,
  while the third author was in residence at the MSRI, Berkeley CA,
  during the Spring 2015 and Fall 2016 semesters.  The second author
  thanks Universit\'e Paris-Sud, Forschungsinstitut f\"ur Mathematik
  of ETH Z\"urich, and Vilho, Yrj\"o ja Kalle V\"ais\"al\"an rahasto
  for their support during the preparation of this work. This research
  was supported by the CNRS PICS $n^0$ 6950 ``Equidistribution et
  comptage en courbure négative et applications arithmétiques''. We
  thank, for interesting discussions on this text, Y.~Benoist,
  J.~Buzzi (for his help in Sections \ref{subsec:codagesimplicial} and
  \ref{subsec:mixingratesimpgraphs}, and for kindly agreeing to insert
  Appendix \ref{appendixBuzzi} used in Section \ref{subsec:varprinc}),
  N.~Curien, S.~Mozes, M.~Pollicott (for his help in Section
  \ref{subsec:mixingratemetgraphs}), R.~Sharp and J.-B. Bost (for his
  help in Sections \ref{subsec:valuedfields} and
  \ref{subsec:mertens}).  We especially thank O.~Sarig for his help in
  Section \ref{subsec:mixingratesimpgraphs}: In a long email, he
  explained to us how to prove Theorem
  \ref{theo:critexpdecaysimpldynsymb}.}

\newpage
\section*{\hfill General notation\hfill}
\addcontentsline{toc}{section}{General notation}
\label{subsec:genenota}

\bigskip\bigskip
In this preamble, we introduce some general notation that will be used
throughout the book. We recommend the use of the List of symbols
(mostly in alphabetical order by the first letter) and of the Index
for easy references to the definitions in the text.

\medskip
\noindent Let $A$ be a subset of a set $E$. We denote by
$\gls{characteristicfunction}:E\ra \{0,1\}$ the {\it characteristic}
(or indicator) {\it function} of $A$: $\mathbbm{1}_A(x)=1$ if $x\in
A$, and $\mathbbm{1}_A(x)=0$ otherwise. We denote by
$\gls{complementary}=E-A$ the complementary subset of $A$ in $E$.

\smallskip\noindent 
We denote by $\lfloor x\rfloor=\sup\{n\in \NN\;:\; n\leq x\}$ the
lower integral part of any $x\in\RR$ and by $\lceil x\rceil=\inf\{n\in
\NN\;:\; x\leq n\}$ its upper integral part.

\smallskip\noindent 
We denote by $\gls{ln}$ the natural logarithm (with $\ln (e)=1$). 

\smallskip\noindent 
We denote by $\card(E)$ or by $|E|$ the order of a finite set $E$.

\smallskip\noindent 
We denote by $\|\mu\|$ the total mass of a finite positive measure
$\mu$. 

\smallskip\noindent 
If $(X,\A)$ and $(Y,\B)$ are measurable spaces, $f:X\ra Y$ a
measurable map, and $\mu$ a measure on $X$, we denote by $f_*\mu$
the image measure of $\mu$ by $f$, with $f_*\mu(B)=\mu(f^{-1}(B))$ for
every $B\in\B$. 

\smallskip\noindent 
If $(X,d)$ is a metric space, then $\gls{closedball}$ is the closed
ball with centre $x\in X$ and radius $r>0$.

\smallskip\noindent 
For every subset $A$ of a metric space and for every $\epsilon>0$, we
denote by $\gls{neighbourhood}$ the closed $\epsilon$-neighbourhood of
$A$, and by convention $\N_0A=\overline{A}$. We denote by
$\gls{minusneighbourhood}$ the set of points of $A$ at distance at least
$\epsilon$ from the complement of $A$.


\smallskip\noindent 
Given a topological space $Z$, we denote by $\gls{espacecontcompsupp}$
the vector space of continuous maps from $Z$ to $\RR$ with compact
support.

\smallskip\noindent 
Given a locally compact topological space $Z$, we denote by
$\gls{weakstarconv}$ the weak-star convergence of (Borel, positive)
measures on $Z$: We have $\mu_n\weakstar \mu$ if and only if
$\lim_{n\ra+\infty} \mu_n(f)= \mu(f)$ for every $f\in\C_c(Z)$.

\smallskip\noindent 
The {\it negative part}\index{negative!part} of a real-valued map $f$ is
$\gls{negativepart}=\max\{0,-f\}$.

\smallskip\noindent 
We denote by $\gls{diracmass}$ the unit Dirac mass at a point $x$ in
any measurable space.

\medskip
\noindent 
Finally, the symbol $\Box$ right at the end of a statement indicates that
this statement will not  be given a proof, either since a reference is
given or since it is an immediate consequence of previous statements.

\part{Geometry and dynamics in negative curvature}
\label{part1}

\chapter{Negatively curved geometry}
\label{sec:negcurv}

\section{Background on $\CAT(-1)$ spaces}
\label{subsec:catmoinsun}

Let $X$ be a geodesically complete proper $\CAT(-1)$ space, let
$x_0\in X$ be an arbitrary basepoint, and let $\Ga$ be a
nonelementary discrete group of isometries of $X$.

We refer for example to \cite{BriHae99} for the relevant terminology,
proofs and complements on these notions. In this Section, we recall some
definitions and notation for the sake of completeness.

A metric space is {\em proper}\index{proper} if its closed balls are
compact. A {\em geodesic}\index{geodesic} in a metric space $X'$ is an
isometric map $c$ from an interval $I$ of $\RR$ into $X'$.\footnote{We
  say that $c$ is a {\em geodesic segment}\index{geodesic!segment} if
  $I$ is compact, a {\em geodesic ray}\index{geodesic!ray} if $I$ is a
  half-infinite interval, and a {\em geodesic
    line}\index{geodesic!line} if $I=\RR$.}  A metric space $X'$ is
{\em geodesic}\index{geodesic} if for all $x,y\in X'$, there exists a
geodesic segment $c:[a,b]\ra X'$ from $x=c(a)$ to $y=c(b)$. A geodesic
metric space $X'$ is {\em geodesically complete}\index{geodesically
  complete} (or has {\em extendible geodesics})\index{extendible
  geodesics} if any isometric map from an interval in $\RR$ to $X'$
extends to at least one isometric map from $\RR$ to $X'$. A {\em
  comparison triangle}\index{comparison triangle} of a triple of
points $(x,y,z)$ in a metric space $X'$ is a (unique up to isometry)
triple of points $(\overline{x},\overline{y},\overline{z})$ in the
real hyperbolic plane $\HH^2_\RR$ such that $d(x,y)= d(\overline{x},
\overline{y})$, $d(y,z)=d(\overline{y},\overline{z})$ and $d(z,x)=
d(\overline{z},\overline{x})$.

\begin{center}
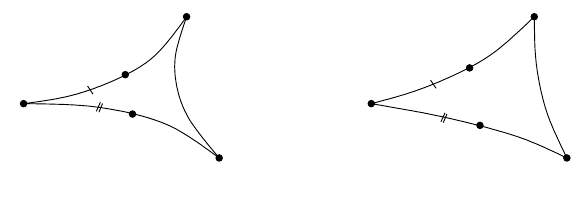
\end{center}

A metric space $X'$ is $\gls{catmoinsun}$ if it is geodesic and if for
every triple of points $(x,y,z)$ in $X'$, for all geodesic segments
$a,b$ respectively from $x$ to $y$ and from $x$ to $z$, and for all
points $p,q$ in the image of $a,b$ respectively, if
$(\overline{x},\overline{y},\overline{z})$ is a comparison triangle of
$(x,y,z)$, if $\overline{p}$ (resp.~$\overline{q}$) is the point on
the geodesic segment from $\overline{x}$ to $\overline{y}$
(resp.~$\overline{z}$) at distance $d(x,p)$ (resp.~$d(x,q)$) from
$\overline{x}$, then $d(p,q)\leq d(\overline{p},\overline{q})$.

We will put a special emphasis on the case when $X$ is a (proper,
geodesically complete) {\em $\RR$-tree}\index{tree@$\RR$-tree}, that
is, a uniquely arcwise connected geodesic metric space.  In the
Introduction, we have denoted by $Y$ the geodesically complete proper
locally $\CAT(-1)$ good orbispace $\Ga\bs X$, see for instance
\cite[Ch.~11]{GhyHar90} for the terminology.

Two geodesic rays $\rho,\rho':[0,+\infty[\;\ra X$ are {\em
      asymptotic}\index{geodesic!ray!asymptotic} if their images are
at finite Hausdorff distance, or equivalently if there exists
$a\in\RR$ such that $\lim_{t\ra+\infty}d(\rho(t),\rho'(t+a))=0$. We
denote by $\gls{boundaryatinfinity}$ the {\em space at
  infinity}\index{space at infinity} of $X$, which consists of the
asymptotic classes of geodesic rays in $X$, and we endow it with the
quotient topology of the compact-open topology. It coincides with
the space of (Freudenthal's) ends of $X$ when $X$ is an $\RR$-tree.

Let $x,y\in X$ and $\xi,\xi'\in\partial_\infty X$. We denote by
$[x,y]=[x,y]$ the unique image of a geodesic segment from $x$ to
$y$. We denote by $[x,\xi[$ the image of the unique geodesic ray
$\rho:[0,+\infty[\ra X$ in the asymptotic class $\xi$ with
$\rho(0)=x$, and we say that $\rho$ {\em starts from $x$} and
{\em ends at $\xi$}.  We denote by $]\xi,\xi'[\;=\;]\xi',\xi[$
the unique image of a geodesic line $\ell: \RR\ra X$ with
$t\mapsto \ell(t)$ and $t\mapsto \ell(-t)$ in the asymptotic
classes $\xi$ and $\xi'$ respectively, and we say that $\ell$
{\em starts from $\xi'$} and {\em ends at $\xi$}.

We endow the disjoint union $X\cup \partial_\infty X$ with the unique
metrisable compact topology (independent of $x_0$) inducing the
above topologies on $X$ and $\partial_\infty X$, such that a
sequence $(y_i)_{i\in\NN}$ in $X$ converges to $\xi\in \partial_\infty
X$ if and only if $\lim_{i\ra+\infty}d(x_0,y_i)=+\infty$ and, with
$c_i:[0,d(x_0,y_i)]\ra X$ the geodesic segment from $x_0$ to $y_i$ and
$\rho:[0,+\infty[\ra X$ the geodesic ray from $x_0$ to $\xi$, we have
$\lim_{i\ra+\infty}d(c_i(t),\rho(t))=0$ for every $t\geq 0$.

We denote by $\gls{isometrygroup}$ the isometry group of $X$, and we
endow it with the compact-open topology. Its action on $X$ uniquely
extends to a continuous action on $X\cup\partial_\infty X$. We say
that a discrete subgroup $\Ga'$ of $\Isom(X)$ is {\em
  nonelementary}\index{nonelementary} if it does not fix a point or an
unordered pair of points in $X\cup\partial_\infty X$.  We denote by
$\gls{limitset}$ the {\em limit set}\index{limit set} of $\Ga$, which
is the set of accumulation points in $\partial_\infty X$ of any orbit
of $\Ga$ in $X$. It is the smallest closed nonempty $\Ga$-invariant
subset of $\partial_\infty X$.

A subset $D$ of $X\cup\partial_\infty X$ is {\it convex}\index{convex}
if for all $u,v\in D$, the image of the unique geodesic segment, ray
or line from $u$ to $v$ is contained in $D$. We denote by
$\gls{convexhull}$ the {\em convex hull}\index{convex hull} in $X$ of
$\Lambda\Ga$, which is the intersection of the closed convex subsets
of $X\cup \partial_\infty X$ containing $\Lambda\Ga$.  When $X$ is an
$\RR$-tree, then a subset $D$ of $X$ is convex if and only if it is
connected, and we will call it a {\em subtree}.\index{subtree} In
particular, if $X$ is an $\RR$-tree, then $\C\Lambda\Ga$ is equal to
the union of the geodesic lines between pairs of distinct points in
$\Lambda\Ga$, since this union is connected and contained in
$\C\Lambda\Ga$.

A point $\xi \in\partial_\infty X$ is called a {\em conical limit
  point}\index{conical limit!point}\index{limit point!conical} of
$\Ga$ if there exists a sequence of orbit points of $x_0$ under $\Ga$
converging to $\xi$ while staying at bounded distance from a geodesic
ray ending at $\xi$.  The set of conical limit points of $\Ga$ is the
{\em conical limit set}\index{conical limit!set}
$\gls{connicallimitset}$ of $\Ga$.

A point $p\in \Lambda\Ga$ is a {\em bounded parabolic limit
  point}\index{bounded parabolic limit point}%
\index{limit point!bounded parabolic} of $\Ga$ if its stabiliser
$\Ga_p$ in $\Ga$ acts properly discontinuously with compact quotient
on $\Lambda\Ga-\{p\}$. The discrete nonelementary group of isometries
$\Ga$ of $X$ is said to be {\em geometrically
  finite}\index{geometrically finite} if every element of $\Lambda\Ga$
is either a conical limit point or a bounded parabolic limit point of
$\Ga$. See for instance \cite{Bowditch95}, as well as \cite{Paulin04b}
when $X$ is an $\RR$-tree, and \cite{DasSimUrb14arxiv} for a very
interesting study of equivalent conditions in an even greater
generality.

For all $x\in X\cup\partial_\infty X$ and $A\subset X$, the {\em
  shadow of $A$ seen from $x$}\index{shadow} is the subset
$\gls{shadow}$ of $\partial_\infty X$ consisting of the endpoints
towards $+\infty$ of the geodesic rays starting at $x$ and meeting $A$
if $x\in X$, and of the geodesic lines starting at $x$ and meeting $A$
if $x\in\partial_\infty X$.

The {\em translation length} \index{translation!length} of an isometry
$\ga\in\Isom(X)$ is
$$
\gls{distancetranslat}=\inf_{x\in X} d(x,\ga x)\,.
$$
An element $\ga\in\Isom(X)$ is {\em elliptic}\index{elliptic} if it
fixes a point in $X$, and then $\lambda(\ga)=0$. An element
$\ga\in\Isom(X)$ is {\em parabolic}\index{parabolic} if it is not
elliptic and fixes a unique point in $\partial_\infty X$, and then
$\lambda(\ga)=0$. An element $\ga\in\Isom(X)$ is {\em
  loxodromic}\index{loxodromic} if $\lambda(\ga)>0$, and then
$$
\Ax_{\ga}=\{x\in X\;:\; d(x,\ga x)= \lambda(\ga)\}
$$
is (the image of) a geodesic line in $X$, called the {\em translation
  axis}\index{translation!axis} of $\ga$. In particular, the
restriction of $\ga$ to $\Ax_{\ga}$ is conjugated, by any isometry
between $\Ax_{\ga}$ and $\RR$, to a translation of $\RR$ of the form
$t\mapsto t\pm\lambda(\ga)$. For all $\beta\in\Isom(X)$ and
$n\in\ZZ-\{0\}$, we have
\begin{equation}\label{eq:propritranslaxis}
\Ax_{\beta\ga\beta^{-1}}=\beta\Ax_\ga,\;\;\;
\lambda(\beta\ga\beta^{-1})=\lambda(\ga),\;\;\;
\Ax_{\ga^n}=\Ax_\ga,\;\;\;\lambda(\ga^n)= |n|\lambda(\ga)\;.
\end{equation}
A loxodromic element $\ga\in\Isom(X)$ has exactly two fixed points
$\ga_-,\ga_+$ in $X\cup\partial_\infty X$, with $\ga_-\in
\partial_\infty X$ its {\em repulsive}\index{repulsive} fixed
point\footnote{For every $x\in X\cup(\partial_\infty X-\{\ga_+\}$, we
  have $\lim_{n\ra+\infty}\ga^{-n}x=\ga_-$ .} and $\ga_+\in
\partial_\infty X$ its {\em attractive}\index{attractive} fixed
point.\footnote{For every $x\in X\cup(\partial_\infty X-\{\ga_-\})$,
  we have $\lim_{n\ra+\infty}\ga^{+n}x=\ga_+$ .}

We will need the following well-known lemma later on. An element of
$\Ga$ is {\em primitive}\index{primitive} in $\Ga$ if there is no
$\ga_0\in\Ga$ and $k\in\NN-\{0,1\}$ such that $\ga=\ga_0^k$. Note that
there might exist a primitive loxodromic element $\ga$ in $\Ga$, whose
translation length is not minimal amongst the translation lengths of
the loxodromic elements $\ga'\in\Ga$ with $\Ax_\ga=
\Ax_{\ga'}$.\footnote{For instance, if $X$ is the real hyperbolic
  plane $\HH^2_\RR$ and if $\Ga$ contains an orientation preserving
  loxodromic element $\ga$ such that the stabiliser in $\Ga$ of
  $\Ax_\ga$ is generated by $\ga$ and by the symmetry $\sigma$ with
  respect to $\Ax_\ga$, then $\ga^{2^n}\sigma$ is primitive for all 
  $n\in\ZZ$.}

\blemm \label{lem:axespastropproches} (1) For every loxodromic element
$\ga\in\Ga$, there exist $k'\in\NN-\{0\}$ and a primitive loxodromic
element $\ga_1\in\Ga$ such that $\ga=\ga_1^{k'}$, and there exist
$k\in\NN-\{0\}$, a primitive loxodromic element $\ga_0\in\Ga$ whose
translation length is minimal amongst the translation lengths of the
loxodromic elements $\ga'\in\Ga$ with $\Ax_\ga= \Ax_{\ga'}$, and an
element $\ga'\in\Ga$ pointwise fixing $\Ax_\ga$ such that
$\ga=\ga_0^k\ga'$.

\medskip
(2) For every compact subset $K$ of $X$, for all $A>0$ and $r>0$,
there exists $L>0$ such that for all loxodromic elements
$\ga,\alpha\in\Ga$, if $\Ax_\ga$ meets $K$, if $\lambda(\alpha)=
\lambda(\ga) \leq A$ and if $\Ax_\ga$ and $\Ax_\alpha$ have segments
of length at least $L$ at Hausdorff distance at most $r$, then
$\Ax_\alpha=\Ax_\ga$.

\medskip
(3) For every compact subset $K$ of $X$ and for every $A>0$, there
exists $N\in\NN$ such that for every loxodromic element $\ga\in\Ga$
whose translation axis meets $K$, the cardinality of the set of
loxodromic elements $\alpha\in\Ga$ with $\Ax_\alpha=\Ax_\ga$ and
$\lambda(\alpha) \leq A$ is at most $N$.
\elemm

\dem (1) If $\ga\in\Ga$ is loxodromic, then the group of restrictions
to $\Ax_\ga$ of the elements of $\Ga$ preserving $\Ax_\ga$ is
conjugated, by any isometry between $\Ax_\ga$ and $\RR$, to a discrete
group of isometries $\Lambda$ of $\RR$. Since replacing $\ga$ by an
element of $\Ga$ having a power at least $2$ equal to $\ga$ strictly
decreases the translation length and by Equation
\eqref{eq:propritranslaxis}, the first claim of (1) holds. The normal
subgroup of $\Lambda$ consisting of translations is isomorphic to
$\ZZ$, generated by the conjugate of the restriction to $\Ax_\ga$ of
an element $\ga_0\in\Ga$. Any such element has minimal translation
length, hence is primitive since if there exist $\ga_1\in\Ga$ and
$n\in\NN-\{0,1\}$ with $\ga_0=\ga_1^n$, by Equation
\eqref{eq:propritranslaxis}, we would have $\Ax_{\ga_1}=\Ax_{\ga_0}$
and $\lambda(\ga_1)=
\frac{1}{n}\,\lambda(\ga_0)<\lambda(\ga_0)$. There exists
$k\in\ZZ-\{0\}$ such that the restrictions of $\ga$ and $\ga_0^k$ to
$\Ax_\ga$ coincide. Hence $\ga'=\ga_0^{-k}\ga$ pointwise fixes
$\Ax_\ga$, and up to replacing $\ga_0$ by its inverse, Assertion (1)
of Lemma \ref{lem:axespastropproches} holds.

\medskip
(2) Since the action of $\Ga$ on $X$ is properly discontinuous, there
exists $N=N(K,A,r)\geq 1$ such that for every loxodromic element
$\ga\in\Ga$ whose translation axis meets $K$ and whose translation
length is at most $A$, for every $x\in \Ax_\ga$, the cardinality of
$\{\beta\in\Ga\;:\; d(x,\beta x)\leq 2r\}$ is at most $N$. Let
$L=AN$. For every loxodromic element $\alpha\in\Ga$ with
$\lambda(\alpha) =\lambda(\ga) \leq A$, assume that $[x,y]$ and
$[x',y']$ are segments in $\Ax_\ga$ and $\Ax_\alpha$ respectively,
with length exactly $L$ such that $d(x,x'),d(y,y')\le r$. We may
assume, up to replacing $\ga$ and $\alpha$ by their inverses, that
$\ga$ translates from $x$ towards $y$ and $\alpha$ translates from
$x'$ towards $y'$. In particular for $k=0,\dots,N$, we have
$d(\alpha^{-k}\ga^kx,x) \leq d(\ga^kx,\alpha^{k}x')+ d(x',x)\leq 2r$
since $\ga^kx$ and $\alpha^{k} x'$ are respectively the points at
distance $k\lambda(\ga) \leq kA\leq L$ from $x$ and $x'$ on the
segments $[x,y]$ and $[x',y']$. Hence by the definition of $N$, there
exists $k\neq k'$ such that $\alpha^{-k}\ga^k =\alpha^{-k'}
\ga^{k'}$. Therefore $\ga^{k-k'} = \alpha^{k-k'}$, which implies by
Equation \eqref{eq:propritranslaxis} that $\Ax_\ga=\Ax_\alpha$.

\medskip
(3) Since the action of $\Ga$ on $X$ is properly discontinuous, there
exists $N'\in\NN$ such that the cardinality of the stabiliser in $\Ga$
of a point of $K$ is at most $N'$, and there exists $\eta>0$ such that
$\lambda(\ga)\geq \eta$ for every loxodromic element $\ga\in\Ga$ whose
translation axis meets $K$. Let us fix such an element $\ga\in\Ga$. We
may assume that its translation length is minimal amongst the
translation lengths of the loxodromic elements $\ga'\in\Ga$ with
$\Ax_\ga= \Ax_{\ga'}$.  Then as seen in the proof of Assertion (1),
for every loxodromic element $\alpha\in\Ga$ with $\Ax_\alpha=\Ax_\ga$,
there exist $k\in \NN-\{0\}$ and $\alpha'\in\Ga$ pointwise fixing
$\Ax_\ga$ such that $\alpha=\ga^k\alpha'$. Thus if $\lambda(\alpha)
\leq A$, then $|k|= \frac{\lambda(\alpha)}{\lambda(\ga)}\leq
\frac{A}{\eta}$, and there are at most
$N=N'(2\lceil\frac{A}{\eta}\rceil +1)$ such elements $\alpha$.
\cqfd

\bigskip
For every $x\in X$, recall that the Gromov-Bourdon {\em visual
  distance}\index{visual distance}\index{distance!visual}
$\gls{visualdist}$ on $\partial_\infty X$ seen from $x$ (see
\cite{Bourdon95}) is defined by
\begin{equation}\label{eq:defidistvis}
d_x(\xi,\eta)=
\lim_{t\ra+\infty} e^{\frac{1}{2}(d(\xi_t,\,\eta_t)-d(x,\,\xi_t)-d(x,\,\eta_t))}\;,
\end{equation}
where $\xi,\eta\in\partial_\infty X$ and $t\mapsto \xi_t,\eta_t$ are
any geodesic rays ending at $\xi,\eta$ respectively. If $X$ is an
$\RR$-tree, if $\xi,\eta\in \partial_\infty X$ are distinct, if $p\in
X$ is such that $[x,p]=\mathopen{[}x,\xi\mathclose{[}
\cap\mathopen{[}x,\eta\mathclose{[}$, then
\begin{equation}\label{eq:distviscastree}
d_x(\xi,\eta)= e^{-d(x,\,p)}\;.
\end{equation}         
For all $x\in X$, $\xi,\eta\in\partial_\infty X$ and $\ga\in\Isom(X)$,
we have
$$
d_{\ga x}(\ga \xi,\ga\eta)=d_x(\xi,\eta)\;.
$$
By the triangle inequality, for all $x,y\in X$ and
$\xi,\eta\in\partial_\infty X$, we have
\begin{equation}\label{eq:lipequivdistvis}
e^{-d(x,\,y)}\leq \frac{d_x(\xi,\eta)}{d_y(\xi,\eta)}\leq e^{d(x,\,y)}\;.
\end{equation} 
In particular, the identity map from $(\partial_\infty X,d_{x})$ to
$(\partial_\infty X,d_{y})$ is a bilipschitz homeomorphism. Under our
assumptions, $(\partial_\infty X,d_{x_0})$ is hence a compact metric
space, on which $\Isom(X)$ acts by bilipschitz homeomorphisms.  The
following well-known result compares shadows of balls to balls for the
visual distance.

\blemm\label{lem:comparshadowball} For every geodesic ray $\rho$ in
$X$, starting from $x\in X$ and ending at $\xi\in\partial_\infty X$,
for all $R\geq 0$ and $t\in\;]R,+\infty[\,$, we have
$$
B_{d_x}(\xi,R\,e^{-t}) \subset \OOO_xB(\rho(t),R)\subset 
B_{d_x}(\xi,e^R\,e^{-t})\;.
$$
\elemm

\dem In order to prove the left inclusion, we adapt the proof of the
left inclusion in \cite[Lem.~3.1]{HerPau01} (which only uses the
$\CAT(-1)$ property). Let $\xi'\in B_{d_x}(\xi,R\,e^{-t})-\{\xi\}$,
let $\rho'$ be the geodesic ray from $x$ to $\xi'$ and let $p$ be the
closest point to $w=\rho(t)$ on the image of $\rho'$. For every $s>t$,
let $y=\rho(s)$ and $z=\rho'(s)$. Let $(\overline{x},\overline{y},
\overline{z})$ be a comparison triangle in $\HH^2_\RR$ of $(x,y,z)$,
and let $\overline{\theta}\in[0,\pi]$ be its angle at $\overline{x}$.
Let $\overline{w}$ be the point on $[\overline{x},\overline{y}]$ at
distance $t$ from $\overline{x}$, and let $\overline{p}$ be its
orthogonal projection on $[\overline{x},\overline{z}]$. By the
$\CAT(-1)$ property, we have
$$
d(w,p)\leq d(\overline{w},\overline{p})\;.
$$

\begin{center}
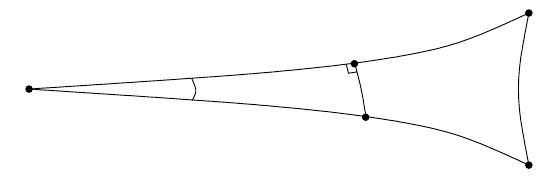
\end{center}

By the hyperbolic sine rule for right angled triangles in $\hdr$, we
have
$$
\sin\,\overline{\theta} =
\frac{\sinh d(\overline{w},\overline{p})}{\sinh t}\;\;\;{\rm and}
\;\;\;
\sin\,\frac{\overline{\theta}}{2}=
\frac{\sinh \frac{1}{2}\,d(\overline{y},\overline{z})}
     {\sinh d(\overline{x},\overline{y})}
=\frac{\sinh \frac{1}{2}\,d(y,z)}{\sinh s}\;.
$$
Hence
$$
d(w,p)\leq d(\overline{w},\overline{p})\leq
\frac{e^{d(\overline{x},\overline{w})}}{2 \sinh d(\overline{x},\overline{w})}
\sinh d(\overline{w},\overline{p})=
\frac{1}{2}e^t\sin\overline{\theta}\leq
e^t\sin\,\frac{\overline{\theta}}{2}\;.
$$
Since
$$
\lim_{s\ra+\infty}\sin\,\frac{\overline{\theta}}{2}=
\lim_{s\ra+\infty}e^{\frac{1}{2}\,d(y,z)-s}=d_x(\xi,\xi')\leq R\,e^{-t}\;,
$$
we hence have $d(w,p)\leq R$, so that $\xi'\in \OOO_xB(w,R)$,
as wanted.

\medskip 
In order to prove the right inclusion in Lemma
\ref{lem:comparshadowball}, let $\xi'\in \OOO_xB(\rho(t),R)$ and let
$\rho'$ be the geodesic ray from $x$ to $\xi'$. The closest point $p$
to $\rho(t)$ on the image of $\rho'$ satisfies $d(p,\rho(t))\leq R$,
hence $d(x,p)\geq t-R$. Therefore, for $t'$ large enough,
$d(\rho'(t'),p)\leq t'-(t-R)$, and
\begin{align*}
d_x(\xi',\xi)& \leq \limsup_{t'\ra+\infty} \;
e^{\frac{1}{2}\big(\,d(\rho(t'),\,\rho(t))+d(\rho(t),\,p)+d(p,\,\rho'(t'))\,\big)-t'}
\\ &\leq \lim_{t'\ra+\infty} 
e^{\frac{1}{2}\big((t'-t) + R+ (t'-t+R)\big)-t'}=e^Re^{-t}\;.
\end{align*} 
Therefore $\xi'\in B_{d_x}(\xi,e^R\,e^{-t})$, as wanted.
\cqfd

\medskip
The {\em Busemann cocycle}\index{Busemann
  cocycle}\index{cocycle!Busemann} of $X$ is the map
$\beta: \partial_{\infty} X\times X\times X\to\RR$ defined by
$$
(\xi,x,y)\mapsto \beta_{\xi}(x,y)=
\lim_{t\to+\infty}d(\xi_t,x)-d(\xi_t,y)\;,
$$
where $t\mapsto \xi_t$ is any geodesic ray ending at $\xi$. If $X$ 
is an $\RR$-tree, if $p\in X$ is such that $\mathopen{[}x,\xi
\mathclose{[} \cap \mathopen{[}y,\xi\mathclose{[}= \mathopen{[}p,
\xi\mathclose{[}\,$, then
\begin{equation}\label{eq:buscocycastree}
\beta_\xi(x,y)= d(x,\,p)-d(y,\,p)\;.
\end{equation}
The triangle inequality gives immediately the upper bound
\begin{equation}\label{eq:majobusemann}
|\,\beta_{\xi}(x,y)\,|\leq d(x,y)\;.
\end{equation}
The {\em horosphere}\index{horosphere} with centre $\xi \in
\partial_{\infty}X$ through $x\in X$ is $\{y\in X\;:\;
\beta_{\xi}(x,y)=0\}$, and $\{y\in X\;:\; \beta_{\xi}(x,y)\leq 0\}$ is
the (closed) {\em horoball}\index{horoball} centred at $\xi$ bounded
by this horosphere. Horoballs are nonempty proper closed (strictly)
convex subsets of $X$. Given a horoball $\H$ and $t\geq 0$, we
denote by $\gls{horoshrink}=\{x\in\H\;:\;d(x,\partial \H)\geq t\}$ the
horoball contained in $\H$ (hence centred at the same point at
infinity as $\H$) whose boundary is at distance $t$ from the boundary
of $\H$.

\section{Generalised geodesic lines}
\label{subsec:lines}

Let $\gls{gengeod}$ be the space of $1$-Lipschitz maps $w:\RR\to X$
which are isometric on a closed interval and locally constant outside
it.\footnote{that is, constant on each complementary component} This
space has been introduced by Bartels and L\"uck in \cite{BarLuc12}, to
which we refer for the following basic properties. The elements of
$\gengeod X$ are called the {\em generalised geodesic
  lines}\index{generalised geodesic!line}%
\index{geodesic!line!generalised} of $X$.  Any geodesic segment or ray
of $X$ will be considered as an element of $\gengeod X$, by extending
it continuously to $\RR$ as locally constant outside its domain of
definition.

We endow $\gengeod X$ with
the distance $\gls{distGchechX}$ defined by
\begin{equation}\label{eq:geodesicmetric}
\forall\;w, w'\in\gengeod X,\;\;\; d(w, w')=
\int_{-\infty}^{+\infty} d(w(t),w'(t))\;e^{-2 |t|}\,dt\,.
\end{equation}
The group $\Isom(X)$ acts isometrically on $\gengeod X$ by
postcomposition. The distance $d$ induces the topology of uniform
convergence on compact subsets on $\gengeod X$, and $\gengeod X$ is a
proper metric space. 

The {\em geodesic flow}\index{geodesic!flow} $\gls{geodesicflow}$ on
$\gengeod X$ is the one-parameter group of homeomorphisms of the space
$\gengeod X$ defined by $\flow t w:s \mapsto w(s+t)$ for all
$w\in\gengeod X$ and $t\in\RR$. It commutes with the action of
$\Isom(X)$.  If $w$ is isometric exactly on the interval $I$, then
$\flow {-t} w$ is isometric exactly on the interval $t+I$.

The {\em footpoint projection}\index{footpoint projection} is the
$\Isom(X)$-equivariant $\frac{1}{2}$-H\"older-continuous\footnote{See
  Section \ref{subsec:holdercont} for the definition of the (locally
  uniform) H\"older-continuity used in this book, and Proposition
  \ref{prop:footpointholdercont} for a proof of this claim.} map
$\gls{footpointprojection}:\gengeod X\to X$ defined by $\pi(w)=w(0)$
for all $w\in\gengeod X$. The {\em antipodal map}\index{antipodal map}
of $\gengeod X$ is the $\Isom(X)$-equivariant isometric map
$\gls{antipodalmap}: \gengeod X\to\gengeod X$ defined by $\iota
w:s\mapsto w(-s)$ for all $w\in\gengeod X$, which satisfies
$\iota\circ\flow{t}=\flow{-t}\circ\iota$ for every $t\in\RR$ and
$\pi\circ\iota=\pi$.

The {\em positive}\index{positive endpoint}\index{endpoint!positive}
and {\em negative}\index{negative!endpoint}\index{endpoint!negative}
endpoint maps are the continuous maps from $\gengeod X$ to
$X\cup\partial_\infty X$ defined by
$$
w\mapsto \gls{endpointmap}=\lim_{t\to\pm\infty}w(t)\,.
$$

The space $\gls{spacegeodlin}$ of geodesic lines in $X$ is the
$\Isom(X)$-invariant closed metric subspace of $\gengeod X$ consisting
of the elements $\ell\in\gengeod X$ with $\gls{endpointmapline} \in
\partial_\infty X$.  Note that the distances on $\G X$ considered in
\cite{BarLuc12} and \cite{PauPolSha15} are topologically equivalent
to, although slightly different from, the restriction to $\G X$ of the
distance defined in Equation \eqref{eq:geodesicmetric}.  The factor
$e^{- 2 |t|}$ in this equation, sufficient in order to deal with
H\"older-continuity issues, is replaced by $e^{-t^2}/\sqrt{\pi}$ in
\cite{PauPolSha15} and by $e^{- |t|}/2$ in \cite{BarLuc12}.

Note that for all $w\in \gengeod X$ and $s\in\RR$, we have
\begin{equation}\label{eq:disttranslatgeod}
d(w, \flow{s}w)\leq |s|\,,
\end{equation}
with equality if $w\in \G X$.

We will also consider the $\Isom(X)$-invariant closed subspaces
$$
\gls{gengeodpm}=\{w\in\gengeod X:w_\pm\in\partial_\infty X\}\,,
$$ 
and their $\Isom(X)$-invariant closed subspaces $\gls{gengeodpmz}$
consisting of the elements $\rho\in\G_\pm X$ which are isometric
exactly on $\pm\mathopen{[}0,+\infty\mathclose{[}$. 

The subspaces $\G X$ and $\G_\pm X$ satisfy $\G_- X\cap\G_+ X=\G X$
and they are invariant under the geodesic flow. The antipodal map
$\iota$ preserves $\G X$, and maps $\G_\pm X$ to $\G_\mp X$ as well as
$\G_{\pm,\,0}X$ to $\G_{\mp,\,0}X$. We denote again by
$\gls{antipodalmapdown}: \Ga\bs \gengeod X\to \Ga\bs \gengeod X$ and by
$\gls{geodesicflowdown}$ with
$\flow t: \Ga\bs \gengeod X\to \Ga\bs \gengeod X$ the
quotient maps of $\iota$ and $\flow t$, for every $t\in\RR$.

Let $w\in\gengeod X$ be isometric exactly on an interval $I$ of $\RR$.
If $I$ is compact then $w$ is a {\em (generalised) geodesic segment}%
\index{generalised geodesic!segment}, and if $I= \mathopen{]} -\infty,
a\mathclose{]}$ or $I=\mathopen{[}a, +\infty\mathclose{[}$ for some
$a\in\RR$, then $w$ is a {\em (generalised) (negative or positive)
  geodesic ray}\index{generalised geodesic!ray}%
\index{ray!generalised geodesic} in $X$. Any geodesic line $\wh w\in\G
X$ such that $\wh w|_I=w|_I$ is an {\em extension}\index{extension} of
$w$.  Note that $\wh w$ is an extension of $w$ if and only if $\ga\wh
w$ is an extension of $\ga w$ for any $\ga\in\Isom(X)$, if and only if
$\iota\wh w$ is an extension of $\iota w$, and if and only if $\flow s
\wh w$ is an extension of $\flow s w$ for any $s\in\RR$.  For any
subset $\Omega'$ of $\gengeod X$ and any subset $A$ of $\RR$, let
$$
\Omega'|_A=\{w|_A:w\in\Omega'\}\,.
$$

\brema Let $(\ell_i)_{i\in\NN}$ be a sequence of generalised geodesic
lines such that $[t_i^-,t_i^+]$ is the maximal segment on which
$\ell_i$ is isometric. Let $(s_i)_{i\in\NN}$ be a sequence in $\RR$
such that $t_i^\pm-s_i\to\pm\infty$ as $i\to+\infty$ and $\ell_i(s_i)$
stays in a compact subset of $X$, then $d(\ell_i,\G
X)\to 0$ as $i\to+\infty$. Furthermore if $(s_i)_{i\in\NN}$ is
bounded, then up to extracting a subsequence, $(\ell_i)_{i\in\NN}$
converges to an element in $\G X$.  
\erema

This conceptually important observation explains how it is conceivable
that long common perpendicular segments may equidistribute towards
measures supported on geodesic lines. See Chapter \ref{sec:equidarcs}
for further developments of these ideas.

\section{The unit tangent bundle}
\label{subsec:unitbundle}

In this book, we define the {\em unit tangent
  bundle}\index{unit!tangent bundle} $\gls{unittangentspace}$ of $X$
as the space of germs at $0$ of the geodesic lines in $X$. It is the
quotient space
$$
T^1X=\G X/\sim
$$ 
where $\ell\sim\ell'$ if and only if there exists $\epsilon>0$ such
that $\ell|_{[-\epsilon,\epsilon]}= \ell'|_{[-\epsilon, \epsilon]}$.  The 
canonical projection from $\G X$ to $T^1X$ will be denoted by
$\ell\mapsto \gls{germgeodline}$. When $X$ is a Riemannian manifold,
the spaces $\G X$ and $T^1X$ canonically identify with the usual unit
tangent bundle of $X$, but in general, the map $\ell\mapsto v_\ell$
has infinite fibers.

We endow $T^1X$ with the quotient distance $\gls{distgerm}$ of the
distance of $\G X$, defined by:
\begin{equation}\label{eq:germdistance}
  \forall\;v, v'\in T^1X,\;\;\; d_{T^1X}(v, v')= 
  \inf_{\ell,\,\ell'\in\G X\;:\;v=v_\ell,\, v'=v_{\ell'}} d(\ell, \ell')\,.
\end{equation}
It is easy to check that this distance is indeed Hausdorff, hence that
$T^1X$ is locally compact, and that it induces on $T^1X$ the quotient
topology of the compact-open topology of $\G X$.  The map $\ell\mapsto
v_\ell$ is $1$-Lipschitz.

The action of $\Isom(X)$ on $\G X$ induces an isometric action of
$\Isom(X)$ on $T^1X$. The antipodal map and the footpoint projection
restricted to $\G X$ respectively induce an $\Isom(X)$-equivariant
isometric map $\iota:T^1X\ra T^1X$ and an $\Isom(X)$-equivariant
$\frac{1}{2}$-H\"older-continuous map $\pi:T^1X\ra X$ called the {\em
  antipodal map}\index{antipodal map} and {\em footpoint projection}%
\index{footpoint projection} of $T^1X$.  The canonical projection from
$\G X$ to $T^1X$ is $\Isom(X)$-equivariant and commutes with the
antipodal map: For all $\ga \in\Isom(X)$ and $\ell\in\G X$, we have
$\ga v_\ell = v_{\ga\ell}$, $\iota v_\ell =v_{\iota\ell}$ and
$\pi(v_\ell) = \pi(\ell)$. We denote again by $\iota: \Ga\bs T^1 X\to
\Ga\bs T^1 X$ the quotient map of $\iota$.

\bigskip 
Let $\gls{boundaryatinfinitybis}$ be the subset of $\partial_\infty
X\times\partial_\infty X$ which consists of pairs of distinct points
at infinity of $X$.  {\em Hopf's parametrisation}%
\index{Hopf parametrisation} of $\G X$ is the homeomorphism which
identifies $\G X$ with $\partial_\infty^2 X \times\RR$, by the map
$\ell \mapsto (\ell_-,\ell_+,t)$, where $t$ is the signed distance
from the closest point to the basepoint $x_0$ on the geodesic line
$\ell$ to $\ell(0)$.\footnote{More precisely, $\ell(t)$ is the closest
  point to $x_0$ on $\ell$.}  We have $\flow{s} (\ell_-,\ell_+,t) =
(\ell_-,\ell_+,t+s)$ for all $s\in\RR$, and for all $ \ga\in \Ga$, we
have $\ga (\ell_-,\ell_+,t) =(\ga \ell_-,\ga \ell_+, t+
t_{\ga,\,\ell_-,\,\ell_+})$ where $t_{\ga,\,\ell_-,\,\,\ell_+} \in\RR$
depends only on $\ga$, $\ell_-$ and $\ell_+$.  In Hopf's
parametrisation, the restriction of the antipodal map to $\G X$ is the
map $(\ell_-,\ell_+,t)\mapsto (\ell_+,\ell_-,-t)$.

\medskip The {\em strong stable leaf}\index{strong!stable
  leaf}\index{leaf!strong stable} of $w\in\G_+ X$ is
$$
\gls{Wss}=
\big\{\ell\in \G X\;:\lim_{t\to+\infty}d(\ell(t),w(t))=0\big\}\,,  
$$
and the {\em strong unstable leaf}\index{strong!unstable
  leaf}\index{leaf!strong unstable} of $w\in\G_-X$ is
$$
\gls{Wsu}=\iota\wss(\iota w)=
\big\{\ell\in  \G X\;:\lim_{t\to-\infty}d(\ell(t),w(t))= 0 \big\}\,. 
$$

\begin{center}
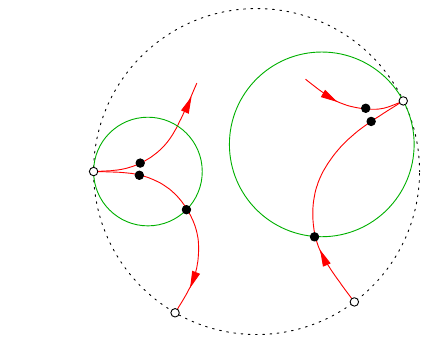
\end{center}

For every $w\in\G_\pm X$, let $\gls{Hamenstadtdistancepm}$ be
{\em Hamenst\"adt's distance}%
\index{Hamenst\"adt's!distance}\index{distance!Hamenst\"adt's} on
$\wssu(w)$ defined as follows:\footnote{See \cite[Appendix]{HerPau97} 
and compare with \cite{Hamenstadt89}.} for all $\ell,\ell'\in
\wssu(w)$, let
$$
d_{\wssu(w)}(\ell,\ell') = 
\lim_{t\ra+\infty} e^{\frac{1}{2}d(\ell(\mp t),\;\ell'(\mp t))-t}\,.
$$
The above limits exist, and Hamenst\"adt's distances are distances
inducing the original topology on $\wssu(w)$. For all $\ell,\ell'\in
\wssu(w)$ and $\ga\in \Isom(X)$, we have 
$$
\ga\wssu(w)=\wssu(\ga w)
$$ 
and
$$
d_{\wssu(\ga w)}(\ga \ell,\ga \ell')= d_{\wssu(w)}(\ell,\ell')=
d_{\wsus (\iota w)}(\iota\ell,\iota\ell')\,.
$$ 
Furthermore, for every $s\in\RR$, we have 
$$
\flow{s}\wssu(w)=\wssu(\flow{s} w)
$$ 
and
for all $\ell,\ell'\in \wssu (w)$
\begin{equation}\label{eq:expandHamdist}
d_{\wssu (\flow sw)}(\flow s\ell,\flow s\ell')=
e^{\mp s}d_{\wssu(w)}(\ell,\ell')\;.
\end{equation}

If $X$ is an $\RR$-tree, for all $w\in\G_+X$ and $\ell,\ell'\in
\wss(w)$, if $\mathopen{[}s,+\infty\mathclose{[}$ is the maximal
interval on which $\ell$ and $\ell'$ agree, then
$d_{\wss(w)}(\ell,\ell')=e^s$.

The following lemma compares the distance in $\G X$ with
Hamenst\"adt's distance for two geodesic lines in the same strong
(un)stable leaf.

\blemm\label{lem:comparddHam} There exists a universal constant $c>0$
such that for all $w\in\G_\pm X$ and $\ell,\ell'\in \wssu (w)$, we
have
$$
d(\ell,\ell')\leq c\;d_{\wssu(w)}(\ell,\ell')\;\;\;{\rm and}\;\;\;
d(\pi(\ell),\pi(\ell'))\leq d_{\wssu(w)}(\ell,\ell')\;.
$$
\elemm

\dem We could refer to \cite[Lem.~3]{ParPau14ETDS} (see also
\cite[Lem.~2.4]{PauPolSha15}) for a proof of the first result. Note that
the distance on $\G X$ considered in loc.~cit.~is slightly different
from the one in this book, hence we give a full proof for the sake of
completeness. We assume that $w\in\G_+ X$, the proof when $w\in\G_- X$
is similar.

\begin{center}
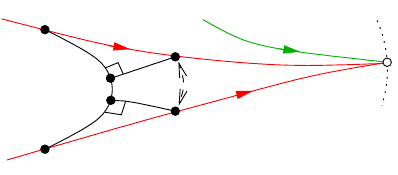
\end{center}

Let $\ell,\ell'\in \wss(w)$. We may assume that $\ell\neq \ell'$. By
the convexity properties of the distance in $X$, the map from $\RR$ to
$\RR$ defined by $t\mapsto d(\ell(t),\ell'(t))$ is decreasing, with
image $\mathopen{]}0,+\infty\mathclose{[}$. Let $S\in\RR$ be such that
$d(\ell(S),\ell'(S)) =1$.  For every $t\le S$, let $p$ and $p'$ be
the closest point projections of $\ell(S)$ and $\ell'(S)$ on the
geodesic segment $[\ell(t),\ell'(t)]$. We have $d(p,\ell(S)),
d(p',\ell'(S)) \leq 1$ by comparison. Hence, by convexity and the
triangle inequality,
\begin{align*}
d(\ell(t),\ell'(t))&\geq d(\ell(t),p)+d(p',\ell'(t))\\& 
\geq d(\ell(t),\ell(S))-1+ d(\ell'(t),\ell'(S))-1=2(S-t-1)\;.
\end{align*}
Thus by the definition of the Hamenst\"adt distance $d_{\wss(w)}$,
we have
\begin{equation}\label{eq:minodistss}
d_{\wss(w)}(\ell,\ell')\geq e^{S-1}\;.
\end{equation}

By the triangle inequality, if $t\leq S$, then
$$
d(\ell(t),\ell'(t))\leq
d(\ell(t),\ell(S))+d(\ell(S),\ell'(S))+ d(\ell'(S),\ell'(t))
=2(S-t)+1\;.
$$
Since $X$ is $\CAT(-1)$, if $t\geq S$, we have by comparison
$$
d(\ell(t),\ell'(t))\leq e^{S-t}\;\sinh d(\ell(S),\ell'(S))=
(\sinh 1)\,e^{S-t}\;.
$$ Therefore, by the definition of the distance $d$ on $\G X$ (see
Equation \eqref{eq:geodesicmetric}),
$$
d(\ell,\ell')\leq \int_{-\infty}^{S}(2(S-t)+1)\,e^{-2|t|}\;dt+(\sinh 1)
\int_{S}^{+\infty}e^{S-t}\,e^{-2|t|}\;dt=\operatorname{O}(e^{S})\;.
$$
The first inequality of Lemma \ref{lem:comparddHam} hence follows
from Equation \eqref{eq:minodistss}.

\medskip
The second one is proved in \cite[Lem.~2.4]{PauPolSha15}, and we again
only give a proof for the sake of completeness.

Let $x=\pi(\ell)$, $x'=\pi(\ell')$ and $\rho= d_{\wssu(w)}(\ell,
\ell')$.  Consider the ideal triangle $\Delta$ with vertices
$\ell_+,\ell'_+$ and $\ell_-=\ell'_-$ (see the picture below on the
left).  Let $p\in\ell(\RR)$, $p'\in\ell'(\RR)$ and $q\in
\mathopen{]}\ell_+, \ell'_+\mathclose{[}$ be the tangency points of
the unique triple of pairwise tangent horospheres centred at the
vertices of $\Delta$: $\beta_{\ell_-}(p,p')=0$, $\beta_{\ell_+}(p,q)
=0$ and $\beta_{\ell'_+}(p',q) =0$.  By the definition of the
Hamenst\"adt distance, we have $p=\ell(-\ln \rho)$.

\begin{center}   
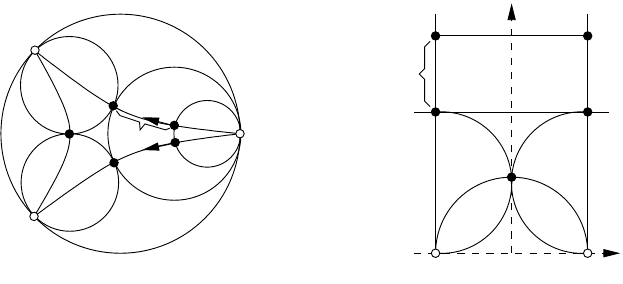
\end{center}

Consider the ideal triangle $\overline{\Delta}$ in the hyperbolic
upper half-plane $\HH^2_\RR$, with vertices $-\frac{1}{2}$,
$\frac{1}{2}$ and $\infty$ (see the above picture on the right). Let
$\overline{p}=(-\frac 12,1)$, $\overline{p}'=(\frac 12,1)$ and
$\overline{q}=(0,\frac{1}{2})$ be the pairwise tangency points of
horospheres centred at the vertices of $\overline{\Delta}$. Let
$\overline{x}$ and $\overline{x}'$ be the point at algebraic
(hyperbolic) distance $-\ln \rho$ from $\overline{p}$ and
$\overline{p}'$, respectively, on the upwards oriented vertical lines
through them. By comparison, we have $d(x,x')\leq
d(\overline{x},\overline{x}')\leq 1/e^{-\ln \rho} =\rho$. \cqfd

\medskip Let $\H$ be a horoball in $X$, centred at $\xi
\in\partial_\infty X$.  The strong stable leaves $\wss(w)$ are equal
for all geodesic rays $w$ starting at time $t=0$ from a point of
$\partial \H$ and converging to $\xi$. Using the homeomorphism
$\ell\mapsto \ell_-$ from $\wss(w)$ to $\partial_\infty X-\{\xi\}$,
Hamenst\"adt's distance on $\wss(w)$ defines a distance
$\gls{Hamenstadtdistanceinfinity}$ on $\partial_\infty X-\{\xi\}$ that
we also call {\em Hamenst\"adt's distance}%
\index{Hamenst\"adt's!distance}\index{distance!Hamenst\"adt's}. For
all $\ell,\ell'\in \wss(w)$, we have
$$
d_\H(\ell_-,\ell_+)=d_{\wss(w)}(\ell,\ell')\;,
$$
and for all $\eta,\eta'\in\partial_\infty X-\{\xi\}$, we have
\begin{equation}\label{eq:defidisthamenbord}
d_\H(\eta,\eta')=\lim_{t\ra+\infty} 
e^{\frac{1}{2}d(\ell_\eta(-t),\;\ell_{\eta'}(- t))-t}\;,
\end{equation}
where $\ell_\eta,\ell_{\eta'}$ are the geodesic lines starting from
$\eta,\eta'$ respectively, ending at $\xi$, and passing through the
boundary of $\H$ at time $t=0$. Note that for every $t\geq 0$, if
$\H[t]$ is the horoball contained in $\H$ whose boundary is
at distance $t$ from the boundary of $\H$, then we have
\begin{equation}\label{eq:shrinkhoroballdist}
d_{\H[t]}=e^{-t}\;d_\H\;.
\end{equation}

\bigskip 
Let $w\in\G_\pm X$ and $\eta'>0$. We define $\gls{ballhamenstadt}$ as
the set of $\ell\in\wssu(w)$ such that there exists an extension $\wh
w\in\G X$ of $w$ with $d_{\wssu(w)}(\ell,\wh w)<\eta'$.  In
particular, $B^\pm(w,\eta')$ contains all the extensions of $w$, and
is the union of the open balls centred at the extensions of $w$, of
radius $\eta'$, for Hamenst\"adt's distance on $\wssu(w)$.

\begin{center}
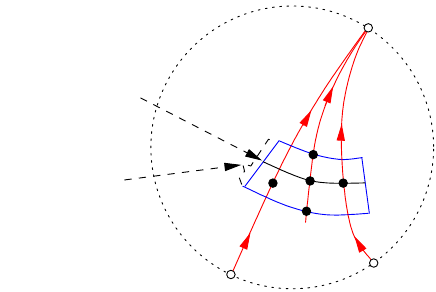
\end{center}

The union over $t\in\RR$ of the images under $\flow t$ of the strong
stable leaf of $w\in \G_+ X$ is the {\em stable
  leaf}\index{stable!leaf}\index{leaf!stable}
$$
\gls{Ws}=\bigcup_{t\in\RR}\flow t \wss(w)
$$ 
of $w$, which consists of the elements $\ell\in \G X$ with
$\ell_+=w_+$. Similarly,  
the {\em unstable leaf}\index{unstable!leaf}\index{leaf!unstable} of
$w\in\G_-X$ 
$$
\gls{Wu}=\bigcup_{t\in\RR}\flow t \wsu(w)\,,
$$ 
consists of the elements $\ell\in \G X$ with $\ell_-=w_-$. Note
that the (strong) (un)stable leaves are subsets of the space of
geodesic lines $\G X$. The (un)stable leaves are invariant under the
geodesic flow, and for all $w\in \G_\pm X$ and $\ga\in\Isom(X)$, we
have
$$
\iota\, \wosu(w) = \wous(\iota\, w)\;\;\;{\rm and}\;\;\;
\ga \wosu(w)=\wosu(\ga w)\;.
$$

\medskip The {\em unstable horosphere}\index{unstable!horosphere}%
\index{horosphere!unstable} $\gls{horosphereunstable}$ of $w\in\G_- X$
is the horosphere in $X$ centred at $w_{-}$ and passing through $\wh
w(0)$ for any extension $\wh w\in\G X$ of $w$ (see the picture above 
the definition of Hamenst\"adt's distance). The {\em stable
  horosphere}\index{stable!horosphere} \index{horosphere!stable}
$\gls{horospherestable}$ of $w\in\G_+ X$ is the horosphere in $X$
centred at $w_{+}$ and passing through $\wh w(0)$ for any extension
$\wh w\in\G X$ of $w$.  These horospheres $H_{\pm}(w)$ do not depend
on the chosen extensions $\wh w$ of $w\in\G_\pm X$. The {\em unstable
  horoball}\index{unstable!horoball}\index{horoball!unstable}
$\gls{horoballunstable}$ of $w\in\G_- X$ and {\em stable
  horoball}\index{stable!horoball}\index{horoball!stable}
$\gls{horoballstable}$ of $w\in\G_+ X$ are the horoballs bounded by
these horospheres. Note that
\begin{equation}\label{eq:footpowssuhpm}
  \pi(\wssu(w))=H_\pm(w)
\end{equation}
for every $w\in\G_\pm X$, and that $w(0)$ belongs to $H_\pm(w)$ if and
only if $w$ is isometric at least on $\pm[0,+\infty[$.

\section{Normal bundles and dynamical neighbourhoods}
\label{subsect:nbhd}

In this Section, adapting \cite[\S 2.2]{ParPau16ETDS} to the present
context, we define spaces of geodesic rays that generalise the unit
normal bundles of submanifolds of negatively curved Riemannian
manifolds. When $X$ is a manifold, these normal bundles are
submanifolds of the unit tangent bundle of $X$, which identifies with
$\G X$. In general and in particular in trees, it is essential to use
geodesic rays to define normal bundles, and not geodesic lines.

\medskip 
Let $D$ be a nonempty {\em proper}\index{proper}\footnote{that is, different
from $X$} closed convex subset in $X$. We denote by $\partial D$ its
boundary in $X$ and by $\partial_\infty D$ its set of points at
infinity.  Let
$$
\gls{closestpointmap}:X\cup(\partial_\infty X-\partial_\infty D)\to D
$$ 
be the (continuous) {\em closest point map to $D$}\index{closest
  point map}, defined on $\xi\in\partial_\infty X -\partial_\infty D$
by setting $P_D(\xi)$ to be the unique point in $D$ that minimises the
function $y\mapsto \beta_\xi(y,x_0)$ from $D$ to $\RR$. The {\em outer
  unit normal bundle}\index{outer unit normal bundle}%
\index{unit!normal bundle!outer} $\normalout D$ of (the boundary of)
$D$ is
$$
\gls{outerunitnormalbundle}= 
\{\rho\in\G_{+,\,0}X\;:\; P_{D}(\rho_+) =\rho(0)\}\,.
$$ 
The {\em inner unit normal bundle}\index{inner unit normal
  bundle}\index{unit!normal bundle!inner} $\normalin D$ of (the
boundary of) $D$ is
$$
\gls{innerunitnormalbundle}=\iota\normalout
D= \{\rho\in\G_{-,\,0}X\;:\; P_{D}(\rho_-) =\rho(0)\}.
$$ 
Note that $\normalout D$ and $\normalin D$ are spaces of geodesic
rays.  If $X$ is a smooth manifold, then these spaces have a natural
identification with subsets of $\G X$ because every geodesic ray is
the restriction of a unique geodesic line. But this does not hold in
general.

\begin{center}
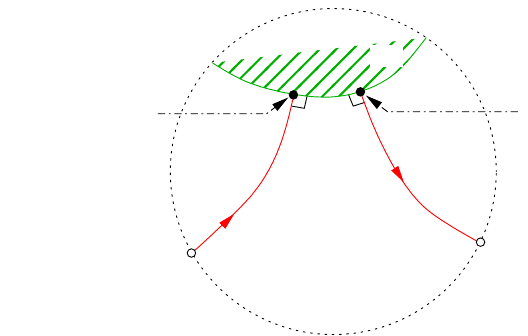
\end{center}

\brema As $X$ is assumed to be proper with extendible geodesics, we
have 
$$
\pi(\normalpm D)=\partial D\;.
$$ 
To see this, let $x\in\partial D$ and let $(x_k)_{k\in\NN}$ be a
sequence of points in the complement of $D$ converging to $x$. For all
$k\in\NN$, let $\rho_k\in\normalout D$ be a geodesic ray with
$\rho_k(0)=P_D(x_k)$ and such that the image of $\rho_k$ contains
$x_k$. As the closest point map does not increase distances, the
sequence $(P_D(x_k))_{k\in\NN}$ converges to $x$. Since $X$ is proper,
the space $\partial_\infty X$ is compact and the sequence
$((\rho_k)_+)_{k\in\NN}$ has a subsequence that converges to a point
$\xi\in\partial_\infty X$.  The claim follows from the continuity of
the closest point map.

The possible failure of this equality when $X$ is not proper is easy
to see. For example, let $X$ be the $\RR$-tree constructed by starting
with the Euclidean line $D=\RR$ and attaching a copy of the halfline
$\mathopen{[}0,+\infty\mathclose{[}$ at each $x\in D$ such that $x>0$. 
Then $0\in\partial D-\pi(\normalpm D)$.  
\erema

The restriction of the endpoint map $\rho\mapsto \rho_\pm$ to
$\normalpm D$ is a homeomorphism to its image $\partial_\infty
X-\partial_\infty D$. We denote its inverse map by
$\gls{closestpointmaptonormal}$ (see the above picture): for every
$\rho\in\normalpm D$, we have
$$
\rho=P^\pm_D(\rho_\pm)\;.
$$
Note that $P_D=\pi\circ P^\pm_D$. For every isometry $\ga$ of $X$, we
have $\normalpm (\ga D)= \ga\,\normalpm D$ and $P^\pm_{\ga
  D}\circ\ga=\ga\circ P^\pm_D$. In particular, $\normalpm D$ is
invariant under the isometries of $X$ that preserve $D$.

For every $w\in\G_\pm X$, we have a canonical homeomorphism
$$
\gls{homeoinftystableaf}: \wssu(w)\to\normalmp H\!B_\pm(w)\;,
$$ 
that associates to each geodesic line $\ell\in\wssu(w)$ the unique
geodesic ray $\rho \in \normalmp H\!B_\pm(w)$ such that
$\ell_\mp=\rho_\mp$, or, equivalently, such that $\ell(t)=\rho(t)$ for
every $t\in\RR$ with $\mp t>0$.  It is easy to check that $N^\pm_{\ga
  w}\circ\ga= \ga\circ N^\pm_w$ for every $\ga\in \Isom(X)$.

\medskip
We define 
\begin{equation}\label{eq:defiUCp}
\gls{calUpmsubD}=\{\ell\in \G X\;:\; \ell_\pm\notin\partial_\infty D\}\;.
\end{equation}
Note that $\U^\pm_D$ is an open subset of $\G X$, invariant under the
geodesic flow.  We have $\U^\pm_{\ga D}=\ga \U^\pm_D$ for every
isometry $\ga$ of $X$ and, in particular, $\U^\pm_D$ is invariant
under the isometries of $X$ preserving $D$. Define 
$$
\gls{fibrationpm}:\U^\pm_D \ra \normalpm D
$$ 
as the composition of the continuous endpoint map $\ell\mapsto
\ell_\pm$ from $\U^\pm_D$ onto $\partial_\infty X -\partial_\infty D$
and the homeomorphism $P^\pm_{D}$ from $\partial_\infty
X-\partial_\infty D$ to $\normalpm D$ (see the picture below). The
continuous map $f^\pm_D$ takes $\ell\in\U^\pm_D$ to the unique element
$\rho\in \normalpm D$ such that $\rho_\pm=\ell_ \pm$.  The fiber of
$\rho\in \normalout D$ for $f^+_D$ is exactly the stable leaf
$\ws(\rho)$, and the fiber of $\rho\in \normalin D$ for $f^-_D$ is the
unstable leaf $\wu(\rho)$.  For all $\ga\in \Isom(X)$ and $t\in\RR$,
we have
\begin{equation}\label{eq:equivfibrationf}
f^\pm_{\ga D} \circ\ga=\ga\circ f^\pm_D\;\;\;{\rm and}\;\;\;
f^\pm_D\circ \flow t=f^\pm_D\,.
\end{equation}

\begin{center}
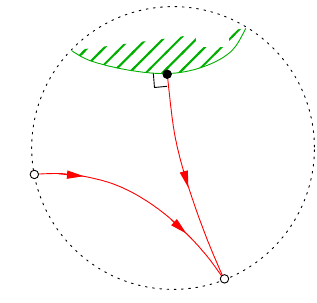
\end{center}

\medskip
Let $w\in\G_\pm X$ and $\eta,\eta'>0$. We define (see the picture
after Equation \eqref{eq:shrinkhoroballdist}) the {\em dynamical 
$(\eta,\eta')$-neighbourhood}%
\index{dynamical!neighbourhood!of a point} of $w$ by
\begin{equation}\label{eq:defidynamicalneighbourhood}
\gls{dynamicalneighbourhood} =
\bigcup_{s\in\mathopen{]}-\eta,\,\eta\,\mathclose{[}}\flow s
B^\pm(w, \eta')\;.
\end{equation}

\bexem If $X$ is an $\RR$-tree, $w\in \G_+X$ and $\eta<\ln \eta'$,
then $V^+_{w,\,\eta,\,\eta'}$, which is the set of $\flow{s}\ell$
where $s\in\mathopen{]}-\eta,\eta\mathclose{[}$ and $\ell\in\G X$ is
such that there exists an extension $\wh w$ of $w$ with $\inf 
\{t\in\RR\;:\; \ell(t)=\wh w (t)\}\leq \ln \eta'$, is as in the 
following picture.

\begin{center}
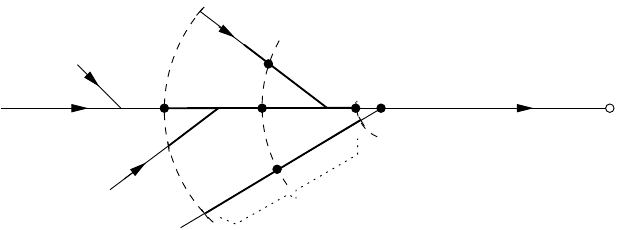
\end{center}
\eexem

\bigskip
Clearly, $B^\pm(w,\eta')=\iota B^\mp(\iota w,\eta')$, and hence we
have $V^\pm_{w,\,\eta,\,\eta'}=\iota V^\mp_{\iota
  w,\,\eta,\,\eta'}$. Furthermore, for every $s\in\RR$,
\begin{equation}\label{eq:flotVetaeta}
\flow s B^\pm(w, \eta') = B^\pm(\flow s w, e^{\mp s} \eta')\,, 
\;\;{\rm hence}\;\; \flow s
V^\pm_{w, \,\eta,\,\eta'} = V^\pm_{\flow s w, \,\eta,\,e^{\mp s} \eta'}\,.
\end{equation} 
For every $\ga\in\Isom(X)$, we have $\ga B^\pm(w, \eta')= B^\pm(\ga w,
\eta') $ and $\ga V^\pm_{w,\,\eta,\,\eta'} = V^\pm_{\ga w,\,\eta,\,\eta'} \,$.  
The map from $\mathopen{]}-\eta, \eta\mathclose{[} \times B^\pm(w,
\eta')$ to $V^\pm_{w,\,\eta,\,\eta'}$ defined by $(s,\ell')
\mapsto \flow s \ell'$ is a homeomorphism.

For all subsets $\Omega^-$ of $\G_+X$ 
and $\Omega^+$ of $\G_-X$, let
\begin{equation}\label{eq:defcalV}
\gls{dynamicalneighbourhoodsubset}=
\bigcup_{w\in\Omega^\mp}  V^\pm_{w,\,\eta,\,\eta'}\,,
\end{equation}
that we call the {\em dynamical neighbourhoods}%
\index{dynamical!neighbourhood} of $\Omega^\mp$. Note that they are
subsets of $\G X$, not of $\G_\pm X$.  The families
$(\V^\pm_{\eta,\,\eta'} (\Omega^{\mp}))_{\eta,\eta'>0}$ are
nondecreasing in $\eta$ and in $\eta'$.  For every $\ga\in \Isom(X)$,
we have $\ga \V^\pm_{\eta,\,\eta'} (\Omega^{\mp})=\V^\pm_{\eta,\,\eta'} 
(\ga \Omega^{\mp})$ and for every $t\geq 0$, we have
\begin{equation}\label{eq:flowbehavdynaneigh}
\flow {\pm t}\V^\pm_{\eta,\,\eta'}(\Omega^{\mp})
=\V^\pm_{\eta,\,e^{-t}\eta'}(\flow{\pm t}\Omega^{\mp})\;. 
\end{equation}

Note that
$$
\bigcup_{\eta\,,\;\eta'>0}\V^\pm_{\eta,\,\eta'}(\partial^1_{\pm} D)
=\U^\pm_D\;,
$$ 
and that $\bigcap_{\eta,\,\eta'>0}\V^\pm_{\eta,\,\eta'}
(\partial^1_{\pm} D)$ is the set of all extensions in $\G X$ of the
elements of $\partial^1_{\pm} D$. Assume that $\Omega^\mp$ is a subset
of $\partial^1_\pm D$.  The restriction of $f^\pm_D$ to
$\V^\pm_{\eta,\,\eta'}(\Omega^{\mp})$ is a continuous map onto
$\Omega^{\mp}$, with fiber over $w\in \Omega^{\mp}$ the open subset
$V^\pm_{w,\,\eta,\,\eta'}$ of $\wosu(w)$.

\bigskip 
We will need the following elementary lemma in Section
\ref{subsec:rateequidtrees}.

\blemm \label{lem:estimdellw}
There exists a universal constant $c'>0$ such that for every
$w\in\G_+ X$ which is isometric on $[s_w,+\infty[$ and every
    $\ell\in V^+_{w,\,\eta,\,\eta'}$, we have
$$
d(\ell,w)\leq c'(\eta+\eta'+e^{s_w})\;.
$$
\elemm

\dem By Equation \eqref{eq:defidynamicalneighbourhood} and by the
definition of $B^+(w,\eta')$ in Section \ref{subsec:lines}, there
exist $s\in\mathopen{]}-\eta,+\eta\mathclose{[}$ and an extension 
$\wh{w}\in\G X$ of $w$ such that $\flow{s}\ell\in \wss(w)$ and
$d_{W^{+}(w)}(\flow{s}\ell, \wh{w})\leq \eta'$. By Equation
\eqref{eq:disttranslatgeod}, we have $d(\ell, \flow{s}\ell)\leq
|s|\leq \eta$. By Lemma \ref{lem:comparddHam}, we have
$d(\flow{s}\ell, \wh{w})\leq c\; d_{W^{+}(w)}(\flow{s}\ell, \wh{w})
\leq c\;\eta'$. By the definition of the distance on $\gengeod X$ (see
Equation \eqref{eq:geodesicmetric}), we have
$$
d(\wh{w},w)\leq \int_{-\infty}^{s_w} |s_w-t|\;e^{-2|t|}\;dt =\bigO(e^{s_w})\;.
$$ 
Therefore the result follows from the triangle inequality
$$
d(\ell,w)\leq d(\ell,\flow{s}\ell)+d(\flow{s}\ell, \wh{w})+d(\wh{w},w)
\;.\;\;\;\Box
$$

\section{Creating common perpendiculars}
\label{subsec:creatcommonperp}

Let $D^-$ and $D^+$ be two nonempty proper closed convex subsets of
$X$, where $X$ is as in the beginning of Section
\ref{subsec:catmoinsun}.  A geodesic arc $\alpha:\mathopen{[}0,
  T\mathclose{]}\to X$, where $T>0$, is a {\em common perpendicular of
  length}\index{common perpendicular}%
\index{length!of common perpendicular} $T$ {\em from} $D^-$ {\em to}
$D^+$ if there exists $w^\mp\in\partial^1_\pm D^\mp$ such that
$w^-|_{[0,\,T]}= \flow{-T}w^+|_{[0,\,T]} =\alpha$. Since $X$ is
$\CAT(-1)$, this geodesic arc $\alpha$ is the unique shortest geodesic
segment from a point of $D^-$ to a point of $D^+$. There exists a common
perpendicular from $D^-$ to $D^+$ if and only if the closures of $D^-$
and $D^+$ in $X\cup \partial_\infty X$ are disjoint.  When $X$ is an
$\RR$-tree, then two closed subtrees of $X$ have a common
perpendicular if and only if they are nonempty and disjoint.

\medskip One of the aims of this book is to count orbits of common
perpendiculars between two equivariant families of closed convex
subsets of $X$. The crucial remark is that two nonempty proper closed
convex subsets $D^-$ and $D^+$ of $X$ have a common perpendicular
$\alpha$ of length a given $T> 0$ if and only if the subsets
$\flow{T/2}\normalout D^-\big|_{[-\frac T2,\,\frac T2]}$ and
$\flow{-T/2}\normalin D^+\big|_{[-\frac T2,\frac T2]}$ of $\gengeod X$
intersect. This intersection then consists of the common perpendicular
from $D^-$ to $D^+$ reparametrised by $[-\frac T2,\frac T2]$. As a
controlled perturbation of this remark, we now give an effective
creation result of common perpendiculars in $\RR$-trees. It has a
version satisfied for $X$ in the generality of Section
\ref{subsec:catmoinsun}, see the end of this Section.

\blemm \label{lem:creationperp} Assume that $X$ is an $\RR$-tree. For
all $R>1$, $\eta\in\mathopen{]}0,1\mathclose{]}$ and $t\ge2\ln R+4$,
for all nonempty closed connected subsets $D^-,D^+$ in $X$, and for
every geodesic line $\ell\in\flow{t/2} \V^+_{\eta,\,R}(\normalout D^-)
\cap \flow{-t/2} \V^-_{\eta,\,R} (\normalin D^+)$, there exist $s\in
\mathopen{]}-2\eta, 2\eta \mathclose{[}$ and a common perpendicular
$\wt c$ from $D^-$ to $D^+$ such that
\begin{itemize}
\item the length of $\wt c$ is $t+s$,
\item the endpoint of $\wt c$ in $D^\mp$ is $w^\mp(0)$ where
  $w^\mp=f^\pm_{D^\mp}(\ell)$,
\item the footpoint $\ell(0)$ of $\ell$ lies on $\wt c$, and 
$$
\max\Big\{\;d(w^-\big(\frac t2\big),\ell(0)),\; 
d(w^+\big(-\frac{t}{2}\big),\ell(0))\;\Big\}\leq \eta \;.
$$
\end{itemize}  
\elemm

\begin{center}
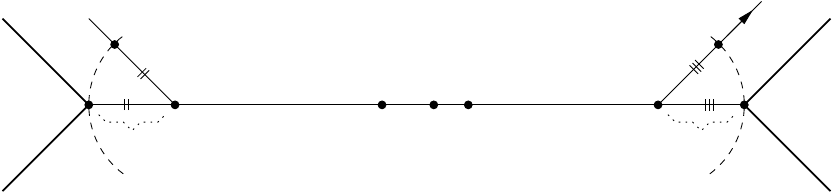
\end{center}

\dem Let $R,\eta,t,D^\pm,\ell$ be as in the statement.  By the
definition of the sets $\V^\mp_{\eta,\,R}(\partial^1_\mp D^\pm)$,
there exist geodesic rays $w^\pm\in\normalmp D^\pm$, geodesic lines
$\wh w^\pm\in\G X$ extending $w^\pm$, and $s^\pm\in\mathopen{]}-\eta,
+\eta\mathclose{[}$, such that $\ell_\pm=(w^\mp)_\pm$ and
$$
 d_{\wssu(w^\mp)}(\flow{\mp\frac{t}{2}\mp s^\mp}\ell,\wh w^\mp)\leq R\;.
$$
Let $x_\pm$ be the closest point to $w^\pm(0)$ on $\ell$. By the
definition of Hamenst\"adt's distances, we have 
$$
d(w^\pm(0),x^\pm)=
d(\ell(\pm\frac t2\pm s^\pm),x^\pm)\le\ln R\,,
$$ 
and in particular $x_\pm=w^\pm(0)$ if and only if $\ell(\pm\frac t2\pm
s^\pm)= w^\pm(0)$. As $t\ge 2\ln R+4$ and $|s^\pm|\leq 2\eta\leq 2$,
the points $\ell(-\frac t2 - s^-)$, $x^-$, $\ell(0)$, $x^+$,
$\ell(\frac t2+s^+)$ are in this order on $\ell$.  In particular, the
segment $[w^-(0),x^-]\cup[x^-,x^+]\cup [x^+, w^+(0)]$ is a nontrivial
geodesic segment from a point of $D^-$ to a point of $D^+$ that meets
$D^\mp$ only at an endpoint.  Hence, $D^-$ and $D^+$ are disjoint, and
$[w^-(0),w^+(0)]$ is the image of the common perpendicular from $D^-$
to $D^+$.

Let $s=s_-+s_+$. The length of $\wt c$ is $(\frac t2+s^+)-(-\frac
t2-s^-)=t+s$. The point $\ell(0)$ lies on $\wt c$, we have $w^\mp=
f^\pm_{D^\mp}(\ell)$ and the endpoints of $\wt c$ are $w^\pm(0)$.
Furthermore,
$$
d(w^\mp(\pm\frac t2),\ell(0))=\big|\;d(w^\mp(\pm\frac t2),w^\mp(0))-
d(\ell(0),\ell(\mp\frac t2\mp s^\mp))\;\big|
=|s^\mp|\le\eta\,.\quad\quad\Box
$$

When $X$ is as in the beginning of Section \ref{subsec:catmoinsun},
the statement and the proof of the following analog of Lemma
\ref{lem:creationperp} is slightly more technical. We refer to
\cite[Lem.~7]{ParPau16ETDS} for a proof in the Riemannian case, and we
leave the extension to the reader, since we will not need it in this
book.

\blemm \label{lem:creationperpRiem} Let $X$ be as in the beginning of
Section \ref{subsec:catmoinsun}. For every $R>0$, there exist
$t_0,c_0>0$ such that for all $\eta\in\;]0,1]$ and all $t\in \;
[t_0,+\infty[\,$, for all nonempty closed convex subsets $D^-,D^+$ in
$X$, and for all $w\in\flow{t/2}\V^+_{\eta,\,R}(\normalout D^-)
\cap \flow{-t/2} \V^-_{\eta,\,R}(\normalin D^+)$, there exist
$s\in\; ]-2\eta,2\eta[$ and a common perpendicular $\wt c$ from
$D^-$ to $D^+$ such that
\begin{itemize}
\item the length of $\wt c$ is contained in
  $[t+s-c_0\,e^{-\frac{t}{2}},t+s+c_0\,e^{-\frac{t}{2}}]$,
\item if $w^\mp=f^\pm_{D^\mp}(w)$ and if $p^\pm$ is the endpoint of
  $\wt c$ in $D^\pm$, then $d(\pi(w^\pm),p^\pm)\leq
  c_0\,e^{-\frac{t}{2}}$,
\item the basepoint $\pi(w)$ of $w$ is at distance at most
$c_0\,e^{-\frac{t}{2}}$ from a point of $\wt c$, and 
$$
\max\{\;d(\pi(\flow{\frac{t}{2}} w^-),\pi(w)),\;
d(\pi(\flow{-\frac{t}{2}} w^+),\pi(w))\;\}\leq 
\eta + c_0\,e^{-\frac{t}{2}}\;.\quad\quad\Box
$$
\end{itemize}  
\elemm

\section{Metric and simplicial trees, and graphs of groups}
\label{subsec:trees}

Metric and simplicial trees and graphs of groups are important
examples throughout this book. In this Section, we
recall the definitions and basic properties of these objects.  

Using Serre's definitions in \cite[\S 2.1]{Serre83}, a {\em
  graph}\index{graph} $\XX$ is the data consisting of two sets
$\gls{setofvertices}$ and $\gls{setofedges}$, called the {\em set of
  vertices} and the {\em set of edges} of $\XX$, of two maps
$\gls{origin}, \gls{terminus}:E\XX\ra V\XX$ and of a fixed point free
involution $e\mapsto \overline{e}$ of $E\XX$, such that
$t(\overline{e})=o(e)$ for every $e\in E\XX$. The elements $o(e)$,
$t(e)$ and $\overline{e}$ are called the {\em initial
  vertex}\index{vertex!initial}, the {\em terminal
  vertex}\index{vertex!terminal} and the {\em opposite
  edge}\index{opposite edge}\index{edge!opposite} of an edge $e\in
E\XX$.  The quotient of $E\XX$ by the involution $e\mapsto
\overline{e}$ is called the set of {\em nonoriented
  edges}\index{edge!nonoriented} of $\XX$. Recall that a connected
graph is {\em bipartite}\index{bipartite}\index{graph!bipartite} if it
is endowed with a partition of its set of vertices into two nonempty
subsets such that any two elements of either subset are not related by
an edge.

The {\em degree}\index{degree} of a vertex $x\in V\XX$ is the
cardinality of the set $\{e\in E\XX:o(e)=x\}$.  For all $j,k\in\NN$, a
graph $\XX$ is $k$-{\em regular}\index{regular}\index{tree!regular} if
the degree of each vertex $x\in V\XX$ is $k$, and it is $(j,k)$-{\em
  biregular}\index{biregular}\index{tree!biregular} if it is bipartite
with the elements of the partition of its vertices into two subsets
having degree $j$ and $k$ respectively.

A {\em metric graph}\index{metric!tree}\index{tree!metric}
$(\XX,\lambda)$ is a pair consisting of a graph $\XX$ and a map
$\lambda:E\XX\ra \mathopen{]}0, +\infty \mathclose{[}$ with a positive
lower bound\footnote{This assumption, though not necessary at this
  stage, will be used repeatedly in this book, hence we prefer to add
  it to the definition.} such that $\lambda(\overline{e})=\lambda(e)$,
called its {\em edge length map}\index{edge!length map}.  A {\em
  simplicial graph}\index{simplicial tree}\index{tree!simplicial}
$\XX$ is a metric graph whose edge length map is constant equal to
$1$.

The {\em topological realisation} of a graph $\XX$ is the topological
space obtained from the family $(I_e)_{e\in E\XX}$ of closed unit
intervals $I_e$ for every $e\in E\XX$ by the finest equivalence
relation that identifies intervals corresponding to an edge and its
opposite edge by the map $t\mapsto 1-t$ and identifies the origins of
the intervals $I_{e_1}$ and $I_{e_2}$ if and only if $o(e_1)=o(e_2)$,
see \cite[Sect.~2.1]{Serre83}.

The {\em geometric realisation}\index{geometric realisation} of a
metric tree $(\XX,\lambda)$ is the topological realisation of $\XX$
endowed with the maximal geodesic metric that gives length
$\lambda(e)$ to the topological realisation of each edge $e\in E\XX$,
and we denote it by $X=\gls{geomrealiz}$. We identify $V\XX$ with its
image in $X$. The metric space $X$ determines $(\XX,\lambda)$ up to
subdivisions of edges, hence we will often not make a strict
distinction between $X$ and $(\XX,\lambda)$. In particular, we will
refer to convex subsets of $(\XX,\lambda)$ as convex subsets of
$X$, etc.

If $\XX$ is a tree, the metric space $X$ is an $\RR$-tree, hence it is
a $\CAT(-1)$ space. Since $\lambda$ is bounded from below by a
positive constant, the $\RR$-tree $X$ is geodesically complete if and
only if $\XX$ is not reduced to one vertex and has no {\em terminal
  vertex}\index{terminal vertex} (that is, no vertex of degree $1$).

We will denote by $\gls{automorphismgroupmetrictree}$, and
$\gls{automorphismgroupsimptree}$ in the simplicial case, the group of
edge-preserving isometries of $X$ that have no inversions.\footnote{An
  automorphism $g$ of a graph has an {\em inversion}\index{inversion}
  if there exists an edge $e$ of the graph such that
  $ge=\overline{e}$. The assumption that the elements of $\Aut(\XX)$
  have no inversion ensures that, for every subgroup $\Ga'$ of
  $\Aut(\XX)$, the quotient map $\Ga'\bs E\XX\ra \Ga'\bs E\XX$ of
  $e\mapsto \overline{e}$ is still a fixed point free involution, so
  that with the quotient maps $\Ga'\bs E\XX\ra \Ga'\bs V\XX$ defined by
  $o$ and $t$, we do have a quotient graph structure $\Ga' \bs \XX$.}
Since the edge length map has a positive lower bound, the metric space
$X$ is proper if and only if $\XX$ is locally finite. In this case,
the nonelementary discrete subgroups $\Ga$ of isometries of $X$ we
will consider will always be edge-preserving and without inversion. If
$\Ga$ is a subgroup of $\Aut(\XX,\lambda)$, we will again denote by
$\lambda:\Ga\bs E\XX\ra \mathopen{]}0, +\infty \mathclose{[}$ the map
    induced by $\lambda:E\XX\ra \mathopen{]}0, +\infty \mathclose{[}$.

A locally finite metric tree $(\XX',\lambda)$ is {\em
  uniform}\index{uniform tree}\index{tree!metric!uniform} if there
exists some discrete subgroup $\Ga'$ of $\Aut(\XX',\lambda)$ such that
$\Ga'\bs\XX'$ is a finite graph. See \cite{BasKul90,BasLub01} for
characterisations of this property in the case of simplicial trees.

\subsection*{Discrete time geodesic flow on trees}

Let $\XX$ be a locally finite simplicial tree. The {\em space of
  generalised discrete geodesic 
  lines}\index{generalised geodesic!line!discrete}%
\index{geodesic!line!generalised discrete} of $\XX$ is the locally
compact space $\gls{gengeoddiscrete}$ of $1$-Lipschitz mappings $w$
from $\RR$ to the geometric realisation $X=|\XX|_1$ which are
isometric on a closed interval with endpoints in $\ZZ\cup\{-\infty,
+\infty\}$ and locally constant outside it, such that $w(0)\in V\XX$
(or equivalently $w(\ZZ)\subset V\XX$). Note that $\gengeod\XX$ is
hence a proper subset of $\gengeod X$, unless $\XX$ is reduced to one
vertex.

By restriction to $\gengeod \XX$, or intersection with $\gengeod \XX$,
of the objects defined in Sections \ref{subsec:lines} and
\ref{subsect:nbhd} for $\gengeod X$, we define the distance $d$ on
$\gengeod\XX$, the subspaces $\G_\pm\XX$, $\G \XX$, $\G_{\pm,0} \XX$,
the strong stable/unstable leaves $\wssu(w)$ of $w\in\G_\pm\XX$ and
their Hamenst\"adt distances $d_{\wssu(w)}$, the stable/unstable
leaves $\wosu(w)$ of $w\in\G_\pm\XX$, the outer and inner unit normal
bundles $\normalpm \DD$ of a nonempty proper simplicial subtree $\DD$
of $\XX$, the dynamical neighbourhoods $\V^\pm_{\eta,\,\eta'}
(\Omega^\mp)$ of subsets $\Omega^\mp$ of $\partial^1_\pm \DD$ as well
as the fibrations 
$$
f^\pm_\DD:\U^\pm_\DD= \{\ell\in \G \XX\;:\; 
\ell_\pm\notin\partial_\infty |\DD|_1\}\ra \normalpm \DD\;,
$$ 
whose fiber over $\rho\in\normalpm \DD$ is $\wosu(\rho)$. Note that some
definitions actually simplify when considering generalised discrete
geodesic lines. For instance, for all $w\in\G_\pm\XX$, $\eta'>0$ and
$\eta\in\,]0,1[\,$, the dynamical neighbourhood
$V^\pm_{w,\,\eta,\,\eta'}$ is equal to $B^\pm(w,\eta')$, and is hence
independent of $\eta\in\,]0,1[\,$.

Besides the map $\pi:\G \XX\ra V\XX$ defined as in the continuous case
by $\ell\mapsto \ell(0)$, we have another natural map $\gls{Tpi}:\G
\XX\ra E\XX$, which associates to $\ell$ the edge $e$ with
$o(e)=\ell(0)$ and $t(e)=\ell(1)$. This map is equivariant under the
group of automorphisms (without inversions) $\Aut(\XX)$ of $\XX$, and
we also denote by $T\pi:\Ga\bs\G \XX\ra \Ga\bs E\XX$ its quotient map,
for every subgroup $\Ga$ of $\Aut(\XX)$.

If $\XX$ has no terminal vertex, for every $e\in E\XX$, let
$$
\gls{boundarytowardse}=\{\ell_+\;:\;\ell\in\G \XX,\;T\pi(\ell)=e\}
$$ 
be the set of points at infinity of the geodesic rays whose initial
(oriented) edge is $e$.

Given $x_0\in V\XX$, the {\em discrete Hopf
  parametrisation}\index{Hopf parametrisation!discrete} now identifies
$\G \XX$ with $\partial^2_\infty X\times\ZZ$ by the map $\ell\mapsto
(\ell_-,\ell_+,t)$ where $t\in\ZZ$ is the signed distance from the
closest vertex to the basepoint $x_0$ on the geodesic line $\ell$ to
the vertex $\ell(0)$.

The {\em discrete  time geodesic flow}\index{geodesic!flow!discrete time}
$\gls{geodesicflowdis}$ on $\gengeod\XX$ is the
one-(discrete-)parameter group of homeomorphisms of $\gengeod\XX$
consisting of (the restriction to $\gengeod\XX$ of) the integral time
maps of the continuous time geodesic flow of the geometric realisation
of $\XX$: we have $\flow t w:s \mapsto w(s+t)$ for all $w\in\gengeod
\XX$ and $t\in\ZZ$.

\subsection*{Crossratios of ends of trees}

Let $\XX$ be a locally finite simplicial tree, with geometric
realisation $X=|\XX|_1$. Recall\footnote{See \cite{Paulin96}, as well as
  \cite{Otal92} in the case of Riemannian manifolds, and note that the
  convention on the order varies in the literature.} that if
$(\xi_1,\xi_2,\xi_3, \xi_4)$ is an ordered quadruple of pairwise
distinct points in $\partial_\infty X$, then their (logarithmic) {\em
  crossratio}\index{crossratio} is
\begin{equation}\label{eq:defcrossratio}
\gls{crossratiotree}=\lim_{x_i\ra\xi_i,\,x_i\in V\XX}\;\;\frac{1}{2}
\big(d(x_1,x_4)-d(x_4,x_3)+d(x_3,x_2)-d(x_2,x_1)\big)\;.
\end{equation}
A similar definition is valid for general $\CAT(-1)$-spaces, but we
will only need the case of simplicial trees in this book.

If $x$ and $y$ are the closest points on the geodesic line
$\mathopen{]}\xi_1,\xi_3\mathclose{[}$ to $\xi_2$ and $\xi_4$
respectively, then
$$
\ldbrack\xi_1,\xi_2,\xi_3,\xi_4\rdbrack=d(x,y)
$$
if $\xi_1,x, y,\xi_3$ are in this order on $\mathopen{]}\xi_1, \xi_3
\mathclose{[}$ and
$$
\ldbrack\xi_1,\xi_2,\xi_3,\xi_4\rdbrack=-d(x,y)
$$ 
otherwise. In particular, $\ldbrack\xi_1,\xi_2,\xi_3,\xi_4\rdbrack=0$
if the geodesic lines $\mathopen{]}\xi_1,\xi_3\mathclose{[}$ and
$\mathopen{]}\xi_2, \xi_4\mathclose{[}$ are disjoint.

\begin{center}
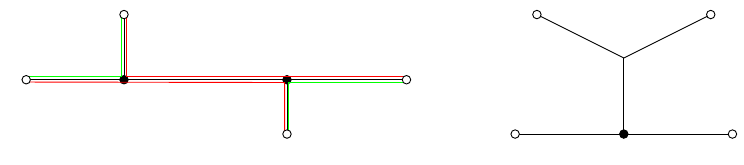
\end{center}

We have the following properties:

$\bullet$~ $\ldbrack\xi_1,\xi_2,\xi_3,\xi_4\rdbrack=
\ldbrack\xi_4,\xi_3,\xi_2,\xi_1\rdbrack
=-\ldbrack\xi_3,\xi_2,\xi_1,\xi_4\rdbrack$.

$\bullet$~ the crossratio is continuous, and even locally constant on
the space of pairwise distinct quadruples of elements of
$\partial_\infty X$,

$\bullet$~ if $\ga$ is a loxodromic element of $\Aut(\XX)$, with
repulsive and attractive fixed points $\ga_-$ and $\ga_+$ in
$\partial_\infty X$ respectively, then for every $\xi\in
\partial_\infty X- \{\ga_-,\ga_+\}$, the translation length of $\ga$
satisfies
$$
  \lambda(\ga)= \ldbrack\ga_-,\xi,\ga\xi,\ga_+\rdbrack\;.
$$

\subsection*{Bass-Serre's graphs of groups}

 Recall (see for instance \cite{Serre83, BasLub01}) that a {\em 
  graph of groups}\index{graph!of groups}
 $(\YY,G_*)$ consists of
\begin{enumerate}
\item[$\bullet$] a graph $\YY$, which is connected unless otherwise
  stated,
\item[$\bullet$] a group $G_v$ for every vertex $v\in V\YY$,
\item[$\bullet$] a group $G_e$ for every edge $e\in E\YY$ such that
  $G_e=G_{\overline{e}}$,
\item[$\bullet$] an injective group morphism $\rho_e:G_e\ra G_{t(e)}$
  for every edge $e\in V\YY$.
\end{enumerate}
We will identify $G_e$ with its image in $G_{t(e)}$ by $\rho_e$,
unless the meaning is not clear (which might be the case for instance
if $o(e)=t(e)$). We refer to op.~cit.\footnote{though see Example
  \ref{exem:quotgraphgroup} for the main example encountered in this
  book}~for the definition of the {\em Bass-Serre tree} $T(\YY,G_*)$
of $(\YY,G_*)$, of its {\em fundamental group} $\pi_1(\YY,G_*)$ when a
basepoint in $V\YY$ is chosen, and of the simplicial action of
$\pi_1(\YY,G_*)$ on $T(\YY,G_*)$. Note that the fundamental group of
$(\YY,G_*)$ does not always act faithfully on its Bass-Serre tree
$T(\YY,G_*)$, that is, the kernel of its action might be nontrivial.

\medskip
A {\it subgraph of subgroups}\index{subgraph of subgroups}
of $(\YY,G_*)$ is a  graph of groups $(\YY',G'_*)$ where 

$\bullet$~ $\YY'$ is a subgraph of $\YY$,

$\bullet$~ for every $v\in V\YY'$, the group $G'_v$ is a subgroup of
$G_v$,

$\bullet$~ for every $e\in E\YY'$, the group $G'_e$ is a subgroup of
$G_e$, 

$\bullet$~ the monomorphism $\rho'_e:G'_e\ra G'_{t(e)}$ is the
restriction to $G'_e$ of the monomorphism $\rho_e:G_e\ra G_{t(e)}$,
and
$$
G'_{t(e)}\cap \rho_e(G_e)=\rho'_e(G'_e)\;.
$$

This condition, first introduced in \cite[Coro.~1.14]{Bass93}, is
equivalent to the injectivity of the natural map $G'_{t(e)}
/\rho'_e(G'_e)\ra G_{t(e)}/\rho_e(G_e)$ for every $e\in E\YY$. It
implies by \cite[2.15]{Bass93} when the underlying basepoint is chosen
in $\YY'$, that

$\bullet$~~the fundamental group $\Ga'=\pi_1(\YY',G'_*)$ of
$(\YY',G'_*)$ injects into the fundamental group $\Ga=\pi_1(\YY,G_*)$
of $(\YY,G_*)$,

$\bullet$~~the Bass-Serre tree $\XX'$ of $(\YY',G'_*)$ injects in an
equivariant way into the Bass-Serre tree $\XX$ of $(\YY,G_*)$ so that
the stabiliser of $\XX'$ in $\Ga$ is $\Ga'$, and

$\bullet$~~the map $(\Ga'\bs \XX') \ra (\Ga\bs\XX)$ induced by the
inclusion map $\XX'\ra\XX$ by taking quotient is injective:
\begin{equation}\label{eq:simplesubtree}
\forall\; \ga\in\Ga ,\;\forall\; z\in V\XX'\cup E\XX',\;\;\;
{\rm if}\;\; \ga z\in V\XX'\cup E\XX', \;{\rm then}
\;\; \exists\; \ga'\in\Ga',\;\;\; \ga' z= \ga z\;.
\end{equation}

\medskip
The {\it edge-indexed
  graph}\index{edge-indexed graph}\index{graph!edge-indexed}
$(\YY,i)$ of the graph of groups $(\YY,G_*)$ is the graph $\YY$
endowed with the map $i:E\YY \ra \NN-\{0\}$ defined by $i(e) =
[G_{o(e)}:G_e]$ (see for instance \cite{BasKul90, BasLub01}).

In Section \ref{subsect:equicountmetricgraphgroup}, we will consider
{\em metric graphs of groups}\index{metric!graph of groups}%
\index{graph!of groups!metric} $(\YY,G_*,\lambda)$ which are graphs of
groups endowed with an edge length function $\lambda:E\YY\to\;
]0,+\infty[$ (with $\lambda(\overline{e})=\lambda(e)$ for every
$e\in E\YY$).

\smallskip 
\bexem\label{exem:quotgraphgroup} The main examples of graphs of
groups that we will consider in this book are the following ones. Let
$\XX$ be a simplicial tree and let $\Ga$ be a subgroup of
$\Aut(\XX)$. The {\em quotient graph of groups}%
\index{graph!of groups!quotient} $\Ga \dbs\XX$ is the following graph
of groups $(\YY,G_*)$.  Its underlying graph $\YY$ is the quotient
graph $\Ga \bs\XX$. Fix a lift $\wt z\in V\XX\cup E\XX$ for every
$z\in V\YY\cup E\YY$. For every $e\in E\YY$, assume that
$\overline{\wt e}= \wt{\overline{e}}$, and fix an element $g_e\in\Ga$
such that $g_e \,\wt{t(e)}= t(\wt e)$. For every $y\in V\YY\cup E\YY$,
take as $G_{y}$ the stabiliser $\Ga_{\wt y}$ in $\Ga$ of the fixed
lift $\wt y$. Take as monomorphism $\rho_e:G_e\ra G_{t(e)}$ the
restriction to $\Ga_{\wt e}$ of the conjugation $\ga\mapsto
g_e^{-1}\ga g_e$ by $g_e^{-1}$.  Note that $\Ga \dbs\XX$ has finite
vertex groups if $\XX$ is locally finite and $\Ga$ is discrete. For
every choice of basepoint in $V\YY$, there exist a group isomorphism
$\theta: \pi_1(\YY,G_*)\ra \Ga$ and a $\theta$-equivariant
simplicial isomorphism from the Bass-Serre tree $T(\YY,G_*)$ to 
$\XX$ (see for instance \cite{Serre83, BasLub01}).
\eexem

The {\em volume form}\index{volume!form!of a graph of groups} of a
graph of finite groups $(\YY,G_*)$ is the measure
$\gls{volformgraphgroup}$ on the discrete set $V\YY$, such that for
every $y\in V\YY$,
$$
\vol_{(\YY,\,G_*)}(\{y\})=\frac{1}{|G_y|}\;,
$$
where $|G_y|$ is the order of the finite group $G_y$. Its total mass,
called the {\em volume}\index{volume!of a graph of groups} of
$(\YY,G_*)$, is
$$
\gls{volumegraphgroup}=\|\vol_{(\YY,\,G_*)}\|=
\sum_{y\in V \YY}\frac{1}{|G_y|}\;.
$$ 
We denote by $\gls{ldeuxgraphgroup}=\Leb^2(V\YY,\vol_{(\YY,G_*)})$ the
complex Hilbert space of square integrable maps $V\YY\ra\CC$ for this
measure $\vol_{(\YY,G_*)}$, and by $f\mapsto {\|f\|}_2$ and $(f,g)
\mapsto {\langle f,g\rangle}_2$ its norm and (antilinear on the right)
scalar product.  Let
$$
\Leb^2_0(\YY,G_*)=\big\{f\in \Leb^2(\YY,G_*)\;:\;\int
f\;d\vol_{(\YY,\,G_*)} =0\big\}\,.
$$ 
When $\Vol(\YY,G_*)$ is finite, $\Leb^2_0(\YY,G_*)$ is 
the orthogonal subspace to the constant functions.

We also consider a {\em (edge-)volume form}%
\index{volume!form!of a graph of groups (edge-)} $\gls{volformgraphgroupedge}$
on the discrete set $E\YY$ such that for every $e\in E\YY$,
$$
\Tvol_{(\YY,\,G_*)}(\{e\})=\frac{1}{|G_e|}\;,
$$
with total mass 
$$
\gls{voltotgraphgroupedge}=\|\Tvol_{(\YY,\,G_*)}\|=
\sum_{e\in E\YY}\frac{1}{|G_e|}\;.
$$
The {\em (edge-)volume form}\index{volume!form!of a metric graph of
  groups} of a metric graph of groups $(\YY,G_*,\lambda)$ is given by
$$
\gls{voltotgraphgroupedgemet}=\frac{ds}{|G_e|}
$$ on each edge $e$ of $\YY$ parameterised by its arclength $s$, so
that its total mass is
$$
\gls{voltotmetgraphgroupedge}=\|\Tvol_{(\YY,\,G_*,\lambda)}\|=
\sum_{e\in E\YY}\frac{\lambda(e)}{|G_e|}\;.
$$  
For $\lambda\equiv 1$, this total mass agrees with that of the
discrete definition above.

\brema \label{rem:TVoltorsionfree}
Note that $\TVol(\YY,\,G_*)=\card(E\YY)$ when the edge groups
are trivial. We have
$$
\TVol(\YY,\,G_*) =\sum_{e\in E \YY}\frac{1}{|G_e|}
=\sum_{y\in V \YY}\frac{1}{|G_y|}\sum_{e\in E \YY,\;o(e)=y}\frac{|G_y|}{|G_e|}
=\sum_{y\in V \YY}\frac{\deg (\wt y)}{|G_y|}\;,
$$
where $\wt y$ is any lift of $y$ in the Bass-Serre tree of
$(\YY,G_*)$. In particular, if $\XX$ is a uniform simplicial tree and
$\Ga$ is discrete subgroup of $\Aut \XX$, then the finiteness of
$\Vol(\Ga\dbs\XX)$ and of $\TVol(\Ga\dbs\XX)$ are equivalent.
Defining the volume form on $V\YY$ by $\{y\}\mapsto \frac{\deg (\wt
  y)}{|G_y|}$ sometimes makes formulas simpler, but we will follow the
convention which occurs in the classical references (see for instance
\cite{BasLub01}). 

If the Bass-Serre tree of $(\YY,\,G_*)$ is $(q+1)$-regular, then
\begin{equation}\label{eq:volTvol}
\pi_*\Tvol_{\YY,\,G_*}= (q+1) \vol_{\YY,\,G_*}\;\;\;{\rm and}\;\;\;
\TVol(\YY,\,G_*)= (q+1) \Vol(\YY,\,G_*)\,.
\end{equation}
\erema

\medskip
We say that a discrete group of (inversion-free) automorphisms $\Ga$
of a locally finite metric or simplicial tree $(\XX,\lambda)$ is a
{\em (tree) lattice}\index{lattice!(tree)} of $(\XX,\lambda)$ if the
quotient graph of groups $\Ga \dbs\XX$ has finite volume:
$$
\Vol(\Ga\dbs\XX)<+\infty\;.
$$
This implies by \cite[Prop.~4.5]{BasKul90}\footnote{using the fact that
$\Aut(\XX,\lambda)$ is a closed subgroup of $\Aut(\XX)$} that $\Ga$ is
a lattice\footnote{Recall that a {\em lattice}\index{lattice} in a
  locally compact group $G$ is a discrete subgroup $\Ga'$ of $G$ such
  that the left quotient space $\Ga'\bs G$ admits a probability
  measure invariant under translations on the right by $G$.} in the
locally compact group $\Aut(\XX,\lambda)$ (hence that
$\Aut(\XX,\lambda)$ is unimodular, see for instance \cite[Chap.~I,
  Rem.~1.9]{Raghunathan72}), the converse being true for instance if
$(\XX,\lambda)$ is uniform.  If $\Ga$ is a {\em uniform
  lattice}\index{lattice!uniform} of $(\XX,\lambda)$, that is, if
$\Ga$ is a discrete subgroup of $\Aut(\XX,\lambda)$ and if the
quotient graph $\Ga \bs\XX$ is finite, then $\Ga$ is clearly a (tree)
lattice of $(\XX,\lambda)$.

\medskip 
A graph of finite groups $(\YY,G_*)$ is a {\em cuspidal
  ray}\index{cuspidal ray}\index{ray!cuspidal} if $\YY$ is a
simplicial ray such that the homomorphisms $G_{e_n}\to G_{o(e_n)}$ are
surjective for its sequence of consecutive edges $(e_i)_{i\in\NN}$
oriented towards the unique end of $\YY$.  By \cite{Paulin02}, a
discrete group $\Ga'$ of $\Aut(\XX)$ (hence of $\Isom(|\XX|_1)$) is
geometrically finite if and only if it is nonelementary and if the
quotient graph of groups by $\Ga'$ of its minimal nonempty invariant
subtree is the union of a finite graph of groups and a finite number
of cuspidal rays attached to the finite graph at their finite
endpoints.

\brema \label{rem:geomfiniimpllattice} If $\XX$ is a locally finite
simplicial tree and if $\Ga'$ is a geometrically finite discrete group
of $\Aut(\XX)$ such that the convex hull of its limit set
$\C\Lambda\Ga'$ is (the geometric realisation of) a uniform tree, then
$\Ga'$ is a lattice of the simplicial tree $\C\Lambda\Ga'$.  
\erema

\dem
Since $\C\Lambda\Ga'$ is uniform, there is a uniform upper bound
on the length of an edge path in $\C\Lambda\Ga'$ which injects in
$\Ga'\bs\C\Lambda\Ga'$ and such that the stabiliser of each edge of this
edge path is equal to the stabilisers of both endpoints of this
edge. It is hence easy to see that the volume of each of the (finitely
many) cuspidal rays in $\Ga'\dbs\C\Lambda\Ga'$ is finite, by a
geometric series argument. Hence the volume of $\Ga'\dbs\C\Lambda\Ga'$
is finite.  
\cqfd

\bigskip
Note that contrarily to the case of Riemannian manifolds, there
are many more (tree) lattices than there are geometrically finite
(tree) lattices, even in regular trees, see for instance
\cite{BasLub01}.

\bigskip
In Part \ref{sect:arithappli} of this book, we will consider
simplicial graphs of groups that arise from the arithmetic of
non-Archimedean local fields.  We say that a discrete group $\Ga$ of
(inversion-free) automorphisms of a simplicial tree $\XX$ is {\em
  algebraic}\index{algebraic lattice}\index{lattice!algebraic} if
there exist a non-Archimedean local field $\wh K$ (a finite extension
of $\QQ_p$ for some prime $p$ or the field of formal Laurent series
over a finite field) and a connected semi-simple algebraic group
$\underline{G}$ with finite centre defined over $\wh K$, of $\wh
K$-rank one, such that $\XX$ identifies with the Bruhat-Tits tree of
$\underline{G}$ in such a way that $\Ga$ identifies with a lattice of
the locally compact group $\underline{G}(\wh K)$.  If $\Ga$ is
algebraic, then $\Ga$ is geometrically finite by
\cite{Lubotzky91}. Note that $\XX$ is then bipartite, see Section 2 of
op.~cit.~for a discussion and references. See Sections
\ref{sec:fieldsandvaluations} and \ref{subsec:BruhatTitstrees} for
more details, and their subsequent Sections for arithmetic
applications arising from algebraic lattices.

\chapter{Potentials, critical exponents and Gibbs 
cocycles}
\label{sect:Potentialmeasure}

Let $X$ be a geodesically complete proper $\CAT(-1)$ space, let
$x_0\in X$ be an arbitrary basepoint, and let $\Ga$ be a
nonelementary discrete group of isometries of $X$. 

In this Chapter, we define potentials on
$T^1X$, which are new data on $X$ in addition to its geometry. We
introduce the fundamental tools associated with potentials, and we
give some of their basic properties. The development follows
\cite{PauPolSha15} with modifications to fit the present more general
context.

In Section \ref{subsec:cond}, given a simplicial or metric tree
$(\XX,\lambda)$, with geometric realisation $X$, we introduce a
natural method to associate a ($\Ga$-invariant) potential $\wt F_c:
T^1X\ra \RR$ to a $\Ga$-invariant function $\wt c:E\XX \ra \RR$
defined on the set of edges of $\XX$, that we call a system of
conductances on $\XX$.  This construction gives a nonsymmetric
generalisation of electric networks.

\section{Background on (uniformly local) H\"older-continuity}
\label{subsec:holdercont}

In this preliminary Section, we recall the notion of H\"older-continuity
we will use in this book, which needs to be defined appropriately in
order to deal with noncompactness issues. The H\"older-continuity will
be used on the one hand for potentials when $X$ is a Riemannian manifold
in Section \ref{subsec:potentials}, and on the other hand for error
term estimates in Chapters \ref{sec:mixingrate},
\ref{sec:skinningwithpot} and \ref{sec:equidarcs}.

\medskip
As in \cite{PauPolSha15}, we will use the following uniformly local
definition of H\"older-continuous maps. Let $E$ and $E'$ be two metric
spaces, and let $\alpha\in \mathopen{]}0, 1\mathclose{]}$. A map
$f:E\ra E'$ is
\begin{itemize}
\item {\em $\alpha$-H\"older-continuous}\index{Holder@H\"older-continuity}
if there exist $c,\epsilon>0$ such that for all $x,y\in E$ with
$d(x,y)\leq \epsilon$, we have
$$
d(f(x),f(y))\leq c\;d(x,y)^\alpha\;.
$$
\item
 {\em locally
  $\alpha$-H\"older-continuous}\index{locally!H\"older-continuous}%
\index{Holder@H\"older-continuity!local} if for every $x\in E$, there
exists a neighbourhood $U$ of $x$ such that the restriction of $f$ to
$U$ is $\alpha$-H\"older-continuous;
\item
 {\em H\"older-continuous}%
 \index{Holder@H\"older-continuity!local}%
 \index{locally!H\"older-continuous} (respectively {\em locally
   H\"older-continuous})\index{Holder@H\"older-continuity!local}%
 \index{locally!H\"older-continuous} if
 there exists $\alpha\in\mathopen{]}0,1\mathclose{]}$ such that $f$ is
$\alpha$-H\"older-continuous (respectively locally
$\alpha$-H\"older-continuous);
\item {\em Lipschitz}\index{Lipschitz} if it is
$1$-H\"older-continuous and {\em locally
  Lipschitz}\index{locally!Lipschitz} if it is locally
$1$-H\"older-continuous.
\end{itemize}

\medskip
Let $E$ be a set. Two distances $d$ and $d'$ on $E$ are (uniformly
locally) {\em H\"older-equivalent}\index{Holder@H\"older-equivalent} if
the identity map from $(E,d)$ to $(E,d')$ and the identity map from
$(E,d')$ to $(E,d)$ are H\"older-continuous. This is an equivalence
relation on the set of distances on $E$, and a {\em H\"older
  structure}\index{Holder@H\"older-equivalent} on $E$ is the choice of
such an equivalence class. For a map between two metric spaces, to be
H\"older-continuous depends only on the Hölder structures on the
source and target spaces.

\medskip
Let $E$ and $E'$ be two metric spaces. We say that a map $f:E\ra E'$ has
\begin{itemize} 
\item {\em at most linear growth}\index{growth!linear}\index{linear
    growth} if there exist $a,b\geq0$ such that $d(f(x),f(y))\leq
  a\,d(x,y)+b$ for all $x,y\in E$,
\item {\em subexponential
  growth}\index{growth!subexponential}\index{subexponential growth} if
for every $a>0$, there exists $b\geq 0$ such that $d(f(x),f(y))\leq
b\,e^{a\,d(x,\,y)}$ for all $x,y\in E$.
\end{itemize}

\brema\label{rem:soriteholder} 
When $E$ is a geodesic space, a consequence of the (uniformly local)
H\"older-continuous property of $f:E\ra E'$ is that $f$ then has at
most linear growth: the definition implies that
$$
d(f(x),f(y))\leq c\,\epsilon^{\alpha-1}\,d(x,y) + c\,\epsilon^{\alpha}
$$ 
for all $x,y$ in $X$, by subdividing the geodesic segment in $E$ from
$x$ to $y$ into $\big\lceil \frac{d(x,y)}{\epsilon}\big\rceil$ segments
of equal lengths at most $\epsilon$ and using the triangle inequality
in $E'$.  
\erema

\medskip
The following result (due to Bartels-Lück \cite{BarLuc12} with a
different distance on $\gengeod X$) proves in particular that the
footpoint projection $\pi:\gengeod X\ra X$ is
$\frac{1}{2}$-H\"older-continuous, as claimed in Section
\ref{subsec:lines}. Recall that $\gengeod X$ is endowed with the
distance $d$ defined by Equation \eqref{eq:geodesicmetric}.

\bprop \label{prop:footpointholdercont}
For every $t\in\RR$, the map from $\gengeod X$ to $X$ defined by
$\ell\mapsto \ell(t)$ is  $\frac{1}{2}$-H\"older-continuous.  
\eprop

\dem
Let $\ell,\ell'\in\gengeod X$ be such that $d(\ell,\ell')\leq 1$.
Assume that $t\geq 0$, otherwise the argument is similar, replacing
$\int_{t}^{t+\epsilon}$ by $\int_{t-\epsilon}^{t}$.  For every
$\epsilon>0$, we have by the triangle inequality
$$
d(\ell,\ell')\geq 
\int_{t}^{t+\epsilon} d(\ell(s),\ell'(s))\;e^{-2 s}\,ds \geq
\big(d(\ell(t),\ell'(t)) -2\epsilon\big)\,\epsilon\; e^{-2 t -2\epsilon}\;.
$$
If $d(\ell(t),\ell'(t))\geq 4$, let $\epsilon=1$, so that
$d(\ell,\ell')\geq \frac{e^{-2t-2}}{2}\;d(\ell(t),\ell'(t))$, hence 
$$
d(\ell(t),\ell'(t))\leq 2\;e^{2t+2}\;d(\ell,\ell')^{\frac{1}{2}}\;.
$$ 
If $d(\ell(t),\ell'(t))\leq 4$, let $\epsilon= \frac{1}{4}
\;d(\ell(t),\ell'(t))\leq 1$. Hence $d(\ell,\ell')\geq \frac{1}{8}
\;d(\ell(t),\ell'(t))^2 \,e^{-2t-2}$, so that
$$
d(\ell(t),\ell'(t))\leq 
2\,\sqrt{2}\;e^{t+1}\;d(\ell,\ell')^{\frac{1}{2}}\;. \;\;\;\Box
$$

\medskip

When $X$ is an $\RR$-tree, the regularity property of the footpoint
projection $\pi:\G X\ra X$ is stronger than the one given by
Proposition \ref{prop:footpointholdercont} (see Lemma
\ref{lem:troisboulesmet} (2)). The results below, that will be needed
in Sections \ref{subsec:rateequidtrees} and
\ref{subsect:errormetricgraphgroup}, say that not only the evaluation
maps $\ell\mapsto \ell(t)$ are $\frac{1}{2}$-H\"older-continuous, but
so are the endpoint maps $\ell\mapsto \ell_\pm$. As we will only need
these facts in the tree case, we prove them only when $X$ is an
$\RR$-tree, and we start by a simplicial version of it.

\blemm \label{lem:troisboulessimp}
Let $\XX$ be a simplicial tree. There are universal constants
$\epsilon_0>0, c_0\geq 1$ such that for all $\epsilon\in
\mathopen{]}0, \epsilon_0\mathclose{[}\,$ and $\ell\in\G\XX$, the ball
$B_d(\ell,\epsilon)$ is contained in
$$
\{\ell'\in\G\XX\;:\; \ell'(0)=\ell(0), \;\ell'_\pm\in
B_{d_{\ell(0)}}(\ell_\pm, c_0\; \sqrt{\epsilon}\,) \}
$$
and contains $\{\ell'\in\G\XX \;:\; \ell'(0) =\ell(0), \;\ell'_\pm
\in B_{d_{\ell(0)}}(\ell_\pm, \frac{1}{c_0}\;\sqrt{\epsilon}\,)\}$.

In particular, the endpoint maps $\ell\mapsto \ell_\pm$ from $\G \XX$
to $\partial_\infty \XX$ are $\frac{1}{2}$-H\"older-continuous.
\elemm

\dem 
If $\ell,\ell'\in \G\XX$ have distinct footpoints, then
$d(\ell(0),\ell'(0))\geq 1$, so that $d(\ell(t),\ell'(t))\geq
\frac{1}{2}$ if $|t|\leq \frac{1}{4}$, so that $d(\ell,\ell')\geq
\int_{-\frac{1}{4}}^{\frac{1}{4}} \frac{1}{2}\; e^{-2|t|}=\epsilon_0>0$.

Conversely, assume that $\ell,\ell'\in \G\XX$ have equal footpoints,
so that they coincide on $[-N,N']$ for some $N,N'\in\NN$. By
the definition of the visual distances (see Equation
\eqref{eq:distviscastree}), we have
$$
d_{\ell(0)}(\ell_+ ,\ell'_+)=e^{-N'}
$$ 
and similarly $d_{\ell(0)}(\ell_-,\ell'_-)=e^{-N}$.  By the definition
of the distance on $\gengeod X$ (see Equation
\eqref{eq:geodesicmetric}), we have, by an easy change of variables,
\begin{align*}
d(\ell,\ell')&=\int_{N'}^{+\infty}2\,|t-N'|\;e^{-2 t}\;dt +
\int_{-\infty}^{-N} 2\,|-N-t|\;e^{2 t}\;dt \\ &
=(e^{-2N'} +e^{-2N})\;\int_{0}^{+\infty} 2\,u\;e^{-2u} \;du=
\frac 12(e^{-2N'} +e^{-2N})\;.
\end{align*}
The result follows.
\cqfd

\medskip
Let us now give a (more technical) version of this lemma for
$\RR$-trees, also proving that the footpoint projection is
Lipschitz. If $a$ and $b$ are positive functions of some parameters,
we write $a\asymp b$ if there exists a universal constant $C>0$ such
that $\frac{1}{C}\;b\leq a \leq C\; b$.

\blemm \label{lem:troisboulesmet} Assume that $X$ is an $\RR$-tree.
\begin{enumerate}
\item
  There exists a universal constant $c_1>0$ such that for all
  $\ell,\ell'\in \G X$, if $d(\ell,\ell')\leq c_1$, then $\ell'(0) \in
  \ell(\RR)$, the intersection $\ell(\RR)\cap \ell'(\RR)$ is not
  reduced to a point, the orientations of $\ell$ and $\ell'$ coincide
  on this intersection, and
$$
d(\ell,\ell')\asymp d_{\ell(0)}(\ell_- ,\ell'_- )^2+
d_{\ell(0)}(\ell_+ ,\ell'_+ )^2+d(\ell(0),\ell'(0))\;.
$$
\item The footpoint map $\pi:\G X\ra X$ defined by $\ell\mapsto
  \ell(0)$ is (uniformly locally) Lipschitz.
\item 
There are universal constants $\epsilon_0>0, c_0\geq 1$ such that for
all $\epsilon\in\;]0,\epsilon_0[\,$ and $\ell\in\G X$, the ball
$B_d(\ell,\epsilon)$ in $\G X$ is contained in
$$
\{\ell'\in\G X\;:\;  \ell'(0)\in\ell(\RR), 
\;d(\ell'(0),\ell(0))\leq c_0\;\epsilon, \;\ell'_\pm \in 
B_{d_{\ell(0)}}(\ell_\pm, c_0\; \sqrt{\epsilon}\,) \}
$$ 
and contains 
$$
\{\ell'\in\G X \;:\; \ell'(0)\in\ell(\RR), 
\;d(\ell'(0),\ell(0))\leq \frac{1}{c_0}\;\epsilon, \;\ell'_\pm
 \in B_{d_{\ell(0)}}(\ell_\pm , \frac{1}{c_0}
\;\sqrt{\epsilon}\,)\}\;.
$$
\item The endpoint maps $\ell\mapsto \ell_\pm$ from $\G X$
to $\partial_\infty X$ are $\frac{1}{2}$-H\"older-continuous.
\end{enumerate} 
\elemm

\dem Note that Assertion (2) follows from Assertion (1) and that
Assertion (4) follows from Assertion (3).

\medskip
(1) Let $\ell,\ell'\in \G X$. If $\ell'(0) \notin \ell(\RR)$, then
$\ell'(t)\notin \ell(\RR)$ for all $t\geq 0$ or $\ell'(t)\notin
\ell(\RR)$ for all $t\leq 0$, since $X$ is an $\RR$-tree. In the first
case, we hence have $d(\ell(t),\ell'(t))\geq t$ for all $t\geq 0$,
thus $d(\ell,\ell')$ is at least $c_2=\int_0^{+\infty} te^{-2\,t}\;dt
= \frac14>0$.  The same estimate holds in the second case. By
symmetry, if $\ell(0) \notin \ell'(\RR)$, then $d(\ell,\ell') \geq
c_2>0$. This argument furthermore shows that if the geodesic segment
(or ray or line) $\ell(\RR)\cap \ell'(\RR)$ is reduced to a point,
then $d(\ell,\ell')$ is at least $c_2>0$.

If $d(\ell'(0),\ell(0))\geq 1$, then $d(\ell(t),\ell'(t))\geq
\frac{1}{2}$ for $|t|\leq \frac{1}{4}$, thus
$$
d(\ell,\ell')\geq \int_{-\frac14}^{\frac14}\;\frac{1}{2}\;e^{-2|t|}\;dt\;,
$$
which is a positive universal constant.

If $d(\ell'(0),\ell(0))\leq 1$, if $\ell(\RR)\cap \ell'(\RR)$ contains
$\ell'(0)$ and is not reduced to a point, but if the orientations of
$\ell$ and $\ell'$ do not coincide on this intersection, then
$$
d(\ell(t),\ell'(t))\geq 2t-d(\ell(0),\ell'(0))\geq 2t-1
$$
for all $t\geq 1$, so that $d(\ell,\ell')$ is at least
$\int_{1}^{+\infty} (2t-1)\,e^{-2t}\;dt$, which is a positive
constant.

\medskip
Assume now that $\ell'(0) \in \ell(\RR)$ and $\ell(0) \in \ell'(\RR)$,
that $d(\ell'(0),\ell(0)) \leq 1$, that $\ell(\RR)\cap \ell'(\RR)$ is
not reduced to a point and that the orientations of $\ell$ and $\ell'$
coincide on this intersection.  Then there exists $s\in\RR$ such that
$\ell'(0)= \ell(s)$, so that $|s|=d(\ell(0), \ell'(0))\leq 1$. Assume
for instance that $s\geq 0$, the other case being treated similarly.
Then there exist $S,S'\geq 0$ maximal such that $\ell'(t)=\ell(t+s)$
for all $t\in[-S,S']$. We use the conventions that $S=+\infty$ if
$\ell'_-= \ell_-$, that $S'=+\infty$ if $\ell'_+= \ell_+$, and that
$e^{-\infty} =0$. Since $\ell(0) \in \ell'(\RR)$, we have $-S+s\leq
0$.

\begin{center}
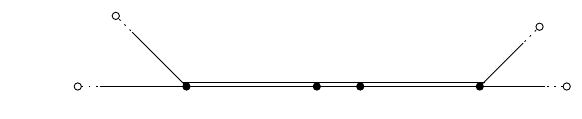
\end{center}

By the definition of the visual distances (see Equation
\eqref{eq:distviscastree}), we have, for $t$ big enough,
$$
d_{\ell(0)}(\ell_+ ,\ell'_+ )=e^{-d(\ell(0),\,\ell(S'+s))}=e^{-S'-s}
\asymp e^{-S'}\;.
$$ 
Similarly $d_{\ell(0)}(\ell_-,\ell'_-)\asymp e^{-S}$.

As can be seen in the above picture, we have
$$
d(\ell(t),\ell'(t))=\left\{\begin{array}{ll}
-2t-2S+s & {\rm if}\;\;t\leq -S\\
s & {\rm if}\;\;-S\leq t\leq S'+s\\
2t-2S'-s & {\rm if}\;\;t\geq S'+s\;.
\end{array}\right.
$$
By the definition of the distance on $\gengeod X$ (see Equation
\eqref{eq:geodesicmetric}), by easy changes of variables, assuming
that at least one of $S,S'$ is at least $1$ for the last line (otherwise
the previous line shows that $d(\ell,\ell')$ is larger than a positive
constant), we have
\begin{align*}
d(\ell,\ell')&=\int_{-\infty}^{-S}\;(-2t-2S+s)\;e^{-2 |t|}\;dt+
\int_{-S}^{S'+s}\;s\;e^{-2 |t|}\;dt+
\int_{S'+s}^{+\infty}(2t-2S'-s)\;e^{-2 t}\;dt \\ &
=\frac{1}{2}\,e^{-2S+s}\int_{s}^{+\infty}\;u\;e^{-u}\;du
+s\int_{-S}^{S'+s}\;e^{-2 |t|}\;dt
+ \frac{1}{2}\,e^{-2S'-s}\int_{s}^{+\infty}\;u\;e^{-u}\;du \\ &
\asymp e^{-2S}+e^{-2S'}+s\;.
\end{align*}
Assertion (1) of Lemma \ref{lem:troisboulesmet}  follows.

\medskip (3) The first inclusion in Assertion (3) follows easily from
Assertion (1). The second inclusion follows from the argument of its
proof, and the fact that if $d_{\ell(0)}(\ell_+ ,\ell'_+ )$,
$d_{\ell(0)}(\ell_-,\ell_-)$ and $d(\ell(0),\ell'(0))$ are at most
some small positive constant, then $\ell(\RR)\cap \ell'(\RR)$ contains
$\ell'(0)$ and is not reduced to a point.  \cqfd

\medskip 
When $X$ is a Riemannian manifold with pinched sectional curvature,
the following result says that the H\"older structure defined by the
distance $d$ given in Equation \eqref{eq:geodesicmetric} on the unit
tangent bundle of $X$, identified with $\G X$ and $T^1X$ as explained
previously, is the usual one.

Recall that {\it Sasaki's metric}\index{Sasaki's metric} on $T^1X$ is
the Riemannian metric induced on the submanifold $T^1X$ by the
canonical Riemannian metric on the tangent bundle $TX$, such that for
every $v\in TX$, if $T_vTX=V_v\oplus H_v$ is the direct sum
decomposition defined by the Levi-Civita connection of the Riemannian
manifold $X$, then
\begin{enumerate}
\item[$\bullet$] the direct sum $V_v\oplus H_v$ is orthogonal for
  Sasaki's metric,
\item[$\bullet$] the canonical isomorphism $V_v\simeq T_{\pi(v)}X$ is
  an isometry, when $V_v$ is endowed with Sasaki's scalar product,
\item[$\bullet$] the restriction of the tangent map of the footpoint
  projection $T\pi : H_v\ra T_{\pi(v)}X$ is an isometry, when $H_v$
  is endowed with Sasaki's scalar product.
\end{enumerate}

\bprop \label{prop:maniouholder}
When $X$ is a Riemannian manifold with pinched sectional
curvature, the following distances on $T^1X$ are H\"older-equivalent:
\begin{enumerate}
\item
  the Riemannian distance on $T^1X$ defined by Sasaki's metric,
\item
  the distance $\delta_1$ on $T^1X$ defined, for all
  $\ell,\ell'\in T^1X$, by
$$
\delta_1(\ell,\ell')=
\exp({-\sup\{t\geq 0\;:\;\sup_{s\in[-t,t]}d(\ell(s),\ell'(s))
  \leq 1\}})\;,
$$ 
with the convention $\delta_1(\ell,\ell')=1$ if $d(\ell(0),\ell'(0))>
1$ and $\delta_1(\ell,\ell')=0$ if $\ell=\ell'$,
\item
  the distance $\delta_2$ on $T^1X$ defined, for all
  $\ell,\ell'\in T^1X$, by
$$
\delta_2(\ell,\ell')=\sup_{s\in[0,1]}d(\ell(s),\ell'(s))\;,
$$
\item
  the distance $d$ defined by Equation \eqref{eq:geodesicmetric}.
\end{enumerate}
\eprop

Note that we could replace the interval $[0,1]$ in the definition of
$\delta_2$ by any interval $[a,b]$ with $a<b$.

\medskip
\dem 
The fact that the first three distances on $T^1X$ are
H\"older-equivalent is already known, see for instance
\cite[p.~70]{Ballmann95}.

Let us hence prove that $d$ and $\delta_2$ are H\"older-equivalent. By
the convexity of the distance function between two geodesic segments
in a $\CAT(-1)$-space, we have
$$
\delta_2(\ell,\ell')=
\max \{d(\ell(0),\ell'(0)),\;d(\ell(1),\ell'(1))\}\;.
$$
Hence by Proposition \ref{prop:footpointholdercont} applied twice
(with $t=0$ and $t=1$), there exists $c>0$ such that if
$d(\ell,\ell')\leq 1$, then $\delta_2(\ell,\ell') \leq
c\;d(\ell,\ell')^{\frac{1}{2}}$. Therefore the identity map from
$(T^1X,d)$ to $(T^1X,\delta_2)$ is H\"older-continuous.

\medskip
Let us now prove conversely that there exist two constants $c'_1\geq
1$ and $c'_2\in\mathopen{]}0,1\mathclose{]}$ such that for all $\ell,
\ell'\in\G X$, if $\delta_2(\ell,\ell')\leq \frac{1}{2}$, then
\begin{equation}\label{eq:holderddeltade}
d(\ell,\ell') \leq c'_1\;\delta_2(\ell,\ell')^{c'_2}\;.
\end{equation}

Let $x=\ell(0)$, $y=\ell'(0)$, $z=\ell(1)$, $z'=\ell'(1)$. Note that
$d(y,z)\geq \frac{1}{2}$ by the triangle inequality, since $d(x,y)\leq
\delta_2(\ell,\ell')\leq \frac{1}{2}$ and $d(x,z)=1$. Let $\ell''$ be
the geodesic line through $y$ and $z$, oriented from $y$ to $z$ with
$\ell''(0)=y$.  Let us prove that $d(\ell,\ell'') \leq
\frac{c'_1}{2}\;d(x,y)^{c'_2}$ for appropriate constants $c'_1$ and
$c'_2$.  A similar argument proves that $d(\ell'',\ell') \leq
\frac{c'_1}{2}\;d(z,z')^{c'_2}$, and the triangle inequality for the
distance $d$ gives Equation \eqref{eq:holderddeltade}.

\begin{center}
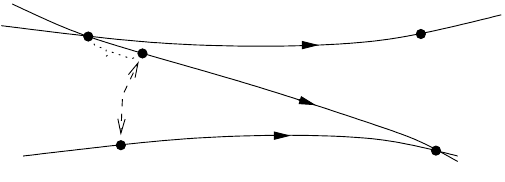
\end{center}

If $d(y,z)\leq d(x,z)$, let $x'\in [x,z]$ be such $d(x',z)=d(y,z) \geq
\frac{1}{2}$.  We have $s=d(x,x')=d(x,z)-d(x',z)= d(x,z)-d(y,z)\leq
d(x,y)$ by the triangle inequality, and $d(x',y)\leq d(x,y)$ by
convexity. By Equation \eqref{eq:disttranslatgeod}, we have
$$
d(\ell,\ell'')\leq d(\ell,\flow{s}\ell)+ d(\flow{s}\ell,\ell'')\leq 
s+ d(\flow{s}\ell,\ell'')\leq d(x,y)+ d(\flow{s}\ell,\ell'')\;.
$$
If $d(y,z)\geq d(x,z)$, let $y'\in [y,z]$ be such $d(y',z) =
d(x,z)=1\geq \frac{1}{2}$ (see the above picture). We similarly have
$d(x,y')\leq d(x,y)$ and $d(\ell,\ell'')\leq d(x,y)+
d(\flow{s}\ell'',\ell'')$ if $s=d(y,y')$.

We hence only have to prove that for all $x,y,z\in X$ with $0<d(x,y)
\leq \frac{1}{2}$ and $d(x,z)=d(y,z) \geq \frac{1}{2}$, if $u$ and $v$
are the unit tangent vectors at $x$ and $y$ respectively pointing to
$z$, then $d(u,v)\leq c''_1\;d(x,y)^{c''_2}$ for appropriate constants
$c''_1$ and $c''_2$. This follows from the following lemma with $t=0$,
since $\frac{1}{2}\geq \frac{1}{2}\;\sqrt{d(x,y)}$. We will need its
more general version in the proof of Lemma \ref{lem:technicholderPPS}.

\blemm\label{lem:holderdeltaded} There exist two constants $c_8\geq 1$
and $c_3\in\mathopen{]}0,1\mathclose{]}$ such that for all $x,y,z\in
X$ with $0<d(x,y)\leq 1$ and $\rho=d(x,z)=d(y,z) \geq
\frac{1}{2}\;\sqrt{d(x,y)}$, for every $t\in\mathopen{[}0,
  \rho\mathclose{]}$, if $u$ and $v$ are the unit tangent vectors at
$x$ and $y$ respectively pointing to $z$, then
\begin{equation}\label{eq:holderdeltaded}
d(\flow{t}u,\flow{t}v)\leq c_8\;d(x,y)^{c_3}(e^{-t}+e^{2t-2\rho})\;.
\end{equation}
\elemm

\dem
For every $s\in\RR$, let $x_s= \pi(\flow{s}u)$ and $y_s=
\pi(\flow{s}v)$, so that $x_0=x$, $y_0=y$ and $x_\rho=y_\rho=z$.  Let
$-r<r'$ in $\RR$ be such that $d(x_{-r},y_{-r})=d(x_{r'},y_{r'}) =1$,
which exist since $x\neq y$.  We have $r> 0$ by convexity since
$d(x_0,y_0)\leq 1$, and $r'>\rho$.

\begin{center}
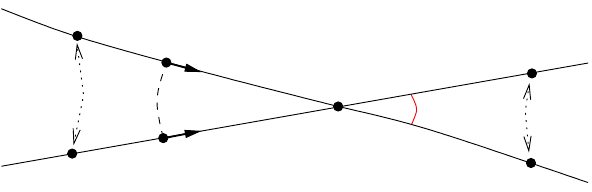
\end{center}

\medskip
\noindent{\bf Claim 1 : } There exists a constant $c_3\in
\mathopen{]}0, 1\mathclose{]}$ depending only on the lower bound of
    the curvature of $X$ such that
$$
e^{-r}\leq 2\; d(x,y)^{c_3}\;\;\;{\rm and}\;\;\;
e^{-r'}\leq d(x,y)^{c_3}\;e^{-\rho}\;.
$$

\medskip
\noindent
\begin{minipage}{9.4cm}
\dem 
Let $a\geq 1$ be such that the sectional curvature of $X$ is at least
$-a^2$. Consider the comparison triangle $(\overline{x_{-r}},
\overline{y_{-r}}, \overline{z})$, in the real hyperbolic space
$\frac{1}{a}\,\hdr$ with constant sectional curvature $-a^2$, to the
triple of points $(x_{-r},y_{-r},z)$ in $X$. Let $\overline{\theta}$
be its angle at $\overline{z}$, and let $\overline{x}, \overline{y}$
be the points corresponding to $x,y$ (so that we have $d(\overline{x}, 
\overline{z})=d(\overline{y}, \overline{z})=\rho$).
\end{minipage}
\begin{minipage}{5.5cm}
\begin{center}
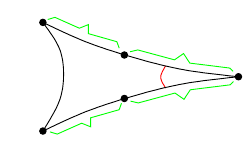
\end{center}
\end{minipage}

\medskip
Since the geodesic triangles in $X$ are less pinched than the geodesic
triangles in $\frac{1}{a} \,\hdr$, we have $d(\overline{x},
\overline{y}) \leq d(x,y)$. By the hyperbolic sine rule in
$\frac{1}{a}\,\hdr$, we have
$$
\frac{\sinh \frac{a}{2}}
     {\sinh (a(r+\rho))}=\frac{\sinh \big(\frac{a}{2}\;d(x_{-r},y_{-r})\big)}
     {\sinh \big(a\;d(x_{-r},z)\big)}=\sin\frac{\overline{\theta}}{2}=
\frac{\sinh \big(\frac{a}{2}\;d(\overline{x},\overline{y})\big)}
     {\sinh \big(a\;d(x,z)\big)}\leq
\frac{\sinh \big(\frac{a}{2}\;d(x,y)\big)}{\sinh (a\,\rho)}\;.
$$

Since the map $t\mapsto (\ln 2)t+\ln(1-\frac{t}{2})$ is nonnegative on
$[0,1]$, and since $\rho\geq\frac{1}{2}\sqrt{d(x,y)}$, we have
$$
e^{2a\rho}\geq e^{\ln 2\sqrt{d(x,y)}}\geq \frac{1}{1-\frac{1}{2}\sqrt{d(x,y)}}\;.
$$ 
Hence $\sinh (a\rho)=\frac{e^{a\rho}-e^{-a\rho}}{2}\geq
\frac{1}{4}\,e^{a\rho}\sqrt{d(x,y)}$. Therefore, since $d(x,y)\leq 1$,
$$
e^{-ar}=\frac{e^{a\rho}}{e^{ar +a\rho}}
\leq \frac{4\sinh (a\rho)\,/\,\sqrt{d(x,y)}}{2\sinh (ar +a\rho)}
\leq
\frac{2\sinh \big(\frac{a}{2}\;d(x,y)\big)}{\sqrt{d(x,y)}\,
\sinh \frac{a}{2}} \leq 2\sqrt{d(x,y)}\;.
$$
So that the first assertion of Claim 1 follows, since
$2^{\frac{1}{a}}\leq 2$, with $c_3=\frac{1}{2a}$.

\medskip
In order to prove the second assertion of Claim 1, let $\theta$ be the
angle at $z$ between the unit tangent vectors $\flow{\rho}u$ and
$\flow{\rho}v$ (see the picture before the statement of Claim 1).
Since the geodesic triangles in $X$ are less pinched than their
comparison triangles in $\frac{1}{a} \,\hdr$, and again by the
hyperbolic sine rule in $\frac{1}{a}\,\hdr$, we have
$$
\sin\frac{\theta}{2}\geq
\frac{\sinh \big(\frac{a}{2}\;d(x_{r'},y_{r'})\big)}
{\sinh \big(a\;d(x_{r'},z)\big)}=
\frac{\sinh \frac{a}{2}}{\sinh (a(r'-\rho))}\;.
$$
Since the geodesic triangles in $X$ are more pinched than their
comparison triangles in $\hdr$, since $d(x,y)\leq 1$ and since $\sinh
\rho\geq \rho\geq\frac{1}{2}\sqrt{d(x,y)}$, we have
$$
\sin\frac{\theta}{2}\leq
\frac{\sinh \big(\frac{1}{2}\;d(x,y)\big)}
{\sinh d(x,z)}\leq \frac{(\sinh \frac{1}{2})d(x,y)}{\sinh \rho}
\leq 2\, (\sinh \frac{1}{2})\,\sqrt{d(x,y)}\;.
$$
Let $c_4=\Big(\frac{\sinh \frac{1}{2}}{\sinh \frac{a}{2}}
\Big)^{\frac{1}{a}}$, which is positive and strictly less than $1$
since $a>1$. Then
\begin{align}
  e^{-ar'}&=
  \frac{e^{-a\rho}}{e^{ar' -a\rho}}\leq\frac{e^{-a\rho}}{2\sinh(a(r' -\rho))}
\leq e^{-a\rho}\;\frac{\sinh \frac{1}{2}}{2\sinh \frac{a}{2}}\,
\frac{d(x,y)}{\sinh \rho}\nonumber\\ &
\leq e^{-a\rho}\;{c_4}^{a}\;
\sqrt{d(x,y)}\label{eq:majorexpmoinsarp}\;.
\end{align}
So that the second assertion of Claim 1 follows, since $c_4\leq 1$,
again with $c_3=\frac{1}{2a}$. \cqfd

\medskip
Let us remark that since $d(x,y)\leq 1$, Equation
\eqref{eq:majorexpmoinsarp} also gives that
$$
r'-\rho\geq -\ln c_4>0\;.
$$
Since $\theta$ is also the angle between the unit tangent vectors
$-\flow{\rho}u$ and $-\flow{\rho}v$, a similar argument gives that
$$
\frac{\sinh \frac{1}{2}}{\sinh (r'-\rho)}=
\frac{\sinh \big(\frac{1}{2}\;d(x_{r'},y_{r'})\big)}
     {\sinh d(x_{r'},z)}\geq \sin\frac{\theta}{2}\geq
     \frac{\sinh \big(\frac{a}{2}\;d(x,y)\big)}
{\sinh \big(a\,d(x,z)\big)}\geq \frac{a\;d(x,y)}{2\;\sinh (a\, \rho)}\;.
$$ 
Let $c_5=\frac{2\,\sinh\frac{1}{2}}{a(1-c_4^2)}>0$. Since the map
$t\mapsto \frac{\sinh(t)}{e^t}$ is nondecreasing on
$\RR$, we hence have
$$
e^{r'-\rho}\leq \frac{e^{-\ln c_4}}{\sinh(-\ln c_4)}\sinh(r'-\rho) \leq
\frac{2}{1-c_4^2}\;\frac{2\,\sinh\frac{1}{2}}{a\;d(x,y)}\;\sinh(a\,
\rho) \leq \frac{c_5}{d(x,y)}\;e^{a\rho}
\;.
$$
Therefore
\begin{equation}\label{eq:controlsuprprime}
1+r'\leq (a+1)\rho + (1+\ln c_5) -\ln d(x,y)\;,
\end{equation}
a formula which be useful later on.

\medskip
Let $t\in\mathopen{[}0,\rho \mathclose{]}$. By the definition of the
distance $d$ on $\G X$ (see Equation \eqref{eq:geodesicmetric}),
we have
\begin{equation}\label{eq:subdivicinq}
d(\flow{t}u,\flow{t}v)=\int_{-\infty}^{+\infty}
d(\pi(\flow{s}\flow{t}u),\pi(\flow{s}\flow{t}v))\;e^{-2|s|}\;ds
=\int_{-\infty}^{+\infty} d(x_s,y_s)\;e^{-2|s-t|}\;ds\;.
\end{equation}
We subdivide the integral $\int_{-\infty}^{+\infty}$ as
$\int_{-\infty}^{-r}+\int_{-r}^{0}+\int_{0}^{\rho}+
\int_{\rho}^{r'}+\int_{r'}^{+\infty}$.

\medskip
\noindent{\bf Claim 2 : } We have
$$
I_1=\int_{-\infty}^{-r}d(x_s,y_s)\;
e^{-2|s-t|}\;ds\leq 2\;e^{-2t}\, d(x,y)^{c_3}\;.
$$

\dem
By the triangle inequality, for every $s\in\mathopen{[}r,+\infty
\mathclose{[}$, we have
$$
d(x_{-s},y_{-s})\leq d(x_{-s},x_{-r})+
d(x_{-r},y_{-r}) +d(y_{-r},y_{-s})\leq 2s+1\;.
$$
Hence
$$
I_1=\int_r^{+\infty}d(x_{-s},y_{-s})\;e^{-2(s+t)}\;ds\leq e^{-2t}
\int_r^{+\infty}(2s+1)\;e^{-2s}\;ds\leq e^{-2t} e^{-r}\;,
$$
and the result follows from the first assertion of Claim 1. 
\cqfd

\medskip
\noindent{\bf Claim 3 : } We have
$$
I_2=\int_{-r}^0 d(x_s,y_s)\;
e^{-2|s-t|}\;ds \leq 2\,(\sinh 1)\,e^{-2t}\, d(x,y)^{c_3}\;.
$$

\dem 
Recall that since $X$ is $\CAT(-1)$ and by an easy exercise in
hyperbolic geometry (see for instance \cite[Lem.~2.5
  (i)]{PauPolSha15}), for all $x',y',z'$ in $X$ such that
$d(x',z')=d(y',z')$, for every $t'\in\mathopen{[}0,
  d(x',z')\mathclose{]}$, if $x'_t$ (respectively $y'_t$) is the point
on $\mathopen{[}x', z'\mathclose{]}$ (respectively $\mathopen{[}y',
  z'\mathclose{]}$) at distance $t'$ from $x'$ (respectively $y'$),
then
\begin{equation}\label{eq:easyexoPPS}
d(x'_t,y'_t)\leq e^{-t'}\sinh d(x',y')\;.
\end{equation}
Hence for all $s\in \mathopen{[}0, r\mathclose{]}$, we have
$d(x_{-s},y_{-s})\leq e^{-(r-s)}\sinh d(x_{-r},y_{-r})= e^{-r+s}\sinh 1$.
Thus
\begin{align*}
I_2&=\int_0^r d(x_{-s},y_{-s})\;e^{-2(s+t)}\;ds\\ &\leq (\sinh 1)\;e^{-2t}
\int_0^{+\infty}e^{-r+s}\;e^{-2s}\;ds= (\sinh 1)\;e^{-2t}\, e^{-r}\;,
\end{align*}
and the result also follows from the first assertion of Claim 1. \cqfd

\medskip
\noindent{\bf Claim 4 : } There exists a universal
constant $c_6>0$ such that
$$
I_3=\int_{0}^{\rho}d(x_s,y_s)\;e^{-2|s-t|}\;ds\leq c_6\;e^{-t}\, d(x,y)\;.
$$

\dem
By Equation \eqref{eq:easyexoPPS} and since $d(x,y)\leq 1$, for
every $s\in \mathopen{[}0, \rho\mathclose{]}$, we have
$$
d(x_s,y_s)\leq e^{-s}\sinh d(x,y)\leq (\sinh 1)\,e^{-s}\, d(x,y)\;.
$$
Therefore
\begin{align*}
I_3&=\int_0^t d(x_{s},y_{s})\;e^{-2(t-s))}\;ds+
\int_t^\rho d(x_{s},y_{s})\;e^{-2(s-t))}\;ds\\ &
\leq (\sinh 1)\, d(x,y)
\Big(e^{-2t}\int_0^t e^s\,ds+e^{2t}\int_t^{+\infty} e^{-3s}\,ds\Big)
\leq\frac{4\,\sinh 1}{3}\, d(x,y)\,e^{-t}\;.\;\;\;\Box
\end{align*}

\medskip
\noindent{\bf Claim 5 : } We have
$$
I_4=\int_{\rho}^{r'}d(x_s,y_s)
\;e^{-2|s-t|}\;ds\leq (\sinh 1)\,e^{2t-2\rho} \,d(x,y)^{c_3}\;.
$$

\dem
By Equation \eqref{eq:easyexoPPS}, for every $s\in \mathopen{[}\rho,
r'\mathclose{]}$, we have
$$
d(x_s,y_s)\leq e^{-(r'-s)}\sinh d(x_{r'},y_{r'})= (\sinh 1)\,e^{-r'+s}\;.
$$
Hence
\begin{align*}
I_4&=\int_\rho^{r'} d(x_{s},y_{s})\;e^{-2(s-t))}\;ds
\leq (\sinh 1)\,e^{2t}\,e^{-r'}\int_\rho^{r'} e^{-s}\,ds
\leq(\sinh 1)\,e^{2t-\rho}\,e^{-r'}\;.
\end{align*}
The result follows from the second assertion of Claim 1.
\cqfd

\medskip
\noindent{\bf Claim 6 : } There exists a constant $c_7>0$ depending
only on the lower bound of the curvature of $X$ such that
$$
I_5=\int_{r'}^{+\infty}d(x_s,y_s) \;e^{-2|s-t|}\;ds\leq
c_7\;e^{2t-2\rho}\, d(x,y)^{c_3}\;.
$$

\dem
By the triangle inequality, for every $s\in\mathopen{[}r',+\infty
\mathclose{[}$, we have
$$
d(x_s,y_s)\leq d(x_{s},x_{r'})+
d(x_{r'},y_{r'}) +d(y_{r'},y_{s})\leq 2s+1\;.
$$
Hence, using the second assertion of Claim 1, we have
$$
I_5=\int_{r'}^{+\infty}d(x_{s},y_{s})\;e^{-2(s-t)}\;ds\leq e^{2t}
\int_{r'}^{+\infty}(2s+1)\;e^{-2s}\;ds= e^{2t} (1+r')e^{-2r'}\;.
$$ 
By Equation \eqref{eq:majorexpmoinsarp}, since $\sinh \rho\geq
\rho\geq\frac{1}{2}\;\sqrt{d(x,y)}$ and $c_3=\frac{1}{2a}$, we have
$$
e^{-2r'}\leq e^{-2\rho}\;c_4^2\;
\Big(\frac{d(x,y)}{2\sinh \rho}\Big)^{2/a}
\leq e^{-2\rho}\;d(x,y)^{2c_3}\;.
$$
Since $c\leq e^c$ for every $c\geq 0$ and $d(x,y)\leq 1$, we have
$-\ln d(x,y)\leq \frac{1}{c_3}\,d(x,y)^{-c_3}$. Hence by Equation
\eqref{eq:controlsuprprime}, we have
$$
I_5\leq e^{2t-2\rho}\Big(
(a+1) c_4^2\; d(x,y)^{2/a}\frac{\rho}{(2\sinh \rho)^{2/a}}
+ (1+\ln c_5)\,d(x,y)^{2c_3}+ \frac{1}{c_3}\,d(x,y)^{c_3}
\Big)\;.
$$ Assuming, as we may, that $a\ge 2$, the map $\rho'\mapsto
\frac{\rho'}{(\sinh \rho')^{2/a}}$ is bounded on $\mathopen{[}0,
+\infty \mathclose{[}$. Since $c_3=\frac{1}{2a}\leq \frac{2}{a}$,
this proves Claim 6.
\cqfd

\medskip
Since $d(x,y)\leq 1$, it follows from Equation \eqref{eq:subdivicinq}
and from Claims 2 to 6 that there exists a constant $c_8>0$ depending
only on the lower bound of the curvature of $X$ such that
$$
d(\flow{t}u,\flow{t}v)=I_1+I_2+I_3+I_4+I_5\leq
c_8\;d(x,y)^{c_3}(e^{-t}+e^{2t-2\rho})\;.
$$ 
This proves Lemma \ref{lem:holderdeltaded}, hence concludes the
proof of Proposition \ref{prop:maniouholder}.  
\cqfd~$\Box$

\medskip
Let $D$ be a nonempty proper closed convex subset of $X$. The
regularity property in the Riemannian manifold case of the fibrations
$f^\pm_D:\U^\pm_D \ra \normalpm D$ defined in Section
\ref{subsect:nbhd}, that will be needed in Section
\ref{subsect:erroterms}, is the H\"older-continuity one, as proved in
the following lemma (see also \cite[Lem.~6]{ParPau16ETDS}).

\blemm \label{lem:fpmholder} 
Assume that $X$ is a Riemannian manifold with pinched negative
curvature. The maps $f^\pm_D$ are Hölder-continuous on any set of
elements $\ell\in \U^\pm_D$ such that $d(\pi(\ell),
\pi(f^\pm_D(\ell)))$ is bounded.  
\elemm

\dem We prove the result for $f^+_D$, the one for $f^-_D$ follows
similarly.  We will use the H\"older-equivalent distances $\delta_1$
and $\delta_2$, defined in the statement of Proposition
\ref{prop:maniouholder}, on the unit tangent bundle of $X$, identified
with $\G X$ as explained previously.

\begin{center}
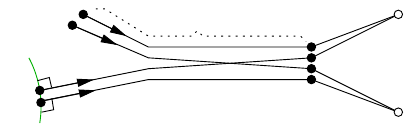
\end{center}

Let $v,v'\in T^1X$ be such that $d(v(0),v'(0))\leq 1$, let
$w=f^+_D(v)$ and $w'=f^+_D(v')$. Let $T=\sup\{t\geq 0\;:\;
\sup_{s\in[0,t]} d(v(s),v'(s)) \leq 1\}$, so that $\delta_1(v,v')\geq
e^{-T}$. We may assume that $T$ is finite, otherwise $v_+=v'_+$, hence
$w=w'$. Let $x=v(T)$ and $x'=v'(T)$, which satisfy $d(x,x')\leq
1$. Let $y$ (respectively $y'$) be the closest point to $x$
(respectively $x'$) on the geodesic ray defined by $w$ (respectively
$w'$). By convexity, since $d(v(0),w(0))$ and $d(v'(0),w'(0))$ are
bounded by a constant $c>0$ and since $v_+=w_+,v'_+=w'_+$, we have
$d(x,y)\leq c$ and $d(x',y')\leq c$. By the triangle inequality, we
have $d(y,y')\leq 2c+1$, $d(y,w(1))\geq T-2c-1$ and $d(y',w'(1))\geq
T-2c-1$. By convexity, and since closest point maps exponentially
decrease the distances, there exists a constant $c'>0$ such that
$$
\delta_2(w,w')=d(w(1),w'(1))\leq c'\,d(y,y')\; e^{-(T-2c-1)}
\leq c'\,(2c+1)\;e^{2c+1}\;\delta_1(v,v')\,.
$$
The result follows.
\cqfd

\medskip
When $X$ is an $\RR$-tree, we have a stronger version of Lemma
\ref{lem:fpmholder}, that will be needed in Section
\ref{subsect:errormetricgraphgroup}.

\blemm\label{lem:fibrationonelip} Assume that $X$ is an $\RR$-tree.
Let $\eta,R>0$ be such that $\eta\leq 1\leq \ln R$, and let $D$ be a
nonempty closed convex subset of $X$. Then the restriction to the
dynamical neighborhood $\V^\pm_{\eta,\,R} (\normalpm D)$ of the
fibration $f^\pm_D$ is (uniformly locally) Lipschitz, with constants
independent of $\eta$.
\elemm

\dem We assume for instance that $\pm=+$. Let $\ell,\ell'\in
\V^+_{\eta,\,R}(\normalout D)$ and let $w=f^+_D(\ell),
w'=f^+_D(\ell')$. 

Since the fiber over $\rho\in\normalout D$ of the restriction to
$\V^+_{\eta,\,R}(\normalout D)$ of $f^+_D$ is $V^+_{\rho,\,\eta,\,R}$, 
\footnote{See the end of Section \ref{subsect:nbhd}.} there exist
$s,s'\in\;]-\eta,\eta[$ such that $\flow{s}\ell\in B^+(w, R)$ and
$\flow{s'}\ell'\in B^+(w',R)$, so that $\flow{s}\ell(t) =w(t)$
and $\flow{s'}\ell'(t)=w'(t)$ for all $t\geq \ln R$ by the
definition of the Hamenstädt balls. Up to permuting $\ell$ and
$\ell'$, we assume that $s'\geq s$.

\begin{center}
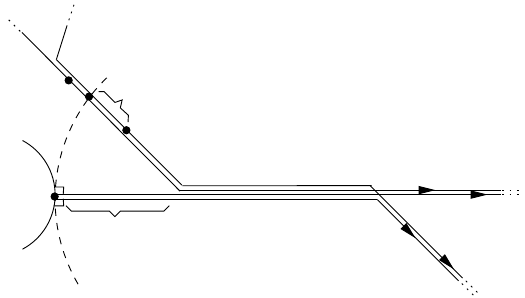
\end{center}

By (the proof of) Lemma \ref{lem:troisboulesmet} (1), there exists a
constant $c_R>0$ depending only on $R$ such that if $d(\ell,\ell')\leq
c_R$ and $s''=d(\ell(0),\ell'(0))$, then $s''=s'-s$ and the geodesic
lines $\flow{s}\ell$ and $\flow{s'}\ell'$ coincide at least on $[-\ln
R-1,\ln R+1]\,$. In particular, we have, since $|s|,|s'|\leq \eta\leq 1$,
$$
w(\ln R)=\ell(s+\ln R)=\ell'(s'+\ln R)=w'(\ln R)\;.
$$ 
Since the origin of $w$ is the closest point on $D$ to any point of
$w([0,+\infty[)$, we hence have that $w(t)=w'(t)$ for all $t\in[0,\ln
R]$. Therefore (using Equation \eqref{eq:disttranslatgeod} for the
last inequality),
\begin{align*}
d(w,w')&=\int_{\ln R}^{+\infty} d(w(t),w'(t))\,e^{-2\,t}\;dt
=\int_{\ln R}^{+\infty} 
d(\flow{s}\ell(t),\flow{s'}\ell'(t))\,e^{-2\,t}\;dt \\ & =
e^{2s}\int_{\ln R+s}^{+\infty}
d(\ell(u),\flow{s''}\ell'(u))\,e^{-2\,u}\;du \leq
e^{2s}\;d(\ell,\flow{s''}\ell')\\ & \leq
e^{2s}\;\big(d(\ell,\ell')+d(\ell',\flow{s''}\ell')\big)\leq
e^{2s}\;(d(\ell,\ell')+s'')\\ & =
e^{2s}\;\big(d(\ell,\ell')+d(\ell(0),\ell'(0))\big)\;,
\end{align*}
so that the result follows from Lemma \ref{lem:troisboulesmet} (2).

\medskip
Note that when $X$ is (the geometric realisation of) a simplicial
tree, we have $s=s'=s''=0$ and the above computations simplify to
give $d(w,w') \leq d(\ell,\ell')$.  
\cqfd

\bigskip
For any metric space $Z$ and $\alpha\in\;]0,1]$, the {\it H\"older
norm}\index{holder@H\"older norm}\index{norm!holder@H\"older} of a
bounded $\alpha$-H\"older-continuous function $f:Z\ra\RR$ is
$$
\gls{normealphaholder}= \|f\|_\infty+ \|f\|'_\alpha\,,
$$ 
where 
$$
\|f\|'_\alpha=\sup_{\substack{x,\,y\,\in \,Z\\ 0<d(x,\,y)\leq 1}}
\frac{|f(x)-f(y)|}{d(x,y)^\alpha}\;.
$$ 
When the diameter of $Z$ is bounded by $1$,\footnote{This is in particular
the case for the sequence spaces of symbolic dynamical systems, see Sections
\ref{subsec:codagesimplicial} and \ref{subsec:mixingratesimpgraphs}.}
this coincides with the usual definition. Note that even if the
constant $\epsilon$ in the definition of a
$\alpha$-H\"older-continuous map is less than $1$, this norm is
finite, since
$$
\sup_{\substack{x,\,y\in Z\\ \epsilon\leq d(x,y)\leq 1}} 
\frac{|f(x)-f(y)|}{d(x,y)^\alpha}\leq 
2\;\epsilon^{-\alpha}\;\|f\|_\infty\;.
$$
Note that for all bounded $\alpha$-H\"older-continuous maps
$f,g:Z\ra\RR$, we have
\begin{equation}\label{eq:prodholdernorm}
\|fg\|_\alpha\leq\|f\|_\alpha\;\|g\|_\infty+\|f\|_\infty\;\|g\|_\alpha\;.
\end{equation}

We denote by $\gls{espacealphaholder}$ (respectively
$\gls{espacealphaholdeb}$) the space of $\alpha$-H\"older-continuous
real-valued functions with compact support (respectively which are
bounded) on $Z$, endowed with this norm. Note that $\C_{\rm b}^\alpha
(Z)$ is a real Banach space.\footnote{The standard proof using
  Arzela-Ascoli's theorem applies with our slightly different
  definition of the H\"older norms.}

\bigskip
The next two lemmas will be needed only in Part III of this book.
The first one is a metric estimate on the extension of geodesic
segments to geodesic rays, with its functional counterpart.

\blemm\label{lem:holderextensiontoray} Let $X$ be a geodesically
complete proper $\CAT(-1)$ space, let $T\geq 1$, and let $\alpha
\in\gengeod X$ be a generalised geodesic line which is isometric
exactly on $[0,T]$. For every generalised geodesic line $\rho\in
\gengeod X$ which is isometric exactly on $[0,+\infty[\,$, such that
$\rho|_{[0,T]}=\alpha|_{[0,T]}$, we have
$$
d(\alpha, \rho) =\frac{e^{-2\,T}}{4}\; < 1\;,
$$ 
and hence for all $\beta\in\;]0,1]$ and $\wt \psi\in
\C^\beta_b(\gengeod X)$,
$$
|\wt \psi(\alpha)-\wt \psi(\rho)|\leq \frac{e^{-2\,\beta\,T}}{4^\beta} \;
\|\wt \psi\,\|_\beta\;.
$$
\elemm

\dem By Equation \eqref{eq:geodesicmetric} defining the distance on
$\gengeod X$, we have, since $d(\alpha(t), \rho(t))=0$ for all
$t\in\;]-\infty,T]$ and $d(\alpha(t), \rho(t))=t-T$ otherwise,
$$ 
d(\alpha, \rho) = \int_T^{+\infty}(t-T)\,e^{-2t}\,dt= 
e^{-2\,T}\;\int_0^{+\infty}u\,e^{-2u}\,du=\frac{e^{-2\,T}}{4}\;.
$$
The result follows.
\cqfd

\medskip The second lemma gives a metric estimate, with its functional
counterpart, on the map which associates to a geodesic ray in an outer
normal unit bundle its point at infinity, emphasising the
$\frac{1}{2}$-H\"older-continuity of endpoints maps (see Lemma
\ref{lem:troisboulesmet} (4)).  We start by giving some definitions.

\medskip
Let $X$ be a geodesically complete proper $\CAT(-1)$ space, and let
$D$ be a nonempty proper closed convex subset of $X$.  The {\em
  distance-like map
$$
d_D:\big(\partial_\infty X-\partial_\infty D\big)^2\to [0,+\infty[
$$ 
associated with $D$}\index{distance-like map} is defined in
\cite[\S 2.2]{HerPau10} as follows. For $\xi,\xi'\in\partial_\infty
X-\partial_\infty D$, let $\xi_t,\xi'_t:[0,+\infty[\;\to X$ be the
geodesic rays starting at the closest points $P_D(\xi),P_D(\xi')$ to
$\xi,\xi'$ on $D$ and converging to $\xi,\xi'$ as $t\to\infty$. Let
\begin{equation}\label{eq:distancelike}
\gls{distancelike}(\xi,\xi')=\lim_{t\to+\infty}e^{\frac{1}{2}d(\xi_t,\,\xi'_t)-t}\;.
\end{equation}
The distance-like map $d_D$ is invariant by the diagonal action of the
isometries of $X$ preserving $D$. If $D$ consists of a single point
$x$, then $d_D$ is the visual distance\footnote{See
Equation \eqref{eq:defidistvis}.} $d_x$ on $\partial_\infty X$ based
at $x$.  If $D$ is a horoball with point at infinity $\xi_0$, then
$d_D$ is Hamenst\"adt's distance\footnote{See
Equation \eqref{eq:defidisthamenbord}.} on $\partial_\infty
X-\{\xi_0\}$.  As seen in \cite[\S 2.2, Ex.~(4)]{HerPau10}, if $X$ is
a metric tree, then
$$
d_D(\xi,\xi')=\begin{cases}
e^{\frac{1}{2}\;d(P_D(\xi),\,P_D(\xi'))} \;>1 & 
{\rm if}\;\; P_D(\xi)\neq P_D(\xi')\\
d_{x}(\xi,\xi')=e^{-d(x,\,y)} \;\leq 1 & 
{\rm if}\,  \begin{array}[t]{ll} 
P_D(\xi)=P_D(\xi')=x\;\;  \\
{\rm and}\;\;
[x,\xi[\;\cap \,[x,\xi'[\;=[x,y]\;.\end{array}
\end{cases}
$$
In particular, although in general it is not an actual distance on its
whole domain $\partial_\infty X-\partial_\infty D$, the map $d_D$ is
locally a distance, and we can define with the standard formula the
$\beta$-H\"older-continuity of maps with values in $(\partial_\infty
X- \partial_\infty D,d_D)$ and the $\beta$-H\"older-norm of a function
defined on $(\partial_\infty X- \partial_\infty D,d_D)$. In the next
result, we endow $\partial_\infty X- \partial_\infty D$ with the
distance-like map $d_D$.

\bprop \label{prop:partialplusholder}
Let $\XX$ be a locally finite simplicial tree without terminal
vertices, and let $\DD$ be a proper nonempty simplicial subtree of
$\XX$. Let $X=|\XX|_1$ and $D=|\DD|_1$ be their geometric
realisations. The homeomorphism $\partial^+:\normalout
\DD\ra(\partial_\infty X - \partial_\infty D)$ defined by $\rho
\mapsto \rho_+$ is $\frac{1}{2}$-H\"older-continuous, and for all
$\beta\in \;]0,1]$ and $\psi\in \C^{\beta}_b(\partial_\infty X-
\partial_\infty D)$, the map $\psi\circ\partial^+:\normalout\DD\ra\RR$ 
is bounded and $\frac{\beta}{2}$-H\"older-continuous, with
$$
\|\psi\circ\partial^+\|_{\frac{\beta}{2}} \leq 
(1+2^{\frac{\beta}{2}+1})\;\|\psi\|_{\beta}\;.
$$
\eprop

\dem
Let us prove that for every $\rho,\rho'\in \normalout \DD$, if
$d(\rho,\rho')\leq 1$, then $\rho(0)=\rho'(0)$, and
\begin{equation}\label{eq:demiholder}
d_D(\rho_+,\rho'_+)=\sqrt{2}\;d(\rho,\rho')^{\frac{1}{2}}\;.
\end{equation}
This proves that the map $\partial^+$ is
$\frac{1}{2}$-H\"older-continuous.  We may assume that $\rho\neq
\rho'$.

Let $\rho,\rho'\in\normalout\DD$. If $\rho(0)\neq \rho'(0)$, then the
images of $\rho$ and $\rho'$ are disjoint and their connecting segment
in the tree $\XX$ joins $\rho(0)$ and $\rho'(0)$; hence for every
$t\in [0, +\infty[\,$, we have
$$
d(\rho(t), \rho'(t))=d(\rho(0), \rho'(0)) + d(\rho(t), \rho(0))
+ d(\rho'(t), \rho'(0))\geq 1+2\,t\;.
$$
Thus
\begin{align*}
d(\rho,\rho')&=
\int_{-\infty}^{0} d(\rho(0),\rho'(0))\;e^{2\;t}\;dt+
\int_{0}^{+\infty} d(\rho(t),\rho'(t))\;e^{-2\;t}\;dt\\ &
\geq \int_{-\infty}^{0} e^{2\;t}\;dt+\int_0^{+\infty}(1+2t)\;e^{-2\,t}\;dt
>2\int_0^{+\infty}e^{-2\,t}\;dt=1\;.
\end{align*}

Assume that $x=\rho(0)=\rho'(0)$ and let $n$ be the length of the
intersection of $\rho$ and $\rho'$. Then
$$
d_D(\rho_+,\rho'_+)= d_x(\rho_+,\rho'_+)=
\lim_{t\ra+\infty} \;e^{\frac{1}{2}\,d(\rho(t),\, \rho'(t))-t} = e^{-n}\;.
$$
Furthermore, since $\rho(t)=\rho'(t)$ for $t\leq n$ and
$d(\rho(t), \rho'(t))=2(t-n)$ otherwise, we have, using the 
change of variables $u= 2(t-n)$,
$$
d(\rho,\rho')=\int_{n}^{+\infty} \;2\,(t-n)\;e^{-2\,t}\;dt=
e^{-2\,n}\int_{0}^{+\infty} \;u\;e^{-u}\;\frac{du}{2}=
\frac{e^{-2\,n}}{2}\;.
$$
This proves Equation \eqref{eq:demiholder}.

Let $\beta\in \;]0,1]$ and $\psi\in \C^{\beta}_b(\partial_\infty X
-\partial_\infty D)$. We have $\|\psi\circ\partial^+\|_\infty
= \|\psi\|_\infty$ since $\partial^+$ is a homeomorphism, and, by
Equation \eqref{eq:demiholder},
\begin{align*}
\|\psi\circ\partial^+\|'_{\frac{\beta}{2}} &=
\sup_{\rho,\,\rho'\in \normalout\DD,\;0<d(\rho,\,\rho')\leq 1}
\frac{|\psi\circ\partial^+(\rho)-\psi\circ\partial^+(\rho')|}
{d(\rho,\rho')^{\frac{\beta}{2}}}\\ & \leq
\sup_{\rho,\,\rho'\in \normalout\DD,\;0<d(\rho,\,\rho')\leq \frac{1}{2}}
\frac{|\psi\circ\partial^+(\rho)-\psi\circ\partial^+(\rho')|}
{d(\rho,\rho')^{\frac{\beta}{2}}}+
\frac{2\,\|\psi\circ\partial^+\|_\infty}{2^{-\frac{\beta}{2}}}\\ & \leq
\sup_{\substack{\xi,\,\xi'\in \partial_\infty X-\partial_\infty D\\
 0< d_D(\xi,\,\xi') \leq 1}}\;\;
\frac{|\psi(\xi)-\psi(\xi')|}{2^{-\frac{\beta}{2}}\;d_D(\xi,\xi')^{\beta}}
+ 2^{\frac{\beta}{2}+1}\,\|\psi\|_\infty\\ &
\leq 2^{\frac{\beta}{2}+1}\,\|\psi\|_\beta
\;.
\end{align*}
Since
$\|\psi\circ\partial^+\|_{\frac{\beta}{2}}= \|\psi\circ\partial^+\|_{\infty}
+\|\psi\circ\partial^+\|'_{\frac{\beta}{2}}$, this proves the last
claim of Proposition \ref{prop:partialplusholder}.
\cqfd

\bigskip
A stronger assumption than the H\"older regularity is the locally
constant regularity, that we now define. Alhough it is only useful for
totally disconnected metric spaces, several error terms estimates in
the literature use this stronger regularity (see for instance
\cite{AthGhoPra12, KemPauSch17} and Part \ref{sect:arithappli} of this
book).

\medskip
Let $\epsilon>0$. For every metric space $E$ and every set $E'$, we say
that a map $f:E\ra E'$ is {\it $\epsilon$-locally
  constant}\index{locally!constant} if $f$ is constant on every closed
ball of radius $\epsilon$ (or equivalently of radius at most
$\epsilon$) in $E$. We say that $f:E\ra E'$ is {\it locally
  constant}\index{locally!constant} if there exists $\epsilon>0$ such
that $f$ is $\epsilon$-locally constant.

Note that if $E$ is a geodesic metric space and $f:E\ra E'$ is locally
constant, then $f$ is constant. But when $E$ is for instance an
ultrametric space, since two distinct closed balls of the same radius
are disjoint, the above definition turns out to be very interesting
(and much used in representation theory in positive characteristic,
for instance). For example, the characteristic function
$\mathbbm{1}_A$ of a subset $A$ of $E$ is $\epsilon$-locally constant
if and only if for every $x\in A$, the closed ball $B(x,\epsilon)$ is
contained in $A$. In particular, the characteristic function of a
closed ball of radius $\epsilon$ in an ultrametric space is
$\epsilon$-locally constant.

\brema\label{rem:locconstholder} Let $E$ and $E'$ be two metric
spaces. If a map $f:E\ra E'$ is $\epsilon$-locally constant, then it
is $\alpha$-H\"older-continuous for every $\alpha\in \mathopen{]}0,
  1\mathclose{]}$. Indeed, for all $x,y\in E$, if $d(x,y)\leq
\epsilon$ then $d(f(x),f(y))=0\leq c\;d(x,y)^\alpha$ for all $c>0$. If
furthermore $E'=\RR$ and $f$ is bounded, then
$$
\sup_{x,y\in E,\,x\neq y}\frac{|f(x)-f(y)|}{d(x,y)^\alpha}=
\sup_{x,y\in E,\;d(x,\,y)>\epsilon}\frac{|f(x)-f(y)|}{d(x,y)^\alpha}
\leq\frac{2}{\epsilon^\alpha}\;\|f\|_\infty\;.
$$
\erema

\medskip 
For all $\epsilon\in\;]0,1]$ and $\beta>0$, we denote by
$\C^{\epsilon\,{\rm lc},\,\beta}_b(E)$ the real vector
space\footnote{Note that a linear combination of $\epsilon$-locally
  constant functions is again a $\epsilon$-locally constant function.}
of $\epsilon$-locally constant functions $f:E\ra\RR$ endowed with the
{\it $\epsilon\,{\rm lc}$-norm of exponent
  $\beta$}\index{elc@$\epsilon\,{\rm lc}$-norm}%
\index{norm!$\epsilon\,{\rm lc}$} defined by
$$
\|f\|_{\epsilon\,{\rm lc},\,\beta}=\;\epsilon^{-\beta}\;\|f\|_{\infty}\;.
$$
The above remark proves that if $\beta\in\;]0,1]$, then
$\|f\|_\beta\le 3\|f\|_{\epsilon\,{\rm lc},\,\beta}$, so that the
inclusion map from $\C^{\epsilon\,{\rm lc},\,\beta}_b(E)$ into
$\C^{\beta}_b(E)$ is continuous.

\section{Potentials}
\label{subsec:potentials}
  
In this book, a {\em potential}\index{potential} for $\Ga$ is a
continuous $\Ga$-invariant function $\gls{potentialup}: T^1X\to\RR$.
The quotient function $\gls{potentialdown}:\Ga\bs T^1X\to\RR$ of $\wt
F$ is called a {\em potential}\index{potential} on $\Ga\backslash
T^1X$. Precomposing by the canonical projection $\G X\ra T^1X$, the
function $\wt F$ defines a continuous $\Ga$-invariant function from
$\G X$ to $\RR$, also denoted by $\wt F$, by $\wt F(\ell)=\wt
F(v_\ell)$ for every $\ell\in \G X$.

\medskip
For all $x,y\in X$, and any geodesic line $\ell\in\G X$   such that
$\ell(0)=x$ and $\ell(d(x,y))=y$, let 
$$
\int_x^y\wt F=\int_0^{d(x,y)}\wt F(v_{\flow{t}\ell})\,dt\,.
$$
Note that for all $t\in\mathopen{]}0,d(x,y)\mathclose{[}$, the germ
$v_{\flow{t}\ell}$ is independent on the choice of such a line $\ell$,
hence $\int_x^y\wt F$ does not depend on the extension $\ell$ of
the geodesic segment $\mathopen{[}x,y\mathclose{]}$.  The following
properties are easy to check using the $\Ga$-invariance of $\wt F$
and the basic properties of integrals: For all $\ga\in\Ga$
\begin{equation}\label{eq:equivint}
\int_{\ga x}^{\ga y} \wt F=\int_x^y\wt F\,,
\end{equation}
for the antipodal map $\iota$ 
\begin{equation}\label{eq:timereversal}
\int_y^x \wt F=\int_x^y\wt F\circ\iota\;,
\end{equation}
 and, for any $z\in[x,y]$,
\begin{equation}\label{eq:additive}
\int_x^y\wt F=\int_x^z\wt F+\int_z^y\wt F\,.
\end{equation}

The {\em period}\index{period} of a loxodromic isometry $\ga$ of $X$
for the potential $\wt F$ is
$$
\per_F(\ga)=\int_x^{\ga x}\wt F
$$
for any $x$ in the translation axis of $\ga$.  Note that, for all
$\alpha\in\Ga$ and $n\in\NN-\{0\}$, we have
\begin{equation}\label{eq:periodinvconj}
\operatorname{Per}_F(\alpha\ga\alpha^{-1})=\operatorname{Per}_F(\ga),
\;\;\;\operatorname{Per}_F(\ga^n)=n\;\operatorname{Per}_F(\ga)
\;\;\;{\rm and}\;\;\; 
\operatorname{Per}_F(\ga^{-1})= \operatorname{Per}_{F\circ\iota}(\ga)\;.
\end{equation}

\bigskip 
In trees, we have the following Lipschitz-type control on the
integrals of the potentials along segments.

\blemm\label{lem:potentialcontinuity} When $\wt F$ is constant or when
$X$ is an $\RR$-tree, for all $x,x',y,y'\in X$, we have
$$
\Big|\int_x^y\wt F-\int_{x'}^{y'}\wt F\big|\le \;
d(x,x')\sup_{\pi^{-1}([x,\,x'])}|\wt F|\;+\;
d(y,y')\sup_{\pi^{-1}([y,\,y'])}|\wt F|\,.
$$
\elemm 

\dem When $\wt F$ is constant, the result follows from the triangle
inequality. 

Assume that $X$ is an $\RR$-tree. Consider the case $x=x'$. Let $z\in
X$ be such that $[x,z]=[x,y]\cap[x,y']$. Using Equation
\eqref{eq:additive} and the fact that $d(y,z)+d(z,y')=d(y,y')$, the
claim follows. The general case follows by combining this case $x=x'$
and a similar estimate for the case $y=y'$.  
\cqfd

\medskip 
Some form of uniform H\"older-type control of the potential, analogous
to the Lipschitz-type one in the previous lemma, will be crucial
throughout the present work. The following Definition
\ref{defi:HCproperty} formalises this (weaker) assumption.

\bdefi\label{defi:HCproperty} The triple $(X,\Ga,\wt F)$ satisfies the
{\em HC-property}\index{property HC} (H\"older-type control) if $\wt
F$ has subexponential growth when $X$ is not an $\RR$-tree and if
there exist $\kappa_1\geq 0$ and $\kappa_2\in\;]0,1]$ such that for
all $x,y,x',y'\in X$ with $d(x,x'),d(y,y')\leq 1$, we have
\begin{align*}\label{eq:HC}
& \Big|\int_x^y   \wt F -\int_{x'}^{y'}\wt F\;\Big|       \tag{HC}\\
\leq \;\;&
\big(\kappa_1+2\max_{\pi^{-1}(B(x,\,1)\cup B(x',\,1))}|\wt F|\,\big)
\;d(x,x')^{\kappa_2}+
\big(\kappa_1+2\max_{\pi^{-1}(B(y,\,1)\cup B(y',\,1))}|\wt F|\,\big)
\;d(y,y')^{\kappa_2}\;.                                                  
\end{align*}
\edefi

By Equation \eqref{eq:timereversal}, $(X,\Ga,\wt F\circ\iota)$
satisfies the \ref{eq:HC}-property if and only if $(X,\Ga,\wt F)$
does. By the triangle inequality $|\,d(x,y)-d(x',y')\,|\leq
d(x,x')+d(y,y')$, for every $\kappa\in\RR$, the triple $(X,\Ga,\wt
F+\kappa)$ satisfies the \ref{eq:HC}-property (up to changing the
constant $\kappa_1$) if and only if $(X,\Ga,\wt F)$ does.

When $X$ is assumed to be a Riemannian manifold with pinched sectional
curvature, requiring the potentials to be H\"older-continuous as in
\cite{PauPolSha15} is sufficient to have the \ref{eq:HC}-property, as
we will see below.  

\bprop\label{prop:exempHCprop}  
The triple $(X,\Ga,\wt F)$ satisfies the \ref{eq:HC}-property if one
of the following conditions is satisfied:
\begin{itemize}
\item $\wt F$ is constant,
\item $X$ is an $\RR$-tree,
\item $X$ is a Riemannian manifold with pinched sectional
  curvature and $\wt F$ is H\"older-conti\-nuous.
\end{itemize}
\eprop

\dem  The first two cases are treated in Lemma
\ref{lem:potentialcontinuity}, and we may take for them $\kappa_1=0$
and $\kappa_2=1$ in the definition of the \ref{eq:HC}-property.

The claim for Riemannian manifolds follows from the property of at
most linear growth of the H\"older-conti\-nu\-ous maps (see Remark
\ref{rem:soriteholder}) and from the following lemma, so that the
constants $\kappa_1>0$ and $\kappa_2\in\;]0,1]$ of the
\ref{eq:HC}-property depend only on the H\"older-continuity constants
of $\wt F$ and on the bounds on the sectional curvature of $X$.

\blemm\label{lem:technicholderPPS} 
If $X$ is a Riemannian manifold with pinched sectional curvature and
$\wt F$ is H\"older-conti\-nuous, there exist two constants $c_1>0$
and $c_2\in\mathopen{]}0,1\mathclose{]}$ such that, for all $x,y,z$ in
$X$ with $d(x,y)\leq 1$, we have
$$
\Big|\;\int_x^{z}\wt F\;-\int_{y}^{z}\wt F\;\Big|\leq 
c_1\,d(x,y)^{c_2}+ 
2\;d(x,y)^{\frac{1}{2}}\max_{\pi^{-1}(B(x,\,1)\cup B(y,\,1))}|\wt F|\;.
$$
The constants $c_1$ and $c_2$ depend only on the H\"older-continuity
constants of $\wt F$ and the bounds on the sectional curvature of $X$.
\elemm

This lemma is similar to the second claim in
\cite[Lem.~3.2]{PauPolSha15}, but the proof of this claim (and more
precisely the proof of \cite[Lem.~2.3]{PauPolSha15} used in the proof
of \cite[Lem.~3.2]{PauPolSha15}), which involves a different distance
$d$ on $\G X$, does not extend with the present definition of $d$.

\medskip\dem By symmetry, we may assume that $d(x,z)\geq d(y,z)$. The
result is true if $x=y$, hence we assume that $x\neq y$. Let $x'$ be
the point on $\mathopen{[}x,z\mathclose{]}$ at distance $d(y,z)$ from
$z$.

\begin{center}
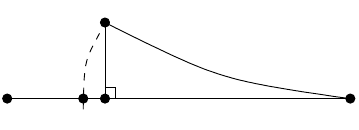
\end{center}

\noindent 
The closest point $p$ of $y$ on $\mathopen{[}x, z\mathclose{]}$ lies
in $\mathopen{[}x',z\mathclose{]}$ by convexity. Hence 
$$
d(x,x')\leq d(x,p)\leq d(x,y)\leq \sqrt{d(x,y)}\leq 1\;,
$$ 
since closest point maps do not increase distances and $d(x,y)\leq
1$. Therefore
\begin{equation}\label{eq:controlintxxprim}
\Big|\int_x^{x'}\wt F\;\Big|\;\leq
d(x,x')\max_{\pi^{-1}([x,\,x'])}|\wt F|\;\leq
\sqrt{d(x,y)}\max_{\pi^{-1}(B(x,\,1))}|\wt F|\;,
\end{equation}
Since $\int_x^{z}\wt F = \int_x^{x'}\wt F + \int_{x'}^{z}\wt F$ (see
Equation \eqref{eq:additive}), we have
\begin{equation}\label{eq:reducxegxprim}
\Big|\;\int_x^{z}\wt F\;-\int_{y}^{z}\wt F\;\Big| 
\leq
\Big|\;\int_{x'}^{z}\wt F\;-\int_{y}^{z}\wt F\;\Big| + 
\Big|\int_x^{x'}\wt F\;\Big|\;.
\end{equation}

Assume first that $d(y,z)=d(x',z)\leq\frac{1}{2}\sqrt{d(x,y)}$. We
have
\begin{align*}
&\;\Big|\;\int_{x}^{z}\wt F\;-\int_{y}^{z}\wt F\;\Big|\leq
\;\Big|\;\int_{x'}^{z}\wt F\;\Big|\;+\;\Big|\;\int_{y}^{z}\wt F
\;\Big|\;+\;\Big|\;\int_{x}^{x'}\wt F\;\Big|\\ \leq & \;
d(x',z)\max_{\pi^{-1}([x',\,z])}|\wt F|+
d(y,z)\max_{\pi^{-1}([y,\,z])}|\wt F|+
\sqrt{d(x,y)}\max_{\pi^{-1}(B(x,\,1))}|\wt F|\\ \leq & \;
2\sqrt{d(x,y)}\;\max_{\pi^{-1}(B(x,\,1)\cup B(y,\,1))}|\wt F|\;,
\end{align*}
and Lemma \ref{lem:technicholderPPS} follows, for any $c_1>0$ and
$c_2\in\mathopen{]}0,1\mathclose{]}$.

\medskip
Now assume that that $d(y,z)\geq\frac{1}{2}\sqrt{d(x,y)}$. Since the
distance function from a given point to a point varying on a geodesic
line is convex, we have $d(x',y)\leq d(x,y)$. By Equations
\eqref{eq:controlintxxprim} and \eqref{eq:reducxegxprim}, we may
therefore assume that $x=x'$ and prove that 
$$
\Big|\;\int_x^{z}\wt F-\int_{y}^{z}\wt F\;\Big|\leq c_1\,d(x,y)^{c_2}
$$ 
for appropriate constants $c_1,c_2$.

Since $X$ is a Riemannian manifold, we identify $\G X$ and $T^1X$ with
the usual unit tangent bundle of $X$ as explained previously. Let $u$
(respectively $v$) be the unit tangent vector at $x$ (respectively
$y$) pointing towards $z$. Let $\rho=d(x,z)=d(y,z)\geq \frac{1}{2}
\sqrt{d(x,y)}$, and  $t\in\mathopen{[}0,\rho \mathclose{]}$. We 
apply Lemma \ref{lem:holderdeltaded}, whose hypotheses are indeed
satisfied.

Since $t\in\mathopen{[}0,\rho \mathclose{]}$ and $d(x,y)\leq 1$, the
term on the right hand side of Equation \eqref{eq:holderdeltaded} is
bounded by $2\,c_8$. Since $\wt F$ is H\"older-continuous, let $c>0$
and $\alpha\in \mathopen{]}0, 1\mathclose{]}$ be the
H\"older-continuity constants such that $|\wt F(u')-\wt F(v')|\leq
c\; d(u',v')^\alpha$ for all $u',v'\in T^1X$ such that
$d(u',v')\leq 2\,c_8$. Then, by Lemma \ref{lem:holderdeltaded},
\begin{align*}
  \Big|\;\int_x^{z}\wt F\;-\int_{y}^{z}\wt F\;\Big|&=
  \Big|\int_0^{\rho}\big(\wt F(\flow{t}u)-\wt F(\flow{t}v)\big)dt\;\Big|
  \leq \int_0^{\rho}c \;d(\flow{t}u,\flow{t}v)^\alpha\,dt\\ &
  \leq c\;{c_8}^\alpha\;
  d(x,y)^{\alpha c_3}\; \Big(\int_0^{+\infty}e^{-\alpha t}\,dt+
  \int_{-\infty}^\rho e^{2\alpha t -2\alpha\rho}\,dt\Big)\\ &=
  \frac{3c}{2\alpha}\;{c_8}^\alpha\;d(x,y)^{\alpha c_3} \;.
\end{align*}
This concludes the proof of Lemma \ref{lem:technicholderPPS} with
$c_2=\alpha c_3$ and $c_1= \frac{3c}{2\alpha}\,c_8^\alpha$, hence
completes the proof of Proposition \ref{prop:exempHCprop} in the
Riemannian manifold case.  \cqfd $\;\Box$

\brema\label{rem:soritesholderbis} (1) If $X=\wt M$ is a Riemannian
manifold, then $T^1X$ is naturally identified with the usual
Riemannian unit tangent bundle of $\wt M$. If the potential $\wt F:
T^1\wt M\to\RR$ is H\"older-continuous for Sasaki's Riemannian metric
on $T^1\wt M$, it is a potential as in \cite{Ruelle81} and
\cite{PauPolSha15}.  Furthermore, the definition of $\int_x^y\wt F$
coincides with the one in these references.

\medskip\noindent (2) The quotient function $F$ is H\"older-continuous
when $\wt F$ is H\"older-continuous. 
\erema

\medskip
Let $\wt F,\wt F^* :T^1X\ra \RR$ be potentials for $\Ga$. We say that
$\wt F^*$ is {\em
  cohomologous}\index{cohomologous}\index{potential!cohomologous} to
$\wt F$ (see for instance \cite{Livsic72}) if there exists a
continuous $\Ga$-invariant function $\wt G :T^1X\ra \RR$, such that, for
every $\ell\in \G X$, the map $t\mapsto \wt G(v_{\flow{t}\ell})$ is
differentiable and
\begin{equation}\label{eq:cohomol}
\wt F^*(v_\ell)-\wt F(v_\ell)=\frac{d}{dt}_{\mid t=0}\wt G(v_{\flow{t}\ell})\;.
\end{equation}

When working with H\"older-continuous potentials, the regularity
requirement is for $\wt G$ to also be H\"older-continuous.  Note that
the right-hand side of Equation \eqref{eq:cohomol} does not depend on
the choice of the representative $\ell$ of its germ $v_\ell$. In
particular, $\operatorname{Per}_F(\ga)=\operatorname{Per}_{F^*}(\ga)$
for any loxodromic isometry $\ga$ if $\wt F$ and $\wt F^*$ are
cohomologous potentials.

A potential $\wt F$ is said to be {\em
  reversible}\index{reversible}\index{potential!reversible} if $\wt
F$ and $\wt F\circ\iota$ are cohomologous.

\section{Poincar\'e series and critical exponents}
\label{subsec:criticexpo}

Let us fix a potential $\wt F: T^1X\to\RR$ for $\Ga$, and $x,y\in X$.

The {\em critical exponent}\index{critical exponent} of $(\Ga,F)$ is
the element $\gls{criticalexponent}$ of the extended real line
$\mathopen{[}-\infty, +\infty\mathclose{]}$ defined by
$$
\delta =\limsup_{n\to+\infty}\frac 1n\ln\sum_{\ga\in\Ga, \,
  n-1<d(x,\ga y)\le n} e^{\int_x^{\ga y}\wt F}\;.
$$
The {\em Poincar\'e series}\index{Poincar\'e!series} of $(\Ga,F)$ is
the map $\gls{poincareseries}: \RR\ra[0,+\infty]$ defined by
$$
Q:s\mapsto\sum_{\ga\in\Ga}e^{\int_x^{\ga y}(\wt F-s)}\,.
$$
If $\delta<+\infty$, we say that $(\Ga,F)$ is {\em of divergence
  type}\index{divergence type} if the series $Q_{\Ga,\,F,\,x,\,y}(\delta)$
diverges, and {\em of convergence type}\index{convergence type}
otherwise. 

When $F=0$, the critical exponent $\delta_{\Ga,\,0}$ is the usual
critical exponent\index{critical exponent}
$\delta_\Ga\in\mathopen{]}0,+\infty\mathclose{]}$ of $\Gamma$, the
Poincar\'e series $Q_{\Ga,\,0,\,x,\,y}$ is the usual Poincar\'e
series\index{critical exponent} of $\Gamma$, and we recover the usual
notion of divergence or convergence type of $\Ga$, see for instance
\cite{Roblin03}.

The Poincar\'e series of $(\Ga,F)$ and its critical exponent make
sense even if $\Ga$ is elementary (see for instance Lemma
\ref{lem:proprielemcritexpo} (10)). The following result collects some
of the basic properties of the critical exponent.  

\blemm \label{lem:proprielemcritexpo} Assume that $(X,\Ga,\wt F)$
satisfies the \ref{eq:HC}-property.  Then
\begin{enumerate} 
\item the critical exponent $\delta_{\Ga,\,F}$ and the divergence or
  convergence of $Q_{\Ga,\,F,\,x,\,y}(s)$ are independent of the
  points $x,y\in X$; they depend only on the cohomology class of $\wt
  F$\,;
\item 
$Q_{\Ga,\,F\circ\iota,\,x,\,y}=Q_{\Ga,\,F,\,y,\,x}$ and
  $\delta_{\Ga,\,F\circ \iota}=\delta_{\Ga,\,F}$\,; in particular,
  $(\Ga,F)$ is of divergence type if and only if $(\Ga,F\circ\iota)$
  is of divergence type \,;
\item the Poincaré series $Q(s)$ diverges if $s<\delta_{\Ga,\,F}$ and
  converges if $s>\delta_{\Ga,\,F}$\,;
\item 
$\delta_{\Ga,\,F+\kappa}=\delta_{\Ga,\,F}+\kappa$ for any
  $\kappa\in\RR$, and $(\Ga,F)$ is of divergence type if and only if
  $(\Ga,F+\kappa)$ is of divergence type\,;
\item if $\Ga'$ is a nonelementary subgroup of $\Ga$, denoting by
  $F':\Ga'\backslash T^1X\ra\RR$ the map induced by $\wt F$, then $
  \delta_{\Ga',\,F'}\leq \delta_{\Ga,\,F} $\,;
\item if $\delta_\Ga<+\infty$, then $\delta_\Ga+ \inf
  \limits_{\pi^{-1}(\C\Lambda\Ga)}\wt F\le\delta_{\Ga,\,F}\le
  \delta_{\Ga}+\sup\limits_{\pi^{-1}(\C \Lambda\Ga)}\wt F$\,;
\item 
$\delta_{\Ga,\,F}>-\infty$\,;
\item the map $\wt F\mapsto \delta_{\Ga,\,F}$ is convex, sub-additive,
  and $1$-Lipschitz for the uniform norm on the vector space of real
  continuous maps on $\pi^{-1}(\C\Lambda\Ga)$\,;\footnote{That is, if
    $\wt F,\wt{F^*}:T^1X\ra \RR$ are potentials for $\Ga$ satisfying
    the \ref{eq:HC}-property, inducing $F,F^*:\Ga\backslash
    T^1X\ra\RR$, and if $\delta_{\Ga,\,F},\delta_{\Ga,\,F^*}<+\infty$,
    then $\delta_{\Ga,\,tF+(1-t)F^*}\leq
    t\,\delta_{\Ga,\,F}+(1-t)\,\delta_{\Ga,\,F^*}$ for every
    $t\in\mathopen{[}0,1\mathclose{]}$,
$$
\delta_{\Ga,\,F+F^*}\leq \delta_{\Ga,\,F}+\delta_{\Ga,\,F^*}\;,
$$
$$
|\;\delta_{\Ga,\,F^*}-\delta_{\Ga,\,F}\;|\leq 
\sup_{v\in \pi^{-1}(\C\Lambda\Ga)} |\;\wt {F^*}(v)-\wt F(v)\;|\;.
$$}
\item if $\Ga''$ is a discrete cocompact group of isometries
of $X$ such that $\wt F$ is $\Ga''$-invariant, denoting by
$F'':\Ga''\backslash T^1X\ra\RR$ the map induced by $\wt F$,
then
$$
\delta_{\Ga,\,F}\leq \delta_{\Ga'',\,F''}\;;
$$
\item
  if $\Ga$ is infinite cyclic, generated by a loxodromic isometry
  $\ga$ of $X$, then $(\Ga,F)$ is of divergence type and
$$
\delta_{\Ga,\,F}=\max\Big\{\frac{\per_F(\ga)}{\lambda(\ga)},
\frac{\per_{F\circ\iota}(\ga)}{\lambda(\ga)}\Big\}\,. 
$$
\end{enumerate}
\elemm

\dem We give details of the proofs of the statements (1)--(7) and (10)
which are the ones used in this book, only for the sake of
completeness. The proofs from \cite[Lem.~3.3]{PauPolSha15} generalise
to the current setting, replacing the use of
\cite[Lem.~3.2]{PauPolSha15} by the following consequence of the
\ref{eq:HC}-property: there exist $\kappa_1\geq 0$ and
$\kappa_2\in\;]0,1]$ such that for every $N\in\NN-\{0\}$, for all
$x,y,x',y'\in X$ with $d(x,x'),d(y,y')\leq N$, we have
\begin{equation}\label{eq:HCN}
\Big|\int_x^y   \wt F  -\int_{x'}^{y'}\wt F\;\Big| 
\leq 2N\big(\kappa_1+
2\max_{\pi^{-1}(B(x,\,N)\cup B(x',\,N)\cup B(y,\,N)\cup B(y',\,N))}|\wt F|\,\big)
\;.                                                  
\end{equation}
This is obtained from Equation \eqref{eq:HC} by subdividing the
segments $[x,x']$ and $[y,y']$ into $N$ subintervals of equal lengths
at most $1$, and by using the triangle inequality.

\medskip
(1) For $x',y'\in X$, if $N=\max\{\lceil d(x,x')\rceil,\lceil
d(y,y')\rceil\}$, by Equation \eqref{eq:HCN} and by the
$\Ga$-invariance of $\wt F$, we have, for every $\ga\in\Ga$,
$$
\Big|\int_x^{\ga y} \wt F - \int_{x'}^{\ga y'} \wt F \;\Big|\leq c=
2\,N\,\big(\kappa_1+
2\max_{\pi^{-1}(B(x,\,N)\cup B(x',N)\cup B(y,\,N)\cup B(y',N))}|\wt F|\,\big)\;,
$$
which is finite since the continuous map $\wt F$ is bounded on
compact subsets of $T^1X$. Hence by the triangle inequality, we have,
for every $s\in\RR$,
$$
e^{-c -s\,d(x,\,x')-s \,d(y,\,y')}\;Q_{\Ga,\,F,\,x,\,y}(s)
\leq Q_{\Ga,\,F,\,x',\,y'}(s)\leq
e^{c+s\,d(x,\,x')+s\,d(y,\,y')}\;Q_{\Ga,\,F,\,x,\,y}(s)\;.
$$
The first claim  of Assertion (1) follows.

\medskip
Let $\wt F^* :T^1X\ra \RR$ be a potential for $\Ga$ which is
cohomologous to $\wt F$. Let $\wt G :T^1X\ra \RR$ be a continuous
$\Ga$-invariant function satisfying Equation \eqref{eq:cohomol}. Let
$\kappa_x= \max_{\pi^{-1}(x)}|\wt G|$, which is finite by continuity.
By $\Ga$-invariance, for every $\ga\in \Ga$, we have $\kappa_{\ga y}=
\kappa_y$. For every $\ga\in \Ga$, with $\ell\in\G X$ any geodesic
line such that $\ell(0)=x$ and $\ell(d(x,\ga y))=\ga y$, we have
\begin{align*}
\Big|\int_x^{\ga y} \wt F^* - \int_{x}^{\ga y} \wt F \;\Big|
&=\Big|\int_0^{d(x,\ga y)}
\frac{d}{dt}_{\mid t=0}\wt G(v_{\flow{t}\flow{s}\ell})\,ds\Big| =
\Big|\int_0^{d(x,\ga y)}
\frac{d}{ds}\wt G(v_{\flow{s}\ell})\,ds\Big| \\ & =
\big|\;\wt G(v_{\ell})-\wt G(v_{\flow{d(x,\ga y)}\ell})\;\big| \leq
\kappa_x+\kappa_y\;.
\end{align*}
Hence by the triangle inequality, we have, for every $s\in\RR$,
$$
e^{-\kappa_x-\kappa_y }\;Q_{\Ga,\,F,\,x,\,y}(s)
\leq Q_{\Ga,\,F^*,\,x,\,y}(s)\leq
e^{\kappa_x+\kappa_y}\;Q_{\Ga,\,F,\,x,\,y}(s)\;.
$$
The second claim of Assertion (1) follows.

\medskip
(2) This assertion follows from Equations \eqref{eq:timereversal} and
\eqref{eq:equivint}, by the change of variable $\ga\mapsto \ga^{-1}$
in the summation of the Poincar\'e series.

\medskip
(3) This assertion is a standard argument of Poincaré series. For
every $s\neq\delta_{\Ga,\,F}$, let $\epsilon= \frac{1}{2}\,
|\delta_{\Ga,\,F} - s|>0$. First assume that $s>\delta_{\Ga,\,F}$. By
the definition of the critical exponent $\delta_{\Ga,\,F}$, there
exists $N\in\NN$ such that for every integer $n\geq N$,
$$
\sum_{\ga\in\Ga, \,n-1<d(x,\ga y)\le n} e^{\int_x^{\ga y}\wt F}
\leq e^{n(\delta_{\Ga,\,F}+\epsilon)}\;.
$$
Hence there exists $c>0$ such that
$$
Q(s)\leq c+\sum_{n\in\NN}e^{n(\delta_{\Ga,\,F}+\epsilon)-s\,n+|s|}
=c+e^{|s|}\sum_{n\in\NN}e^{-\epsilon n}<+\infty\;.
$$

Now assume that $s<\delta_{\Ga,\,F}$. By the definition of
$\delta_{\Ga,\,F}$, there exists an increasing sequence
$(n_k)_{k\in\NN}$ in $\NN$ such that for every $k\in\NN$, we have
$$
\sum_{\ga\in\Ga, \,n_k-1<d(x,\ga y)\le n_k} e^{\int_x^{\ga y}\wt F}
\geq e^{n_k(\delta_{\Ga,\,F}-\epsilon)}\;.
$$
Hence 
$$
Q(s)\geq \sum_{k\in\NN}e^{n_k(\delta_{\Ga,\,F}-\epsilon)-s\,n_k-|s|)}
=e^{-|s|}\sum_{k\in\NN}e^{\epsilon\, n_k}=+\infty\;.
$$
This proves Assertion (3).

\medskip
Assertions (4) and (5) are immediate by Assertion (3), since with
$\kappa$ and $\Ga',F'$ as in these assertions, for every $s\in \RR$,
we have
$$
Q_{\Ga,\,F+\kappa,\,x,\,y}(s)=Q_{\Ga,\,F,\,x,\,y}(s-\kappa)
\;\;\;{\rm and}\;\;\;
Q_{\Ga',\,F',\,x,\,y}(s)\leq Q_{\Ga,\,F,\,x,\,y}(s)\;.
$$

\medskip
(6) If $x$ is a point in the convex hull $\C\Lambda\Gamma$ of the
limit set of $\Ga$, then, for every $\ga\in\Ga$, the geodesic segment
between $x$ and $\ga x$ is contained in $\C\Lambda\Gamma$. Hence
$$
d(x,\ga x)\;\Big(\inf_{\pi^{-1}(\C\Lambda\Ga)}\wt F-s\Big)\leq
{\int_x^{\ga x} (\wt F-s)}\leq d(x,\ga x)\;
\Big(\sup_{\pi^{-1}(\C\Lambda\Ga)}\wt F -s\Big)\;.
$$
This proves Assertion (6) by taking the exponential, summing over
$\ga\in\Ga$ with $n-1<d(x,\ga y)\le n$, taking the logarithm, dividing
by $n$ and taking the upper limit as $n$ tends to $+\infty$.

\medskip
(7) Let $\Ga'$ be a nonelementary convex-cocompact subgroup of $\Ga$
(for instance a Schottky subgroup of $\Ga$, which exists since $\Ga$
is nonelementary). Denote by $F':\Ga'\bs T^1X\ra\RR$ the
map induced by $\wt F$. Since $|\wt F\,|$ is $\Ga$-invariant and bounded
on compact subsets of $T^1X$, by Assertion (6), we have
$\delta_{\Ga',\,F'}>-\infty$ as $\delta_\Gamma\ge 0$. Assertion (7) then follows from
Assertion (5).

\medskip
(10) For every $s\in\RR$, if $x$ belongs to the translation axis of
$\ga$, we have by Equation \eqref{eq:periodinvconj}
\begin{align*}
\sum_{\alpha\in\Ga} \;
e^{\int_x^{\alpha x} (\wt F-s)}& =\sum_{n\in\NN} \;
e^{\int_x^{\ga^n x} (\wt F-s)}+\sum_{n\in\NN-\{0\}} \;
e^{\int_x^{\ga^{-n} x} (\wt F-s)}\\ &=\sum_{n\in\NN} \;
e^{n(\operatorname{Per}_F(\ga)-s\,\ell(\ga))}+\sum_{n\in\NN-\{0\}} \;
e^{n(\operatorname{Per}_{F\circ\iota}(\ga)-s\,\ell(\ga))}\,.
\end{align*}
Hence $Q_{\Ga,\,F,\,x,\,x}(s)$ converges if and only if
$\operatorname{Per}_F(\ga)-s\,\ell(\ga)<0$ and
$\operatorname{Per}_{F\circ\iota}(\ga)-s\,\ell(\ga)<0$. Let
$\overline{\delta}= \max\big\{\frac{\per_F(\ga)}{\lambda(\ga)},
\frac{\per_{F\circ\iota}(\ga)}{\lambda(\ga)}\big\}$. By Assertion (3),
letting $s$ tend to $\delta_{\Ga,\,F}$ on the right gives that $
\delta_{\Ga,\,F}\geq \overline{\delta}$, and letting $s$ tend to
$\overline{\delta}$ on the left gives the reverse inequality.  The
above computation also gives that $Q_{\Ga,\,F,\,x,\,x}(\,\overline{\delta}\,)$
diverges, which proves Assertion (10) and concludes the proof of Lemma
\ref{lem:proprielemcritexpo}.  \cqfd

\medskip
\bexems (1) If $\delta_\Ga$ is finite and $\wt F$ is bounded, then the
critical exponent $\delta$ is finite by Lemma
\ref{lem:proprielemcritexpo} (6).

\smallskip\noindent(2) If $X$ is a Riemannian manifold with pinched
negative curvature or when $X$ has a compact quotient, then
$\delta_\Ga$ is finite. See for instance \cite{Bourdon95}.

\smallskip\noindent(3) There are examples of $(X,\Ga)$ with
$\delta_\Ga=+\infty$ (and hence $\delta=+\infty$ if $\wt F$ is
constant), for instance when $X$ is the complete ideal hyperbolic
triangle complex with $3$ ideal triangles along each edge, see
\cite{GabPau01}, and $\Ga$ its isometry group.  Hence the finiteness
assumption of the critical exponent is nonempty in general. For the
type of results treated in this book, it is however natural and
essential.  
\eexems

\medskip
We may replace upper limits by limits in the definition of the
critical exponents, as follows.

\btheo Assume that $(X,\Ga,\wt F)$ satisfies the
\ref{eq:HC}-property.  If $c>0$ is large enough, then
$$
\delta =\lim_{n\to+\infty}\;\frac 1n\;\ln\sum_{\ga\in\Ga, \,
  n-c<d(x,\ga y)\le n} e^{\int_x^{\ga y}\wt F}\;.
$$
If $\delta>0$, then
$$
\delta =\lim_{n\to+\infty}\;\frac 1n\;\ln\sum_{\ga\in\Ga, \,
  d(x,\ga y)\le n} e^{\int_x^{\ga y}\wt F}\;.
$$
\etheo

\dem 
The proofs of \cite[Theo.~4.2 and Theo.~4.3]{PauPolSha15}, either
using the original arguments of \cite{Roblin02} valid when $\wt F$ is
constant, or the super-multiplicativity arguments of
\cite{DalPeiSam12}, extend, using the \ref{eq:HC}-property (see
Definition \ref{defi:HCproperty}) instead of
\cite[Lem.~3.2]{PauPolSha15}.  \cqfd

\bigskip 
In what follows, we fix a potential $\wt F$ for $\Ga$
such that $(X,\Ga,\wt F)$ satisfies the HC-property.  We
define $\wt F^+=\wt F$ and $\wt F^-=\wt F\circ \iota$, we denote by
$F^\pm:\Ga\bs T^1X\ra\RR$ their induced maps, and we assume that
$\delta= \delta_{\Ga,\,F^+}=\delta_{\Ga,\,F^-}$ is finite.

\section{Gibbs cocycles}
\label{subsec:Gibbscocycle}

The {\em (normalised) Gibbs cocycle}\index{Gibbs!cocycle}%
\index{cocycle!Gibbs} associated with the group $\Ga$ and the
potential $\wt F^\pm$ is the map $C^\pm=\gls{cocycle}
: \partial_\infty X\times X\times X\to\RR$ defined by
$$
(\xi,x,y)\mapsto C^\pm_\xi(x,y) =\lim_{t\to+\infty}\int_y^{\xi_t}(\wt
F^\pm-\delta)-\int_x^{\xi_t}(\wt F^\pm-\delta)\;,
$$
where $t\mapsto\xi_t$ is any geodesic ray with endpoint
$\xi\in\partial_\infty X$. 

We will prove in Proposition \ref{prop:continuGibbscocycle} below that
this map is well defined, that is, the above limits exist for all
$(\xi,x,y)\in\partial_\infty X\times X\times X$ and they are
independent of the choice of the geodesic rays $t\mapsto \xi_t$.  If
$\wt F^\pm=0$, then $C^-=C^+=\delta_{\Ga} \beta$, where $\beta$ is the
Busemann cocycle.  If $X$ is an $\RR$-tree, then
\begin{equation}\label{eq:cocycletreecase}
C^\pm_\xi(x,y)=\int_y^p(\wt F^\pm-\delta)-\int_x^{p}(\wt F^\pm-\delta)\;, 
\end{equation}
where $p\in X$ is the point for which $\mathopen{[}\,p, \xi
\mathclose{[} = \mathopen{[}x,\xi\mathclose{[}\, \cap\, \mathopen{[}y,
\xi\mathclose{[}\,$; in particular, the map $\xi\mapsto
C^\pm_\xi(x,y)$ is locally constant on the totally discontinuous space
$\partial_\infty X$.

The Gibbs cocycles satisfy the following equivariance and cocycle
properties: For all $\xi\in\partial_\infty X$ and $x,y,z\in X$, and
for every isometry $\ga$ of $X$, we have
\begin{equation}\label{eq:cocycle}
 C^\pm_{\ga \xi}(\ga x,\ga y)= C^\pm_\xi(x,y)\; \textrm{ and }\;
 C^\pm_\xi(x,z)+ C^\pm_\xi(z,y)= C^\pm_\xi(x,y)\,.
\end{equation}
For every $\ell\in\G X$, for all $x$ and $y$ on the image of the
geodesic line $\ell$, if $\ell_-,x,y,\ell_+$ are in this order on
$\ell$, we have
\begin{equation}\label{eq:changemoinsplus}
  C^-_{\ell_-}(x,y)=C^+_{\ell_+}(y,x)=-C^+_{\ell_+}(x,y)=
\int_x^y(\wt F^+-\delta)\;.
\end{equation}

\bprop\label{prop:continuGibbscocycle} Assume that $(X,\Ga,\wt F)$
satisfies the \ref{eq:HC}-property and that $\delta<+\infty$.
\begin{enumerate}
\item The maps $C^\pm:\partial_\infty X\times X\times X\to\RR$ are
  well-defined.
\item\label{eq:cocycleLip} With the constants $\kappa_1,\kappa_2$ of
  the \ref{eq:HC}-property, for all $x,y\in X$ and $\xi\in
  \partial_\infty X$, if we assume that $d(x,y)\leq 1$, then
\begin{equation*}
|\,C^\pm_\xi(x,y)\,|\le (\kappa_1+2\,|\delta|+
2 \max_{\pi^{-1}(B(x,\,1)\cup B(y,\,1))}|\wt F|)\,d(x,y)^{\kappa_2}\;,
\end{equation*}
and, in general, if $N=\lceil d(x,y)\rceil$, then 
$$
|\,C^\pm_\xi(x,y)\,|\le N\,(\kappa_1+2\,|\delta|+
2 \max_{\pi^{-1}(B(x,\,N)\cup B(y,\,N))}|\wt F|)\;.
$$
If $X$ is an $\RR$-tree, then for all $x,y\in X$ and
$\xi\in \partial_\infty X$, we have
\begin{equation*}
|\,C^\pm_\xi(x,y)\,|\le d(x,y)\max_{\pi^{-1}([x,\,y])}|\wt F^\pm-\delta|\,.
\end{equation*}
\item\label{eq:holdercontpotential} The maps $C^\pm:\partial_\infty
  X\times X\times X\to\RR$ are locally H\"older-continuous (and locally
  Lipschitz when $X$ is an $\RR$-tree). In particular, they are
  continuous.
\item \label{eq:cocycleombre} For all $r>0$, $x,y\in X$ and
  $\xi\in\partial_\infty X$, if $\xi$ belongs to the shadow
  $\OOO_xB(y,r)$ of the ball $B(y,r)$ seen from $x$, then with the
  constants $\kappa_1,\kappa_2$ of the \ref{eq:HC}-property, if $r\leq
  1$, we have
\begin{equation*}
\Big|\,C^\pm_\xi(x,y)+\int_x^y(\wt F^\pm-\delta)\,\Big|\le 2(\kappa_1+
2\,|\delta|+ 2\max_{\pi^{-1}(B(y,\,2))}|\wt F|)\,r^{\kappa_2}\;,
\end{equation*} and in general
$$
\Big|\,C^\pm_\xi(x,y)+\int_x^y(\wt F^\pm-\delta)\,\Big|\le 
2\, \lceil r\rceil\,(\kappa_1+
2\,|\delta|+2\max_{\pi^{-1}(B(y,\,2\, \lceil r\rceil))}|\wt F|)\;.
$$
If $X$ is an $\RR$-tree, then 
\begin{equation*}
\Big|\,C^\pm_\xi(x,y)+\int_x^y(\wt F^\pm-\delta)\,\Big|\le 
2r \max_{\pi^{-1}(B(y,\,r))}|\wt F^\pm-\delta|\,.
\end{equation*}
\end{enumerate}
\eprop

\dem 
(1) The fact that $C^\pm_\xi(x,y)$ is well defined when $X$ is an
$\RR$-tree follows from Equation \eqref{eq:cocycletreecase}.

When $X$ is not an $\RR$-tree, let $\rho:t\mapsto\xi_t$ be any
geodesic ray with endpoint $\xi\in\partial_\infty X$, let $t\mapsto
x_t$ (respectively $t\mapsto y_t$) be the geodesic ray from $x$
(respectively $y$) to $\xi$. Let $t_x=\beta_\xi(x,\xi_0)$ and
$t_y=\beta_\xi(y,\xi_0)$, so that the quantity $\beta=t_y-t_x$ is
equal to $\beta_\xi(y,x)$ (which is independent of $\rho$), and for
every $t$ big enough, we have $\beta_\xi(\xi_t,x_{t+t_x})=
\beta_\xi(\xi_t,y_{t+t_y})= 0$. 

\begin{center}
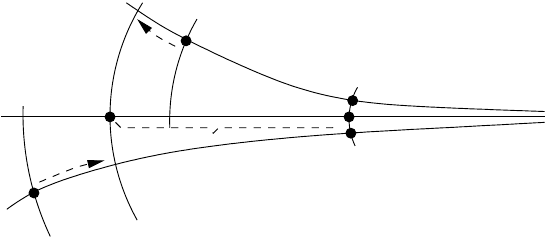
\end{center}

Since $X$ is $\CAT(-1)$, if $t$ is big enough, then the distances
$d(\xi_t,x_{t+t_x})$ and $d(\xi_t,y_{t+t_y})$ are at most one, and
converge, in a nonincreasing way, exponentially fast to $0$ as
$t\ra+\infty$. For $s\geq 0$, let $a_s=\int_y^{y_s}(\wt F^\pm-\delta)
-\int_x^{x_{s-\beta}}(\wt F^\pm-\delta)$ (which is independent of
$\rho$).  We have, using Equation \eqref{eq:HC}, and the fact that
$B(x',1)\cup B(y',1)\subset B(x',2)$ if $d(x',y')\leq 1$,
\begin{align*}
&\Big|\Big(\int_y^{\xi_t}(\wt F^\pm-\delta)-
\int_x^{\xi_t}(\wt F^\pm-\delta)\Big)
-a_{t+t_y}\Big|\\=\;& \Big|\Big( \int_y^{\xi_t}
(\wt F^\pm-\delta)-\int_y^{y_{t+t_y}}(\wt F^\pm-\delta)\Big)+ 
\Big(\int_x^{x_{t+t_x}}(\wt F^\pm-\delta)-
\int_x^{\xi_t}(\wt F^\pm-\delta)\Big)\Big|\\ 
\leq\;& 2(\kappa_1 +2\max_{\pi^{-1}(B(\xi_t,\,2))}|\wt F-\delta|\,)\;
\max\{d(\xi_t,x_{t+t_x}),d(\xi_t,y_{t+t_y})\}^{\kappa_2}\;,
\end{align*}
which converges to $0$ since $\wt F$ has subexponential growth by the
assumptions of the \ref{eq:HC}-property.  Hence in order to prove
Assertion (1), we only have to prove that $\lim_{s\ra+\infty} a_s$
exists.

For all $s\geq t\geq |\beta|$, we have, by the additivity of the
integral along geodesics (see Equation \eqref{eq:additive}) and by
using again Equation \eqref{eq:HC},
\begin{align*}
  &|a_t-a_s| = \Big|\int_{y_t}^{y_s}(\wt F^\pm-\delta)-
  \int_{x_{t-\beta}}^{x_{s-\beta}}(\wt F^\pm-\delta)\Big|\\  \leq\;&
  (\kappa_1 + 2\max_{\pi^{-1}(B(x_{t-\beta},\,2))}|\wt
  F-\delta|)\,d(y_t,x_{t-\beta})^{\kappa_2}+(\kappa_1+
   2\max_{\pi^{-1}(B(y_{s},\,2))}|\wt F-\delta|\,)\,
  d(y_s,x_{s+\beta})^{\kappa_2}\;.
\end{align*}
Again by the subexponential growth of $\wt F$, the above expression
converges (exponentially fast, for future use) to $0$ as $t\ra
+\infty$ uniformly in $s$, hence $\lim_{s\ra+\infty} a_s$ exists by a
Cauchy type argument.

\medskip\noindent (2) Let $(\xi,x,y) \in \partial_\infty X\times X
\times X$.  Assertion \eqref{eq:cocycleLip} of Proposition
\ref{prop:continuGibbscocycle} follows from Equation
\eqref{eq:cocycletreecase} when $X$ is an $\RR$-tree, since $[x,y]=
[x,p]\cup[p,y]$ where $p$ is the closest point to $y$ on $[x,\xi[$.
When $X$ is not an $\RR$-tree, the first claim of Assertion
\eqref{eq:cocycleLip} follows immediately from the
\ref{eq:HC}-property of $(X,\Ga,\wt F^\pm-\delta)$, and the second
claim from this \ref{eq:HC}-property and the same subdivision argument
of the geodesic segment $[x,y]$ into $\lceil d(x,y)\rceil$ 
subintervals of equal lengths (at most $1$), as in the 
proof of Equation \eqref{eq:HCN}.

\medskip\noindent (3) Let $(\xi,x,y),(\xi',x',y') \in \partial_\infty
X\times X \times X$. By the cocycle property \eqref{eq:cocycle}, we
have
\begin{equation}\label{eq:cocycluseforcont}
|C^\pm_\xi(x,y)-C^\pm_{\xi'}(x',y')|\leq |C^\pm_\xi(x,y)-C^\pm_{\xi'}(x,y)|
+|C^\pm_{\xi'}(x,x')|+|C^\pm_{\xi'}(y,y')|\;.
\end{equation}
First assume that $X$ is an $\RR$-tree. Let $K$ be a compact subset of
$X$, and let $$\epsilon_K=\inf_{x,y\in K} e^{-d(x,\,x_0)-d(x,\,y)}>0\,.$$
\begin{center}
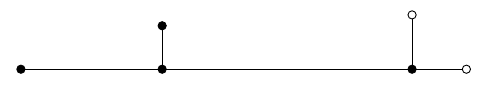
\end{center}
Let $p,q$ be the points in $X$ such that $[x,\xi[ \;\cap\,[y,\xi[\;=
\,[p,\xi[$ and $[x,\xi[ \;\cap\,[x,\xi'[\;=\,[x,q]$. If
$d_{x_0}(\xi,\xi')\leq \epsilon_K$, then by the definition of the
visual distance (see Equation \eqref{eq:distviscastree}) and by
Equation \eqref{eq:lipequivdistvis}, we have
$$
e^{-d(x,\,q)}= d_x(\xi,\xi')  \leq e^{d(x,\,x_0)}\;d_{x_0}(\xi,\xi')
\leq e^{-d(x,\,y)}\leq e^{-d(x,\,p)}\,.
$$
In particular $q\in [p,\xi[$, so that $[x,\xi'[ \;\cap\,[y,\xi'[\;
=\,[p,\xi'[$. Thus by Equation \eqref{eq:cocycletreecase}, we have
$$
|\,C^\pm_\xi(x,y) - C^\pm_{\xi'}(x,y)\,| = 0\;.
$$ 
Therefore, by Equation \eqref{eq:cocycluseforcont} and by the
$\RR$-tree case of Assertion (2), if $d_{x_0}(\xi, \xi') \leq
\epsilon_K$, if $x,y\in K$ and $d(x,x'),d(y,y') \leq 1$, then
\begin{equation*}
|C^\pm_\xi(x,y)-C^\pm_{\xi'}(x',y')|\leq 
d(x,x')\max_{\pi^{-1}([x,\,x'])}|\wt F^\pm-\delta|+
d(y,y')\max_{\pi^{-1}([y,\,y'])}|\wt F^\pm-\delta|\;.
\end{equation*}
Since $\wt F$ is bounded on compact subsets of $T^1X$, this proves
that $C^\pm$ is locally Lipschitz.

\medskip 
Let us now consider the case when $X$ is general. For all distinct
$\xi,\xi' \in \partial_\infty X$, let $t\mapsto\xi_t$ and $t\mapsto
\xi'_t$ be the geodesic rays from $x_0$ to $\xi$ and $\xi'$
respectively.  By the end of the proof of Assertion (1), for every
compact subset $K$ of $X$, there exist $a_1,a_2>0$ such that for every
$x,y\in K$, we have for all $\eta\in\{\xi,\xi'\}$ and $t\geq 0$,
$$
\Big|C^\pm_\eta(x,y)-\Big(\int_y^{\eta_t}(\wt F^\pm-\delta)-
\int_x^{\eta_t}(\wt F^\pm-\delta)\Big)\Big|\leq a_1\,e^{-a_2t}\;.
$$ 
Let $T=-\frac{1}{2}\,\ln d_{x_0}(\xi,\xi')$. If $T\geq 0$, by the
properties of $\CAT(-1)$-spaces (see Equation \eqref{eq:easyexoPPS}
for the second inequality), there exist two constants $a_3,a_4>0$ such
that $d(\xi_{2T},\xi'_{2T})\leq a_3$ and $d(\xi_{T},\xi'_{T}) \leq
a_4\,e^{-T}$.  Hence if $d_{x_0}(\xi,\xi')\leq \min\big\{\frac{1}{a_4^2},
1\big\}$ (so that $T\geq 0$ and $d(\xi_{T}, \xi'_{T}) \leq 1$), we have,
using Equation \eqref{eq:HC} for the last inequality,
\begin{align*}
&|\,C^\pm_\xi(x,y)-C^\pm_{\xi'}(x,y)\,|\\ \leq\;& 
\Big|\int_y^{\xi_T}(\wt F^\pm-\delta)-
\int_y^{\xi'_T}(\wt F^\pm-\delta)\Big|
+\Big|\int_x^{\xi_T}(\wt F^\pm-\delta)-
\int_x^{\xi'_T}(\wt F^\pm-\delta)\Big|+2 \,a_1\,e^{-a_2T}\\ \leq &
2\,(\kappa_1+ 2\max_{\pi^{-1}(B(\xi_T,\,2))}|\wt F^\pm-\delta|)\, 
d(\xi_T,\xi'_T)^{\kappa_2}+ 2\,a_1\,e^{-a_2T}\;.
\end{align*}
By the subexponential growth of $\wt F$, there exists $a_5>0$ such
that
$$
|\,C^\pm_\xi(x,y)-C^\pm_{\xi'}(x,y)\,| \leq 
a_5\,e^{-\frac{\kappa_2}{2}\,T}+ 2\,a_1\,e^{-a_2T}\leq 
(a_5+2\,a_1)\,d_{x_0}(\xi,\xi')^{\min\{\frac{\kappa_2}{4},\frac{a_2}{2}\}}\;.
$$
We now conclude from Equation \eqref{eq:cocycluseforcont} and
Assertion (2) as in the end of the above tree case that $C^\pm$ is
locally H\"older-continuous.

\medskip\noindent 
(4) Let $r>0$, $x,y\in X$ and $\xi\in\partial_\infty X$ be such that
$\xi\in\OOO_xB(y,r)$. Let $p$ be the closest point to $y$ on
$[x,\xi[\,$, so that $d(p,y)\leq r$. By Equations
\eqref{eq:changemoinsplus} and \eqref{eq:cocycle}, we have
\begin{align}
\Big|\,C^\pm_\xi(x,y)+\int_x^y(\wt F^\pm-\delta)\,\Big| &= 
\Big|\,C^\pm_\xi(x,y)-C^\pm_\xi(x,p)-\int_x^p(\wt F^\pm-\delta)+
\int_x^y(\wt F^\pm-\delta)\,\Big| \nonumber\\ & \leq
|\,C^\pm_\xi(p,y)|+\Big|\int_x^p(\wt F^\pm-\delta)-
\int_x^y(\wt F^\pm-\delta)\,\Big|\;.\label{eq:assertion4cocyle}
\end{align}

First assume that $X$ is an $\RR$-tree.  Then by Assertion (2) and by
Lemma \ref{lem:potentialcontinuity}, we deduce from Equation
\eqref{eq:assertion4cocyle} that
\begin{align*}
\Big|\,C^\pm_\xi(x,y)+\int_x^y(\wt F^\pm-\delta)\,\Big|
 \leq 2\,d(y,p)\max_{\pi^{-1}([y,p])}|\wt F^\pm-\delta|\leq 
2r \max_{\pi^{-1}(B(y,\,r))}|\wt F^\pm-\delta|\,. 
\end{align*}

\medskip In the general case, the result then follows similarly from
Equation \eqref{eq:assertion4cocyle} by using Assertion (2) and the
\ref{eq:HC}-property if $r\leq 1$ of Equation \eqref{eq:HCN} in
general.  
\cqfd

\section{Systems of conductances on trees and generalised 
electrical networks}
\label{subsec:cond}

Let $(\XX,\lambda)$ be a locally finite metric tree without terminal
vertices, let $X=|\XX|_\lambda$ be its geometric realisation, and let
$\Ga$ be a nonelementary discrete subgroup of $\Isom(\XX,\lambda)$.

Let $\wt c:E\XX\to\RR$ be a $\Ga$-invariant function, called a {\em
  system of (logarithmic) conductances}\index{conductance}%
\index{system of conductances} for $\Ga$. We denote by $c:\Ga\bs
E\XX\to\RR$ the function induced by $\wt c$, which we also
call a {\em system of conductances}\index{conductance}%
\index{system of conductances} on $\Ga\bs\XX$.

Classically, an {\em electric network}\footnote{A potential in this
  work is not the analog of a potential in an electric network, we
  follow the dynamical systems terminology as in for example
  \cite{PauPolSha15}.}  (without sources or reactive elements) is a
pair $(G,e^c)$, where $G$ is a graph and $c:EG\ra \RR$ a function,
such that $c$ is {\em reversible}\index{conductance!reversible}%
\index{system of conductances!reversible}: $c(e)= c(\overline{e})$ for
all $e\in EG$, see for example \cite{Nash-Williams59},
\cite{Zemanian91}. In this text, we do not assume our system of
conductances $\wt c$ to be reversible. In Chapter \ref{sec:laplacian},
we will even sometimes assume that the system of conductances is {\em
  antireversible}%
\index{antireversible}\index{system of conductances!antireversible},
that is, satisfying $c(\overline{e})=-\,c(e)$ for every $e\in E\XX$.

Two systems of conductances $\wt c,\wt {c'}:E\XX\to\RR$ are said
to be {\em cohomologous},\index{cohomologous}%
\index{system of conductances!cohomologous} if there exists a
$\Ga$-invariant map $f:V\XX\ra \RR$ such that 
$$
\wt {c'}-\wt c=df\;,
$$
where for all $e\in E\XX$, we have
$$
df(e)=\frac{f(t(e))-f(o(e))}{\lambda(e)}\;.
$$

\bprop \label{prop:integpotconduct} Let $\wt c:E\XX\to\RR$ be a system
of conductances for $\Ga$.  There exists a potential $\wt F$ on $T^1X$
for $\Ga$ such that for all $x,y\in V\XX$, if $(e_1,\dots,e_n)$ is the
edge path in $\XX$ without backtracking such that $x=o(e_1)$ and
$y=t(e_n)$, then
$$
\int_x^y\wt F=\sum_{i=1}^{n}\wt c(e_i)\,\lambda(e_i)\,.
$$
\eprop

\dem Any germ $v\in T^1X$ determines a unique edge $e_v$ of the tree
$\XX$, the first one into which it enters: if $\ell$ is any geodesic
line whose class in $T^1 X$ is $v$, the edge $e_v$ is the unique edge
of $\XX$ containing $\pi(v)$ whose terminal vertex is the first vertex
of $\XX$ encountered at a positive time by $\ell$.  The function
$\wt F:T^1X\to\RR$ defined by
\begin{equation}\label{eq:defpotfromconduc}
\wt F(v)=\frac{4\, \wt c(e_v)}{\lambda(e_v)}
\min\big\{d(\pi(v),\,o(e_v)),d(\pi(v),t(e_v))\big\}
\end{equation}
is a (indeed $\Ga$-invariant) potential on the $\RR$-tree $X$, with
$\wt F(v)=0$ if $\pi(v)\in V\XX$.

Let us now compute $\int_x^y\wt F$, for all $x,y\in X$.  For every
$\lambda>0$, let $\psi_\lambda: [0,\lambda] \ra \RR$ be the continuous
map defined by $\psi_\lambda(t)= \frac{t^2}2$ if $t\in [0,\frac\lambda
  2]$ and $\psi_\lambda(t)= \frac{\lambda^2}4- \frac{(\lambda-t)^2}2$
if $t \in[\frac\lambda 2, \lambda]$. Let $(e_0,e_1, \dots, e_n)$ be
the edge path in $\XX$ without backtracking such that $x\in e_0
-\{t(e_0)\}$ and $y\in e_n- \{o(e_n)\}$. An easy computation shows
that
$$
\int_x^y\wt F\;=\;
\sum_{i=0}^{n-1}\wt c(e_i)\,\lambda(e_i)+
\frac{4\, \wt c(e_n)}
{\lambda(e_n)}\,\psi_{\lambda(e_n)}\big(d(y,o(e_n))\big)-
\frac{4\, \wt c(e_0)}
{\lambda(e_0)}\,\psi_{\lambda(e_0)}\big(d(x,o(e_0))\big)\,.
$$
If $x$ and $y$ are vertices, the expression simplifies to the
sum  of the lengths of the edges weighted by the conductances. 
\cqfd

\medskip 
We denote by $\gls{potentialconduc}$ the potential defined by Equation
\eqref{eq:defpotfromconduc} in the above proof, and by
$\gls{potentialconducdown}:\Ga\bs T^1X\ra\RR$ the induced
potential. Note that $F_c$ is bounded if $c$ is bounded. We call $\wt
F_c$ and $F_c$ the potentials {\em associated
  with}\index{potential!associated with} the system of conductances
$\wt c$ and $c$.  This is by no means the unique potential with the
property required in Proposition \ref{prop:integpotconduct}. The
following result proves that the choice is unimportant.

Given a potential $\wt F:T^1X\ra \RR$ for $\Ga$, let us define a map
$\gls{conducpotential}:E\XX\ra\RR$ by
\begin{equation}\label{eq:deficsubF}
\wt c_F:e\mapsto \wt c_F(e)=
\frac{1}{\lambda(e)}\;\int_{o(e)}^{t(e)}\wt F\;.
\end{equation}
Note that $\wt c_F$ is $\Ga$-invariant, hence it is a system of
conductances for $\Ga$. We denote by $\gls{conducpotentialdown}:\Ga\bs
E\XX\to\RR$ the function induced by $\wt c_F:E\XX\to\RR$. Note that
$\wt c_{F+\kappa}=\wt c_{F}+\kappa$ for every constant $\kappa\in\RR$,
that $\wt c_F$ is bounded if $\wt F$ is bounded, and that $c_{F_c}=c$
by the above proposition.

\bprop\label{prop:relatpotentialconductance}\mbox{}  
\begin{enumerate}
\item Every potential (resp.~bounded potential) for $\Ga$ is
  cohomologous to a potential (resp.~bounded potential) associated
  with a system of conductances for $\Ga$.
\item If two systems of conductances $\wt{c'}$ and $\wt c$ are
  cohomologous, then their associated potentials $\wt F_{c'}$ and
  $\wt F_c$ are cohomologous.
\item If $\XX$ has no vertex of degree $2$, if two potentials $\wt
  F^*$ and $\wt F$ for $\Ga$ are cohomologous, then the systems of
  conductances $\wt c_{F^*}$ and $\wt c_F$ for $\Ga$ are cohomologous.
\end{enumerate}
\eprop

Hence if $\XX$ has no vertex of degree $2$, the map $[F]\mapsto [c_F]$
from the set of cohomology classes of potentials for $\Ga$ to the set
of cohomology classes of systems of conductances for $\Ga$ is
bijective, with inverse $[c]\mapsto [F_c]$.

\medskip \dem (1) Let $\wt F$ be a potential for $\Ga$, and let $\wt
F^*=\wt F_{c_F}$ be the potential associated with the system of
conductances $\wt c_F$. Note that if $\wt F$ is bounded, so is $\wt
c_F$ by Equation \eqref{eq:deficsubF}, hence $\wt F^*$ is bounded by
Equation \eqref{eq:defpotfromconduc}. For all $e\in E\XX$ and $t\in
\;]0,\lambda(e)[\,$, let $v_{e,\,t}\in T^1X$ be the germ of any
geodesic line passing at time $0$ through the point of $e$ at distance
$t$ from $o(e)$. Let $\wt G:T^1X\ra \RR$ be the map defined by
$\wt G(v)=0$ if $\pi(v)\in V\XX$ and such that for all $e\in E\XX$ and
$t\in \;]0,\lambda(e)[\,$,
$$
\wt G(v_{e,\,t}) = \int_0^t(\wt F^*(v_{e,\,s})-\wt F(v_{e,\,s}))\;ds\;.
$$ 
Since $\int_{o(e)}^{t(e)}\wt F=\lambda(e)\;\wt c_F(e)$ by the
construction of $\wt c_F$ and $\lambda(e)\;\wt c_F(e)=
\int_{o(e)}^{t(e)}\wt F^*$ by Proposition \ref{prop:integpotconduct},
the map $\wt G:T^1X\ra \RR$ is continuous.  Let $\ell$ be a geodesic
line. The map $t\mapsto \wt G(v_{\flow{t}\ell})$ is obviously
differentiable at time $t=0$ if $\pi(\ell)\notin V\XX$, with
derivative $\wt F^*(v_{\ell})-\wt F(v_{\ell})$. By considering right
and left derivatives of this map at $t=0$, by using the fact that
$\int_{o(e)}^{t(e)}(\wt F-\wt F^*)=0$ for every $e\in E\XX$, and by
the continuity of $\wt F$ and $\wt F^*$ at such a point, this is still
true if $\pi(\ell)\in V\XX$. Hence $\wt F^*$ and $\wt F$ are
cohomologous, and this proves the first claim.

\medskip\noindent (2) Assume that $\wt{c'}$ and $\wt c$ are
cohomologous systems of conductances for $\Ga$, and let $f:V\XX\ra\RR$
be a $\Ga$-invariant function such that $\wt{c'}-\wt c=df$.  Define
$\wt G(v)=f(\pi(v))$ if $\pi(v)\in V\XX$. For all $e\in E\XX$ and
$t\in \;]0,\lambda(e)[\,$, define
$$
\wt G(v_{e,\,t}) =
\int_0^t(\wt F_{c'}(v_{e,\,s})-\wt F_c(v_{e,\,s}))\;ds +f(o(e))\;,
$$ 
which is $\Ga$-invariant. Its limit as $t\ra 0$ is $f(o(e))$
(independent of the edge $e$ with given origin), and its limit as
$t\ra \lambda(e)$ is, by the construction of $\wt F_c$ and $\wt
F_{c'}$,
$$
\lambda(e)\big(\wt{c'}(e)-\wt c(e)\big)+f(o(e))=
\lambda(e)\;df(e)+f(o(e))=f(t(e))
$$ 
(independent of the edge $e$ with given extremity). This proves that
$\wt G$ is continuous. One checks as in Assertion (1) that $\wt F_{c'}$
and $\wt F_c$ are cohomologous.

\medskip\noindent 
(3) In order to prove the third claim, assume that $\wt F^*$
and $\wt F$ are two cohomologous potentials for $\Ga$, and let $\wt G
:T^1X\ra \RR$ be as in the definition of cohomologous potentials, see
Equation \eqref{eq:cohomol}. By the continuity of $\wt G$, for all
elements $v$ and $v'$ in $T^1X$ such that $\pi(v)=\pi(v')\in V\XX$, we
have $\wt G(v)= \wt G(v')$, since (by the assumption on the degrees of
vertices) the two edges (possibly equal) into which $v$ and $v'$ enter
can be extended to geodesic lines with a common negative subray. Hence
for every $x\in V\XX$, the value $f(x)=\wt G(v_x)$ for every $v_x\in
T^1X$ such that $\pi(v_x)=x$ does not depend on the choice of
$v_x$. The map $f:V\XX\ra \RR$ thus defined is $\Ga$-invariant.  With
the above notation and by Equation \eqref{eq:deficsubF}, we hence
have, for every $e\in E\XX$,
\begin{align*}
\wt c_{F^*}(e)-\wt c_F(e)&=\frac{1}{\lambda(e)}\;
\int_0^{\lambda(e)} \big(\wt F^*(v_{e,\,t})-\wt F(v_{e,\,t})\big)\,dt
=\frac{1}{\lambda(e)}\;
\int_0^{\lambda(e)} \frac{d}{dt} \wt G(v_{e,\,t})\,dt\\ &=
\frac{1}{\lambda(e)}\;\big(\wt G(v_{t(e)})-
\wt G(v_{o(e)})\big)=\frac{f(t(e))-f(o(e))}{\lambda(e)}=df(e)\;.
\end{align*}
Hence $\wt c_{F^*}$ and $\wt c_F$ are cohomologous.
\cqfd 

\medskip
Given a metric tree $(\XX,\lambda)$, we define the {\em critical
  exponent}\index{critical exponent} of a $\Ga$-invariant system of
conductances $\wt c:E\XX\ra\RR$ (or of the induced system of
conductances $c:\Ga\bs E\XX\ra\RR$) as the critical exponent of
$(\Ga,\wt F_c)$ where $\wt F_c$ is the potential for $\Ga$ associated
with $\wt c$~:
$$
\delta_c=\delta_{\Ga,\,F_c}\;.
$$
By Proposition \ref{prop:relatpotentialconductance} (2) and Lemma
\ref{lem:proprielemcritexpo} (1), this does not depend on the choice
of a potential $\wt F_c$ satisfying Proposition
\ref{prop:integpotconduct}.

\chapter{Patterson-Sullivan and Bowen-\-Margulis 
measures with potential on $\CAT(-1)$ spaces}
\label{sect:measures}

Let $X,x_0,\Ga$ be as in the beginning of Section
\ref{subsec:catmoinsun},\footnote{That is, $X$ is a geodesically
  complete proper $\CAT(-1)$ space, $x_0\in X$ is a basepoint, and
  $\Ga$ is a nonelementary discrete group of isometries of $X$.} and
let $\wt F$ be a potential for $\Ga$. From now on, we assume that the
triple $(X,\Ga,\wt F)$ satisfies the \ref{eq:HC}-property of
Definition \ref{defi:HCproperty} and that the critical exponent
$\delta= \delta_{\Ga,\,F^\pm}$ is finite.

In this chapter, we discuss geometrically and dynamically relevant
measures on the boundary at infinity of $X$ and on the space of
geodesic lines $\G X$. We extend the theory of Gibbs measures from the
case of manifolds with pinched negative sectional curvature treated in
\cite{PauPolSha15} \footnote{See also the previous works
  \cite{Ledrappier95b, Hamenstadt97, Coudene03, Mohsen07}.} to
$\CAT(-1)$ spaces with the \ref{eq:HC}-property.

\section{Patterson densities}
\label{subsec:Pattersondens}

A family $\gls{Pattersondensity}$ of finite nonzero (positive Borel)
measures on $\partial_\infty X$, whose support is $\Lambda\Ga$, is a
{\em (normalised) Patterson density}\index{Patterson
  density}\index{measure!Patterson} for the pair $(\Ga,\wt F^\pm)$ if
\begin{equation}\label{eq:equivarPatdens}
\ga_*\mu^\pm_x=\mu^\pm_{\ga x}
\end{equation}
for all $\ga\in\Ga$ and $x\in X$, and if the following Radon-Nikodym
derivatives exist for all $x,y\in X$ and satisfy for (almost) all
$\xi\in\partial_\infty X$
\begin{equation}\label{eq:quasinivarPatdens}
\frac{d\mu^\pm_x}{d\mu^\pm_y}(\xi) = e^{-C^\pm_\xi(x,\,y)}\,.
\end{equation}
In particular, the measures $\mu^\pm_{x}$ are in the same measure class
for all $x\in X$, and, by Proposition \ref{prop:continuGibbscocycle},
they depend continuously on $x$ for the weak-star convergence of
measures. Note that a Patterson density for $(\Ga,F^\pm)$ is also a
Patterson density for $(\Ga,F^\pm+s)$ for every $s\in\RR$, since the
definition involves only the normalised potential $\wt F^\pm- \delta$.
If $F=0$, we get the usual notion of a Patterson-Sullivan density (of
dimension $\delta_\Ga$) for the group $\Ga$, see for instance
\cite{Patterson76, Sullivan79, Nicholls89, Coornaert93, Bourdon95,
  Roblin03}.

\bprop\label{prop:existPatterson} There exists at least one Patterson
density for the pair $(\Ga,\wt F^\pm)$. 
\eprop

\dem The Patterson construction (see \cite{Patterson75},
\cite{Coornaert93}) modified as in \cite{Ledrappier95b} (with a
multiplicative rather than additive parameter $s$), \cite{Mohsen07}
and \cite[Section 3.6]{PauPolSha15} (all in the Riemannian manifold
case) gives the result, and we give a proof only for the sake of
completeness.

We start by an independent lemma, generalizing
\cite[Lem.~3.1]{Patterson76} with a similar proof.

\blemm \label{lem:pattersonh} 
Let $\delta'\in\RR$. Let $(a_n)_{n\in\NN}$ and $(b_n)_{n\in\NN}$ be
sequences of positive real numbers such that $\lim_{n\ra+\infty}a_n=
+\infty$ and the generalised Dirichlet series $\sum_{n\in\NN} b_n
\,a_n^{-s}$ converges if $s> \delta'$ and diverges if $s<\delta'$.
Then there exists a positive nondecreasing map $h$ on
$\mathopen{]}0,+\infty\mathclose{[}$ such that

$\bullet$~ for every $\epsilon>0$, there exists $r'_\epsilon> 0$
    such that $h (t'r')\leq {t'}^{\,\epsilon}\,h (r')$ for all
    $t'> 1$ and $r'\geq r'_\epsilon$;

$\bullet$~ the series $\sum_{n\in\NN} \;b_n\; a_n^{-s} \;h(a_n)$
    converges if and only if $s>\delta'$.
\elemm

\dem Let $t_0=0$, $t_1=1$ and $h_1:\mathopen{]}0,1\mathclose{]}\ra
  \mathopen{[}1, +\infty \mathclose{[}$ the constant map $1$. Let us
define by induction on $n\in\NN-\{0\}$ a positive real number $t_n$
and a continous map $h_n:\mathopen{]}t_{n-1},t_n\mathclose{]} \ra
\mathopen{[}1, +\infty\mathclose{[}\,$. If $t_n$ and $h_n$ are
constructed, let $t_{n+1}\in\RR$ be such that
$$
\frac{h_n(t_n)}{t_n^{1/n}}\; 
\sum_{k\in \NN \;:\; t_n<a_k\leq t_{n+1}} b_k \,a_k^{-\delta'+1/n} \geq 1\;,
$$ 
which exists since the generalised Dirichlet series diverges at
$s=\delta'-1/n$. For every $t\in \mathopen{]}t_{n},t_{n+1}
\mathclose{]}$, let $h_{n+1}(t)= h_{n}(t_n)\;\big(\frac{t}{t_n}
\big)^{1/n}$.

Note that the sequence $(t_n)_{n\in \NN}$ is increasing. Let $h:
\mathopen{]}0,+\infty\mathclose{[}\ra \mathopen{[}1, +\infty
\mathclose{[}$ be the map equal to $h_n$ on 
$\mathopen{]} t_{n-1},t_n\mathclose{]}$ for every $n\in\NN-\{0\}$.  
The map $h$ is positive, continuous and nondecreasing. For every
$\epsilon>0$, let $n=\lceil\frac{1}{\epsilon}\rceil$. Since $\ln h
(t)$ is continuous and piecewise affine in $\ln t$ with slopes at most
$\frac{1}{n}$ on $\mathopen{[}t_n,+\infty\mathclose{[}$, we have
$\ln(h(t'r')) -\ln h(r')\leq \frac{1}{n} \ln t'\leq \epsilon \ln t'$
if $t'>1$ and $r'\geq t_n$, which proves the first claim on $h$.

We have 
$$
\sum_{n\in\NN} \;b_n\; a_n^{-\delta'} \;h(a_n) =
\sum_{n\in\NN} \;\sum_{\substack{k\in \NN\\t_n<a_k\leq t_{n+1}}}
\;b_k\; a_k^{-\delta'} \;h(t_n)
\big(\frac{a_k}{t_n}\big)^{1/n}\geq \sum_{n\in\NN} 1=+\infty\;.
$$
If $s>\delta'$, let $\epsilon=\frac{s-\delta'}{2}>0$. Since by
construction $h(t)=\bigO(t^\epsilon)$ as $t\ra+\infty$, there exists a
constant $C>0$ such that $h(a_n)\leq C \,a_n^\epsilon$ for every $n\in
\NN$, and the convergence of the series $\sum_{n\in\NN} \;b_n\;
a_n^{-s}\;h(a_n)$ follows from the convergence of the generalised
Dirichlet series $\sum_{n\in\NN} \;b_n\; a_n^{-\delta-\epsilon}$, thus
proving the second claim on $h$.  
\cqfd

\medskip
Now, for every $z\in X$, recall that $\Delta_z$
denotes the unit Dirac mass at $z$. Let
$h^\pm :\mathopen{[}0,+\infty\mathclose{[}\;\ra
    \mathopen{]}0,+\infty\mathclose{[}$ be a nondecreasing map such
    that

$\bullet$~ for every $\epsilon>0$, there exists $r_\epsilon\geq 0$ such
that $h^\pm (t+r)\leq e^{\epsilon\, t}h^\pm (r)$ for all $t\geq 0$ and $r\geq
r_\epsilon$;

$\bullet$~ if $\overline{Q}_{x}(s)=\sum_{\ga\in\Ga} \;
e^{\int_x^{\ga x_0} (\wt F^\pm -s)}\;h^\pm (d(x,\ga x_0))$, then $\overline{Q}_{x}(s)$
diverges if and only if the inequality $s\leq \delta$ holds.

\noindent 
If $(\Ga, F^\pm )$ is of divergence type, we may take $h^\pm =1$
constant. Otherwise, the existence of $h^\pm$ follows from Lemma
\ref{lem:pattersonh}.\footnote{Let $(\ga_n)_{n\in\NN}$ be an
  enumeration of the elements of $\Ga$, take $a_n=e^{d(x,\ga_n x_0)}$
  and $b_n=e^{\int_x^{\ga_n x_0}\wt F^\pm}$, and then take $h^\pm
  =h\circ\exp$ for $h$ the map given by Lemma \ref{lem:pattersonh}.}
For all $s>\delta$ and $x\in X$, define the measure
$$
\mu^\pm_{x,\,s}= \frac{1}{\overline{Q}_{x_0}(s)}\sum_{\ga\in\Ga} \;
e^{\int_x^{\ga x_0} (\wt F^\pm -s)}\;h^\pm (d(x,\ga x_0))\;\Delta_{\ga x_0}
$$ on $X$. By compactness for the weak-star topology of the space of
probability measures on the compact space $X\cup \partial_\infty X$,
there exists a sequence $(s_k)_{k\in\NN}$ in $\mathopen{]}\delta,
  +\infty \mathclose{[}$ converging to $\delta$ such that the sequence
of probability measures $(\mu_{x_0,\,s_k})_{k\in\NN}$ weak-star
converges to a probability measure $\mu_{x_0}$ on $X \cup
\partial_\infty X$.  Since $\overline{Q}_{x_0} (\delta)= +\infty$ and
since the support of $\mu_{x_0,\,s}$ in $X\cup \partial_\infty X$ is
equal to $\overline{\Ga x_0}=\Ga x_0\,\cup \,\Lambda\Ga$, the support
of $\mu_{x_0}$ is contained in $\Lambda\Ga$, hence equal to
$\Lambda\Ga$ by minimality.  The Radom-Nikodym derivative
$\frac{d\mu^\pm_{x,\,s}}{d\mu^\pm_{x_0,\,s}}$ is the map with support
$\Ga x_0$ defined by
\begin{equation}\label{eq:construcPS}
\frac{d\mu^\pm_{x,\,s}}{d\mu^\pm_{x_0,\,s}}(\ga x_0)= e^{\int_x^{\ga x_0} 
(\wt F^\pm-s)-\int_{x_0}^{\ga x_0}(\wt F^\pm-s)}\;
\frac{h^\pm (d(x,\ga x_0))}{h^\pm (d(x_0,\ga x_0))}\;.
\end{equation}
For every $k\in\NN$, if $d(x_0,\ga x_0)$ is large enough, then 
$$
h^\pm (d(x,\ga x_0))\leq
h^\pm (d(x,x_0)+d(x_0,\ga x_0))\leq e^{(s_k-\delta)d(x,\,x_0)}\;h^\pm
(d(x_0,\ga x_0))\;.
$$ 
By the \ref{eq:HC}-property, as $k\ra+ \infty$, the right hand side
of Equation \eqref{eq:construcPS} with $s=s_k$ converges to
$e^{-C^\pm_\xi(x,\,x_0)}$ uniformly in $\xi\in \Lambda\Ga$ as $\ga
x_0$ tends to $\xi$.  Therefore, as $k\ra+ \infty$, the measures
$\mu^\pm_{x,\,s_k}$ converge to a (finite nonzero) measure
$\mu^\pm_{x}$ with support $\Lambda\Ga$ such that
$\frac{d\mu^\pm_{x}}{d\mu^\pm_{x_0}} (\xi)= e^{-C^\pm_\xi(x,\,x_0)}$.
Let $\ga\in\Ga$. Since $\ga_*\Delta_z=\Delta_{\ga \,z}$ for every
$z\in X$, a change of variable in the summation defining
$\mu^\pm_{x,\,s}$ gives that $\ga_*\mu^\pm_{x,\,s}=\mu^\pm_{\ga
  x,\,s}$.  By the continuity of pushforwards of measures, we have
$\ga_*\mu^\pm_{x}=\mu^\pm_{\ga x}$.  By the cocycle properties of the
Radom-Nikodym derivatives and of the Gibbs cocycles, the family
$(\mu^\pm_{x})_{x\in X}$ is a (normalised) Patterson density for
$(\Ga,\wt F^\pm)$.  
\cqfd

\medskip
We refer to Theorem \ref{theo:HTSR} for the uniqueness up to scalar
multiple of the Patterson density when $(\Ga,F^\pm)$ is of divergence
type and to \cite[Coro.~17.1.8]{DasSimUrb14arxiv} for a
characterisation of the uniqueness when $F=0$.

\bigskip The Patterson densities satisfy the following extension of
the classical Sullivan shadow lemma (which gives the claim when $\wt
F$ is constant, see \cite{Roblin03}), and its corollaries. 

If $\mu$ is a positive Borel measure on a metric space $(X,d)$, the
triple $(X,d,\mu)$ is called a {\em metric measure
  space}\index{measured metric space}. A metric measure space
$(X,d,\mu)$ is {\em doubling}%
\index{doubling measure}\index{measure!doubling} if there exists
$c\geq 1$ such that, for all $x\in X$ and $r>0$,
$$
\mu(B(x,2\,r))\leq c\,\mu(B(x,r))\,.
$$
Note that, up to changing $c$, the number $2$ may be replaced by any
constant larger than $1$.  See for instance \cite{Heinonen01} for
more details on doubling metric measure spaces.  We refer for instance
to \cite[Ex.~17.4.12]{DasSimUrb14arxiv} for examples of nondoubling
Patterson(-Sullivan) measures, and to
\cite[Prop.~17.4.4]{DasSimUrb14arxiv} for a characterisation of the
doubling property of the Patterson measures when $\Ga$ is
geometrically finite and $F=0$.

A family $((X,\mu_i,d_i))_{i\in I}$ of positive Borel measures $\mu_i$
and distances $d_i$ on a common set $X$ is called {\it uniformly
  doubling}\index{uniformly doubling}%
\index{doubling measure!uniformly} if there exists $c\geq 1$ such 
that, for all $i\in I$, $x\in X$ and $r>0$,
$$
\mu_i(B_{d_i}(x,2\,r)) \leq c\,\mu_i(B_{d_i}(x,r))\,.
$$

\blemm\label{lem:shadowlemma}%
\index{shadow!lemma} Let $(\mu^\pm_{x})_{x\in X}$ be a Patterson
density for the pair $(\Ga,F^\pm)$, and let $K$ be a compact subset of
$X$.
\begin{enumerate}
\item {\rm [Mohsen's shadow lemma]} If $R$ is large enough, there
exists $C>0$ such that for all $\ga\in\Ga$ and $x,y\in K$,
$$
\frac{1}{C}\;e^{\int_x^{\ga y} (\wt F^\pm-\delta)}\leq 
\mu^\pm_{x}\big(\OOO_xB(\ga y,R)\big)
\leq C\;e^{\int_x^{\ga y} (\wt F^\pm-\delta)}\;.
$$
\item For all $x,y\in X$, there exists $c>0$ such that for
every $n\in\NN$ 
$$
\sum_{\ga\in\Ga\;:\;n-1< d(x,\,\ga y)\leq n}\;
e^{\int_x^{\ga y} (\wt F^\pm-\delta)}\leq c\;.
$$
\item For every $R>0$ large enough, there exists $C=C(R)>0$ such that
for all $\ga \in \Ga$ and all $x,y\in K$ 
$$
\mu^\pm_x(\OOO_{x}B(\ga y,5R))\le C\,\mu^\pm_x(\OOO_xB(\ga y,R))\;.
$$
\item If $\Ga$ is convex-cocompact, then the metric measure space
  $(\Lambda\Ga, d_x, \mu^\pm_x)$ is doubling for every $x\ in X$, and
  the family of metric measure spaces $((\Lambda\Ga, d_x,
  \mu^\pm_x))_{x\in C\Lambda \Ga}$ is uniformly doubling.
 \end{enumerate}
\elemm

\dem For the first assertion, the proof of
\cite[Lem.~3.10]{PauPolSha15} (see also \cite[Lem.~4]{Coudene03} with
the multiplicative rather than additive convention, as well as
\cite{Mohsen07}) extends, using Proposition
\ref{prop:continuGibbscocycle} \eqref{eq:cocycleLip},
\eqref{eq:cocycleombre} instead of \cite[Lem.~3.4 (i),
  (ii)]{PauPolSha15}. The second assertion is similar to the one of
\cite[Lem.~3.11 (i)]{PauPolSha15}, and the proof of the last two
assertions is similar to the one of \cite[Prop.~3.12]{PauPolSha15},
using Lemma \ref{lem:comparshadowball} instead of
\cite[Lem.~2.1]{HerPau10}. The uniformity in the last assertion
follows from the compactness of $\Ga\bs C\Lambda \Ga$ and the
invariance and continuity properties of the Patterson densities.  
\cqfd

\section{Gibbs measures}
\label{subsec:gibbsmeasure}

We fix from now on two Patterson densities $(\mu^\pm_x)_{x\in X}$ for
the pairs $(\Ga,F^\pm)$.

The {\em Gibbs measure}\index{Gibbs!measure}\index{measure!Gibbs}
$\gls{gibbsmeasure}$ on $\G X$ (associated with this ordered pair of
Patterson densities) is the measure $\wt m_{F}$ on $\G X$ given by the
density
\begin{equation}\label{eq:defigibbs}
d\wt m_{F}(\ell)=
e^{C_{\ell_-}^-(x_0,\,\ell(0))\,+\,C^+_{\ell_+}(x_0,\,\ell(0))}\;
d\mu_{x_0}^-(\ell_-)\,d\mu^+_{x_0}(\ell_+)\,dt
\end{equation} 
in Hopf's parametrisation with respect to the basepoint $x_0$.  The
Gibbs measure $\wt m_{F}$ is independent of $x_0$ by Equations
\eqref{eq:quasinivarPatdens} and \eqref{eq:cocycle}. Hence it is
invariant under the action of $\Ga$ by Equations
\eqref{eq:equivarPatdens} and \eqref{eq:cocycle}. It is invariant
under the geodesic flow by construction (and the invariance of
Lebesgue's measure under translation). Thus,\footnote{See for instance
  \cite[\S 2.6]{PauPolSha15} for the precautions in order to push
  locally forward an invariant measure by an orbifold covering, since
  the group $\Ga$ does not necessarily act freely on $\G X$.}  it
defines a measure $\gls{gibbsmeasuredown}$ on $\Ga\backslash\G X$
which is invariant under the quotient geodesic flow, called the {\em
  Gibbs measure}\index{Gibbs!measure}\index{measure!Gibbs} on
$\Ga\backslash\G X$ (associated with the above ordered pair of
Patterson densities). If $F=0$ and the Patterson densities
$(\mu^+_x)_{x\in X}$ and $(\mu^-_x)_{x\in X}$ coincide, the Gibbs
measure $m_{F}$ coincides with the {\em Bowen-Margulis
  measure}\index{Bowen-Margulis measure}%
\index{measure!Bowen-Margulis} $m_{\rm BM}$ on $\Ga\bs \G X$
(associated with this Patterson density), see for instance
\cite{Roblin03}.

\medskip 
\brema\label{rem:iiiiii} 
(i) The (positive Borel) measure given by the density
\begin{equation}\label{eq:geodesicurrent}
d\lambda(\xi,\eta)= 
e^{C_{\xi}^-(x_0,\,p)\,+\,C^+_{\eta}(x_0,\,p)}\;
d\mu^-_{x_0}(\xi)\,d\mu^+_{x_0}(\eta)
\end{equation}
on $\partial^2_\infty X$ is (by the same arguments as above)
independent of $p\in\,]\xi,\eta[\,$, locally finite and invariant
under the diagonal action of $\Ga$ on $\partial^2_\infty X$. It is a
{\em geodesic current}\index{geodesic!current} for the action of $\Ga$
on the Gromov-hyperbolic proper metric space $X$ in the sense of
Ruelle-Sullivan-Bonahon, see for instance \cite{Bonahon91} and
references therein.

\medskip
(ii) Another parametrisation of $\G X$ also depending on the choice of
the basepoint $x_0$ in $X$, is the map from $\G X$ to
$\partial_\infty^2 X \times \RR$ sending $\ell$ to $(\ell_-,
\ell_+,s)$ where now $s= \beta_{\ell_-}(\pi(\ell),x_0)$ (one may also
use the different time parameter $s= \beta_{\ell_+}(x_0,\pi(\ell))$).
For every $(\eta,\xi)$ in $\partial_\infty^2 X$, let
$p_{\eta,\,\xi}$ be the closest point to $x_0$ on the geodesic line
between $\eta$ and $\xi$.
\begin{center}
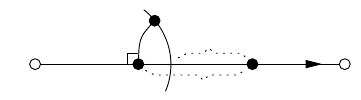
\end{center}
For every $\ell\in \G X$, with $(\ell_-, \ell_+,t)$ the original
Hopf parametrisation, since
$$
s-t=\beta_{\ell_-}(\ell(0),x_0)-\beta_{\ell_-}(\ell(0),p_{\ell_-,\ell_+})=
\beta_{\ell_-}(p_{\ell_-,\ell_+},x_0)
$$ 
depends only on $\ell_-$ and $\ell_+$, the measures
$d\mu^-_{x_0}(\ell_-) \,d\mu^+_{x_0}(\ell_+)\,dt$ and
$d\mu^-_{x_0}(\ell_-)\, d\mu^+_{x_0}(\ell_+) \, ds$ are equal. Hence using
the above variant of Hopf's parametrisation does not change the Gibbs
measures $\wt m_F$ and $m_F$.

\medskip
(iii) Since the time reversal map $\iota$ is $(\ell_-,\ell_+,t)\mapsto
(\ell_+,\ell_-, -t)$ in Hopf's coordinates, the measure $\iota_*\wt m_F$
is the Gibbs measure on $\G X$ associated with the switched pair of
Patterson densities $\big((\mu^+_{x})_{x\in X},(\mu^-_{x})_{x\in X}
\big)$ (and similarly on $\Ga\bs\G X$).  
\erema

\subsection*{The Gibbs property of Gibbs measures}
\addcontentsline{toc}{subsection}{The Gibbs property of Gibbs
  measures} 

Let us now indicate why the terminology of Gibbs measures is indeed
appropriate. This explanation will be the aim of Proposition
\ref{prop:gibbsgibbs}, but we need to give some definitions first.

\medskip
For all $\ell\in \G X$ and $r>0$, $T,T'\geq 0$, the {\em dynamical (or
  Bowen) ball}%
\index{dynamical!ball}\index{Bowen ball} around $\ell$ is
$$
\gls{dynamball}=\big\{\ell'\in \G X\;:\;\; 
\sup_{t\in\mathopen{[}-T',T\mathclose{]}}\; d(\ell(t),\ell'(t)) <r\big\}\;.
$$
Bowen balls have the following invariance properties: for all
$s\in[-T',T]$ and $\ga\in\Ga$,
$$
\flow{s} B(\ell;T,T',r)= B(\flow{s} \ell;T-s,T'+s,r)
\;\;\;{\rm and}\;\;\;
\ga B(\ell;T,T',r)= B(\ga \ell;T,T',r)\;.
$$ 
The following inclusion properties of the dynamical balls are
immediate: If $r\leq s$, $T\geq S$, $T'\geq S'$, then $B(\ell;T,T',r)$
is contained in $B(\ell;S,S',s)$.  The dynamical balls are almost
independent on $r$~: For all $r'\geq r>0$, there exists
$T_{r,\,r'}\geq 0$ such that for all $\ell\in \G X$ and $T,T'\geq 0$,
the dynamical ball $B(\ell; T+T_{r,\,r'},T'+T_{r,\,r'},r')$ is
contained in $B(\ell;T,T',r)$. This follows from the properties of
long geodesic segments with endpoints at bounded distance in a
$\CAT(-1)$-space.

For every $\ell\in \Ga\bs\G X$, let us define $B(\ell;T,T',r')$ as the
image by the canonical projection $\G X\ra\Ga\bs\G X$ of $B(\wt
\ell;T,T',r')$, for any preimage $\wt \ell$ of $\ell$ in $\G X$.

\medskip 
A $(\flow t)_{t\in\RR}$-invariant measure $m'$ on $\Ga\bs \G
X$ {\em satisfies the Gibbs property}%
\index{Gibbs!property}\index{measure!satisfying the Gibbs property}
for the potential $F$ with constant $c(F)\in\RR$ if for every compact
subset $K$ of $\Ga\bs \G X$, there exist $r>0$ and $c_{K,\,r}\geq 1$
such that for all large enough $T,T'\geq 0$, for every $\ell \in
\Ga\bs \G X$ with $\flow{-T'}\ell,\flow{T}\ell\in K$, we have
$$
\frac{1}{c_{K,\,r}}\leq \frac{m'(B(\ell;T,T',r))}
{e^{\int_{-T'}^T (F( v_{\flow{t}\ell})-c(F))\,dt}}\leq c_{K,\,r}\;.
$$

We refer to \cite[Sect.~3.8]{PauPolSha15} for equivalent variations on
the definition of the Gibbs property $m_F$. The following result shows
that the Gibbs measures indeed satisfy the Gibbs property on the
dynamical balls of the geodesic flow, thereby justifying the name.
We refer for instance to \cite[Sect.~3.8]{PauPolSha15} for the
explanations of the connection with symbolic dynamics mentioned in the
introduction.  See also Proposition \ref{prop:Gibbsproptrees} for a
discussion of the case when $X$ is a simplicial tree -- here the
correspondence with symbolic dynamics is particularly clear.

\bprop \label{prop:gibbsgibbs} Let $m_F$ be the Gibbs measure on
$\Ga\bs \G X$ associated with a pair of Patterson densities
$(\mu^\pm_{x})_{x\in X}$ for $(\Ga,\wt F^\pm)$. Then $m_F$
satisfies the Gibbs property for the potential $F$, with constant
$c(F)= \delta$.  
\eprop

\dem The proof is similar to the one of \cite[Prop.~3.16]{PauPolSha15}
(in which the key Lemma 3.17 uses only $\CAT(-1)$ arguments), up to
replacing \cite[Lem.~3.4 (1)]{PauPolSha15} by Proposition
\ref{prop:continuGibbscocycle} (2).  
\cqfd

\subsection*{The Hopf-Tsuji-Sullivan-Roblin theorem}
\addcontentsline{toc}{subsection}{The Hopf-Tsuji-Sullivan-Roblin theorem}

The basic ergodic properties of the Gibbs measures are summarised in
the following result. The case when $\wt F$ is constant is due to
\cite{Roblin03}, see also \cite[\S 6]{BurMoz96}.

\btheo[Hopf-Tsuji-Sullivan-Roblin]\label{theo:HTSR}
\index{theorem@Theorem!of Hopf-Tsuji-Sullivan-Roblin}
\index{Hopf-Tsuji-Sullivan-Roblin theorem} The following conditions
are equivalent
\begin{enumerate}
\item[{\rm\;\;(i)\;\;\;}] The pair $(\Ga,F)$ is of divergence type.
\item[{\rm (ii)}$^-$] The conical limit set of $\Ga$ has positive
  measure with respect to $\mu^-_x$ for some (equivalently every)
  $x\in X$.
\item[{\rm(ii)}$^+$] The conical limit set of $\Ga$ has positive
  measure with respect to $\mu^+_x$ for some (equivalently every)
  $x\in X$.
\item[{\rm(iii)}\;] The dynamical system $(\partial^2_\infty X, \Ga,
  (\mu_x^-\otimes \mu_x^+)_{\mid \partial^2_\infty X})$ is ergodic for
  some (equivalently every) $x\in X$.
\item[{\rm(iv)}\;] The dynamical system $(\partial^2_\infty X, \Ga,
  (\mu_x^-\otimes \mu_x^+)_{\mid \partial^2_\infty X})$ is
  conservative for some (equivalently every) $x\in X$.
\item[{\rm (v)}\;] The dynamical system $(\Ga\bs\G X, (\flow
  t)_{t\in\RR}, m_{F})$ is ergodic.
\item[{\rm (vi)}\;] The dynamical system $(\Ga\bs\G X, (\flow
  t)_{t\in\RR}, m_{F})$ is conservative.
\end{enumerate}
 
\medskip
\noindent If one of the above conditions is satisfied, then 
\begin{enumerate}
\item the measures $\mu^\pm_x$ have no atoms for any  $x\in X$,
\item the diagonal of $\partial_\infty X\times\partial_\infty X$ has
  measure $0$ for $\mu^-_x \otimes \mu^+_x$,
\item the Patterson densities $(\mu^\pm_x)_{x\in X}$ are unique up to
  a scalar multiple, and
\item for all $x,y\in X$, as $n\ra+\infty$,
$$
\max_{\ga\in\Ga,\;n-1<d(x,\ga y)\leq n} e^{\int_x^{\ga y} \wt F^\pm} = 
\smallo(e^{\delta \,n})\;.
$$ 
\end{enumerate}  
\etheo

\dem The proof\footnote{The proof occupies about 16 pages in
  \cite{PauPolSha15}, hence we cannot reproduce it in this book.} of
  the equivalence claim is similar to the one of
  \cite[Theo.~5.4]{PauPolSha15}, using

$\bullet$~ Proposition \ref{prop:continuGibbscocycle}
\eqref{eq:cocycleLip}, \eqref{eq:cocycleombre} instead of
\cite[Lem.~3.4]{PauPolSha15},

$\bullet$~ the \ref{eq:HC}-property instead of
\cite[Lem.~3.2]{PauPolSha15},

$\bullet$~ Lemma \ref{lem:proprielemcritexpo} (2) instead of
\cite[Lem.~3.3 (ii)]{PauPolSha15},

$\bullet$~ Lemma \ref{lem:shadowlemma} (2) instead of \cite[Coro.~3.11
  (i)]{PauPolSha15},

$\bullet$~ Lemma \ref{lem:shadowlemma} (3) instead of \cite[Coro.~3.12
  (i)]{PauPolSha15},

$\bullet$~ Remark \ref{rem:iiiiii}  (ii) instead of
\cite[Rem.~(ii), \S 3.7]{PauPolSha15},

$\bullet$~ Lemma \ref{lem:shadowlemma} (1) instead of
\cite[Lem.~3.10]{PauPolSha15}.  

\noindent Claims (1) and (4) are proved as in
\cite[Prop.~5.13]{PauPolSha15}, Claim (2) is proved as in
\cite[Prop.~5.5 (c)]{PauPolSha15}, and Claim (3) is proved as in
\cite[Coro.~5.12]{PauPolSha15}.  \cqfd

\medskip The following corollary follows immediately from Poincaré's
recurrence theorem and the Hopf-Tsuji-Sullivan-Roblin theorem, see
\cite[Coro.~5.15, Theo.~5.4 (ii')-(iii')]{PauPolSha15} for the
arguments written for the manifold case, which extend.

\bcoro \label{coro:finitudeGibbsdivuniq} If $m_{F}$ is finite, then
\begin{enumerate}
\item the pair $(\Ga,F^\pm)$ is of divergence type,
\item the Patterson densities $(\mu^\pm_x)_{x\in X}$ are unique up to
  a multiplicative constant and the Gibbs measure $m_{F}$ is uniquely
  defined up to a multiplicative constant.
\item the Gibbs measure $m_{F}$ gives full measure to the image
  $\gls{recurrentset}$ of
$$
\wt\Omega_{\rm c}=\{\ell\in\G X\;:\;\ell_\pm\in\Lambda_{\rm c}\Ga\}
$$  
  in $\Ga\bs\G X$, and
\item the geodesic flow is ergodic for $m_F$. \cqfd 
\end{enumerate}
\ecoro

\subsection*{On the finiteness of Gibbs measures}
\addcontentsline{toc}{subsection}{On the finiteness of Gibbs measures}

As the finiteness of the Gibbs measures will be a standing
hypothesis in many of the following results, we now give criteria for
Gibbs measures to be finite.  Recall\footnote{See Section
\ref{subsec:catmoinsun}.} that the discrete nonelementary group of
isometries $\Ga$ of $X$ is geometrically finite if every element of
$\Lambda\Ga$ is either a conical limit point or a bounded parabolic
limit point of $\Ga$.

\btheo \label{theo:DOPB}%
Assume that  $\Ga$ is a geometrically finite discrete group of
isometries of $X$.
\begin{enumerate}
\item If $(\Ga,F^\pm)$ is of divergence type, then the Gibbs measure
  $m_{F}$ is finite if and only if for every bounded parabolic limit
  point $p$ of $\Ga$, the series
$$
\sum_{\alpha\in\Ga_p} \; 
d(x,\alpha y) \;e^{\int_x^{\alpha y} (\wt F^\pm-\delta)}
$$ 
converges, where $\Ga_p$ is the stabiliser of $p$ in $\Ga$.  
\item If we have $\delta_{\Ga_p,\,F_p^\pm}< \delta$, for every bounded
  parabolic limit point $p$ of $\Ga$ with stabiliser $\Ga_p$ in $\Ga$
  and with $F_p^\pm:\Ga_p\bs X\ra\RR$ the map induced by $\wt F^\pm$,
  then $(\Ga,F)$ is of divergence type.  In particular, $m_{F}$ is
  finite.
\end{enumerate}
\etheo

When $X$ is a manifold, this result is due to
\cite[Theo.~B]{DalOtaPei00} for the case $F=0$, and to
\cite{Coudene03} and \cite[Theo.~8.3, 8.4]{PauPolSha15} for the
general case of H\"older-continuous potentials. When $F=0$ but on much
more general assumptions on $X$ with optimal generality, this result
is due to \cite[Theo.~17.1.2]{DasSimUrb14arxiv}.

\medskip \dem The proof is similar to the manifold case in
\cite{PauPolSha15}, which follows closely the proof of
\cite{DalOtaPei00}.  Note that the convergence or divergence of the
above series does not depend on the choice of the sign $\pm$.

Let $\Par_\Ga$ be the set of bounded parabolic limit points of $\Ga$.
By \cite[Lem.~1.9]{Roblin03}\footnote{See also \cite{Paulin04b} for
  the case of simplicial trees and
  \cite[Theo.~12.4.5]{DasSimUrb14arxiv} for a greater generality on
  $X$.}, there exists a $\Ga$-equivariant family
$(\H_p)_{p\in\Par_\Ga}$ of pairwise disjoint closed horoballs, with
$\H_p$ centred at $p$, such that the quotient
$$
M_0=\Ga\bs\big(\C\Lambda\Ga-\bigcup_{p\in\Par_\Ga} \H_p\big)
$$ 
is compact. Using Proposition \ref{prop:continuGibbscocycle},
Theorem \ref{theo:HTSR} and Equation \eqref{eq:HCN} instead of
respectively \cite[Lem.~3.4, Coro.~5.10,Lem.~3.2]{PauPolSha15}, the
\ref{eq:HC}-property, the proofs of \cite[Theo.~8.3,
  8.4]{PauPolSha15} then extend to our situation.  
\cqfd

\medskip Recall that the {\em length spectrum}\index{length!spectrum}
of $\Ga$ on $X$ is the additive subgroup of $\RR$ generated by the
translation lengths in $X$ of the elements of $\Ga$.

Recall that a continuous-time one-parameter group $(h^t)_{t\in\RR}$ of
homeomorphisms of a topological space $Z$ is {\it topologically
  mixing}\index{topologically!mixing} if for all nonempty open subsets
$U,V$ of $Z$, there exists $t_0\in\RR$ such that for every $t\ge t_0$,
we have $U\cap h^t V\neq \emptyset$.

We have the following result, due to \cite[Theo.~1]{Babillot02b} in the
manifold case, with developments by \cite{Roblin03} when $\wt F=0$,
and by \cite[Theo.~8.1]{PauPolSha15} for manifolds with pinched negative
curvature.

\btheo\label{theo:mixing} If the Gibbs measure is finite, then the
following assertions are equivalent :
\begin{enumerate}
\item
the geodesic flow of $\Ga\backslash X$ is mixing for the Gibbs measure,
\item
the geodesic flow of $\Ga\backslash X$ is topologically mixing on its
nonwandering set, which is the quotient under $\Ga$ of the space of
geodesic lines in $X$ both of whose endpoints belong to $\Lambda\Ga$.
\item
the length spectrum of $\Ga$ on $X$ is not contained in a discrete
subgroup of $\RR$. \cqfd
\end{enumerate}
\etheo

\medskip In the manifold case, the third assertion of Theorem
\ref{theo:mixing} is satisfied, for example, if $\Ga$ has a parabolic
element, if $\Lambda\Ga$ is not totally disconnected (hence if
$\Ga\backslash X$ is compact), or if $X$ is a surface or a (rank-one)
symmetric space, see for instance \cite{Dalbo99,Dalbo00}.  

Error terms for the mixing property will be described in Chapter
\ref{sec:mixingrate}. The above result holds for the continuous time
geodesic flow when $X$ is a metric tree. See Theorem
\ref{theo:uniflatmBMmixing} for a version of this theorem for the
discrete time geodesic flow on simplicial trees.  At least when $X$ is
an $\RR$-tree and $\Ga$ is a uniform lattice (so that $\Ga\bs X$ is a
finite metric graph), we have a stronger result under additional
regularity assumptions, see Section \ref{subsec:mixingratesimpgraphs}.

\subsection*{Bowen-Margulis measure computations in locally symmetric 
spaces}\label{subsec:BMcomputlocsym}
\addcontentsline{toc}{subsection}{Bowen-Margulis measure computations in 
locally symmetric spaces}

Assume in this subsection that the potential $\wt F$ is zero. The next
result, Proposition \ref{prop:computexplimBMhomogen}, gathers
computations done in \cite{ParPau16LMS,ParPau16MA} of the
Bowen-Margulis measures $m_{\rm BM}$ when $X$ is a real or complex
hyperbolic space, and  $\Ga$ is a lattice. We start by giving the
notation necessary in order to state Proposition
\ref{prop:computexplimBMhomogen}.

When $X$ is a complete simply connected Riemannian manifold with
dimension $m\geq 2$ and sectional curvature at most $-1$, we endow the
usual unit tangent bundle $T^1X$ with Sasaki's Riemannian
metric.\footnote{See Section \ref{subsec:holdercont}.}  Its Riemannian
measure $\dvol_{T^1X}$, called {\em Liouville's
  measure},\index{Liouville's measure}\index{measure!Liouville}
disintegrates (equivariantly with respect to $\Ga$) under the
fibration $\pi:T^1X\ra X$ over the Riemannian measure $\dvol_{X}$ of
$X$, as
$$
\dvol_{T^1X}=\int_{x\in X} \;\dvol_{T^1_xX}\;\dvol_{X}(x)\;,
$$
where $\dvol_{T^1_xX}$ is the spherical measure on the fiber $T^1_xX$
of $\pi$ above $x\in X$. In particular,
$$ 
\Vol(T^1(\Ga\bs X))=\Vol(\SSS^{m-1})\Vol(\Ga\bs X)\;.
$$ 

Assume furthermore that $X$ is a symmetric space, and that $\Ga$ is
a lattice in $\Isom(X)$. In particular, $\Ga$ is geometrically finite
and the critical exponent of the stabiliser in $\Ga$ of every bounded
parabolic fixed point of $\Ga$ is strictly less than the critical
exponent of $\Ga$, see for instance \cite{Dalbo99,Dalbo00}. Then the
Patterson density is independent of $\Ga$ and uniquely defined up to a
multiplicative constant by Theorem \ref{theo:DOPB} (2) and Lemma
\ref{coro:finitudeGibbsdivuniq} (2). We take $\mu^-_x=\mu^+_x$ for
every $x\in\hnc$ and we will denote this measure simply by $\mu_x$.
By homogeneity reasons, the Bowen-Margulis measure of $\Ga$ on $X$ is
proportional to the Liouville measure $\vol_{T^1X}$, and the main
point of Proposition \ref{prop:computexplimBMhomogen} is to compute
the proportionality constant.

\medskip
Let $n\geq 2$. We endow $\RR^{n-1}$ with its usual Euclidean norm
$\|\cdot\|$ and its usual Lebesgue measure $\lambda_{n-1}$. We use the
upper halfspace model for the real hyperbolic space $\hnr$ of
dimension $n$, that is, $\hnr=\{(x_1,\cdots, x_n)\in\RR^n\;:\;x_n>0\}$
endowed with the Riemannian metric
$$ 
ds^2_{\,\hnr}= \frac{1}{x_n^2}\;(dx_1^2+\cdots+dx_n^2)\;.
$$ 
We identify $\RR^{n-1}$ with $\RR^{n-1}\times\{0\}$ in $\RR^n$, and
again denote by $\|\cdot\|$ the usual Euclidean norm on $\RR^n$.  The
boundary at infinity of $\hnr$ is $\partial_\infty \hnr=\RR^{n-1}\cup
\{\infty\}$. Assuming that $\Ga$ is a lattice in $\Isom(\hnr)$, we
normalise its Patterson density so that in the ball model of $\hnr$
with center $0$, the measure $\mu_0$ is the spherical measure on the
space at infinity $\SSS^{n-1}$.

\medskip
Let $n\geq 2$. We refer to \cite{Goldman99} and \cite[\S
  3]{ParPau16MA} for background on the complex hyperbolic $n$-space
$\hnc$. We denote by $\zeta\cdot\overline{\zeta'}= \sum_{i=1}^{n-1}
\zeta_i\, \overline{\zeta'_i}$ the standard Hermitian product and by
$d\zeta$ the standard Lebesgue measure on $\CC^{n-1}$. We denote by
$\gls{heisenberggroup}$ the {\it Heisenberg group}\index{Heisenberg group} of
dimension $2n-1$, which is the real Lie group structure on
$\CC^{n-1}\times \RR$ with law
$$
(\zeta,u)(\zeta',u')=
(\zeta+\zeta',u+u'+2\,\Im\;\zeta\cdot\overline{\zeta'})\;.
$$ 
We endow $\Heis_{2n-1}$ with the usual left-invariant Haar measure
$d\lambda_{2n-1}(\zeta,u) = d\zeta \,du$ and with the {\it Cygan
  distance}\index{distance!Cygan }\index{Cygan distance}\footnote{See
  \cite[page 160]{Goldman99}. It is called the {\it Kor\'anyi
    distance} by many people working in sub-Riemannian geometry,
  though Kor\'anyi \cite{Koranyi85} does attribute it to Cygan
  \cite{Cygan78}.} $d_{\rm Cyg}$ which is the unique left-invariant
distance on $\Heis_{2n-1}$ with $d_{ \rm Cyg} ((\zeta,u),(0,0))=
(|\zeta|^4 +u^2)^{\frac{1}{4}}$.  

We use the {\it Siegel domain}\index{Siegel domain} model of the
complex hyperbolic space $\gls{complexhyperbolicspace}$ of dimension
$n$, normalised to have maximum constant sectional curvature $-1$,
hence to be $\CAT(-1)$.  This is the manifold $\Heis_{2n-1}\times\,
]0,+\infty[$ endowed with the Riemannian metric given, in the {\it
 horospherical coordinates}\index{horospherical coordinates}
$(\zeta,u,t)\in \CC^{n-1} \times \RR \times\,]0,+\infty[$, by
$$
ds^2_{\,\hnc}= \frac{1}{4\,t^2}\big(dt^2
+(du+2\,\Im\;d\zeta\cdot\overline{\zeta}\,)^2+
4\,t\,d\zeta\cdot\overline{d\zeta} \,\big)\,,
$$
so that its  volume form is
$$
d\operatorname{vol}_{\hnc}(\zeta,u,t)=
\frac{1}{4\,t^{n+1}}\;d\zeta\,du\,dt\;.
$$
Note that the action of $\Heis_{2n-1}$ on $\hnc= \Heis_{2n-1}
\times\, ]0,+\infty[$ by left translations on the first factor,
preserving the second one, is isometric. We identify $\Heis_{2n-1}$
with $\Heis_{2n-1}\times\{0\}$ and we endow $\Heis_{2n-1}\times
[0,+\infty[$ with the distance $d_{\rm Cyg}$ extending the Cygan
distance on $\Heis_{2n-1}$, defined by
$$
d_{\rm Cyg}((\zeta,u,t),(\zeta',u',t'))=
\Big|\,|\zeta-\zeta'|^2+|t-t'|+
i(u-u'+2\,\Im\;\zeta\cdot\overline{\zeta'})\Big|^{1/2}\,.
$$
The space at infinity $\partial_\infty \hnc$ of $\hnc$ is the
Alexandrov compactification $\Heis_{2n-1} \cup \{\infty\}$, so that
the extension at infinity of the isometric action of $\Heis_{2n-1}$ on
$\hnc$ fixes $\infty$ and is the left translation on $\Heis_{2n-1}$.
We denote by $\H_\infty=\Heis_{2n-1} \times [1,+\infty[$ the horoball
of $\hnc$ centered at $\infty$ whose boundary contains the point
$(0,0,1)$.

Assuming that $\Ga$ is a lattice in $\Isom(\hnc)$, we normalise its
Patterson density $(\mu_x)_{x\in\hnc}$ as follows. If
$\mu_{\H_\infty}$ is the measure defined in Proposition
\ref{prop:HPun} associated with the horoball $\H_\infty$, then
$$
d\mu_{\H_\infty}(\zeta,u)= d\lambda_{2n-1}(\zeta,u) =d\zeta\,du\;.
$$ 
This is possible since $\mu_{\H_\infty}$ is invariant under the
isometric action of $\Heis_{2n-1}$ on $\hnc$, which preserves
$\H_\infty$, hence is a left-invariant Haar measure on
$\partial_\infty \hnc-\{\infty\} = \Heis_{2n-1}$.

As the arguments of the following result are purely computational and
rather long, we do not copy them in this book, but we refer respectively
to the proofs of \cite[Eq.~(19), Eq.~(21), Prop.~10]{ParPau16LMS} and
\cite[Lem.~12 (i), (ii), (iii)]{ParPau16MA}.  Analogous computations can
be done when $X$ is the quaternionic hyperbolic $n$-space $\HH_{\HH}^n$.

\bprop\label{prop:computexplimBMhomogen} (1) Let $\Ga$ be a lattice in
$\Isom(\hnr)$, with Patterson density $(\mu_x)_{x\in\hnr}$ normalised
as above.  For all $x=(x_0,\cdots,x_n)$ in $\hnr$, $\xi\in
\partial_\infty \hnr- \{\infty\}$ and $v\in T^1\hnr$ such that
$v_\pm\neq\infty$, we have

\smallskip\noindent
(i)~ 
$d\mu_x(\xi)= \big(\frac {2 x_n}{\|x-\xi\|^2}\big)^{n-1}\;
d\lambda_{n-1}(\xi)$,

\medskip\noindent
(ii)~ using a Hopf parametrisation $v\mapsto (v_-,v_+,s)$,
$$
d\wt m_{\rm BM}(v)=\frac{2^{2(n-1)}d\lambda_{n-1}(v_{-})\;
  d\lambda_{n-1}(v_{+})\;dt}{\|v_+-v_-\|^{2(n-1)}}\;,
$$

\noindent(iii)
$$
\wt m_{\rm BM}=2^{n-1}\,\Vol_{T^1\hnr}\;,
$$ 
and in particular,
$$
\|m_{\rm BM}\|=2^{n-1}\Vol(\SSS^{n-1})\Vol(\Ga\bs\hnr)\;.
$$

\medskip
(2) Let $\Ga$ be a lattice in $\Isom(\hnc)$, with Patterson density
$(\mu_x)_{x\in\hnc}$ normalised as above.  For all $x=(\zeta,u,t)$ in
$\hnc$, $(\xi,r)\in \partial_\infty \hnc-\{\infty\}$ and $v\in
T^1\hnc$ such that $v_\pm\neq\infty$, we have

\smallskip\noindent
(i)~ $\displaystyle d\mu_x(\xi,r)=
\frac{t^{n}}{d_{\rm Cyg}(x,(\xi,r))^{4n}}\;d\xi\,dr\;;$

\medskip\noindent
(ii) using a Hopf parametrisation $v\mapsto (v_-,v_+,s)$,
$$
d\wt m_{\rm BM}(v)=
\frac{d\lambda_{2n-1}(v_-)\,d\lambda_{2n-1}(v_+)\,ds} 
{d_{\rm  Cyg}(v_-,v_+)^{4n}}\,;
$$

\noindent(iii)
$$
\wt m_{\rm BM}= \frac{1}{2^{2n-2}}\; \vol_{T^1\hnc}\,,
$$
and in particular 
$$
\|m_{\rm BM}\|= \frac{\pi^{n}}{2^{2n-3}\,(n-1)!}\;\Vol(\Ga\bs\hnc)\;.\;\;\;\Box
$$
\eprop

\subsection*{On the cohomological invariance of Gibbs measures}
\addcontentsline{toc}{subsection}{On the cohomological invariance 
of Gibbs measures}

We end this Section by an elementary remark on the independence of
Gibbs measures upon replacement of the potential $F$ by a cohomologous
one.

\brema\label{rem:cohomologuustuff}{\rm
Let $\wt F^*:T^1X\ra\RR$ be a potential for $\Ga$ cohomologous to $\wt
F$ and satisfying the \ref{eq:HC}-property.  As usual, let $\wt
F^{*+}=\wt F^*$ and $\wt F^{*-}=\wt F^*\circ \iota$, and let
$F^*:\Ga\bs T^1X\ra\RR$ be the induced map.  Let $\wt G :T^1X\ra \RR$
be a continuous $\Ga$-invariant function such that, for every $\ell\in
\G X$, the map $t\mapsto \wt G(v_{\flow{t}\ell})$ is differentiable
and $\wt F^*(v_\ell)-\wt F(v_\ell)=\frac{d}{dt}_{\mid t=0}\wt
G(v_{\flow{t}\ell})$.  Furthermore assume that $X$ is an $\RR$-tree or
that $G$ is uniformly continuous (for instance Hölder-continuous).

For all $x\in X$ and $\xi\in\partial_\infty X$, let $\ell_{x,\,\xi}$
be any geodesic line with footpoint $\ell_{x,\,\xi}(0)=x$ and positive
endpoint $(\ell_{x,\,\xi})_+=\xi$, and let $\ell_{\xi,\,x}$ be any
geodesic line with $\ell_{\xi,\,x}(0)=x$ and origin
$(\ell_{\xi,\,x})_-=\xi$.  Note that the value $\wt
G(v_{\ell_{x,\,\xi}})$ is independent of the choice of
$\ell_{x,\,\xi}$, by the continuity of $\wt G$, and similarly for $\wt
G(v_{\ell_{\xi,\,x}})$. In particular, for all $\ga\in\Ga$, by the
$\Ga$-invariance of $\wt G$, we have
\begin{equation}\label{eq:equivwtG}
\wt G(v_{\ell_{x,\,\ga^{-1}\xi}})=\wt G(v_{\ell_{\ga x,\,\xi}})
\;\;\;{\rm and}\;\;\;
\wt G\circ\iota \,(v_{\ell_{x,\,\xi}})=\wt G (v_{\ell_{\xi,\,x}})\;.
\end{equation}

Note that $\wt F^{*-}=\wt F^*\circ \iota$ and $\wt F^-=\wt F\circ
\iota$ are cohomologous, since if $\wt G^*=-\,\wt G\circ\iota$, for
every $\ell\in\G X$, we have
\begin{align*}
\wt F^*\circ \iota(v_\ell)-\wt F\circ \iota(v_\ell)&=
\wt F^*(v_{\iota\ell})-\wt F(v_{\iota\ell})=
\frac{d}{dt}_{\mid t=0}\wt G(v_{\flow{t}\iota\ell})\\ & =
\frac{d}{dt}_{\mid t=0}\wt G(\iota v_{\flow{-t}\ell})=
\frac{d}{dt}_{\mid t=0}\wt G^*(v_{\flow{t}\ell})\;.
\end{align*}

As already seen in Lemma \ref{lem:proprielemcritexpo} (1) and (2), the
critical exponent $\delta_{\Ga,\,F^{*\pm}}$ is equal to the critical
exponent $\delta_{\Ga,\,F^\pm}$, and independent of the choice of
$\pm$, and we denote by $\delta$ the common value in the definition of
the Gibbs cocycle.

Let us prove that if $C^{*\pm}=C^\pm_{\Ga,\,F^{*\pm}}$ is the Gibbs
cocycle associated with $(\Ga,F^{*\pm})$, then $C^{*\pm}$ and $C^\pm$
are {\em cohomologous}\index{cocycle!cohomologous}\index{cohomologous}:
\begin{equation}\label{eq:cohomologuecocycleplus}
C^{*+}_{\xi}(x,y)-C^+_{\xi}(x,y)=\wt G(v_{\ell_{x,\,\xi}})-\wt G(v_{\ell_{y,\,\xi}})\;,
\end{equation}
and similarly
\begin{equation}\label{eq:cohomologuecocyclemoins}
C^{*-}_{\xi}(x,y)-C^-_{\xi}(x,y)=
\wt G^*(v_{\ell_{x,\,\xi}})-\wt G^*(v_{\ell_{y,\,\xi}})\;.
\end{equation}
We only prove the first equality, the second one follows similarly,
noting that $\wt G^*$ is uniformly continuous if $\wt G$ is, since
$\iota$ is isometric.  For all $x,y$ in $X$ and $\xi\in\partial_\infty
X$, let $t\mapsto \xi_t$ be a geodesic ray with point at infinity
$\xi$, let $a_t=d(x,\xi_t)$, let $b_t=d(y,\xi_t)$, and for $z=x,y$,
let $\ell_{z,\,\xi_t}$ be any geodesic line with footpoint $z$ passing
through $\xi_t$.  Then
\begin{align*}
&\Big(\int_y^{\xi_t}(\wt F^{*+}-\delta)-\int_x^{\xi_t}(\wt F^{*+}-\delta)\Big)
-\Big(\int_y^{\xi_t}(\wt F^+-\delta)-\int_x^{\xi_t}(\wt F^+-\delta)\Big) \\=\;&
  \int_y^{\xi_t}(\wt F^{*+}-\wt F^+)-\int_x^{\xi_t}(\wt F^{*+}-\wt F^+)\\=\;&
\int_0^{b_t}\frac{d}{dt}\wt G(v_{\flow{s}\ell_{y,\,\xi_t}})\;ds -
  \int_0^{a_t}\frac{d}{dt}\wt G(v_{\flow{s}\ell_{x,\,\xi_t}})\;ds\\=\;&
  \wt G(v_{\ell_{x,\,\xi_t}})- \wt G(v_{\ell_{y,\,\xi_t}}) +
  \big(\wt G(v_{\flow{b_t}\ell_{y,\,\xi_t}})- \wt G(v_{\flow{a_t}\ell_{x,\,\xi_t}})\big)\;.
\end{align*}
When $t$ goes to $+\infty$, the first term of this series of
equalities converges to $C^{*+}_{\xi}(x,y)-C^+_{\xi}(x,y)$ by the
definition of the Gibbs cocycle (see Section
\ref{subsec:Gibbscocycle}). By continuity, $\wt
G(v_{\ell_{y,\,\xi_t}})$ and $\wt G(v_{\ell_{x,\,\xi_t}})$ converge to
$\wt G(v_{\ell_{y,\,\xi}})$ and $\wt G(v_{\ell_{x,\,\xi}})$
respectively. If $X$ is an $\RR$-tree, then if $t$ is large enough, we
have $v_{\flow{b_t}\ell_{y,\,\xi_t}} =v_{\flow{a_t}\ell_{x,\,\xi_t}}$,
hence Equation \eqref{eq:cohomologuecocycleplus} follows. Otherwise,
by the uniform continuity of $\wt G$, since
$v_{\flow{b_t}\ell_{y,\,\xi_t}}$ and $v_{\flow{a_t}\ell_{x,\,\xi_t}}$
are uniformly arbitrarily close as $t$ tends to $0$ by the $\CAT(-1)$
property, the result also follows.

\medskip
Let $(\mu^\pm_{x})_{x\in X}$ be a Patterson density for $(\Ga,
F^\pm)$.  In order to simplify the notation, let $\wt G^+=\wt G$ and
$\wt G^-=\wt G^*$.  The family of measures $(\mu^{*\pm}_{x})_{x\in X}$
defined by setting, for all $x\in X$ and $\xi\in\partial_\infty X$,
\begin{equation}\label{eq:defiPattdenscohomol}
d\mu^{*\pm}_x(\xi)=e^{-\wt G^\pm(v_{\ell_{x,\,\xi}})}\;d\mu^{\pm}_x(\xi)\;,
\end{equation} 
is also a Patterson density for $(\Ga,F^{*\pm})$. Indeed, the
equivariance property \eqref{eq:equivarPatdens} for
$(\mu^{*\pm}_{x})_{x\in X}$ follows from the one for
$(\mu^{\pm}_{x})_{x\in X}$ and from Equation \eqref{eq:equivwtG}. The
absolutely continuous property \eqref{eq:quasinivarPatdens} for
$(\mu^{*\pm}_{x})_{x\in X}$ follows from the one for
$(\mu^{\pm}_{x})_{x\in X}$ and Equations
\eqref{eq:cohomologuecocycleplus} and
\eqref{eq:cohomologuecocyclemoins}.

Assume that the Patterson density for $(\Ga,F^{*\pm})$ defined by
Equation \eqref{eq:defiPattdenscohomol} is chosen in order to
construct the Gibbs measure $\wt m_{F^*}$ for $(\Ga,F^*)$ on $\G X$.
Then using 

\smallskip
$\bullet$~ Hopf's parametrisation with respect to the base point $x_0$
and Equation \eqref{eq:defigibbs} with $F$ replaced by $F^*$ for the
first equality,

\smallskip
$\bullet$~ Equations \eqref{eq:cohomologuecocycleplus},
\eqref{eq:cohomologuecocyclemoins}, \eqref{eq:defiPattdenscohomol} and
cancellations for the second equality,

\smallskip
$\bullet$~ the definition of $\wt G^*=-\,\wt G\circ\iota$ and again
Equation \eqref{eq:defigibbs} for the third equality,

\smallskip
$\bullet$~ Equation \eqref{eq:equivwtG} and the fact that we may
choose $\ell_{\ell_-,\,\ell(0)}=\ell$ and $\ell_{\ell(0),\,\ell_+}=
\ell$ for the last equality,

\smallskip
\noindent
we have 
\begin{align*}
d\wt m_{F^*}(\ell)&=
e^{C_{\ell_-}^{*-}(x_0,\,\ell(0))\,+\,C^{*+}_{\ell_+}(x_0,\,\ell(0))}\;
d\mu_{x_0}^{*-}(\ell_-)\,d\mu^{*+}_{x_0}(\ell_+)\,dt\\ 
& =
e^{C_{\ell_-}^{-}(x_0,\,\ell(0))\,-\wt G^*(v_{\ell_{\ell(0),\,\ell_-}})\,-
\,C^{+}_{\ell_+}(x_0,\,\ell(0))\,+\wt G(v_{\ell_{\ell(0),\,\ell_+}})}\;
d\mu^-_{x_0}(\ell_-)\,d\mu^+_{x_0}(\ell_+)\,dt\\ 
& =
e^{\wt G\circ\iota (v_{\ell_{\ell(0),\,\ell_-}})\,-\wt G(v_{\ell_{\ell(0),\,\ell_-}})}\;
d\wt m_{F}(\ell)\\ 
& = 
e^{\wt G(v_\ell)-\wt G(v_\ell)}\;d\wt m_{F}(\ell)\;,
\end{align*}
hence $\wt m_{F^*}=\wt m_{F}$.

In particular, since the Gibbs measure, when finite, is independent up
to a multiplicative constant on the choice of the Patterson densities
by Corollary \ref{coro:finitudeGibbsdivuniq}, we have that $m_F$ is finite
if and only if $m_{F^*}$ is finite, and then
\begin{equation}\label{eq:equalgibbsmeascohomolog}
\frac{m_{F^*}}{\|m_{F^*}\|} = \frac{m_F}{\|m_F\|}\;.
\end{equation} 

}\erema

\section{Patterson densities for simplicial trees}
\label{subsec:pattersongibbstrees}

In this Section and the following one, we specialise and modify the
general framework of the previous sections to treat simplicial trees.
Recall\footnote{See Section \ref{subsec:trees}.} that a simplicial
tree $\XX$ is a metric tree whose edge length map is constant equal to
$1$.  The time $1$ map of the geodesic flow $(\flow{t})_{t\in\RR}$ on
the space $\gengeod X$ of all generalised geodesic lines of the
geometric realisation $X=|\XX|_1$ of $\XX$ preserves for instance its
subset of generalised geodesic lines whose footpoints are at distance
at most $1/4$ from vertices. Since both this subset and its complement
have nonempty interior in $\gengeod X$, the geodesic flow on $\gengeod
X$ has no mixing or ergodic measure with full support.  This is why we
considered the discrete  time geodesic flow $(\flow{t})_{t\in\ZZ}$ on
$\gengeod \XX$ in Section \ref{subsec:trees}.

\medskip
Let $\XX$ be a locally finite simplicial tree without terminal
vertices, and let $X=|\XX|_1$ be its geometric realisation.  Let $\Ga$
be a nonelementary discrete subgroup of $\Aut(\XX)$. Let $\wt
F:T^1X\ra \RR$ be a potential for $\Ga$, and let $\wt F^+=\wt F$, $\wt
F^-= \wt F\circ\iota$. Let $\delta=\delta_{\Ga,F^\pm}$ be the critical
exponent of $(\Ga,F^\pm)$. Let $C^\pm:\partial_\infty X\times X\times
X\ra \RR$ be the associated (normalised) Gibbs cocycles.  Let
$(\mu^\pm_x)_{x\in X}$ be two Patterson densities on $\partial_\infty
X$ for the pairs $(\Ga,F^\pm)$.

Note that only the restrictions of the cocycles $C^\pm$ to
$\partial_\infty X\times V\XX\times V\XX$ are useful and that it is
often convenient and always sufficient to replace the cocycles by
finite sums involving a system of conductances (as defined in Section
\ref{subsec:cond}), see below.  Furthermore, only the restriction
$(\mu^\pm_x)_{x\in V\XX}$ of the family of Patterson densities to the
set of vertices of $\XX$ is useful.

\bexem Let $\XX$ be a simplicial tree with geometric realisation $X$
and let $\wt c:E\XX\ra \RR$ be a system of conductances on $\XX$.  For
all $x,y$ in $V\XX$ and $\xi\in \partial_\infty X$, with the usual
convention on the empty sums, let
$$
c^+_\xi(x,y)=\sum_{i=1}^m \;\wt c(e_i)-\sum_{j=1}^n \;\wt c(f_j)
$$ 
and 
$$
c^-_\xi(x,y)=\sum_{i=1}^m \;\wt c(\overline{e_i})-
\sum_{j=1}^n \;\wt c(\overline{f_j})\,,
$$ 
where, if $p\in V\XX$ is such that $[p,\xi[\;=[x,\xi[\;\cap\,
[y,\xi[\,$, then $(e_1,e_2,\dots, e_m)$ is the geodesic edge
path in $\XX$ from $x=o(e_1)$ to $p=t(e_m)$ and $(f_1,f_2,\dots, f_n)$
is the geodesic edge path in $\XX$ from $v=o(f_1)$ to $p=t(f_n)$.

\begin{center}
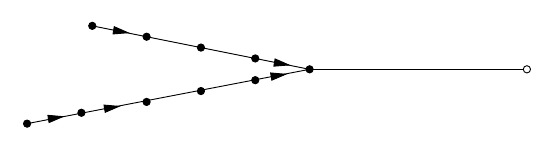
\end{center}

With $\delta_c$ defined in the end of Section \ref{subsec:cond} and
with $C^\pm$ the Gibbs cocycles for $(\Ga,\wt F_c)$, by Equation
\eqref{eq:cocycletreecase} and by Proposition
\ref{prop:integpotconduct}, we have, for all $\xi\in
\partial_\infty X$ and $x,y\in V\XX$,
$$
C^\pm_\xi(x,y)=-\,c^\pm_\xi(x,y)+\delta_c\,\beta_\xi(x,y)\,,
$$
and Equation \eqref{eq:quasinivarPatdens} gives
$$
d\mu^\pm_x(\xi)=e^{c^\pm_\xi(x,\,y)-\delta_c\,\beta_\xi(x,\,y)}\;d\mu^\pm_y(\xi)\,.
$$
\eexem

\medskip
Using the particular structure of trees, we can prove a version of the
Shadow Lemma \ref{lem:shadowlemma} where one can take the radius $R$
to be $0$. When $F=0$, this result is due to Coornaert
\cite{Coornaert93}.

\blemm[Mohsen's shadow lemma for trees]\label{lem:shadowlemmatrees}%
\index{shadow!lemma!for trees} Let $K$ be a finite subset of $V\XX$.
There exists $C>0$ such that for all $\ga\in\Ga$ and $x,y\in K$ with
$y\in\C\Lambda\Ga$, we have
$$
\frac{1}{C}\;e^{\int_x^{\ga y} (\wt F^\pm-\delta)}\leq 
\mu^\pm_{x}\big(\OOO_x\{\ga y\}\big)
\leq C\;e^{\int_x^{\ga y} (\wt F^\pm-\delta)}\;.
$$
\elemm

\dem The structure of the proof is the same one as for Lemma
\ref{lem:shadowlemma} (1) (that is, the one of
\cite[Lem.~3.10]{PauPolSha15}) with differences towards the end of the
argument.  Let us prove that there exists $C=C_K>0$ such that for all
$\ga\in\Ga$ and $x,y\in K$ with $\ga y\in \C\Lambda\Ga$, we have
\begin{equation}\label{eq:ombrerreduc}
\frac{1}{C}\leq
\mu^\pm_{\ga y}\big(\OOO_x\{\ga y\}\big)
\leq C\;.
\end{equation}
Assuming this, let us prove Lemma \ref{lem:shadowlemmatrees}.  By
Equation \eqref{eq:quasinivarPatdens}, we have
$$
\mu^\pm_{x}\big(\OOO_x\{\ga y\}\big)=\int_{\xi\in\OOO_x\{\ga y\}}
e^{-C_{\xi}^\pm(x,\,\ga y)}\;d\mu_{\ga y}^\pm(\xi)\;.
$$
Note that $C^\pm_\xi(x,\ga y)+\int_x^{\ga y} (\wt
F^\pm-\delta)=0$ if $\xi\in \OOO_x\{\ga y\}$ (that is, if $\ga y\in
[x,\xi[\,$), by Equation \eqref{eq:changemoinsplus}.
Hence
$$
 \mu_{x}^\pm(\OOO_x\{\ga y\}) =
e^{\int_x^{\ga y} (\wt F^\pm-\delta)}\;\mu_{\ga y}^\pm(\OOO_x\{\ga y\})\;,
$$
and Lemma \ref{lem:shadowlemmatrees} follows from Equation
\eqref{eq:ombrerreduc}.

\medskip
Let us now prove the upper bound in Equation \eqref{eq:ombrerreduc}.
Fix $z_0\in K$, and let
$$
C'=\sup_{x,\,y\in K,\;\xi\in \partial_\infty X}\;|C_{\xi}^\pm(x,y)|\;,
$$ 
which is finite by Proposition \ref{prop:continuGibbscocycle} (2),
since $K$ is compact and $\wt F^\pm$ continuous. Then, using Equation
\eqref{eq:equivarPatdens} for the equality and Equation
\eqref{eq:quasinivarPatdens} for the last inequality, we have
$$
\mu_{\ga y}^\pm(\OOO_x\{\ga y\})
\leq \|\mu_{\ga y}^\pm\|=\|\mu_{y}^\pm\| \leq e^{C'}\;\|\mu_{z_0}^\pm\|\;,
$$ 
and the upper bound holds if $C\geq e^{C'}\;\|\mu_{z_0}^\pm\|$.

\medskip 
Finally, in order to prove the lower bound in Equation
\eqref{eq:ombrerreduc}, we assume for a contradiction that there exist
sequences $(x_i)_{i\in\NN}$, $(y_i)_{i\in\NN}$ in $K$ with
$y_i\in\C\Lambda\Ga$ and $(\ga_i)_{i\in\NN}$ in $\Ga$ such that
$\mu^\pm_{\ga_i y_i}( \OOO_{x_i}\{\ga_i y_i\})$ converges to $0$ as
$i\to+\infty$.  Up to extracting a subsequence, since $K$ is finite,
we may assume that the sequences $(x_i)_{i\in\NN}$ and
$(y_i)_{i\in\NN}$ are constant, say with value $x$ and $y$
respectively. Since $y\in\C\Lambda\Ga$, as every point in
$\C\Lambda\Ga$ belongs to at least one geodesic line between two limit
points of $\Ga$, the geodesic segment from $x$ to $\ga_i y$ may be
extended to a geodesic ray from $x$ to a limit point of $\Ga$. Since
the support of $\mu^\pm_{z}$ is equal to $\Lambda\Ga$ for any $z\in
X$, we have $\mu^\pm_{\ga_i y}(\OOO_{x}\{\ga_i y\})>0$ for all
$i\in\NN$. Thus, up to taking a subsequence, we can assume that
$\ga_i^{-1} x$ converges to $\xi\in\Lambda\Ga$ (otherwise by
discreteness, we may extract a subsequence so that $(\ga_i)_{i\in\NN}$
is constant, and $\mu^\pm_{\ga_i y}( \OOO_{x}\{\ga_i y\})$ cannot
converge to $0$).

Since $\XX$ is a tree, there exists a positive integer $N$ such that
$\OOO_{\ga^{-1}_ix}\{y\}=\OOO_{\ga^{-1}_Nx}\{y\}=\OOO_{\xi}\{y\}$ for
all $i\ge N$. As above, $\OOO_{\xi}\{y\}$ meets $\Lambda\Ga$ since
$y\in\C\Lambda\Ga$, thus $\mu^\pm_y(\OOO_{\xi}\{y\})>0$. But for every
$i\geq N$,
$$
\mu^\pm_y\big(\OOO_{\xi}\{y\}\big)=
(\ga_i^{-1})_*\mu^\pm_{\ga_i y}\big(\OOO_{\ga_i^{-1} x}\{y\}\big)
=\mu^\pm_{\ga_i y}\big(\OOO_x\{\ga_i y\}\big)
$$ 
tends to $0$ as $i\ra+\infty$, a contradiction.  
\cqfd

\medskip
Let $\gls{totalmassPatterson}:V\XX\ra \mathopen{[}0, +\infty
  \mathclose{[}$ be the {\em total mass functions of the Patterson
      densities\,}:
$$
\wt \phi_{\mu^\pm}(x)=\|\mu_x^\pm\|
$$ 
for every $x\in V\XX$. These maps are $\Ga$-invariant by Equation
\eqref{eq:equivarPatdens}, hence they induce maps $\phi_{\mu^\pm}:\Ga\bs
V\XX\ra [0,+\infty[$. In the case of real hyperbolic manifolds and
vanishing potentials, the total mass functions have important links to
the spectrum of the hyperbolic Laplacian (see \cite{Sullivan87}). See
also \cite{CooPap97,CooPap99} for the case of simplicial trees and the
discrete Laplacian, Section \ref{subsec:quantumgraph} for a
generalisation of the result of Coornaert and Papadopoulos, and for
instance \cite{BerKuc13} for developments in the field of quantum
graphs.

\section{Gibbs measures for metric and simplicial trees}
\label{subsec:ergodictrees}

Let $(\XX,\lambda)$ be a locally finite metric tree without terminal
vertices, let $X=|\XX|_\lambda$ be its geometric realisation, and let
$x_0\in V\XX$ be a basepoint.  Let $\Ga$ be a nonelementary discrete
subgroup of $\Aut(\XX,\lambda)$. Let $\wt F:T^1X\ra \RR$ be a
potential for $\Ga$, and let $\wt F^+=\wt F$, $\wt F^-= \wt
F\circ\iota$.  Let $\delta=\delta_{\Ga,F^\pm}$ be the critical
exponent of $(\Ga,F^\pm)$, assumed to be finite.  Let
$(\mu^\pm_x)_{x\in V\XX}$ be two (normalised) Patterson densities on
$\partial_\infty X$ for the pairs $(\Ga,F^\pm)$.

 The {\em Gibbs measure
  $\gls{gibbsmeasuresimplicial}$}\index{Gibbs!measure} on the space of
discrete geodesic lines $\G\XX$ of $\XX$, invariant under $\Ga$ and
under the discrete  time geodesic flow $(\flow{t})_{t\in\ZZ}$ of
$\gengeod\XX$, is defined analogously with the continuous time case,
using the discrete Hopf parametrisation for any basepoint $x_0\in
V\XX$, by
\begin{equation}\label{eq:defigibbsdis}
d\wt m_{F}(\ell)=
e^{C_{\ell_-}^-(x_0,\,\ell(0))\,+\,C^+_{\ell_+}(x_0,\,\ell(0))}\;
d\mu_{x_0}^-(\ell_-)\,d\mu^+_{x_0}(\ell_+)\,dt\;,
\end{equation}
where now $dt$ is the counting measure on $\ZZ$. We again denote by
$\gls{gibbsmeasuresimplicialdown}$ the measure that $\wt m_{F}$
induces on $\Ga\bs \G\XX$.

In this Section, we prove that the Gibbs measures in the case of trees
satisfy a Gibbs property even closer to the one in symbolic dynamics,
we give an analytic finiteness criterion of the Gibbs measures for
metric trees, and we recall the ergodic properties of tree lattices.

\medskip 
As recalled in the introduction, Gibbs measures were first introduced
in statistical mechanics and consequently in symbolic dynamics, see
for example \cite{Bowen75}, \cite{ParPol90}, \cite{PauPolSha15}. In
order to motivate the terminology used in this book, we recall the
definition of a Gibbs measure for the full two-sided shift on a finite
alphabet:\footnote{See Section \ref{subsec:TMS} for the appropriate
  definition when the alphabet is countable.} Let $\Sigma_n = \{1,2,
\dots,n\}^\ZZ$ be the product space of sequences $x= (x_n)_{n\in\ZZ}$
indexed by $\ZZ$ in the finite discrete set $\{1,2,\dots,n\}$, and let
$\sigma:\Sigma_n\to\Sigma_n$ be the {\em shift map}\index{shift}
defined by $\sigma((x_n)_{n\in\ZZ})= (x_{n+1})_{n\in\ZZ}$.  A
shift-invariant probability measure $\mu$ on $\Sigma_n$ {\it satisfies
  the Gibbs property}\index{Gibbs!property}%
\index{measure!satisfying the Gibbs property} for an energy function
$\phi:\Sigma_n\to\RR$ if
$$
\frac 1C\le\frac{\mu([a_{-m_-}, a_{-m_-+1},\dots,a_{m_+-1}, a_{m_+}])}
{e^{-P(m_-+m_++1)+\sum_{k=-m_-}^{m_+}\phi(\sigma^kx)}} \le C
$$ 
for some constants $C\geq1$ and $P\in\RR$ (called the {\em
  pressure}\index{pressure}) and for all $m_\pm$ in $\NN$ with $m_-
\leq m_+$ and $x$ in the cylinder $[a_{-m_-}, a_{-m_-+1}, \dots,
  a_{m_+-1}, a_{m_+}]$ that consists of those $x\in\Sigma_n$ for which
$x_k=a_k$ for all $k\in[-m_-,m_+]$.
 
Let $x_-,x_+\in V\XX$ and let $x_0\in V\XX\cap [x_-,x_+]$. Let us
define the {\em tree cylinder}\index{tree cylinder}\index{cylinder!tree} of
the triple $(x_-,x_0,x_+)$ by
$$
\cyl{x_-,x_0,x_+}=\{\ell\in\G\XX:
\ell_\pm\in\OOO_{x_0}\{x_\pm\},\ \ell(0)=x_0\}\,.
$$
These cylinders are close to the dynamical balls that have been
introduced in Section \ref{subsec:gibbsmeasure}, and the parallel with
the symbolic case is obvious, as this cylinder is the set of geodesic
lines which coincides on $[-m_-,m_+]$, where $m_\pm=d(x_0,x_\pm)$,
with a given geodesic line passing through $x_\pm$ and through $x_0$
at time $t=0$. The Gibbs measure $\wt m_F$ on the space of discrete
geodesic lines $\G \XX$ satisfies a variant of the Gibbs property
which is even closer to the one in symbolic dynamics than the general
case described in Proposition \ref{prop:gibbsgibbs}.


\bprop[Gibbs property]\label{prop:Gibbsproptrees} Let $K$ be a finite
subset of $V\XX\cap\C\Lambda\Ga$.  There exists $C'>1$ such that for
all $x_\pm\in \Ga K$ and $x_0\in V\XX\cap[x_-,x_+]$,
$$
\frac 1{C'}\le \frac{\wt m_F(\cyl{x_-,x_0,x_+})}
{e^{-\delta\, d(x_-,\,x_+)+\int_{x_-}^{x_+}\wt F}} \le C'\;.
$$
\eprop

\dem The result is immediate if $d(x_-,x_+)$ is bounded, since the
above denominator and numerator take only finitely many values by the
finiteness of $K$ and by $\Ga$-invariance, and the numerator is
nonzero since $x_\pm\in\C\Lambda\Ga$, hence the tree cylinder
$\cyl{x_-,x_0,x_+}$ meets the support of $\wt m_F$. We may hence
assume that $d(x_-,x_+)\geq 2$. Using the invariance of $\wt m_F$
under the discrete time geodesic flow, we may thus assume that
$x_0\neq x_-,x_+$. By $\Ga$-invariance, we may assume that $x_-$
varies in the finite set $K$, and that $x_0$ is at distance $1$ from
$x_-$, hence also varies in a finite set.

Using the discrete Hopf parametrisation with respect to the vertex
$x_0$, we have, by Lemma \ref{lem:shadowlemmatrees}, for some $C>0$,
\begin{align*}
\wt m_F(\cyl{x_-,x_0,x_+}) & 
=\mu^-_{x_0}(\OOO_{x_0}\{x_-\})\,\mu^+_{x_0}(\OOO_{x_0}\{x_+\})\\
&\leq C^2 \,e^{\int_{x_0}^{x_-}(\wt F\circ\iota -\delta)}\,
e^{\int_{x_0}^{x_+}(\wt F-\delta)}= C^2 e^{\int_{x_-}^{x_+}(\wt F-\delta)}\;.
\end{align*}
This proves the upper bound in Proposition \ref{prop:Gibbsproptrees}
with $C'=C^2$ and the lower bound follows similarly.  
\cqfd

\medskip
Next, we give a finiteness criterion of the Gibbs measure for metric
trees in terms of the total mass functions of the Patterson densities,
extending the case when $\Ga$ is torsion free and $\wt F=0$, due to
\cite[Theo.~1.1]{CooPap97b}.

\bprop \label{prop:majoGibbsL2} Let $(\XX,\lambda,\Ga,\wt F)$ be as in
the beginning of this Section.
\begin{enumerate}
\item If $(\XX,\lambda)$ is simplicial and $\|\cdot\|_2$ is the
  Hilbert norm of $\Leb^2(\Ga\bs V\XX, \vol_{\Ga\dbs\XX})$, we
  have\footnote{The maps $\phi_{\mu^\pm}$ are defined at the end of
    Section \ref{subsec:pattersongibbstrees}.}
$$
\|m_{F}\|\leq {\|\phi_{\mu^+}\|}_2\;{\|\phi_{\mu^-}\|}_2\;.
$$
\item In general, with $\|\cdot\|_2$ the Hilbert norm of
  $\Leb^2(\Ga\bs E\XX, \Tvol_{(\Ga\dbs\XX,\lambda)})$, we
  have\footnote{Recall that $o:\Ga\bs E\XX\to\Ga\bs V\XX$ is the
    initial vertex map, see Section \ref{subsec:trees}.}
$$
\|m_{F}\|\leq {\|\phi_{\mu^+}\circ o\|}_2\;{\|\phi_{\mu^-}\circ o\|}_2\;.
$$
\end{enumerate}
\eprop

\dem 
(1) The simplicial assumption on $(\XX,\lambda)$ means that all edges
have length $1$.  The space $\Ga\bs\G \XX$ is the disjoint union of
the subsets $\{\ell\in \Ga \bs\G \XX\;:\; \pi(\ell)= \ell(0)=[x]\}$ as
the orbit $[x]=\Ga x$ of $x\in V\XX$ ranges over $\Ga\bs V \XX$. By
Equation \eqref{eq:defigibbsdis}, using the discrete Hopf
decomposition with respect to the basepoint $x$, we have
$$
d(\wt m_{F})\big|_{\{\ell\in\G \XX \;:\;
  \ell(0)=x\}}(\ell)= d\mu_{x}^-(\ell_-)\, d\mu^+_{x}(\ell_+)\;.
$$ 
Let $\Dirac_{[x]}$ be the unit Dirac mass at $[x]$. By ramified
covering arguments, we hence have the following equality of measures
on the discrete set $\Ga\bs V \XX$:
\begin{equation}\label{eq:decompGibbsbas}
\pi_*m_F=\sum_{[x]\in \Ga\bs V \XX} \frac{1}{|\Ga_x|}\;
\big(\mu^-_x\times\mu^+_x\big)\big(\{(\ell_-,\ell_+)\in\partial_\infty^2 X
\;:\;x\in\;]\ell_-,\ell_+[\}\big)\;\Dirac_{[x]}\:.
\end{equation}
Thus, using the Cauchy-Schwarz inequality and by the definition of the
measure $\vol_{\Ga\dbs \XX}$,
\begin{align*}
\|m_{F}\|=\|\pi_*m_{F}\|&\leq
\sum_{[x]\in \Ga\bs V \XX} \frac{1}{|\Ga_x|}\;\|\mu^-_x\times\mu^+_x\|
={\langle \phi_{\mu^-},\phi_{\mu^+}\rangle}_2 \\ &\leq 
{\|\phi_{\mu^-}\|}_2{\|\phi_{\mu^+}\|}_2\;.
\end{align*}
This proves Assertion (1) of Proposition \ref{prop:majoGibbsL2}.

\medskip\noindent 
(2) The argument is similar to the previous one. Since the singletons
in $\RR$ have zero Lebesgue measure, the space $\Ga\bs\G X$ is, up to
a measure zero subset for $m_F$, the disjoint union for $[e]\in\Ga\bs
E\XX$ of the sets $A_{[e]}$ consisting of the elements $\ell\in
\Ga\bs\G X$ such that $\ell(0)$ belongs to the interior of the edge
$[e]$ and the orientations of $\ell$ and $e$ coincide on $e$. We fix a
representative $e$ of each $[e]\in\Ga\bs E\XX$. For every $t\in[0,
  \lambda(e)]$, let $e_t$ be the point of $e$ at distance $t$ from
$o(e)$. By Equation \eqref{eq:defigibbs}, using Hopf's decomposition
with respect to the basepoint $o(e)$ in $A_{[e]}$, we have as above
\begin{align*}
\|m_{F}\|=\sum_{[e]\in\Ga\bs E\XX}\frac{1}{|\Ga_e|}\; & 
\int_{\ell_-\in\partial_{\ov{e}} \XX} \int_{\ell_+\in\partial_{e} \XX} 
\int_{0}^{\lambda(e)}  \\  &
e^{C_{\ell_-}^-(o(e),\,e_t)+
C^+_{\ell_+}(o(e),\,e_t)}\;d\mu_{o(e)}^-(\ell_-)\;d\mu^+_{o(e)}(\ell_+)
\;dt \;.
\end{align*}
Since $\ell_-, o(e), e_t,\ell_+$ are in this order on the geodesic
line $\ell$ with $\ell_-\in \partial_{\ov{e}} \XX$ and $\ell_+\in
\partial_{e} \XX$, we have $C_{\ell_-}^-(o(e),\,e_t)+
C^+_{\ell_+}(o(e),\,e_t)=0$ by Equation \eqref{eq:changemoinsplus}.
Hence, by the definition of the measure
$\Tvol_{(\Ga\dbs\XX,\lambda)}$,\footnote{See Section
  \ref{subsec:trees}.}
\begin{align*}
\|m_{F}\|&=\sum_{[e]\in\Ga\bs E\XX} \frac{\lambda(e)}{|\Ga_e|}\; 
\mu_{o(e)}^-(\partial_{\ov{e}} \XX)\;\mu^+_{o(e)}(\partial_{e} \XX)\\ &\leq
\sum_{[e]\in\Ga\bs E\XX} \frac{\lambda(e)}{|\Ga_e|}\; 
\|\mu_{o(e)}^-\|\;\|\mu^+_{o(e)}\|=
\langle\phi_{\mu^-}\circ o,\phi_{\mu^+}\circ o\rangle_2\\ &\leq
{\|\phi_{\mu^-}\circ o\|}_2\;{\|\phi_{\mu^+}\circ o\|}_2\;,
\end{align*}
which finishes the proof.
\cqfd

\medskip Let us give some corollaries of this proposition in
the case of simplicial trees. It follows from Assertion (1) of
Proposition \ref{prop:majoGibbsL2} that if the $\Leb_2$-norms of the
total mass of the Patterson densities are finite, then the Gibbs
measure $m_{F}$ is finite. Taking $\wt F=0$ and $(\mu_x^+)_{x\in V\XX}
= (\mu_x^-)_{x\in V\XX}$, so that the Gibbs measure $m_F$ is the
Bowen-Margulis measure $m_{\rm BM}$, it follows from this proposition
that
\begin{equation}\label{eq:majoBM}
\|m_{\rm BM}\|\leq {{\|\phi_{\mu^\pm}\|}_2}^2\leq 
\Vol(\Ga\dbs\XX)\sup_{x\in V\XX}\|\mu^\pm_x\|^2\;.
\end{equation}
In particular, if $\Ga$ is a (tree) lattice\footnote{See Section
  \ref{subsec:trees}.}  of $\XX$ and if the total mass of the
Patterson density is bounded, then the Bowen-Margulis measure $m_{\rm
  BM}$ is finite.

\bigskip 
The following statement summarises the basic ergodic properties of the
lattices of $(\XX,\lambda)$ when $F=0$.

\bprop \label{prop:uniflatmBMfinie} Let $(\XX,\lambda)$ be a metric or
simplicial tree, with geometric realisation $X$.  Assume that
$(\XX,\lambda)$ is uniform and that $\Ga$ is a lattice in
$\Aut(\XX,\lambda)$.  Then
\begin{enumerate}
\item $\Ga$ is of divergence type, and its critical exponent
  $\delta_\Ga$ is the Hausdorff dimension of any visual distance $d_x$
  on $\partial_\infty X=\Lambda\Ga$;
\item the Patterson density $(\mu_x)_{x\in X}$ coincides, up to a
  scalar multiple, with the family of Hausdorff measures
  $\gls{hausdorffmeasure}$ of dimension $\delta_\Ga$ of the visual
  distances $(\partial_\infty X, d_x)$; in particular, it is
  $\Aut(\XX,\lambda)$-equivariant: for all $\ga\in\Aut(\XX,\lambda)$
  and $x\in X$, we have $\ga_*\mu_x=\mu_{\ga x}$;
\item the Bowen-Margulis measure $\wt m_{\rm BM}$ of $\Ga$ on $\G X$ is
$\Aut(\XX,\lambda)$-invariant, and the Bowen-Margulis measure $m_{\rm
  BM}$ of $\Ga$ on $\Ga\bs\G X$ is finite.  
\end{enumerate}
\eprop

\dem 
Let $\Ga'$ be any uniform lattice of $(\XX,\lambda)$, which exists
since the metric tree $(\XX,\lambda)$ is uniform.  It is well-known
(see for instance \cite{Bourdon95}) that the critical exponent
$\delta_{\Ga'}$ of $\Ga'$ is finite and equal to the Hausdorff
dimension of any visual distance $(\partial_\infty X, d_x)$, and that
the family $(\mu_x^{\rm Haus})_{x\in VX}$ of Hausdorff measures of
the visual distances $(\partial_\infty X, d_x)$ is a Patterson density
for any discrete nonelementary subgroup of $\Aut(\XX,\lambda)$ with
critical exponent equal to $\delta_{\Ga'}$.

By \cite[Coro.~6.5(2)]{BurMoz96}, the lattice $\Ga$ in $\Aut(\XX,
\lambda)$ is of divergence type and $\delta_\Ga= \delta_{\Ga'}$. By
the uniqueness property of the Patterson densities when $\Ga$ is of
divergence type (see Theorem \ref{theo:HTSR} (3)), the family
$(\mu_x)_{x\in VX}$ coincides, up to a scalar multiple, with
$(\mu_x^{\rm Haus})_{x\in VX}$.

As the graph $\Ga'\bs\XX$ is compact, the total mass function of the
Hausdorff measures of the visual distances is bounded, hence so is
$(\|\mu_x\|)_{x\in VX}$. By Proposition \ref{prop:majoGibbsL2},
since $\Ga$ is a tree lattice of $(\XX,\lambda)$, hence of $\XX$, this
implies that the Bowen-Margulis measure $m_{\rm BM}$ of $\Ga$ is
finite.  
\cqfd

\medskip Note that as in \cite{DalOtaPei00}, when $(\XX, \lambda)$ (or
its minimal nonempty $\Ga$-invariant subtree) is not assumed to be
uniform, there are examples of $\Ga$ that are lattices (or are
geometrically finite) whose Bowen-Margulis measure $m_{\rm BM}$ is
infinite, see Section \ref{subsec:exampleinfiniteBM} for more
details.

\bigskip 
Assume till the end of this Section that $(\XX,\lambda)$ is
simplicial, that is, that $\lambda\equiv 1$. Let us now discuss the
mixing properties of the discrete time geodesic flow on $\Ga\bs \G\XX$
for the Gibbs measure $m_F$.

\medskip Let $\gls{lengthspectrum}$ be the {\em length
  spectrum}\index{length!spectrum} of $\Ga$, which is, in the present
simplicial case, the subgroup of $\ZZ$ generated by the translation
lengths in $\XX$ of the elements of $\Ga$.

Recall that $x_0\in V\XX$ is a fixed basepoint. Let
$\gls{setofevenvertices}= \{x\in V\XX\;:\; d(x,x_0) = 0\mod 2\}$ be
the set of vertices of $\XX$ at an even distance from the basepoint
$x_0$, and let $\gls{setofoddvertices}= V\XX-\Veven\XX$. Let
$\gls{evendiscgeodflowspace}$ (respectively
$\gls{evendiscgenegeodflowspace}$) be the subset of $\G\XX$
(respectively $\gengeod\XX$) that consists of the geodesic lines
(respectively generalised geodesic lines) in $\XX$ whose origin is in
$\Veven\XX$.

Recall\footnote{See above Theorem \ref{theo:mixing} for the
  continuous-time version} that a discrete time one-parameter group
$(h^n)_{n\in\ZZ}$ of homeomorphisms of a topological space $Z$ is {\it
  topologically mixing}\index{topologically!mixing} if for all
nonempty open subsets $U,V$ of $Z$, there exists $n_0\in\NN$ such that
for all $n\ge n_0$, we have $U\cap h^n (V)\neq \emptyset$.

Recall that given a measured space $(Z,m)$, with $m$ nonzero and
finite, endowed with a discrete time one-parameter group
$(h^n)_{n\in\ZZ}$ of measure-preserving transformations, the measure
$m$ is {\it mixing}\index{mixing} under $(h^n)_{n\in\ZZ}$ if for all
$f,g\in \LL^2(Z,m)$, we have
$$
\lim_{n\ra +\infty}\int_Z f\;g\circ h^n\;dm=
\frac{1}{\|m\|}\int_Z f\;dm\;\int_Z g\;dm\;,
$$
or equivalently if for every $g\in \LL^2(Z,m)$, the functions
$g\circ h^n$ weakly converge in the Hilbert space $\LL^2(Z,m)$ to the
constant function $\frac{1}{\|m\|}\int_Z g\;dm$ as $n\ra+\infty$.

\btheo \label{theo:uniflatmBMmixing} Assume that the smallest nonempty
$\Ga$-invariant simplicial subtree of $\XX$ is uniform, without
vertices of degree $2$, and that $m_F$ is finite. Then the following
assertions are equivalent:
\begin{itemize}
\item the length spectrum of $\Ga$ satisfies
  $\loops{\Ga}=\ZZ$;
\item the discrete time geodesic flow on $\Ga\bs \G\XX$ is
  topologically mixing on its nonwandering set;
\item the quotient graph $\Ga\bs\XX$ is not bipartite;
\item the Gibbs measure $m_{F}$ is mixing under the discrete time
  geodesic flow $(\flow t)_{t\in\ZZ}$ on $\Ga\bs\G\XX$.
\end{itemize}
Otherwise $\loops{\Ga}=2\ZZ$, and the square of the discrete time
geodesic flow $(\flow{2t})_{t\in\ZZ}$ is topologically mixing on the
nonwandering subset of $\Ga\bs\Geven\XX$ and the restriction of the
Gibbs measure $m_{F}$ to $\Ga\bs\Geven\XX$ is mixing under
$(\flow{2t})_{t\in\ZZ}$.  
\etheo

\dem The nonwandering set of $(\flow t)_{t\in\ZZ}$ on $\Ga\bs\G\XX$ is
$\Ga\bs \{\ell\in\G\XX\;:\;\ell_\pm\in \Lambda\Ga\}$, and the
nonwandering set of $(\flow {2t})_{t\in\ZZ}$ on $\Ga\bs\Geven\XX$ is
$$
\Omega_{\rm even}=\Ga\bs \{\ell\in\Geven\XX\;:\;\ell_\pm\in
\Lambda\Ga\}\;.
$$
Since the translation axis of any loxodromic element of $\Ga$ is
contained in the convex hull of the limit set, we may hence assume
that the geometric realisation of $\XX$ is equal to $\C\Lambda\Ga$.

\blemm\label{lem:GabLev}
If $\XX$ is a locally finite tree without vertices of degree $2$, if
$\Ga$ is a nonelementary discrete sugbroup of $\Aut(\XX)$ such that
$\XX$ is {\rm tree-minimal} (that is, does not contain a
$\Ga$-invariant proper nonempty subtree), then the length spectrum
$\loops{\Ga}$ of $\Ga$ is equal either to $\ZZ$ or to $2\ZZ$, and
equal to $2\ZZ$ if and only if the quotient graph $\Ga\bs\XX$ is
bipartite.
\elemm

\dem This lemma is essentially due to \cite{GabLev95}. By for instance
\cite[Lem.~4.3]{Paulin89a}, since $\XX$ is tree-minimal, every
geodesic segment (and in particular any two consecutive edges) is
contained in the translation axis of a loxodromic element of
$\Ga$. Hence if $x$ and $y$ are the two endpoints of any edge $e$ of
$\XX$, since they have degree at least $3$, there exist at least two
loxodromic elements $\alpha$ and $\beta$ of $\Ga$ such that the
translation axes $\Ax_\alpha$ and $\Ax_\beta$ contain $x$ and $y$
respectively, but do not meet the interior of the edge $e$, so that
$d(x,y)=d(\Ax_\alpha,\Ax_\beta)$. In particular $\Ax_\alpha$ and
$\Ax_\beta$ are disjoint, which implies by for instance
\cite[Prop.~1.6]{Paulin89a} that the translation lengths of $\alpha$,
$\beta$ and $\alpha\beta$ satisfy
$$
\lambda(\alpha\beta)=
\lambda(\alpha)+\lambda(\beta)+2\;d(\Ax_\alpha,\Ax_\beta)\;.
$$
Hence $2=2\,d(x,y)=\lambda(\alpha\beta)-\lambda(\alpha)-\lambda(\beta)
\in \loops{\Ga}$. Therefore $2\ZZ\subset \loops{\Ga}\subset \ZZ$, and
either $\loops{\Ga}= \ZZ$ or $\loops{\Ga}= 2\ZZ$.

Note that for all vertices $x,y,z$ in a simplicial tree, if $d(x,y)$
and $d(y,z)$ are both even or both odd, then
\begin{equation}\label{eq:pairtriangle}
d(x,z)=d(x,y)+d(y,z)-2\,d(y,[x,z])=0\mod 2\;.
\end{equation}
Note that for all $x\in V\XX$ and $\ga\in\Ga$, we have
\begin{equation}\label{eq:pairdlamb}
d(x,\ga x) = \lambda(\ga) \mod 2\;.
\end{equation}
Indeed, if $\ga$ is loxodromic, then $d(x,\ga x) = \lambda(\ga) +2\;
d(x,\Ax_\ga)$ and otherwise, $d(x,\ga x) = 2\;d(x,
\operatorname{Fix}(\ga))$ where $\operatorname{Fix}(\ga)$ is the set
of fixed points of the elliptic element $\ga$.  For future use, this
proves that the following assertions are equivalent :
\begin{align}
&(1)\;\;\; \loops{\Ga}\subset 2\ZZ\nonumber\\
&(2)\;\;\; \forall\;x\in X,\;\;\forall\; \ga\in\Ga,\;\;\;\;
d(x,\ga x)\in 2\ZZ\;.\label{eq:equivLamGaevendisteven}
\end{align}

Assume that $\loops{\Ga} = 2\ZZ$. Then $\Veven\XX$ (hence $\Vodd\XX$) is
$\Ga$-invariant, since for all $x\in \Veven\XX$, the distance $d(x,\ga x)$
is even by Equation \eqref{eq:pairdlamb}, and $d(\ga x,\ga x_0)=d(x,
x_0)$ is even, so that $d(\ga x,x_0)$ is even by Equation
\eqref{eq:pairtriangle}.  Since no edge of $\XX$ has both endpoints in
$\Veven\XX$, this proves that $\Ga\bs \XX$ is bipartite, with partition of
its set of vertices $(\Ga\bs \Veven\XX)\sqcup (\Ga\bs \Vodd\XX)$.

Assume conversely that $\Ga \bs \XX$ is bipartite. The set $\Veven\XX$,
which is the lift of one of the two elements of the partition of its
vertices by the canonical projection $V\XX\ra \Ga\bs V\XX$, is
$\Ga$-invariant. By Equation \eqref{eq:pairdlamb}, this proves that
$\loops{\Ga}\subset 2\ZZ$, hence that $\loops{\Ga}= 2\ZZ$.
\cqfd

\medskip
The equivalence of the first, second and fourth claims in the
statement of Theorem \ref{theo:uniflatmBMmixing} follows from a
discrete time version with potential of \cite[Theo.~3.1]{Roblin03} or
a discrete time version of \cite[Theo.~1]{Babillot02b} (which can be
extended to $\CAT(-1)$ spaces by the remark in \cite[page
  70]{Babillot02b}).  It can also be recovered from the following
arguments when $\loops{\Ga} =2\ZZ$, and we prefer to concentrate on
this case, since it requires a lot of modifications and is stated
with almost no proof in \cite[Prop.~3.3]{BroPau07Tou}, and
only when $\wt F=0$.

Assume from now on that $\loops{\Ga}=2\ZZ$. Since $\Veven\XX$ is
$\Ga$-invariant as seen above, and since $d(\ell(0),\flow{2s}\ell(0))
=2|s|$ is even for all $\ell\in\G \XX$ and $s\in\ZZ$, it follows from
the definition of $\Geven\XX=\{\ell\in\G \XX\;:\;\pi(\ell)\in \Veven
\XX\}$ and from Equation \eqref{eq:pairtriangle} that $\Geven\XX$ is
invariant under the even discrete time geodesic flow
$(\flow{2s})_{s\in\ZZ}$ and under $\Ga$. Note that the discrete Hopf
parametrisation of $\G\XX$ gives a homeomorphism from $\Geven\XX$ to
$$
\big\{(\xi,\eta,t)\in\partial_\infty^2\XX\times\ZZ\;:\;
t=d(x_0,\mathopen{]}\xi,\eta\mathclose{[}\,)\mod 2\big\}\;.
$$
The restriction of the Gibbs measure $\wt m_F$ to $\Geven\XX$, that we
will again denote by $\wt m_F$, disintegrates by the projection on the
first factor $\partial_\infty^2\XX\times\ZZ \ra \partial_\infty^2\XX$
over the geodesic current $\wh m_F$ where, for every $(\xi,\eta)\in
\partial_\infty^2\XX$ and (any) $x\in\mathopen{]}\xi,\eta\mathclose{[}$,
\begin{equation}\label{eq:geodcurrquasiprod}
d\wh m_F(\xi,\eta)=e^{C_{\xi}^-(x_0,\,x)\,+\,C^+_{\eta}(x_0,\,x)}
d\mu_{x_0}^-(\xi)\,d\mu^+_{x_0}(\eta)\;,
\end{equation}
with conditional measure on the fiber over $(\xi,\eta)$ the counting
measure on the discrete set $\{t\in\ZZ\;:\; t=d(x_0,\mathopen{]}\xi,
\eta\mathclose{[}\,)\mod 2\}$. Since $m_F$ is finite and invariant
under the discrete time geodesic flow, it is conservative by
Poincaré's recurrence theorem. Hence the measure quasi-preserving
action of $\Ga$ on the measured space $(\partial_\infty^2\XX,\wh
m_F)$ is ergodic by (the discrete time version of) Theorem
\ref{theo:HTSR}.

Since the distance between two points in a horosphere of a simplicial
tree is even, and again by Equation \eqref{eq:pairtriangle}, every
horosphere of $\XX$ is either entirely contained in $\Veven \XX$ or
entirely contained in $\Vodd \XX$. For every $\ell\in\Geven\XX$, its
strong stable/unstable leaf
$$
W^\pm(\ell)=\{\ell'\in\G\XX\;:\;\lim_{t\ra\pm\infty} d(\ell(t),
\ell'(t)) =0\}
$$
is contained in $\Geven\XX$, since the image by the footpoint
projection of a strong stable/instable leaf is a horosphere (see
Equation \eqref{eq:footpowssuhpm}). Thus $\Geven\XX$ is saturated by
the partition into strong stable/instable leaves of $\G\XX$.

We now follow rather closely the arguments of
\cite[Theo.~1]{Babillot02b} in order to prove the last claims of
Theorem \ref{theo:uniflatmBMmixing}, the main point being that the
geodesic current $\wh m_F$ is a quasi-product measure.

The following lemmas are particular cases of respectively Lemma 1 and
Fact page 64 of \cite {Babillot02b}, valid for general finite measure
preserving dynamical systems, applied, with the notation of loc.~cit,
to $(T_t)_{t\in A}=(\flow{2t})_{t\in\ZZ}$.

\blemm\label{lem:bablem1} Let $f\in\LL^2(\Ga\bs\Geven\XX,m_F)$ be such
that $\int f\,dm_F=0$. If there exists an increasing sequence
$(t_n)_{n\in\NN}$ in $\ZZ$ such that $f\circ \flow{2t_n}$ does not
converge to $0$ for the weak topology on $\LL^2(\Ga\bs\Geven\XX,m_F)$,
then there exist an increasing sequence $(s_n)_{n\in\NN}$ in $\ZZ$ and
a nonconstant\footnote{Recall that an element of $\LL^2(Z,m)$ is
  nonconstant if any representative function is not almost everywhere
  constant.}  element $f^*\in\LL^2(\Ga\bs\Geven\XX,m_F)$ such that
$f\circ \flow{2s_n}$ and $f\circ \flow{-2s_n}$ both converge to $f^*$
for the weak topology on $\LL^2(\Ga\bs\Geven\XX,m_F)$.  \cqfd 
\elemm

\blemm\label{lem:babfact}
If $(f_n)_{n\in\NN}$ is a sequence in $\LL^2(\Ga\bs\Geven\XX,m_F)$
weakly converging to $f^*$ in $\LL^2(\Ga\bs\Geven\XX,m_F)$, then there
exists a subsequence $(f_{n_k})_{k\in\NN}$ such that the Ces\`aro averages
$\frac{1}{N^2}\sum_{k=0}^{N^2-1} f_{n_k}$ converge pointwise almost everywhere
to $f^*$ as $N\ra+\infty$.
\cqfd \elemm

Recall that the support of $m_F$ is the nonwandering set $\Omega_{\rm
  even}$.  Assume for a contradiction that the restriction of $m_F$ to
$\Ga\bs\Geven\XX$ is not mixing under the even discrete time geodesic
flow. Then there exists a continuous function $f$ with compact support
on $\Omega_{\rm even}$ such that $\int f\;dm_F=0$ and $(f\circ
\flow{2n})_{n\in\NN}$ does not weakly converge to $0$ in
$\LL^2(\Ga\bs\Geven\XX,m_F)$. By Lemmas \ref{lem:bablem1} and
\ref{lem:babfact}, there exist a nonconstant element $f^*\in
\LL^2(\Ga\bs\Geven\XX,m_F)$ and increasing sequences
$(n_k^\pm)_{k\in\NN}$ in $\NN$ such that $\frac{1}{N^2}
\sum_{k=0}^{N^2-1} f\circ \flow{\pm 2n^\pm_k}$ pointwise almost
everywhere converges to $f^*$ as $N\ra+\infty$.

Let $\wt f^*=f^*\circ p_{\rm even}$, where $p_{\rm even}:\Geven\XX\ra
\Ga\bs\Geven\XX$ is the canonical projection, be the lift of $f^*$ to
$\Geven\XX$. Since the conditional measures for the disintegration of
$\wt m_F$ over $\wh m_F$ are counting measures on countable sets,
there exists a full $\wh m_F$-measure subset $E_0$ of
$\partial_\infty^2 \XX$ such that, for every $\ell\in \Geven\XX$ with
$(\ell_-,\ell_+)\in E_0$, the above convergences hold after lifting to
$\Geven\XX$ at the points $\flow{2n}\ell$ for all $n\in\ZZ$.

For every $\ell\in \Geven\XX$, the subgroup $A_\ell$ of $2\ZZ$ given
by the periods of the map $2n\mapsto \wt f^*(\flow{2n}\ell)$ only
depends on $(\ell_-,\ell_+)$.  Thus, we have a measurable map from
$E_0$ into the (discrete) set of subgroups of $2\ZZ$, which is
$\Ga$-invariant, hence is constant $\wh m_F$-almost everywhere by the
ergodicity of $\wh m_F$ under $\Ga$.

Assume for a contradiction that this almost everywhere constant
subgroup is $2\ZZ$, that is, that the values of $\wt f^*$ almost
everywhere do not depend on the time parameter in the discrete Hopf
parametrisation of $\Geven\XX$. Then $\wt f^*$ defines a
$\Ga$-invariant measurable function on $\partial_\infty^2\XX$.  Again
by ergodicity, this function is almost everywhere constant,
contradicting the fact that $f^*$ is not almost everywhere constant.

Hence there exist a full $\wh m_F$-measure subset $E_1$ of $E_0$ and
$\kappa\in\NN-\{0,1\}$ such that $A_\ell=2\,\kappa\,\ZZ$ for every
$\ell\in \Geven\XX$ with $(\ell_-,\ell_+)\in E_1$.  Let us finally
prove that $\loops{\Ga}$ is contained in $2\,\kappa\,\ZZ$, which
contradicts the original assumption that $\loops{\Ga}=2\ZZ$.

Let $\wt f^\pm=\limsup_{N\ra+\infty}\;\frac{1}{N^2} \sum_{k=0}^{N^2}
f\circ \flow{\pm 2 n^\pm_k}\circ p_{\rm even}$, so that the set
$$
E=\{(\xi,\eta)\in E_1\;:\;
\forall\;\ell\in\Geven\XX,\;{\rm if}\;\ell_-=\xi\;{\rm and}\;
\ell_+=\eta,\;{\rm then} \;\wt f^+(\ell)=\wt f^-(\ell)=\wt f^*(\ell)\}
$$
has full $\wh m_F$-measure. By the hyperbolicity of the geodesic flow
(see Equation \eqref{eq:expandHamdist}) and the uniform continuity of
$f$, the map $\wt f^+$ is constant along any strong stable leaf of
$\Geven\XX$ and $\wt f^-$ is constant along any strong unstable
leaf. Let
$$
E^-= \{\xi\in\Lambda\Ga\;:\;(\xi,\eta')\in E
\;{\rm for}\; \mu_{x_0}^+{\text -{\rm almost~every}}\;\eta'\in\Lambda\Ga\}
$$
and
$$
E^+= \{\eta\in\Lambda \Ga \;:\;(\xi',\eta)\in E
\;{\rm for}\; \mu_{x_0}^- {\text -{\rm almost~every}} \;\xi'\in\Lambda\Ga\}\;.
$$

Since $\wh m_F$ is in the same measure class as the product measure
$\mu^-_{x_0}\otimes \mu^+_{x_0}$ (see Equation
\eqref{eq:geodcurrquasiprod}), and by Fubini's theorem, we have
$\mu^-_{x_0}(\;^c E^-)=0$ and $\mu^+_{x_0}(\;^c E^+)=0$, and the set
$E^-\times E^+$ has full $\wh m_F$-measure.

Let $\ga$ be a loxodromic element of $\Ga$, and let $x$ be any vertex
of $\XX$ on the translation axis of $\ga$. Since $d(x,\ga x)$ is even
(see Equation \eqref{eq:equivLamGaevendisteven}), the midpoint $y$ of
the geodesic segment $[x,\ga x]$ is a vertex of $\XX$. Since $x$ and
$y$ have degree at least $3$, there exist $\xi_x$ and $\xi_y$ in
$\partial_\infty X$ whose closest points on the translation axis of
$\ga$ are respectively $x$ and $y$. Note that $\xi_x,\xi_y \in
\Lambda\Ga$ as $X=\C\Lambda\Ga$. Since $y$ and $x$ are the closest
points to $\ga_+$ and $\ga_-$ on the geodesic line $\mathopen{]}
\xi_x, \xi_y\mathclose{[}$, we have\footnote{See Section
\ref{subsec:trees} for the definition of the crossratio of an
ordered quadruple of pairwise distinct points in
$\partial_\infty X$.} $\ldbrack\xi_x,\ga_+,\xi_y,\ga_-\rdbrack=d(x,y)$.

\begin{center}
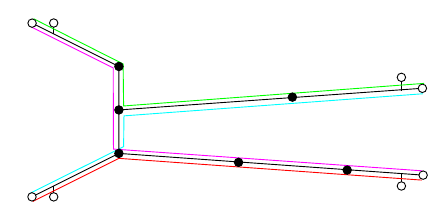
\end{center}

Since $E\cap(E^-\times E^+)$ has full $\wh m_F$-measure, there exists
$(\xi,\eta)\in E\cap (E^-\times E^+)$ arbitrarily close to $(\ga_-,
\xi_x)$. Since the set 
$$
\{(\xi',\eta')\in \partial^2_\infty\Ga\;:\;
(\xi,\eta'), (\xi',\eta), (\xi',\eta') \in E\}
$$ 
has full $\wh m_F$-measure, there exists such a $(\xi',\eta')$
arbitrarily close to $(\ga_+, \xi_y)$.  Let $\ell^0\in\Geven\XX$ be
such that $\ell^0_-=\xi$ and $\ell^0_+=\eta$. Let $\ell^1\in
W^+(\ell^0)$ be such that $\ell^1_-=\xi'$. Let $\ell^2\in W^-(\ell^1)$
be such that $\ell^2_+=\eta'$. Let $\ell^3\in W^+(\ell^2)$ be such
that $\ell^3_-=\xi$. Finally, let $\ell^4\in W^-(\ell^3)$ be such that
$\ell^4_+=\eta$. Then $\ell^4=\flow{2s}\ell^0$ for some $s\in\ZZ$ with
$2|s|= d(\ell^0(0),\ell^4(0))$.

By the definition of $E$ and since $\wt f^+$ (resp.~$\wt f^-$) is
constant along the strong stable (resp.~unstable) leaves, we have
\begin{align*}
\wt f^*(\ell^0)&=\wt f^+(\ell^0)=\wt f^+(\ell^1)=\wt f^-(\ell^1)=
\wt f^-(\ell^2)=\wt f^+(\ell^2)=\wt f^+(\ell^3)\\ &=\wt f^-(\ell^3)
=\wt f^-(\ell^4)=\wt f^*(\ell^4)=\wt f^*(\flow{2\,s}\ell^0)\;.
\end{align*}
Hence $2\,|s|$ is a period in $A_{\ell^0}$, thus is contained in
$2\,\kappa\,\ZZ$.

If $t>0$ is large enough, we have
$$
\ell^0(t)=\ell^1(t), \;\;\;\ell^1(-t)= \ell^2(-t), \;\;\;
\ell^2(t)=\ell^3(t), \;\;\;\ell^3(-t)= \ell^4(-t)\;,
$$
which respectively tend to $\eta$, $\xi'$, $\eta'$, $\xi$ as
$t\ra+\infty$.  Since the crossratio is locally constant and by its
properties, in particular its definition in Equation
\eqref{eq:defcrossratio}, we have
\begin{align*}
\ell(\ga)&=d(x,\ga x)=2\,d(x,y)= 2\,\ldbrack\xi_x,\ga_+,\xi_y,\ga_-\rdbrack=
2\,\ldbrack\eta,\xi',\eta',\xi\rdbrack \\&=\lim_{t\ra+\infty}
d(\ell^0(t),\ell^3(-t))-d(\ell^3(-t),\ell^2(t))+d(\ell^2(t),\ell^1(-t))-
d(\ell^1(-t),\ell^0(t))\\&= \lim_{t\ra+\infty}
d(\ell^0(t),\ell^3(-t)) -d(\ell^1(-t),\ell^1(t)) = \lim_{t\ra+\infty}
d(\ell^0(t),\ell^4(-t))-2\,t\\&=d(\ell^0(0),\ell^4(0))=2\,|s|\in
2\,\kappa\,\ZZ\;.
\end{align*}
Thus $\loops{\Ga}\subset 2\,\kappa\,\ZZ$, which contradicts the fact
that $\loops{\Ga}=2\ZZ$.
\cqfd

\medskip
By Proposition \ref{prop:uniflatmBMfinie}, the general assumptions of
Theorem \ref{theo:uniflatmBMmixing} are satisfied if $\XX$ is
uniform, without vertices of degree $2$, $\Ga$ is a lattice of $\XX$
and $\wt F=0$. Thus, if we assume furthermore that $\Ga\bs\XX$ is not
bipartite, then the Bowen-Margulis measure $m_{\rm BM}$ of $\Ga$ is
mixing under the discrete time geodesic flow on $\Ga\bs\G\XX$.

\chapter{Symbolic dynamics of geodesic flows on trees}
\label{sec:symbdyngeodflowarbre}

In this Chapter, we give a coding of the discrete time geodesic flow
on the nonwandering sets of quotients of locally finite simplicial
trees $\XX$ without terminal vertices by nonelementary discrete
subgroups of $\Aut(\XX)$ by a subshift of finite type on a countable
alphabet. Similarly we give a coding of the continuous time geodesic
flow on the nonwandering sets of quotients of locally finite metric
trees $(\XX,\lambda)$ without terminal vertices by nonelementary
discrete subgroups of $\Aut(\XX,\lambda)$ by suspensions of such
subshifts. These codings are used in Section \ref{subsec:varprinc} to
prove the variational principle in both contexts, and in Sections
\ref{subsec:mixingratesimpgraphs} and \ref{subsec:mixingratemetgraphs}
to obtain rates of mixing of the flows.

\section{Two-sided topological Markov shifts}
\label{subsec:TMS}

In this short and independent Section, that will be used in Sections
\ref{subsec:codagesimplicial}, \ref{subsec:codagemetric},
\ref{subsec:varprinc}, \ref{subsec:mixingratesimpgraphs} and
\ref{subsec:mixingratemetgraphs}, we recall some definitions
concerning symbolic dynamics on countable alphabets.\footnote{See for
  instance \cite{Kitchens98, Sarig15}.}

\medskip
A (two-sided) topological {\it Markov
  shift}\index{Markov!shift!two-sided}%
\index{topological Markov shift!two-sided}\footnote{Note that the
  terminology could be misleading, a topological Markov shift comes a
  prori without a measure, and many probability measures invariant
  under the shift do not satisfy the Markov chain property that the
  probability to pass from one state to another depends only on the
  previous state, not of all past states.} is a topological dynamical 
  system $(\Sigma,\sigma)$ constructed from a countable
discrete {\it alphabet} $\A$ and a {\it transition matrix} $A=
(A_{i,\,j})_{i,\,j\in\A}\in \{0,1\}^{\A\times\A}$, where $\Sigma$ is
the closed subset of the topological product space $\A^\ZZ$ defined by
$$
\Sigma=\big\{x=(x_n)_{n\in\ZZ}\in \A^\ZZ\;:\; \forall \;n\in\ZZ,\;\;\;
A_{x_n,x_{n+1}}=1\}\;,
$$ 
and $\sigma: \Sigma\ra \Sigma$ is the (two-sided) {\it
  shift}\index{shift} defined by
$$
(\sigma(x))_n=x_{n+1}
$$ 
for all $x\in \Sigma$ and $n\in\ZZ$. Note that to be given $(\A,A)$
is equivalent to be given an oriented graph with countable set of
vertices $\A$ (and set of oriented edges a subset of $\A\times \A$)
and with incidence matrix $A$ such that $A_{i,\,j}=1$ if there is an
oriented edge from the vertex $i$ to the vertex $j$ and $A_{i,\,j}=0$
otherwise.

For all $p\leq q$ in $\ZZ$, a finite sequence $(a_n)_{p\leq n\leq q}
\in\A^{\{p,\dots,q\}}$ is {\it admissible}\index{admissible} (or
{\it $A$-admissible} when we need to make $A$ precise) if
$A_{a_n,\,a_{n+1}}=1$ for all $n\in\{p,\dots, q-1\}$. A topological
Markov shift is {\it transitive}\index{transitive}%
\index{Markov!shift!transitive}%
\index{topological Markov shift!transitive} if for all $x,y\in\A$,
there exists an admissible finite sequence $(a_n)_{p\leq n\leq q}$
with $a_p=x$ and $a_q=y$. This is equivalent to requiring the dynamical
system $(\Sigma,\sigma)$ to be {\it topologically
  transitive}\index{topologically!transitive}: for all nonempty open
subsets $U,V$ in $\Sigma$, there exists $n\in\ZZ$ such that $U\cap
\sigma^n(V)\neq \emptyset$.

Note that the product space $\A^\ZZ$ is not locally compact when $\A$
is infinite. When the matrix $A$ has only finitely many nonzero
entries on each line and each colum, then $(\Sigma,\sigma)$ is also
called a {\it subshift of finite type (on a countable
  alphabet)}\index{subshift of finite type}. The topological space
$\Sigma$ is then locally compact: By diagonal extraction, for all $p
\leq q$ in $\ZZ$ and $a_{p},a_{p+1}, \dots, a_{q-1}, a_q$ in $\A$,
every {\em cylinder}\index{cylinder}
$$
[a_{p},a_{p+1}, \dots, a_{q-1}, a_q]=\big\{(x_n)_{n\in\ZZ}\in \Sigma\;:\; 
\forall \;n\in \{p,\dots, q\},\;\; x_n=a_n\big\}
$$
is a compact open subset of $\Sigma$.

Given a continuous map $F_{\rm symb}:\Sigma\ra\RR$ and a constant
$c_{F_{\rm symb}}\in\RR$, we say that a measure $\PP$ on $\Sigma$,
invariant under the shift $\sigma$, {\it satisfies the Gibbs
  property\index{Gibbs!property}%
\index{measure!satisfying the Gibbs property}\footnote{Note that some
  references have a stronger notion of Gibbs measure (see for instance
  \cite{Sarig03}), with the constant $C$ independent of $E$.} with
Gibbs constant\index{Gibbs!constant} $c_{F_{\rm symb}}$ for the
potential $F_{\rm symb}$} if for every finite subset $E$ of the
alphabet $\A$, there exists $C_E\geq 1$ such that for all $p\leq q$ in
$\ZZ$ and for every $x=(x_n)_{n\in\ZZ}\in\Sigma$ such that $x_p,x_q\in
E$, we have
\begin{equation}\label{eq:defiGibbspropertyZ}
\frac {1}{C_E}\le\frac{\PP([x_{p}, x_{p+1},\dots,x_{q-1}, x_{q}])}
{e^{-c_{F_{\rm symb}}(q-p+1)+\sum_{n=p}^{q}F_{\rm symb}(\sigma^n x)}} \le C_E\;.
\end{equation}

Two continuous maps $F_{\rm symb},F'_{\rm symb}:\Sigma\ra \RR$ are
{\em cohomologous}\index{cohomologous} if there exists a continuous map
$G:\Sigma\ra \RR$ such that 
$$
F'_{\rm symb}-F_{\rm symb}=G\circ\sigma - G\;.
$$

\section{Coding discrete time geodesic flows on simplicial trees}
\label{subsec:codagesimplicial}

Let $\XX$ be a locally finite simplicial tree without terminal
vertices, with $X=|\XX|_1$ its geometric realisation. Let $\Ga$
be a nonelementary discrete subgroup of $\Aut(\XX)$, and let
$\wt F:T^1X\ra \RR$ be a potential for $\Ga$.

In this Section, we give a coding of the discrete time geodesic flow
$(\flow{t})_{t\in\ZZ}$ on the nonwandering subset of $\Ga\bs\G \XX$ by
a locally compact transitive (two-sided) topological Markov shift.
This explicit construction will be useful later on in order to study
the variational principle (see Section \ref{subsec:varprinc}) and
rates of mixing (see Section \ref{subsec:mixingratesimpgraphs}).

The main technical aspect of this construction, building on \cite[\S
6]{BroPau07Tou}, is to allow the case when $\Ga$ has torsion. When
$\Ga$ is torsion free and $\Ga/\XX$ is finite, the construction is
well-known, we refer for instance to \cite{CooPap02} for a more
general setting when the potential is $0$. In order to consider
for instance nonuniform tree lattices, it is important to allow
torsion in $\Ga$. Our direct approach also
avoids the assumption that the discrete subgroup $\Ga$ is {\it
  full}\index{full}, that is, equal to the subgroup consisting of the
elements $g\in\Aut(\XX)$ such that $p\circ g=p$ where
$p:\XX\ra\Ga\bs\XX$ is the canonical projection, as in \cite{Kwon15}
(building on \cite[7.3]{BurMoz96}).

\medskip
Let $\XX'$ be the minimal nonempty $\Ga$-invariant simplicial subtree
of $\XX$, whose geometric realisation is $\C\Lambda\Ga$. Since we are
only interested in the support of the Gibbs measures, we will only
code the geodesic flow on the nonwandering subset $\Ga\bs \XX'$ of
$\Ga\bs\G \XX$. The same construction works with the full space
$\Ga\bs \G \XX$, but the resulting Markov shift is then not necessarily
transitive.

Let $(\YY,G_*)= \Ga\dbs \XX'$ be the quotient graph of groups of
$\XX'$ by $\Ga$ (see for instance Example \ref{exem:quotgraphgroup}),
and let $p:\XX'\ra\YY=\Ga\bs \XX'$ be the canonical projection. We
denote by $[1]=H$ the trivial double coset in any double coset set
$H\bs G/ H$ of a group $G$ by a subgroup $H$.

\medskip
We consider the alphabet $\A$ consisting of the triples $(e^-,h,e^+)$
where 

\smallskip $\bullet$~ $e^\pm\in E\YY$ satisfy $t(e^-)=o(e^+)$ and 

\smallskip $\bullet$~ $h\in \rho_{e^-}(G_{e^-})\bs G_{o(e^+)}
/\rho_{\,\overline{e^+}}(G_{e^+})$ satisfy $h\neq [1]$ if
$\overline{e^+}=e^-$.

\medskip\noindent This set is countable (and finite if and only if the
quotient graph $\Ga\bs\XX'$ is finite), we endow it with the discrete
topology.  We consider the (two-sided) topological Markov shift with
alphabet $\A$ and transition matrix $A_{(e^-,\,h,\,e^+),\,
  ({e'}^-,\,h',\,{e'}^+)} = 1$ if $e^+={e'}^-$ and $0$ otherwise.
Note that this matrix $A= (A_{i,j})_{i,j\in\A}$ has only finitely many
nonzero entries on each line and each column, since $\XX'$ is locally
finite and $\Ga$ has finite vertex stabilisers in $\XX'$.  We consider
the subspace
$$
\Sigma=\big\{(e_i^-,h_i,e^+_i)_{i\in\ZZ}\in\A^\ZZ\;:\; 
\forall\;i\in\ZZ,\;\;e^+_{i-1}=e_i^-\big\}
$$ 
of the product space $\A^\ZZ$, and the shift $\sigma:\Sigma\ra
\Sigma$ defined by $(\sigma(x))_i=x_{i+1}$ for all $(x_i)_{i\in\ZZ}$
in $\Sigma$ and $i$ in $\ZZ$. As seen above, $\Sigma$ is locally compact.

\medskip
Let us now construct a natural coding map $\Theta$ from
$\Ga\bs\G \XX'$ to $\Sigma$, by slightly modifying the construction of
\cite[\S 6]{BroPau07Tou}.

\begin{center}
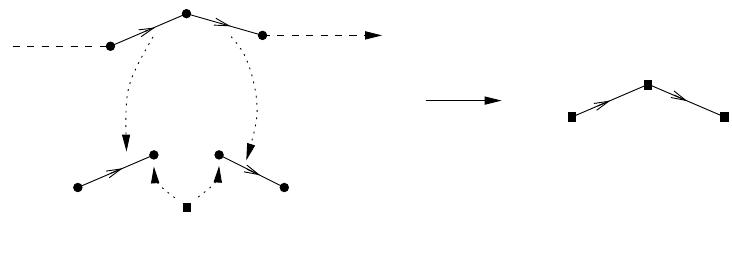
\end{center}

For every discrete geodesic line $\ell\in \G \XX'$, for every
$i\in\ZZ$, let $f_i=f_i(\ell)$ be the edge of $\XX'$ whose geometric
realisation is $\ell([i,i+1])$ with origin $f(i)$ and endpoint
$f(i+1)$, and let $e_i=p(f_i)$, which is an edge in $\YY$. Let us use
the notation of Example \ref{exem:quotgraphgroup}: we fix lifts $\wt
e$ and $\wt v$ of every edge $e$ and vertex $v$ of $\YY$ in $\XX'$
such that $\overline{\wt e}= \wt{\overline{e}}$, and elements
$g_e\in\Ga$ such that $g_e\,\wt{t(e)} =t(\wt e)$. Since
$p(\wt{e_i})=e_i=p(f_i)$, there exists $\ga_i= \ga_i(\ell)\in\Ga$,
well defined up to multiplication on the left by an element of
$G_{e_i}=\Ga_{\wt{e_i}}$, such that $\ga_i f_i= \wt{e_i}$ for all
$i\in\ZZ$.

We define $e_{i+1}^-(\ell)=e_{i}$, $e_{i+1}^+(\ell)=e_{i+1}$, and 
\begin{equation}\label{eq:defigi}
h_{i+1}(\ell)= g_{e_{i+1}^-(\ell)}^{\;\;-1}\,\ga_{i}(\ell)\,
\ga_{i+1}(\ell)^{-1}\,g_{\,\overline{e_{i+1}^+}(\ell)}\;.
\end{equation} 
Since for every edge $e$ of $\YY$ the structural monomorphism
$$
\rho_e:G_e=\Ga_{\wt e}\;\;\longrightarrow\;\; G_{t(e)}=\Ga_{\wt {t(e)}}
$$ 
is the map $g\mapsto g_e^{-1}gg_e$, the double coset of $h_i(\ell)$
in $\rho_{e_i^-(\ell)}(G_{e_i^-(\ell)}) \bs G_{o(e_i^+(\ell))} /
\rho_{\,\overline{e_i^+}(\ell)}(G_{e_i^+(\ell)})$ does not depend on the choice of
the $\ga_i$'s, and we again denote it by $h_i(\ell)$.

\medskip The next result shows that, assuming only that 
$\Ga$ is discrete and nonelementary, the time-one discrete
geodesic flow $\flow{1}$ on its nonwandering subset of $\Ga\bs\G
\XX$ is topologically conjugate to a locally compact
transitive (two-sided) topological Markov shift.

\btheo\label{theo:codingdisgeodflo} With $\XX'=\C\Lambda\Ga$, the map
$\Theta:\Ga\bs\G \XX' \ra\Sigma$ defined by 
$$
\Ga\ell\mapsto
(e_i^-(\ell),h_i(\ell), e^+_i(\ell))_{i\in\ZZ}
$$ 
is a homeomorphism which conjugates the time-one discrete geodesic
flow $\flow{1}$ and the shift $\sigma$, that is, the following diagram
commutes
$$
\begin{CD}
\Ga\bs\G\XX' @>\flow{1}>>\Ga\bs\G \XX'\\
@V\Theta VV @VV\Theta V\\
\Sigma @>\sigma>> \;\,\Sigma\,,
\end{CD}
$$
and the topological Markov
shift $(\Sigma,\sigma)$ is locally compact and transitive.

Furthermore, if we endow $\Ga\bs\G \XX'$ with the quotient distance
of
$$
d(\ell,\ell')=
e^{-\sup\big\{n\in\NN\;:\;\;\ell|_{[-n,n]}\;=\;\ell'|_{[-n,n]}\big\}}
$$ 
on $\G\XX'$ and $\Sigma$ with the
distance
$$
d(x,x')= e^{-\sup\big\{n\in\NN\;:\;\;
\forall\, i\,\in\,\{-n,\dots,n\},\;\;x_i\;=\;x'_i\big\}}\;,
$$
then $\Theta$ is a bilipschitz homeomorphism.

Finally, if $\XX'$ is a uniform tree without vertices of degree at
most $2$, if the Gibbs measure $m_F$ of $\Ga$ is finite, and if the
length spectrum $L_\Ga$ of $\Ga$ is equal to $\ZZ$, then the
topological Markov shift $(\Sigma,\sigma)$ is topologically mixing.
\etheo

Note that when $\YY$ is finite (or equivalently when $\Ga$ is
cocompact), the alphabet $\A$ is finite (hence $(\Sigma,\sigma)$ is a
standard subshift of finite type). When furthermore the vertex groups
of $(\YY,G_*)$ are trivial (or equivalently when $\Ga$ acts freely,
and in particular is a finitely generated free group), this result is
well-known, but it is new if the vertex groups are not
trivial. Compare with the construction of \cite{CooPap02}, whose
techniques might be applied since $\Ga$ is word-hyperbolic if $\YY$ is
finite, up to replacing Gromov's (continuous time) geodesic flow of
$\Ga$ by the (discrete time) geodesic flow on $\G\XX'$, thus avoiding
the suspension part (see also the end of op.~cit.~when $\Ga$ is a
free group).

\medskip \dem 
For all $\ell\in\G \XX'$ and $\ga\in\Ga$, we can take $\ga_i(\ga\ell)
= \ga_i(\ell)\ga^{-1}$, and since $p(\ga f_i)=p(f_i)$, we have
$e_i^\pm(\ga\ell)=e_i^\pm(\ell)$ and $h_i(\ga \ell)= h_i(\ell)$, hence
the map $\Theta$ is well defined. By construction, the map $\Theta$ is
equivariant for the actions of $\flow{1}$ on $\Ga\bs\G \XX'$ and
$\sigma$.

%
With the distances indicated in the statement of Theorem
\ref{theo:codingdisgeodflo}, if $\ell,\ell'\in\G\XX'$ satisfy
$\ell_{\mid[-n,\,n]}\;=\;\ell'_{\mid[-n,\,n]}$ for some $n\in\NN$, then we
have $e_i^\pm(\ell)=e_i^\pm(\ell')$ for $-n \leq i \leq n-1$, and we
may take $\ga_i(\ell)=\ga_i(\ell')$ for $-n \leq i \leq n-1$, so that
$h_i(\ell)= h_i(\ell')$ for $-n \leq i \leq n-1$. Therefore, we have
$$
d(\,\Theta(\Ga\ell),\Theta(\Ga\ell')\,)\leq e\;d(\Ga\ell,\Ga\ell')\;,
$$
and $\Theta$ is Lipschitz (hence continuous).

\medskip 
Let us construct an inverse $\Psi:\Sigma \ra\Ga\bs\G \XX'$ of
$\Theta$, by a more general construction that will be useful later
on. Let $I$ be a nonempty interval of consecutive integers in $\ZZ$,
either finite or equal to $\ZZ$ (the definition of the inverse of
$\Theta$ only requires the second case $I=\ZZ$). For all $e^-,e^+\in
E\YY$ such that $t(e^-)=o(e^+)$, we fix once and for all a
representative of every double coset in $\rho_{e^-}(G_{e^-}) \bs
G_{o(e^+)} / \rho_{\,\overline{e^+}} (G_{e^+})$, and we will denote
this double coset by its representative.

Let $w= (e_i^-, h_i, e^+_i)_{i\in I}$ be a sequence indexed by $I$ in
the alphabet $\A$ such that for all $i\in I$ such that $i-1\in I$, we
have $e^+_{i-1}= e_i^-$ (when $I$ is finite, this means that
$w$ is an $A$-admissible sequence in $\A$, and when $I=\ZZ$, this
means that $w\in \Sigma$). In particular, the element $h_i\in
G_{o(e_i^+)}= \Ga_{\wt{o(e_i^+)}}$ is the chosen representative of its
double coset $\rho_{e_i^-}(G_{e_i^-}) \;h_i\;
\rho_{\,\overline{e_i^+}}(G_{e_i^+})$.

For every $i\in I$, note that
$$
o(\,h_i \;g_{\,\overline{e_i^+}}^{\;\;\;\;-1}\;
\wt{e_i^+}\;)=o(\,g_{\,\overline{e_i^+}}^{\;\;\;\;-1}\;
\wt{e_i^+}\;)=\wt{o(e_i^+)}=\wt{t(e_i^-)}= t(\,g_{e^-_{i}}^{\;\;-1}
\;\wt{e_{i}^-}\;)\;.
$$ 
But $h_i \;g_{\,\overline{e_i^+}}^{\;\;\;\;-1}\; \wt{e_i^+}$ is not
the opposite edge of the edge $g_{e^-_{i}}^{\;\;-1} \;\wt{e_{i}^-}$,
since the double coset of $h_i$ is not the trivial one $[1]$ when
$e^+_i= \overline{e^-_i}$, hence $h_i$ does not fix 
$g_{e^-_{i}}^{\;\;-1} \;\wt{e_{i}^-}$. Therefore the length $2$ edge path
(see the picture below)
$$
(\;g_{e^-_{i}}^{\;\;-1}\;\wt{e_{i}^-}, \;h_i \;
g_{\,\overline{e_i^+}}^{\;\;\;\;-1} \; \wt{e_i^+}\;)
$$ 
is geodesic.

\begin{center}
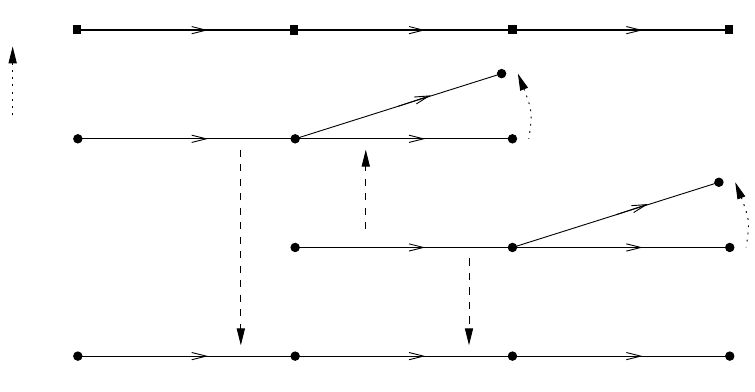
\end{center}

Let us construct by induction a geodesic segment $\wt w$ in $\XX'$
(which will be a discrete geodesic line if $I=\ZZ$), well defined up
to the action of $\Ga$, as follows.

We fix $i_0\in I$ (for instance $i_0=0$ if $I=\ZZ$ or $i_0=\min I$ if
$I$ is finite), and $\alpha_{i_0}\in \Ga$.
Let us define
$$
f_{i_0}=f_{i_0}(w)=\alpha_{i_0}\;g_{e^+_{i_0}}^{\;\;-1}\;\wt{e^+_{i_0}}\;.
$$
Let us then define 
$$
\alpha_{i_0-1} =\alpha_{i_0-1}(w)=
\alpha_{i_0}\;g_{e^+_{i_0}}^{\;\;\;-1}\;g_{\,\overline{e^+_{i_0}}}^{\,}
\;h_{i_0}^{\;-1} \;\;\;{\rm and}\;\;\; f_{i_0-1}=f_{i_0-1}(w)=
\alpha_{i_0-1} \;g_{e^-_{i_0}}^{\;\;\;\;\;-1}\;\wt{e_{i_0}^-}\;.
$$
We have $\alpha_{i_0-1}\;h_{i_0}\;g_{\,\overline{e^+_{i_0}}}^{\;\;\;-1}\;
\wt{e_{i_0}^+} =f_{i_0}$ and $(f_{i_0-1},f_{i_0})$ is a geodesic edge
path of length $2$ (as the image by $\alpha_{i_0-1}$ of such a path).

Let $i-1,i'\in I$ be such that $i'\leq i_0\leq i-1$. Assume by
increasing induction on $i$ and decreasing induction on $i'$ that a
geodesic edge path $(f_{i'-1}=f_{i'-1}(w),\dots, f_{i-1}=f_{i-1}(w))$
in $\XX'$ and a sequence $(\alpha_{i'-1}=\alpha_{i'-1}(w), \dots,
\alpha_{i-1}=\alpha_{i-1}(w))$ in $\Ga$ have been constructed such
that
$$
f_{j}=\alpha_{j}\;g_{e^+_{j}}^{\;\;\;-1}\;\wt{e_{j}^+} 
\;\;\;{\rm and}\;\;\; 
\alpha_{j}= 
\alpha_{j-1}\;h_{j}\;g_{\,\overline{e^+_{j}}}^{\;\;-1}\;g_{e^+_{j}}
$$
for every $j\in\NN$ such that $i'-1\leq j\leq i-1$, with besides
$j\geq i'$ for the equality on the right.

If $i$ does not belong to $I$, we stop the construction on the right
hand side at $i-1$.  If on the contrary $i\in I$, let us define (see
the above picture)
$$
\alpha_{i}= 
\alpha_{i-1}\;h_{i}\;g_{\,\overline{e^+_{i}}}^{\;\;\;\;-1}\;g_{e^+_{i}}
\;\;\;{\rm and}\;\;\;
f_{i}=f_{i}(w)=\alpha_{i}\;g_{e^+_{i}}^{\;\;\;-1} \;\wt{e_{i}^+}\;.
$$
Then 
$$
(f_{i-1},f_{i})=\big(\alpha_{i-1}\;g_{e^-_{i}}^{\;\;\;-1} \;\wt{e_{i}^-},\;
\alpha_{i-1}\;h_{i}\;
g_{\,\overline{e^+_{i}}}^{\;\;\;\;-1} \;\wt{e_{i}^+}\big)\;,
$$ 
is a geodesic edge path of length $2$ (as the image by
$\alpha_{i-1}$ of such a path). As an edge path is geodesic if and
only if it has no back-and-forth, $(f_{i'},\dots, f_{i})$ is a
geodesic edge path in $\XX'$. Thus the construction holds at rank $i$
on the right.

If $i'-1$ does not belong to $I$, we stop the construction on the left
side at $i'$. Otherwise we proceed as for the construction of
$\alpha_{i_0-1}$ and $f_{i_0-1}$ in order to construct $\alpha_{i'-2}$ and
$f_{i'-2}$ with the required properties.

If $I=[p,q]\cap\ZZ$ with $p\leq q$ in $\ZZ$, let $I'=[p-1,q]\cap\ZZ$. If
$I=\ZZ$, let $I'=\ZZ$. We have thus constructed a geodesic edge path
\begin{equation}\label{eq:canonicliftw}
(f_i)_{i\in I'}= (f_i(w))_{i\in I'}
\end{equation}
in $\XX'$. We denote by $\wt w$ its parametrisation
by $\RR$ if $I=\ZZ$ and by $[p-1,q+1]$ if $I=[p,q]\cap\ZZ$, in such a
way that $\wt w(i)=o(f_i)$ for all $i\in I$. In particular, $f_i=\wt
w([i,i+1])$ for all $i\in I'$. When $I=[p,q]\cap\ZZ$, we consider $\wt
w$ as a generalised discrete geodesic line, by extending it to a
constant on $]-\infty, p-1]$ and on $[q+1,+\infty[\,$.

The orbit $\Ga \wt w$ of $\wt w$ does not depend on the choice of
$\alpha_{i_0}$, since replacing $\alpha_{i_0}$ by $\alpha'_{i_0}$
replaces $f_i$ by $\alpha'_{i_0} \alpha_{i_0}^{-1} f_i$ for all $i\in
I'$, hence replaces $\wt w$ by $\alpha'_{i_0}\alpha_{i_0}^{-1}\wt
w$. This also implies that $\Ga \wt w$ does not depend on the choice
of $i_0\in I$.

\medskip Assume from now on that $I=\ZZ$, and define $\Psi:\Sigma
\ra\Ga\bs\G \XX'$ by
$$
\Psi(w)=\Ga \wt w\;.
$$ 
Let $w=(e_i^-, h_i, e^+_i)_{i\in \ZZ}$ and $w'=({e'_i}^-, {h'}_i,
{e'_i}^+)_{i\in I}$ in $\Sigma$ satisfy $e_i^\pm={e'_i}^\pm$ and $h_i=
{h'}_i$ for all $i\in \{-n,\dots, n\}$ for some $n\in\NN$. Then we may
take the same $i_0=0$ and $\alpha_{i_0}$ in the construction of $\wt
w$ and $\wt{w'}$. We thus have $\alpha_i(w)=\alpha_i(w')$ and $f_i(w)=
f_i(w')$ for $-n \leq i \leq n$. Therefore, with the distances
indicated in the statement of Theorem \ref{theo:codingdisgeodflo}, we
have
$$
d(\,\Psi(w),\Psi(w')\,)\leq d(w,w')\;,
$$
and $\Psi$ is Lipschitz.

Let us prove that $\Psi$ is indeed the inverse of $\Theta$. As in the
construction of $\Theta$, for all $\ell\in \G \XX'$ and $i\in\ZZ$, we
define $f_i=\ell([i,i+1])$, $e_i^+=p(f_{i})$ and $e_i^-=e_{i-1}^+$.
We denote by $\ga'_i\in\Ga$ an element sending $f_i$ to
$g_{e_i^+}^{\;\;-1}\; \wt{e_i^+}$ for all $i\in\ZZ$ (see the picture
below): with the notation above the statement of Theorem
\ref{theo:codingdisgeodflo}, we have
$\ga'_i=g_{e_i^+}^{\;\;-1}\;\ga_i(\ell)$.  \smallskip

\begin{center}
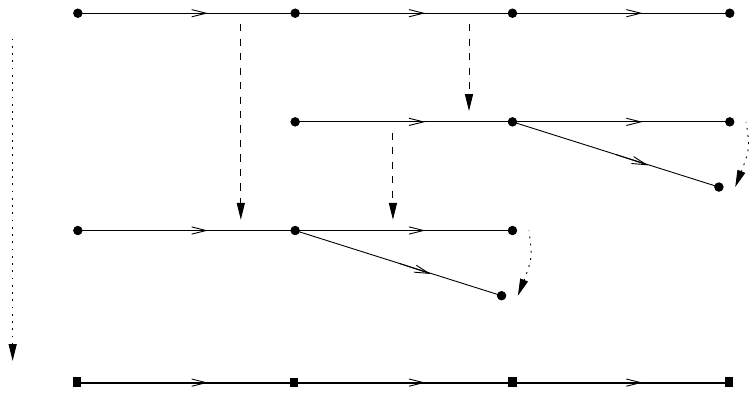
\end{center}

\smallskip
\noindent 
Then $\ga'_i$ is well defined up to multiplication on the left by an
element of $\Ga_{g_{e_i^+}^{\;\;-1} \;\wt{e_i^+}}=\rho_{e_i^+}(G_{e_i^+})$.  
Let $h'_i$ be an element in $G_{o(e_i^+)}$ sending
$g_{\overline{e_i^+}}^{\;\;\;-1}\; \wt{e_i^+}$ to $\ga'_{i-1}f_i$. It
exists since these two edges have the same origin $\wt{o(e_i^+)}$, and
same image by $p$:
$$
p(\ga'_{i-1}f_i)=p(f_i)=e_i^+=p(\wt{e_i^+})=
p(\,g_{\overline{e_i^+}}^{\;\;\;-1}\; \wt{e_i^+}\,)\;.
$$ 
Furthermore, it is well defined up to multiplication on the right by
an element of $\Ga_{g_{\overline{e_i^+}}^{\;\;\;-1} \;\wt{e_i^+}} =
\rho_{\overline{e_i^+}}(G_{e_i^+})$, and we have (see the above
picture)
$$
\ga'_{i-1}\,{\ga'_i}^{-1}\;g_{e_i^-}^{\;\;-1}\,g_{\,\overline{e_i^+}}
\;\in \;h'_i\;\rho_{\overline{e_i^+}}(G_{e_i^+})
$$

Using $\ga_j(\ell)= g_{e_j^+}\ga'_j$ for $j=i,i-1$ in Equation
\eqref{eq:defigi} gives $h_i(\ell)= \ga'_{i-1}{\ga'_i}^{-1}
g_{e_i^+}^{\;-1} g_{\overline{e_i^+}}$. Hence by the construction of
$\Theta$ (see with $\ga'_i=g_{e_i^-}^{\;\;-1}\;\ga_i(\ell)$ for all
$i\in\ZZ$), we have
$$
\Theta(\ell)=(e_i^-,\;\rho_{e_i^-}(G_{e_i^-})\; h'_i\;
\rho_{\overline{e_i^+}}(G_{e_i^+}),\; e^+_i)_{i\in\ZZ}\;.
$$
Let $h_i$ be the chosen representative of the double coset
$\rho_{e_i^-}(G_{e_i^-})\; h'_i\; \rho_{\overline{e_i^+}}
(G_{e_i^+})$~: there exist $\alpha\in \rho_{e_i^-}(G_{e_i^-})
=\rho_{e_{i-1}^+}(G_{e_{i-1}^+})$ and $\beta\in
\rho_{\overline{e_i^+}}(G_{e_i^+})$ such that $h_i=\alpha h'_i
\beta$. Up to replacing $\ga'_{i-1}$ by $\alpha^{-1}\ga'_{i-1}$ and $h'_i$ by
$h'_i\beta$, we then may have $h'_i=h_i$. By taking
$\alpha_{i_0}={\ga'_{i_0}}^{-1}$, we have $\alpha_{i}={\ga'_i}^{-1}$
for all $i\in\ZZ$, and an inspection of the above two constructions
gives that $\Theta\circ\Psi=\id$ and $\Psi\circ\Theta=\id$.

Since the discrete time geodesic flow is topologically transitive on
its nonwandering subset and by conjugation, the topological Markov
shift $(\Sigma,\sigma)$ is topologically transitive.

If $\XX'$ is a uniform tree without vertices of degree at most $2$,
if the length spectrum of $\Ga$ is equal to $\ZZ$ and if the Gibbs
measure $m_F$ for $\Ga$ is finite, then by Theorem
\ref{theo:uniflatmBMmixing}, the discrete time geodesic flow on
$\Ga\bs \G\XX'$ is topologically mixing, hence by conjugation by
$\Theta$, the topological Markov shift $(\Sigma,\sigma)$ is
topologically mixing. This concludes the proof of Theorem
\ref{theo:codingdisgeodflo}.  
\cqfd

\bigskip 
When the length spectrum $L_\Ga$ of $\Ga$ is different from $\ZZ$,
the topological Markov shift $(\Sigma,\sigma)$ constructed above is
not always topologically mixing. We now modify the above construction
in order to take care of this problem.

Recall that $\XX'=\C\Lambda\Ga$ and that $\Geven\XX'$ is the space of
geodesic lines $\ell\in\G\XX'$ whose origin $\ell(0)$ is at even
distance from the basepoint $x_0$ (we assume that $x_0\in \XX'$),
which is invariant under the time-two discrete geodesic flow
$\flow{2}$ and, when $L_\Ga=2\ZZ$, under $\Ga$, as seen in the proof
of Theorem \ref{theo:uniflatmBMmixing}.

Consider $\A_{\rm even}$ the alphabet consisting of the quintuples
$(f^-,h^-,f^0,h^+,f^+)$ where the triples $(f^-,h^-,f^0)$ and
$(f^0,h^+,f^+)$ belong to $\A$ and $o(f^0)$ is at even distance from
the image in $\YY=\Ga\bs\XX'$ of the basepoint $x_0$. Let $A_{\rm
  even}=(A_{{\rm even},\,i,j})_{i,j\in\A_{\rm even}}$ be the
transition matrix with line and column indices in $\A_{\rm even}$ such
that for all $i= (f^-,h^-,f^0, h^+,f^+)$ and
$j=(f^-_*,h^-_*,f^0_*,h^+_*,f^+_*)$, we have $A_{{\rm even},\,i,j}=1$
if and only if $f^+=f^-_*$. We denote by $(\Sigma_{\rm
  even},\sigma_{\rm even})$ the associated topological Markov
shift. We endow $\Sigma_{\rm even}$ with the slightly modified
distance
$$
d_{\rm even}(x,x')= e^{-2\sup\big\{n\in\NN\;:\;\;
\forall\, k\,\in\,\{-n,\,\dots,\,n\},\;\;x_k\;=\;x'_k\big\}}\;,
$$ 
where $x=(x_k)_{k\in\ZZ}$ and $x'=(x'_k)_{k\in\ZZ}$ are in
$\Sigma_{\rm even}$.

We have a canonical injection $\inj:\Sigma_{\rm even}\ra \Sigma$
sending the sequence $(f^-_n,h^-_n,f^0_n,h^+_n,f^+_n)_{n\in\ZZ}$ to 
$(e_n^-, h_n, e^+_n)_{n\in\ZZ}$ with, for every $n\in\ZZ$,
$$
e_{2n}^-=f^-_n,\;h_{2n}=h^-_n,\;e_{2n}^+=f^0_n,\;
e_{2n+1}^-=f^0_n,\;h_{2n+1}=h^+_n,\;e_{2n+1}^+=f^+_n,\;.
$$ 
By construction, $\inj$ is clearly a homeomorphism onto its image,
and
$$
\Theta(\Ga\bs\Geven\XX')=\inj(\Sigma_{\rm  even})\;.
$$ 
If two sequences in $\Sigma_{\rm even}$ coincide between $-n$ and $n$,
then their images by $\inj$ coincide between $-2n$ and
$2n$. Conversely, if the images by $\inj$ of two sequences in
$\Sigma_{\rm even}$ coincide between $-2n-1$ and $2n+1$, then these
sequences coincide between $-n$ and $n$. Hence $\inj$ is bilipschitz,
for the above distances.

Let us define $\Theta_{\rm even}=\inj^{-1}\circ\,
\Theta|_{\Ga\bs\Geven\XX'} :\Ga\bs\Geven\XX' \ra\Sigma_{\rm
  even}$. The following diagram hence commutes
$$
\begin{CD}
\Ga\bs\Geven \XX' @>\Theta_{\rm even}>> \Sigma_{\rm even}\\
@VVV @VV\inj V\\
\Ga\bs\G \XX' @>\Theta>> \;\,\Sigma\,,
\end{CD}
$$
where the vertical map on the left hand side is the inclusion map.

\btheo\label{theo:codingeven} Assume that $\XX'=\C\Lambda\Ga$ is a
uniform tree without vertices of degree at most $2$, that the Gibbs
measure $m_{F}$ of $\Ga$ is finite, and that the length spectrum
$L_\Ga$ of $\Ga$ is equal to $2\ZZ$.  Then the map $\Theta_{\rm even}
: \Ga\bs\Geven\XX' \ra\Sigma_{\rm even}$ is a bilipschitz homeomorphism
which conjugates the time-two discrete geodesic flow $\flow{2}$ and
the shift $\sigma_{\rm even}$, and the topological Markov shift
$(\Sigma_{\rm even},\sigma_{\rm even})$ is locally compact and
topologically mixing.  
\etheo

\dem The only claims that remains to be proven is the last one, which
follows from Theorem \ref{theo:uniflatmBMmixing}, by
conjugation. \cqfd

\bigskip 
Let us now study the properties of the image by the coding map
$\Theta$ of finite Gibbs measures on $\Ga\bs \G\XX$.

Let $\delta=\delta_{\Ga,F^\pm}$ be the critical exponent of
$(\Ga,F^\pm)$. Let $(\mu^\pm_x)_{x\in V\XX}$ be two (normalised)
Patterson densities on $\partial_\infty X$ for the pairs
$(\Ga,F^\pm)$, where as previously $\wt F^+=\wt F$, $\wt F^-= \wt
F\circ\iota$. Assume that the associated Gibbs measure $m_F$ on
$\Ga\bs\G \XX$ (using the convention for discrete time of Section
\ref{subsec:pattersongibbstrees}) is finite.

Let us define
\begin{equation}\label{eq:defPP}
\PP=\frac{1}{\|m_F\|}\;\Theta_*m_F
\end{equation}
as the image of the Gibbs measure $m_F$ (whose support is $\Ga\bs\G
\XX'$) by the homeomorphism $\Theta$, normalised to be a
probability measure. It is a probability measure on $\Sigma$,
invariant under the shift $\sigma$. 

Let $(Z_n)_{n\in\ZZ}$ be the random process classically associated with
the full shift $\sigma$ on $\Sigma$: it is the random process on the
Borel space $\Sigma$ indexed by $\ZZ$ with values in the discrete
alphabet $\A$, where $Z_n:\Sigma\ra\A$ is the (continuous hence
measurable) $n$-th projection $(x_k)_{k\in\NN}\mapsto x_n$ for all
$n\in\ZZ$.

\medskip 
The following Proposition \ref{prop:propriPP} summarises the
properties of the probability measure $\PP$. We start by recalling and
giving some notation used in this proposition.

For every admissible finite sequence $w=(a_p,\dots, a_q)$ in
$\A$, where $p\leq q$ in $\ZZ$, we denote

\smallskip
$\bullet$~ by $[w]=[a_p,\dots, a_q]=\{(x_n)_{n\in\ZZ}\in\Sigma\;:\;
\forall \;n\in \{p,\dots, q\},\;\; x_n=a_n\big\}$ the associated
cylinder in $\Sigma$,

\smallskip
$\bullet$~ by $\wh w$ the associated geodesic edge path in $\XX'$ with
length $q-p+2$ constructed in the proof of Theorem
\ref{theo:codingdisgeodflo} (see Equation \eqref{eq:canonicliftw}),
with origin $\wh w_-$ and endpoint $\wh w_+$.

\medskip
For every geodesic edge path $\alpha=(f_{p-1},\dots, f_q)$ in $\XX'$,
we define (See Section \ref{subsec:trees} for the notation, and the
picture below)
$$
\partial^+_\alpha\XX'=\partial_{f_q}\XX'\;\;\;{\rm and}\;\;\;
\partial^-_\alpha\XX'=\partial_{\,\overline{f_{p-1}}}\XX'\;,
$$ 
and 
$$
\G_\alpha\XX=\{\ell\in\G\XX\;:\; \ell(p-1)=o(f_{p-1})
\;\;\;{\rm and}\;\;\; \ell(q+1)=t(f_{q})\}\;.
$$
\begin{center}
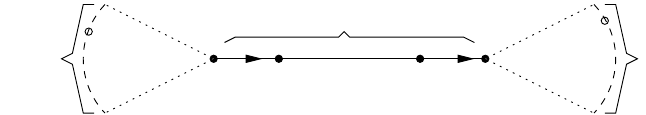
\end{center}

We define a map $F_{\rm symb}:\Sigma\ra \RR$ by
\begin{equation}\label{eq:defipotentsymbol}
F_{\rm symb}(x)= \int_{o(e^+_0)}^{t(e^+_0)} F
\end{equation}
if $x=(x_i)_{i\in\ZZ}$ with $x_0=(e_0^-,h_0,e^+_0)$.  Note that for
all $(x_n)_{n\in\ZZ}, (y_n)_{n\in\ZZ} \in\Sigma$, if $x_0=y_0$, then
$F_{\rm symb}(x)=F_{\rm symb} (y)$, so that $F_{\rm symb}$ is locally
constant (constant on each cylinder of length $1$ at time $0$), hence
continuous.

For instance, if $F=F_c$ is the potential associated with a system of
conductances $c:\Ga\bs E\XX'\ra\RR$ (see Section \ref{subsec:cond}),
then
$$ 
F_{\rm symb}(x)= c(e^+_0)\;.
$$ 
Note that if $c,c':\Ga\bs E\XX'\ra\RR$ are cohomologous systems of
conductances on $\Ga\bs E\XX'$, then the corresponding maps $F_{\rm
  symb}, F'_{\rm symb}:\Sigma\ra \RR$ are cohomologous. Indeed if $f:
\Ga\bs V\XX\ra \RR$ is a map such that $c'(e)-c(e)=f(t(e))-f(o(e))$
for every $e\in \Ga\bs E\XX$, with $G:\Sigma\ra \RR$ the map defined by
$G(x)=f(o(e^+_0))$ if $x=(x_i)_{i\in\ZZ}$ with $x_0=(e_0^-,h_0,e^+_0)$, 
then $G$ is locally constant, hence continuous, and since
$t(e^+_0)=o(e^+_1)$, we have, for every $x\in\Sigma$,
$$
F'_{\rm symb}(x)-F_{\rm  symb}(x)=G(\sigma x)-G(x)\;.
$$

\bdefi\label{def:Markovgood} {\rm Let $\XX''$ be a locally finite
  simplicial tree. A nonelementary discrete subgroup $\Ga'$ of
  $\Aut(\XX'')$ is {\it Markov-good}\index{Markov-good} 
if for every $n\in\NN-\{0\}$ and every geodesic edge path $(e_0,\dots,
e_{n+1})$ in $\C\Lambda\Ga'$, we have
\begin{equation}\label{eq:Markov-good}
|\Ga'_{e_0}\cap \dots \cap \Ga'_{e_n}|\; 
|\Ga'_{e_{n-1}}\cap \Ga'_{e_n}\cap \Ga'_{e_{n+1}}|
=|\Ga'_{e_0}\cap \dots \cap \Ga'_{e_{n+1}}|\;
|\Ga'_{e_{n-1}}\cap \Ga'_{e_n}|\;.
\end{equation}}
\edefi

\brema\label{rem:notMarkovgood}
(1) Note that Equation \eqref{eq:Markov-good} is automatically
satisfied if $n=1$ and that $\Ga'$ is Markov-good if $\Ga'$ acts
freely on $\XX''$.

\medskip
\noindent (2) A group action on a simplicial tree is {\it
  $2$-acylindrical}~\index{acylindrical}\footnote{See for instance
  \cite{Sela97,GuiLev11}, which require other minor hypotheses that
  are not relevant here.} if the stabiliser of any geodesic edge path
of length $2$ is trivial.
If $\Ga'$ is $2$-acylindrical on $\XX$, then $\Ga'$ is Markov-good,
since all groups appearing in Equation \eqref{eq:Markov-good} are
trivial.

\medskip
\noindent(3) If $\XX''$ has degrees at least $3$ and if $\Ga'$ is a
noncocompact geometrically finite lattice of $\XX''$ with abelian edge
stabilisers, then $\Ga'$ is not Markov-good.  \erema

\dem (3) Since the quotient graph $\Ga'\bs \XX''$ is infinite, the
graph of groups $\Ga'\dbs\XX''$ contains at least one cuspidal
ray. Consider a geodesic ray in $\XX''$ with consecutive edges
$(f_n)_{n\in\NN}$ mapping injectively onto this cuspiday ray, pointing
towards its end. The stabilisers of the edges $f_n$ in $\Ga'$ are
hence nondecreasing in $n$: we have $\Ga'_{f_n} \subset
\Ga'_{f_{n+1}}$ for all $n\in\NN$. By the finiteness of the volume,
there exists $n\geq 3$ such that $\Ga'_{f_{n-2}}$ is strictly
contained in $\Ga'_{f_{n-1}}$. Since $\XX''$ has degrees at least $3$,
there exists $\ga\in\Ga'$ fixing $t(f_{n-1})$ but not fixing
$f_{n-1}$. Let $e_0=f_0,\dots, e_{n-1}= f_{n-1}$, $e_{n}=
\ga\;\overline{f_{n-1}}$ and $e_{n+1}=\ga\; \overline{f_{n-2}}$.
Since the stabilisers of $f_{n-1}$ and $e_n$ are conjugated by $\ga$
within the abelian stabiliser of $f_n$, they are equal. Then $(e_0,
\dots, e_{n+1})$ is a geodesic edge path in the simpliciak tree
$\XX''$ (whose geometric realisation is equal to $\C\Lambda\Ga'$ since
$\Ga'$ is a lattice). Since $\Ga'_{e_0}\cap \dots \cap
\Ga'_{e_n}=\Ga'_{f_0}$, $\Ga'_{e_{n-1}}\cap \Ga'_{e_n}\cap
\Ga'_{e_{n+1}} =\Ga'_{f_{n-2}}$, $\Ga'_{e_0}\cap \dots \cap
\Ga'_{e_{n+1}} = \Ga'_{f_0}$, $\Ga'_{e_{n-1}}\cap \Ga'_{e_n}=
\Ga'_{f_{n-1}}$ and $|\Ga'_{f_{n-2}}|\neq |\Ga'_{f_{n-1}}|$, the
subgroup $\Ga'$ is not Markov-good.  
\cqfd

\bigskip
Recall that a random process $(Z'_n)_{n\in\ZZ}$ on $(\Sigma,\PP)$ is a
{\it Markov chain}\index{Markov!chain} if and only if for all $p\leq
q$ in $\ZZ$ and $a_p,\dots, a_q,a_{q+1}$ in $\A$, we have, when defined,
\begin{equation}\label{eq:defirandchainZ}
\PP(Z'_{q+1}=a_{q+1}\,|\,Z'_{q}=a_{q},\dots, Z'_{p}=a_{p})=
\PP(Z'_{q+1}=a_{q+1}\,|\,Z'_{q}=a_{q})\;.
\end{equation}

\bprop\label{prop:propriPP} (1) For every admissible finite sequence
$w$ in $\A$, we have
$$
\PP([w])=\frac{\mu^-_{\wh w_-}(\partial^-_{\wh w}\XX')\;
\mu^+_{\wh w_+}(\partial^+_{\wh w}\XX')\;e^{\int_{\wh w_-}^{\wh w_+}(\wt F-\delta)}}
{|\Ga_{\wh w}|\;\|m_F\|}\;.
$$

\medskip\noindent (2) The random process $(Z_n)_{n\in\ZZ}$ on
$(\Sigma,\PP)$ is a Markov chain if and only if $\Ga$ is Markov-good.

\medskip\noindent (3) The measure $\PP$ on the topological Markov
shift $\Sigma$ satisfies the Gibbs property with Gibbs constant
$\delta$ for the potential $F_{\rm symb}$.  
\eprop

It follows from the above Assertion (2) and from Remark
\ref{rem:notMarkovgood} that when $\XX$ has degrees at least $3$ and
$\Ga$ is a noncocompact geometrically finite lattice of $\XX$ with
abelian edge stabilisers (and more generally, this is not a necessary
assumption), then $(Z_n)_{n\in\ZZ}$ is not a Markov chain. The fact
that codings of discrete time geodesic flows on trees might not
satisfy the Markov chain property had been noticed by Burger and Mozes
around the time the paper \cite{BurMoz96} was
published.\footnote{Personal communication.} When proving the
variational principle in Section \ref{subsec:varprinc} and the
exponential decay of correlations in Section
\ref{subsec:mixingratesimpgraphs}, we will hence have to use tools
that are not using the Markov chain property.

\medskip
\dem (1) Let $w=(a_p,\dots, a_q)$, with $p\leq q$ in $\ZZ$, be an
admissible finite sequence in $\A$. Recall that $[w]=\{x\in\Sigma
\;:\;\forall\, i\in\{p,\dots,q\},\;\; x_i=a_i\}$. By the construction of
$\Theta$, the preimage $\Theta^{-1}([w])$ is equal to the image
$\Ga\G_{\wh w} \XX'$ of $\G_{\wh w}\XX'$ in $\Ga\bs\G\XX'$.  Hence,
since $\Ga_{\wh w}$ is the stabiliser of $\G_{\wh w}\XX'$ in $\Ga$,
$$
\PP([w])=\frac{1}{\|m_F\|}\;m_F(\Ga\G_{\wh w}\XX')=
\frac{1}{|\Ga_{\wh w}|\;\|m_F\|}\;\wt m_F(\G_{\wh w}\XX')\;.
$$ 
In the expression of $\wt m_F$ given by Equation
\eqref{eq:defigibbsdis}, let us use as basepoint  the origin $\wh
w_-$ of the edge path $\wh w$, and note that all elements of $\G_{\wh
  w} \XX'$ pass through $\wh w_-$ at time $t=p-1$, so that by the
invariance of $\wt m_F$ under the discrete time geodesic flow, we have
\begin{align*}
\wt m_F(\G_{\wh w}\XX') &=\int_{\ell\in \G_{\wh w}\XX'} d\wt m_F(\flow{1-p}\ell)
=\int_{\ell_-\in\partial^-_{\wh w}\XX'}\int_{\ell_+\in\partial^+_{\wh w}\XX'}
d\mu^-_{\wh w_-}(\ell_-)d\mu^+_{\wh w_-}(\ell_+)\\ &
=\mu^-_{\wh w_-}(\partial^-_{\wh w}\XX')\;
\mu^+_{\wh w_-}(\partial^+_{\wh w}\XX')
=\mu^-_{\wh w_-}(\partial^-_{\wh w}\XX')\;
\mu^+_{\wh w_+}(\partial^+_{\wh w}\XX')\;e^{\int_{\wh w_-}^{\wh w_+}(\wt F-\delta)}\;,
\end{align*}
where this last equality follows by Equations
\eqref{eq:quasinivarPatdens} and \eqref{eq:changemoinsplus} with
$x=\wh w_-$ and $y=\wh w_+$, since for every $\ell_+ \in
\partial^+_{\wh w}\XX'$, we have $\wh w_+\in[\wh w_+, \ell_+[\,$.

\medskip\noindent (2) Let us fix $p\leq q$ in $\ZZ$ and $a_p,\dots,
a_q,a_{q+1}$ in $\A$, and let us try to verify Equation
\eqref{eq:defirandchainZ} for $(Z'_n)_{n\in\NN}=(Z_n)_{n\in\NN}$. Let
$\alpha_*=(a_p,\dots, a_q)$, which is an admissible sequence, since we
assumed the conditional probability $\PP(Z_{q+1}=a_{q+1}\,|\,
Z_{q}=a_{q}, \dots, Z_{p}=a_{p})$ to be well defined.  We may assume
that $\alpha=(a_p,\dots, a_q,a_{q+1})$ is an admissible sequence,
otherwise both sides of Equation \eqref{eq:defirandchainZ} are $0$.
Let us consider
$$
Q_\alpha=\frac{\PP(Z_{q+1}=a_{q+1}\,|\,Z_{q}=a_{q},\dots, Z_{p}=a_{p})}
{\PP(Z_{q+1}=a_{q+1}\,|\,Z_{q}=a_{q})}=
\frac{\PP([a_p,\dots , a_{q+1}])\;\PP([a_q])}
{\PP([a_p,\dots , a_{q}])\;\PP([a_q , a_{q+1}])}\;.
$$

Let us replace each one of the four terms in this ratio by its value
given by Assertion (1). Since $\partial^-_{\wh \alpha}\XX'=
\partial^-_{\wh {\alpha_*}}\XX'$, $\partial^+_{\wh \alpha}\XX'=
\partial^+_{\widehat{a_q,a_{q+1}}}\XX'$, $\partial^+_{\wh {\alpha_*}}
\XX' = \partial^+_{\widehat{a_q}}\XX'$ and $\partial^-_{\wh{a_q}}\XX'=
\partial^-_{\widehat{a_q,a_{q+1}}}\XX'$, all Patterson measure terms
cancel. Denoting by $y_1$ the common origin of $\wh \alpha$ and
$\wh{\alpha_*}$, by $y_2$ the common origin of $\wh a_q$ and
$\widehat{a_q,a_{q+1}}$, by $y_3$ the common terminal point of $\wh
a_q$ and $\wh {\alpha_*}$, and by $y_4$ the common terminal point of
$\widehat{a_q,a_{q+1}}$ and $\wh \alpha$, we thus have by Assertion
(1)
$$
Q_\alpha=\frac{|\Ga_{\wh {\alpha_*}}|\;|\Ga_{\widehat{a_q,a_{q+1}}}|}
{|\Ga_{\wh \alpha}|\;|\Ga_{\widehat{a_q}}|} \;\;
\frac{e^{\int_{y_1}^{y_4}(\wt F-\delta)}\;e^{\int_{y_2}^{y_3}(\wt F-\delta)}}
{e^{\int_{y_1}^{y_3}(\wt F-\delta)}\;e^{\int_{y_2}^{y_4}(\wt F-\delta)}}\;.
$$
Since $y_1,y_2,y_3,y_4$ are in this order on $[y_1,y_4]$, we have
$$
Q_\alpha = \frac{|\Ga_{\wh {\alpha_*}}|\;|\Ga_{\widehat{a_q,a_{q+1}}}|}
{|\Ga_{\wh \alpha}|\;|\Ga_{\widehat{a_q}}|} \;.
$$

Since every geodesic edge path of length $n+1$ at least $3$ in $\XX'$
defines an admissible sequence of length $n$ at least $2$ in $\A$, by
Equation \eqref{eq:Markov-good}, we have $Q_\alpha=1$ for every
admissible sequence $\alpha$ in $\A$ if and only if $\Ga$ is
Markov-good.

\medskip\noindent (3) Let $E$ be a finite subset of the alphabet $\A$,
and let $w=(a_p,\dots, a_q)$ with $p\leq q$ in $\ZZ$ be an admissible
sequence in $\A$ such that $a_p,a_q\in E$. By Assertion (1), we have
$$
\PP([w])=\frac{\mu^-_{\wh w_-}(\partial^-_{\wh w}\XX')\;
\mu^+_{\wh w_+}(\partial^+_{\wh w}\XX')\;e^{\int_{\wh w_-}^{\wh w_+}(\wt F-\delta)}}
{|\Ga_{\wh w}|\;\|m_F\|}\;.
$$ 
Since $a_p,a_q$ are varying in the finite subset $E$ of $\A$, the
first and last edges of $\wh w$ vary amongst the images under elements
of $\Ga$ of finitely many edges of $\XX'$. Since $w$ is admissible, the
sets $\partial^\pm_{\wh w}\XX'$ are nonempty open subsets of $\Lambda
\Ga$, hence they have positive Patterson measures. Furthermore, the
quantities $\mu^\pm_{\wh w_\pm} (\partial^\pm_{\wh w}\XX')$ are invariant 
under the action of $\Ga$ on the first/last edge of $\wh w$.  Hence
there exists $c_1\geq 1$ depending only on $E$ such that $1 \leq
|\Ga_{\wh w}| \leq |\Ga_{\wh w_-}| \leq c_1$ and $\frac1{c_1}\leq
\mu^\pm_{\wh w_\pm} (\partial^\pm_{\wh w}\XX') \leq c_1$.

Note that the length of $\wh w$ is equal to $q-p+2$. Therefore
$$
\frac{e^{-\delta}}{c_1^3\;\|m_F\|}\;e^{-\delta(q-p+1)+\int_{\wh w_-}^{\wh w_+} \wt F} 
\leq \; \PP([w])\; \leq 
\frac{e^{-\delta}\;c_1^2}{\|m_F\|}\;e^{-\delta(q-p+1)+\int_{\wh w_-}^{\wh w_+} \wt F}\;.
$$

If $\wh w= (f_{p-1},f_p,\dots, f_{q})$ and $x\in [w]$, we have by the
definition of $F_{\rm symb}$
$$
\int_{\wh w_-}^{\wh w_+} \wt F=\sum_{i=p-1}^q \int_{o(f_i)}^{t(f_i)} \wt F=
\int_{o(f_{p-1})}^{t(f_{p-1})} \wt F+ \sum_{i=p}^q F_{\rm symb}(\sigma^i(x))\;.
$$ 

Since $\wt F$ is continuous and $\Ga$-invariant, and since
$o(f_{p-1})$ remains in the image under $\Ga$ of a finite subset of
$V\XX'$, there exists $c_2>0$ depending only on $E$ such that $|\wt
F(v)|\leq c_2$ for every $v\in T^1X$ with $\pi(v)\in [o(f_{p-1}),
  t(f_{p-1})]$. Hence $|\int_{o(f_{p-1})}^{t(f_{p-1})} \wt F|\leq
c_2$, and Assertion (3) of Proposition \ref{prop:propriPP} follows
(see Equation \eqref{eq:defiGibbspropertyZ} for the definition of the
Gibbs property).  
\cqfd

\bigskip
Again in order to consider the case when the length spectrum $L_\Ga$
of $\Ga$ is $2\ZZ$, we define
$$
\PP_{\rm even}=\frac{1}{\|(m_F)_{\mid\;  \Ga\bs\Geven\XX'}\|}\;
(\Theta_{\rm even})_*\big((m_F)_{\mid\; \Ga\bs\Geven\XX'}\big)\;,
$$ 
and $(Z_{{\rm even},\,n})_{n\in\ZZ}$ the random process associated
with the full shift $\sigma_{\rm even}$ on $\Sigma_{\rm even}$, with
$Z_{{\rm even},\,n}: \Sigma_{\rm even}\ra\A_{\rm even}$ the $n$-th
projection for every $n\in\ZZ$.

By a proof similar to the one of Proposition \ref{prop:propriPP}, we
have the following result. We define a map $F_{\rm symb,~even}:
\Sigma_{\rm even}\ra \RR$ by
\begin{equation}\label{eq:defipotentsymboleven}
F_{\rm symb,~even}(x)= \int_{o(f^0_0)}^{t(f^+_0)} F
\end{equation}
if $x=(x_i)_{i\in\ZZ}$ with $x_0=(f_0^-,h_0^-,f_0^0,h_0^+,f_0^+)$.  As
previously, $F_{\rm symb,~even}$ is locally constant, hence
continuous.

\bprop\label{prop:propriPPeven} The measure $\PP_{\rm even}$ on the
topological Markov shift $\Sigma_{\rm even}$ satisfies the Gibbs
property with Gibbs constant $\delta$ for the potential $F_{\rm
  symb,~even}$.  \cqfd 
\eprop

Again, if $\Ga$ is a noncocompact geometrically finite lattice of
$\XX$ with abelian edge stabilisers and $\XX'$ has degrees at least
$3$, then $(Z_{{\rm even} , \,n})_{n\in\ZZ}$ is not a Markov chain.

\section{Coding continuous time geodesic flows on metric trees}
\label{subsec:codagemetric}

Let $(\XX,\lambda)$ be a locally finite metric tree without terminal
vertices, with $X=|\XX|_\lambda$ its geometric realisation. Let $\Ga$
be a nonelementary discrete subgroup of $\Aut(\XX,\lambda)$, and let
$\wt F:T^1X\ra \RR$ be a potential for $\Ga$. Let $X'=\C\Lambda\Ga$,
which is the geometric realisation $|\XX'|_\lambda$ of a metric
subtree $(\XX', \lambda)$. Let $\delta=\delta_{\Ga,F^\pm}$ be the
critical exponent of $(\Ga,F^\pm)$.  Let $(\mu^\pm_x)_{x\in V\XX}$ be
the (normalised) Patterson densities on $\partial_\infty X$ for the
pairs $(\Ga,F^\pm)$, and assume that the associated Gibbs measure
$m_F$ is finite.  We also assume in this Section that the lengths of
the edges of $(\XX',\lambda)$ have a finite upper bound (which is in
particular the case if $(\XX',\lambda)$ is uniform). They have a
positive lower bound by definition (see Section \ref{subsec:trees}).

In this Section, we prove that the continuous time geodesic flow on
$\Ga\bs\G X'$ is isomorphic to a suspension of a transitive
(two-sided) topological Markov shift on a countable alphabet, by an
explicit construction that will be useful later on in order to study
the variational principle (see Section \ref{subsec:varprinc}) and
rates of mixing (see Section \ref{subsec:mixingratemetgraphs}). Since
we are only interested in the support of the Gibbs measures, we will
only give such a description for the geodesic flow on the nonwandering
subset $\Ga\bs \G X'$ of $\Ga\bs\G X$. The same construction works
with the full space $\Ga\bs \G X$, but the resulting Markov shift is
then not necessarily transitive.

\medskip 
We start by recalling (see for instance \cite[\S 1.11]{BriStu02}) the
definitions of the suspension of an invertible discrete time dynamical
system and of the first return map on a cross-section of a continuous
time dynamical system, which allow to pass from transformations to
flows and back, respectively.

Let $(Z,\mu, T)$ be a metric space $Z$ endowed with a homeomorphism
$T$ and a $T$-invariant (positive Borel) measure $\mu$. Let $r:Z\ra
\;]0,+\infty[$ be a continuous map, such that for all $z\in Z$, the
subset $\{r(T^nz)\;:\;n\in\NN\} \cup \{-\,r(T^{-(n+1)}z) \;:\;
n\in\NN\}$ is discrete in $\RR$. Then the {\it
  suspension}\index{suspension} (or also {\it special
  flow}\index{special flow}) {\it over $(Z,\mu, T)$ with roof
  function\index{roof function} $r$} is the following continuous time
dynamical system $(Z_r,\mu_r, (T^t_r)_{t\in\RR})$~:

\smallskip
$\bullet$~ The space $Z_r$ is the quotient topological space
$(Z\times\RR)/\!\sim$ where $\sim$ is the equivalence relation on
$Z\times\RR$ generated by $(z,s+r(z))\sim (Tz,s)$ for all $(z,s)\in
Z\times\RR$. We denote by $[z,s]$ the equivalence class of
$(z,s)$. Note that 
$$
\F=\{(z,s)\;:\; z\in Z, \;\;0\leq s <r(z)\}
$$ 
is a measurable strict fundamental domain for this equivalence
relation. We endow $Z_r$ with the {\it Bowen-Walters
  distance}\index{Bowen-Walters
  distance}\index{distance!Bowen-Walters}, see \cite{BowWal72} and
particularly the appendix in \cite{BarSau00}.

\smallskip
$\bullet$~ For every $t\in\RR$, the map $T^t_r:Z_r\ra Z_r$ is the map
$[z,s]\mapsto [z,s+t]$. Equivalently, when $(z,s)\in\F$ and
  $t\geq 0$, then $T^t_r([z,s])=[T^nz,s']$ where $n\in\NN$ and
  $s'\in\RR$ are such that 
$$
t+s=\sum_{i=0}^{n-1}r(T^iz) + s'\;\;\;{\rm and}\;\;\;
 0\leq s'< r(T^nz)\;.
$$

\smallskip
$\bullet$~ Denoting by $ds$ the Lebesgue measure on $\RR$, the measure
$\mu_r$ is the pushforward of the restriction to $\F$ of the product
measure $d\mu\,ds$ by the restriction to $\F$ of the canonical
projection $(Z\times\RR)\ra Z_r$.

\medskip
Note that $(T^t_r)_{t\in\RR}$ is indeed a continuous one-parameter
group of homeomorphisms of $Z_r$, preserving the measure $\mu_r$. The
measure $\mu_r$ is finite if and only if $\int_Zr\,d\mu$ is finite,
since
$$
\|\mu_r\|\;=\;\int_\F \;d\mu ds=\int_Zr\,d\mu\;.
$$ 
We will denote by $(Z,\mu, T)_r$ the continuous time dynamical
system $(Z_r,\mu_r, (T^t_r)_{t\in\RR})$ thus constructed.

\medskip
Conversely, let $(Z,\mu, (\phi_t)_{t\in\RR})$ be a metric space $Z$
endowed with a continuous one-parameter group of homeomorphisms
$(\phi_t)_{t\in\RR}$, preserving a (positive Borel) measure
$\mu$. Let $Y$ be a {\it cross-section}\index{cross-section} of
$(\phi_t)_{t\in\RR}$, that is a closed subspace of $Z$ such that for
every $z\in Z$, the set $\{t\in\RR\;:\; \phi_t(z)\in Y\}$ is infinite
and discrete. Let $\tau:Y\ra\;]0,+\infty[$ be the (continuous) {\it
first return time}\index{first return!time} on the cross-section
$Y$: for every $y\in Y$,
$$
\tau(y)=\min\{t>0\;:\; \phi_t(y)\in Y\}\;.
$$
Let $\phi_Y:Y\ra Y$ be the (homeomorphic) {\it first return map}%
\index{first return!map} to  (or {\it Poincaré 
map}\index{Poincar\'e!map} of) the cross-section $Y$, defined by
$$
\phi_Y:y\mapsto \phi_{\tau(y)}(y)\;.
$$ 
By the invariance of $\mu$ under the flow $(\phi_t)_{t\in\RR}$, the
restriction of $\mu$ to 
$$
\{\phi_t(y)\;:\; y\in Y,\; 0\leq t<\tau(y)\}
$$
disintegrates\footnote{with conditional measure on the fiber
  $\{\phi_t(y)\;:\; 0\leq t<\tau(y)\}$ over $y\in Y$ the image of the
  Lebesgue measure on $[0,\tau(y)[$ by $t\mapsto \phi_t(y)$} by the
(well-defined) map $\phi_t(y)\mapsto y$ over a measure $\mu_Y$ on $Y$,
which is invariant under the first return map $\phi_Y$:
$$
d\mu(\phi_t(y))\;=\;dt\;d\mu_Y(y)\;.
$$ 
Note that if $\tau$ has a positive lower bound and if $\mu$ is finite,
then $\mu_Y$ is finite, since 
$$
\|\mu\|\geq \|\mu_Y\|\;\inf \tau\;,
$$
and
 $(Y,\mu_Y,\phi_Y)$ is a discrete time dynamical system.

\medskip
Recall that an {\it isomorphism}\index{isomorphism} from a
continuous time dynamical system $(Z,\mu, (\phi_t)_{t\in\RR})$ to
another one $(Z',\mu', (\phi'_t)_{t\in\RR})$ is a homeomorphism
between the underlying spaces preserving the underlying measures and
commuting with the underlying flows.

\bexem\label{exem:isomsuspension}
  {\rm If $(Z,\mu, T)$ and $(Z',\mu', T')$ are (invertible)
  discrete time dynamical systems, endowed with roof functions $r:Z\ra
  \;]0,+\infty[$ and $r':Z'\ra \;]0,+\infty[$ respectively, if
  $\theta: Z\ra Z'$ is a measure preserving homeomorphism commuting
  with the transformations $T$ and $T'$ (that is, $\theta_*\mu=\mu'$,
  $\theta\circ T=T'\circ\theta$) and such that 
$$
r'\circ\theta=r\;,
$$ 
then the map $\widecheck{\theta}: Z_r\ra Z'_{r'}$ defined by
$[z,s]\mapsto [\theta(z),s]$ is an isomorphim between the suspensions
$(Z,\mu, T)_r$ and $(Z',\mu', T')_{r'}$.}
\eexem

\medskip
It is well known (see for instance \cite[\S 1.11]{BriStu02}) that the
above two constructions are inverses one to another, up to
isomorphism. In particular, we have the following result.

\bprop\label{prop:inverseinduiresuspendre} With the general notation
above Example \ref{exem:isomsuspension}, the suspension
$(Y,\mu_Y,\phi_Y)_\tau$ over $(Y,\mu_Y,\phi_Y)$ with roof function
$\tau$ is isomorphic to $(Z,\mu, (\phi_t)_{t\in\RR})$ by the map
$f_Y:[y,s]\mapsto \phi_sy$.
\cqfd \eprop

\medskip
In order to describe the continuous time dynamical system $\big(\Ga\bs
\G X',\frac{m_F}{\|m_F\|},(\flow{t})_{t\in\RR}\big)$ as a suspension
over a topological Markov shift, we will start by describing it as a
suspension of the discrete time geodesic flow on $\Ga\bs \G
\XX'$. Note that the Patterson densities and Gibbs measures depend not
only on the potential, but also on the lengths of the
edges.\footnote{The fact that the Patterson densities could be
  singular one with respect to another when the metric varies is a
  well known phenomenon, even when the potential vanishes. For
  instance, let $\Sigma=\Ga\bs\hdr$ and $\Sigma'=\Ga'\bs\hdr$ be two
  closed connected hyperbolic surfaces, uniformised by the real
  hyperbolic plane $(\hdr,ds^2_{\rm hyp})$ endowed with torsion free
  cocompact Fuchsian groups $\Ga$ and $\Ga'$. Let
  $\phi:\Sigma\ra\Sigma'$ be a diffeomorphism, with lift
  $\wt\phi:\hdr\ra \hdr$. Then $\Ga$ is a discrete group of isometries
  for the two $\CAT(-1)$-metrics $ds^2_{\rm hyp}$ and
  $\wt\phi^*ds^2_{\rm hyp}$. Kuusalo's theorem \cite{Kuusalo73} says
  that the corresponding two Patterson densities of $\Ga$ are
  absolutely continuous one with respect to the other if and only if
  $\phi$ is isotopic to the identity. See an extension of this result
  in \cite{HerPau97}. See also the result of \cite{KapNag07} which
  parametrises the Culler-Vogtmann space using Patterson densities for
  cocompact and free actions of free groups on metric trees.}  We
hence need to relate precisely the continuous time and discrete time
situations, and we will use in this Section the left exponent $\sharp$
to indicate a discrete time object whenever needed.

For instance, we set\footnote{See Section \ref{subsec:trees} for the
  definition of the geometric realisation $|\YY|_1$ of a simplicial
  tree $\YY$.} $^\sharp X'=|\XX'|_1$ and we denote by $(^\sharp
\flow{t})_{t\in\ZZ}$ the discrete time geodesic flow on
$\Ga\bs\G\XX'$.  Note that $X'$ and $^\sharp X'$ are equal as
topological spaces (but not as metric spaces). The boundaries at
infinity of $X'$ and $^\sharp X'$, which coincide with their spaces of
ends as topological spaces (by the assumption on the lengths of the
edges), are hence equal and denoted by $\partial_\infty\XX$.

We may assume by Section \ref{subsec:cond} that the potential $\wt
F:T^1X\ra\RR$ is the potential $\wt F_c$ associated with a system of
conductances $\wt c$ on the metric tree $(\XX,\lambda)$ for $\Ga$. Let
$\delta_c=\delta_{F_c}$. We denote by $^\sharp \wt c:E\XX\ra\RR$ the
$\Ga$-invariant system of conductances
\begin{equation}\label{eq:defisharpc}
^\sharp \wt c : e\mapsto (\wt c(e)-\delta_c)\lambda(e)
\end{equation}
on the simplicial tree $\XX$ for $\Ga$, by $\wt F_{^\sharp c}:
T^1(\,^\sharp X)\ra\RR$ its associated potential, and by $^\sharp
c:\Ga\bs E\XX\ra\RR$ and $F_{^\sharp c}: \Ga\bs T^1(\,^\sharp X)\ra\RR$
their quotient maps.

Note that the inclusion morphism $\Aut(\XX,\lambda)\ra\Aut(\XX)$ is a
homeomorphism onto its image (for the compact-open topologies), by the
assumption of a positive lower bound on the lengths of the edges,
hence that $\Ga$ is also a nonelementary discrete subgroup of
$\Aut(\XX)$.

Now, let $(\Sigma,\sigma,\PP)$ be the (two-sided) topological Markov
shift conjugated to the discrete time geodesic flow $\big( \Ga\bs
\G\XX',\;^\sharp\flow{1},\frac{m_{F_{^\sharp c}}}{\|m_{F_{^\sharp c}}\|} \big)$ 
by the bilipschitz homeomorphism $\Theta: \Ga\bs\G\XX'\ra\Sigma$ of
Theorem \ref{theo:codingdisgeodflo} (where the potential $F$ is
replaced by $F_{^\sharp c}$). Let $r: \Sigma \ra\; ]0,+\infty[$ be the
map
\begin{equation}\label{eq:defirrrrr}
r:x\mapsto \lambda(e^+_0)
\end{equation}
if $x=(x_n)_{n\in\ZZ}\in\Sigma$ and $x_0=(e_0^-,h_0,e_0^+)\in\A$.
This map is locally constant, hence continuous on $\Sigma$, and has a
positive lower bound, since the lengths of the edges of $(\XX',
\lambda)$ have a positive lower bound.

\btheo \label{theo:mainsuspension} Assume that the lengths of the
edges of $(\XX',\lambda)$ have a finite
upper bound, and that the Gibbs measure $m_{F}$ is finite.  Then there
exists $a>0$ such that the continuous time dynamical system $\big(
\Ga\bs \G X', \frac{m_{F}}{\|m_{F}\|}, (\flow{t})_{t\in\RR} \big)$ is
isomorphic to the suspension $(\Sigma,\sigma,a\,\PP)_r$ over
$(\Sigma,\sigma,a\,\PP)$ with roof function $r$, by a bilipschitz
homeomorphism $\Theta_r:\Ga\bs \G X'\ra\Sigma_r$.  
\etheo

\dem Let
$$
Y=\{\ell\in\Ga\bs \G X'\;:\;\ell(0)\in \Ga\bs V\XX\}\;.
$$

Then the (closed) subset $Y$ of $\Ga\bs \G X'$ is a cross-section of
the continuous time geodesic flow $(\flow{t})_{t\in\RR}$, since every
orbit meets $Y$ infinitely many times  in a discrete set of times and
since the lengths of the edges of $(\XX', \lambda)$ have a positive
lower bound. Let $\tau:Y\ra\mathopen{]}0,+\infty\mathclose{[}$ be
the first return time, let $\mu_Y$ be the measure on $Y$ (obtained by
disintegrating $\frac{m_{F}}{\|m_{F}\|}$), and let $\flow{}_Y: Y\ra Y$
be the first return map associated with this cross-section $Y$.

We have a natural reparametrisation map $R:Y\ra \Ga\bs \G \XX'$,
defined by $\ell\mapsto \;^\sharp\ell$, where $^\sharp\ell(n)=
(\flow{}_Y^n\ell)(0)$ is the $n$-th passage of $\ell$ in $V\XX$, for
every $n\in\ZZ$. Since there exist $m,M>0$ such that $\lambda(E\XX)
\subset [m,M]$, the map $R$ is a bilipschitz homeomorphism. It
conjugates the first return map $\flow{}_Y$ and the discrete time
geodesic flow on $\Ga\bs \G \XX'$:
$$
R\circ\flow{}_Y=\;^\sharp\flow{1}\circ R\;.
$$

The main point of this proof is the following result relating the
measures $\mu_Y$ and $m_{F_{{}^\sharp c}}$.

\blemm\label{lem:samePatdenssR}\mbox{}\begin{enumerate}
\item[(1)] The family $(\mu^\pm_x)_{x\in V\XX}$ is a Patterson density
  for $(\Ga,F_{{}^\sharp c})$ on the boundary at infinity of the
  simplicial tree $\XX'$, and the critical exponent $\delta_{{}^\sharp
    c}$ of ${}^\sharp c$ is equal to $0$.
\item[(2)] We have $R_*\frac{\mu_Y}{\|\mu_Y\|} = \frac{m_{F_{{}^\sharp
      c}}}{\|m_{F_{{}^\sharp c}}\|}$.
\end{enumerate}
\elemm

\dem (1) By the definition of the potential associated with a system
of conductances,\footnote{See Section \ref{subsec:cond}.} for every
$x,y\in V\XX'$, if $(e_1,\dots, e_n)$ is the edge path in $\XX$ with
$o(e_1)=x$ and $t(e_n)=y$, then (noting that the integrals along paths
depend on the lengths of the edges, the first one below being in $X'$,
the second one in $^\sharp X'$)
\begin{equation}\label{eq:calintegralmetsimp}
\int_x^y(\wt{F_c}-\delta_c)=
\sum_{i=1}^n (\wt c(e_i)\lambda(e_i)-\delta_c\lambda(e_i))
=\int_x^y\;\wt{F_{^\sharp c}}\;.
\end{equation}
Let us denote\footnote{See Section \ref{subsec:criticexpo}.} by 
$$
^\sharp Q(s)= Q_{\Ga,\,F_{^\sharp c},\,x,\,y}(s)=
\sum_{\ga\in\Ga} \; e^{\int_{x}^{\ga y} (\wt F_{^\sharp c}-s)} 
$$ 
and
$$
Q(s) = Q_{\Ga,\,F_{c}-\delta_c,\,x,\,y}(s) = \sum_{\ga\in\Ga}\;
e^{\int_{x}^{\ga y} (\wt F_{c}-\delta_c-s)}
$$
the Poincaré series for the simplicial tree with potential $\wt
F_{^\sharp c}$ and for the metric tree with normalised potential $\wt
F_{c}-\delta_c$, respectively. With $M$ an upper bound on the lengths
of the edges of $(\XX,\lambda)$, the distances $d_{{}^\sharp X'}$ on
${}^\sharp X'$ and $d_{X'}$ on $X'$ satisfy $d_{{}^\sharp
  X'}\geq\frac{1}{M}\;d_{X'}$ on the pairs of vertices of $\XX'$. We
hence have $^\sharp Q(s)\leq Q(\frac{s}{M})<+\infty$ if $s>0$ and
$^\sharp Q(s)\geq Q(\frac{s}{M})=+\infty$ if $s<0$. Thus the critical
exponent $\delta_{^\sharp c}$ of $(\Ga,\,F_{^\sharp c})$ for the
simplicial tree $\XX'$ is equal to $0$, hence $F_{^\sharp c}$ is a
normalised potential.

By the definition\footnote{See Section \ref{subsec:Gibbscocycle}.} of
the Gibbs cocycles (which uses the normalised potential), Equation
\eqref{eq:calintegralmetsimp} also implies that the Gibbs cocycles
$C^\pm$ and $^\sharp C^\pm$ for $(\Ga,F_c)$ and $(\Ga,F_{^\sharp c})$
respectively coincide on $\partial_\infty\XX\times V\XX\times V\XX$.
Thus by Equations \eqref{eq:equivarPatdens} and
\eqref{eq:quasinivarPatdens}, the family $(\mu^\pm_x)_{x\in V\XX}$ is
indeed a Patterson density for $(\Ga,F_{^\sharp c})$: for all
$\ga\in\Ga$ and $x,y \in V\XX$, and for (almost) all
$\xi\in\partial_\infty \XX$,
$$
\ga_*\mu^\pm_x=\mu^\pm_{\ga x} \;\;\;{\rm and}\;\;\;
\frac{d\mu^\pm_x}{d\mu^\pm_y}(\xi) = e^{-\,^\sharp C^\pm_\xi(x,\,y)}\,.
$$

\medskip
\noindent (2) We may hence choose these families $(\mu^\pm_x)_{x\in
  V\XX}$ in order to define the Gibbs measure $m_{F_{^\sharp c}}$
associated with the potential $F_{^\sharp c}$ on $\Ga\bs\G\XX$. Note
that since we will prove that $m_{F_{^\sharp c}}$ is finite, the
normalised measure $\frac{m_{F_{^\sharp c}}} {\|m_{F_{^\sharp c}}\|}$ is
independent of this choice (see Corollary
\ref{coro:finitudeGibbsdivuniq}).

Let $\wt Y=\{\ell\in \G X'\;:\;\ell(0)\in V\XX'\}$ be the
($\Ga$-invariant) lift of the cross-section $Y$ to $\G X'$, let $\wt
R:\wt Y\ra \G \XX'$ be the lift of $R$, mapping a geodesic line
$\ell\in\wt Y$ to a discrete geodesic line $^\sharp\ell$ with same
footpoint obtained by reparametrisation, and let $\wt{\mu_Y}$ be the
measure on $\wt Y$ whose induced measure on $Y=\Ga\bs \wt Y$ is
$\mu_Y$. We have a partition of $\wt Y$ into the closed-open subsets
$\wt Y_{x}=\{\ell\in \G X'\;:\;\ell(0)=x\}$ as $x$ varies in $V\XX'$.

Let us fix $x\in V\XX$. By the definition of $\mu_Y$ as a
disintegration of $\frac{m_{F_c}}{\|m_{F_c}\|}$ with respect to the
continuous time, by lifting to $\G X'$, by using Hopf's
parametrisation with respect to $x$ and Equation \eqref{eq:defigibbs}
with $x_0=x$ , we have for every $\ell\in \wt Y_{x}$,
$$
d\wt{\mu_Y}(\ell)=
\frac{1}{\|m_{F_c}\|}\; d\mu^-_{x}(\ell_-)\;d\mu^+_{x}(\ell_+)\;.
$$ 
Note that $\ell(0)=\;^\sharp\ell(0)$, $\ell_-=\;^\sharp\ell_-$,
$\ell_+= \;^\sharp\ell_+$ since the reparametrisation does not change
the origin or the two points at infinity.  Hence by Assertion (1), we
have
$$ 
\wt R_*(\wt{\mu_Y}) = \frac{1}{\|m_{F_c}\|}\;\wt m_{F_{^\sharp c}}\;.
$$ 
As $\mu_Y$ is a finite measure since $\tau$ has a positive lower
bound, this implies that $m_{F_{^\sharp c}}$ is finite. By
renormalizing as probability measures, this proves Assertion (2).
\cqfd

\bigskip
Let $a=\|\mu_Y\|>0$, so that by Lemma \ref{lem:samePatdenssR} (2) we
have $R_*\mu_Y=\frac{a\;m_{F_{^\sharp c}}}{\|m_{F_{^\sharp c}}\|}$.
Let $\widecheck{r} :\Ga\bs\G \XX'\ra\;]0,+\infty]$ be the map
$$
\widecheck{r}:\Ga\ell\mapsto \lambda\big(\ell([0,1])\big)
$$ 
given by the length for $\lambda$ of the first edge followed by a
discrete geodesic line $\ell\in \G \XX'$. Note that $\widecheck{r}$ is
locally constant, hence continuous, and that $\widecheck{r}$ is a roof
function for the discrete time dynamical system $(\Ga\bs\G
\XX',\;^\sharp\flow{1})$. Also note that
$$
\widecheck{r}\circ R = \tau\;\;\;\;{\rm and} \;\;\;\;
r\circ\Theta=\widecheck{r}
$$
by the definitions of $\tau$ and $r$.

Let us finally define $\Theta_r: \Ga\bs\G X'\ra\Sigma_r$
as the compositions of the following three maps
{\small
\begin{equation}\label{eq:defiThetar}
(\Ga\bs\G X',\frac{m_{F}}{\|m_{F}\|},(\flow{t})_{t\in\RR})
\;\stackrel{f_Y^{\;-1}}{\longrightarrow}\;
(Y,\mu_Y,g_Y)_\tau
\;\stackrel{\widecheck{R}}{\longrightarrow}\;
(\Ga\bs\G \XX',\frac{a\;m_{F_{^\sharp c}}}{\|m_{F_{^\sharp c}}\|},
\;^\sharp\flow{1})_{\widecheck{r}}
\;\stackrel{\widecheck{\Theta}}{\longrightarrow}\;
(\Sigma,a\,\PP,\sigma)_r\;,
\end{equation}
}

\noindent
where the first one is the inverse of the tautological isomorphism
given by Proposition \ref{prop:inverseinduiresuspendre} and the last
two ones, given by Example \ref{exem:isomsuspension}, are the
isomorphisms $\widecheck{R}$ and $\widecheck{\Theta}$ of continuous
time dynamical systems obtained by suspensions of the isomorphisms $R$
and $\Theta$ of discrete time dynamical systems. It is easy to check
that $\Theta_r$ is a bilipschitz homeomorphism, using the following
description of the Bowen-Walters distance, see for instance
\cite[Appendix]{BarSau00}.

\bprop \label{prop:defiBowenWalters} Let $(Z,\mu, T)_r$ be the
suspension over an invertible dynamical system such that $T$ is a
bilipschitz homeomorphism, with roof function $r$ having a positive
lower bound and a finite upper bound. Let $d_{\rm BW} : Z_r\times
Z_r\ra\RR$ be the map\footnote{The map $d_{\rm BW}$ is actually not a
  distance, but this proposition says that it may replace the
  Bowen-Walters true distance when working up to multiplicative
  constants or bilipschitz homeomorphisms.} defined (using the
canonical representatives) by
\begin{multline*}
d_{\rm BW}([x,s],[x',s'])=\\ \min\{d(x,x')+|s-s'|,
\;\;d(Tx,x') + r(x) - s + s',\;\;d(x,Tx') + r(x') + s - s'\}
\end{multline*}
Then there exists a constant $C_{\rm BW}>0$ such that the
Bowen-Walters distance $d$ on $\Sigma_r$ satisfies
$$
\frac{1}{C_{\rm BW}}\;d_{\rm BW}\;\leq \;d\;\leq \;C_{\rm BW}\;d_{\rm BW}\;.
\;\;\;\Box
$$
\eprop

\medskip
This concludes the proof of Theorem \ref{theo:mainsuspension}.  \cqfd

\section{The variational principle for metric and 
simplicial trees} 
\label{subsec:varprinc}

In this Section, we assume that $X$ is the geometric realisation of a
locally finite metric tree without terminal vertices $(\XX,\lambda)$
(respectively of a locally finite simplicial tree $\XX$ without
terminal vertices).  Let $\Ga$ be a nonelementary discrete subgroup of
$\Aut(\XX,\lambda)$ (respectively $\Aut(\XX)$). 

We relate in this Section the Gibbs measures\footnote{See the
  definition in Sections \ref{subsec:gibbsmeasure} and
  \ref{subsec:pattersongibbstrees}.} to the equilibrium
states\footnote{See the definitions below.} for the continuous time
geodesic flow on $\Ga\bs\G X$ (respectively for the discrete time
geodesic flow on $\Ga\bs\G\XX$).

When $X$ is a Riemannian manifold with pinched negative curvature such
that the derivatives of the sectional curvature are uniformly bounded,
and when the potential is H\"older-continuous, the analogs of the
results of this Section are due to \cite[Theo.~6.1]{PauPolSha15}. Their
proofs generalise the proofs of Theorems 1 and 2 of \cite{OtaPei04},
with ideas and techniques going back to \cite{LedStr82}. When $Y$ is a
compact locally $\CAT(-1)$-space, a complete statement about
existence, uniqueness and Gibbs property of equilibrium states for any
H\"older-continuous potential is given in \cite{ConLafTho16}.

The proof of the metric tree case will rely strongly (via the
suspension process described in Section \ref{subsec:codagemetric}) upon
the proof of the simplicial tree case, hence we start by the latter.

\bigskip
\noindent{\bf The simplicial tree case. }

\bigskip
Let $\XX$ be a locally finite simplicial tree without terminal
vertices, with geometric realisation $X=|X|_1$. Let $\Ga$ be a
nonelementary discrete subgroup of $\Aut(\XX)$. Let $\wt c: E\XX
\ra\RR$ be a system of conductances for $\Ga$ on $\XX$ and $c: \Ga\bs
E\XX \ra\RR$ its quotient map. Let $\wt F_c: T^1X \ra \RR$ be its
associated potential, with quotient map $F_c: \Ga\bs T^1X \ra \RR$,
and let $\delta_c$ be the critical exponent of $c$.

We define a  map $\widetilde{\FF_c}:\G\XX\ra \RR$ by
$$
\widetilde{\FF_c}(\ell)=\wt c\,(e^+_0(\ell))=
\int_{o(e^+_0(\ell))}^{t(e^+_0(\ell))}\wt F_c
$$ 
for all $\ell\in \G\XX$, where $e^+_0(\ell)$ is the edge of $\XX$
in which $\ell$ enters at time $t=0$. This map is locally constant,
hence continuous, and it is $\Ga$-invariant, hence it induces a
continuous map $\FF_c:\Ga\bs\G\XX\ra\RR$ which is also called a {\em
  potential}.\footnote{See after the proof of Theorem
  \ref{theo:varprinarbsimp} for a comment on cohomology classes.}

The following result proves that the Gibbs measure of $(\Ga,F_c)$ for
the discrete time geodesic flow on $\Ga\bs\G\XX$ is an equilibrium
state for the potential $\FF_c$. We start by recalling the definition
of an equilibrium state,\footnote{This definition is given for
  transformations and not for flows, and for possibly unbounded
  potentials, contrarily to the one the Introduction.} see also
\cite{Bowen75,Ruelle04}.

Let $Z$ be a locally compact topological space, let $T:Z\ra Z$ be a
homeomorphism, and let $\phi:Z\ra\RR$ be a continuous map.  Let
$\M_\phi$ be the set of Borel probability measures $m$ on $Z$,
invariant under the transformation $T$, such that the negative part
$\phi^-=\max\{0,-\phi\}$ of $\phi$ is $m$-integrable. Let
$\gls{metricentropydisc}$ be the (metric) entropy of the
transformation $T$ with respect to $m\in\M_\phi$ (see for instance
\cite{BriStu02}).  The {\em metric
  pressure}\index{metric!pressure}\index{pressure!metric} for the
potential $\phi$ of a measure $m$ in $\M_\phi$ is
$$
\gls{metricpressuredisc}=h_m(T)+\int_{Z}\phi\,dm\;.
$$
The fact that the negative part of $\phi$ is $m$-integrable, which is
in particular satisfied if $\phi$ is bounded, implies that
$\pressure{\phi}(m)$ is well defined in $\RR\cup\{+\infty\}$. The {\em
  pressure}\index{pressure} of the potential $\phi$ is the element of
$\RR\cup\{+\infty\}$ defined by
$$
\gls{pressuredisc}=\sup_{m\in\M_\phi}\;\pressure{\phi}(m)\;.
$$ A measure $m_0$ in $\M_\phi$ is an {\em equilibrium state}%
\index{equilibrium state} for the potential $\phi$ if
$\pressure{\phi}(m_0)= \pressure{\phi}$.

\btheo[The variational principle for simplicial
  trees]\label{theo:varprinarbsimp} Let $\XX,\Ga,\wt c$ be as
above. Assume that $\delta_c<+\infty$ and that there exists a finite
Gibbs measure $m_c$ for $F_c$ such that the negative part of the
potential $\FF_c$ is $m_c$-integrable. Then $\frac{m_c}{\|m_c\|}$ is
the unique equilibrium state for the potential $\FF_c$ under the
discrete time geodesic flow on $\Ga\bs\G\XX$, and the pressure of
$\FF_c$ coincides with the critical exponent $\delta_c$ of $c$~:
$$
\pressure{\FF_c}=\delta_c\;.
$$
\etheo

In order to prove this result, using the coding of the discrete time
geodesic flow given in Section \ref{subsec:codagesimplicial}, the main
tool is the following result of J.~Buzzi in symbolic dynamics,
building on works of Sarig and Buzzi-Sarig, whose proof is given in
Appendix \ref{appendixBuzzi}.

Let $\sigma:\Sigma\ra\Sigma$ be a two-sided topological Markov
shift\footnote{See Section \ref{subsec:TMS} for definitions.}  with
(countable) alphabet $\A$ and transition matrix $A$, and let
$\phi:\Sigma\ra\RR$ be a continuous map.

For every $n\in\NN$, we denote\footnote{with a shift of indices
  compared with the notation of Appendix \ref{appendixBuzzi}} by
$$
\operatorname{var}_n\phi=\sup_{\substack{ x,\,y\,\in\,\Sigma\\
\forall\;i\,\in\,\{-n,\,\dots,\,n\},\;x_i=y_i}} |\phi(x)-\phi(y)|
$$ 
the {\it $n$-variation}\index{variation@$n$-variation} of $\phi$.  For
instance, if $\phi(x)$ depends only on $x_0$ where $x=
(x_i)_{i\in\ZZ}$, then $\operatorname{var}_n\phi =0$ for every $n\in\NN$
(and hence $\sum_{n\in\NN}\; (n+1)\operatorname{var}_n\phi=0$
converges).

A {\it weak Gibbs measure}\index{weak Gibbs measure}%
\index{Gibbs!measure!weak} for $\phi$ with Gibbs constant $C(m)\in\RR$
is a $\sigma$-invariant (positive Borel) measure $m$ on $\Sigma$ such
that for every $a\in\A$, there exists $c_a\geq 1$ such that for all
$n\in\NN-\{0\}$ and $x\in [a]$ such that $\sigma^n(x)=x$, we have
\begin{equation}\label{eq:defweakgibbs}
\frac{1}{c_a}\leq \frac{m([x_{0},\dots, x_{n-1}])}
{e^{-C(m)\,n}\;e^{\,\sum_{i=0}^{n-1}\phi(\sigma^i x)}}\leq c_a\;.
\end{equation}

\btheo [J.~Buzzi, see Corollary
  \ref{coro:mainappendix_twosided}]\label{theo:buzzi} Let
$(\Sigma,\sigma)$ be a two-sided transitive topological Markov shift
on a countable alphabet and let $\phi:\Sigma\ra\RR$ be a
continuous map such that $\sum_{n\in\NN}\; (n+1)\operatorname{var}_n\phi$
converges. Let $m$ be a weak Gibbs measure for $\phi$ on $\Sigma$ with
Gibbs constant $C(m)$, such that $\int \phi^-\;dm<+\infty$. Then the
pressure of $\phi$ is finite, equal to $C(m)$, and $m$ is the unique
equilibrium state.  
\cqfd 
\etheo

\medskip \noindent{\bf Proof of Theorem \ref{theo:varprinarbsimp}. }
In Section \ref{subsec:codagesimplicial}, we constructed a transitive
topological Markov shift $(\Sigma,\sigma)$ on a countable alphabet
$\A$ and a homeomorphism $\Theta:\Ga\bs\G \XX' \ra\Sigma$ which
conjugates the time-one discrete geodesic flow $\flow{1}$ on the
nonwandering subset $\Ga\bs\G \XX'$ of $\Ga\bs\G \XX$ and the shift
$\sigma$ on $\Sigma$ (see Theorem \ref{theo:codingdisgeodflo}). 
Let us define a potential $\FF_{c,\,{\rm symb}}:\Sigma\ra \RR$ by
\begin{equation}\label{eq:defipotentsymbolbis}
\FF_{c,\,{\rm symb}}(x)= c(e^+_0)
\end{equation}
if $x=(x_i)_{i\in\ZZ}\in\Sigma$ with $x_0=(e_0^-,h_0,e^+_0)$.  Note
that this potential is the one denoted by $F_{\rm symb}$ in Equation
\eqref{eq:defipotentsymbol}, when the potential $F$ on $T^1X$ is
replaced by $F_c$. By the construction of $\Theta$ and the definition
of $\FF_c$, we have
\begin{equation}\label{FsymbThetaFFc}
\FF_{c,\,{\rm symb}}\circ\Theta=\FF_c\;.
\end{equation} 
Note that all probability measures on $\Ga\bs\G\XX$ invariant under
the discrete  time geodesic flow are supported on the nonwandering set
$\Ga\bs\G\XX'$. The pushforward of measures $\Theta_*$ hence gives a
bijection from the space $\M_{\FF_c}$ of $\flow{1}$-invariant
probability measures on $\Ga\bs\G\XX$ for which the negative part of
$\FF_c$ is integrable to the space $\M_{\FF_{c,\,{\rm symb}}}$ of
$\sigma$-invariant probability measures on $\Sigma$ for which the
negative part of $\FF_{c,\,{\rm symb}}$ is integrable.  This bijection
induces a bijection between the subsets of equilibrium states.

Note that $\FF_{c,\,{\rm symb}}(x)$ depends only on $x_0$ for every
$x=(x_i)_{i\in\ZZ}\in\Sigma$. Hence as seen above, the series
$\sum_{n\in\NN}\; (n+1)\operatorname{var}_n \FF_{c,\,{\rm symb}}$
converges.

By definition,\footnote{See Equation \eqref{eq:defPP}.} the measure
$\PP$ is the pushforward of $\frac{m_c}{\|m_c\|}$ by $\Theta$. Hence,
$\PP$ is a $\sigma$-invariant probability measure on $\Sigma$ for
which the negative part of $\FF_{c,\,{\rm symb}}$ is integrable, by
the assumption of Theorem \ref{theo:varprinarbsimp}. By Proposition
\ref{prop:propriPP} (3), the measure $\PP$ on $\Sigma$ satisfies the
Gibbs property with Gibbs constant $\delta_c$ for the potential
$\FF_{c,\,{\rm symb}}$, hence\footnote{For every $a\in\A$, for the
  constant $c_a$ required by the definition of the weak Gibbs property
  in Equation \eqref{eq:defweakgibbs}, take the constant $C_{E}$ given
  by the definition (see Equation \eqref{eq:defiGibbspropertyZ}) of
  the Gibbs property with $E=\{a\}$.}  satisfies the weak Gibbs
property with Gibbs constant $\delta_c$.  Theorem
\ref{theo:varprinarbsimp} then follows from Theorem \ref{theo:buzzi}.
\cqfd

\medskip
\rem It follows from Equation \eqref{FsymbThetaFFc}, from the remark
above Definition \ref{def:Markovgood} and from the fact that
$\Theta\circ \flow{1}=\sigma\circ\Theta$, that if $c,c':\Ga\bs
E\XX'\ra\RR$ are cohomologous systems of conductances on $\Ga\bs
E\XX'$, then the corresponding maps $\FF_c, \FF_{c'}:\Ga\bs
\G\XX'\ra\RR$ are cohomologous: there exists a continuous map
$\GG:\Ga\bs \G\XX'\ra \RR$ such that for every $\ell\in\Ga\bs \G\XX'$,
$$
\FF_{c'}(\ell)-\FF_c(\ell)=\GG(\flow{1} \ell)-\GG(\ell)\;.
$$

Note that given two bounded cohomologous continuous potentials on a
topological dynamical system $(Z,T)$ as above, one has finite pressure
(resp.~admits an equilibrium state) if and only if the other one does,
and they have the same pressure and same set of equilibrium states.

\bigskip
\noindent{\bf The metric tree case. }

\bigskip
Let $(\XX,\lambda)$ be a locally finite metric tree without terminal
vertices with geometric realisation $X=|\XX|_\lambda$, let $\Ga$ be a
nonelementary discrete subgroup of $\Aut(\XX,\lambda)$ and let $\wt
c:E\XX\ra\RR$ be a system of conductances for $\Ga$ on $\XX$. Let $\wt
F_c:T^1X \ra \RR$ be its associated potential (see Section
\ref{subsec:cond}), and let $\delta_c=\delta_{F_c}$ be the critical
exponent of $c$.

Recall\footnote{See Section \ref{subsec:unitbundle}.} that we have a
canonical projection $\G X\ra T^1 X$ which associates to a geodesic
line $\ell$ its germ $v_\ell$ at its footpoint $\ell(0)$. Let $\wt
\FF^\natural_c:\G X\ra \RR$ be the $\Ga$-invariant map obtained by
precomposing the potential $\wt F_c:T^1 X\ra \RR$ with this canonical
projection:
$$
\wt \FF^\natural_c : \ell\mapsto \;\wt F_c(v_\ell)\;.
$$ 
Let $\FF^\natural_c:\Ga\bs \G X\ra \RR$ be its quotient map, which is
continuous as a composition of continuous maps.

The following result proves that the Gibbs measure of $(\Ga,F_c)$ for
the continuous time geodesic flow on $\Ga\bs\G X$, once renormalised
to be a probability measure, is an equilibrium state for the potential
$\FF^\natural_c$. We start by recalling the definition of an
equilibrium state for a possibly unbounded potential under a
flow.\footnote{This requires only minor modifications to the
  definition given in the Introduction for bounded potentials.}

Given $(Z,(\phi_t)_{t\in\RR})$ a topological space endowed with a
continuous one-parameter group of homeomorphisms and $\psi:Z\ra\RR$ a
continuous map (called a {\em potential}\index{potential}), let
$\M_\psi$ be the set of Borel probability measures $m$ on $Z$
invariant under the flow $(\phi_t)_{t\in\RR}$, such that the negative
part $\psi^-=\max\{0,-\psi\}$ of $\psi$ is $m$-integrable.  Let
$\gls{metricentropy}$ be the (metric) entropy of the geodesic flow
with respect to $m\in\M_\psi$ (see for instance \cite{BriStu02}).  The
{\em metric pressure}\index{metric!pressure}\index{pressure!metric}
for $\psi$ of a measure $m\in\M_{\psi}$ is
$$
\gls{metricpressure}=h_m(\phi_1)+\int_{Z}\,\psi\;dm\;.
$$
The fact that the negative part of $\psi$ is $m$-integrable, which
is in particular satisfied if $\psi$ is bounded, implies that
$\pressure{\psi}(m)$ is well defined in $\RR\cup\{+\infty\}$. The {\em
  pressure}\index{pressure} of the potential $\psi$ is the element of
$\RR\cup\{+\infty\}$ defined by
$$
\gls{pressure}=\sup_{m\in\M_{\psi}}\;\pressure{\psi}(m)\;.
$$
Note that $\pressure{\psi+c}=\pressure{\psi}+c$ for every constant
$c\in\RR$. An element $m\in\M_\psi$ is an {\em equilibrium
  state}\index{equilibrium state} for $\psi$ if the least upper bound
defining $\pressure{\psi}$ is attained on $m$.

Note that if $\psi'$ is another potential {\it
  cohomologous}\index{cohomologous} to $\psi$, that is, if there
exists a continuous map $G:Z\ra\RR$, differentiable along every orbit
of the flow, such that $\psi'(x)-\psi(x)= \frac{d}{dt}_{\mid t=0}
G(\phi_t(x))$, if $\psi$ and $\psi'$ are bounded (so that $\M_{\psi'}=
\M_\psi$), then, for every $m\in\M_\psi$, we have $\pressure{\psi'}(m)
= \pressure{\psi}(m)$, $\pressure{\psi'}= \pressure{\psi}$ and the
equilibrium states for $\psi'$ are exactly the equilibrium states for
$\psi$.

\btheo[The variational principle for metric
  trees]\label{theo:varprinarbmet} Let $(\XX,\lambda), \Ga,\wt c$ be
as above. Assume that the lengths of the edges of $(\XX,\lambda)$ have
a finite upper bound.\footnote{They have a positive lower bound by
  definition, see Section \ref{subsec:trees}.} Assume that
$\delta_c<+\infty$ and that there exists a finite Gibbs measure $m_c$
for $F_c$ such that the negative part of the potential
$\FF^\natural_c$ is $m_c$-integrable.  Then $\frac{m_c}{\|m_c\|}$ is
the unique equilibrium state for the potential $\FF^\natural_c$ under
the continuous time geodesic flow on $\Ga\bs\G X$, and the pressure
 of $\FF^\natural_c$ coincides with the
critical exponent $\delta_c$ of $c$~:
$$
\pressure{\FF^\natural_c}=\delta_c\;.
$$
\etheo

Using the description of the continuous time dynamical system
$\big(\Ga\bs \G X',\frac{m_c}{\|m_c\|},(\flow{t})_{t\in\RR}\big)$ as a
suspension over a topological Markov shift (see Theorem
\ref{theo:mainsuspension}), this statement reduces to well-known
techniques in the thermodynamic formalism of suspension flows, see for
instance \cite{IomJorTod15}, as well as \cite{BarIom06, Kempton11,
  IomJor13, JaeKesLam14}. Our situation is greatly simplified by the
fact that our roof function has a positive lower bound and a finite
upper bound, and that our symbolic potential is constant on the
$1$-cylinders $\{x=(x_i)_{i\in\ZZ}\in \Sigma\;:\;x_0=a\}$ for $a$ in
the alphabet.

\medskip
\dem 
Since finite measures invariant under the geodesic flow on
$\Ga\bs\G X$ are supported on its nonwandering set, up to replacing
$X$ by $X'=\C\Lambda \Ga$, we assume that $X=X'$.

Since equilibrium states are unchanged up to adding a constant to the
potential, under the assumptions of Theorem \ref{theo:varprinarbmet},
let us prove that $\frac{m_c}{\|m_c\|}$ is the unique equilibrium
state for the potential $\FF^\natural_c-\delta_c$ under the continuous
time geodesic flow on $\Ga\bs\G X$, and that the pressure of
$\FF^\natural_c-\delta_c$ vanishes. The last claim of Theorem
\ref{theo:varprinarbmet} follows, since
$$
\pressure{\FF^\natural_c}-\delta_c=\pressure{\FF^\natural_c-\delta_c}\;.
$$

We refer to the paragraphs before the statement of Theorem
\ref{theo:mainsuspension} for the definitions of

\smallskip
$\bullet$~ the system of conductances $^\sharp \wt c$ for $\Ga$ on the
simplicial tree $\XX$,

\smallskip
$\bullet$~ the (two-sided) topological Markov shift $(\Sigma,
\sigma,\PP)$ on the alphabet $\A$, conjugated to the discrete time
geodesic flow $\big( \Ga\bs \G\XX,\;^\sharp\flow{1},
\frac{m_{F_{^\sharp c}}}{\|m_{F_{^\sharp c}}\|} \big)$ by the
homeomorphism $\Theta: \Ga\bs\G\XX\ra\Sigma$,

\smallskip
$\bullet$~ the roof function $r: \Sigma \ra\; ]0,+\infty[$, and

\smallskip
$\bullet$~ the suspension $(\Sigma,\sigma, a\,\PP)_r = (\Sigma_r,
(\sigma^t_r)_{t\in\RR},a\,\PP_r)$ over $(\Sigma,\sigma,a\,\PP)$ with
roof function $r$, conjugated to the continuous time geodesic flow
$\big(\Ga\bs \G X, \frac{m_c}{\|m_c\|}, (\flow{t})_{t\in\RR}\big)$ by
the homeomorphism $\Theta_r:\Ga\bs \G X\ra\Sigma_r$ defined at the end
of the proof of Theorem \ref{theo:mainsuspension}. We will always
(uniquely) represent the elements of $\Sigma_r$ as $[x,s]$ with $x\in
\Sigma$ and $0\leq s < r(x)$.

\medskip
We denote by $\FF^\natural_{c,\,{\rm symb}}:\Sigma_r\ra\RR$ the
potential defined by
\begin{equation}\label{eq:defiFFnaturalcsymb}
\FF^\natural_{c,\,{\rm symb}}= \FF^\natural_c\circ \Theta_r^{-1}\;, 
\end{equation}
which is continuous as a composition of continuous maps. The key
technical observation in this proof is the following one.

\blemm\label{lem:relatavecsansnatural} 
For every $x\in\Sigma$, we have $\FF_{^\sharp c,\,{\rm symb}}(x)
=\int_{0}^{r(x)}(\FF^\natural_{c,\,{\rm symb}}-\delta_c)([x,s])\;ds$.
For every $x\in\Sigma$, the sign of $\FF^\natural_{c,\,{\rm symb}}
([x,s])$ is constant on $s \in \mathopen{[}0, r(x)\mathclose{[}$.
\elemm

\dem 
Let $x=(x_n)_{n\in\ZZ}\in\Sigma$ and $x_0=(e_0^-,h_0,e_0^+)\in\A$.  By
the definition of the first return time $r$ in Equation
\eqref{eq:defirrrrr}, we have in particular
$$
r(x)=\lambda(e^+_0)\;.
$$ 
By Equation \eqref{eq:defipotentsymbolbis} and by the definition of
$^\sharp \wt c$ in Equation \eqref{eq:defisharpc}, we have
$$
\FF_{^\sharp c,\,{\rm symb}}(x)=\;^\sharp c(e^+_0)= 
(c(e^+_0)-\delta_c)\lambda(e^+_0)\;.
$$ 
Using in the following sequence of equalities respectively

$\bullet$~ the two definitions of the potential
$\FF^\natural_{c,\,{\rm symb}}$ in Equation
\eqref{eq:defiFFnaturalcsymb} and of the suspension flow
$(\sigma^t_r)_{t\in\RR}$,

$\bullet$~ the fact that the suspension flow $(\sigma^t_r)_{t\in\RR}$
is conjugated to the continuous time geodesic flow by $\Theta_r$,

$\bullet$~ the definition of $\Theta_r$ using the reparametrisation
map\footnote{See its definition above Lemma \ref{lem:samePatdenssR}.}
$R$ of continuous time geodesic lines with origin on vertices to
discrete time geodesic lines,\footnote{See Equation
  \eqref{eq:defiThetar}, with the notation of Example
  \eqref{exem:isomsuspension} and Proposition
  \ref{prop:inverseinduiresuspendre}.}

$\bullet$~ the definition of the potential $\FF^\natural_c $,

$\bullet$~ the fact that $e^+_0$ is the first edge followed by the
discrete time geodesic line $\Theta^{-1}x$, hence by the continuous
time geodesic line $R^{-1}\Theta^{-1}x$, and the relation between $c$ and the
potential $F_c$ associated with $c$ (see Proposition
\ref{prop:integpotconduct}),

\noindent
we have
\begin{align*}
\int_{0}^{r(x)}\FF^\natural_{c,\,{\rm symb}}([x,s])\;ds&= 
\int_{0}^{r(x)}\FF^\natural_c\big(\Theta_r^{-1}\sigma^s_r[x,0]\big)\;ds= 
\int_{0}^{r(x)}\FF^\natural_c\big(\flow{s}\Theta_r^{-1}[x,0]\big)\;ds\\  &
= \int_{0}^{r(x)}\FF^\natural_c\big(\flow{s}R^{-1}\Theta^{-1}x\big)\;ds= 
\int_{0}^{\lambda(e^+_0)} F_c\big(v_{\flow{s}R^{-1}\Theta^{-1}x}\big)\;ds\\  & =
c(e^+_0)\lambda(e^+_0)\;.
\end{align*}

Since $\int_{0}^{r(x)}\delta_c\;ds=\delta_c\;\lambda(e^+_0)$, the
first claim of Lemma \ref{lem:relatavecsansnatural} follows. The
second claim follows by the definition of the potential $F_c$
associated with $c$, see Equation \eqref{eq:defpotfromconduc}. 
 \cqfd

\medskip
By Equation \eqref{eq:defiFFnaturalcsymb}, the pushforwards of
measures by the homeomorphism $\Theta_r$, which conjugates the flows
$(\flow{t})_{t\in\RR}$ and $(\sigma^t_r)_{t\in\RR}$, is a bijection
from $\M_{\FF^\natural_c}$ to $\M_{\FF^\natural_{c,\,{\rm symb}}}$,
such that
$$
\pressure{\FF^\natural_{c,\,{\rm symb}}} ((\Theta_r)_*m)=
\pressure{\FF^\natural_c}(m)
$$ 
for every $m\in \M_{\FF^\natural_c}$. In particular, we only have to
prove that $(\Theta_r)_*\frac{m_{c}}{\|m_c\|}=a\,\PP_r$ is the unique
equilibrium state for the potential $\FF^\natural_{c,\,{\rm symb}}
-\delta_c$ under the suspension flow $(\sigma^t_r)_{t\in\RR}$, and
that the pressure of $\FF^\natural_{c,\,{\rm symb}}-\delta_c$
vanishes. 

\medskip
The uniqueness follows for instance from
\cite[Theo.~3.5]{IomJorTod15},\footnote{Note that a topological Markov
  shift which is (incorrectly) called topologically mixing in
  \cite[page 551]{IomJorTod15} is actually (topologically) transitive
  with the definition in this book, Section \ref{subsec:TMS}.} since
the roof function $r$ is locally constant and the potential
$g=\FF^\natural_{c,\,{\rm symb}}$ is such that the
map\footnote{denoted by $\Delta_{g}$ in loc.~cit.} from $\Sigma$ to
$\RR$ defined by $x\mapsto \int_0^{r(x)} g([x,s])\,ds$ is locally
H\"older-continuous by Lemma \ref{lem:relatavecsansnatural} and since
$\FF_{^\sharp c,\,{\rm symb}}$ is locally constant.

\medskip
Let us now relate the $\sigma$-invariant measures on $\Sigma$ with the
$(\sigma_r^t)_{t\in\RR}$-invariant measures on $\Sigma_r$. Recall that
we denote the Lebesgue measure on $\RR$ by $ds$ and the points in
$\Sigma_r$ by $[x,s]$ with $x\in\Sigma$ and $0\leq s< r(x)$.

\blemm \label{lem:construcSm}
The map $S:\M_{\FF_{^\sharp c,\,{\rm symb}}}\ra
\M_{\FF^\natural_{c,\,{\rm symb}}}$, which associates to any measure
$m$ in $\M_{\FF_{^\sharp c,\,{\rm symb}}}$ on $\Sigma$ the measure
$$
d\,S(m) ([x,s])=\frac{1}{\int_\Sigma\; r\;dm}\;d\mu(x)\,ds
$$
on $\Sigma_r$, is a bijection, such that, for every $m\in
\M_{\FF_{^\sharp c,\,{\rm symb}}}$,
$$
P_{\FF^\natural_{c,\,{\rm symb}}-\delta_c}(S(m))=
\frac{P_{\FF_{^\sharp  c,\,{\rm symb}}}(m)}{\int_\Sigma\; r\;dm}\;.
$$
\elemm

\dem Note that $\int_\Sigma \;r\;dm$ is the total mass of the measure
$d\mu_r([x,s])= d\mu(x)\,ds$ on $\Sigma_r$. In particular, $S(m)$ is
indeed a probability measure on $\Sigma_r$.

Since $r$ has a positive lower bound and a finite upper bound, it is
well known since \cite{AmbKak42}, see also \cite[\S 2.4]{IomJorTod15},
that the map $S$ defined above\footnote{and denoted by $R$ in \cite[\S
2.4]{IomJorTod15}} is a bijection from the set of $\sigma$-invariant
probability measures $m$ on $\Sigma$ to the set of
$(\sigma_r^t)_{t\in\RR}$-invariant probability measures on $\Sigma_r$.

Furthermore, for every $\sigma$-invariant probability measure $m$ on
$\Sigma$, we have the following {\it Kac formula}\index{Kac formula},
by the definition of the probability measure $S(m)$ and by Lemma
\ref{lem:relatavecsansnatural},
\begin{align}
\int_{\Sigma_r}\FF^\natural_{c,\,{\rm symb}}\;d\,S(m)-\delta_c &=
\int_{\Sigma_r}(\FF^\natural_{c,\,{\rm symb}}-\delta_c)\;d\,S(m)\nonumber\\ & =
\frac{1}{\int_\Sigma \;r\;dm}\; \int_{x\in \Sigma}\int_{0}^{r(x)}
(\FF^\natural_{c,\,{\rm symb}}-\delta_c)([x,s])\;dm(x)\,ds\nonumber\\ & 
=\frac{1}{\int_\Sigma \;r\;dm}\; \int_{\Sigma}
\FF_{^\sharp c,\,{\rm symb}}\;dm\;.\label{eq:integpotensymb}
\end{align}
By the comment on the signs at the end of Lemma
\ref{lem:relatavecsansnatural}, this computation also proves that the
negative part of $\FF^\natural_{c,\,{\rm symb}}$ is integrable for
$S(m)$ if and only if the negative part of $\FF_{^\sharp c,\,{\rm
    symb}}$ is integrable for $m$. Hence $S$ is indeed a bijection
from $\M_{\FF_{^\sharp c,\,{\rm symb}}}$ to $\M_{\FF^\natural_{c,\,{\rm symb}}}$. 

By Abramov's formula \cite{Abramov59}, see also
\cite[Prop.~2.14]{IomJorTod15}, we have
\begin{equation}\label{eq:Abramov}
h_{S(m)}(\sigma^1_r)=\frac{h_m(\sigma)}{\int_\Sigma\; r\;dm}\;.
\end{equation}
The last claim of Lemma \ref{lem:construcSm} follows by summation from
Equations \eqref{eq:integpotensymb} and \eqref{eq:Abramov}.  
\cqfd

\medskip 
By the proof of Theorem \ref{theo:varprinarbsimp} (replacing the
potential $c$ by $^\sharp c$), the pressure of the potential
$\FF_{^\sharp c,\,{\rm symb}}$ is equal to the critical exponent
$\delta_{^\sharp c}$ of the potential $^\sharp c$, and by Lemma
\ref{lem:samePatdenssR} (1), we have $\delta_{^\sharp c} =0$. Hence
for every $m\in \M_{\FF_{^\sharp c,\,{\rm symb}}}$, we have
$$
P_{\FF^\natural_{c,\,{\rm symb}}-\delta_c}(S(m))=
\frac{P_{\FF_{^\sharp  c,\,{\rm symb}}}(m)}{\int_\Sigma \;r\;dm}\leq 
\frac{P_{\FF_{^\sharp c,\,{\rm symb}}}}{\int_\Sigma \;r\;dm}=
\frac{\delta_{^\sharp c}}{\int_\Sigma \;r\;dm}=0\;.
$$ 
In particular, the pressure of the potential $\FF^\natural_{c,\,{\rm
    symb}} -\delta_c$ is at most $0$, since $S$ is a bijection.  By
the proof of Theorem \ref{theo:varprinarbsimp} (replacing the
potential $c$ by $^\sharp c$), we know that $\PP$ is an equilibrium
state for the potential $\FF_{^\sharp  c,\,{\rm symb}}$. Hence
$$ 
P_{\FF^\natural_{c,\,{\rm symb}}-\delta_c}(S(\PP))=
\frac{P_{\FF_{^\sharp c,\,{\rm symb}}}(\PP)}{\int_\Sigma \;r\;d\PP}=
0\;.
$$ 
Therefore, $S(\PP)$ is an equilibrium state of the potential
$\FF^\natural_{c,\,{\rm symb}}-\delta_c$, with pressure $0$.  But
$a\PP_r$, which is equal to $\frac{\PP_r}{\|\PP_r\|}$ since $a\PP_r$
is a probability measure, is by construction equal to $S(\PP)$.  The
result follows.  
\cqfd

\bigskip
With slightly different notation, this result implies Theorem
\ref{theo:varprintreeintro} in the Introduction. 

\medskip
\noindent{\bf Proof of Theorem \ref{theo:varprintreeintro}. } Any
bounded potential $\wt F$ for $\Ga$ on $T^1X$ is cohomologous to a
bounded potential $\wt F_c$ associated with a system of conductances
(see Proposition \ref{prop:relatpotentialconductance}). If two
potentials $\wt F$ and $\wt F'$ for $\Ga$ on $T^1X$ are
cohomologous\footnote{See the definition at the end of Section
  \ref{subsec:potentials}.} then the potentials $\ell\mapsto \wt
F(v_\ell)$ and $\ell\mapsto \wt F'(v_\ell)$ for $\Ga$ on $\G X$ are
cohomologous for the definition given before the statement of Theorem
\ref{theo:varprinarbmet}. Since the existence and uniqueness of an
equilibrium state depends only on the cohomology class of the
bounded\footnote{hence with integrable negative part for any
  probability measure} potentials on $\G X$, the result follows.
\cqfd


%% file: meridonII.tex
\chapter{Random walks on weighted graphs of groups}
\label{sec:laplacian}

Let $\XX$ be a locally finite simplicial tree without terminal
vertices, and let $X=|\XX|_1$ be its geometric realisation.  Let $\Ga$
be a nonelementary discrete subgroup of $\Aut(\XX)$.

In Section \ref{subsec:quantumgraph}, given a (logarithmic) system
of conductances $c:\Ga\bs E\XX\ra\RR$, we define an operator ${\bf
  \Delta}_c$ on the functions defined on the set of vertices of the
quotient graph of groups $\Ga\dbs\XX$.
This operator is the infinitesimal generator\footnote{More precisely,
  it is the infinitesimal generator of the continuous time random
  process on the graph $\Ga\bs\XX$ whose co-called ``discrete
  skeleton'' or ``jump chain'' is the aforementioned random walk. The
  process waits an exponentially distributed time with parameter $1$
  at a vertex $x$, then instantaneously jumps along an edge $e$
  starting from $x$ with probability $i(e)e^{c(e)}/\deg_c(x)$. See for
  instance \cite[\S 2.1.2]{AldFil14}.}  of the random walk on
$\Ga\dbs\XX$ associated with the (normalised) exponential of this
system of conductances.  When $\Ga$ is torsion free and the system of
conductances vanishes, the construction recovers the standard Laplace
operator on the graph $\Ga\bs\XX$.

Under appropriate antireversibility assumptions on the system of
conductances $c$, using techniques of Sullivan and
Coornaert-Papadopoulos, we prove that the total mass of the Patterson
densities is a positive eigenvector for the operator ${\bf \Delta}_c$
associated with $c$.

In Section \ref{subsec:harmonicmeasure}, we study the nonsymmetric
nearest neighbour random walks on $V\XX$ associated with
antireversible systems of conductivities, and we show that the
Patterson densities are the harmonic measures of these random walks.

\section{Laplacian operators on  weighted graphs of groups}
\label{subsec:quantumgraph}

Let $\XX$ be a locally finite simplicial tree without terminal
vertices, and let $X=|\XX|_1$ be its geometric realisation.  Let $\Ga$
be a nonelementary discrete subgroup of $\Aut(\XX)$. Let $\wt c:E\XX
\ra \RR$ be a ($\Ga$-invariant) system of conductances for $\Ga$.  

We define $\wt c^{\,+}=\wt c$ and $\wt c^{\,-}:e\mapsto \wt
c(\overline{e})$, which is another system of conductances for
$\Ga$. Recall (see Section \ref{subsec:cond}) that $\wt c$ is
reversible (respectively antireversible) if $\wt c^{\,-}=\,\wt
c^{\,+}$ (respectively $\wt c^{\,-}=-\,\wt c^{\,+}$). For every $x\in
V\XX$, we define
$$
\deg_{\wt c^{\,\pm}}(x)= \sum_{e\in E\XX,\;o(e)=x} e^{\wt c^{\,\pm} (e)}\;.
$$
The quotient graph of groups $\Ga\dbs\XX$ is endowed with the quotient
maps $c^\pm: \Ga\bs E\XX\ra\RR$ of $\wt c^{\,\pm}$. Also note that the
quantity $\deg_{\wt c^{\,\pm}}(x)$ is constant on the $\Ga$-orbit of
$x$. Hence, it defines a map $\deg_{c^\pm}: \Ga\bs V\XX\ra
\;]0,+\infty[\,$.

On the vector space $\CC^{V\XX}$ of maps from $V\XX$ to $\CC$, we
consider the operator $\gls{laplacianconductanceup}$, called the {\it
  (weighted) Laplace operator}\index{Laplacian} of
$(\XX,c^\pm)$,\footnote{or on $\XX$ associated with the system of
  conductances $\wt c^{\,\pm}$} defined by setting, for all $f\in
\CC^{V\XX}$ and $x\in V\XX$,
\begin{equation}\label{eq:defiLaplacianup}
\Laplacian_{\wt c^{\,\pm}} f(x)=\frac{1}{\deg_{\wt c^{\,\pm}}(x)}\; 
\sum_{e\in E\XX,\;o(e)=x} e^{\wt c^{\,\pm}(e)} \big(f(x)-f(t(e))\big)\;.
\end{equation}
This is the standard Laplace operator\footnote{See for example
  \cite{Cartier72} with the opposite choice of the sign, or
  \cite{ChuGriYau00}.} of a weighted graph for the weight function
$e\mapsto e^{\wt c^{\,\pm}(e)}$, except that usually one requires that
$\wt c(e)=\wt c(\bar e)$. Note that $p^\pm(e)=\frac{e^{\wt c^{\,\pm}(e)}}
{\deg_{\wt c^{\,\pm}}(o(e))}$ is a Markov transition kernel on the tree $\XX$,
see the following Section \ref{subsec:harmonicmeasure}.

The weighted Laplace operator $\Laplacian_{\wt c^{\,\pm}}$ is invariant under
$\Ga$: for all $f\in\CC^{V\XX}$ and $\ga\in\Ga$, we have
$$
\Laplacian_{\wt c^{\,\pm}} (f\circ\ga)=(\Laplacian_{\wt c^{\,\pm}} f)\circ\ga\;.
$$
In particular, this operator induces an operator on functions defined
on the quotient graph $\Ga\bs\XX$, as follows.

Let $(\YY,G_*)$ be a graph of finite groups. We denote\footnote{See
  Section \ref{subsec:trees} for the definition of the measure
  $\vol_{(\YY,G_*)}$ on the discrete set $V\YY$.} by $\LL^2(V\YY,
\vol_{(\YY,G_*)})$ the Hilbert space of maps $f:V\YY\ra \CC$ with
finite norm $\|f\|_{\vol}$ for the following scalar product:
$$
\langle f,g\rangle_{\vol}= 
\sum_{x\in V\YY}\;\frac{1}{|G_x|}\; f(x)\;\overline{g(x)}\;.
$$ 
We denote\footnote{See Section \ref{subsec:trees} for the definition
  of the measure $\Tvol_{(\YY,G_*)}$ on the discrete set $E\YY$.} by
$\LL^2(E\YY, \Tvol_{(\YY,G_*)})$ the Hilbert space of maps
$\phi:E\YY\ra \CC$ with finite norm $\|\phi\|_{\Tvol}$ for the
following scalar product:
$$
\langle \phi,\psi\rangle_{\Tvol}= \frac{1}{2}\;
\sum_{e\in E\YY}\;\frac{1}{|G_e|}\; \phi(e)\;\overline{\psi(e)}\;.
$$

Let $i:E\YY \ra
\NN-\{0\}$ be the index map $i(e) = [G_{o(e)}:G_e]$. For every
function $c:E\YY\ra\RR$, let $\deg_c:V\YY\ra\RR$ be the positive function
defined by
$$
\deg_c(x)=\sum_{e\in E\YY,\;o(e)=x} i(e)\;e^{c(e)}\;.
$$
The {\em Laplace operator}\index{Laplacian}\footnote{See for instance
  \cite{Morgenstern94SIAM} when $c=0$.}  of $(\YY,G_*,c)$ is the
operator $\gls{laplacianconductancedown}=\Laplacian_{\YY,G_*,c}$ on
$\Leb^2(V\YY,\vol_{\YY,G_*})$ defined by
$$
\Laplacian_c f:x\mapsto\frac{1}{\deg_c(x)}\; 
\sum_{e\in E\YY,\;o(e)=x} i(e)\,e^{c(e)} \big(f(x)-f(t(e))\big)\;.
$$

\brema \medskip (1) Let $(\YY,G_*)=\Ga\dbs\XX$ be a graph of finite
groups with $p: V\XX\ra V\YY=\Ga\bs V\XX$ the canonical
projection. Let $\wt c:E\XX\to\RR$ be a potential for $\Ga$ and let
$c:E\YY=\Ga\bs E\XX\to\RR$ be the map induced by $\wt c$.  An easy
computation shows that for all $f\in\CC^{V\YY}$ and $x\in V\YY$, we
have
$$
\Laplacian_c f(x)=\Laplacian_c \wt f(\wt x)\;
$$
if $\wt f=f\circ p:V\XX\ra\CC$ and $\wt x\in V\XX$ satisfies $p(\wt
x)=x$.

\medskip\noindent (2) For every $x\in V\YY$, let 
$$
i(x)=\sum_{e\in E\YY,\;o(e)=x} i(e)\;.
$$ 
Then $i(x)$ is the degree of any vertex of any universal cover of
$(\YY,G_*)$ above $x$. In particular, the map $i:V\YY\ra\RR$ is
bounded if and only if the universal cover of $(\YY,G_*)$ has
uniformly bounded degrees. When $c=0$, we denote the Laplace operator
by $\Laplacian=\Laplacian_{\YY,\,G_*}$ and for every $x\in V\YY$, we
have
$$
\Laplacian f(x)=\frac{1}{i(x)}\; 
\sum_{e\in E\YY,\;o(e)=x} i(e) \big(f(x)-f(t(e))\big)\;.
$$
We thus recover the Laplace operator of \cite{Morgenstern94SIAM} on
the edge-indexed graph $(\YY,i)$.
\erema

\bprop \label{prop:propriLaplacian}
Let $(\YY,G_*)$ be a graph of finite groups, whose map
$i:V\YY\ra\RR$ is bounded. Let $c:E\YY\ra\RR$ be a system of
conductances on $\YY$, and let $$p(e)= \frac{e^{c(e)}}{\deg_c(o(e))}$$ for
every $e\in E\YY$. The following properties hold.
\begin{enumerate}
\item The Laplace operator $\Laplacian_c:\LL^2(V\YY,\vol_{(\YY,G_*)})
  \ra \LL^2(V\YY,\vol_{(\YY,G_*)})$ is linear and bounded.
\item The map $d_c: \LL^2(V\YY,
  \vol_{(\YY,G_*)}) \ra \LL^2(E\YY,\Tvol_{(\YY,G_*)})$ defined by
$$
d_c(f):e\mapsto \sqrt{p(e)}\,\big(f(t(e))-f(o(e))\big)
$$
is linear and bounded, and its dual operator 
$$
d_c^*:\LL^2(E\YY,\Tvol_{(\YY,G_*)}) \ra \LL^2(V\YY,\vol_{(\YY,G_*)})
$$ 
is given by
$$
d_c^*(\phi):x\mapsto \sum_{e\in
  E\YY,\;o(e)=x}\;\frac{i(e)}{2}\Big(\sqrt{p(\ov e)}\;\phi(\ov e)
-\sqrt{p(e)}\;\phi(e)\Big)\;.
$$
\item
Assume that $c$ is reversible and that the map $\deg_c:V\YY\ra\RR$ is
constant. Then
$$
\Laplacian_c=d_c^*d_c\;.
$$
In particular, $\Laplacian_c$ is self-adjoint and nonnegative.
\end{enumerate}
\eprop

\dem By the assumptions, there exists $M\in\NN$ such that $i(x)\leq M$
for every $x\in V\YY$, and hence $i(e)\leq M$ for every $e\in
E\YY$. Note that $i(e)= \frac{|G_{o(e)}|}{|G_e|}$, that
$G_{e}=G_{\overline{e}}$ and that $p(e)\leq 1$ .

\medskip\noindent (1) For every $f\in \LL^2(V\YY,\vol_{(\YY,G_*)})$,
using in the following computations

$\bullet$~ the Cauchy-Schwarz inequality for the first inequality,

$\bullet$~ the fact that for every $x\in V\YY$, we have
$$
\sum_{e\in E\YY,\;o(e)=x} i(e)^2\,p(e)^2
\leq M^2\sum_{e\in E\YY,\;o(e)=x} p(e)
\leq M^2\sum_{e\in E\YY,\;o(e)=x} i(e)p(e)=M^2
$$
for the second inequality,

$\bullet$~ the fact that $|G_{o(e)}|\geq |G_e|$ for every $e\in E\YY$
for the third inequality, and

$\bullet$~ the change of variable $e\mapsto \ov e$ in $\sum_{e\in
  E\YY} \frac{1}{|G_e|}|f(t(e))|^2$ (since $G_{\ov e}=G_e$) 
for the first equality on the fifth line of the computations,

\noindent we have
\begin{align*}
\|\Laplacian_c f\|_{\vol}^2&=\sum_{x\in V\YY}\;\frac{1}{|G_x|}\; 
\Big|\sum_{e\in E\YY,\;o(e)=x} i(e)\,p(e) \big(f(x)-f(t(e))\big)\Big|^2
\\ & \leq \sum_{x\in V\YY}\;\frac{1}{|G_x|}\; 
\Big(\sum_{e\in E\YY,\;o(e)=x} i(e)^2\,p(e)^2\Big)
\Big(\sum_{e\in E\YY,\;o(e)=x}  \big|f(x)-f(t(e))\big|^2\Big)
\\ & \leq 2\,M^2 \sum_{x\in V\YY}\;\frac{1}{|G_x|}\; 
\Big(\sum_{e\in E\YY,\;o(e)=x} \big(\big|f(x)\big|^2+\big|f(t(e))\big|^2\big)\Big)
\\ & \leq 2\,M^2 \; 
\sum_{e\in E\YY} \frac{1}{|G_e|}\Big(\big|f(o(e))\big|^2+\big|f(t(e))\big|^2\Big)
\\ & = 4\,M^2 \; 
\sum_{e\in E\YY} \frac{1}{|G_e|}\big|f(o(e))\big|^2 = 
4\,M^2 \; \sum_{x\in V\YY}\;\frac{1}{|G_x|}\; 
\sum_{e\in E\YY,\;o(e)=x} i(e)\,|f(x)|^2
\\ & = 4\,M^2 \; \sum_{x\in V\YY}\;\frac{i(x)}{|G_x|}\;|f(x)|^2
\leq 4\,M^3 \; \sum_{x\in V\YY}\;\frac{1}{|G_x|}\; 
|f(x)|^2=4\,M^3 \;\| f\|_{\vol}^2\;.
\end{align*}
Hence the linear operator $\Laplacian_c$ is bounded.

\medskip\noindent (2) For every $f\in \LL^2(V\YY,\vol_{(\YY,G_*)})$,
we have
\begin{align*}
\|d_c f\|_{\Tvol}^2&=\frac{1}{2}\;\sum_{e\in E\YY}\;
\frac{p(e)}{|G_e|}\; \big|f(t(e))-f(o(e))\big|^2
\\ & \leq\sum_{e\in E\YY}\;
\frac{1}{|G_e|}\; \Big(\big|f(t(e))\big|^2+\big|f(o(e))\big|^2\Big)
=2\sum_{e\in E\YY}\;
\frac{1}{|G_e|}\;\big|f(o(e))\big|^2
\\ & =2\sum_{e\in E\YY}\;
\frac{i(e)}{|G_{o(e)}|}\;\big|f(o(e))\big|^2=2\sum_{x\in V\YY}\;
\frac{i(x)}{|G_{x}|}\;|f(x)|^2\leq 2\,M\|f\|_{\vol}^2\;.
\end{align*}
Hence the linear operator $d_c$ is bounded.

For all $f\in \LL^2(V\YY,\vol_{(\YY,G_*)})$ and $\phi\in
\LL^2(E\YY,\Tvol_{(\YY,G_*)})$, using again the change of variable
$e\mapsto \ov e$, we have
\begin{align*}
\langle \phi , d_c f \rangle_{\Tvol} & = \frac{1}{2}
\sum_{e\in E\YY}\; \frac{1}{|G_{e}|}\;
\sqrt{p(e)}\; \phi(e) \;\overline{\big(f(t(e))-f(o(e)\big)}
\\ & =\frac{1}{2}
\bigg(\sum_{e\in E\YY}\; \frac{\sqrt{p(\ov e)}}{|G_{e}|}\;
\phi(\ov e) \;\overline{f(o(e))}-
\sum_{e\in E\YY}\; \frac{\sqrt{p(e)}}{|G_{e}|}\;
\phi(e) \;\overline{f(o(e))}\bigg)
\\ & =\sum_{x\in V\YY}\;\frac{1}{|G_x|}\;
\sum_{e\in E\YY,\;o(e)=x}\;\frac{i(e)}{2}\Big(\sqrt{p(\ov e)}\;\phi(\ov e)
-\sqrt{p(e)}\;\phi(e)\Big)\;\overline{f(x)}\;.
\end{align*}
This gives the formula for $d_c^*$.

\medskip\noindent (3) Let $f,g\in \LL^2(V\YY,\vol_{(\YY,G_*)})$.  Note
that $p(e)=p(\overline{e})$\footnote{This is the usual reversibility
  requirement for the corresponding Markov chain.} by the
reversibility of $c$ and the fact that $\deg_c$ is constant. Hence, by
developping the products in the first line and by making the change of
variable $e\mapsto \ov e$ in half the values, we have
\begin{align*}
\langle d_cf,d_cg\rangle_{\Tvol} & = 
\frac{1}{2}\;\sum_{e\in E\YY}\; \frac{1}{|G_{e}|}\;
p(e) \big(f(t(e)) - f(o(e))\big)\big(\,\overline{g(t(e))-g(o(e))}\,\big)
\\ & = \sum_{e\in E\YY}\; \frac{i(e)}{|G_{o(e)}|}\;
p(e) \big(f(o(e))\;\overline{g(o(e))}-f(t(e))\;\overline{g(o(e))}\,\big)
\\ & = \sum_{x\in V\YY}\; \frac{1}{|G_x|}\; \sum_{e\in E\YY,\;o(e)=x} 
i(e)\;p(e) \big(f(x)-f(t(e))\big)\;\overline{g(x)}
\\ & = \langle\Laplacian_c f , g \rangle_{\vol}\;.
\end{align*}
This proves the last claim in Proposition \ref{prop:propriLaplacian}.
\cqfd

\medskip The following result is an extension to antireversible
systems of conductances of \cite[Prop.~3.3]{CooPap96} (who treated the
case of zero conductances), which is a discrete version of Sullivan's
analogous result for hyperbolic manifolds (see \cite{Sullivan87}).
Let $\wt F_c:T^1X\ra \RR$ be the potential for $\Ga$ associated with
$\wt c$,\footnote{See Section \ref{subsec:cond}.} so that $(\wt
F_c)^\pm=\wt F_{c^\pm}$, and let $\delta_c$ be their common critical
exponent. Let $C^\pm: \partial_\infty X\times V\XX\times V\XX\ra \RR$
be the associated Gibbs cocycles.  Let $(\mu^\pm_x)_{x\in V\XX}$ be
two Patterson densities on $\partial_\infty X$ for the pairs
$(\Ga,F_{c^\pm})$.

\bprop \label{prop:lambdaharm} Assume that $\XX$ is $(q+1)$-regular,
that the system of conductances $\wt c$ is antireversible and that
the map $\deg_{\wt c^{\,\pm}}:V\XX\ra\RR$ is constant with value
$\kappa^\pm$. Then the total mass $\phi_{\mu^\pm} :x\mapsto
\|\mu_x^\pm\|$ of the Patterson density is a positive eigenvector
associated with the eigenvalue
$$
1-\frac{e^{\delta_c}+q e^{-\delta_c}}{\kappa^\pm}\,.
$$ 
for the Laplace operator $\Laplacian_{\wt c^{\,\pm}}$ on $\CC^{V\XX}$.
\eprop

\dem Note that the function $\wt c^{\,\pm}:E\XX\ra\RR$ is bounded,
since $e^{\wt c^{\,\pm}(e)}\leq \deg_{\wt c^{\,\pm}}(o(e))=\kappa^\pm$
for every $e\in E\XX$. Hence $(\wt F_c)^\pm=\wt F_{c^\pm}$ is bounded
by its definition in Section \ref{subsec:cond}. Since $\XX$ is
$(q+1)$-regular, the critical exponent $\delta_\Ga$ is finite and
hence the critical exponent $\delta_c=\delta_{\Ga,\,F_{c^\pm}}$ is
finite by Lemma \ref{lem:proprielemcritexpo} (6).  Since
$$
\phi_{\mu^\pm}(x)=\int_{\partial_\infty X}d\mu_x^\pm=
\int_{\partial_\infty X}\;e^{-C^\pm_\xi(x,\,x_0)}\;d\mu_{x_0}^\pm\,,
$$
by Equation \eqref{eq:quasinivarPatdens} and by linearity, we only
have to prove that for every fixed $\xi\in \partial_\infty X$  the
map 
$$
f:x\mapsto e^{-C^\pm_\xi(x,\,x_0)}
$$ 
is an eigenvector with eigenvalue $1-\frac{e^{\delta_c}+q
  e^{-\delta_c}}{\kappa^\pm}$ for $\Laplacian_{\wt c^{\,\pm}}$.

For every $e\in E\XX$, recall\footnote{See Section
  \ref{subsec:trees}.} that $\partial_e\XX$ is the set of points at
infinity of the geodesic rays in $\XX$ whose initial edge is $e$.  By
Equation \eqref{eq:changemoinsplus} and by the definition of the
potential associated with a system of conductances\footnote{See
  Proposition \ref{prop:integpotconduct} with the edge length map
  $\lambda$ constant equal to $1$.}, for all $e\in E\XX$ and $\eta\in
\partial_e\XX$, since $t(e)\in [o(e),\eta[$ (independently of the
choice of sign $\pm$), we have
$$
C^\pm_\eta(t(e),o(e))=\int_{o(e)}^{t(e)} (\wt F_{c^\pm}-\delta_c) 
=\wt c^{\,\pm}(e)-\delta_c\;.
$$
Thus if $\xi\in\partial_e\XX$, we have
$$
f(t(e))=e^{-C^\pm_\xi(t(e),\,o(e))-C^\pm_\xi(o(e),\,x_0)}=
e^{-\wt c^{\,\pm}(e)+\delta_c}\;f(o(e))\;,
$$
and otherwise
$$
f(t(e))=e^{C^\pm_\xi(t(\ov{e}),\,o(\ov{e}))-C^\pm_\xi(o(e),\,x_0)}=
e^{\wt c^{\,\pm}(\ov{e})-\delta_c}\;f(o(e))\;.
$$
For every $x\in V\XX$, let $e_\xi$ be the unique edge of $\XX$ with
origin $x$ such that $\xi\in\partial_{e_\xi}\XX$.  Then, 
\begin{align*}
\Laplacian_{\wt c^{\,\pm}}f (x)&=
f(x)-\frac{1}{\deg_{\wt c^{\,\pm}}(x)} \sum_{o(e)=x}e^{\wt c^{\,\pm}(e)} f(t(e))\\ & =
f(x)-\frac{1}{\kappa^\pm} e^{\wt c^{\,\pm}(e_\xi)} f(t(e_\xi)) -
\frac{1}{\kappa^\pm} \sum_{e\neq e_\xi,\,o(e)=x}e^{\wt c^{\,\pm}} f(t(e)) \\ &=
\Big(1-\frac{e^{\delta_c}}{\kappa^\pm}  -
\frac{q\,e^{-\delta_c}}{\kappa^\pm}  \Big)f(x) \;.
\end{align*}
This proves the result.
\cqfd 

\medskip Note that the antireversibility of the potential is used in
an essential way in order to get the last equation in the proof of
Proposition \ref{prop:lambdaharm}.

\section{Patterson densities as harmonic measures for simplicial \\ trees}
\label{subsec:harmonicmeasure}

In this Section, we define and study a Markov chain on the set of
vertices of a simplicial tree endowed with a discrete group of
automorphisms and with an appropriate system of conductances, such
that the associated (nonsymmetric, nearest neighbour) random walk
converges almost surely to points in the boundary of the tree, and we
prove that the Patterson densities, once normalised, are the
corresponding harmonic measures.  We thereby generalise the zero
potential case treated in \cite{CooPap96}, which is also a
special case of \cite{ConMuc07a} when $X$ is a tree under the
additional restriction that the discrete group is cocompact. For other
connections between harmonic measures and Patterson measures, we refer
for instance to \cite{ConMuc07a,BlaHaiMat11,Tanaka14,GouMatMau15} and
their references.

\medskip Let 
$\XX$ be a $(q+1)$-regular simplicial tree, with $q\geq 2$. Let $\Ga$
be a nonelementary discrete subgroup of $\Aut(\XX)$. Let $\wt c:E\XX
\ra \RR$ be an antireversible system of conductances for $\Ga$, such
that the associated map $\deg_{\wt c}:V\XX \ra \RR$ on the vertices of $\XX$
is constant. Let $(\mu_x)_{x\in V\XX}$ be a Patterson density for
$(\Ga,F_c)$, where $F_c$ is the potential associated with $c$. We
denote by $\phi_\mu:x\mapsto \|\mu_{x}\|$ the associated total mass
function on $V\XX$.

\medskip We start this Section by recalling a few facts about discrete
Markov chains, for which we refer for instance to
\cite{Revuz84,Woess94}. A {\em state space}\index{state space} is a
discrete and countable set $I$. A {\em transition
  kernel}\index{transition kernel} on $I$ is a map $p:I\times I\ra
[0,1]$ (considered as a square matrix with coefficients in $[0,1]$,
and with row and column indices in $I$) such that for every $x\in I$,
$$
\sum_{y\in I}\;p(x,y)=1\;.
$$ 
Let $\lambda$ be a probability measure on $I$. A (discrete) {\em
  Markov chain}\index{Markov!chain} on a state space $I$ with initial
distribution $\lambda$ and transition kernel $p$ is a sequence
$(Z_n)_{n\in\NN}$ of random variables with values in $I$ such that for
all $n\in\NN$ and $x_0,\dots, x_{n+1}\in I$, the probability of events
$\PP$ satisfies
\begin{enumerate}
\item $\PP[Z_0=x_0]=\lambda(\{x_0\})$,
\item $\PP[Z_{n+1}=x_{n+1} \mid Z_0=x_0,Z_1=x_1,\dots, Z_{n}=x_{n}]=$\\
  $\PP[Z_{n+1}=x_{n+1} \mid Z_{n}=x_{n}]=p(x_n,x_{n+1})$.
\end{enumerate}
The associated {\it random walk}\index{random walk} consists in
choosing a point $x_0$ in $I$ with law $\lambda$, and by induction,
once $x_n$ is constructed, in choosing $x_{n+1}$ in $I$ with
probability $p(x_n,x_{n+1})$.  Note that
$$
\PP[Z_0=x_0,Z_1=x_1,\dots, Z_{n}=x_{n}]
=\lambda(\{x_0\})\,p(x_0,x_1)\dots p(x_{n-1},x_n)\;.
$$
When the initial distribution $\lambda$ is the unit Dirac mass
$\Delta_x$ at $x\in I$, the Markov chain is then uniquely determined
by its transition kernel $p$ and by $x$, and is denoted by
$(Z^x_n)_{n\in\NN}$.

For every $n\in\NN$, we denote by $p^{(n)}$ the iterated matrix
product of the transition kernel $p$: we define $p^{(0)}(x,y)$ to be
the Kronecker symbol $\delta_{x,y}$ for all $x,y\in I$, and by
induction $p^{(n+1)}=p \;p^{(n)}$, that is, for all $x,z\in I$,
$$
p^{(n+1)}(x,z)=\sum_{y\in I} p(x,y)p^{(n)}(y,z)\;.
$$
Note that
$$
p^{(n)}(x,y)=\PP[Z^x_n=y]
$$
is the probability for the random walk starting at time $0$ from $x$
of being at time $n$ at the point $y$.  The {\em Green
  kernel}\index{Green!kernel} of $p$ is the map $G_p$ from $I\times I$
to $[0,+\infty ]$ defined by
$$
(x,y)\mapsto G_p(x,y)=\sum_{n\in\NN}\; p^{(n)}(x,y)\;,
$$
and its {\em Green function}\index{Green!function} is the following
power series in the complex variable $z$~:
$$
G_p(x,y\mid z)= \sum_{n\in\NN} \;p^{(n)}(x,y) \;z^n\;.
$$
Recall that if $G_p(x,y)\neq 0$ for all $x,y\in I$, then the random
walk is {\em recurrent}\index{recurrent}\footnote{that is,
  $\card\{n\in\NN\;:\; Z_n^x=y\}=\infty$ for every $y\in I$ (or
  equivalently, there exists $y\in I$ such that $\card\{n\in\NN\;:\;
  Z_n^x=y\}=\infty$)} if $G_p(x,y)=\infty$ for any (hence all)
$(x,y)\in I\times I$, and {\em transient}\index{transient} otherwise.
Note that, using again matrix products of $I\times I$ matrices,
\begin{equation}\label{eq:matrixproduct}
G_p=\operatorname{Id} +p\; G_p\;.
\end{equation}

\bigskip 
We will from now on consider as state space the set $V\XX$ of vertices
of $\XX$. If a Markov chain $(Z^x_n)_{n\in\NN}$ starting at time $0$
from $x$ converges almost surely in $V\XX \cup \partial_\infty X$ to
a random variable $Z^x_\infty$, the law of $Z^x_\infty$ is called the
{\em harmonic measure}\index{harmonic measure} (or {\em hitting
  measure} on the boundary) associated with this Markov chain, and is
denoted by
$$
\nu_x=(Z^x_\infty)_*(\PP)\;.
$$
Note that $\nu_x$ is a probability measure on $\partial_\infty X$.

For instance, the transition kernel of the simple nearest neighbour
random walk on $\XX$ is defined by taking as transition kernel the map
$\underline{p}$ where
$$
\underline{p}(x,y)=\frac 1{q+1}\;A(x,y)
$$ 
for all $x,y\in V\XX$, with $A:V\XX\times V\XX\ra \{0,1\}$  the
{\em adjacency matrix}\index{adjacency matrix} of the tree $\XX$,
defined by $A(x,y)=1$ for any two vertices $x,y$ of $\XX$ that are
joined by an edge in $\XX$ and $A(x,y)=0$ otherwise. We denote by
$$
\underline{\green{}}(x,y\mid z)=
\sum_{k\in\NN}\;\underline{p}^{(n)}(x,y)\;z^n\,.
$$ 
the Green function of $\underline{p}$, whose radius of convergence is
$\underline r=\frac{q+1}{2\sqrt q}$ and which diverges at
$z=\underline r$, see for example \cite{Woess94},
\cite[Ex.~9.82]{Woess09}, \cite[\S 6.3]{LyoPerBook}.

\medskip 
The antireversible system of conductances $\wt c:E\XX\to\RR$
defines a cocycle on the set of vertices of $\XX$, as follows.  For
all $u,v\in V\XX$, let $c(u,v)=0$ if $u=v$ and otherwise
let $$c(x,y)=\sum_{i=1}^n \wt c(e_i)\,,$$ where $(e_1,e_2,\dots, e_m)$
is the geodesic edge path in $\XX$ from $u=o(e_1)$ to $v=t(e_n)$.

\blemm \label{lem:propricocyclec}\mbox{}
\begin{enumerate}
\item For every edge path $(e'_1,e'_2,\dots, e'_{n'})$ from $u$ to
  $v$, we have
$$
c(u,v)= \sum_{i=1}^{n'} \wt c(e'_i)\;.
$$
\item The map $c:V\XX\times V\XX\ra \RR$ has the following
  cocycle property: for all $u,v,w\in V\XX$,
$$
c(u,v)+c(v,w)= c(u,w) \;\;\;{\rm and~hence}\;\;\;c(v,u)=-c(u,v)\;.
$$
\item We have $c(u,v)=\int_u^v\wt F_c$.
\item For all $\xi\in\partial_\infty X$ and $u,v\in V\XX$, if
  $C^c_\cdot(\cdot, \cdot)$ is the Gibbs cocycle associated with $\wt
  F_c$, we have
$$
C^c_\xi(u,v)= c(v,u)+\delta_c\beta_\xi(u,v)\;.
$$
\end{enumerate}
\elemm

\dem (1) Since $\XX$ is a simplicial tree, any nongeodesic edge path
from $u$ to $v$ has a back-and-forth on some edge, which contributes
to $0$ in the sum defining $c(x,y)$ by the antireversibility
assumption on the system of conductances. Therefore, by induction, the
sum in Assertion (1) indeed does not depend on the choice of the edge
path from $u$ to $v$.

Assertion (2) is immediate from Assertion (1).  Assertion (3) follows
from the definition of $c(\cdot,\cdot)$ by Proposition
\ref{prop:integpotconduct}.

\medskip\noindent (4) For every $\xi\in\partial_\infty X$, if $p\in
V\XX$ is such that $[u,\xi[\;\cap \,[v,\xi[\;=[p,\xi[\,$, then using
Equation \eqref{eq:cocycletreecase} and Assertions (3) and (2), we
have
\begin{align*}
C^c_\xi(u,v)&=\int_v^p(\wt F_c-\delta_c)-\int_u^p(\wt F_c-\delta_c)
= c(v,p)-c(u,p)+\delta_c\,\beta_\xi(u,v)\\ &=
c(v,u)+\delta_c\,\beta_\xi(u,v)\;.\;\;\;\; \Box
\end{align*}

\bigskip Let 
$$
\kappa_c=\frac{q+1}{e^{\delta_c}+q\,e^{-\delta_c}}\;,
$$
which belongs to $]0,\frac{q+1}{2\sqrt q}]\,$, with $\kappa_c=
\frac{q+1}{2\sqrt q}$ if and only if $e^{\delta_c}=\sqrt{q}$. Note
that this constant $\kappa_c$ is less than the radius of convergence
$\underline{r}= \frac{q+1}{2\sqrt q}$ of the Green function
$\underline{\green{}}(x,y\mid z)$ if and only if $\delta_c\neq
\frac{1}{2}\,\ln q$. The computation (due to Kesten) of the Green
function of $\underline{p}$ is well known, and gives the following
formula, see for instance \cite[Prop.~3.1]{CooPap96}:
If $\delta_c\neq \frac{1}{2}\,\ln q$, then there exists $\alpha>0$
such that\footnote{We actually have $\alpha= \frac{q+1} {e^{\delta_c}+
    (q-1)\,e^{-\delta_c}}$.} for all $x,y\in V\XX$
\begin{equation}\label{eq:calcgreensimpleRW}
\underline{\green{}}(x,y\mid \kappa_c)=
\alpha\; e^{-\delta_c\,d(x,\,y)}\;.
\end{equation}

We now define the transition kernel $p_c$ {\em associated
  with\footnote{The transition kernel also depends on the choice of
    the Patterson density if $\Ga$ is not of divergence type.} the
  (logarithmic) system of conductances} $c$ by, for all $x,y\in V\XX$,
$$
p_c(x,y)= \kappa_c\;
\frac{\phi_\mu(y)}{\phi_\mu(x)}\;e^{c(x,\,y)}\;\underline{p}(x,y)\;.
$$

From now on, we denote by $(Z_n^x)_{n\in\NN}$ the Markov chain with
initial distribution $\Delta_x$ and transition kernel $p_c$.

\blemm \label{lem:propripotentieldec}\mbox{}
\begin{enumerate}
\item The map $p_c$ is a transition kernel on $V\XX$.
\item The Green kernel $\green_c=G_{p_c}$ of $p_c$ is
\begin{equation}\label{eq:greenconduc}
\green_c(x,y) = e^{c(x,\,y)}\;\frac{\phi_\mu(y)}{\phi_\mu(x)}\;\;
\underline{\green{}}(x,y\mid \kappa_c)\;.
\end{equation}
In particular, the Green kernel of $p_c$ is finite if $\delta_c\neq
\frac{1}{2}\,\ln q$.
\item Assume that $\delta_c\neq \frac{1}{2}\,\ln q$. For all $x,y,z\in
  V\XX$, we have
$$
\frac{\phi_\mu(y)\;\green_c(y,z)}{\phi_\mu(x)\;\green_c(x,z)}=
e^{c(y,\,x)+\delta_c(d(x,\,z)-d(y,\,z))}\;.
$$

If furthermore $z\notin [x,y[\,$, then, for every $\xi\in
\OOO_x(z)$,\footnote{Recall that given $x,z\in V\XX$, the shadow
$\OOO_x(z)$ of $z$ seen from $x$ is the set of points at
infinity of the geodesic rays from $x$ through $z$.}

$$
\frac{\phi_\mu(y)\;\green_c(y,z)}{\phi_\mu(x)\;\green_c(x,z)}=
e^{C^c_\xi(x,\,y)}\;.
$$
\end{enumerate}
\elemm

\dem (1) By the proof of Proposition \ref{prop:lambdaharm}, the
positive function $\phi_\mu$ is an eigenvector with eigenvalue
$e^{\delta_c}+ q\, e^{-\delta_c}$ for the operator $$f\mapsto
\{x\mapsto \sum_{e\in E\XX,\;o(e)=x} e^{\wt c(e)} \,f(t(e))\}\,.$$
Since $\underline{p}(o(e),t(e))=\frac{1}{q+1}$ for every $e\in E\XX$,
we hence have
\begin{align*}
&\sum_{y\in V\XX}p_c(x,y) \;= \sum_{e\in E\XX,\;o(e)=x} p_c(x,t(e))
\\ =\;&\frac{1+q}{(e^{\delta_c}+q\,e^{-\delta_c})\;\phi_\mu(x)}
\sum_{e\in E\XX,\;o(e)=x} e^{\wt c(e)}\,\phi_\mu(t(e))\;
\underline{p}(x,t(e)) =1\;.
\end{align*}

\noindent(2) Let us first prove that for all $x,y\in V\XX$ and
$n\in\NN$, we have
\begin{equation}\label{eq:probatimesn}
p_c^{(n)}(x,y)= 
(\kappa_c)^n\;
\frac{\phi_\mu(y)}{\phi_\mu(x)}\;e^{c(x,\,y)}\;
\underline{p}^{(n)}(x,y)\;.
\end{equation}
Indeed, by the cocycle property of $c(\cdot,\cdot)$ and by a
telescopic cancellation argument, we have
\begin{align*}
 p_c^{(n)}&(x,y) \\
=\;&  \sum_{x_1,\dots,\, x_{n-1}\in V\XX}
  p_c(x,x_1)\;p_c(x_1,x_2)\dots p_c(x_{n-2},x_{n-1})\;p_c(x_{n-1},y)\\ 
=\;& (\kappa_c)^n\;\frac{\phi_\mu(y)}{\phi_\mu(x)}\;e^{c(x,\,y)}
\sum_{x_1,\dots, \,x_{n-1}\in V\XX}
\underline{p}(x,x_1)\;\underline{p}(x_1,x_2)
\dots \;\underline{p}(x_{n-2},x_{n-1})\;\underline{p}(x_{n-1},y)\\ =\;&
(\kappa_c)^n\;\frac{\phi_\mu(y)}{\phi_\mu(x)}\;e^{c(x,\,y)}\;
\underline{p}^{(n)}(x,y)\;. 
\end{align*}

\medskip Equation \eqref{eq:greenconduc} follows from Equation
\eqref{eq:probatimesn} by summation on $n$. As we have already seen,
$\kappa_c<\underline{r}$ if and only if $\delta_c\neq \frac{1}{2}\,\ln
q$. The last claim of Assertion (2) follows.

\medskip\noindent (3) Let $x,y,z\in V\XX$. Using (twice) Assertion
(2), the cocycle property of $c$ and (twice) Equation
\eqref{eq:calcgreensimpleRW}, we have
\begin{align*}
\frac{\green_c(y,z)}{\green_c(x,z)}&=
\frac{e^{c(y,\,z)}\;\phi_\mu(z)\;\phi_\mu(x)\;\;
\underline{\green{}}(y,z\mid \kappa_c)}
{e^{c(x,\,z)}\;\phi_\mu(y)\;\phi_\mu(z)\;\;
\underline{\green{}}(x,z\mid \kappa_c)} =
e^{c(y,\,x)}\;\frac{\phi_\mu(x)}{\phi_\mu(y)}\;
\frac{\alpha\; e^{-\delta_c\,d(y,\,z)}}{\alpha\; e^{-\delta_c\,d(x,\,z)}}
\\ &=
\frac{\phi_\mu(x)}{\phi_\mu(y)}\;e^{c(y,\,x)+\delta_c(d(x,\,z)-d(y,\,z))}\;.
\end{align*}
This proves the first claim of Assertion (3).  Under the additional
assumptions on $x,y,z,\xi$, we have
$$
\beta_\xi(x,y)=d(x,z)-d(y,z)\;.
$$
The last claim of Assertion (3) hence follows from Lemma
\ref{lem:propricocyclec} (4).  \cqfd

\medskip Using the criterion that the random walk starting from a
given vertex of $\XX$ with transition probabilities $p_c$ is transient
if and only if the Green kernel $\green_c(x,y)$ of $p_c$ is finite
(for any, hence for all, $x,y\in V \XX$), Lemma
\ref{lem:propripotentieldec} (2) implies that if $\delta_c\neq
\frac{1}{2}\,\ln q$, then $(Z_n^x)_{n\in\NN}$ almost surely leaves
every finite subset of $V\XX$. The following result strengthens this
remark.

\bprop \label{prop:harmmesexist}
If $\delta_c\neq \frac{1}{2}\,\ln q$, then for every $x\in
V\XX$, the Markov chain $(Z_n^x)_{n\in\NN}$ (with initial distribution
$\Delta_x$ and transition kernel $p_c$) converges almost surely in
$V\XX\cup \partial_\infty X$ to a random variable with values in
$\partial_\infty X$. In particular the harmonic measure $\nu_x$ of
$(Z_n^x)_{n\in\NN}$ is well defined if $\delta_c\neq \frac{1}{2}\,\ln q$.
\eprop

\dem Since $\XX$ is a tree, if $(x_n)_{n\in\NN}$ is a sequence in
$V\XX$ such that $d(x_n,x_{n+1})=1$ for every $n\in\NN$ and which does
not converge to a point in $\partial_\infty X$, then there exists a
point $y$ such that this sequence passes infinitely often through $y$,
that is, $\{n\in\NN\;:\;x_n=y\}$ is infinite. The result then follows
from the fact that the Markov chain $(Z_n^x)_{n\in\NN}$ is transient
since $\delta_c\neq \frac{1}{2}\,\ln q$.  
\cqfd

\medskip
The following result, generalising \cite[Theo.~4.5]{CooPap96} when
$\wt c=0$, says that the Patterson measures associated with the system of
conductances $\wt c$, once renormalised to probability measures, are
exactly the harmonic measures for the random walk with transition
probabilities $p_c$.

\btheo \label{theo:pattersonharmonic} Let $(\XX,\Ga,\wt c,
(\mu_x)_{x\in V\XX})$ be as in the beginning of Section
\ref{subsec:harmonicmeasure}. If $\delta_c\neq \frac{1}{2}\,\ln q$,
then for every $x\in V \XX$, the harmonic measure of the Markov chain
$(Z_n^x)_{n\in\NN}$ is
$$
\nu_x=\frac{\mu_x}{\|\mu_x\|}\;.
$$
\etheo

\dem We fix $x\in V\XX$. For every $n\in\NN$, we denote by $S(x,n)$
and $B(x,n)$ the sphere and (closed) ball of centre $x$ and radius $n$
in $V\XX$, and we define two maps $f_1,f_2: V\XX\ra \RR$ with finite
support by
$$
f_1(z)=\frac{\mu_x(\OOO_x(z))}{\|\mu_x\|\;\green_c(x,z)}
\;\;\;{\rm and}\;\;\;
f_2(z)=\frac{\nu_x(\OOO_x(z))}{\green_c(x,z)}
$$
if $z\in S(x,n)$, and $f_1(z)=f_2(z)=0$ otherwise. Let us prove that
$f_1=f_2$ for every $n\in\NN$. Since $\{\OOO_x(z)\;:\;z\in V\XX\}$
generates the Borel $\sigma$-algebra of $\partial _\infty X$, this
proves that the Borel measures $\nu_x$ and $\frac{\mu_x}{\|\mu_x\|}$
coincide.

We will use the following criterion. For all maps $G:V\XX\times
V\XX\ra \RR$ and $f:V\XX\ra \RR$ such that $f$ has finite support, let
us again denote by $G\; f:V\XX\ra \RR$ the matrix product of the
square matrix $G$ and the column matrix $f$, defined by, for every
$y\in V\XX$,
$$
G\; f(y)=\sum_{z\in V\XX} G(y,z)f(z)\;.
$$

\blemm\label{lem:greenagree} For all $f,f':V\XX\to\RR$ with finite
support, if $\green_c\; f=\green_c\; f'$, then $f=f'$.  
\elemm

\dem
By Equation \eqref{eq:matrixproduct}, we have
$$
f'=\green_c\; f'-p_c\;\green_c\; f'
=\green_c\; f-p_c\green_c\; f=f\;.
\;\;\;\Box
$$

Let us hence fix $n\in\NN$ and prove that $\green_c\; f_1=
\green_c\; f_2$. Theorem \ref{theo:pattersonharmonic} then follows.

\bigskip \noindent {\bf Step 1: } For every $y\in B(x,n)$, since
$\{\OOO_x(z)\;:\;z\in S(x,n)\}$ is a Borel partition of $\partial
_\infty X$, by Equation \eqref{eq:quasinivarPatdens}, since
$z\notin[x,y[$ if $z\in S(x,n)$ and $y\in B(x,n)$, and by the second
claim of Lemma \ref{lem:propripotentieldec} (3), we have
\begin{align}
1 & =\frac{1}{\phi_\mu(y)}\;\int_{\partial _\infty X}d\mu_y=
\frac{1}{\phi_\mu(y)}\;\sum_{z\in S(x,n)} \int_{\OOO_x(z)}  
e^{-C^c_{\xi}(y,\,x)} \;d\mu_x(\xi)\nonumber\\ &=\frac{1}{\phi_\mu(y)}\;
\sum_{z\in S(x,n)}  \int_{\OOO_x(z)}  \;
\frac{\phi_\mu(y)\;\green_c(y,z)}{\phi_\mu(x)\;\green_c(x,z)}\;d\mu_x(\xi)
\nonumber\\ &= \sum_{z\in S(x,n)}  
\green_c(y,z)\;\frac{\mu_x(\OOO_x(z))}{\|\mu_x\|\;\green_c(x,z)}
=(\green_c \; f_1)(y)\;.\label{eq:computGcfun}
\end{align}

\bigskip \noindent {\bf Step 2: } For all $y,z \in V\XX$ such that
$z\notin[x,y[\,$, any random walk starting at time $0$ from $y$ and
converging to a point in $\OOO_x(z)$ goes through $z$. Let us denote
by $C_x(z)$ the set of vertices different from $z$ on the geodesic
rays from $z$ to the points in $\OOO_x(z)$. Partioning by the last
time the random walk passes through $z$, using the Markov property
saying that what happens before the random walk arrives at $z$ and
after it leaves $z$ are independent, we have
$$
\nu_y(\OOO_x(z))=\PP[Z^y_\infty\in \OOO_x(z)]=
\green_c(y,z)\;\PP[\forall\,n>0,\;Z^z_n\in C_x(z)]\;,
$$
so that 
\begin{equation}\label{eq:variatharmmes}
\frac{\nu_y(\OOO_x(z))}{\nu_x(\OOO_x(z))}=
\frac{\green_c(y,z)}{\green_c(x,z)}\;.
\end{equation}

\bigskip \noindent {\bf Step 3: } For every $y\in B(x,n)$, again since
$\{\OOO_x(z)\;:\;z\in S(x,n)\}$ is a Borel partition of $\partial
_\infty X$, and by Equation \eqref{eq:variatharmmes}, we have
\begin{equation}\label{eq:computGcfdeux}
1  =\|\nu_y\|= \sum_{z\in S(x,n)}  \nu_y(\OOO_x(z)) = \sum_{z\in S(x,n)}  
\green_c(y,z)\;\frac{\nu_x(\OOO_x(z))}{\green_c(x,z)}
=(\green_c \; f_2)(y)\;. 
\end{equation}

\bigskip
\noindent
{\bf Step 4: } By Steps $1$ and $3$, we have $(\green_c \; f_1)(y)=
(\green_c \; f_2)(y)$ for every $y\in B(x,n)$. Let now $y\in V\XX-
B(x,n)$. Define $y'\in S(x,n)$ as the point at distance $n$ from $x$
on the geodesic segment $[x,y]$. For every $z\in S(x,n)$, we have
$d(y',z)-d(y,z)=-d(y,y')$, which is independent of $z$.

\begin{center}
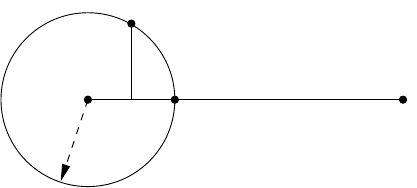
\end{center}

\noindent
Since $y'\in B(x,n)$, we have, as just said, $(\green_c \; f_1)(y')=
(\green_c \; f_2)(y')$.  Hence by the first claim of Lemma
\ref{lem:propripotentieldec} (3), we have
\begin{align*}
(\green_c \; f_1)(y) & = \sum_{z\in S(x,n)} \green_c(y,z)\;f_1(z)
\\ &= \sum_{z\in S(x,n)}  e^{c(y,\,y')+\delta_c(d(y',\,z)-d(y,\,z))}\;
\frac{\phi_\mu(y')}{\phi_\mu(y)}\;
\green_c(y',z)\;f_1(z)\\ &= e^{c(y,\,y')-\delta_cd(y,\,y')}\;
\frac{\phi_\mu(y')}{\phi_\mu(y)}\;(\green_c \; f_1)(y')
\\ &= e^{c(y,\,y')-\delta_cd(y,\,y')}\;\frac{\phi_\mu(y')}{\phi_\mu(y)}
\;(\green_c \; f_2)(y')= (\green_c \; f_2)(y)\;.
\end{align*}
This proves that $\green_c\; f_1= \green_c\; f_2$, thereby
concluding the proof of Theorem \ref{theo:pattersonharmonic}.  
\cqfd

\chapter{Skinning measures with potential on $\CAT(-1)$ spaces}
\label{sect:skinning}

In this Chapter, we introduce skinning measures as weighted
pushforwards of the Patterson-Sullivan densities associated with a
potential to the unit normal bundles of convex subsets of a $\CAT(-1)$
space. The development follows \cite{ParPau14ETDS} with modifications
to fit the present context.
 
Let $X,x_0,\Ga,\wt F$ be as in the beginning of Chapter
\ref{sect:measures}, and $\wt F^\pm$, $F^\pm$, $\delta=
\delta_{\Ga,\,F^\pm}<+\infty$ the associated notation.  Let
$(\mu^\pm_x)_{x\in X}$ be (normalised) Patterson densities on
$\partial_\infty X$ for the pairs $(\Ga,F^\pm)$.\footnote{See Section
  \ref{subsec:Pattersondens}.}

\section{Skinning measures}
\label{subsec:skinningmeasures}

Let $D$ be a nonempty proper closed convex subset of $X$.  The {\em
  outer skinning measure}%
\index{skinning measure}\index{measure!skinning} $\wt\sigma^+_D$ on
the outer normal bundle $\normalout D$ of $D$ and the {\em inner
  skinning measure}%
\index{skinning measure} $\wt\sigma^-_D$ on the inner normal bundle
$\normalin D$ of $D$ associated with the Patterson densities
$(\mu^\pm_x)_{x\in X}$ for $(\Ga,\wt F^\pm)$ are the measures
$\gls{skinningpm}=\wt\sigma^\pm_{D,\,F^\pm}$ on $\partial^1_\pm D$ 
defined by
\begin{equation}\label{eq:defiskin}
d\wt\sigma^\pm_D(\rho)  = 
e^{C^\pm_{\rho_\pm}(x_0,\,\rho(0))}\;d\mu^\pm_{x_0}(\rho_{\pm}) \;,
\end{equation}
where $\rho\in\normalpm D$, using the endpoint homeomorphisms
$\rho\mapsto \rho_\pm$ from $\normalpm D$ to $\partial_\infty X
-\partial_\infty D$, and noting that $\rho(0)=P_D(\rho_\pm)$ depends
continuously on $\rho_\pm$. 

When $\wt F=0$, the skinning measure has been defined by Oh and Shah
\cite{OhSha13} for the outer unit normal bundles of spheres,
horospheres and totally geodesic subspaces in real hyperbolic spaces. 
The definition was generalised in \cite{ParPau14ETDS} to the outer
unit normal bundles of nonempty proper closed convex sets in
Riemannian manifolds with variable negative curvature.

Note that the Gibbs measure is defined on the space $\G X$ of geodesic
lines, the potential is defined on the space $T^1X$ of germs at time
$t=0$ of geodesic lines, and since $\normalpm D$ is contained in
$\G_{\pm,\,0} X$ (see Section \ref{subsect:nbhd}), the skinning
measures are defined on the spaces $\G_{\pm,\,0} X$ of (generalised)
geodesic rays. In the manifold case, all the above spaces are
canonically identified with the standard unit tangent bundle, but in general,
the natural restriction maps $\G X\ra T^1X$ and $\G X\ra\G_{\pm,\,0}
X$ have infinite (though compact) fibers.

\brema\label{rem:skinrem}
(1) If $D=\{x\}$ is a singleton, then
\begin{equation}\label{eq:skinnpoint}
d\wt\sigma^\pm_D(\rho)=\;d\mu^\pm_{x}(\rho_{\pm})\;
\end{equation}
where $\rho$ is a geodesic ray starting (at time $t=0$) from $x$.

\medskip
\noindent 
(2) When the potential $\wt F$ is equal to $\wt F\circ \iota$ (in
particular when $F=0$), we have $C^-=C^+$, and we may (and we will)
take $\mu_x^-=\mu_x^+$ for all $x\in X$, hence $\iota_*\wt m_F=\wt
m_F$ and
$$
\wt \sigma^-_D=\iota_*\wt \sigma^+_D\;.
$$ More generally, if $\wt F$ is reversible, let $\wt G :T^1X\ra \RR$
be a continuous $\Ga$-invariant function such that, for every $\ell\in
\G X$, the map $t\mapsto \wt G(v_{\flow{t}\ell})$ is differentiable
and $\wt F^*(v_\ell)-\wt F(v_\ell)=\frac{d}{dt}_{\mid t=0}\wt
G(v_{\flow{t}\ell})$.  Furthermore assume that $X$ is an $\RR$-tree or
that $G$ is uniformly continuous (for instance
H\"older-continuous). Then we have, for every $\rho\in\partial_+^1D$,
denoting by $\wh\rho\in\G X$ any extension of $\rho$ to a geodesic
line in $X$,
$$
d\iota_*\wt \sigma^-_D(\rho)=
e^{-\wt G(v_{\wh\rho})}\;d\wt \sigma^+_D(\rho)\;.
$$
Indeed, for all $x,y\in X$ and $\xi\in\partial_\infty X$, let
$\ell_{x,\,\xi}$ be any geodesic line with footpoint
$\ell_{x,\,\xi}(0) =x$ and positive endpoint $(\ell_{x,\,\xi})_+=\xi$.
Then by Remark \ref{rem:cohomologuustuff}, we have
$$
C^{-}_{\xi}(x,y)-C^+_{\xi}(x,y)=
\wt G(v_{\ell_{x,\,\xi}})-\wt G(v_{\ell_{y,\,\xi}})\;,
$$
and we may (and we will) take
$$
d\mu^{-}_x(\xi)=e^{-\wt G(v_{\ell_{x,\,\xi}})}\;d\mu^+_x(\xi)\;.
$$
Hence for every $\rho\in\partial_+^1D$, we have
\begin{align*}
  d\wt \sigma^-_D(\iota\rho) &=
  e^{C^-_{(\iota\rho)_-}(x_0,\,(\iota\rho)(0))}\;d\mu^-_{x_0}((\iota\rho)_-)
  =e^{C^-_{\rho_+}(x_0,\,\rho(0))}\;d\mu^-_{x_0}(\rho_+)\\ & =
  e^{C^+_{\rho_+}(x_0,\,\rho(0)) -\wt G(v_{\ell_{\ell(0),\,\rho_+}})}\;d\mu^+_{x_0}(\rho_+)
  =e^{-\wt G(v_{\ell_{\rho(0),\,\rho_+}})}\;d\wt \sigma^+_D(\rho)\;.
\end{align*}

\medskip
\noindent (3) The (normalised) Gibbs cocycle being unchanged when the
potential $F$ is replaced by the potential $F+\sigma$ for any constant
$\sigma$, we may (and will) take the Patterson densities, hence the
Gibbs measure and the skinning measures, to be unchanged by such a
replacement.
\erema

\medskip
When $D$ is a horoball in $X$, let us now relate the skinning
measures of $D$ with previously known measures on $\partial_\infty X$,
constructed using techniques due to Hamenstädt.  

Let $\H$ be a horoball centred at a point $\xi\in\partial_\infty
X$. Recall that $P_\H: \partial_\infty X -\{\xi\} \ra\partial \H$ is
the closest point map on $\H$, mapping $\eta \neq \xi$ to the
intersection with the boundary of $\H$ of the geodesic line from
$\eta$ to $\xi$.  The following result is proved in \cite[\S
  2.3]{HerPau04} when $F=0$.

\bprop \label{prop:HPun} 
Let $\rho:[0,+\infty[\,\ra X$ be the geodesic ray starting from any
point of the boundary of $\H$ and converging to $\xi$.  The following
weak-star limit of measures on $\partial_\infty X-\{\xi\}$
$$
d\mu^\pm_\H(\eta)=\lim_{t\ra+\infty} \;
e^{-\int_{\rho(t)}^{P_\H(\eta)}(\wt F^\pm-\delta)}\,d\mu^\pm_{\rho(t)}(\eta)
$$
exists, and it does not depend on the choice of $\rho$. The measure
$\mu^\pm_\H$ is invariant under the elements of $\Ga$ preserving $\H$,
and it satisfies, for every $x\in X$ and (almost) every $\eta
\in \partial_\infty X-\{\xi\}$,
$$
\frac{d\mu^\pm_\H}{d\mu^\pm_{x}}(\eta)=e^{-C^\pm_\eta(P_\H(\eta),\,x)}\;.
$$
\eprop

\dem We prove all three assertions simultaneously. Let us fix $x\in
X$. For all $t\geq 0$ and $\eta\in \partial_\infty X-\{\xi\}$, let
$z_t$ be the closest point to $P_\H(\eta)$ on the geodesic ray from
$\rho(t)$ to $\eta$.

\begin{center}
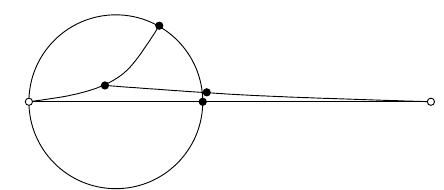
\end{center}

Using Equation \eqref{eq:quasinivarPatdens} with $x$ replaced by
$\rho(t)$ and $y$ by the present $x$, by the cocycle equation
\eqref{eq:cocycle} and by Equation \eqref{eq:changemoinsplus} as
$z_t\in[\rho(t),\eta[\,$, we have
\begin{align*}
e^{-\int_{\rho(t)}^{P_\H(\eta)}(\wt F^\pm-\delta)}\,d\mu^\pm_{\rho(t)}(\eta) &
= e^{-\int_{\rho(t)}^{P_\H(\eta)}(\wt F^\pm-\delta)}\,e^{-C^\pm_\eta(\rho(t),\,x)}
\,d\mu^\pm_{x}(\eta)\\
& = e^{-\int_{\rho(t)}^{P_\H(\eta)}(\wt F^\pm-\delta)}\,e^{-C^\pm_\eta(\rho(t),\,z_t)}
\,e^{-C^\pm_\eta(z_t,\,x)}\,d\mu^\pm_{x}(\eta)\\
& = e^{-\int_{\rho(t)}^{P_\H(\eta)}(\wt F^\pm-\delta)+\int_{\rho(t)}^{z_t}(\wt F^\pm-\delta)}
\,e^{-C^\pm_\eta(z_t,\,x)}\,d\mu^\pm_{x}(\eta)\;.
\end{align*}
As $t\ra+\infty$, note that $z_t$ converges to $P_\H(\eta)$ and that
by the \ref{eq:HC}-property (and since $\wt F$ is bounded on any
compact neighbourhood of $P_\H(\eta)$), we have
$$
\Big|\int_{\rho(t)}^{P_\H(\eta)}(\wt F^\pm-\delta)-
\int_{\rho(t)}^{z_t}(\wt F^\pm-\delta)\Big| \ra 0\;.
$$
The result then follows by the continuity of the Gibbs cocycle (see
Proposition \ref{prop:continuGibbscocycle} (3)).  
\cqfd

\medskip
Using this proposition and the cocycle property of $C^\pm$ in the
definition \eqref{eq:defigibbs} of the Gibbs measure, we obtain, for
every $\ell\in\G X$ such that $\ell_\pm\neq \xi$,
\begin{equation}\label{eq:relatGibbsHamen}
d\wt m_{F}(\ell)= 
e^{C_{\ell_-}^-(P_\H(\ell_-),\,\ell(0))\,+\,C^+_{\ell_+}(P_\H(\ell_+),\,\ell(0))}\;
d\mu_{\H}^-(\ell_-)\,d\mu^+_{\H}(\ell_+)\,dt\;.
\end{equation} 
Note that it is easy to see that for every $\rho\in\normalpm \H$, we
have
\begin{equation}\label{eq:relatskinhoroballHerPau}
d\wt\sigma^\pm_\H(\rho)=d\mu^\pm_\H(\rho_\pm)\;.
\end{equation}

When $F=0$, we obtain {\em Hamenst\"adt's
  measure}\index{Hamenst\"adt's!measure}\index{measure!Hamenst\"adt's}
\begin{equation}\label{eq:defiHamensmeasure}
\mu_{\H}=\lim_{t\ra+\infty} e^{\delta_\Ga\, t}\mu_{\rho(t)}
\end{equation}
on $\partial_\infty X-\{\xi\}$ associated with the horoball $\H$,
which is independent of the choice of the geodesic ray $\rho$ starting
from a point of the horosphere $\partial \H$ and converging to $\xi$.
Note that for every $t\geq 0$, if $\H[t]$ is the horoball contained in
$\H$ whose boundary is at distance $t$ from the boundary of $\H$, we
then have
\begin{equation}\label{eq:shrinkhoroballmeas}
\mu_{\H[t]}=e^{-\delta_\Ga\, t}\;\mu_\H\;.
\end{equation}

\bigskip
Assume till the end of Proposition \ref{prop:computskinninglocsym}
that the potential $\wt F$ is zero. The next result gathers
computations done in \cite{ParPau16LMS,ParPau16MA} of the skinning
measures of horoballs and some totally geodesic subspaces, when $X$ is
a real or complex hyperbolic space and $\Ga$ is a lattice. We
consider the notation $\hnr$, $\hnc$, $\H_\infty$, $\Heis_{2n-1}$,
$\lambda_{2n-1}$ introduced in Section \ref{subsec:BMcomputlocsym},
and we again endow $T^1\hnr$ and $T^1\hnc$ with their Sasaki's
Riemannian metric. Recall that a {\it complex hyperbolic
  line}\index{complex geodesic line} in $\hnc$ is a totally geodesic
plane with constant sectional curvature $-4$.

As the arguments of the following result are purely computational and
rather long, we do not copy them in this book, but we refer respectively
to the proofs of \cite[Prop.~11 (1), (2)]{ParPau16LMS} and
\cite[Lem.~12 (iv), (v), (vi)]{ParPau16MA}.  Analogous computations can
be done when $X$ is the quaternionic hyperbolic $n$-space $\HH_{\HH}^n$.

\bprop\label{prop:computskinninglocsym} (1) Let $\Ga$ be a lattice in
$\Isom(\hnr)$, with Patterson density $(\mu_x)_{x\in\hnr}$ normalised
as in Section \ref{subsec:BMcomputlocsym}.

\smallskip\noindent (i) If $D$ is a horoball in $\hnr$, if
$\vol_{\normalpm D}$ is the Riemannian measure of the submanifold
$\normalpm D$ in $T^1\hnr$, then
$$
\wt \sigma^\pm_{D}= 2^{n-1}\vol_{\normalpm D}\;.
$$
If the point at infinity of $D$ is a parabolic fixed point of $\Ga$, with
stabiliser $\Ga_{D}$ in $\Ga$, then\footnote{See for instance
\cite[p.~473]{Hersonsky93} for the last equality.}
$$
\|\sigma^\pm_{D}\|=2^{n-1}\Vol(\Ga_{D}\bs\normalpm D)
=2^{n-1}\Vol(\Ga_{D}\bs\partial D)
=2^{n-1}(n-1)\Vol(\Ga_{D}\bs D)\;.
$$

\medskip\noindent (ii) If $D$ is a totally geodesic hyperbolic
subspace of dimension $k\in\{1,\dots,n-1\}$ in $\hnr$, if
$\vol_{\normalpm D}$ is the Riemannian measure of the submanifold
$\normalpm D$ in $T^1\hnr$, then
$$
\wt \sigma^\pm_{D}= \vol_{\normalpm D}\,.
$$ 
With $\Ga_{D}$ the stabiliser of $D$ in $\Ga$ and $m$ the order of
the pointwise stabiliser of $D$ in $\Ga$, if $\Ga_D\bs D$ has finite
volume, then
$$
\|\sigma^\pm_{D}\|=\frac{\Vol(\SSS^{n-k-1})}{m}\;\Vol(\Ga_D\bs D)\;.
$$

\medskip
(2) Let $\Ga$ be a lattice in $\Isom(\hnc)$, with Patterson density
$(\mu_x)_{x\in\hnc}$ normalised as in Section
\ref{subsec:BMcomputlocsym}.  

\medskip\noindent(i) Using the homeomorphism $v\mapsto v_\pm$ from
$\normalpm\H_\infty$ to $\partial_\infty \hnc-\{\infty\}=
\Heis_{2n-1}$, we have
$$
d\wt\sigma^\pm_{\H_\infty}(v)=d\lambda_{2n-1}(v_\pm)\,.
$$
For every horoball $D$ in $\hnc$, if $\vol_{\partial D}$ is the
Riemannian measure of the hypersurface $\partial D$ in $\hnc$, then
$$
\pi_*\wt\sigma^\pm_{D}=2\,\vol_{\partial D}\,.
$$ 
If the point at infinity of $D$ is a parabolic fixed point of $\Ga$,
with stabiliser $\Ga_{D}$ in $\Ga$, then
$$
\|\sigma^\pm_{D}\|= 4n\, \Vol(\Ga_{D}\bs D)\,.
$$

\medskip\noindent (ii) For every geodesic line $D$ in $\hnc$, if
$\vol_{\normalpm D}$ is the Riemannian measure of the submanifold
$\normalpm D$ in $T^1\hnc$, we have
$$
d\pi_*\wt \sigma^\pm_{D}=\frac{n}{4^{n-1}\,(2n-1)}\;
d\pi_*\vol_{\normalpm D}\,.
$$
With $\Ga_{D}$ the stabiliser of $D$ in $\Ga$ and $m$ the order of the
pointwise stabiliser of $D$ in $\Ga$, if $\Ga_D\bs D$ has finite
length, then
$$
\|\sigma^\pm_{D}\|= 
\frac{2\,\pi^{n-1}\, n!}{m\,(2n-1)!}\;\Vol(\Ga_{D}\bs D)\,.
$$

\medskip\noindent (iii) For every complex geodesic line $D$ in
$\hnc$, if $\vol_{\normalout D}$ is the Riemannian measure of the
submanifold $\normalout D$ in $T^1\hnc$, we have
$$
d\pi_*\wt \sigma^+_{D}= \frac{1}{2^{2n-3}}\;d\pi_*\vol_{\normalout D}\,.
$$
With $\Ga_{D}$ the stabiliser of $D$ in $\Ga$ and $m$ the order of the
pointwise stabiliser of $D$ in $\Ga$, if $\Ga_D\bs D$ has finite
area, then
$$
\|\sigma^+_{D}\|= 
\frac{\pi^{n-1}}{m\,4^{n-2}\,(n-2)!}\;\Vol(\Ga_{D}\bs D)\;.\;\;\;\Box
$$  
\eprop

\medskip
The following results give the basic properties of the skinning
measures analogous to those in \cite[Sect.~3]{ParPau14ETDS} when the
potential is zero.

\bprop\label{prop:basics} Let $D$ be a nonempty proper closed convex
subset of $X$, and let $\wt\sigma^\pm_{D}$ be the skinning
measures on $\normalpm D$ for the potential $\wt F$. 

\noindent {\rm (i)} The skinning measures $\wt\sigma^\pm_{D}$ are
independent of $x_{0}$.

\noindent {\rm (ii)} For all $\ga\in\Ga$, we have
$\ga_*\wt\sigma^\pm_D= \wt\sigma^\pm_{\ga D}$.  In particular, the
measures $\wt\sigma^\pm_D$ are invariant under the stabiliser of $D$
in $\Ga$.

\noindent {\rm (iii)} For all $s\geq 0$ and $w\in\partial^1_{\pm} D$,
we have\footnote{denoting by $(\flow {\pm s} w)_{\mid\pm[0,+\infty[}$
 the element of $\G_{\pm,\,0}$ which coincides with $\flow {\pm s} w$ 
on $\pm[0,+\infty[$}
\begin{align*}
d\,\wt\sigma^\pm_{\N_sD}((\flow {\pm s} w)\big|_{\pm[0,+\infty[})
&=e^{C^\pm_{w_\pm}(\pi(w),\,\pi(\flow{\pm s}w))}\;d\wt\sigma^\pm_{D}(w)\\ &=
e^{-\int_{\pi(w)}^{\pi(\flow{\pm s}w)}(\wt F^\pm-\delta)}\;
d\wt\sigma^\pm_{D}(w)\;.
\end{align*}

\noindent {\rm (iv)} The support of $\wt\sigma^\pm_{D}$ is
$$
\{v\in\partial_{\pm}^1 D:v_{\pm}\in\Lambda\Gamma\}= P^\pm_D(\Lambda
\Ga -(\Lambda\Ga\cap\partial_\infty D))\;.
$$
In particular, $\wt \sigma^\pm_D$ is the zero measure if and only if
$\Lambda\Ga$ is contained in $\partial_\infty D$.  
\eprop

For future use, the version\footnote{contained in
  \cite[Prop.~4]{ParPau14ETDS}} of Assertion (iii) when $F=0$ is
\begin{equation}\label{eq:skinscale}
\frac{d(\flow{\pm s})_{*}\wt\sigma^\pm_{D}}
{d\,\wt\sigma^\pm_{\N_sD}}(\flow {\pm s} w)= e^{-\delta_\Gamma\,s}\;,
\end{equation}
where $w\in\partial^1_\pm D$ and we again denote by $\flow{\pm s}$ the
map from $\partial^1_\pm D$ to $\partial^1_\pm \N_1D$ defined by
$w\mapsto (\flow {\pm s} w)\big|_{\pm[0,+\infty[}$.

As another particular case of Assertion (iii) for future use, consider
the case when $X= |\XX|_\lambda$ is the geometric realisation of a
metric tree $(\XX, \lambda)$ and when $\wt F=\wt F_c$ is the potential
associated with a system of conductances $\wt c$ on $\XX$ for a
subgroup $\Ga$ of $\Aut(\XX)$ (see Equation
\eqref{eq:defpotfromconduc} and Proposition
\ref{prop:relatpotentialconductance}). Then for all
$w\in\partial^1_{+} D$ (respectively $w\in\partial^1_{-} D$), if $e_w$
is the first (respectively the last) edge followed by $w$, with length
$\lambda(e_w)$, then
$$
\int_{\pi(w)}^{\pi(\flow{\pm \lambda(e_w)}w)}\wt F^\pm=\wt c(e_w)\;\lambda(e_w)
$$
by Proposition \ref{prop:integpotconduct}, so that
\begin{equation}\label{eq:pousstempunconduct}
d\,\wt\sigma^\pm_{\N_sD}
((\flow {\pm \lambda(e_w)} w)_{\mid\pm[0,+\infty[})
= e^{-(\wt c(e_w)+\delta)\lambda(e_w)}\;d\wt\sigma^\pm_{D}(w)
\;.
\end{equation}

\medskip 
\dem The proofs of the claims are straightforward modifications of
those for zero potential in \cite[Prop.~4]{ParPau14ETDS}.  We give
details of the proofs for the measure $\wt\sigma^+_{D}$, the case of
$\wt\sigma^-_{D}$ being similar.

\smallskip
\noindent(i) The claim follows from Equation
\eqref{eq:quasinivarPatdens} and the cocycle property
\eqref{eq:cocycle}.

\smallskip
\noindent(ii) The claim follows from Equation
\eqref{eq:equivarPatdens}, the first part of Equation
\eqref{eq:cocycle} and Claim (i).

\smallskip\noindent(iii) Since $\big((\flow {s}
w)_{\mid[0,+\infty[}\big)_+=w_+$ and since $w\in \normalout D$ if and
only if $(\flow {s} w)_{\mid[0,+\infty[}\in \normalout \N_s D$, we
have, using the definition \eqref{eq:defiskin} of the skinning measure
and the cocycle property \eqref{eq:cocycle}, for all $s\geq 0$,
$$
  d\,\wt\sigma^+_{\N_{s}D}((\flow {s}w)_{\mid [0,+\infty[})=
  e^{C^+_{w_{+}}(x_0,\,\pi(\flow sw))}\,d\mu^+_{x_{0}}(w_{+})=
  e^{C^+_{w_{+}}(\pi(w),\,\pi(\flow s w))}\,d\,\wt\sigma^+_{D}( w)\,.
$$
This proves Claim (iii) for $\wt\sigma^+_{D}$, using Equation
\eqref{eq:changemoinsplus}. 

\smallskip
\noindent (iv) The claims follow from the fact that the support of any
Patterson measure is $\Lambda\Ga$, see Subsection
\ref{subsec:Pattersondens}.  
\cqfd

\medskip
Given two nonempty closed convex subsets $D$ and $D'$ of $X$, let
$$
A_{D,D'}=\partial_{\infty}X-(\partial_{\infty} D \cup 
\partial_{\infty}D')
$$ 
and let $h^\pm_{D,D'}: P^\pm_{D}(A_{D,D'}) \to P^\pm_{D'}(A_{D,D'})$
be the restriction of $P^\pm_{D'}\circ (P^\pm_{D})^{-1}$ to
$P^\pm_{D}(A_{D,D'})$. It is a homeomorphism between open subsets of
$\normalpm D$ and $\normalpm D'$, associating to the element $w$ in
the domain the unique element $w'$ in the range with $w'_\pm=w_\pm$.

\bprop\label{prop:abscontskinmeas} Let $D$ and $D'$ be nonempty closed
convex subsets of $X$ and let $h^\pm=h^\pm_{D,D'}$.  The measures
$(h^\pm)_{*}\;\wt\sigma^\pm_{D}$ and $\wt\sigma^\pm_{D'}$ on $
P^\pm_{D'}(A_{D,D'})$ are absolutely continuous one with respect to
the other, with
$$
\frac{d(h^\pm)_{*}\;\wt\sigma^\pm_{D}}{d\wt\sigma^\pm_{D'}}(w')=
e^{-C^\pm_{w_{\pm}}(\pi(w),\,\pi(w'))},
$$
for (almost) all $w\in P^\pm_{D}(A_{D,D'})$ and $w'=h^\pm(w)$. 
\eprop
\dem
As $w'_\pm=w_\pm$, we have
\begin{align*}
d\wt\sigma^\pm_{D'}(w')  & = 
e^{C^\pm_{w'_\pm}(x_0,\,\pi(w'))}\;d\mu^\pm_{x_0}(w'_{\pm}) =
e^{C^\pm_{w'_\pm}(x_0,\,\pi(w))}\;e^{C^\pm_{w'_\pm}(\pi(w),\,\pi(w'))}\;
d\mu^\pm_{x_0}(w'_{\pm})\\
& =
e^{C^\pm_{w_\pm}(\pi(w),\,\pi(w'))}\;d\wt\sigma^\pm_{D}(w) 
\end{align*}
using the definition \eqref{eq:defiskin} of the skinning measure and
the cocycle property \eqref{eq:cocycle}.  
\cqfd
 
\medskip
Let $w\in\G_\pm X$. With $N^\pm_w: \wssu(w)\to\normalmp H\!B_\pm(w)$
the canonical homeomorphism defined in Section \ref{subsect:nbhd}, we
define the {\em skinning measures}%
\index{skinning measure}\index{measure!skinning} $\mussu{w}$ on
the strong stable or strong unstable leaves $\wssu(w)$ by
$$
\gls{strongmeasure}=((N^\pm_w)^{-1})_*\wt \sigma^\mp_{H\!B_{\pm}(w)} \;,
$$
so that 
\begin{equation}\label{eq:condtionelless}
d\mussu{w}(\ell) =e^{C^\mp_{\ell_\mp}(x_0,\,\ell(0))} \;
d\mu^\mp_{x_0}(\ell_{\mp})
\end{equation}
for every $\ell\in\wssu(w)$. By Proposition \ref{prop:basics} (ii) and
the naturality of $N^\pm_w$, for every $\ga\in\Ga$, we have
\begin{equation}\label{eq:mussequiv}
\ga_*\mussu{w}=\mussu{\ga w}\;.
\end{equation}
By Proposition \ref{prop:basics} (iv), the support of $\mussu{w}$ is
$\{\ell\in \wssu(w)\;:\;\ell_{\mp} \in \Lambda \Gamma\}$.  For all
$t\in\RR$ and $\ell\in \wssu(w)$, we have, using Equations
\eqref{eq:condtionelless}, \eqref{eq:cocycle} and
\eqref{eq:changemoinsplus}, and since $\ell_\pm=w_\pm$,
\begin{equation}\label{eq:mussscaling}
\frac{d(\flow{-t})_{*}\mussu{w}}{d\,\mussu{\flow{t}w}}(\flow{t}\ell)
=e^{C^\mp_{\ell_\mp}(\ell(t),\,\ell(0))}=e^{C^\pm_{w_\pm}(\ell(0),\,\ell(t))}\;.
\end{equation}

Let $w\in\G_\pm X$. The homeomorphisms $\wssu(w)\times\RR\to \wosu(w)$,
defined by 
$$
(\ell,s)\mapsto \ell'=\flow s\ell\;,
$$ 
conjugate the actions of $\RR$ by translation on the second factor of
the domain and by the geodesic flow on the range, and the actions of
$\Ga$ (trivial on the second factor of the domain).  Let us consider
the measures $\gls{conditionalweak}$ on $\wosu(w)$ given, using the
above homeomorphism, by
\begin{equation}\label{eq:definumoins}
d\nu^\mp_w(\ell')=
e^{C^\pm_{w_\pm}(w(0),\,\ell(0))}\;d\mussu{w}(\ell)\;ds\;.
\end{equation}
They satisfy $(\flow t)_*\nu^\pm_w=\nu^\pm_w$ for all $t\in\RR$ (since
if $\ell'=\flow{s}\ell$, then $\flow{-t}\ell'=\flow{s-t}\ell$, and by
invariance under translations of the Lebesgue measure on $\RR$).
Furthermore, $\ga_* \nu^\pm_w=\nu^\pm_{\ga w}$ for all $ \ga\in\Ga$.
In general, they depend on $w$, not only on $W^\pm(w)$. Furthermore,
the support of $\nu^\pm_w$ is $\{\ell'\in \wosu(w)\;:\;\ell'_{\mp} \in
\Lambda \Gamma\}$.  These properties follow easily from the properties
of the skinning measures on the strong stable or strong unstable
leaves.

\blemm \label{lem:defiR} (i) For every nonempty proper closed convex
subset $D'$ in $X$, there exists $R_0>0$ such that for all $R\geq
R_0$, $\eta>0$, and $w\in\partial^1_\pm D'$, we have
$\nu^\mp_w(V^\pm_{w,\,\eta,\,R})>0$.\footnote{See Section
  \ref{subsect:nbhd} for the definition of $V^\pm_{w,\,\eta,\,R}$.}

\smallskip\noindent (ii) For all $w\in\G_\pm X$ and $t\in\RR$, the
measures $\nu^\mp_{\flow{t}w}$ and $\nu^\mp_w$ are proportional:
$$
\nu^\mp_{\flow{t}w}=e^{C^\pm_{w_\pm}(w(t),\,w(0))}\;\nu^\mp_w\,.
$$ 
\elemm

\dem (i) Let us show, as in \cite[Lem.~7]{ParPau14ETDS}, that there exists
$R_0>0$ (depending only on $D'$ and on the Patterson densities) such
that for all $R\geq R_0$, $w\in\partial^1_+D'$ and
$w'\in\partial^1_-D'$, we have $\muss{w} (B^+(w,R))>0$ and
$\musu{w'}(B^-(w',R))>0$.  The result follows from this by the definitions
of $\nu^\mp_w$ and $V^\pm_{w,\,\eta,\,R}$.

\medskip\noindent
\begin{minipage}{9.5cm} 
~~~ We give the proof of the claim on $B^+(w,R)$, the proof of the
  claim on $B^-(w,R)$ is similar.  For all $w\in \normal D'$ and
  $\xi'\in D'\cup \partial_{\infty} D'$, by a standard comparison and
  convexity argument applied to the geodesic triangle with vertices
  $\pi(w), w_+,\xi'$, the point $\pi(w)$ is at distance at most
  $2\ln(\frac{1+\sqrt 5}{2})$ from the intersection between the
  stable horosphere $H_+(w)$ and the geodesic ray or line between
  $\xi'$ and $w_+$.
\end{minipage}
\begin{minipage}{4.4cm}
\begin{center}
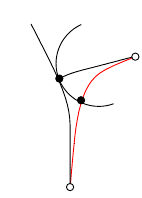
\end{center}
\end{minipage}

\medskip
The triangle inequality and the definition of Hamenstädt distances
imply that, for all $\ell,\ell'\in \wss(w)$,
\begin{equation}\label{eq:majodistHamen}
d_{\wss(w)}(\ell,\ell')\leq e^{\frac{1}{2}d(\pi(\ell),\,\pi(\ell'))}\;.
\end{equation}
Hence, for every $\xi'\in \partial_{\infty} D'$, for every extension
$\wh w\in\G X$ of $w$, if $\ell'$ is the element of $\wss(w)$ such
that $\ell'_-=\xi'$, we have
$$
d_{\wss(w)}(\wh w, \ell')\leq \frac{1+\sqrt 5}{2}\;.
$$ 
Thus, if $\partial_{\infty} D'\cap\Lambda\Ga\ne\emptyset$, then we
may take $R_0 = 2 >\frac{1+\sqrt 5}{2}$, since by Proposition
\ref{prop:basics} (iv), the support of $\muss{w}$ is $\{\ell\in
\wss(w)\;:\; \ell_-\in\Lambda\Ga\}$.

Assume now that $\partial_{\infty} D'\cap\Lambda\Ga=\emptyset$. For a
contradiction, assume that, for all $n\in\NN$, there exists
$w_n\in\normal D'$ such that $\muss{w_n}(B^+(w_n,n))=0$.  Assume first
that $(w_n)_{n\in\NN}$ has a convergent subsequence with limit $w\in
\normal D'$. Since the measure $\muss{w'}$ depends continuously on
$w'\in\G_+X$, for every compact subset $K$ of $\wss(w)$, we have
$\muss{w}(K)=0$. By Proposition \ref{prop:basics} (iv) and by Equation
\eqref{eq:condtionelless}, this implies that the support of the
Patterson measure $\mu^-_{x_0}$, which is the limit set of $\Ga$, is
contained in $\{w_{+}\}$. This is impossible, since $\Ga$ is
nonelementary.

In the remaining case, the points $\pi(w_n)$ in $D'$ converge, up to
extracting a subsequence, to a point $\xi$ in $\partial_{\infty}D'$.
By definition of the map $P^+_{D'}$ and of $\normal D'$, the points at
infinity $(w_n)_{+}$ converge to $\xi$. For every $\eta$ in
$\partial_\infty X$ different from $\xi$, the geodesic lines from
$\eta$ to $(w_n)_{+}$ converge to the geodesic line from $\eta$ to
$\xi$. 

\medskip\noindent
\begin{minipage}{9cm} 
  ~~~ By convexity, if $n$ is large enough, the geodesic line
  $]\eta,(w_n)_{+}[$ meets $\N_1D'$, hence passes at distance at most
  $2$ from $\pi(w_n)$. This implies by Equation 
  \eqref{eq:majodistHamen}  that if $n$ is large enough, then
  there exists $\ell\in B^+(w_n, \,n)$ such that $\eta = \ell_-$.
\end{minipage}
\begin{minipage}{4.9cm}
\begin{center}
~~~~~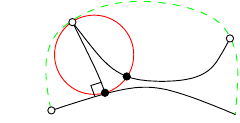
\end{center}
\end{minipage}

\medskip Since we assumed that $\muss{w_n}(B^+(w_n,n))=0$ for
all $n\in\NN$, Proposition \ref{prop:basics} (iv) implies that we
have $\eta\notin\Lambda\Ga$. Hence $\Lambda\Ga$ is contained in
$\{\xi\}$, a contradiction since $\Ga$ is nonelementary. 

\smallskip\noindent (ii) For all $w\in\G_\pm X$, $s,t\in\RR$ and
$\ell\in \wssu(w)$, we have by Equations \eqref{eq:definumoins} and
\eqref{eq:mussscaling}, and by the cocycle property \eqref{eq:cocycle}
of $C^\pm$,
\begin{align*}
d\nu^\mp_{\flow{t}w}(\flow{s}\ell)&=
d\nu^\mp_{\flow{t}w}(\flow{s-t}\flow{t}\ell)=
e^{C^\pm_{(\flow{t}w)_\pm}(\flow{t}w(0),\,\flow{t}\ell(0))}
\;d\mussu{\flow{t}w}(\flow{t}\ell)\;d(s-t)\\ & =
e^{C^\pm_{w_\pm}(w(t),\,\ell(t))}\;
e^{-C^\pm_{w_\pm}(\ell(0),\,\ell(t))}\;d\mussu{w}(\ell)\;ds\\ &=
e^{C^\pm_{w_\pm}(w(t),\,\ell(0))}\;
e^{-C^\pm_{w_\pm}(w(0),\,\ell(0))}\;d\nu^\mp_w(\flow{s}\ell)\\
&= e^{C^\pm_{w_\pm}(w(t),\,w(0))}d\nu^\mp_w(\flow{s}\ell)\;. \;\;\;\Box
\end{align*}

\medskip The following disintegration result of the Gibbs measure over
the skinning measures of any closed convex subset is a crucial tool
for our equidistribution and counting results. Recall the definition
in Equation \eqref{eq:defiUCp} of the flow-invariant open sets
$\U^\pm_D$ and the definition of the fibrations $f^\pm_D:\U^\pm_D\to
\normalpm D$ from Section \ref{subsect:nbhd}.

\bprop\label{prop:disintegration} Let $D$ be a nonempty proper closed
convex subset of $X$.  The restriction to $\U^\pm_D$ of the Gibbs
measure $\wt m_{F}$ disintegrates by the fibration $f^\pm_D: \U^\pm_D
\to \partial^1_{\pm} D$ over the skinning measure $\wt\sigma^\pm_D$
of $D$, with conditional measure $\nu^\mp_\rho$ on the fiber
$(f^\pm_D)^{-1}(\rho)=W^{0\pm}(\rho)$ of $\rho\in \partial^1_{\pm} D$:
when $\ell$ ranges over $\U^\pm_D$, we have
$$
d\wt m_{F}{}_{\mid \U^\pm_D}(\ell)=\int_{\rho\in \partial^1_{\pm} D}
d\nu^\mp_\rho(\ell)\,
d\wt\sigma^\pm_{D}(\rho)\;.\
$$
\eprop

\dem In order to prove the claim for the fibration $f^+_{D}$, let
$\phi\in\C_{\rm c}(\U^+_D)$. Using in the various steps below:
\begin{itemize}
\item Hopf's parametrisation with time parameter $t$ and the
  definitions of $\wt m_F$ (see Equation \eqref{eq:defigibbs}) and of
  $\U^+_D$ (see Equation \eqref{eq:defiUCp}),

\item the positive endpoint homeomorphism $w\mapsto w_+$ from
  $\normalout D$ to $\partial_\infty X-\partial_\infty D$, and the
  negative endpoint homeomorphism $\ell'\mapsto \ell'_-$ from
  $\wss(w)$ to $\partial_\infty X-\{w_+\}$, with $s\in\RR$ the real
  parameter such that $\ell'=\flow{-s}\ell\in\wss(w)$ where $\ell\in
  \ws(w)$, noting that $t-s$ depends only on $\ell_{+}=w_{+}$ and
  $\ell_{-}=\ell'_{-}$,

\item the definitions Equation \eqref{eq:condtionelless} and
  \eqref{eq:defiskin} of the measures $\muss{w}$ and
  $\wt\sigma^+_{D}$, and the cocycle property \eqref{eq:cocycle} of
  $C^\pm$,

\item Equation \eqref{eq:changemoinsplus} and the cocycle property
  \eqref{eq:cocycle} of $C^+$,
\end{itemize}
\noindent
we have {\small
\begin{align*}
   &\int_{\ell\in \U^+_D} \phi(\ell)\,d\wt m_F(\ell)  \\
  = &\int_{\ell_+\in\partial_\infty X-\partial_\infty D}
  \int_{\ell_-\in\partial_\infty X-\{\ell_+\}}\int_{t\in\RR}
  \phi(\ell)\;e^{C^-_{\ell_-}(x_0,\,\pi(\ell))+C^+_{\ell_+}(x_0,\,\pi(\ell))}\,
dt\,d\mu^-_{x_0}(\ell_-)\,
  d\mu^+_{x_0}(\ell_+)\\
  = & \int_{w\in\normalout D} \int_{\ell'\in\wss(w)}\int_{s\in\RR}
  \phi(\flow s \ell')\;e^{C^-_{\ell'_-}(x_0,\,\pi(\flow s
    \ell'))+C^+_{w_+}(x_0,\,\pi(\flow s \ell'))} 
 \,ds\,d\mu^-_{x_0}(\ell'_-)\,d\mu^+_{x_0}(w_+)\\
  = & \int_{w\in\normalout D} \int_{\ell'\in \wss(w)}\int_{s\in\RR}
  \phi(\flow s \ell')\;e^{C^-_{\ell'_-}(\pi(\ell'),\,\pi(\flow s \ell'))
    +C^+_{w_+}(\pi(w),\,\pi(\flow s \ell'))}\,
  ds\,d\muss{w}(\ell')\,d\wt\sigma^+_{D}(w)\\
  = & \int_{w\in\normalout D} \int_{\ell'\in \wss(w)}\int_{s\in\RR}
  \phi(\flow s \ell')\;e^{C^+_{w_+}(\pi(w),\,\pi(\ell'))}\,ds\,
  d\muss{w}(\ell')\,d\wt\sigma^+_{D}(w)\,,
\end{align*}}%
which implies the claim for the fibration $f^+_{D}$, by the definition
\eqref{eq:definumoins} of the measure $\nu^-_w$. The proof for the
fibration $f^-_{D}$ is similar.  
\cqfd

\bigskip For every $u\in\G_-X$,  if $D=H\!B_-(u)$, we have 
$\normalout D=N^-_u(\wsu(u))$ and
$$
\U^+_D=\G X-\ws(\iota u)=\bigcup_{w\in\wsu(u)} \ws(w)\;.
$$ 
Applying the above proposition and a change of variable, 
the restriction to $\G X-\ws(\iota u)$ of the Gibbs measure $\wt m_F$
disintegrates over the strong unstable measure $\musu{u}=
((N^-_u)^{-1})_* \wt\sigma^+_D$, with conditional measure on the fiber
$\ws(w)$ of $w\in \wsu(u)$ the measure $\nu^-_w=\nu^-_{N^-_u(w)}$: for every
$\phi\in\C_{\rm c}(\G X-\ws(\iota u))$, we have
\begin{align}
  &\int_{\ell\in\G X-\ws(\iota u)}\phi(\ell)\,d\wt m_F(\ell) =  \nonumber
\\ &
 \int_{w\in \wsu(u)} \int_{\ell'\in \wss(w)}\int_{s\in\RR}
  \phi(\flow s \ell')\;e^{C^+_{w_+}(\pi(w),\,\pi(\ell'))}\,
ds\,d\mu_{\wss(w)}(\ell')\,d\musu{u}(w)\,.\label{eq:desintegrbox}
\end{align}
Note that if the Patterson densities have no atoms, then the stable
and unstable leaves have measure zero for the associated Gibbs
measure. This happens for instance if the Gibbs measure $m_F$ is
finite, see Corollary \ref{coro:finitudeGibbsdivuniq} and Theorem
\ref{theo:HTSR}.

\section{Equivariant families of convex subsets and 
their skinning measures}
\label{subsec:equivfammult}

Let $I$ be an index set endowed with a left action of $\Ga$. A family
$\D=(D_i)_{i\in I}$ of subsets of $X$ or of $\gengeod X$ indexed by $I$
is {\em $\Ga$-equivariant}\index{equivariant family} if $\ga
D_i=D_{\ga i}$ for all $\ga\in\Ga$ and $i\in I$. We will denote by
$$
\gls{equivrelatequivfamil}
$$ 
the equivalence relation on $I$ defined by $i\sim j$ if and only if
$D_i=D_j$ and there exists $\ga\in\Ga$ such that $j=\ga i$. This
equivalence relation is $\Ga$-equivariant: for all $i,j\in I$ and
$\ga\in\Ga$, we have $\ga i\sim\ga j$ if and only if $i\sim j$. We say
that $\D$ is {\em locally finite}%
\index{equivariant family!locally finite}\index{locally!finite} if for
every compact subset $K$ in $X$ or in $\gengeod X$, the quotient set
$\{i\in I: D_i\cap K\ne\emptyset\}/_\sim$ is finite.

\medskip
\noindent{\bf Examples.} (1) Fixing a nonempty proper closed convex
subset $D$ of $X$, taking $I=\Ga$ with the left action by translations
$(\ga,i)\mapsto \ga i$, and setting $D_i=i D$ for every $i\in\Ga$
gives a \mbox{$\Ga$-equivariant} family $\D=(D_i)_{i\in I}$. In this case, we
have $i\sim j$ if and only if $i^{-1}j$ belongs to the stabiliser
$\Ga_D$ of $D$ in $\Ga$, and $I/_\sim\; =\Ga/\Ga_D$.  Note that $\ga
D$ depends only on the class $[\ga]$ of $\ga$ in $\Ga/\Ga_D$. We could
also take $I'=\Ga/\Ga_D$ with the left action by translations $(\ga,
[\ga'])\mapsto [\ga \ga']$, and $\D'=(\ga D)_{[\ga]\in I'}$, so that
for all $i,j\in I'$, we have $i\sim_{\D'} j$ if and only if $i=j$, and
besides, $\D'$ is locally finite if and only if $\D$ is locally
finite.  The following choices of $D$ yield equivariant families with
different characteristics:
\begin{enumerate}
\item[(a)] Let $\ga_0\in\Ga$ be a loxodromic element with translation
  axis $D=\Ax_{\ga_0}$.  The family $(\ga D)_{\ga\in \Ga}$ is locally
  finite and $\Ga$-equivariant. Indeed, by Lemma
  \ref{lem:axespastropproches}, only finitely many elements of the
  family $(\ga D)_{\ga\in \Ga/\Ga_D}$ meet any given bounded subset of
  $X$.
\item[(b)] Let $\ell\in\G X$ be a geodesic line whose image under the
  canonical map $\G X\ra\Ga\backslash\G X$ has a dense orbit in
  $\Ga\backslash\G X$ under the geodesic flow, and let $D=\ell(\RR)$
  be its image. Then the $\Ga$-equivariant family $(\ga
  D)_{\ga\in\Ga}$ is not locally finite.
\item[(c)] More generally, let $D$ be a convex subset such that
  $\Ga_D\bs D$ is compact. Then the family $(\ga D)_{\ga\in \Ga}$ is a
  locally finite $\Ga$-equivariant family.
\item[(d)] Let $\xi\in \partial_\infty X$ be a bounded parabolic limit
  point of $\Ga$, and let $\H$ be any horoball in $X$ centred at
  $\xi$. Then the family $(\ga\H)_{\ga\in \Ga}$ is a locally finite
  $\Ga$-equivariant family.
\end{enumerate}

\smallskip\noindent 
(2) More generally, let $(D^\alpha)_{\alpha\in A}$ be a finite family
of nonempty proper closed convex subsets of $X$, and for every
$\alpha\in A$, let $F_\alpha$ be a finite set. Define $I=
\bigcup_{\alpha\in A} \Ga\times\{\alpha\}\times F_\alpha$ with the
action of $\Ga$ by left translation on the first factor, and for every
$i=(\ga,\alpha,x)\in I$, let $D_i=\ga D^\alpha$. Then $I/_\sim\;
=\bigcup_{\alpha\in A} \Ga/\Ga_{D^\alpha}\times\{\alpha\}\times
F_\alpha$ and the $\Ga$-equivariant family $\D=(D_i)_{i\in I}$ is
locally finite if and only if the family $(\ga D^\alpha)_{\ga\in \Ga}$
is locally finite for every $\alpha\in A$. The cardinalities of
$F_\alpha$ for $\alpha\in A$ contribute to the multiplicities (see
Section \ref{subsec:downstairs}).

\bigskip 
Let $\D=(D_i)_{i\in I}$ be a locally finite $\Ga$-equivariant family
of nonempty proper closed convex subsets of $X$.  Let $\Omega=
(\Omega_i)_{i\in I}$ be a $\Ga$-equivariant family of subsets of
$\gengeod X$, where $\Omega_i$ is a measurable subset of $\normalpm
D_i$ for all $i\in I$ (the sign $\pm$ being constant), such that
$\Omega_i=\Omega_j$ if $i\sim_\D j$.  Then
$$
\gls{skinningpmOmega}=
\sum_{i\in I/_\sim}\wt\sigma^\pm_{D_i}|_{\Omega_i}\;,
$$ 
is a well-defined $\Ga$-invariant locally finite measure on
$\gengeod X$, whose support is contained in $\G_{\pm,\,0}X$.  Hence,
the measure $\wt\sigma^\pm_{\Omega}$ induces\footnote{See for instance
  \cite[\S 2.6]{PauPolSha15} and the beginning of Chapter
  \ref{sec:equidcountdownstairs} for details on the definition of the
  induced measure when $\Ga$ may have torsion, hence does not
  necessarily acts freely on $\gengeod X$.} a locally finite measure
on $\Ga \bs \gengeod X$, denoted by $\sigma^\pm_{\Omega}$.
When $\Omega=\normalpm \D=(\normalpm D_i)_{i\in I}$, the measure
$\wt\sigma^\pm_{\Omega}$ is denoted by
$$
\gls{skinningpmcalD}=\sum_{i\in I/_\sim} \wt \sigma^\pm_{D_i}\,.
$$
The measures $\wt\sigma^+_{\D}$ and $\wt\sigma^-_{\D}$ are
respectively called the {\em outer and inner skinning
  measures}\index{skinning measure!inner}%
\index{skinning measure!outer} of $\D$ on $\gengeod X$, and their
induced measures $\gls{skinningpcalDdown}$ and
$\gls{skinningmcalDdown}$ on $\Ga\bs \gengeod X$ are the {\em outer
  and inner skinning measures}%
\index{skinning measure!inner}\index{skinning measure!outer} of 
$\D$ on $\Ga\bs \gengeod X$.

\medskip
\noindent{\bf Example.}  Consider the $\Ga$-equivariant family $\D=
(\ga D)_{\ga\in \Ga/\Ga_{x}}$ with $D=\{x\}$ a singleton in $X$.  With
$\pi^\pm = (P^\pm_D)^{-1}:\normalpm D\ra\partial_\infty X$ the
homeomorphism $\rho\mapsto \rho_\pm$, we have $(\pi^\pm)_*
\wt\sigma^\pm_{D}=\mu^\pm_x$ by Remark \ref{rem:skinrem} (1), and\footnote{See also the beginning of Chapter \ref{sec:equidcountdownstairs}.}
\begin{equation}\label{eq:skinning singleton}
\|\sigma^\pm_{\D}\|\;=\;\frac{\|\mu^\pm_x\|}{|\Ga_x|}\;.
\end{equation}

\chapter{Explicit measure computations for 
simplicial trees and graphs of groups}
\label{sect:explicitcomputdiscret}

In this Chapter, we compute skinning measures and Bowen-Margulis
measures for some highly symmetric simplicial trees $\XX$ endowed with
a nonelementary discrete subgroup $\Ga$ of $\Aut(\XX)$. These
computations are parallel to the ones given in Section
\ref{subsec:BMcomputlocsym} when $X$ is a rank one symmetric space.
The potentials $F$ are supposed to be $0$ in this Chapter, and we
assume that the Patterson densities $(\mu^+_x)_{x\in V\XX}$ and
$(\mu^-_x)_{x\in V\XX}$ of $\Ga$ are equal, denoted by $(\mu_x)_{x\in
  V\XX}$. As the study of geometrically finite discrete subgroups of
$\Aut(\XX)$ mostly reduces to the study of particular (tree) lattices
(see Remark \ref{rem:geomfiniimpllattice}), we will assume that $\Ga$
is a lattice in this Chapter.

The results of these computations will be useful when we state special
cases of the equidistribution and counting results in regular and
biregular trees and, in particular, in the arithmetic applications in
Part \ref{sect:arithappli}.  The reader only interested in the
continuous time case may skip directly to Chapter
\ref{sec:mixingrate}.

\medskip
A rooted simplicial tree $(\XX,x_0)$ is {\em spherically
  symmetric}\index{spherically symmetric} if $\XX$ is not reduced to
$x_0$ and has no terminal vertex, and if the stabiliser of $x_0$ in
$\Aut(\XX)$ acts transitively on each sphere of centre $x_0$. The set
of isomorphism classes of spherically symmetric rooted simplicial
trees $(\XX,x_0)$ is in bijection with the set of sequences
$(p_n)_{n\in\NN}$ in $\NN-\{0\}$, where $p_n+1$ is the degree of any
vertex of $\XX$ at distance $n$ from $x_0$.

\begin{center}
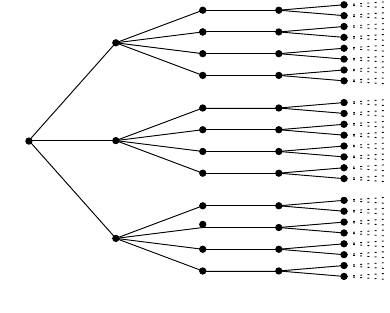
\end{center}

If $(\XX,x_0)$ is spherically symmetric, it is easy to check that the
simplicial tree $\XX$ is uniform\footnote{See Section
  \ref{subsec:trees} for the terminology concerning simplicial trees.}
if and only if the sequence $(p_n)_{n\in\NN}$ is periodic with
palindromic period in the sense that there exists $N\in \NN -\{0\}$
such that $p_{n+N}=p_n$ for every $n\in\NN$ and $p_{N-n}=p_n$ for
every $n\in \NN$ such that $n\leq N$.  If $N=1$, then $\XX=\XX_{p_0}$
is the regular tree of degree $p_0+1$, and if $N=2$, then $\XX=
\XX_{p_0,\,p_1}$ is the biregular tree of degrees $p_0+1$ and $p_1+1$.

We denote by $X=|\XX|_1$ the geometric realisation of $\XX$. The
Hausdorff dimension $h_\XX$ of $\partial_\infty X$ for any visual
distance is then
\begin{equation}\label{eq:lyonshdim}
h_\XX=\frac{1}{N}\ln(p_0\dots p_{N-1})\;,
\end{equation}
see for example \cite[p.~935]{Lyons90}.

\section{Computations of Bowen-Margulis measures 
for simplicial trees}
\label{subsec:BowMargmeas}

The next result gives examples of computations of the total mass of
Bowen-Margulis measures for lattices of simplicial trees having some
regularity properties.

Analogous computations can be performed for Riemannian manifolds
having appropriate regularity properties.  We refer for instance to
\cite[Prop.~10]{ParPau16LMS} and \cite[Prop.~20 (1)]{ParPau16ETDS} for
computations of Bowen-Margulis measures for lattices in the isometry
group of the real hyperbolic spaces, and to \cite[Lem.~12
  (iii)]{ParPau16MA} for the computation in the complex hyperbolic
case. In both cases, the main point is the computation of the
proportionality constant between the Bowen-Margulis measure and
Sasaki's Riemannian volume of the unit tangent bundle. When dealing
now with simplicial trees, similar consequences of homogeneity
properties will appear below.

We refer to Section \ref{subsec:trees} for the definitions of $\vol$,
$\Vol$, $T\pi$, $\Tvol$, $\TVol$ appearing in the following result.

\bprop \label{prop:computBM} Let $(\XX,x_0)$ be a spherically
symmetric rooted simplicial tree, with associated sequence
$(p_n)_{n\in\NN}$, such that $\XX$ is uniform, and let $\Ga$ be a
lattice of $\XX$. 

\smallskip\noindent 
(1) For every $x\in V\XX$, let $r_x=d(x,\Aut(\XX)x_0)$, and let 
$$
c_x= \frac{(p_{r_x}-1)
  e^{2\,r_x\,h_\XX}} {(p_0+1)^2p_1^2\,\dots\, p_{r_x-1}^2p_{r_x}} +
\frac{2p_0}{(p_0+1)^2}
$$ 
if $r_x\neq 0$ and $c_x= \frac{p_0}{p_0+1}$ if $r_x=0$. Then
\begin{align}
\|m_{\rm BM}\| &= \sum_{[x]\in \Ga\bs V \XX} \frac 1{|\Ga_x|}
\big(\|\mu_x\|^2-\sum_{e\in E\XX\;:\; o(e)=x}\mu_x(\partial_{e}\XX)^2\big) 
\nonumber \\ & = \|\mu_{x_0}\|^2
\sum_{[x]\in \Ga\bs V \XX} \frac{c_x}{|\Ga_x|}\;.\label{eq:formulgenemBM}
\end{align}

\smallskip\noindent (2) If $\XX=\XX_{p,\,q}$ is the biregular tree of
degrees $p+1$ and $q+1$, with $V\XX=V_p\XX\,\sqcup\, V_q\XX$ the
corresponding partition of the set of vertices of $\XX$, if the
Patterson density $(\mu_x)_{x\in V \XX}$ of $\Ga$ is normalised so
that $\|\mu_{x}\|= \frac{p+1}{\sqrt{p}}$ for all $x\in V_p\XX$, then
$$
(T\pi)_*m_{\rm BM}=\Tvol_{\Ga\dbs\XX}
$$
and
\begin{equation}\label{eq:relatmBMvolbihomog}
\|m_{\rm BM}\|=\TVol(\Ga\dbs\XX)= 
\sum_{[x]\in \Ga\bs V_p \XX} \frac{p+1}{|\Ga_x|}\;\;+\;\; 
\sum_{[x]\in \Ga\bs V_q \XX} \frac{q+1}{|\Ga_x|}\;.
\end{equation}

\noindent (3) If $\XX=\XX_{q}$ is the regular tree of degree $q+1$, if
the Patterson density $(\mu_x)_{x\in V \XX}$ of $\Ga$ is normalised to
be a family of probability measures, then
$$
\pi_*m_{\rm BM} = \frac{q}{q+1}\;\vol_{\Ga\dbs\XX}
$$ 
and in particular
\begin{equation}\label{eq:relatmBMvolhomog}
\|m_{\rm BM}\|=\frac{q}{q+1} \Vol(\Ga\dbs\XX)\;.
\end{equation}
\eprop

\medskip \dem We start by proving the first equality of Assertion (1).
For every $x\in V\XX$, we may partition the set of geodesic lines
$\ell\in\G\XX$ with $\ell(0)=x$ according to the two edges starting
from $x$ contained in the image of $\ell$. The only restriction for
the edges is that they are required to be distinct.

For every $e\in E\XX$, recall from Section \ref{subsec:trees} 
that $\partial_e\XX$ is the set of points at infinity of the geodesic
rays whose initial edge is $e$. For all $e\in E\XX$ and $x\in V\XX$,
say that $e$ {\em points away}\index{pointing!away} from $x$ if
$o(e)\in [x,t(e)]$, and that $e$ {\em points
  towards}\index{pointing!towards} $x$ otherwise. In particular, all
edges with origin $x$ point away from $x$.
Hence by Equation
\eqref{eq:decompGibbsbas}, and since $\mu_x=\mu^-_x=\mu^+_x$, we have
\begin{align}
\pi_*m_{\rm BM}&=\sum_{[x]\in \Ga\bs V \XX} \frac{1}{|\Ga_x|}\;
\sum_{e,\,e'\in E\XX\;:\; o(e)=o(e')=x,\;e\neq e'}
\mu^-_x(\partial_e\XX)\;\mu^+_x(\partial_{e'}\XX)\;\Dirac_{[x]}
\label{eq:calcpartielmBMprequel}\\ &=
\sum_{[x]\in \Ga\bs V \XX} \frac{1}{|\Ga_x|}\;\Big(
\big(\sum_{e\in E\XX\;:\; o(e)=x} \mu_x(\partial_e\XX)\;\big)^2-
\sum_{e\in E\XX\;:\; o(e)=x}\mu_x(\partial_{e}\XX)^2\Big)\;\Dirac_{[x]}
\label{eq:calcpartielmBM}\:.
\end{align}
This gives the first equality  of Assertion (1).

\medskip Let us prove the second equality of Assertion (1).  By
homogeneity, we assume that $\|\mu_{x_0}\|=1$ and we will prove that
$$
\|m_{\rm BM}\|= \sum_{[x]\in \Ga\bs V \XX} \frac{c_x}{|\Ga_x|}\,.
$$ 
Let $N\in \NN -\{0\}$ be such that $p_{n+N}=p_n$ for every $n\in\NN$
and $p_{N-n}=p_n$ for every $n\in \{0,\dots,N\}$, which exists since
$\XX$ is assumed to be uniform.  Then the automorphism group
$\Aut(\XX)$ of the simplicial tree $\XX$ acts transitively on the set
of vertices at distance a multiple of $N$ from $x_0$. Hence for every
$x\in V\XX$, the distance $r_x = d(x,\Aut(\XX)x_0)$ belongs to
$\big\{0,1,\dots,\lfloor \frac{N}{2} \rfloor\big\}$, and there exist
$\ga_x,\ga'_x\in \Aut(\XX)$ such that 
$$
d(x, \ga_x x_0)= r_x, \;\;\;x\in
[\ga_x x_0,\ga'_x x_0]\;\;\;{\rm and}\;\;\; d(\ga_x x_0,\ga'_x x_0)=N\;.
$$ 
The map $x\mapsto r_x$ is constant on the orbits of $\Ga$ in
$V\XX$,\footnote{In fact, it is constant on the orbits of $\Aut(\XX)$.}
hence so is the map $x\mapsto c_x$, and thus the right hand side of
Equation \eqref{eq:formulgenemBM} is well defined.

Since the family $(\mu_x^{\rm Haus})_{x\in V\XX}$ of Hausdorff
measures of the visual distances $(\partial_\infty X, d_x)$ is
invariant under any element of $\Aut(\XX)$, since $\Ga$ is a lattice
and by Proposition \ref{prop:uniflatmBMfinie}, we have $\delta_\Ga =
h_\XX$ and $\ga_*\mu_x=\mu_{\ga x}$ for all $x\in V \XX$ and
$\ga\in\Aut(\XX)$.

Since $(\XX,x_0)$ is spherically symmetric, and since $\mu_{x_0}$ is a
probability measure, we have by induction, for every $e\in E\XX$
pointing away from $x_0$ with $d(x_0,o(e))=n$,
\begin{equation}\label{eq:massatroot}
\mu_{x_0}(\partial_e\XX)=\frac{1}{(p_0+1)\,p_1\,\dots\,p_n}
\end{equation}
if $n\neq 0$, and $\mu_{x_0}(\partial_e\XX)=\frac{1}{p_0+1}$
otherwise.

For every fixed $x\in V\XX$, let us now compute $\mu_x(\partial_e\XX)$
for every edge $e$ of $\XX$ with origin $x$.  Let $\ga=\ga_x,
\ga'=\ga'_x \in \Aut(\XX)$ be as above. By the spherical transitivity,
we may assume that $e$ or $\overline{e}$ belongs to the edge path from
$\ga x_0$ to $\ga'x_0$.

\begin{center}
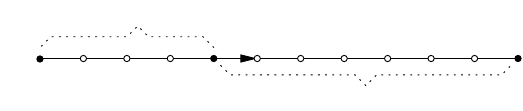
\end{center}

There are two cases to consider.

\medskip
\noindent{\bf Case 1: } Assume first that $e$ points away from $\ga
x_0$.  There are $p_0+1$ such edges starting from $x$ if $r_x=0$, and
$p_{r_x}$ otherwise. By Equation \eqref{eq:massatroot} and by invariance
under $\Aut(\XX)$ of $(\mu_x^{\rm Haus})_{x\in V\XX}$, we have
$$
\mu_{\ga x_0}(\partial_e\XX)=\frac{1}{(p_0+1)\,p_1\,\dots\,p_{r_x}}\;,
$$
with the convention that the denominator is $p_0+1$ if $r_x=0$. Since
the map $\xi\mapsto \beta_\xi(x,\ga x_0)$ is constant with value $-r_x$
on $\partial_e\XX$, and by the quasi-invariance property of the
Patterson density (see Equation \eqref{eq:quasinivarPatdens}), we have
$$
\mu_{x}(\partial_e\XX)=e^{-\delta_\Ga(-r_x)}\mu_{\ga x_0}(\partial_e\XX)=
\frac{e^{r_x h_\XX}}{(p_0+1)\,p_1\,\dots\,p_{r_x}}\;,
$$
with the same convention as above.

\medskip
\noindent
{\bf Case 2: } Assume now  that $e$ points towards $\ga
x_0$. This implies that $r_x\geq 1$, and there is one and only one
such edge starting from $x$. Then as above we have
$$
\mu_{\ga' x_0}(\partial_e\XX)=\frac{1}{(p_N+1)\,p_{N-1}\,\dots\,p_{r_x}}\;,
$$
and
$$
\mu_{x}(\partial_e\XX)=e^{-\delta_\Ga(-(N-r_x))}\mu_{\ga' x_0}(\partial_e\XX)=
\frac{e^{(N-r_x) h_\XX}}{(p_N+1)\,p_{N-1}\,\dots\,p_{r_x}}\;.
$$

For every $x\in V\XX$, let 
\begin{equation}\label{eq:Cx}
C_x=\Big(\sum_{e\in E\XX\;:\; o(e)=x} \mu_x(\partial_e\XX)\Big)^2-
\sum_{e\in E\XX\;:\; o(e)=x}\mu_x(\partial_{e}\XX)^2\;.
\end{equation}
If $r_x\neq 0$, since the stabiliser of $\ga x_0$ in $\Aut(\XX)$ acts
transitively on the $p_{r_x}$ edges with origin $x$ pointing away from
$\ga x_0$, since $e^{Nh_\XX}=p_0p_1\dots p_{N-1}$ and $p_N=p_0$, we
have
\begin{align*}
C_x=\;\;&
\Big(\,p_{r_x}\frac{e^{r_x h_\XX}}{(p_0+1)\,p_1\,\dots\,p_{r_x}}
+\frac{e^{(N-r_x) h_\XX}}{(p_N+1)\,p_{N-1}\,\dots\,p_{r_x}}\;\Big)^2
\\ &-\Big(\,p_{r_x}\big(\frac{e^{r_x h_\XX}}{(p_0+1)\,p_1\,\dots\,p_{r_x}}\big)^2
+\big(\frac{e^{(N-r_x) h_\XX}}{(p_N+1)\,p_{N-1}\,\dots\,p_{r_x}}\big)^2\;\Big)
\\ =\;\;&
\frac{{(p_{r_x}}^2-p_{r_x})\,e^{2\,r_x h_\XX}}{(p_0+1)^2\,{p_1}^2\,\dots\,{p_{r_x}}^2}+
\frac{2\,p_{r_x}e^{N\, h_\XX}}{(p_0+1)\,p_1\,\dots\,p_{r_x}p_{r_x}\dots p_{N-1}(p_N+1)}
\\ =\;\;&
\frac{{(p_{r_x}}-1)\,e^{2\,r_x h_\XX}}
{(p_0+1)^2\,{p_1}^2\,\dots\,{p_{r_x-1}}^2p_{r_x}}+
\frac{2\,p_{0}}{(p_0+1)^2}=c_x\;.
\end{align*}
If $r_x=0$, we have
$$
C_x=\Big((p_0+1)\frac{1}{p_0+1}\Big)^2-(p_0+1)\Big(\frac{1}{p_0+1}\Big)^2
=\frac{p_0}{p_0+1}=c_x\;.
$$ 
Assertion (1) now follows from Equation \eqref{eq:calcpartielmBM}.

\medskip 
Let us prove Assertion (2).  Note that $\XX=\XX_{p,\,q}$ is
spherically symmetric with respect to any vertex of $\XX$, and that,
by Equation \eqref{eq:lyonshdim},
$$
h_\XX=\frac{1}{2}\ln(pq)\;.
$$ 
Let $e$ be an edge of $\XX$, with $x=o(e)\in V_p\XX$ and $y=t(e)\in
V_q\XX$. For every $z\in V\XX$, we define $C_z$ as in Equation
\eqref{eq:Cx}.

Note that by homogeneity, we have $C_z=C_x$ and $\|\mu_z\|=\|\mu_x\|$
for all $z\in V_p\XX$, as well as $C_z=C_y$ and $\|\mu_z\|=\|\mu_y\|$
for all $z\in V_q\XX$.  Hence the normalisation of the Patterson
density as in the statement of Assertion (2) is possible.  By the
spherical symmetry at $x$, and the normalisation of the measure, we
have $\mu_x(\partial_e\XX)= \frac{1}{\sqrt{p}}$ and
$\mu_x(\partial_{\overline{e}}\XX) =\sqrt{p}$.  Therefore
\begin{align*}
\|\mu_y\| &=\mu_y(\partial_e\XX)+\mu_y(\partial_{\overline{e}}\XX) 
=e^{h_\XX}\mu_x(\partial_e\XX)+e^{-h_\XX}\mu_x(\partial_{\overline{e}}\XX)
\\ & =\sqrt{pq}\,\frac{1}{\sqrt{p}}+\frac{1}{\sqrt{pq}}\sqrt{p}=
\frac{q+1}{\sqrt{q}}\;.
\end{align*}
This symmetry in the values of $\|\mu_y\|$ and $\|\mu_x\|$ explains
the choice of our normalisa\-tion of the Patterson density. We have
$$
C_x=\|\mu_x\|^2-(p+1)\Big(\frac{\|\mu_x\|}{p+1}\Big)^2
=\frac{p}{p+1}\|\mu_x\|^2=p+1
$$ 
and similarly $C_y=\frac{q}{q+1}\|\mu_y\|^2=q+1$. This proves the
second equality in Equation \eqref{eq:relatmBMvolbihomog}, by the
first equation of Assertion (1).

In order to prove that $(T\pi)_*m_{\rm BM}= \Tvol_{\Ga\dbs\XX}$,
we now partition $\Ga\bs \G \XX$ as
$$
\bigcup_{[e]\in \Ga\bs E \XX}\;\Ga\bs\big\{\ell \in  \G \XX\;:\;
\ell(0)=\pi(o(e)),\;\ell(1)=\pi(t(e))\big\}\;.
$$ 
Using on every element of this partition Hopf's decomposition with
respect to the basepoint $o(e)$, we have, by a ramified covering
argument already used in the proof of the second part of Proposition
\ref{prop:majoGibbsL2},
$$
(T\pi)_*m_{\rm BM}=\sum_{[e]\in \Ga\bs E \XX} \frac{1}{|\Ga_e|}\;
\mu_{o(e)}(\partial_\infty X-\partial_e\XX)\;
\mu_{o(e)}(\partial_e\XX) \;\Dirac_{[e]}\:.
$$
Since  $\mu_{o(e)}(\partial_e\XX)= e^{-h_\XX}
\mu_{t(e)}(\partial_e\XX)$ and by homogeneity, we have
\begin{align*}
(T\pi)_*m_{\rm BM}&=\sum_{[e]\in \Ga\bs E \XX} \frac{1}{|\Ga_e|}\;
\frac{\deg o(e)-1}{\deg o(e)}\;\|\mu_{o(e)}\|\;
\frac{\deg t(e)-1}{\deg t(e)}\;\|\mu_{t(e)}\|\;e^{-h_\XX}
 \;\Dirac_{[e]}\\ &
=\sum_{[e]\in \Ga\bs E \XX} \frac{1}{|\Ga_e|}\;
\frac{\|\mu_{o(e)}\|\;\|\mu_{t(e)}\|\;\sqrt{pq}}{(p+1)(q+1)}
 \;\Dirac_{[e]}=\Tvol_{\Ga\dbs\XX}\:.
\end{align*}
The first equality of Equation \eqref{eq:relatmBMvolbihomog} follows,
since pushforwards of measures preserve the total mass.

\medskip Finally, the last claim of Assertion (3) of Proposition
\ref{prop:computBM} follows from Equation \eqref{eq:formulgenemBM},
since $c_x=\frac{q}{q+1}$ for every $x\in V\XX_q$ (or by taking $q=p$
in Equation \eqref{eq:relatmBMvolbihomog} and by renormalising). The
first claim of Assertion (3) follows from the first claim of Assertion
(2), by using Equation \eqref{eq:volTvol} and renormalising.  
\cqfd

\brema (1) In particular, if $\XX=\XX_{q}$ is regular, if the Patterson
density is normalised to be a family of probability measures and if $\Ga$
is torsion free, then $\pi_*m_{\rm BM}$ is $\frac{q}{q+1}$ times the
counting measure on $\Ga\bs V\XX$. In this case, Equation
\eqref{eq:relatmBMvolhomog} is given by \cite[Rem.~2]{CooPap97b}.

\smallskip
\noindent(2) If $\XX=\XX_{p,\,q}$ is biregular with $p\neq q$, then
$\pi_*m_{\rm BM}$ is not proportional to $\vol_{\Ga\dbs\XX}$.  In
particular, if $\Ga$ is torsion free and if the Patterson density is
normalised to be a family of probability measures, then $\pi_*m_{\rm
  BM}$ is the sum of $\frac{p}{p+1}$ times the counting measure on
$\Ga\bs V_p\XX$ and $\frac{q}{q+1}$ times the counting measure on
$\Ga\bs V_q\XX$.

This statement is coherent with the well-known fact that in pinched but
variable curvature, the Bowen-Margulis measure is generally not
absolutely continuous with respect to Sasaki's Riemannian measure on
the unit tangent bundle (they would then be proportional by ergodicity
of the geodesic flow in the lattice case).  
\erema

\section{Computations of skinning measures for 
simplicial trees}
\label{subsec:skinmeascomputtree}

We now give examples of computations of the total mass of skinning
measures (for zero potentials), after introducing some notation. Let
$\XX$ be a locally finite simplicial tree without terminal vertices,
and let $\Ga$ be a discrete subgroup of $\Aut(\XX)$.

\medskip
For every simplicial subtree $\DD$ of $\XX$, we define the {\it
  boundary}\index{boundary} $\partial V\DD$ of $V\DD$ in $\XX$ as
$$   
\gls{boundaryVD}=\{x\in V\DD\;:\;\exists\; e\in E\XX,\;\; o(e)=x, \;
t(e)\notin V\DD\}\;.
$$ 
The {\em boundary}\index{boundary} $\gls{boundaryD}$ of $\DD$ is
the maximal subgraph (which might be non connected) of $\XX$ with set
of vertices $\partial V\DD$. It is contained in $\DD$. The stabiliser
$\Ga_\DD$ of $\DD$ in $\Ga$ acts discretely on $\partial \DD$.

For every $x\in V \XX$, we define the {\em codegree}\index{codegree}
of $x$ in $\DD$ as $\codeg_\DD(x)=0$ if $x\notin\DD$ and otherwise
$$
\gls{codegre} = \deg_\XX(x) - \deg_\DD(x)\;.
$$ 
Note that $\codeg_\DD(x)=0$ if $x\notin \partial V\DD$, and that
the codegree $\codeg_{\N_1\DD} (x)$ of $x\in V\XX$ in the
$1$-neighbourhood $\N_1\DD$ of $\DD$ is equal to $0$ unless $x$ lies
in the boundary of $\N_1\DD$, in which case it is equal to
$\deg_\XX(x)-1$.

Let $\D=(\DD_i)_{i\in I}$ be a locally finite $\Ga$-equivariant family
of simplicial subtrees of $\XX$, and let $x\in V\XX$. We define the
{\em multiplicity}\index{multiplicity}\footnote{See Section
  \ref{subsec:downstairs} for explanations on the terminology.} of $x$
in (the boundary of) $\D$ as (see Section \ref{subsec:equivfammult}
for the definition of $\sim_\D$)
$$
\gls{multiplicitysimplicial}=
\frac{\card\{i\in I/_{\sim_\D}\;:\; x\in \partial V\DD_i\}}{|\Ga_x|}\;.
$$ 
The numerator and the denominator are finite by the local
finiteness of the family $\D$ and the discreteness of $\Ga$, and they
depend only on the orbit of $x$ under $\Ga$.  Note that if $\DD$ is a
simplicial subtree of $\XX$ which is {\em precisely invariant}%
\index{precisely invariant} under $\Ga$ (that is, whenever $\ga\in\Ga$
is such that $\DD\cap\ga\DD$ is nonempty, then $\ga$ belongs to the
stabiliser $\Ga_\DD$ of $\DD$ in $\Ga$), if $\D=(\ga\DD)_{\ga \in
  \Ga/\Ga_\DD}$, and if $x\in \partial V\DD$, then
$$
m_{\D}(x)=\frac{1}{|\Ga_x|}\;.
$$   
In particular, if furthermore $\Ga$ is torsion free, then
$m_{\D}(x)=1$ if $x\in \partial V\DD$, and $m_{\D}(x)=0$ otherwise.

\bexem 
Let $\G$ be a connected graph without vertices of degree $\leq 2$ and
let $\XX$ be its universal cover, with covering group $\Ga$. If $C$
is a cycle in $\G=\Ga\bs\XX$ and if $\D$ is the family of geodesic
lines in $\XX$ lifting $C$, then $m_{\D}(x)=1$ for all $x\in V\XX$
whose image in $\G=\Ga\bs\XX$ belongs to $C$ if $C$ is a {\it
  simple}\index{simple} cycle (that is, if $C$ passes through no
vertex twice).  
\eexem

We define the {\em codegree}\index{codegree} of $x$ in $\D$ as
$$
\gls{codegrecalD} = \sum_{i\in I/\sim_\D} \codeg_{\DD_i}(x)\;,
$$
which is well defined as $\codeg_{\DD_i}(x)$ depends only on the class
of $i\in I$ modulo $\sim_\D$. Note that
\begin{equation}\label{eq:codegk}
\codeg_\D(x)= (\deg_\XX x -k)\,|\Ga_x|\,m_{\D}(x)
\end{equation}
if $\deg_{\DD_i}(x)=k$ for every $x\in \partial V\DD_i$ and $i\in
I$. If every vertex of $\XX$ has degree at least $3$, this is in
particular the case with $k=2$ if $\DD_i$ is a line for all $i\in I$
and with $k=1$ if $\DD_i$ is a horoball for all $i\in I$.

We will say that a simplicial subtree $\DD$ of $\XX$, with stabiliser
$\Ga_\DD$ in $\Ga$, is {\em almost precisely invariant}\index{almost
  precisely invariant}\index{precisely invariant!almost} if there
exists $N\in\NN$ such that for every $x\in V\DD$, the number
of $\ga\in\Ga/\Ga_\DD$ such that $x\in \ga V\DD$ is at most
$N$. It follows from this property that if $\D=(\ga\DD)_{\ga\in\Ga}$,
then $\D$ is locally finite and $\codeg_\D(x)\leq N\codeg_\DD(x)$ for
every $x\in V\XX$.

\bprop \label{prop:computskin} Assume that $\XX$ is a regular or
biregular simplicial tree with degrees at least $3$, and that $\Ga$ is
a lattice of $\XX$.

\smallskip
\noindent(1) For every simplicial subtree $\DD$ of $\XX$, we have
$$
\pi_*\wt \sigma^\pm_\DD= \sum_{x\in V\XX} \;\;
\frac{\|\mu_x\|\;\codeg_\DD(x)}{\deg_\XX(x)}\;\;\Dirac_{x}\;.
$$

\smallskip\noindent(2) If $\D=(\DD_i)_{i\in I}$ is a locally finite
$\Ga$-equivariant family of simplicial subtrees of $\XX$, then
$$
\pi_* \sigma^\pm_\D= \sum_{[x]\in \Ga\bs V\XX}\;
\frac{\|\mu_x\|\;\codeg_\D(x)}{|\Ga_x|\;\deg_\XX(x)}\;\;\Dirac_{[x]}\;.
$$

\smallskip\noindent(3) Let $k\in\NN$ and let $\DD$ be a simplicial
subtree of $\XX$ such that $\deg_{\DD}(x)=k$ for every $x\in \partial
V\DD$ and the $\Ga$-equivariant family $\D=(\ga\DD)_{\Ga/\Ga_\DD}$ is
locally finite. Then
$$
\pi_* \sigma^\pm_\D= \sum_{\Ga_\DD y\,\in \,\Ga_\DD\bs \partial V\DD}\;
\frac{\|\mu_y\|\;(\deg_\XX(y)-k)}{|(\Ga_\DD)_y|\;\deg_\XX(y)}\;\;
\Dirac_{\Ga y}\;.
$$

\smallskip\noindent(4) If $\DD$ is a simplicial subtree of $\XX$ such
that the $\Ga$-equivariant family $\D=(\ga\DD)_{\ga\in\Ga}$ is
locally finite, then the skinning measure $\sigma^\pm_\D$ is finite if
and only if the graph of groups $\Ga_\DD\dbs\partial\DD $ has
finite volume.  
\eprop

Before proving Proposition \ref{prop:computskin}, let us give some
immediate consequences of its Assertion (3).  If $\XX= \XX_{p,\,q}$ is
biregular of degrees $p+1$ and $q+1$, let $V\XX=V_p\XX\,\sqcup\,
V_q\XX$ be the corresponding partition of the set of vertices of $\XX$
and, for $r\in\{p,q\}$, let $\partial_r\DD$ be the edgeless graph with
set of vertices $\partial V\DD\cap V_r\XX$.

\bcoro\label{coro:biregularskinningmasses} 
Assume that $(\XX,\Ga)$ is as in Proposition
\ref{prop:computskin}. Let $\DD$ be a simplicial subtree of $\XX$ such
that the $\Ga$-equivariant family $\D=(\ga\DD)_{\ga\in\Ga/\Ga_\DD}$ is
locally finite.

\smallskip\noindent(1) If $\XX= \XX_{p,\,q}$ is biregular of degrees
$p+1$ and $q+1$ and if the Patterson density $(\mu_x)_{x\in V \XX}$ of
$\Ga$ is normalised so that $\|\mu_{x}\|=
\frac{\deg_\XX(x)}{\sqrt{\deg_\XX(x)-1}}$ for all $x\in V\XX$, then

\smallskip
\noindent
$\bullet$~ if $\DD$ is a horoball,
$$
\|\sigma^\pm_\D\|=\sqrt{p}\;\Vol(\Ga_\DD\dbs \partial_p\DD)+
\sqrt{q}\;\Vol(\Ga_\DD\dbs \partial_q\DD)\;,
$$
$\bullet$~ if $\DD$ is a line,
\begin{equation}\label{eq:massskinlinebiregular}
\|\sigma^\pm_\D\|=\frac{p-1}{\sqrt{p}}\;\Vol(\Ga_\DD\dbs \partial_p\DD)+
\frac{q-1}{\sqrt{q}}\;\Vol(\Ga_\DD\dbs \partial_q\DD)\;.
\end{equation}

\smallskip\noindent(2) If $\XX=\XX_q$ is the regular tree of degree
$q+1$ and if the Patterson measures $(\mu_x)_{x\in V\XX}$ are
normalised to be probability measures, then

\smallskip
\noindent
$\bullet$~ if $\DD$ is a horoball,
\begin{equation}\label{eq:massskinhorosph}
\|\sigma^\pm_\D\|= \frac{q}{q+1}\Vol(\Ga_\DD\dbs\partial \DD)
\end{equation}
$\bullet$~ if $\DD$ is a line,
\begin{equation}\label{eq:massskinline}
\|\sigma^\pm_\D\|= \frac{q-1}{q+1}\Vol(\Ga_\DD\dbs\DD)\;.\;\;\; \Box
\end{equation}
\ecoro

\medskip
\noindent{\bf Proof of Proposition \ref{prop:computskin}. } 
(1) We may partition the outer/inner unit normal bundle $\normalpm
\DD$ of $\DD$ according to the first/last edge of the elements in
$\normalpm \DD$. On each of the elements of this partition, for the
computation of the skinning measures using its definition and its
independence of the basepoint (see Section
\ref{subsec:skinningmeasures}), we take as basepoint the
initial/terminal point of the corresponding edge. Since $\DD$ is a
simplicial tree, note that for every $e\in E\XX$ such that $o(e)\in
V\DD$, we have $e\in E\DD$ if and only if $t(e)\in V\DD$. Thus, we
have
\begin{align*}
\pi_*\wt \sigma^+_\DD&=\sum_{e\in E\XX\;:\; o(e)\in V\DD,\; t(e)\notin V\DD}
\mu_{o(e)}(\partial_e\XX) \;\Dirac_{o(e)}\\ &
=\sum_{x\in \partial V\DD} \Big( \sum_{e\in E\XX\;:\; o(e)=x,\; t(e)\notin V\DD}
\mu_{x}(\partial_e\XX)\Big) \;\Dirac_{x}\;.
\end{align*}
and similarly
\begin{align*}
\pi_*\wt \sigma^-_\DD&=\sum_{e\in E\XX\;:\; t(e)\in V\DD,\; o(e)\notin V\DD}
\mu_{t(e)}(\partial_{\overline{e}}\XX) \;\Dirac_{t(e)}\\ &
=\sum_{x\in \partial V\DD} \Big( \sum_{e\in E\XX\;:\; o(e)=x,\; t(e)\notin V\DD}
\mu_{x}(\partial_e\XX)\Big) \;\Dirac_{x}\;.
\end{align*}
As in the proof of Proposition \ref{prop:computBM} (2), since $\XX$ is
spherically homogeneous around each point and since $\Ga$ is a lattice
(so that the Patterson density is $\Aut(\XX)$-equivariant, see
Proposition \ref{prop:uniflatmBMfinie}), we have $\mu_x(\partial_e\XX)
= \frac{\|\mu_x\|} {\deg_\XX(x)}$ for all $x\in V\XX$ and $e\in E\XX$
with $o(e)=x$.  Assertion (1) 
follows, since $\sum_{e\in E\XX\;:\; o(e)=x,\; t(e)\notin V\DD} \;1\;=
\codeg_\DD(x)$ if $x\in \partial V\DD$ and $\codeg_\DD(x)=0$
otherwise.

\medskip\noindent (2) By the definition\footnote{See Section
\ref{subsec:equivfammult}.} of the skinning measures associated with
$\Ga$-equi\-va\-riant families, we have $\wt\sigma^\pm_{\D}=\sum_{i\in
  I/_\sim} \wt \sigma^+_{\DD_i}$, where $\sim\;=\;\sim_\D$. Hence by
Assertion (1)
\begin{align*}
\pi_*\wt\sigma^\pm_{\D}&=\sum_{i\in I/_\sim} \wt \pi_*\sigma^+_{\DD_i}
=\sum_{i\in I/\sim }\;\sum_{x\in V\XX} \;\;
\frac{\|\mu_x\|\;\codeg_{\DD_i}(x)}{\deg_\XX(x)}\;\;\Dirac_{x}\\ &
=\sum_{x\in V\XX} \Big(\sum_{i\in I/\sim }\codeg_{\DD_i}(x)\Big)\;
\frac{\|\mu_x\|}{\deg_\XX(x)}\;\;\Dirac_{x}
=\sum_{x\in V\XX}\;
\frac{\|\mu_x\|\;\codeg_\D(x)}{\deg_\XX(x)}\;\;\Dirac_{x}\;.
\end{align*}
By the definition of the measure induced in $\Ga\bs V\XX$ when $\Ga$
may have torsion (see for instance \cite[\S 2.6]{PauPolSha15} and the
beginning of Chapter \ref{sec:equidcountdownstairs}), Assertion (2)
follows.

\medskip\noindent  (3) It follows from Assertion (2) and from Equation
\eqref{eq:codegk} that
$$
\pi_* \sigma^\pm_\D= \sum_{[x]\in \Ga\bs V\XX}\;
\frac{\deg_\XX(x)-k}{\deg_\XX(x)}\;\|\mu_x\|\;m_{\D}(x)\;\;\Dirac_{[x]}\;.
$$
For every $x\in V\XX$, by the definition of $m_{\D}(x)$, we have, by
partitioning $\partial V\DD$ into its orbits under $\Ga_\DD$,
\begin{align*}
m_{\D}(x)&=\frac{1}{|\Ga_x|}
\card\{\ga\in\Ga_\DD\bs\Ga\;:\; \ga x\in\partial V\DD\}
\\ & =\frac{1}{|\Ga_x|}\sum_{\Ga_\DD y\,\in\,\Ga_\DD\bs \partial V\DD}
\card\{\ga\in\Ga_\DD\bs\Ga\;:\; \Ga_\DD \ga x= \Ga_\DD y\}
\\ & =\frac{1}{|\Ga_x|}\sum_{\Ga_\DD y\,\in\,\Ga_\DD\bs \partial V\DD,\;\Ga x=\Ga y}
\card\{\ga\in\Ga_\DD\bs\Ga\;:\; \Ga_\DD \ga y= \Ga_\DD y\}
\\ & =\frac{1}{|\Ga_x|}\sum_{\Ga_\DD y\,\in\,\Ga_\DD\bs \partial V\DD,\;\Ga x=\Ga y}
[\Ga_y:(\Ga_\DD)_y]=\sum_{\Ga_\DD y\,\in\,\Ga_\DD\bs \partial V\DD,\;\Ga x=\Ga y}
\frac{1}{|(\Ga_\DD)_y|}\;.
\end{align*}
This proves Assertion (3), since $\sum_{[x]\in\Ga\bs V\XX,\;\Ga x=\Ga
  y} \Delta_{[x]}=\Delta_{\Ga y}$.

\medskip\noindent (4) It follows from Assertion (2) that 
$$
\|\sigma^\pm_\D\|=\sum_{[x]\in \Ga\bs V\XX}\;
\frac{\|\mu_x\|\;\codeg_\D(x)}{|\Ga_x|\;\deg_\XX(x)}\;.
$$
Note that for every $x\in\partial V\DD$, we have
$$
|\Ga_x|\,m_\D(x)\leq \codeg_\D(x)\leq \deg_\XX(x)\;|\Ga_x|\,m_\D(x)\;.
$$
Let $m=\min_{x\in V\XX} \|\mu_x\|$ and $M=\max_{x\in V\XX} \|\mu_x\|$,
which are positive and finite, as the total mass of the Patterson
measures takes at most two values, since $\Ga$ is a lattice and $\XX$
is biregular.  By arguments similar to those in the proof of
Assertion (3), we hence have
$$
\frac{m}{\min_{x\in V\XX}\deg_\XX(x)}\;\Vol(\Ga_\DD\dbs\partial\DD)
\leq\;\|\sigma^\pm_\D\|\;\leq
M\;\Vol(\Ga_\DD\dbs\partial\DD)
$$
The result follows. \cqfd

\bigskip 
We now give a formula for the skinning measure (with zero potential)
of a geodesic line in the simplicial tree $\XX$, using\footnote{See
  the definitions of Hamenst\"adt's distance and measure in Sections
  \ref{subsec:lines} and \ref{subsec:skinningmeasures} respectively.}
Hamenst\"adt's distance $d_\H$ and measure $\mu_\H$ associated with a
fixed horoball $\H$ in $\XX$. This expression for the skinning measure
will be useful in Part \ref{sect:arithappli}.

\blemm\label{lem:hamline} Let $\H$ be a horoball in $\XX$ centred at
a point $\xi\in\partial_\infty X$.  Let $L$ be a geodesic line in
$\XX$ with endpoints $L_\pm\in\partial_\infty X-\{\xi\}$.  Then for
all $\rho\in\normalout L$ such that $\rho_+\ne\xi$,
$$
d\wt\sigma^+_L(\rho)=\frac{d_{\H}(L_+,L_-)^{\delta_\Ga}}
{d_\H(\rho_+,L_-)^{\delta_\Ga}\,d_\H(\rho_+,L_+)^{\delta_\Ga}}\,d\mu_\H(\rho_+)\,.
$$
\elemm

\dem 
By Equations \eqref{eq:shrinkhoroballdist} and
\eqref{eq:shrinkhoroballmeas}, the power ${d_{\H}}^{\delta_\Ga}$ of
the distance and the measure $\mu_\H$ scale by the same factor when
the horoball is replaced by another one centred at the same
point. Thus, we can assume in the proof that $L$ does not intersect
the interior of $\H$.

Fix $\rho\in\normalout L$ such that $\rho_+\ne\xi$. Let $y$ be the
closest point to $\xi$ on $L$, let $x_0$ be the closest point to $L$
on $\H$, and let $z$ be the closest point to $\xi$ on
$\rho([0,+\infty[)$.  Let $t\mapsto x_t$ be the geodesic ray starting
from $x_0$ at time $t=0$ and converging to $\xi$.  When $t$ is large
enough, the points $\rho_+$, $z$, $x_t$ and $\xi$ are in this order on
the geodesic line $]\rho_+,\xi[$.

We have, by the definition in Equation \eqref{eq:defiskin} of the
skinning measure, the cocycle property of the Busemann function,
Equation \eqref{eq:changemoinsplus} and the definition of $z$,
\begin{align*}
  d\wt\sigma^+_L(\rho)&
  =e^{\delta_\Ga\,\beta_{\rho_+}(x_t,\,\rho(0))}d\mu_{x_t}(\rho_+)
  =e^{\delta_\Ga\,\beta_{\rho_+}(x_t,\,z)+
    \delta_\Ga\,\beta_{\rho_+}(z,\,\rho(0))}d\mu_{x_t} (\rho_+) \\
  &=e^{-\delta_\Ga\,\beta_{\xi}(x_t,\,z)-
    \delta_\Ga\,d(z,\,\rho(0))}d\mu_{x_t}(\rho_+)\\
  &= e^{\delta_\Ga\,t-\delta_\Ga\,\beta_{\xi}(x_0,\,z)-
    \delta_\Ga\,d(z,\,\rho(0))}d\mu_{x_t}(\rho_+)\;,
\end{align*}
and, by the definition of Hamenst\"adt's measure $\mu_\H$ (see Equation
\eqref{eq:defiHamensmeasure}),
$$
d\mu_\H(\rho_+)=e^{\delta_\Ga\,t}d\mu_{x_t}(\rho_+)\,.
$$

\begin{center}
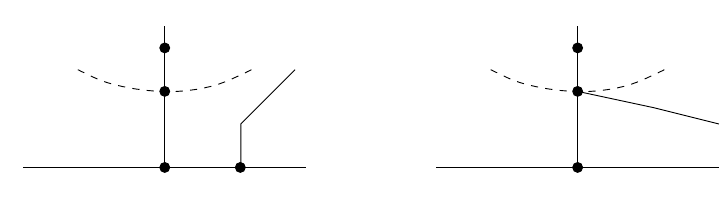
\end{center}

\noindent {\bf Case 1: } Assume first that $\rho(0)\ne y$.  We may
assume that $\rho(0)\in [y,L_+[$. Then $z=\rho(0)$ and $z$ is the
closest point to $\H$ on the geodesic line $]L_+,\rho_+[$.  Thus
$d_\H(L_-,L_+)= d_\H(L_-,\rho_+)$ and $d_\H(L_+,\rho_+)=
e^{-d(z,\,x_0)}= e^{\beta_{\xi}(x_0,\,z)}$, and the claim follows.

\medskip
\noindent {\bf Case 2: } Assume now that $y=\rho(0)$. Then $[y,z]=
[y,\xi[\;\cap\, [y,\rho_+[$, and we may assume that $x_0=z$ up to
adjusting the horoball $\H$ while keeping its point at infinity. Thus
$d_\H(L_-,L_+)=e^{-d(y,\,x_0)}= e^{-d(z,\,\rho(0))}$ and
$d_\H(L_-,\rho_+)=d_\H(L_+,\rho_+)=1$, and the claim follows.  
\cqfd

\chapter{Rate of mixing for the geodesic flow} 
\label{sec:mixingrate}

Let $X$, $x_0$, $\Ga$, $\wt F$, $(\mu^\pm_x)_{x\in X}$ be as in the
beginning of Chapter \ref{sect:skinning}, and $\wt F^\pm$, $F^\pm$,
$\delta= \delta_{\Ga,\,F^\pm}<\infty$, $\wt m_F$, $m_F$ the associated
notation. In this Chapter, we start by collecting in Section
\ref{subsec:mixingratemanifolds} known results on the rate of mixing
of the geodesic flow for manifolds. The main part of the Chapter then
consists in proving analogous bounds for the discrete time and
continuous time geodesic flow for quotient spaces of simplicial and
metric trees respectively.

We define $\gls{gibbsmeasuredownrenorm}= \frac{m_F}{\|m_F\|}$ when the
Gibbs measure is finite. Recall that this measure is nonzero since $\Ga$ is
nonelementary.

Let $\alpha\in\;]0,1]$.\footnote{We refer to Section
\ref{subsec:holdercont} for the definition of the Banach space
$\C_{\rm b}^\alpha (Z)$ of bounded $\alpha$-H\"older-continuous
functions on a metric space $Z$.} We will say that the (continuous
time) geodesic flow on $\Ga\backslash \G X$ is {\em
exponentially mixing for the $\alpha$-H\"older regularity}%
\index{exponentially mixing!Holder@H\"older} or that it has {\em
  exponential decay of $\alpha$-H\"older
  correlations}\index{exponential decay!Holder@H\"older}%
\index{decay of correlations!exponential} for the potential $F$ if
there exist $C,\kappa >0$ such that for all $\phi,\psi\in \C_{\rm
  b}^\alpha (\Ga\backslash \G X)$ and $t\in\RR$, we have
\begin{equation}\label{eq:defiholderexponentialdecay}
\Big|\int_{\Ga\backslash \G X}\phi\circ\flow{-t}\;\psi\;d\overline{m_{F}}-
\int_{\Ga\backslash \G X}\phi\; d\overline{m_{F}}
\int_{\Ga\backslash \G X}\psi\;d\overline{m_{F}}\;\Big|
\le C\;e^{-\kappa|t|}\;\|\phi\|_\alpha\;\|\psi\|_\alpha\,,
\end{equation}
and that it is {\em polynomially mixing}\index{polynomially mixing}
for the $\alpha$-H\"older regularity or has {\em polynomial decay of
  $\alpha$-H\"older correlations}%
\index{decay of correlations!polynomial} if there exist $C>0$ and
$n\in\NN-\{0\}$ such that for all $\phi,\psi\in \C_{\rm b}^\alpha
(\Ga\backslash \G X)$ and $t\in\RR$, we have
$$
\Big|\int_{\Ga\backslash \G X}\phi\circ\flow{-t}\;\psi\;d\overline{m_{F}}-
\int_{\Ga\backslash \G X}\phi\; d\overline{m_{F}}
\int_{\Ga\backslash \G X}\psi\;d\overline{m_{F}}\;\Big|
\le C\;(1+|t|)^{-n}\;\|\phi\|_\alpha\;\|\psi\|_\alpha\;.
$$

\section{Rate of mixing for Riemannian manifolds}
\label{subsec:mixingratemanifolds}

When $X=\wt M$ is a complete simply connected Riemannian manifold with
pinched negative sectional curvature with bounded derivatives, then
the boundary at infinity of $\wt M$, the strong unstable, unstable,
stable, and strong stable foliations of $T^1\wt M$ are H\"older-smooth
and only H\"older-smooth in general.\footnote{See for instance
  \cite{Brin95} when $\wt M$ has a compact quotient (a result first
  proved by Anosov), and \cite[Theo.~7.3]{PauPolSha15}.}  Hence the
assumption of H\"older regularity on functions on $T^1\wt M$ is
appropriate for these manifolds.

The geodesic flow is known to have exponential decay of H\"older
correlations for compact manifolds $M=\Ga\bs\wt M$ when

\begin{itemize}
\item $M$ is two-dimensional and $F$ is any H\"older-continous
  potential by \cite{Dolgopyat98},

\item $M$ is $1/9$-pinched and $F=0$ by
\cite[Coro.~2.7]{GiuLivPol13},

\item $m_F$ is the Liouville measure by \cite{Liverani04}, see also
  \cite{Tsujii10}, \cite[Coro.~5]{NonZwo15} who give more precise
  estimates,
 
\item $M$ is locally symmetric and $F$ is any H\"older-continuous
  potential by \cite{Stoyanov11}, see also \cite{MohOh15}.
\end{itemize}

\bigskip 
When $\wt M$ is a symmetric space, then the boundary at infinity of
$\wt M$, the strong unstable, unstable, stable, and strong stable
foliations of $T^1\wt M$ are smooth.  Hence talking about leafwise
$\C^\ell$-smooth functions on $T^1\wt M$ makes sense. For every
$\ell\in \NN$, we will denote by $\gls{espacellsobolev}$ the vector
space of real-valued $\C^\ell$-smooth functions on the orbifold
$T^1M=\Ga\bs T^1\wt M$ with compact support in $T^1M$, and by
$\gls{normellsobolev}$ the Sobolev $W^{\ell,2}$-norm of any $\psi\in
\C_{\rm c}^\ell(T^1M)$.  This space consists of functions induced on
$T^1M$ by $\C^\ell$-smooth $\Ga$-invariant functions with compact
support on $T^1\wt M$.

Given $\ell\in \NN$, we will say that the geodesic flow on $T^1M$ is
{\em exponentially mixing for the $\ell$-Sobolev
  regularity}\index{exponentially mixing!Sobolev} (or that it has {\em
  exponential decay of $\ell$-Sobolev correlations})%
\index{exponential decay!Sobolev} for the potential $F$ if there exist
$c,\kappa>0$ such that for all $\phi,\psi\in \C_{\rm c}^\ell(T^1M)$ and all
$t\in\RR$, we have
$$
\Big|\int_{T^1M} \phi\circ \flow{-t}\;\psi\;d\overline{m_{F}}-\int_{T^1M}
\phi\;d\overline{m_{F}}\int_{T^1M} \psi\;d\overline{m_{F}}\;\Big|\leq
c\,e^{-\kappa |t|}\;\|\psi\|_\ell\;\|\phi\|_\ell\;.
$$
When $F=0$ and $\Ga$ is an arithmetic lattice in the isometry group of
$\wt M$ (the Gibbs measure then coincides, up to a multiplicative
constant, with the Liouville measure), this property, for some
$\ell\in\NN$, follows from \cite[Theorem~2.4.5]{KleMar96}, with the
help of \cite[Theorem 3.1]{Clozel03} to check its spectral gap
property, and of \cite[Lemma~3.1]{KleMar99} to deal with finite cover
problems.

\section{Rate of mixing for simplicial trees}
\label{subsec:mixingratesimpgraphs}

Let $\XX$ be a locally finite simplicial tree without terminal
vertices, with geometric realisation $X=|\XX|_1$, and $x_0\in
V\XX$. Let $\Ga$ be a nonelementary discrete subgroup of $\Aut(\XX)$
and let $\wt c: E\XX \ra\RR$ be a system of conductances for $\Ga$ on
$\XX$. Let $\wt F_c: T^1X \ra \RR$ be the associated potential of $\wt
c$, and $c: \Ga\bs E\XX \ra\RR$, $F_c: \Ga\bs T^1X \ra \RR$ the
quotient functions. Let $\delta_c= \delta_{\Ga,\,F_c}$ be the critical
exponent of $c$, assumed to be finite.\footnote{That is, to be
  $<+\infty$, since The critical exponentas it is $>-\infty$ by Lemma
  \ref{lem:proprielemcritexpo} (7).}  Let $(\mu^\pm_x)_{x\in V\XX}$ be
(normalised) Patterson densities on $\partial_\infty X$ for the pairs
$(\Ga,F^\pm_c)$, and let $\wt m_c=\wt m_{F_c}$ and $m_c=m_{F_c}$ be
the associated Gibbs measures on $\G\XX$ and $\Ga\bs \G\XX$.

In this Section, building on the end of Section
\ref{subsec:ergodictrees} concerning the mixing properties
themselves,\footnote{See Theorem \ref{theo:uniflatmBMmixing}.} we now
study the rate of mixing of the discrete time geodesic flow on $\Ga\bs
\G\XX$ for the Gibbs measure $m_c=m_{F_c}$, when it is mixing.

\medskip
Let $(Z,m,T)$ be a dynamical system with $(Z,m)$ a probability space
and $T:Z\ra Z$ a (not necessarily invertible) measure preserving
map. For all $n\in\NN$ and $\phi,\psi\in\LL^2(m)$, the (well-defined)
$n$-th {\it correlation coefficient}\index{correlation coefficient} of
$\phi,\psi$ is
$$
\gls{correlationcoeff}=
\int_{Z}\phi\circ T^n\;\psi\;dm-\int_{Z}\phi\; dm\;\int_{Z}\psi\;dm\;.
$$ 
Let $\alpha\in\;]0,1]$ and assume that $Z$ is a metric space
(endowed with its Borel $\sigma$-algebra). Similarly as for the case
of flows in the beginning of Chapter \ref{sec:mixingrate}, we
will say that the dynamical system $(Z,m,T)$ is {\em exponentially
  mixing for the $\alpha$-H\"older regularity}%
\index{exponentially mixing!Holder@H\"older} or that it has {\em
  exponential decay of $\alpha$-H\"older
  correlations}\index{exponential decay!Holder@H\"older} if there
exist $C,\kappa >0$ such that for all $\phi,\psi\in \C_{\rm b}
^\alpha(Z)$ and $n\in\NN$, we have
$$
|\operatorname{cov}_{m,\,n}(\phi,\psi)|
\le C\;e^{-\kappa\, n }\;\|\phi\|_\alpha\;\|\psi\|_\alpha\,.
$$ 
Note that this property is invariant under measure preserving
conjugations of dynamical systems by bilipschitz homeomorphisms.

\medskip
The main result of this Section is a simple criterion for the
exponential decay of correlation of the discrete time geodesic
flow on $\Ga\bs \G\XX$. 

We define $\gls{gibbsmeasuredownrenormc}= \frac{m_c}{\|m_c\|}$ when
the Gibbs measure $m_c$ on $\Ga\backslash \G \XX$ is finite, and we
use the dynamical system $(\Ga\backslash \G \XX,\,\overline{m_c},
\flow{1})$ in the definition of the correlation coefficients.

Given a finite subset $E$ of $\Ga\bs V\XX$, we denote by $\tau_E:
\Ga\bs \G\XX\ra\NN\cup\{+\infty\}$ the first return\footnote{Actually,
  a more precise terminology is ``first positive passage time'', but
  we use the shorter one. If $\ell\in \Ga\bs \G\XX$ is such that
$\ell(0)\in E$, then ``return'' is appropriate.} time to $E$ of the
discrete time geodesic flow:
$$
\tau_E(\ell)=\inf\{n\in\NN-\{0\}\;:\;\flow{n}\ell(0)\in E\}\;,
$$
with the usual convention that $\inf\emptyset = +\infty$.

\btheo \label{theo:critexpdecaysimpl} Let $\XX,\Ga,\wt c$ be as above,
with $\delta_c$ finite. Assume that the Gibbs measure $m_c$ is finite
and mixing for the discrete time geodesic flow on $\Ga\backslash \G
\XX$. Assume moreover that there exist a finite subset $E$ of
$\Ga\bs V\XX$ and $C',\kappa' >0$ such that for all $n\in\NN$, we have
\begin{equation}\label{eq:expodecaycusp}
m_c\big(\{\ell\in \Ga\bs \G\XX\;:\;\ell(0)\in E \;{\rm and}\;
\tau_E(\ell)\geq n\}\big)\leq C'\;e^{-\kappa' n}\;.
\end{equation}
Then the discrete time geodesic flow on $\Ga\backslash \G \XX$ has
exponential decay of $\alpha$-H\"older correlations for the system of
conductances $c$.  
\etheo

A similar statement holds for the square of the discrete time geodesic
flow on $\Ga\bs \Geven \XX$ when $m_c$ is finite, $\C\Lambda\Ga$ is a
uniform simplicial tree with degrees at least $3$ and $L_\Ga=2\ZZ$.

Note that the crucial Hypothesis \eqref{eq:expodecaycusp} of Theorem
\ref{theo:critexpdecaysimpl} is in particular satisfied if $\Ga\bs\XX$
is finite, by taking $E=\Ga\bs V\XX$.  In this case, the result is quite
well-known: when $\Ga$ is torsion free, it follows from
Bowen's result \cite[1.26]{Bowen75} that a mixing subshift of finite
type is exponentially mixing.

\medskip
\dem Let $\XX'=\C\Lambda\Ga$. Using the coding introduced in Section
\ref{subsec:codagesimplicial}, we first reduce this statement to a
statement in symbolic dynamics.

\subsection*{Step 1 : Reduction to two-sided symbolic dynamics.} 
Let $(\Sigma,\sigma)$ be the (two-sided) topological Markov shift with
alphabet $\A$ and transition matrix $A$ constructed in Section
\ref{subsec:codagesimplicial}, which is conjugated to $(\Ga\bs\G \XX',
\flow{1})$ by the homeomorphism $\Theta:\Ga\bs\G \XX' \ra\Sigma$ (see
Theorem \ref{theo:codingdisgeodflo}). Let $\PP=\Theta_*\;
\frac{m_c}{\|m_c\|}$, which is a mixing $\sigma$-invariant probability
measure on $\Sigma$ with full support, since $\Ga\bs\G\XX'$ is the
support of $m_c$. Let
$$
\E=\{(e^-,h,e^+)\in\A : t(e^-)=o(e^+)\in E\}\,.
$$ 
The set $\E$ is finite since the degrees and the vertex stabilisers
of $\XX$ are finite. For all $x\in\Sigma$ and $k\in\ZZ$, we denote by
$x_k$ the $k$-th component of $x=(x_n)_{n\in\ZZ}$. Let
$$
\tau_\E(x)=\inf\{n\in\NN-\{0\}\;:\;x_n\in \E\}\;
$$ 
be the first return\footnote{See the previous footnote.} time to $\E$
of $x$ under iteration of the shift $\sigma$.

Let $\gls{naturalextension}:\Sigma\ra\A^\NN$ be the {\em natural
  extension}\index{natural extension}\footnote{The authors are not
  responsible for this questionable terminology, rather standard in
  symbolic dynamics.} $(x_n)_{n\in\ZZ} \mapsto (x_n)_{n\in\NN}$.
Theorem \ref{theo:critexpdecaysimpl} will follow from the following
two-sided symbolic dynamics result.\footnote{Assumption (1) of Theorem
  \ref{theo:critexpdecaysimpldynsymb} is far from being optimal, but
  will be sufficient for our purpose.}

\btheo \label{theo:critexpdecaysimpldynsymb} Let $(\Sigma,\sigma)$ be
a locally compact transitive two-sided topological Markov shift with
alphabet $\A$ and transition matrix $A$, and let $\PP$ be a mixing
$\sigma$-invariant probability measure with full support on
$\Sigma$. Assume that 
\begin{enumerate}
\item[(1)] for every $A$-admissible finite sequence $w=
  (w_0,\dots,w_n)$ in $\A$, the Jacobian of the map $f_w$ from
  $\{(x_k)_{k\in\NN}\in\pi_+(\Sigma)\;:\;x_0=w_n\}$ to
  $\{(y_k)_{k\in\NN} \in \pi_+(\Sigma)\;:\;y_0=w_0,\dots, y_n=w_n\}$
  defined by $(x_0,x_1,x_2,\dots)\mapsto (w_0,\dots,w_n,
  x_1,x_2,\dots)$, with respect to the restrictions of the pushforward
  measure $(\pi_+)_*\PP$, is constant;
\item[(2)] there exist a finite subset $\E$ of $\A$ and $C',\kappa'
  >0$ such that for all $n\in\NN$, we have
\begin{equation}\label{eq:expodecaycuspdynsymb}
\PP\big(\{x\in \Sigma\;:\;x_0\in \E \;{\rm and}\;
\tau_\E(x)\geq n\}\big)\leq C'\;e^{-\kappa' n}\;.
\end{equation}
\end{enumerate}
Then $(\Sigma,\PP,\sigma)$ has exponential decay of $\alpha$-H\"older
correlations.  
\etheo

\noindent
{\bf Proof that Theorem \ref{theo:critexpdecaysimpldynsymb} implies
  Theorem \ref{theo:critexpdecaysimpl}.} Since $\wt m_c$ is supported
on $\G\XX'$, up to replacing $\XX$ by $\XX'$, we may assume that
$\partial_\infty X=\Lambda\Ga$.

\medskip
By the construction of $\Theta$ just before the statement of Theorem
\ref{theo:codingdisgeodflo}, for every $\ell=\Ga\wt\ell \in \Ga\bs\G
\XX$, we have $ (\Theta\ell)_0=(e^-_0(\wt \ell\,),h_0(\wt \ell\,),
e^+_0(\wt \ell\,))$ with $e^+_0(\wt \ell\,)=p(\wt\ell([0,1]))$ where
$p:\XX\ra\Ga\bs\XX$ is the canonical projection, so that
$o\big(e^+_0(\wt \ell\,)) \big)=\ell(0)$. Since $\Theta$ conjugates
$\flow{1}$ to $\sigma$, for every $n\in\NN$, we have
$$
(\Theta\ell)_n=(\sigma^n(\Theta\ell))_0=(\Theta(\flow{n}\ell))_0\;.
$$ 
Thus $(\Theta\ell)_n\in\E$ if and only if $\flow{n}\ell(0)\in E$, and
$$
\tau_\E(\Theta\ell)=\tau_E(\ell)\;. 
$$

Therefore Theorem \ref{theo:critexpdecaysimpl} will follow from
Theorem \ref{theo:critexpdecaysimpldynsymb} by conjugation since
$\Theta$ is bilipschitz, once we have proved that Hypothesis (1) of
Theorem \ref{theo:critexpdecaysimpldynsymb} is satisfied for the
two-sided topological Markov shift $(\Sigma,\sigma)$ conjugated by
$\Theta$ to $(\Ga\bs\G \XX, \flow{1})$, which is the main point in
this proof.

\medskip
We hence fix an $A$-admissible finite sequence $w=(w_0,\dots,w_n)$ in
$\A$. We consider the (one-sided) {\it cylinders}\index{cylinder}
$$
[w_n]=\{(x_k)_{k\in\NN}\in\pi_+(\Sigma)\;:\;x_0=w_n\}\;,
$$ 
$$
[w]=\{(y_k)_{k\in\NN} \in \pi_+(\Sigma)\;:\;y_0=w_0,\dots, y_n=w_n\}
$$ 
and the map $f_w:[w_n]\ra [w]$ with $(x_0,x_1,x_2,\dots)\mapsto (w_0,\dots,
w_n, x_1,x_2,\dots)$ that appear in
Hypothesis (1). We denote by $\wt w$ and $\wt{w_n}$ the discrete
generalised geodesic lines in $\XX$ associated with $w$ and $w_n$ (see
the proof of Theorem \ref{theo:codingdisgeodflo} just after Equation
\eqref{eq:canonicliftw}).  Since $w$ ends with $w_n$, by the
construction of $\Theta$, there exists $\ga\in\Ga$ sending the two
consecutive edges of $\wt w_n$ to the last two consecutive edges of
$w$.  We denote by $x=\wt w(0)$ and $y=\wt{w_n}(0)$ the footpoints of
$\wt w$ and $\wt{w_n}$ respectively.

\begin{center}
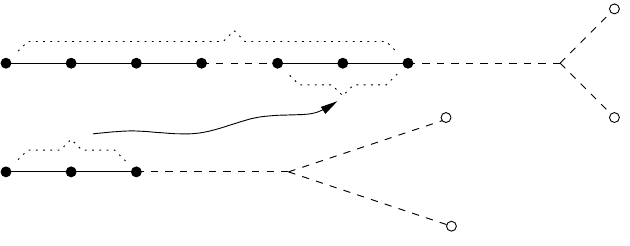
\end{center}

For every discrete generalised geodesic line $\omega\in\gengeod\XX$
which is isometric exactly on an interval $I$ containing $0$ in its
interior (as for $\omega=\wt w,\wt{w_n}$), let
$$
\G_\omega\XX=\{\ell\in\G\XX\;:\; \ell_{\mid I}=\omega_{\mid I}\}
$$ 
be the space of extensions of $\omega_{\mid I}$ to geodesic lines.
With $\partial^\pm_{\omega}\XX=\{\ell_\pm\;:\;\ell\in \G_\omega\XX\}$
its set of points at $\pm\infty$, we have a homeomorphism $\G_\omega
\XX\ra (\partial^-_{\omega}\XX\times \partial^+_{\omega}\XX)$ defined
by $\ell\mapsto (\ell_-,\ell_+)$, using Hopf's parametrisation with
respect to the point $\omega(0)$, since all the geodesic lines in
$\G_\omega\XX$ are at the point $\omega(0)$ at time $t=0$. Using as
basepoint $x_0=\omega(0)$ in the definition of the Gibbs measure (see
Equation \eqref{eq:defigibbsdis}), this homeomorphism sends the
restriction to $\G_\omega \XX$ of the Gibbs measure $d\,\wt m_c(\ell)$
to the product measure $d\mu^-_{\omega(0)}(\ell_-)\;
d\mu^+_{\omega(0)} (\ell_+)$.  Hence the pushforward of $\wt m_c
{}_{\mid \G_\omega\XX}$ by the positive endpoint map $e_+:\ell\mapsto
\ell_+$ is $\mu^-_{\omega(0)}(\partial^-_{\omega}\XX)
\;d\mu^+_{\omega(0)} (\ell_+)$, and note that $\mu^-_{\omega(0)}
(\partial^-_{\omega}\XX)$ is a positive constant.

Let $\wt p:\G\XX\ra\Ga\bs\G\XX$ be the canonical projection. Since
$\pi_+:\Sigma\ra\Sigma_+$ is the map which forgets about the past,
there exist measurable maps $u_w:\partial^+_{\wt w}\XX\ra [w]$ and
$u_{w_n}:\partial^+_{\wt{w_n}}\XX\ra[w_n]$ such that the following
diagrams commute:
$$
\xymatrix{
\G_{\wt{w}} \XX\;\; \ar[rr]^{\pi_+\,\circ\,\Theta\,\circ\, \wt p}\ar[dr]_{e_+}
& &\;\;[w]\\ & \partial^+_{\wt{w}}\XX \ar[ur]_{u_w} &
         }\;\;\;\;\;\;{\rm and}\;\;\;\;\;\; 
\xymatrix{
\G_{\wt{w_n}} \XX\;\; \ar[rr]^{\pi_+\,\circ\,\Theta\,\circ\, \wt p}\ar[dr]_{e_+}
& &\;\;[w_n]\\ & \partial^+_{\wt{w}}\XX \ar[ur]_{u_{w_n}} &
         }\;.
$$
Furthermore, the map $u_w$ (respectively $u_{w_n}$) is surjective,
and has constant finite order fibers given by the orbits of the
finite stabiliser $\Ga_{\wt{w}}$ (respectively $\Ga_{\wt{w_n}}$).
Since $\PP=\Theta_*\frac{m_c}{\|m_c\|}$, the pushforward by the map
$u_w$ (respectively $u_{w_n}$) of the measure $\mu^+_{x}$
(respectively $\mu^+_{y}$) is a constant time the restriction of
$(\pi_+)_*\PP$ to $[w]$ (respectively $[w_n]$). Finally, by the
construction of the (inverse of the) coding in the proof of Theorem
\ref{theo:codingdisgeodflo}, the following diagram is commutative:
$$
\xymatrix{
\partial^+_{\wt{w_n}}\XX\;\; \ar[r]^{u_{w_n}}\ar[d]_{\ga}& 
\;\;[w_n]\;\;\ar[d]^{f_w}\\ 
\partial^+_{\wt{w}}\XX\;\;\ar[r]^{u_w}&\;\;[w]\;. \;\;
         }
$$

Recall that the pushforwards of measures $\mu,\nu$, which are
absolutely continuous one with respect to the other, by a measurable
map $f$ are again absolutely continuous one with respect to the other,
and satisfy (almost everywhere)
$$
\frac{d\,f_*\mu}{d\,f_*\nu}\;\circ\, f \;=\;\;\frac{d\,\mu}{d\,\nu}\;.
$$ 
Hence in order to prove that Hypothesis (1) in the statement of
Theorem \ref{theo:critexpdecaysimpldynsymb} is satisfied, we only have
to prove that the map $\ga: \partial^+_{\wt{w_n}}\XX\ra
\partial^+_{\wt{w}}\XX$ has a constant Jacobian for the measures
$\mu^+_{y}$ on $\partial^+_{\wt{w_n}}\XX$ and $\mu^+_{x}$ on
$\partial^+_{\wt{w}}\XX$ respectively.

For all $\xi,\eta \in \partial^+_{\wt{w}}\XX$, by the properties of
the Patterson densities (see Equations \eqref{eq:equivarPatdens} and
\eqref{eq:quasinivarPatdens}), since $\ga y$ belongs to the geodesic
ray from $x$ to $\xi$ and $\eta$ (see the above picture), and by
Equation \eqref{eq:changemoinsplus}, we have
$$
\frac{\frac{d\,\ga_*\mu^+_{y}}{d\,\mu^+_{x}}(\xi)}
{\frac{d\,\ga_*\mu^+_{y}}{d\,\mu^+_{x}}(\eta)}=
\frac{\frac{d\,\mu^+_{\ga y}}{d\,\mu^+_{x}}(\xi)}
{\frac{d\,\mu^+_{\ga y}}{d\,\mu^+_{x}}(\eta)}=
\frac{e^{-C^+_\xi(\ga y,\,x)}}{e^{-C^+_\eta(\ga y,\,x)}}=
\frac{e^{-\int_x^{\ga y}(\wt F_c^+-\delta_c)}}{e^{-\int_x^{\ga y}(\wt F_c^+-\delta_c)}}=1\;.
$$

\medskip
This proves that Hypothesis (1) in Theorem
\ref{theo:critexpdecaysimpldynsymb} is satisfied, and concludes the
proof of Theorem \ref{theo:critexpdecaysimpl}.  
\cqfd

\medskip
 We now indicate how to pass from a one-sided version of Theorem
 \ref{theo:critexpdecaysimpldynsymb} to the two-sided one, as was
 communicated to us by J.~Buzzi.

\subsection*{Step 2 : Reduction to one-sided symbolic dynamics.}

Let $(\Sigma_+,\sigma_+)$ be the {\it one-sided topological Markov
  shift}\index{Markov!shift!one-sided}%
\index{topological Markov shift!one-sided} with alphabet $\A$ and
transition matrix $A$, that is, $\Sigma_+$ is the closed subset of 
the topological product space $\A^\NN$ defined by
$$
\Sigma_+=\big\{x=(x_n)_{n\in\NN}\in \A^\NN\;:\; \forall \;n\in\NN,\;\;\;
A_{x_n,x_{n+1}}=1\}\;,
$$
and $\gls{onesidedshift}: \Sigma_+\ra \Sigma_+$ is the (one-sided) {\it
  shift}\index{shift}\footnote{Although it is standard to denote the
  one-sided shift by $\sigma$ in the same way as the two-sided shift,
  we use $\sigma_+$ for readability.} defined by
$$
(\sigma_+(x))_n=x_{n+1}
$$ 
for all $x\in \Sigma_+$ and $n\in\NN$. We endow $\Sigma_+$ with the
distance
$$
d(x,x')= e^{-\max\big\{n\in\NN\;:\;\;
\forall\, i\,\in\,\{0,\dots,n\},\;\;x_i\;=\;x'_i\big\}}\;.
$$ 
Note that the distances on $\Sigma$ and $\Sigma_+$ are bounded by $1$.

The natural extension $\pi_+:(x_n)_{n\in\ZZ} \mapsto (x_n)_{n\in\NN}$
maps $\Sigma$ to $\Sigma_+$. It satisfies $\pi_+\circ \sigma =
\sigma_+\circ\pi_+$ and is $1$-Lipschitz. Note that $\Sigma$ is
transitive (respectively locally compact) if and only if $\Sigma_+$ is
transitive (respectively locally compact).

In the one-sided case, we always assume that the cylinders start at
time $t=0$: Given an admissible sequence $w=(w_0,w_1,\dots, w_{n-1})$,
we will say that the cylinder
$$
[w]=[w_0,\dots, w_{n-1}]=\{(x_n)_{n\in\NN}\in\Sigma_+\;:\; 
\forall\, i\,\in\,\{0,\dots,n-1\},\;\;x_i\;=\;w_i\}
$$ 
defined by $w$ has length $\gls{cylinderlength}=n$.

\medskip
We first explain how to relate the decay of correlations for the
two-sided and one-sided systems. This is well-known since the works of
Sinai \cite[\S 3]{Sinai72} and Bowen \cite[Lem.~1.6]{Bowen75}, see for
instance \cite[\S 4]{Young98}, and the following proof has been
communicated to us by J.~Buzzi. We fix $\alpha\in\;]0,1]$. For all
metric spaces $Z$ and bounded $\alpha$-H\"older-continuous functions
$f:Z\ra\RR$, recall\footnote{See Section \ref{subsec:holdercont} for the
  definition of the H\"older norm $\|\cdot\|_\alpha$.} that
$$
\|f\|'_\alpha=\sup_{\substack{x,\,y\,\in \,Z\\ 0<d(x,\,y)\leq 1}}
\frac{|f(x)-f(y)|}{d(x,y)^\alpha}\;\;\;{\rm and}\;\;\;
\|f\|_\alpha =\|f\|_\infty+ \|f\|'_\alpha\;.
$$

\medskip
 For every $a\in\A$, let us fix $z^a\in\Sigma$ such that
$(z^a)_0=a$.

\blemm\label{lem:phinholder} Let $\phi:\Sigma\ra\RR$ be a bounded
$\alpha$-H\"older-continuous map and $N\in\NN$.  For all $x\in
\Sigma_+$, let $y^{(N)}(x)=(y_i)_{i\in\NN}\in\Sigma$, where $y_i=
x_{i+N}$ if $i\geq-N$ and $y_i= (z^{x_0})_{i+N}$ otherwise.  Define
$\phi^{(N)}:\Sigma_+\ra \RR$ by $\phi^{(N)}(x)= \phi(y^{(N)}(x))$.
Then $\phi^{(N)}$ is bounded and $\alpha$-H\"older-continuous on
$\Sigma_+$, with
$$
|\phi\circ\sigma^{N}-\phi^{(N)}\circ\pi_+| \leq 
\|\phi\|'_\alpha\; e^{-\alpha\, N}\;.
$$
Moreover, 
$$
\|\phi^{(N)}\|'_\alpha\leq e^{\alpha \,N}\,\|\phi\|'_\alpha \; \text{ and }\;
\|\phi^{(N)}\|_\infty \leq \|\phi\|_\infty\;.
$$
\elemm

\dem For every $x=(x_n)_{n\in\ZZ}\in\Sigma$, if $y=y^{(N)}(\pi_+(x))$,
we have $(\sigma^{N}(x)\big)_n= x_{n+N}=y_n$ if $|n|\leq N$. Hence,
\begin{align*}
|\phi\circ\sigma^{N}(x)-\phi^{(N)}(\pi_+(x))|&= 
|\phi(\sigma^{N}(x))-\phi(y)|
\\ &\leq \|\phi\|'_\alpha \;d(\sigma^{N}(x),y)^\alpha \leq 
\|\phi\|'_\alpha\; e^{-\alpha\, N}\;.
\end{align*} 
Moreover, for all $x=(x_n)_{n\in\NN}$, $x'=(x'_n)_{n\in\NN}$ in
$\Sigma_+$, if $y=y^{(N)}(x)$ and $y'=y^{(N)}(x')$, then $d(y,y')=
e^{N}\,d(x,x')$ if $d(x,x')<e^{-N}$ and otherwise $d(y,y')\le 1\le
e^Nd(x,x')$, so that
$$
|\phi^{(N)}(x)-\phi^{(N)}(x')| = |\phi(y)-\phi(y')| \leq 
\|\phi\|'_\alpha\; d(y,y')^\alpha\leq 
\|\phi\|'_\alpha\;  e^{\alpha \,N}\,d(x,x')^\alpha\;.\;\;\;\Box
$$

\bprop \label{prop:decaycorrpassonesid} Let $\mu$ be a
$\sigma$-invariant probability measure on $\Sigma$. Assume that the
dynamical system $(\Sigma_+,\sigma_+,(\pi_+)_*\mu)$ has exponential
decay of $\alpha$-H\"older correlations. Then $(\Sigma,\sigma,\mu)$ has
exponential decay of $\alpha$-H\"older correlations.  
\eprop

\dem 
Let $C,\kappa>0$ be such that for all bounded
$\alpha$-H\"older-continuous maps $\phi',\psi':\Sigma_+\ra\RR$ and
$n\in\NN$, we have
$$
|\operatorname{cov}_{(\pi_+)_*\mu,\,n}(\phi',\psi')| \leq 
C\;\|\phi'\|_\alpha\;\|\psi'\|_\alpha\; e^{-\kappa n}\;.
$$ 
Let $\phi,\psi:\Sigma\ra\RR$ be bounded $\alpha$-H\"older-continuous
maps and $n\in\NN$. Denoting by $\pm \,t$ any value in $[-t,t]$ for
every $t\geq 0$, we have, by the first part of Lemma
\ref{lem:phinholder} and for every $N\in \NN$,\footnote{to be chosen
  appropriately below}
\begin{align*}
\int_\Sigma \phi\circ \sigma^n\;\;\psi\; d\mu & = 
\int_\Sigma \phi\circ\sigma^{n+N}\;\;\psi\circ \sigma^{N}\; d\mu \\ & = 
\int_\Sigma (\phi^{(N)}\circ\pi_+ \pm \|\phi\|'_\alpha\;e^{-\alpha\, N})\circ
\sigma^{n}\;\; (\psi^{(N)}\circ\pi_+ \pm \|\psi\|'_\alpha\;e^{-\alpha\, N})\; 
d\mu\\  & = 
\int_{\Sigma_+} \phi^{(N)}\circ\sigma_+^n \;\;\psi^{(N)}\; d(\pi_+)_*\mu
     \;\pm\;\|\phi\|_\alpha\;\|\psi\|_\alpha\;e^{-\alpha \,N}\;.
\end{align*}
A similar estimate holds for the second term in the definition of the
correlation coefficients. Hence, by the second part of Lemma
\ref{lem:phinholder},
\begin{align*}
    &|\operatorname{cov}_{\mu,\,n}(\phi,\psi)| \\
     & \leq |\operatorname{cov}_{(\pi_+)_*\mu,\,n}(\phi^{(N)},\psi^{(N)})| + 
       2\;\|\phi\|_\alpha\;\|\psi\|_\alpha\;e^{-\alpha\, N} \\
     &\leq C \; (\|\phi\|_\infty + \|\phi\|'_\alpha\; e^{\alpha\, N}) \,
      (\|\psi\|_\infty+\|\psi\|'_\alpha\; e^{\alpha\, N})\; e^{-\kappa\, n}  + 
           2\;\|\phi\|_\alpha\;\|\psi\|_\alpha\;e^{-\alpha\, N} \\
     &\leq \|\phi\|_\alpha\;\|\psi\|_\alpha(C\;e^{2\alpha\, N-\kappa\, n}  
      +2\;e^{-\alpha\, N})\;.
\end{align*}
Taking $N=\lfloor \frac{\kappa \,n}{4\,\alpha}\rfloor$,
$C'=C+2\,e^{\alpha}$ and $\kappa'=\frac{\kappa}{4}$, we have
$$
|\operatorname{cov}_{\mu,\,n}(\phi,\psi)| \leq 
C'\;\|\phi\|_\alpha\;\|\psi\|_\alpha\;e^{-\kappa' n}\;,
$$
and the result follows. 
\cqfd

\medskip In order to conclude Step 2, we now state the one-sided
version of Theorem \ref{theo:critexpdecaysimpldynsymb} and prove how
it implies Theorem
\ref{theo:critexpdecaysimpldynsymb}.\footnote{Assumption (1) of
  Theorem \ref{theo:critexpdecaysimpldynsymbonesid} is far from being
  optimal, but will be sufficient for our purpose.} For every finite
subset $\E$ of $\A$, let
$$
\tau_\xi(x)=\inf\{n\in\NN-\{0\}\;:\;x_n\in\E\}
$$ 
be the first return time\footnote{or rather the first positive passage
  time} of $x\in \Sigma_+$ under iteration of the one-sided shift.

\btheo \label{theo:critexpdecaysimpldynsymbonesid} Let
$(\Sigma_+,\sigma_+)$ be a locally compact transitive one-sided
topological Markov shift with alphabet $\A$ and transition matrix $A$,
and let $\PP_+$ be a mixing $\sigma_+$-invariant probability measure
with full support on $\Sigma_+$. Assume that 
\begin{enumerate}
\item[(1)] for every $A$-admissible finite sequence
  $w=(w_0,\dots,w_n)$ in $\A$, the Jacobian of the map from $[w_n]$ to
  $[w]$ defined by $(w_n,x_1,x_2,\dots)\mapsto (w_0,\dots,w_n,
  x_1,x_2,\dots)$ with respect to the restrictions of the measure
  $\PP_+$ is constant;
\item[(2)] 
there exist a finite subset $\E$ of $\A$ and $C',\kappa'
  >0$ such that for all $n\in\NN$, we have
\begin{equation}\label{eq:expodecaycuspdynsymbonesid}
\PP_+\big(\{x\in \Sigma_+\;:\;x_0\in \E \;{\rm and}\;
\tau_\E(x)\geq n\}\big)\leq C'\;e^{-\kappa' n}\;.
\end{equation}
\end{enumerate}
Then $(\Sigma_+,\PP_+,\sigma_+)$ has exponential decay of
$\alpha$-H\"older correlations.  
\etheo

\noindent
{\bf Proof that Theorem \ref{theo:critexpdecaysimpldynsymbonesid}
  implies Theorem \ref{theo:critexpdecaysimpldynsymb}.}  
Let $(\Sigma,\sigma,\PP,\E)$ be as in the statement of Theorem
\ref{theo:critexpdecaysimpldynsymb}.  Let $\PP_+= (\pi_+)_*\PP$, which
is a mixing $\sigma_+$-invariant probability measure on
$\Sigma_+$ with full support. Note that Hypothesis (1) in Theorem
\ref{theo:critexpdecaysimpldynsymbonesid} follows from Hypothesis (1)
of Theorem \ref{theo:critexpdecaysimpldynsymb}. Similarly, Equation
\eqref{eq:expodecaycuspdynsymbonesid} follows from Equation
\eqref{eq:expodecaycuspdynsymb}.  Hence Theorem
\ref{theo:critexpdecaysimpldynsymb} follows from Theorem
\ref{theo:critexpdecaysimpldynsymbonesid} and Proposition
\ref{prop:decaycorrpassonesid}. \cqfd

\medskip
Let us now consider Theorem
\ref{theo:critexpdecaysimpldynsymbonesid}. The scheme of its
proof, using inducing and Young tower arguments, was communicated to
us by O.~Sarig.

\subsection*{Step 3 : Proof of Theorem
  \ref{theo:critexpdecaysimpldynsymbonesid}.}

In this final Step, using inducing of the dynamical system
$(\Sigma_+,\sigma_+)$ on the subspace $\{x\in\Sigma_+\;:\; x_0\in\E\}
=\bigcup_{a\in\E} [a]$ (a finite union of $1$-cylinders), we present
$(\Sigma_+,\sigma_+)$ as a Young tower to which we will apply the
results of \cite{Young99}.

Note that since $\sigma_+$ is mixing and $\PP_+$ has full support,
there exists a $\sigma_+$-invariant measurable subset of full measure
$\Delta$ of $\Sigma_+$ such that the orbit under $\sigma_+$ of every
element of $\Delta$ passes infinitely many times inside the nonempty
open subset $\bigcup_{a\in\E} [a]$. We again denote by
$\tau_\E:\Delta\ra \NN -\{0\}$ the restriction to $\Delta$ of the
first return time in $\bigcup_{a\in\E} [a]$, so that if
$$
\Delta_0=\{x\in\Delta\;:\;x_0\in\E\}=\bigcup_{a\in\E}\;\Delta\cap [a]\;,
$$ 
then $\tau_\E(x)=\min\{n\in\NN-\{0\}\;:\;\sigma_+^n(x)\in\Delta_0\}$ for
all $x\in\Delta$.  We denote by $F:\Delta\ra\Delta_0$ the first return
map to $\Delta_0$ under iteration of the one-sided shift, that is
$$
F:x\mapsto\sigma_+^{\tau_\E(x)}(x)\;.
$$ 
Let $W$ be the set of admissible sequences $w$ of length $|w|$ at
least $2$ such that if $w=(w_0,\dots, w_n)$ with $n=|w|-1$ then
$$
w_0,w_n\in\E\;\;\;{\rm and}\;\;\;w_1,\dots, w_{n-1}\notin\E\;.
$$
We have the following properties:

$\bullet$~ the sets $\Delta_a=\Delta\cap [a]$ for $a\in\E$ form a
finite measurable partition of $\Delta_0$ and for every $a\in\E$, the
sets $\Delta_w=\Delta\cap [w]$ for $w\in W$ and $w_0=a$ form a
countable measurable partition of $\Delta_a$;

$\bullet$~ for every $w\in W$, the first return time $\tau_\E$ is
constant (equal to $|w|-1$) on $\Delta_w$, and if $w_{|w|-1}=b$,
then the first return map $F$ is a bijection from $\Delta_w$ to
$\Delta_{b}$;

$\bullet$~ for all $w\in W$ and $x,y\in \Delta_w$, since $x,y$ have
the same $|w|$ first components, we have
$$
d(F(x),F(y))=d(\sigma_+^{|w|-1}x,\sigma_+^{|w|-1}y)=e^{|w|-1}\;d(x,y)
\geq e\;d(x,y)\;;
$$

$\bullet$~ for all $w\in W$, $n\in\{0,\dots,|w|-2\}$ and $x,y\in
\Delta_w$, we have
$$
d(\sigma_+^nx,\sigma_+^ny)= e^n\;d(x,y)\leq e^{|w|-2}\;d(x,y)
< d(F(x),F(y))\;;
$$

$\bullet$~ for every $w\in W$, the Jacobian of the first return map
$F:\Delta_w\ra \Delta_{w_{|w|-1}}$ for the restrictions to $\Delta_w$
and $\Delta_{w_{|w|-1}}$ of $\PP_+$ is constant.\footnote{Actually,
  only a much weaker assumption is required, such as a
  H\"older-continuity property of this Jacobian, see \cite{Young99}.}

\medskip
By an easy adaptation of \cite[Theo.~3]{Young99} (see also \cite[\S
2.1]{Melbourne07}) which considers the case when $\E$ is a
singleton, we have the following noneffective\footnote{Actually,
  there is in \cite{Young99} (see also \cite{CyrSar09}) a control on
  the constant in terms of some norms of the test functions, but these
  norms are not the ones we are interested in.} exponential decay of
correlation: there exists $\kappa>0$ such that for every $\phi,\psi\in
\C^\alpha_{\rm b}(\Sigma_+)$, there exists a constant $C_{\phi,\psi}>0$ such
that
$$
|\operatorname{cov}_{\PP^+,\,n}(\phi,\psi)| \leq 
C_{\phi,\psi}\;e^{-\kappa n}
$$ 
By an elegant argument using the Principle of Uniform Boundedness,
it is proved in \cite[Appendix B]{ChaColSch05} that this implies that 
there exist $C,\kappa>0$ such that for every $\phi,\psi\in
\C^\alpha_{\rm b}(\Sigma_+)$,  we have
$$
|\operatorname{cov}_{\PP^+,\,n}(\phi,\psi)| \leq 
 C\;\|\phi\|_\alpha\;\|\psi\|_\alpha\;e^{-\kappa\, n }\;.
$$ 
This concludes the proof of Theorem
\ref{theo:critexpdecaysimpldynsymbonesid}, hence the proof of Theorem
\ref{theo:critexpdecaysimpl}.  
\cqfd \cqfd

\bigskip
The next result gives examples of applications of Theorem
\ref{theo:critexpdecaysimpl} when $\Ga\bs\XX$ is infinite. It
strengthens \cite[Theo.~2.1]{AthGhoPra12} that applies only to
arithmetic lattices and only for the locally constant regularity (see
Section \ref{subsec:locconst}), see also \cite{BekLub11} for an
approach using spectral gaps. It was claimed in \cite{Kwon15}, but was
retracted by the author.

\bcoro\label{coro:expdecaygeomfinisimpl} Let $\XX$ be a locally finite
simplicial tree without terminal vertices. Let $\Ga$ be a
geometrically finite subgroup of $\Aut(\XX)$ such that the smallest
nonempty $\Ga$-invariant subtree of $\XX$ is uniform without vertices
of degree $2$. Let $\alpha\in\;]0,1]$.

\smallskip
(1) If $L_\Ga=\ZZ$, then the discrete time geodesic flow on $\Ga\bs \G
\XX$ has exponential decay of $\alpha$-H\"older correlations for the
zero system of conductances.

\smallskip
(2) If $L_\Ga=2\ZZ$, then the square of the discrete time geodesic
flow on $\Ga\bs \Geven \XX$ has exponential decay of $\alpha$-H\"older
correlations for the zero system of conductances, that is, there exist
$C,\kappa >0$ such that for all $\phi,\psi\in \C_{\rm b}^\alpha
(\Ga\backslash \Geven \XX)$ and $n\in\ZZ$, we have
\begin{align*}
&\Big|\int_{\Ga\backslash \Geven \XX}\phi\circ\flow{-2n}\;\psi\; d\, m_{\rm BM}
- \frac{1}{{\text{\scriptsize $m_{\rm BM}(\Ga\backslash \Geven \XX)$}}}
\int_{\Ga\backslash \Geven \XX} \phi\; d\, m_{\rm BM}
\int_{\Ga\backslash \Geven \XX} \psi\;d\, m_{\rm BM}\;\Big| \\ &\le\; 
C\;e^{-\kappa|n|}\;\|\phi\|_\alpha\;\|\psi\|_\alpha\,.
\end{align*}
\ecoro

The main point in order to obtain this corollary is to prove the
exponential decay of volumes of geodesic lines going high in the
cuspidal rays of $\Ga\bs\XX$, stated as Assumption
\eqref{eq:expodecaycusp} in Theorem
\ref{theo:critexpdecaysimpl}. There is a long history of similar
results, starting from the exponential decay of volumes of small cusp
neighbourhoods in noncompact finite volume hyperbolic manifolds (based
on the description of their ends) used by Sullivan to deduce
Diophantine approximation results (see \cite[\S
  9]{Sullivan82}).\footnote{and by probabilists in order to study the
  statistics of cusp excursions (see for instance \cite{EnrFra02})}
These results were extended to the case of locally symmetric
Riemannian manifolds by Kleinbock-Margulis \cite{KleMar99} (based on
the description of their ends using Siegel sets). Note that the
geometrically finite lattice assumption on $\Ga$ is here in order to
obtain similar descriptions of the ends of $\Ga\bs\XX$.

\medskip
\dem Up to replacing $\XX$ by $\C\Lambda\Ga$, we assume that $\XX$ is
a uniform simplicial tree with degrees at least $3$ and that $\Ga$ is
a geometrically finite lattice of $\XX$. We use the zero system of
conductances.

\medskip\noindent(1) By \cite{Paulin04b},\footnote{See also Section
  \ref{subsec:trees}.} the graph $\Ga\bs \XX$ is the union of a finite
graph $\YY$ and finitely many cuspidal rays $R_i$ for $i\in\{1,\dots,
k\}$.  If $(x_{i,\,n})_{n\in\NN}$ is the sequence of vertices in
increasing order along $R_i$ for $i=1,\dots, k$, then the vertex group
$G_{x_{i,\,n}}$ of $x_{i,\,n}$ in the quotient graph of groups $\Ga
\dbs \XX$ satisfies $G_{x_{i,\,n}}\subset G_{x_{i,\,n+1}}$ for every
$n\in\NN$, and the edge group of the edge $e_{i,n}$ with origin
$x_{i,n}$ and endpoint $x_{i,\,n+1}$ is equal to
$G_{x_{i,\,n}}$.\footnote{identifying the edge group of an edge $e$
  with its image by the structural map $G_e\ra G_{o(e)}$} Note
that since the degrees of the vertices of $\XX$ are at least $3$, we
have $[G_{x_{i,\,n+1}}:G_{x_{i,n}}]\geq 2$ and $|G_{x_{i,\,0}}|\geq
1$, so that, for every $n\in\NN$,
\begin{equation}\label{eq:expogrowthcusp}
|G_{x_{i,\,n}}|\geq 2^{n}\;.
\end{equation}

Let $E$ be the (finite) set of vertices $V\YY$ of $\YY$. Note that for
all $n\in\NN-\{0\}$ and $\ell\in \Ga\bs \G\XX$, if $\ell(0)\in E$ and
$\tau_E(\ell)\geq 2n$, then $\ell$ leaves $\YY$ after time $0$
and it travels (geodesically) inside some cuspidal ray for a time at
least $n$, so that there exists $i\in\{1,\dots,k\}$ such that
$\ell(n)=x_{i,\,n}$. Hence for all $n\in\NN$, using 

$\bullet$~ the invariance of $m_{\rm BM}$ under the discrete time
geodesic flow in order to get the third term,

$\bullet$~ Equation \eqref{eq:calcpartielmBMprequel} where $\wt
x_{i,\,n}$ is a fixed lift of $x_{i,\,n}$ in $V\XX$ for the fifth
term, and

$\bullet$~ Equation \eqref{eq:expogrowthcusp} since $|\Ga_{\wt
  x_{i,\,n}}|=|G_{x_{i,\,n}}|$, and the facts that the degrees of the
uniform simplicial tree $\XX$ are uniformly bounded and that the total
mass of the Patterson measures of the lattice $\Ga$ are uniformly
bounded (see Proposition \ref{prop:uniflatmBMfinie}) for the last
term,

\noindent we have
\begin{align*}
& m_{\rm BM}\big(\{\ell\in \Ga\bs \G\XX\;:\;\ell(0)\in E \;{\rm and}\;
\tau_E(\ell)\geq 2n\}\big)\\ \leq\; &
\sum_{i=1}^km_{\rm BM}\big(\{\ell\in\Ga\bs\G\XX\;:\;\ell(n)=x_{i,n}\}\big)
\\ =\; &
\sum_{i=1}^km_{\rm BM}\big(\{\ell\in\Ga\bs\G\XX\;:\;\ell(0)=x_{i,n}\}\big)
=
\sum_{i=1}^k\pi_*m_{\rm BM}(\{x_{i,n}\})
\\ =\; &
\sum_{i=1}^k \;\frac{1}{|\Ga_{\wt x_{i,\,n}}|}\; 
\sum_{e,\,e'\in E\XX\;:\; o(e)=o(e')=\,\wt x_{i,\,n},\;e\neq e'}
\mu_{\wt x_{i,\,n}}(\partial_e\XX)\;\mu_{\wt x_{i,\,n}}(\partial_{e'}\XX)
\\ \leq\; & k\;\frac{1}{2^n}\;
\max_{x\in V\XX }\deg (x)^2\;\max_{x\in V\XX }\|\mu_x\|^2\;.
\end{align*}

The result then follows from Theorem \ref{theo:critexpdecaysimpl}
using the above finite set $E$ which satisfies Assumption
\eqref{eq:expodecaycusp} as we just proved, and using Proposition
\ref{prop:uniflatmBMfinie} and Theorem \ref{theo:uniflatmBMmixing} in
order to check that under the assumption that $L_\Ga=\ZZ$, the
Bowen-Margulis measure $m_{\rm BM}$ of $\Ga$ is finite and mixing
under the discrete time geodesic flow on $\Ga\bs\G\XX$.

\medskip
\noindent (2) The proof of Assertion (2) is similar to the one
of Assertion (1).
 \cqfd

\bigskip
\rem 
The techniques introduced in the above proof in order to check the
main hypothesis of Theorem \ref{theo:critexpdecaysimpl} may be applied
to numerous other examples.  For instance, let $\XX$ be a locally
finite simplicial tree without terminal vertices. Let $\Ga$ be a
nonelementary discrete subgroup of $\Aut(\XX)$ such that the smallest
nonempty $\Ga$-invariant subtree of $\XX$ is uniform without vertices
of degree $2$, and such that $L_\Ga=\ZZ$. Let $\alpha\in\;]0,1]$.
Assume that $\Ga\bs \XX$ is the union of a finite graph $A$ and
finitely many trees $\TT_1,\dots, \TT_k$ meeting $A$ in one and
exactly one vertex $*_1,\dots, *_n$ such that for every edge $e$ in
$\TT_i$ pointing away from the root $*_i$ of $\TT_i$, the canonical
morphism $G_e\ra G_{o(e)}$ between edge and vertex groups of the
quotient graph of groups $\Ga\dbs\XX=(\Ga\bs\XX,G_*)$ is an
isomorphism. Assume that there exist $C,\kappa>0$ such that for all
$n\in\NN$,
$$
\sum_{i=1,\dots,k,\;x\in V\TT_i\;:\; d(x,*_i)=n}\frac{1}{|G_x|}\leq Ce^{-\kappa \,n}\;.
$$ 
Then the discrete time geodesic flow on $\Ga\bs \G \XX$ has
exponential decay of $\alpha$-H\"older correlations for the zero
system of conductances.

This is in particular the case for every $k,q\in\NN$ such that $k\geq
2$, $q> 2k+1$ and $q-k$ is odd, when the quotient graph of groups
$\Ga\dbs\XX$ has underlying edge-indexed graph\footnote{See definition
  in Section \ref{subsec:trees}.} a loop-edge with both indices equal
to $\frac{q-k+1}{2}$ glued to the root of a regular $k$-ary rooted
tree, with indices $1$ for the edges pointing towards the root and
$q-k+1$ for the edges pointing away from the root (see the picture
below with $k=2$). Note that $\XX$ is then the $(q+1)$-regular tree,
and that the loop edge is here in order to ensure that $L_\Ga=\ZZ$.
For instance, the vertex group of a point at distance $n$ from the
root may be chosen to be $\ZZ/(\frac{q-k+1}{2})\ZZ \times
\big(\ZZ/(q-k+1)\ZZ\big)^n$.

\begin{center}
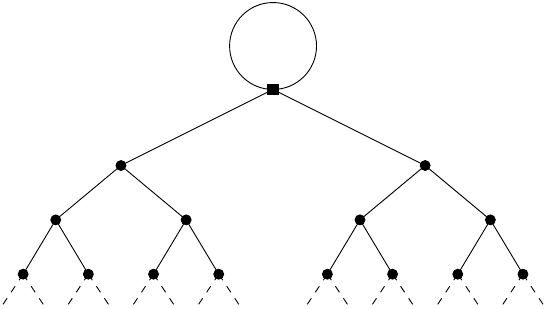
\end{center}

\section{Rate of mixing for metric  trees}
\label{subsec:mixingratemetgraphs}

Let $(\XX,\lambda)$, $X$, $\Ga$, $\wt F$, $\wt F^\pm$, $\delta=
\delta_{\Ga,\,F^\pm}<+\infty$ and $(\mu^\pm_x)_{x\in VX}$ be as in
the beginning of Section \ref{subsec:ergodictrees}. Let $\wt m_F$ and
$m_F$ be the associated Gibbs measures on $\G X$ and $\Ga\bs\G X$. The
aim of this Section is to study the problem of finding conditions on
these data under which the (continuous time) geodesic flow on
$\Ga\bs\G X$ is polynomially mixing for the Gibbs measure $m_F$.

We will actually prove a stronger property, though it applies only to
observables which are smooth enough along the flow. Let us fix
$\alpha\in \;]0,1]$. Let $(Z,\mu,(\phi_t)_{t\in\RR})$ be a
topological space $Z$ endowed with a continuous one-parameter group
$(\phi_t)_{t\in\RR}$ of homeomorphisms preserving a (Borel)
probability measure $\mu$ on $Z$. For all $k\in\NN$, let
$\gls{espacealphaholdebk}$ be the real vector space of maps $f:Z\ra\RR$
such that for all $z\in Z$, the map $t\mapsto f(\phi_tz)$ is
$\C^k$-smooth, and such that the maps $\partial^i_tf:Z\ra\RR$ defined
by $z\mapsto \frac{d^i}{dt^i}\big|_{ t=0} f(\phi_tz)$ for $0\leq i\leq
k$ are bounded and $\alpha$-H\"older-continuous. It is a Banach space
when endowed with the norm
$$
\|f\|_{k,\,\alpha} = \sum_{i=0}^k\; \|\partial^i_tf\|_{\alpha}\;,
$$ 
and it is contained in $\LL^2(Z,\mu)$ by the finiteness of $\mu$. We
denote by $\gls{espacealphaholderkc}$ the vector subspace of elements
of $\C^{k,\,\alpha}_{\rm b} (Z)$ with compact support.

For all $\psi,\psi'\in \LL^2(Z,\mu)$ and $t\in\RR$, let
$$
\gls{correlationcoeffcont}=\int_{Z} \psi\circ\phi_{t}\;\psi'\;d\mu
- \int_{Z}\psi\; d\mu\;\int_{Z}\psi'\;d\mu
$$ 
be the (well-defined) {\it correlation coefficient}%
\index{correlation coefficient} of the observables $\psi,\psi'$ at
time $t$ under the flow $(\phi_t)_{t\in\RR}$ for the measure $\mu$. We
say\footnote{See \cite{Dolgopyat98ETDS}, and more precisely
  \cite[Def.~2.2]{Melbourne07} whose definition is slightly different
  but implies the one given in this paper by the Principle of Uniform
  Boundedness argument of \cite[Appendix B]{ChaColSch05} already used
  in Section \ref{subsec:mixingratesimpgraphs}.} that the (continuous
time) dynamical system $(Z,\mu, (\phi_t)_{t\in\RR})$ has {\it
  superpolynomial decay of $\alpha$-H\"older correlations}%
\index{decay of correlations!superpolynomial} if for every $n\in\NN$
there exist $C=C_n >0$ and $k=k_n\in\NN$ such that for all $\psi,\psi'
\in \C^{k,\,\alpha}_{\rm b} (Z)$ and $t\in\RR$, we have
$$
|\operatorname{cov}_{\mu,\,t}(\psi,\psi')| \le 
C\;(1+|t|)^{-n}\;\|\psi\|_{k,\,\alpha}\;\|\psi'\|_{k,\,\alpha}\;.
$$ 
Following Dolgopyat, we say that the dynamical system $(Z,\mu,
(\phi_t)_{t\in\RR})$ is {\em rapidly mixing}\index{rapid
  mixing} if there exists $\alpha>0$ such that
$(Z,\mu, (\phi_t)_{t\in\RR})$ has superpolynomial decay of
$\alpha$-H\"older correlations

\medskip
We will use the following two assumptions on our data, introduced
respectively in \cite{Dolgopyat98ETDS} and \cite{Melbourne07}. Recall
that the Gibbs measure $m_F$, when finite, is mixing if and only the
length spectrum $L_\Ga$ is dense in $\RR$ (see Theorem
\ref{theo:mixing}). The rapidly mixing property will require stronger
assumptions on $L_\Ga$.

We say that the length spectrum $L_\Ga$ of $\Ga$ is {\it
  $2$-Diophantine}\index{Diophantine@$2$-Diophantine} if there exists
a ratio of two translation lengths of elements of $\Ga$ which is
Diophantine. Recall that a real number $x$ is {\it
  Diophantine}\index{Diophantine} if there exist $\alpha,\beta>0$ such
that
$$
\big|x-\frac{p}{q}\big|\geq \alpha\, q^{-\beta}
$$ 
for all $p,q\in\ZZ$ with $q> 0$.

Let $E$ be a finite subset of vertices of $\Ga\bs\XX$, and let $\wt E$
be the set of vertices of $\XX$ mapping to $E$.  We denote by $T_E$
the set of triples $(\lambda(\ga), d(\gamma), q(\gamma,p))$ where
$\ga\in\Ga$ has translation length $\lambda(\ga)>0$, has $d(\gamma)$
vertices on its translation axis $\Ax(\ga)$ modulo $\ga^\ZZ$ and if
the first return time of a vertex $p$ in $\wt E\cap\Ax(\ga)$ under the
discrete time geodesic flow along the translation axis has period
$q(\gamma,p)$. We say that the length spectrum $L_\Ga$ of $\Ga$ is
{\it $4$-Diophantine}\index{Diophantine@$4$-Diophantine} with respect
to $E$ if for all sequences $(b_k)_{k\in\NN}$ in $[1+\infty[$
converging to $+\infty$ and $(\omega_k)_{k\in\NN}$,
$(\varphi_k)_{k\in\NN}$ in $[0,2\pi[\,$, there exists $N\in\NN$ such
that for all $a\geq N$ and $C,\beta\geq 1$, there exist $k\geq 1$ and
$(\tau,d,q)\in T_E$ such that
$$
d\big((b_k\tau+\omega_kd)\lfloor \beta\ln b_k\rfloor+
q\varphi_k,2\pi\ZZ)\geq C\,q\,b_k^{-a}\;.
$$

We define the first return time after time $\epsilon$ on a finite
subset $E$ of vertices of $\Ga\bs\XX$ as the map $\tau^{> \epsilon} _E
:\Ga\bs \G X\ra \;[0,+\infty]$ defined by $\tau^{> \epsilon}_E
(\ell)=\inf\{t>\epsilon\;:\;\ell(t)\in E\}$.

\btheo\label{theo:mixingtreelattice} Assume that the Gibbs measure
$m_F$ is finite and mixing for the (continuous time) geodesic flow,
and that the lengths of the edges of $(\XX,\lambda)$ have a finite
upper bound.\footnote{They have a positive lower bound by definition,
  see Section \ref{subsec:trees}.} Furthermore assume that
\begin{enumerate} 
\item[(a)] either $\Ga\bs X$ is compact and the length spectrum of
  $\Ga$ is $2$-Diophantine,
\item[(b)] or there exists a finite subset $E$ of vertices of
  $\Ga\bs\XX$ satisfying the following properties:
\begin{enumerate} 
\item[(1)] there exist $C,\kappa>0$ and $\epsilon\in\;]0,\min\lambda[$
        such that for all $t\geq 0$,
$$
m_F(\{\ell\in\Ga\bs \G X\;:\;d(\ell(0),E)\leq \epsilon\;\;\;{\rm and}\;\;\;
\tau^{> \epsilon}_E(\ell)\geq t\})\leq C\;e^{-\kappa\, t}\;,
$$
\item[(2)] the length spectrum of $\Ga$ is $4$-Diophantine with
  respect to $E$.
\end{enumerate} 
\end{enumerate} 
Then the (continuous time) geodesic flow on $\Ga\bs\G X$ has
superpolynomial decay of $\alpha$-H\"older correlations for the
normalised Gibbs measure $\frac{m_F}{\|m_F\|}$.  
\etheo

Note that the existence of $E$ satisfying the exponentially small tail
Hypothesis (1) is in particular satisfied if $\Ga$ is geometrically
finite with $E$ the set of vertices of a finite subgraph of
$\Ga\bs\XX$ whose complement in $\Ga\bs\XX$ is the underlying graph of
a union of cuspidal rays in $\Ga\dbs\XX$~: see the proof of Corollary
\ref{coro:expdecaygeomfinisimpl} and use the hypothesis on the lengths
of edges.

Note that the exponentially small tail Hypothesis (1) might be
weakened to a superpolynomially small tail hypothesis while keeping
the same conclusion, see \cite{Melbourne09}. Since the former is
easier to check than the latter, we prefer to state Theorem
\ref{theo:mixingtreelattice} as it is.

We will follow a scheme of proof analogous to the one in Section
\ref{subsec:mixingratesimpgraphs} for simplicial trees, by reducing
the study to a problem of suspensions of Young towers, and then apply
results of \cite{Dolgopyat98ETDS} and \cite{Melbourne07} for the rapid
mixing property of suspensions of hyperbolic and nonuniformly
hyperbolic dynamical systems.

\medskip
\dem 
Since the Gibbs measure normalised to be a probability measure depends
only on the cohomology class of the potential (see Equation
\eqref{eq:equalgibbsmeascohomolog}), we may assume by Proposition
\ref{prop:relatpotentialconductance} that $F=F_c$ is the potential on
$\Ga\bs T^1X$ associated with a system of conductances $\wt c:E\XX\ra
\RR$ for $\Ga$. We denote by $\delta_c$ the critical exponent of
$(\Ga,F_c)$, and by $m_c$ the Gibbs measure $m_{F_c}$. Up to replacing
$\XX$ by its minimal nonempty $\Ga$-invariant subtree, we assume that
$X=\C\Lambda\Ga$.

\subsection*{Step 1 : Reduction to a suspension of a two-sided 
symbolic dynamics.}

We refer to the paragraphs before the statement of Theorem
\ref{theo:mainsuspension} and at the beginning of Section
\ref{subsec:codagemetric} for the definitions of

\smallskip
$\bullet$~ the system of conductances $^\sharp c$ for $\Ga$ on the
simplicial tree $\XX$,

\smallskip
$\bullet$~ the (two-sided) topological Markov shift $(\Sigma,
\sigma,\PP)$ on the alphabet $\A$, conjugated to the discrete time
geodesic flow $\big( \Ga\bs \G\XX,\;^\sharp\flow{1},
\frac{m_{^\sharp c}}{\|m_{^\sharp c}\|} \big)$ by the
homeomorphism $\Theta: \Ga\bs\G\XX\ra\Sigma$,

\smallskip
$\bullet$~ the roof function $r: \Sigma \ra\; ]0,+\infty[\,$,

\smallskip
$\bullet$~ the suspension $(\Sigma,\sigma, a\,\PP)_r = (\Sigma_r,
(\sigma^t_r)_{t\in\RR},a\,\PP_r)$ over $(\Sigma,\sigma,a\,\PP)$ with
roof function $r$, where $a=\frac{1}{\|\PP_r\|}$. 

The suspension $(\Sigma,\sigma, a\,\PP)_r$ is
conjugated to the 
geodesic flow $\big(\Ga\bs \G X,
\frac{m_c}{\|m_c\|}, (\flow{t})_{t\in\RR}\big)$ by the bilipschitz
homeomorphism $\Theta_r:\Ga\bs \G X\ra\Sigma_r$ defined at the end of
the proof of Theorem \ref{theo:mainsuspension}. We will always
(uniquely) represent the elements of $\Sigma_r$ as $[x,s]$ with $x\in
\Sigma$ and $0\leq s < r(x)$.

\medskip
Note that since $\Theta_r^{-1}$ conjugates $(\sigma^t_r)_{t\in\RR}$
and $(\flow{t})_{t\in\RR}$, we have for all $f:\Ga\bs \G X\ra \RR$ and
$x\in\Sigma_r$, when defined,
$$
\partial^i_t(f\circ\Theta_r^{-1})(x)=
\frac{d^i}{dt^i}\big|_{ t=0} f\circ\Theta_r^{-1}(\sigma^t_r x)=
\frac{d^i}{dt^i}\big|_{ t=0} f(\flow{t}\Theta_r^{-1}(x))=
(\partial^i_tf)\circ\Theta_r^{-1}(x)\;.
$$ 
Hence if $f:\Ga\bs \G X\ra \RR$ is $C^k$-smooth along the orbits of
$(\flow{t})_{t\in\RR}$, then $f\circ\Theta_r^{-1}$ is $C^k$-smooth
along the orbits of $(\sigma^t_r)_{t\in\RR}$. Furthermore, since
$\Theta_r$ is bilipschitz, the precomposition map by $\Theta_r^{-1}$
is a continuous linear isomorphism from $\C^{k,\,\alpha}_{\rm b}(\Ga\bs \G X)$
to $\C^{k,\,\alpha}_{\rm b}(\Sigma_r)$.

Note that since $\Theta_r$ conjugates $(\flow{t})_{t\in\RR}$ and
$(\sigma^t_r)_{t\in\RR}$, and sends $\frac{m_c}{\|m_c\|}$ to
$\frac{\PP_r}{\|\PP_r\|}$, we have, for all $\psi,\psi'
\in \LL^2(\Ga\bs \G X)$ and $t\in\RR$,
$$
\operatorname{cov}_{\frac{m_c}{\|m_c\|},\,t}(\psi,\psi')=
\operatorname{cov}_{\frac{\PP_r}{\|\PP_r\|},\,t}
(\psi\circ\Theta_r^{-1},\psi'\circ\Theta_r^{-1})\;.
$$ 
Therefore we only have to prove that the suspension $(\Sigma_r,
(\sigma^t_r)_{t\in\RR},\frac{\PP_r}{\|\PP_r\|})$ is rapidly mixing
under one of the assumptions (a) and (b).

\subsection*{Step 2 : Reduction to a suspension of a one-sided 
symbolic dynamics.}

In this Step, we explain the rather standard reduction concerning
mixing rates from suspensions of two-sided topological Markov shifts
to suspensions of one-sided topological Markov shifts. We use the
obvious modifications of the notation and constructions concerning the
suspension of a noninvertible transformation to a semiflow, given for
invertible transformations at the beginning of Section
\ref{subsec:codagemetric}.

We consider the one-sided topological Markov shift $(\Sigma_+,
\sigma_+, \PP_+)$ over the alphabet $\A$ constructed at the beginning
of Step 2 of the proof of Theorem \ref{theo:critexpdecaysimpl}, where
the system of conductances $c$ is now replaced by $^\sharp c$, so that
$\PP=\Theta_*(m_{^\sharp c}/\|m_{^\sharp c}\|)$. Let $\pi_+:\Sigma \to
\Sigma_+$ be the natural extension so that $\PP_+=(\pi_+)_*\PP$ and
$\pi_+\circ\sigma=\sigma_+\circ\pi_+$.

We are going to construct in Step 2, as the suspension of $(\Sigma_+,
\sigma_+, \PP_+)$ with an appropriate roof function $r_+$, a semiflow
$((\Sigma_+)_{r_+}, \big((\sigma_+)^t_{r_+}\big)_{t\geq 0},
(\PP_+)_{r_+})$, and prove that the flow $\big(\Sigma_r,
(\sigma^t_r)_{t\in\RR}, \frac{\PP_r}{\|\PP_r\|}\big)$ is rapidly mixing if
the semiflow $((\Sigma_+)_{r_+}, \big((\sigma_+)^t_{r_+} \big)_{t\geq
  0}, \frac{(\PP_+)_{r_+}}{\|(\PP_+)_{r_+}\|})$ is rapidly mixing.
  
Let $r_+: \Sigma_+ \ra\; ]0,+\infty[$ be the map defined by
\begin{equation}\label{eq:defirrrrrplus}
r_+:x\mapsto \lambda(e^+_0)\;,
\end{equation}
where if $x=(x_n)_{n\in\NN}$, the edge $e_0^+$ is such that $x_0=
(e_0^-,h_0,e_0^+)$.  Note that this map has a positive lower bound,
and a finite upper bound, and that it is locally constant (and even
constant on the $1$-cylinders of $\Sigma_+$). By Equation
\eqref{eq:defirrrrr}, we have
\begin{equation}\label{eq:rplusr}
r_+\circ\pi_+=r\;.
\end{equation}
We denote by $((\Sigma_+)_{r_+}, \big((\sigma_+)^t_{r_+}\big)_{t\geq
  0}, (\PP_+)_{r_+})$ the suspension semiflow over $(\Sigma_+,
\sigma_+, \PP_+)$ with roof function $r_+$. We (uniquely) represent
the points of the suspension space $(\Sigma_+)_{r_+}$ as $[x,s]$ for
$x\in \Sigma_+$ and $0\leq s<r_+(x)$. For all $t\geq 0$, we have
$(\sigma_+)^t_{r_+} ([x,s]) = [\sigma_+^nx,s']$ where $n\in\NN$ and
$s'\in\RR$ are such that $t+s=\sum_{i=0}^{n-1}r_+(\sigma_+^ix) + s'$
and $0\leq s'< r_+(\sigma_+^nx)$.

We define the {\it suspended natural 
extension}\index{natural extension!suspended} as the map
$\pi_+^{\;r}:\Sigma_r\ra (\Sigma_+)_{r_+}$ by
$$
\pi_+^{\;r}:[x,s]\mapsto [\pi_+(x),s]\;,
$$ which is well defined by Equation \eqref{eq:rplusr}.  Note that
$\pi_+^{\;r}$ is $1$-Lipschitz for the Bowen-Walters distances on
$\Sigma_r$ and $(\Sigma_+)_{r_+}$ (see Proposition
\ref{prop:defiBowenWalters}).\footnote{Note that Proposition
  \ref{prop:defiBowenWalters} is stated for suspensions of invertible
  maps, but the roof function is constant on $1$-cylinders, and the
  branches of the inverse of $\sigma_+$ on these $1$-cylinders are
  uniformely Lipschitz, hence the proof of \cite[Appendix]{BarSau00}
  extends.}

For all $\psi:\Sigma_r\ra\RR$ and $T\geq 0$, let us construct a
function $\psi^{(T)}:(\Sigma_+)_{r_+}\ra\RR$ as follows.  For every
$[x,s]\in(\Sigma_+)_{r_+}$, let $N\in\NN$ and $s'\geq 0$ be such that
$(\sigma_+)_{r_+}^T[x,s]=[\sigma_+^Nx,s']$, with 
$$
0\leq s'< r_+(\sigma_+^Nx)\;\;\;{\rm and}\;\;\;
s+T=\sum_{i=0}^{N-1} r_+(\sigma_+^ix)\;+ s'\;.
$$
Let
$$
\psi^{(T)}([x,s])=\psi([y,s'])
$$ 
where $y=(y_n)_{n\in\ZZ}$ is such that $y_i=x_{i+N}$ if $i\geq -N$
and $y_i=(z^{x_0})_{i+N}$ otherwise. Note that $y_0=x_N$, hence $r(y)=
r_+(\sigma_+^N(x))$ and $0\leq s' <r(y)$, so that the above map is
well defined.

Finally, for every $\psi\in\C^{k,\,\alpha}_{\rm b}(\Sigma_r)$ or
$\psi\in\C^{k,\,\alpha}_{\rm b}((\Sigma_+)_{r_+})$, let
$$
\|\psi\|_{k,\,\infty} =\sum_{i=0}^k\;\|\partial^i_t\psi\|_{\infty}
\;\;\;{\rm and}\;\;\;
\|\psi\|'_{k,\,\alpha} =\sum_{i=0}^k\;\|\partial^i_t\psi\|'_{\alpha}\;,
$$
so that
$$
\|\psi\|_{k,\,\alpha}=\|\psi\|_{k,\,\infty}+\|\psi\|'_{k,\,\alpha}\;.
$$

\blemm\label{lem:passageauT}
Let $T\geq 0$ and $\psi\in\C^{k,\,\alpha}_{\rm b}(\Sigma_r)$.
\begin{enumerate} 
\item[(1)] 
  For all $t\geq 0$, we have $(\sigma_+)^t_{r_+}\circ
  \pi_+^{\;r} = \pi_+^{\;r}\circ\sigma_r^t$.
\item[(2)] With $\alpha'=\frac{\alpha}{\sup\lambda}$, there exists a
  constant $C_1\geq 1$ (independent of $k$, $T$ and $\psi$) such that
  $|\psi\circ\sigma_r^T- \psi^{(T)}\circ\pi_+^{\;r}| \leq
  C_1\;\|\psi\|'_\alpha\; e^{-\alpha'\, T}$.
\item[(3)] We have $\psi^{(T)}\in\C^{k,\,\alpha}_{\rm b}
  ((\Sigma_+)_{r_+})$ and $\|\psi^{(T)}\|_{k,\,\infty}\leq \|\psi\|_{k,\,\infty}$.
  With $\alpha''=\frac{\alpha}{\inf\lambda}$, there exists a constant
  $C_2\geq 1$ (independent of $k$, $T$ and $\psi$) such that
\begin{equation}\label{eq:controlckalphanorm}
\|\psi^{(T)}\|'_{k,\,\alpha} \leq C_2\; e^{\alpha''\,T}\;\|\psi\|'_{k,\,\alpha}\;.
\end{equation}
\end{enumerate}
\elemm

\dem (1) For all $t\geq 0$ and $[x,s]\in\Sigma_r$, let $n\in\NN$ and
$s'\geq 0$ be such that
$$
t+s=\sum_{i=0}^{n-1}r_+(\sigma_+^i\pi_+(x)) + s'\;\;\;{\rm and}\;\;\;
0\leq s'<r_+(\sigma_+^n\pi_+(x))\;.
$$ 
Since $r_+\circ\sigma_+^i\circ \pi_+=r\circ\sigma^i$ for every
$i\in\NN$, these two conditions are equivalent to
$$
t+s=\sum_{i=0}^{n-1}r(\sigma^ix) + s'\;\;\;{\rm and}\;\;\;0\leq
s'< r(\sigma^nx)\,.
$$ 
Hence
$$
(\sigma_+)^t_{r_+}\circ \pi_+^{\;r}([x,s])
= (\sigma_+)^t_{r_+}([\pi_+(x),s])=[\sigma_+^n\pi_+(x),s']
$$
and
$$
\pi_+^{\;r}\circ\sigma_r^t([x,s])=\pi_+^{\;r}([\sigma^nx,s'])
=[\pi_+(\sigma^nx),s']\;.
$$
This proves Assertion (1) since $\pi_+\circ\sigma=\sigma_+\circ\pi_+$.

\medskip\noindent 
(2) By Proposition \ref{prop:defiBowenWalters}, we may assume that, in
the formula of the H\"older norms, the Bowen-Walters distance is replaced by
the function $d_{\rm BW}$, as this will only change $C_1$ by $C_{\rm
  BW}^\alpha \,C_1$.

For every $[x,s]\in\Sigma_r$, with $[y,s']$ and $N$ associated with
$\pi_+^{\;r}([x,s])=[\pi_+(x),s]$ as in the definition of
$\psi^{(T)}([\pi_+(x),s])$, we have
$$
d_{\rm BW}(\sigma_r^T[x,s],[y,s'])= d_{\rm BW}([\sigma^Nx,s'],[y,s'])\leq
d(\sigma^Nx,y)\leq e^{-N}\;.
$$ 
Since the positive roof function $r$ is bounded from above by the
least upper bound $\sup \lambda$ of the lengths of the edges, we have
$$
N\geq\; \sum_{i=0}^{N-1}\;\frac{r(\sigma^ix)}{\sup \lambda}\;=
\frac{1}{\sup \lambda}\;(s+T-s')\geq \frac{T}{\sup \lambda}-1\;.
$$
Hence
\begin{align*}
|\psi\circ\sigma_r^T([x,s])-\psi^{(T)}(\pi_+^{\;r}(x))| & = 
|\psi(\sigma_r^T[x,s])-\psi([y,s'])|
\\ & \leq \|\psi\|'_\alpha \;d_{\rm BW}(\sigma_r^T[x,s],[y,s'])^\alpha \leq 
\|\psi\|'_\alpha\; e^{-\frac{\alpha}{\sup\lambda}\, T +\alpha}\;.
\end{align*}

\medskip
\noindent 
(3) Let us prove that $\psi^{(T)}$ is $\C^k$ along semiflow lines. Fix
$[x,s]\in (\Sigma_+)_{r_+}$.  With $[y,s']$ and $N$ as in the
construction of $\psi^{(T)}([x,s])$, let us consider $\epsilon>0$
small enough, so that $\epsilon<r_+(x)-s$ and $\epsilon <
r_+(\sigma_+^N(x))-s'=r(y)-s'$. Then
$$
\psi^{(T)}\circ(\sigma_+)_{r_+}^\epsilon([x,s])
=\psi^{(T)}([x,s+\epsilon])=\psi([y,s'+\epsilon])=
\psi\circ\sigma_r^\epsilon([y,s'])\;.
$$ 
Therefore by taking derivatives with respect to $\epsilon$ in this
formula, $\psi^{(T)}$ is indeed $\C^k$ along semiflow lines, and, for
$i=0,\dots,k$, we have
\begin{equation}\label{eq:commutTpartiali}
\partial^i_t\big(\psi^{(T)}\big)=\big(\partial^i_t\psi\big)^{(T)}\;.
\end{equation}

The inequality $\|\psi^{(T)}\|_\infty\leq \|\psi\|_\infty$ is
immediate by construction. Using the above centred equation, we have
$\|\psi^{(T)}\|_{k,\,\infty}\leq \|\psi\|_{k,\,\infty}$.

Let us prove that there exists a constant $C_2\geq 1$ (independent of
$T$ and $\psi$) such that
\begin{equation}\label{eq:majonormTalpha}
\|\psi^{(T)}\|'_\alpha\leq C_2\;e^{\alpha''\,T}\;\|\psi\|'_\alpha\;.
\end{equation}

Let $[x,s]\in(\Sigma_+)_{r_+}$, take $[y,s']$ and $N$ as in the
definition of $\psi^{(T)}([x,s])$. Let $[\underline{x},
  \underline{s}\,]\in(\Sigma_+)_{r_+}$, take $[\underline{y},
  \underline{s}']$ and $\underline{N}$ as in the definition of
$\psi^{(T)}([\underline{x}, \underline{s}\,])$. Up to exchanging
$[x,s]$ and $[\underline{x}, \underline{s}\,]$, we assume that
$\underline{N}\geq N$.

By Proposition \ref{prop:defiBowenWalters},\footnote{See the previous
  footnote.} we may assume that, in the H\"older norms formulas, the
Bowen-Walters distance is replaced by the function $d_{\rm BW}$, as
this will only change $C_2$ by $C_{\rm BW}^{2\,\alpha} \,C_2$.  Let
$$
C_3=\min\{e^{-1},\; \inf\lambda\}\;.
$$

We have
\begin{equation}\label{eq:majopassTun}
\big|\psi^{(T)}([x,s])-\psi^{(T)}([\underline{x}, \underline{s}\,])\big|
=\big|\psi([y,s'])-\psi([\underline{y}, \underline{s'}])\big|
\leq \|\psi\|'_\alpha \;
d_{\rm BW}([y,s'],[\underline{y}, \underline{s'}])^\alpha\;.
\end{equation}
Note that the map $d_{\rm BW}$ on $\Sigma_r\times \Sigma_r$ is
bounded from above by $1 + \sup \lambda$, since the distance on
$\Sigma$ is at most $1$ and since the roof function $r$ is bounded
from above by $\sup \lambda$. 
If $d_{\rm BW}([x,s],[\underline{x}, \underline{s}\,])\geq
C_3\;e^{-\frac{T}{\inf\lambda}}$, then
$$ 
d_{\rm BW}([y,s'],[\underline{y}, \underline{s'}]) \leq 1 + \sup
\lambda \leq \frac{1 + \sup \lambda}{C_3}\;e^{\frac{T}{\inf\lambda}}\;
d_{\rm BW}([x,s],[\underline{x}, \underline{s}\,])\;.
$$
Therefore Equation \eqref{eq:majonormTalpha} follows from Equation
\eqref{eq:majopassTun} whenever $C_2\geq \frac{(1 + \sup
  \lambda)^\alpha}{C_3^\alpha}$.

Conversely, suppose that $d_{\rm BW}([x,s],[\underline{x},
  \underline{s}\,])< C_3\;e^{-\frac{T}{\inf\lambda}}$.  Assume that
$$
d_{\rm BW}([x,s],[\underline{x}, \underline{s}\,])= d(x,\underline{x})
+|s-\underline{s}|\;,
$$ 
the other possibilities are treated similarly.  Since 
$$
s+T=\sum_{i=0}^{N-1} r_+(\sigma_+^ix)\;+ s'
$$ 
and since the roof function $r_+$ is bounded from below by $\inf
\lambda$, we have $T\geq N\inf\lambda - \inf\lambda$, or equivalently
$N\leq \frac{T}{\inf\lambda}+1$.  Hence $d(x,\underline{x}) < C_3
\;e^{-\frac{T}{\inf\lambda}}\leq e^{-N}$ by the definition of $C_3$.
In particular the sequences $x$ and $\underline{x}$ indexed by $\NN$
have the same $N+1$ first coefficients. Since $r_+(z)$ depends only on
$z_0$ for all $z\in\Sigma_+$, we thus have $r_+(\sigma_+^ix)=
r_+(\sigma_+^i\,\underline{x})$ for $i=0,\dots, N$. Note that we have
$$
\underline{s}+T= \sum_{i=0}^{\underline{N}-1}
r_+(\sigma_+^i\,\underline{x}) \;+ \underline{s'}\;.
$$ 

If $\underline{N}=N$, then by taking the difference of the last two
centred equations, we have $\underline{s}-s =\underline{s'}-s'$, and
by construction, the sequences $y$ and $\underline{y}$ indexed by
$\ZZ$ satisfy $y_i=(\underline{y})_i$ if $i\leq 0$ and if $0\leq i\leq
-\ln d(x,\underline{x}) - N$. Therefore $d(y,\underline{y}) \leq
e^N\;d(x,\underline{x})$ and
\begin{align*}
d_{\rm BW}([y,s'],[\underline{y}, \underline{s'}])&\leq
d(y,\underline{y})+|s'-\underline{s'}|\leq 
e^N\;d(x,\underline{x})+|\underline{s}-s|\\ & \leq
e^N\;d_{\rm BW}([x,s],[\underline{x}, \underline{s}\,])
\leq e^{\frac{T}{\inf\lambda}+1}\;d_{\rm BW}([x,s],
[\underline{x}, \underline{s}\,]) \;.
\end{align*}
Therefore Equation \eqref{eq:majonormTalpha} follows from Equation
\eqref{eq:majopassTun} whenever $C_2\geq e^\alpha$.

If  $\underline{N}>N$, then again by difference
$$
\underline{s}-s=\sum_{i=N+1}^{\underline{N}-1}r_+(\sigma_+^i\,\underline{x})
+ r_+(\sigma_+^Nx) - s' + \underline{s'}.
$$ 
Note that $\underline{s'}\geq 0$, that $r_+(\sigma_+^Nx) - s'\geq
0$, and that $|\underline{s}-s|< C_3 \;e^{-\frac{T}{\inf\lambda}}\leq
\inf\lambda$ by the definition of $C_3$. Hence we have
$\underline{N} =N+1$ and $\underline{s}-s=r_+(\sigma_+^Nx) - s' +
\underline{s'}$.  By construction, the sequences $\sigma y$ and
$\underline{y}$ indexed by $\ZZ$ satisfy $(\sigma
y)_i=(\underline{y})_i$ if $i\leq 0$ and if $0\leq i\leq -\ln
d(x,\underline{x}) - N-1$. Hence by the definition of $d_{\rm BW}$ and
since $r(y)=r_+(\sigma_+^Nx)$, we have
\begin{align*}
d_{\rm BW}([y,s'],[\underline{y}, \underline{s'}])&\leq
d(\sigma y,\underline{y})+r(y)- s' + \underline{s'}\leq 
e^{N+1}\;d(x,\underline{x})+|\underline{s}-s|\\ & \leq
e^{N+1}\;d_{\rm BW}([x,s],[\underline{x}, \underline{s}\,])
\leq e^{\frac{T}{\inf\lambda}+2}\;d_{\rm BW}([x,s],
[\underline{x}, \underline{s}\,]) \;.
\end{align*}
Therefore Equation \eqref{eq:majonormTalpha} follows from Equation
\eqref{eq:majopassTun} whenever $C_2\geq e^{2\alpha}$. This ends the proof of 
 Equation \eqref{eq:majonormTalpha}.

\medskip
Now note that Equations \eqref{eq:commutTpartiali} and
\eqref{eq:majonormTalpha} imply Equation \eqref{eq:controlckalphanorm}
by summation (using the independence of $C_2$ on $\psi$), thus
concluding the proof of Lemma \ref{lem:passageauT}.  
\cqfd

\bprop\label{prop:decaycorrpassonesidsusp} 
Let $\mu$ be a $(\sigma_r^t)_{t\in\RR}$-invariant probability measure
on $\Sigma_r$.  Assume that the dynamical system $\big((\Sigma_+)_{r_+},
((\sigma_+)^t_{r_+})_{t\in\RR},(\pi_+^{\;r})_*\mu\big)$ has
superpolynomial decay of $\alpha$-H\"older correlations. Then
$(\Sigma_r, (\sigma_r^t)_{t\in\RR},\mu)$ has superpolynomial decay of
$\alpha$-H\"older correlations.  
\eprop

\dem 
We fix $n\in\NN$. Let $N=1+2\lceil\frac{\sup\lambda}{\inf\lambda}
\rceil$.  Let $k\in\NN$ and $C_4>0$ (depending on $n$) be such that
for all $\psi,\psi'\in\C^{k,\,\alpha}_{\rm b}((\Sigma_+)_{r_+})$, we have for
all $t\geq 1$
\begin{equation}\label{eq:covsusponesid}
\big|\operatorname{cov}_{(\pi_+^{\;r})_*\mu,\,t}(\psi,\psi')\big| \leq 
C_4\;\|\psi\|_{k,\,\alpha}\;\|\psi'\|_{k,\,\alpha}\; t^{N\, n}\;.
\end{equation}

Now let $\psi,\psi'\in\C^{k,\,\alpha}_{\rm b}(\Sigma_r)$. We again
denote by $\pm \,a$ any value in $[-a,a]$ for every $a\geq 0$.  By
invariance of $\mu$ under $(\sigma_r^t)_{t\in\RR}$, by Lemma
\ref{lem:passageauT} (2) and by Lemma \ref{lem:passageauT} (1), we
have, for every $T\geq 0$ (to be chosen appropriately later on),
\begin{align*}
& \int_{\Sigma_r} \psi\circ \sigma^t_r\;\;\psi'\; d\mu = 
\int_{\Sigma_r} \psi\circ\sigma_r^{T+t}\;\;
\psi'\circ \sigma_r^{T}\; d\mu \\  = \;&
\int_{\Sigma_r} (\psi^{(T)}\circ\pi_+^{\;r} \pm 
C_1\;\|\psi\|'_\alpha\;e^{-\alpha'\,T})\circ \sigma_r^{t}\;\; 
({\psi'}^{(T)}\circ\pi_+^{\;r} \pm 
C_1\;\|\psi'\|'_\alpha\;e^{-\alpha'\,T})\; d\mu\\   = \; & 
\int_{(\Sigma_+)_{r_+}} \psi^{(T)}\circ(\sigma_+)_{r_+}^{t} \;\;{\psi'}^{(T)}\; 
d(\pi_+^{\;r})_*\mu\;\;
\pm\;C_1^2\;\|\psi\|_\alpha\;\|\psi'\|_\alpha\;e^{-\alpha'\,T}\;.
\end{align*}
A similar estimate holds for the second term in the definition of the
correlation coefficients.  Hence, applying Equation
\eqref{eq:covsusponesid} to the observables $\psi^{(T)}$ and
${\psi'}^{(T)}$, by Lemma \ref{lem:passageauT} (3), we have since
$C_2\, e^{\alpha''\,T}\ge 1$,
\begin{align*}
& |\operatorname{cov}_{\mu,\,t}(\psi,\psi')| \leq 
|\operatorname{cov}_{(\pi_+^{\;r})_*\mu,\,t}(\psi^{(T)},{\psi'}^{(T)})| 
+  2\;C_1^2\;\|\psi\|_\alpha\;\|\psi'\|_\alpha\;e^{-\alpha'\, T} \\ \leq \;&
C_4\;(\|\psi\|_{k,\,\infty}+ C_2\;\|\psi\|'_{k,\,\alpha}\; e^{\alpha''\, T})\,
(\|\psi'\|_{k,\,\infty}+C_2\;\|\psi'\|'_{k,\,\alpha}\; e^{\alpha''\, T})\; 
t^{-N\,n} \\  \;& + 2\;C_1^2\;\|\psi\|_\alpha\;\|\psi'\|_\alpha\;e^{-\alpha'\,T} 
\\ \leq \;& \|\psi\|_{k,\,\alpha}\;\|\psi'\|_{k,\,\alpha}\;
(C_4\;C_2^2\;e^{2\alpha'' T}\; t^{-N\,n} +2\;C_1^2\;e^{-\alpha'\,T})\;.
\end{align*}
Take $T=\frac{n}{\alpha'}\;\ln t\geq 0$. Since $N=1+2\lceil
\frac{\alpha''}{\alpha'}\rceil$, we have $2\alpha''\;
\frac{n}{\alpha'}-N\,n\leq -n$. Hence with $C_5=C_4\;C_2^2
\;+\;2\;C_1^2$, we have for all $t\geq 1$
$$
|\operatorname{cov}_{\mu,\,t}(\psi,\psi')| \leq 
C_5\;\|\psi\|_{k,\,\alpha}\;\|\psi'\|_{k,\,\alpha}\;t^{-n}\;.
$$
This concludes the proof of Proposition \ref{prop:decaycorrpassonesidsusp}.
\cqfd

\subsection*{Step 3 : Conclusion of the proof of Theorem
  \ref{theo:mixingtreelattice}.}

In this Step, we prove that the semiflow $\big((\Sigma_+)_{r_+},
\big((\sigma_+)^t_{r_+} \big)_{t\geq 0},
\frac{(\PP_+)_{r_+}}{\|(\PP_+)_{r_+}\|}\big)$ is rapidly mixing, which
concludes the proof of Theorem \ref{theo:mixingtreelattice}, using
Proposition \ref{prop:decaycorrpassonesidsusp} with $\mu=
\frac{\PP_r}{\|\PP_r\|}$.

Recall\footnote{See the proof of Theorem \ref{theo:mainsuspension}.}
that $Y= \{\ell\in \Ga\bs\G X\;:\;\ell(0)\in V\XX\}$ is a
cross-section of the geodesic flow on $\Ga\bs\G X$, and that if
$R:Y\ra \Ga\bs\G \XX$ is the reparametrisation map of $\ell\in Y$ to a
discrete geodesic line $^\sharp\ell\in \Ga\bs\G \XX$ with the same
origin, then the measure $\mu_Y$, induced by the Gibbs measure $m_c$
on the cross-section $Y$ by disintegration along the flow, maps by
$R_*$ to a constant multiple of $m_{^\sharp c}=m_{F_{^\sharp c}}$ (see
Lemma \ref{lem:samePatdenssR} (2)).  Hence for all $n\in\NN-\{0\}$ and
$\epsilon\in\mathopen{]}0,\frac{1}{2}\inf\lambda\mathclose{[}\,$, by
Assumption (b) (1) in the statement of Theorem
\ref{theo:mixingtreelattice}, we have
\begin{align*}
m_{^\sharp c} & \Big(\Big\{\,^\sharp\ell\in\Ga\bs\G \XX\;:\; 
\begin{array}{c}\,^\sharp\ell(0)\in E \\
\forall\;k\in\{1,\dots,n-1\},\;\,^\sharp\ell(k)\notin E
\end{array}\Big\}\Big) \\ & \\ 
\leq \; & \frac{\|m_{^\sharp c}\|}{\|\mu_Y\|}\;
\mu_Y\Big(\Big\{\,R^{-1}({}^\sharp\ell)\in Y\;:\; 
\begin{array}{c}R^{-1}({}^\sharp\ell)(0)\in E \\
\forall\;t\in\;]0,n\inf\lambda[ \,,\;   R^{-1}({}^\sharp\ell)(t)\notin E
\end{array}\Big\}\Big) \\ & \\ 
\leq \; & \frac{\|m_{^\sharp c}\|}{\epsilon\,\|\mu_Y\|}\;
m_c\Big(\Big\{\,\flow{s}R^{-1}({}^\sharp\ell)\in \Ga\bs\G X\;:\; 
\begin{array}{c}0\leq s\leq\epsilon,\;\;
d(\flow{s}R^{-1}({}^\sharp\ell)(0),E)\leq\epsilon \\
\forall\;t\in\;]\epsilon,n\inf\lambda-\epsilon[ \,,\;   
\flow{s}R^{-1}({}^\sharp\ell)(t)\notin E\end{array}\Big\}\Big)\\ & \\ 
\leq \; &  \frac{\|m_{^\sharp c}\|}{\epsilon\,\|\mu_Y\|}\;
C\;e^{-\kappa\,(\inf\lambda)\,n+\kappa\,\epsilon}\;.
\end{align*}
Therefore Equation \eqref{eq:expodecaycusp} (where $c$ is replaced by
$\sharp c$) is satisfied, with $C'= \frac{\|m_{^\sharp c}\| \,C\,
  e^{\kappa\,\epsilon}}{\epsilon\,\|\mu_Y\|}$ and $\kappa'= \kappa\,
\inf\lambda$. As seen in the proof of Theorem
\ref{theo:critexpdecaysimpl}, this implies that there exists a finite
subset $\E$ of the alphabet $\A$ such that Equation
\eqref{eq:expodecaycuspdynsymbonesid} is satisfied.

We now apply \cite[Theo.~2.3]{Melbourne07} with the dynamical system
$(X,m_0,T)=(\Sigma_+,\PP_+,\sigma_+)$ (using the system of
conductances $^\sharp c$) and the roof function $h=r_+$. This
dynamical system is presented as a Young tower in Step 3 of the proof
of Theorem \ref{theo:critexpdecaysimpl}.  Equation
\eqref{eq:expodecaycuspdynsymbonesid} for the first return map
$\tau_\E$ and the $4$-Diophantine hypothesis are exactly the
hypothesis needed in order to apply \cite[Theo.~2.3]{Melbourne07}.
Thus the semiflow $((\Sigma_+)_{r_+}, \big((\sigma_+)^t_{r_+}
\big)_{t\geq 0}, \frac{(\PP_+)_{r_+}}{\|(\PP_+)_{r_+}\|})$ has
superpolynomial decay of $\alpha$-H\"older correlations.

When $\Ga\bs X$ is compact, the alphabet $\A$ is finite and
$(\Sigma_+,\sigma_+,\PP_+)$ is a (one-sided) subshift of finite type,
hence we do not need the exponentially small tail assumption, but only
the $2$-Diophantine hypothesis, and we may apply \cite{Dolgopyat98ETDS}.
\cqfd

\bcoro\label{coro:expdecaygeomfinimet} Assume that the Gibbs measure
$m_F$ is finite and mixing for the (continuous time) geodesic flow,
that the lengths of the edges of $(\XX,\lambda)$ have a finite upper
bound, and that $\Ga$ is geometrically finite. There exists a full
measure subset $A$ of $\RR^4$ (for the Lebesgue measure) such that if
$\Ga$ has a quadruple of translation lengths in $A$, or if the length
spectrum is $4$-Diophantine, then the (continuous time) geodesic flow
on $\Ga\bs\G X$ has superpolynomial decay of $\alpha$-H\"older
correlations for the Bowen-Margulis measure $m_{\rm BM}$.  \ecoro

\dem The exponentially small tail Assumption (b) (1) is checked as in
the proof of Corollary \ref{coro:expdecaygeomfinisimpl}. The deduction
of Corollary \ref{coro:expdecaygeomfinimet} from Theorem
\ref{theo:mixingtreelattice} then proceeds, by an argument going back
in part to Dolgopyat, as for the deduction of Corollary 2.4 from Theorem 2.3
in \cite{Melbourne07}.  
\cqfd

\medskip 
Note that under the general assumptions of Theorem \ref{theo:mixing},
the geodesic flow on $\Ga\bs\G X$ might not be exponentially mixing,
see for instance \cite[page 162]{Pollicott86} or \cite{Ruelle83} for
analogous behaviour.

\part{Geometric equidistribution and counting}
\label{part:equid}

\chapter{Equidistribution of equidistant level 
sets to Gibbs measures}
\label{sec:skinningwithpot}

Let $X$ be a geodesically complete proper $\CAT(-1)$ space, let $\Ga$
be a nonelementary discrete group of isometries of $X$, let $\wt F$ be
a continuous $\Ga$-invariant map on $T^1X$ such that
$\delta=\delta_{\Ga,\,F^\pm}$ is finite and positive and that the
triple $(X,\Ga,\wt F)$ satisfies the \ref{eq:HC}-property,\footnote{See
Definition \ref{defi:HCproperty}.} and let $(\mu^\pm_x)_{x\in X}$ be
Patterson densities for the pairs $(\Ga,F^\pm)$.

In this Chapter, we prove that the skinning measure on (any nontrivial
piece of) the outer unit normal bundle of any properly immersed
nonempty proper closed convex subset of $X$, pushed a long time by the
geodesic flow, equidistributes towards the Gibbs measure, under
finiteness and mixing assumptions. This result gives four important
extensions of \cite[Theo.~1]{ParPau14ETDS}, one for general $\CAT(-1)$
spaces with constant potentials, one for Riemannian manifolds with
pinched negative curvature and H\"older-continuous potentials,  one for
$\RR$-trees with general potentials, and one for simplicial trees.

\section{A general equidistribution result}
\label{subsec:skinningwithpot}

Before stating this equidistribution result, we start by a technical
construction which will also be useful in the following Chapter
\ref{sec:equidarcs}.  We refer to Section \ref{subsect:nbhd} for the
notation concerning the dynamical neighbourhoods (including
$V^\mp_{w,\,\eta',\,R}$) and to Chapter \ref{sect:skinning} for the
notation concerning the skinning measures (including $\nu^\pm_w$).

\bigskip\noindent{\bf Technical construction of bump functions.}  Let
$D^\pm$ be nonempty proper closed convex subsets of $X$, and let $R>0$
be such that $\nu^\pm_w (V^\mp_{w,\,\eta'',\,R})>0$ for all $\eta''>0$
and $w\in \normalmp D^\pm$. Let $\eta>0$ and let $\Omega^\pm$ be
measurable subsets of $\normalmp D^\pm$. We now construct functions
$$
\phi^\mp_{\eta,\,R,\,\Omega^\pm} : \G X\to[0,+\infty[
$$ 
whose supports are contained in dynamical neigbourhoods of
$\Omega^\pm$.  If $X=\wt M$ is a Riemanian manifold and $\wt F=0$, we
recover the same bump functions as in \cite{ParPau14ETDS} after the
standard identifications.

For all $\eta'>0$, let $h^\pm_{\eta,\,\eta'}: \G_\mp X\to [0,+\infty[$
be the $\Ga$-invariant measurable maps defined by
\begin{equation}\label{eq:defihpm}
h^\pm_{\eta,\,\eta'}(w)=\frac{1}{\nu^\pm_w(V^\mp_{w,\,\eta,\,\eta'})}
\end{equation}
if $\nu^\pm_w(V^\mp_{w,\,\eta,\,\eta'})>0$ (which is for instance
satisfied if $w_\pm\in\Lambda\Ga$ and for every $w\in\normalmp D^\pm$
if $\eta'=R$ by the choice of $R$) and $h^\pm_{\eta,\,\eta'}(w) = 0$
otherwise.

These functions $h^\pm_{\eta,\,\eta'}$ have the following behaviour
under precomposition by the geodesic flow. By Lemma \ref{lem:defiR}
(2), by Equation \eqref{eq:flotVetaeta}, and by the invariance of
$\nu^\pm_w$ under the geodesic flow, we have, for all $t\in\RR$ and
$w\in\G_\pm X$,
\begin{equation}\label{eq:commutflowh}
  h^\mp_{\eta,\,\eta'}(\flow{\mp t}w)=
  e^{C^\pm_{w_\pm}(w(0),\,w(\mp t))}\; h^\mp_{\eta,\,e^{-t}\eta'}(w)\,.
\end{equation}

Let us also describe the behaviour of $h^\pm_{\eta,\,\eta'}$ when
$\eta'$ is small. Let $w\in\G_\pm X$ be such that $w$ is isometric at
least on $\pm[0,+\infty[$, which is for instance the case if $w\in 
\normalpm D^\mp$. For all $\eta'>0$ and $\ell\in B^\pm(w,\eta')$, let
$\wh w$ be an extension of $w$ such that $d_{W^\pm(w)}(\ell,\wh w) <
\eta'$. Then $\wh w(0)=w(0)$ by the assumption on $w$, and using Lemma
\ref{lem:comparddHam}, we have
$$
d(\ell(0),w(0))=d(\ell(0),\wh w(0))\leq d_{W^\pm(w)}(\ell,\wh w) 
<\eta'\;.
$$
Hence, with $\kappa_1$ and $\kappa_2$ the constants in Definition
\ref{defi:HCproperty}, if $\eta'\leq 1$ and
$$
c_1=\kappa_1+2\,\delta+2
\sup_{\pi^{-1}(B(w(0),\,2))} |\wt F|\;,
$$
we have, by Proposition \ref{prop:continuGibbscocycle}
\eqref{eq:cocycleLip},
$$
|\;C^\mp_{w_\mp}(w(0),\ell(0))\;| \leq c_1\,(\eta')^{\kappa_2}\;.
$$ 
Using the defining Equation \eqref{eq:definumoins} of $\nu^\mp_w$,
for all $s\in\RR$, $\eta'\in\mathopen{]}0,1\mathclose{]}$ and $\ell\in
  B^\pm(w,\eta')$, we have
$$
e^{-c_1(\eta')^{\kappa_2}}\;ds\,d\mu_{\wssu(w)}(\ell)
\leq d\nu^\mp_w(\flow{s}\ell)\leq 
e^{c_1(\eta')^{\kappa_2}}\;ds\,d\mu_{\wssu(w)}(\ell)\,.
$$
It follows that for all $\eta'\in\;]0,1]$ and $w\in\normalpm D^\mp$
such that $w_\pm\in\Lambda\Ga$, we have the following control of
$h^\mp_{\eta,\,\eta'}(w)$:
\begin{equation}\label{eq:estimhsmaletaprimp}
\frac{e^{-c_1(\eta')^{\kappa_2}}}{2\eta\;\mu_{\wssu(w)}(B^\pm(w,\eta'))}
\leq h^\mp_{\eta,\,\eta'}(w)\leq
\frac{e^{c_1(\eta')^{\kappa_2}}}{2\eta\,\mu_{\wssu(w)}(B^\pm(w,\eta'))}\,.
\end{equation}
Note that when $X$ is an $\RR$-tree, we may take $\kappa_2=1$ and
$c_1=\sup_{\pi^{-1}(B(w(0),1))}|\wt F-\delta|$ in this equation, as
said in the last claim of Proposition \ref{prop:continuGibbscocycle}
\eqref{eq:cocycleLip}.  Note that $c_1$ is bounded when $w$ ranges
over any compact subset of $\G_\pm X$, and is uniformly bounded when
$\wt F$ is bounded.

\medskip
Recall that $\mathbbm{1}_A$ denotes the characteristic function of a
subset $A$. We now define the test functions
$\phi^\mp_{\eta,\,R,\,\Omega^\pm}: \G X\to[0,+\infty[$ with support in
a dynamical neighbourhood of $\Omega^\mp$ by
\begin{equation}\label{eq:defiphi}
\phi^\mp_{\eta,\,R,\,\Omega^\mp}=\;h^\mp_{\eta,\,R}\circ f^\pm_{D^\mp}\;\;
\mathbbm{1}_{\V^\pm_{\eta,\,R}(\Omega^\mp)}\,,
\end{equation}
where $\V^\pm_{\eta,\,R}(\Omega^\mp)$ and $f^\pm_{D^\mp}$ are as in
Section \ref{subsect:nbhd}.  Note that if $\ell\in \V^\pm_{\eta, \,R}
(\Omega^\mp)$, then $\ell_\pm \notin \partial_\infty D^\mp$ by
convexity. Thus,  $\ell$ belongs to the domain of definition
$\U^\pm_{D^\mp}$ of $f^\pm_{D^\mp}$. Hence $\phi^\mp_{\eta,\,R,\,
  \Omega^\mp} (\ell)= h^\mp_{\eta, \,R} \circ f^\pm_{D^\mp}(\ell)$ is
well defined. By convention, $\phi^\mp_{\eta,\,R,\,\Omega^\mp}(\ell)
=0$ if $\ell\notin \V^\pm_{\eta,\,R} (\Omega^\mp)$.

\medskip
We now globalise these test functions in order to apply them to
equivariant families of supports. 

Let $\eta>0$.  Let $\D=(D_i)_{i\in I}$ be a locally finite
$\Ga$-equivariant family of nonempty proper closed convex subsets of
$X$ with $\Ga\bs I$ finite, and $\sim=\sim_\D$. Let $R>0$ be such that
$\nu_w^\pm(V^\mp_{w,\,\eta'',\,R})>0$ for all $\eta''>0$, $i\in I$ and
$w\in\normalpm D_i$.  Let $\Omega= (\Omega_i)_{i\in I}$ be a locally
finite $\Ga$-equivariant family of measurable subsets of $\gengeod X$,
with $\Omega_i\subset \normalpm D_i$ for all $i\in I$ and $\Omega_i=
\Omega_j$ if $i\sim j$. We define the global test functions
$\wt\Phi^\mp_\eta:\G X\to[0,+\infty[$ by
\begin{equation}\label{eq:defglobaltestfunction}
\wt\Phi^\mp_\eta=\sum_{i\in I/\sim} \phi^\mp_{\eta,\,R,\,\Omega_i}
=\sum_{i\in I/\sim}\;h^\mp_{\eta,\,R}\circ f^\pm_{D^\mp}\;\;
\mathbbm{1}_{\V^\pm_{\eta,\,R}(\Omega^\mp)}\;.
\end{equation}

A subset $\Delta_\Ga$ of $\gengeod X$ is a {\em fundamental domain}
for the action of $\Ga$ if the interiors of its translates are
disjoint and any compact subset of $\gengeod X$ meets only finitely
many translates of $\Delta_\Ga$: If $m_F$ is finite, a fundamental
domain with boundary of zero measure exists by \cite[p.~13]{Roblin03},
using the fact that $\wt m_F$ has no atoms according to Corollary
\ref{coro:finitudeGibbsdivuniq} (1) and to Theorem \ref{theo:HTSR}.

The following properties of the bump functions are proved as in
\cite[Prop.~18]{ParPau14ETDS}.

\blemm \label{lem:integrable} (1) For every $\eta>0$, the functions
$\phi^\mp_{\eta,\,R,\,\Omega^\pm}$ are measurable, nonnegative and
satisfy
$$
\int_{\G X} \phi^\mp_{\eta,\,R,\,\Omega^\pm}\;d\wt m_{F}=\wt
\sigma^\pm_{D^\mp}(\Omega^\mp)\,.
$$
(2) For every $\eta>0$, the function $\wt\Phi^\mp_\eta$ is well
defined, measurable and $\Ga$-invariant. It defines, by passing to the
quotient, a measurable function $\Phi^\mp_\eta: \Ga\bs\G
X\to[0,+\infty[\,$ such that
\begin{equation}\label{eq:inttestskinmass}
\int_{\Ga\bs\G X}\Phi^\mp_\eta\;dm_{F}= \|\sigma^\pm_\Omega\|\,.
\end{equation}
\elemm 

\dem (1) Recall that the fiber of the restriction of $f^\pm_{D^\mp}$
to $\V^\pm_{\eta,\,R}(\Omega^{\mp})$ over $w\in \Omega^{\mp}$ is the
open subset $V^\pm_{w,\,\eta,\,\eta'}$ of $\wosu(w)$.  By the
disintegration result of Proposition \ref{prop:disintegration}, by the
definition of $h^\mp_{\eta,\,R}$ and by the choice of $R$, we have
\begin{align*}
\int_{\G X} \phi^\mp_{\eta,\,R,\,\Omega^\pm}\;d\wt m_F& = 
\int_{\ell\in \V^\pm_{\eta,\,R}(\Omega^{\mp})}
h^\mp_{\eta,\,R}\circ f^\pm_{D^\mp}(\ell)\;d\wt m_F(\ell)\\ &
=\int_{w\in \Omega^{\mp}}h^\mp_{\eta,\,R}(w)\;\int_{\ell\in V^\pm_{w,\,\eta,\,R}} \;
d\nu^\mp_w(\ell)\;d\wt\sigma^\pm_{D^\mp}(w)
=\wt \sigma^\pm_{D^\mp}(\Omega^\mp)\,.
\end{align*}

\smallskip
\noindent 
(2) The function $\wt\Phi^\mp_\eta$ is well defined, since $\Omega_i =
\Omega_j$ and thus $\V^\pm_{\eta,\,R}(\Omega_i)= \V^\pm_{\eta,\,R}
(\Omega_j)$ if $i\sim j$, since $h^\mp_{\eta,\,R} \circ f^\pm_{D_i}
(\ell)$ is finite if $\ell\in \V^\pm_{\eta,\,R} (\Omega_i)$ (by the
definition of $R$), and since the sum defining $\wt\Phi^\mp_\eta$ has
only finitely many nonzero terms, by the local finiteness of the
family $\Omega$ (given $\ell\in\G X$, the summation over $I/\!\!\sim$
giving $\wt\Phi^\mp_\eta(\ell)$ may be replaced by a summation over
the finite set $\{i\in I: \ell\in \V^\pm_{\eta,\,R}(\Omega_i)\}/\!\!
\sim$).

The function $\wt\Phi^\mp_\eta$ is $\Ga$-invariant since
$$
\mathbbm{1}_{\V^\pm_{\eta,\,R}(\Omega_i)}\circ \ga=
\mathbbm{1}_{\ga^{-1}\V^\pm_{\eta,\,R}(\Omega_i)}=
\mathbbm{1}_{\V^\pm_{\eta,\,R}(\Omega_{\ga^{-1}i})}
$$ 
and\footnote{See Equation \eqref{eq:equivfibrationf}.}
$$
h^\mp_{\eta,\,R} \circ f^\pm_{D_i} \circ \ga=
h^\mp_{\eta,\,R} \circ \ga \circ f^\pm_{\ga^{-1}D_i} =
h^\mp_{\eta,\,R}\circ f^\pm_{D_{\ga^{-1}i}}
$$ 
and by a change of index in Equation \eqref{eq:defglobaltestfunction}.
If $\Delta_\Ga$ is a fundamental domain for the action of $\Ga$, we
have by Assertion (1)
\begin{align*}
\int_{\Ga\bs\G X}\Phi^\mp_\eta\;dm_{F}=
\int_{\Delta_\Ga}\wt\Phi^\mp_\eta\;d\wt m_{F} & =
\sum_{i\in I/\sim} \int_{\G X} \phi^\mp_{\eta,\,R,\,\Omega_i\cap\Delta_\Ga}\;d\wt m_{F}\\
& =
\sum_{i\in I/\sim} \wt\sigma^\pm_{D_i}(\Delta_\Ga\cap \Omega_i)
= \wt\sigma^\pm_\Omega(\Delta_\Ga)=\|\sigma^\pm_\Omega\|\,. \;\;\;\Box
\end{align*}

\medskip

We now state and prove the aforementioned equidistribution result.
Note that as the elements of the outer unit normal bundles are 
only geodesic rays on $[0,+\infty[$, their pushforwards by the
geodesic flow at time $t$ are geodesic rays on $[-t,+\infty[$ and
the convergence towards geodesic lines (defined on $]-\infty,
+\infty[\,$) does take place in the full space of generalised geodesic
lines $\gengeod X$. This explains why it is important not to forget to
consider the negative times in order for the skinning measures,
supported on geodesic rays, when pushed by the geodesic flow, to have
a chance to weak-star converge to Gibbs measures, supported on
geodesic lines, up to renormalisation.

The proof of the following result has similarities with that of
\cite[Theo.~1]{ParPau14ETDS}, but the computations do not apply in the
present context because the proof in loc.~cit.~does not keep track of
the past: here we can no longer reduce our study to the outer unit
normal bundle of the $t$-neighbourhood of the elements of $\D$.

\btheo\label{theo:equid} Let $(X, \Ga,\wt F)$ be as in the beginning
of Chapter \ref{sec:skinningwithpot}. Assume that the Gibbs measure
$m_{F}$ on $\Ga\bs\G X$ is finite and mixing for the geodesic flow.
Let $\D=(D_i)_{i\in I}$ be a locally finite $\Ga$-equivariant family
of nonempty proper closed convex subsets of $X$. Let $\Omega=
(\Omega_i)_{i\in I}$ be a locally finite $\Ga$-equivariant family of
measurable subsets of $\gengeod X$, with $\Omega_i\subset \normalpm
D_i$ for all $i\in I$ and $\Omega_i=\Omega_j$ if $i\sim_\D j$.  Assume
that $\sigma^\pm_{\Omega}$ is finite and nonzero. Then, as $t\ra +
\infty$, for the weak-star convergence of measures on $\Ga\bs\gengeod
X$,
$$
\frac{1}{\|(\flow{\pm t})_*\sigma^\pm_{\Omega}\|}\,
(\flow{\pm t})_*\sigma^\pm_{\Omega}\;\;
\stackrel{*}{\rightharpoonup}
\;\;\frac{1}{\|m_{F}\|}\,m_{F}\;.
$$
\etheo

\dem 
We only give the proof when $\pm=+$, the other case is treated
similarly. Given three numbers $a,b,c$ (depending on some parameters),
we write $a=b \pm c$ if $|a-b|\leq c$.

Let $\eta\in\;]0,1]$. We may assume that $\Ga\bs I$ is finite, since
for every $\epsilon>0$, there exists a $\Ga$-invariant partition
$I=I'\cup I''$ with $\Ga\bs I'$ finite such that if $\Omega'=
(\Omega_i)_{i\in I'}$ and $\Omega''= (\Omega_i)_{i\in I''}$, then
$\sigma^+_{\Omega}=\sigma^+_{\Omega'}+\sigma^+_{\Omega''}$ with
$\|(\flow{-t})_*\sigma^+_{\Omega''}\|=\|\sigma^+_{\Omega''}\|<\epsilon$.
Hence, using Lemma \ref{lem:defiR} (i), we may fix $R>0$ such that
$\nu^-_w (V^+_{w,\, \eta,\,R}) >0$ for all $i\in I$ and
$w\in\normalout D_i$.

Fix $\psi\in\C_{\rm c}(\Ga\bs \gengeod X)$, a continuous function with
compact support on $\Ga\bs \gengeod X$.  Let us prove that
$$
\lim_{t\ra+\infty}\;\frac{1}{\|(\flow{t})_*\sigma^+_{\Omega}\|}\;
\int_{\Ga\bs\G X}\psi\;d(\flow{t})_*\sigma^+_{\Omega}=\frac{1}{\|m_{F}\|}\;
\int_{\Ga\bs\G X}\psi\;dm_{F}\;.
$$ 
Let $\Delta_\Ga$ be a fundamental domain for the action of $\Ga$ on
$\gengeod X$, such that the boundary of $\Delta_\Ga$ has zero measure.
By a standard argument of finite partition of unity and up to
modifying $\Delta_\Ga$, we may assume that there exists a function
$\wt\psi: \gengeod X\ra\RR$ whose support has a small neighbourhood
contained in $\Delta_\Ga$ such that $\wt \psi=\psi\circ p$ on this
neighbourhood, where $p:\gengeod X\ra \Ga\bs\gengeod X$ is the
canonical projection (which is Lipschitz). Fix $\epsilon>0$. Since
$\wt\psi$ is uniformly continuous, for every $\eta>0$ small enough and
for every $t\geq 0$ large enough, for all $w\in\G_+ X$ isometric on
$[-t,+\infty[$ and $\ell\in V^+_{w,\,\eta,\,e^{-t}R}$, we have
\begin{equation}\label{eq:unifcont}
\wt\psi(\ell)= \wt\psi(w)\pm \frac{\epsilon}{2}\;.
\end{equation}

\noindent If $t$ is large enough and $\eta$ small enough, we have, using
respectively
\begin{itemize}
\item the definition of the global test function $\wt \Phi_\eta=\wt \Phi^-_\eta$,
  since the support of $\wt \psi$ is contained in $\Delta_\Ga$ and the
  support of $\phi^-_{\eta,\,R,\,\Omega_i}$ is contained in
  $\U^+_{D_i}$, for the second equality,
\item the disintegration property of $f^+_{D_i}$ in Proposition
  \ref{prop:disintegration} for the third
  equality,
\item the fact that if $\ell$ is in the support of $\nu^-_{\rho}$,
  then $f^+_{D_i}(\flow{-t}\ell)=f^+_{D_i}(\ell)=\rho$ and the change
  of variables by the geodesic flow $w=\flow{t}\rho$ for the fourth
  equality,
\item the fact that the support of $\nu^-_{\flow{-t}w}$ is contained
  in $\ws(\flow{-t} w)$, and that 
$$
\ws(\flow{-t} w)\cap \flow{t}V^+_{\eta,\, R}(\Omega_i)=
\flow{t}\big(\ws(\flow{-t} w)\cap V^+_{\eta,\, R}(\Omega_i))=
\flow{t}V^+_{\flow{-t}w,\,\eta,\,R}
$$ 
  for the fifth equality,
\item Equation \eqref{eq:unifcont} for the sixth equality, and
\item the definition of $h^-$, the invariance of the measure
  $\nu^-_{\flow{-t}w}$ and the Gibbs measure $\wt m_F$ under the
  geodesic flow, and the definition of the measure $\sigma^+_{\Omega}$
  for the last two equalities:
\end{itemize}
\begin{align}
\;&\int_{\Ga\bs\G X} \psi\;\;\Phi_{\eta}\circ\flow{-t}\;dm_{F}=
\int_{\Delta_\Ga\cap\, \G X}
\wt\psi\;\;\wt\Phi_{\eta}\circ\flow{-t}\; d\wt m_{F}\nonumber\\
=\;&
\sum_{i\in I/\sim}\;
\int_{\ell\in \U^+_{D_i}}
\;\wt\psi(\ell)\;\phi^-_{\eta,\,R,\,\Omega_i}(\flow{-t}\ell)\;
d\wt m_{F}(\ell) \nonumber\\
=\;&\sum_{i\in I/\sim}\;
\int_{\rho\in \normalout D_i}
\int_{\ell\in \U^+_{D_i}}
\wt\psi(\ell)\;h^-_{\eta,\,R}(f^+_{D_i}(\flow{-t}\ell))
\mathbbm{1}_{\V^+_{\eta,R}(\Omega_i)}(\flow{-t}\ell)\;
d\nu^{-}_{\rho}(\ell)\;d\wt\sigma^+_{D_i}(\rho)\nonumber\\
=\;&\sum_{i\in I/\sim}\;
\int_{w\in \flow{t}\normalout D_i}
\int_{\ell\in \flow{t}\V^+_{\eta,R}(\Omega_i)}
\wt\psi(\ell)\;h^-_{\eta,\,R}(\flow{-t}w))\;
d\nu^{-}_{\flow{-t}w}(\ell)\;d(\flow t)_*\wt\sigma^+_{D_i}(w)\nonumber\\
=\;&\sum_{i\in I/\sim}\;
\int_{w\in \flow{t}\normalout D_i}
\int_{\ell\in \flow{t}V^+_{\flow{-t}w,\,\eta,\,R}}
\wt\psi(\ell)\;h^-_{\eta,\,R}(\flow{-t}w)\;
d\nu^{-}_{\flow{-t}w}(\ell)\;d(\flow t)_*\wt\sigma^+_{D_i}(w)\nonumber\\
=\;&
\sum_{i\in I/\sim}\;
\int_{w\in \flow{t}\normalout D_i} \wt\psi(w)\;h^-_{\eta,\,R}(\flow{-t}w)\;
\nu^{-}_{\flow{-t}w}(\flow{t}V^+_{\flow{-t}w,\,\eta,\,R})\;
d(\flow t)_*\wt\sigma^+_{D_i}(w)\nonumber\\
& \;\;\pm\frac\epsilon 2
\int_{\Delta_\Ga\cap\, \G X}\wt\Phi_\eta\circ\flow{-t}\,d\wt m_{F}\nonumber\\
=\;&
\sum_{i\in I/\sim}\;\int_{\gengeod X}\wt\psi\;d(\flow t)_*\wt\sigma^+_{D_i}
\pm\frac\epsilon 2\int_{\Ga\bs\G X}\Phi_\eta\circ\flow{-t}\;dm_{F}\nonumber\\
=\;&
\int_{\Ga\bs\gengeod X}\psi\;d(\flow t)_*\sigma^+_{\Omega}
\pm\frac\epsilon 2\int_{\Ga\bs\G X}\Phi_\eta\;dm_{F}\;.\label{ninelineeq}
\end{align}
We then conclude as in the end of the proof of
\cite[Theo.~19]{ParPau14ETDS}. By Equation \eqref{eq:inttestskinmass}, we
have $\|(\flow{t})_*\sigma^+_{\Omega}\|=\|\sigma^+_{\Omega}\|=
\int_{\Ga\bs\G X} \Phi_\eta\; d m_{F}$. By the mixing property of the
geodesic flow on $\Ga\bs\G X$ for the Gibbs measure $m_{F}$, for
$t\geq 0$ large enough (while $\eta$ is small but fixed), we hence
have
$$
\frac{\int_{\Ga\bs\gengeod X} \psi\;d(\flow{t})_*\sigma^+_{\Omega}}
{\|(\flow{t})_*\sigma^+_{\Omega}\|} =
\frac{\int_{\Ga\bs\G X}
\Phi_{\eta}\circ\flow{-t}\;\psi\;dm_{F}}
{\int_{\Ga\bs\G X} \Phi_\eta\;dm_{F}} \pm\frac{\epsilon}{2}=
\frac{\int_{\Ga\bs\G X}\psi\;dm_{F}}
{\|m_{F}\|} \pm\epsilon\,.
$$
This proves the result.
\cqfd

\medskip
Recall that by Proposition \ref{prop:exempHCprop}, Theorem
\ref{theo:equid} applies to Riemannian manifolds with pinched negative
curvature and for $\RR$-trees for which the geodesic flow is mixing
and which satisfy the finiteness requirements of the Theorem.

Since pushforwards of measures are weak-star continuous and preserve
total mass, we have, under the assumptions of Theorem
\ref{theo:equid}, the following equidistribution result in $X$ of the
immersed $t$-neighbourhood of a properly immersed nonempty proper
closed convex subset of $X$: as $t\ra+\infty$,
\begin{equation}\label{eq:equidisequidisbas}
\frac{1}{\|\sigma^+_{\Omega}\|}\,
\pi_*(\flow{t})_*\sigma^+_{\Omega}\;\;
\stackrel{*}{\rightharpoonup}
\;\;\frac{1}{\|m_{F}\|}\,\pi_*m_{F}\;.
\end{equation}

\section{Rate of equidistribution of equidistant level 
sets for manifolds}
\label{subsect:expmixingmanifolds}

If $X=\wt M$ is a Riemannian manifold and if the geodesic flow of
$\Ga\bs\wt M$ is mixing with exponentially decaying correlations, we
get a version of Theorem \ref{theo:equid} with error bounds.  See
Section \ref{subsec:mixingratemanifolds} for conditions on $\Ga$ and
$\wt F$ that imply the exponential mixing.

\btheo\label{theo:expratesequid} Let $\wt M$ be a complete simply
connected Riemannian manifold with negative sectional curvature. Let
$\Ga$ be a nonelementary discrete group of isometries of $\wt M$.  Let
$\wt F:T^1\wt M\to\RR $ be a bounded $\Ga$-invariant
H\"older-continuous function with critical exponent $\delta=
\delta_{\Ga,\,F}$.  Let $\D=(D_i)_{i\in I}$ be a locally finite
$\Ga$-equivariant family of nonempty proper closed convex subsets of
$\wt M$, with finite nonzero skinning measure $\sigma_{\D}$.  Let
$M=\Ga\backslash\wt M$ and let $F:T^1M\ra \RR$ be the potential
induced by $\wt F$.

\smallskip\noindent (i) If $M$ is compact and if the geodesic flow on
$T^1M$ is mixing with exponential speed for the H\"older regularity
for the potential $F$, then there exist $\alpha\in\;]0,1]$ and
$\kappa''>0$ such that for all $\psi\in \C_{\rm c}^\alpha(T^1M)$, we
have, as $t\ra+\infty$,
$$
\frac{1}{\|\sigma_{\D}\|}\int \psi \;d(\flow t)_*\sigma_{\D}=
\frac{1}{\|m_F\|} \int \psi\;dm_F
+\operatorname{O}(e^{-\kappa''t}\;\|\psi\|_\alpha)\;.
$$

\smallskip\noindent(ii) If $\wt M$ is a symmetric space, if $D_i$ has
smooth boundary for every $i\in I$, if $m_F$ is finite and smooth, and
if the geodesic flow on $T^1M$ is mixing with exponential speed for
the Sobolev regularity for the potential $F$, then there exist
$\ell\in\NN$ and $\kappa''>0$ such that for all $\psi\in
\C_{\rm c}^\ell(T^1M)$, we have, as $t\ra+\infty$,
$$
\frac{1}{\|\sigma_{\D}\|}\int \psi \;d(\flow t)_*\sigma_{\D}=
\frac{1}{\|m_F\|} \int \psi\;dm_F
+\operatorname{O}(e^{-\kappa''t}\;\|\psi\|_\ell)\;.
$$
\etheo

Note that if $\wt M$ is a symmetric space and $M$ has finite volume,
then $M$ is geometrically finite.  Theorem \ref{theo:DOPB} implies
that $m_F$ is finite if $F$ is small enough. The maps
$\operatorname{O}(\cdot)$ depend on $\wt M,\Ga,F,\D,$ and the speeds
of mixing.

\medskip \dem Up to rescaling, we may assume under the assumptions of Claim
(i) or (ii) that the sectional curvature is bounded from above by
$-1$. The critical exponent $\delta$ and the Gibbs measure $m_F$ are
finite under the assumptions of the theorem.


Let us consider Claim (i).  Under its assumptions, there exists
$\alpha\in\;]0,1[\,$ such that the geodesic flow on $T^1M$ is
exponentially mixing for the H\"older regularity $\alpha$ and such
that the strong stable foliation of $T^1\wt M$ is
$\alpha$-H\"older.\footnote{See Section
  \ref{subsec:mixingratemanifolds}.}

First assume that $\Ga\bs I$ is finite. Fix $R>0$ large enough and,
for every $\eta>0$, let us consider the test function $\Phi_\eta$ as
in the proof of Theorem \ref{theo:equid}. Up to replacing $D_i$ by
$\N_1D_i$, we may assume that the boundary of $D_i$ is
$\C^{1,1}$-smooth, for every $i\in I$, see for instance
\cite{Walter76}.

Fix $\psi\in \C_{\rm c}^\alpha(T^1M)$. We may assume as in the proof
of Theorem \ref{theo:equid} that there exists a lift $\wt\psi: T^1\wt
M\ra\RR$ of $\psi$ whose support is contained in a given fundamental
domain $\Delta_\Ga$ for the action of $\Ga$ on $T^1\wt M$. There exist
$\eta_0>0$ and $t_0\geq 0$ such that for every $\eta\in\;]0,\eta_0]$,
and for every $t\in [t_0,+\infty[$, for every $w\in T^1\wt M$ and
    $v\in V^+_{w,\,\eta,\,e^{-t}R}$, we have
\begin{equation}\label{eq:holdercont}
\wt\psi(v)= \wt\psi(w)+
\operatorname{O}\big((\eta+e^{-t})^\alpha\|\psi\|_\alpha\big)\;, 
\end{equation}
since $d(v,w)=\operatorname{O}(\eta+e^{-t})$ by Equation
\eqref{eq:disttranslatgeod} and Lemma \ref{lem:comparddHam}.

Let $\overline m_F=\frac{m_F}{\|m_F\|}$ be the normalisation of the
Gibbs measure $m_F$ to a probability measure. As in the proof of
Theorem \ref{theo:equid} using Equation \eqref{eq:holdercont} instead
of Equation \eqref{eq:unifcont}, we have
$$
\frac{\int_{T^1M} \psi\;d(\flow{t})_*\sigma^+_{\Omega}}
{\|(\flow{t})_*\sigma^+_{\Omega}\|}=\frac{\int_{T^1M}
\Phi_\eta\;\psi\circ\flow t\;d\overline{m}_{F}}
{\int_{T^1M}\Phi_\eta \;d\overline{m}_{F}}
+\operatorname{O}\big((\eta+e^{-t})^\alpha\|\psi\|_\alpha\big)\;.
$$ 
As $M$ is compact, the Patterson densities and the Bowen-Margulis
measure are doubling measures by Lemma \ref{lem:shadowlemma}
(4).\footnote{See also \cite[Prop.~3.12]{PauPolSha15}.} Using discrete
convolution approximation,\footnote{See for instance
  \cite[p. 290-292]{Semmes96} or \cite{KinKorShaTuo12}.} there exist
$\kappa'>0$ and, for every $\eta>0$, a nonnegative function
$\reg\Phi_\eta \in \operatorname{C}_{\rm c}^\alpha(T^1M)$ such that
\begin{enumerate}
\item[$\bullet$]
  $\int_{T^1M}\reg\Phi_\eta\;d\overline{m}_{F}=\int_{T^1M}\Phi_\eta
  \;d\overline{m}_{F}$,
\item[$\bullet$] $\int_{T^1M}|\reg\Phi_\eta-\Phi_\eta|
  \;d\overline{m}_{F}= \operatorname{O}(\eta\int_{T^1M}\Phi_\eta
  \;d\overline{m}_{F})$,
\item[$\bullet$] $\|\reg\Phi_\eta\|_\alpha =
  \operatorname{O}(\eta^{-\kappa'} \int_{T^1M}\Phi_\eta
  \;d\overline{m}_{F})$.
\end{enumerate}
Hence, applying the exponential mixing of the geodesic flow, with
$\kappa>0$ as in its definition \eqref{eq:defiholderexponentialdecay},
we have, for $\eta\in\;]0,\eta_0]$ and $t\in [t_0,+\infty[$,
\begin{align*}
&\frac{\int_{T^1M} \psi\;d(\flow{t})_*\sigma^+_{\Omega}}
{\|(\flow{t})_*\sigma^+_{\Omega}\|}=\\=&
\;\frac{\int_{T^1M}
\reg\Phi_\eta\;\psi\circ\flow t\;d\overline{m}_{F}}{\int_{T^1M}\Phi_\eta
\;d\overline{m}_{F}}
+ \operatorname{O}\big(\eta\; \|\psi\|_\infty+ 
(\eta+e^{-t})^\alpha\|\psi\|_\alpha\big)\\=\; & 
\frac{\int_{T^1M}  \reg\Phi_\eta\;d\overline{m}_{F}}{\int_{T^1M}\Phi_\eta
  \;d\overline{m}_{F}}
\;\int_{T^1M}  \psi\;d\overline{m}_{F}
+ \operatorname{O}\big(e^{-\kappa t}\|\reg\Phi_\eta\|_\alpha\|\psi\|_\alpha+ 
\eta\; \|\psi\|_\infty+(\eta+e^{-t})^\alpha\|\psi\|_\alpha\big)\\=\; & 
\int_{T^1M}  \psi\;d\overline{m}_{F}
+ \operatorname{O}\big((e^{-\kappa t}\eta^{-\kappa'}+ 
\eta+(\eta+e^{-t})^\alpha)\|\psi\|_\alpha\big)\;.
\end{align*}
Taking $\eta=e^{-t\lambda}$ for $\lambda$ small enough (for instance
$\lambda =\kappa/(2\kappa')$), the result follows (for instance with
$\kappa''= \min\{\kappa/2,\; \kappa/(2\kappa'), \;\alpha \min\{1,
\kappa/(2\kappa')\}\}$), when $\Ga\bs I$ is finite.  As the implied
constants do not depend on the family $\D$, the result holds in
general.

For Claim (ii), the required smoothness of $m_F$ (that is, the fact
that $m_F$ is absolutely continuous with respect to the Lebesgue
measure with smooth Radon-Nikodym derivative) allows to use the
standard convolution approximation described for instance in \cite[\S
  1.6]{Ziemer89}, instead of the operator $\reg$ as above, and the
proof proceeds similarly. \cqfd

\section{Equidistribution  of  equidistant level 
sets on simplicial \\ graphs and random walks on 
graphs of groups}
\label{subsec:randomwalk}

Let $\XX$, $X$, $\Ga$, $\wt c$, $c$, $\wt F_c$, $F_c$,
$\delta_c<+\infty$, $(\mu^\pm_x)_{x\in V\XX}$, $\wt m_c=\wt m_{F_c}$,
$m_c=m_{F_c}$ be as in the beginning of Section
\ref{subsec:mixingratesimpgraphs}.

In this Section, we state an equidistribution result analogous to
Theorem \ref{theo:equid}, which now holds in the space of generalised
discrete geodesic lines $\Ga\bs\gengeod \XX$, but whose proof is
completely analogous.

\btheo\label{theo:equidsimplicial} Let $\XX,\Ga,\wt c,
(\mu^\pm_x)_{x\in V\XX}$ be as above. Assume that the Gibbs measure
$m_{c}$ on $\Ga\bs\G \XX$ is finite and mixing for the discrete time
geodesic flow.  Let $\D=(\DD_i)_{i\in I}$ be a locally finite
$\Ga$-equivariant family of nonempty proper simplicial subtrees of
$\XX$, and $D_i=|\DD_i|_1$. Let $\Omega= (\Omega_i)_{i\in I}$ be a
locally finite $\Ga$-equivariant family of measurable subsets of
$\gengeod \XX$, with $\Omega_i\subset \normalout D_i$ for all $i\in I$
and $\Omega_i=\Omega_j$ if $i\sim_\D j$.  Assume that
$\sigma^+_{\Omega}$ is finite and nonzero. Then, as $n\ra+\infty$, for
the weak-star convergence of measures on $\Ga\bs\gengeod \XX$,
$$
\frac{1}{\|(\flow{n})_*\sigma^+_{\Omega}\|}\,
(\flow{n})_*\sigma^+_{\Omega}\;\;
\stackrel{*}{\rightharpoonup}
\;\;\frac{1}{\|m_{F}\|}\,m_{F}\;. \;\;\;\Box
$$
\etheo

We leave to the reader the analog of this result when the
restriction to $\Ga\bs \Geven\XX$ of the Gibbs measure is finite and
mixing for the square of the discrete time geodesic flow.

Using Proposition \ref{prop:uniflatmBMfinie} and Theorem
\ref{theo:uniflatmBMmixing} in order to check that the Bowen-Margulis
measure $m_{\rm BM}$ on $\Ga\bs\G \XX$ is finite and mixing, we have
the following consequence of Theorem \ref{theo:equidsimplicial}, using
the system of conductances $\wt c=0$.

\bcoro\label{coro:equidsimpliciallattice} Let $\XX$ be a uniform
simplicial tree. Let $\Ga$ be a lattice of $\XX$ such that the graph
$\Ga\bs\XX$ is not bipartite.  Let $\D=(\DD_i)_{i\in I}$ be a locally
finite $\Ga$-equivariant family of nonempty proper simplicial subtrees
of $\XX$ and $D_i=|\DD_i|_1$. Let $\Omega= (\Omega_i)_{i\in I}$ be a
locally finite $\Ga$-equivariant family of measurable subsets of
$\gengeod \XX$, with $\Omega_i\subset \normalout D_i$ for all $i\in I$
and $\Omega_i=\Omega_j$ if $i\sim_\D j$.  Assume that the skinning
measure $\sigma^+_{\Omega}$ (with vanishing potential) is finite and
nonzero. Then, as $n\ra+\infty$, for the weak-star convergence of
measures on $\Ga\bs\gengeod \XX$,
$$
\frac{1}{\|(\flow{n})_*\sigma^+_{\Omega}\|}\,
(\flow{n})_*\sigma^+_{\Omega}\;\;
\stackrel{*}{\rightharpoonup}
\;\;\frac{1}{\|m_{\rm BM}\|}\,m_{\rm BM}\;. \;\;\;\Box
$$
\ecoro

When furthermore $\XX$ is regular, we have the following corollary,
using Proposition \ref{prop:computBM} (3).

\bcoro\label{coro:equidsimplicialregular} Let $\XX$ be a regular
simplicial tree of degree at least $3$. Let $\Ga$ be a lattice of
$\XX$ such that the graph $\Ga\bs\XX$ is not bipartite. Let
$\D=(\DD_i)_{i\in I}$ be a locally finite $\Ga$-equivariant family of
nonempty proper simplicial subtrees of $\XX$ and $D_i=|\DD_i|_1$. Let
$\Omega= (\Omega_i)_{i\in I}$ be a locally finite $\Ga$-equivariant
family of measurable subsets of $\gengeod \XX$, with $\Omega_i\subset
\normalout D_i$ for all $i\in I$ and $\Omega_i=\Omega_j$ if $i\sim_\D
j$.  Assume that the skinning measure $\sigma^+_{\Omega}$ (with
vanishing potential) is finite and nonzero. Then, as $n\ra+\infty$,
for the weak-star convergence of measures on $\Ga\bs V\XX$,
$$
\frac{1}{\|(\flow{n})_*\sigma^+_{\Omega}\|}\,
\pi_*(\flow{n})_*\sigma^+_{\Omega}\;\;
\stackrel{*}{\rightharpoonup}
\;\;\frac{1}{\Vol(\Ga \dbs \XX)}\,\vol_{\Ga\dbs\XX}\;. \;\;\;\Box
$$
\ecoro

\medskip Let us give an application of Corollary
\ref{coro:equidsimplicialregular} in terms of random walks on graphs
of groups, which might also be deduced from general result on random
walks, as indicated to the third author by M.~Burger and S.~Mozes.

Let $(\YY,G_*)$ be a graph of finite groups with finite volume, and
let $(\YY',G'_*)$ be a connected subgraph of subgroups.\footnote{See
  Section \ref{subsec:trees} for definitions and background.} Note
that $(\YY',G'_*)$ also has finite volume, less than or equal to the
volume of $(\YY,G_*)$. We say that $(\YY,G_*)$ is {\em locally
  homogeneous}\index{homogeneous!locally} if $\sum_{e\in E\YY,\;
  o(e)=x}\frac{|G_x|}{|G_e|}$ is constant at least $3$ for all $x\in
V\YY$. We say that a graph of groups is {\it
  $2$-acylindrical}\index{acylindrical} if the action of its
fundamental group on its Bass-Serre tree is $2$-acylindrical (see
Remark \ref{rem:notMarkovgood}). In particular, this action is
faithful if the graph has at least two edges.

The {\em non-backtracking simple random walk on $(\YY,G_*)$ starting
  transversally to $(\YY',G'_*)$}\index{non-backtracking}%
\index{random walk!non-backtracking} is the following Markovian random
process $(X_n=(f_n,\ga_n))_{n\in\NN}$ where $f_n\in E\YY$ and $\ga_n$
is a double coset or right coset of $G_{o(f_n)}$ for all
$n\in\NN$. Choose at random a vertex $y_0$ of $\YY'$ for the
probability measure $\frac{1}{\Vol(\YY',G'_*)} \; \vol_{\YY',G'_*}$
(we will call $y_0$ the {\it origin}\index{origin} of the random
path). Then choose uniformly at random $X_0= (f_0, \ga_0)$ where
$f_0\in E\YY$ is such that $o(f_0)=y_0$ and $\ga_0$ is a double coset
in $G'_{y_0}\bs G_{y_0}/ \rho_{\,\overline{f_0}}(G_{f_0})$ such that
if $f_0\in E\YY'$ then $\ga_0\notin G'_{y_0}\,
\rho_{\,\overline{f_0}}(G_{f_0})$.\footnote{This last condition says
  that $\ga_0$ is not the double coset of the trivial element.}
Assuming $X_n=(f_n,\ga_n)$ constructed, choose uniformly at random
$X_{n+1}=(f_{n+1},\ga_{n+1})$ where $f_{n+1} \in E\YY$ is such that
$o(f_{n+1})=t(f_n)$ and $\ga_{n+1}\in G_{o(f_{n+1})}/
\rho_{\,\overline{f_{n+1}}}(G_{f_{n+1}})$ is such that if $f_{n+1}=
\overline{f_n}$ then $\ga_{n+1}\notin
\rho_{\,\overline{f_{n+1}}}(G_{f_{n+1}})$.  The {\it $n$-th
  vertex}\index{vertex@$n$-th vertex of a random walk} of
$(X_n=(f_n,\ga_n))_{n\in\NN}$ is $o(f_n)$.

\bcoro 
Let $(\YY,G_*)$ be a locally homogeneous $2$-acylindrical
nonbipartite graph of finite groups with finite volume, and
let $(\YY',G'_*)$ be a locally homogeneous nonempty proper connected
subgraph of subgroups. Then the $n$-th vertex of the non-backtracking
simple random walk on $(\YY,G_*)$ starting transversally to
$(\YY',G'_*)$ converges in distribution to $\frac{1}{\Vol(\YY,G_*)}\;
\vol_{\YY,G_*}$ as $n\ra +\infty$.  
\ecoro

\dem Let $\Ga$ be the fundamental group of $(\YY,G_*)$ (with respect
to a choice of basepoint in $V\YY'$), which is a lattice of the
Bass-Serre tree $\XX$ of $(\YY,G_*)$, since $\Ga$ acts faithfully on
$\XX$ and $(\YY,G_*)$ has finite volume. Note that $\XX$ is regular
since $(\YY,G_*)$ is locally homogeneous.  Let $p:\XX\ra\YY=\Ga\bs\XX$
be the canonical projection.

Let $\Ga'$ be the fundamental group of $(\YY',G_*')$ (with respect to
the same choice of basepoint). As seen in Section \ref{subsec:trees},
there exists a simplicial subtree $\XX'$ whose stabiliser in $\Ga$ is
$\Ga'$, such that the quotient graph of groups $\Ga'\dbs\XX'$
identifies with $(\YY',G_*')$ and the map $(\Ga'\bs \XX')\ra
(\Ga\bs\XX)$ is injective. Similarly, $\XX'$ is regular since
$(\YY',G_*')$ is locally homogeneous. Let $\D=(\ga\XX')_{\ga\in\Ga}$,
which is a locally finite $\Ga$-equivariant family of nonempty proper
simplicial subtrees of $\XX$.

Using the notation of Example \ref{exem:quotgraphgroup} for the graph
of groups $\Ga\dbs\XX$ (which identifies with $(\YY,G_*)$), we fix
lifts $\wt f$ and $\wt y$ in $\XX$ by $p$ of every edge $f$ and vertex
$y$ of $\YY$ such that $\overline{\wt f}= \wt{\overline{f}}$, and
elements $g_f\in\Ga$ such that $g_f\,\wt{t(f)} =t(\wt f\,)$. We may
assume that $\wt f\in E\XX'$ if $f\in E\YY'$, that $\wt y\in V\XX'$ if
$y\in V\YY'$, and that $g_f\in\Ga'$ if $f\in E\YY'$, which is possible
by Equation \eqref{eq:simplesubtree}.

Let $(\Omega,\PP)$ be the (canonically constructed) probability space
of the random walk $(X_n= (f_n,\ga_n))_{n\in\NN}$. For all $n\in\NN$,
let $y_n=o(f_n)$ be the random variable (with values in the discrete
space $\YY=\Ga\bs\XX$) of the $n$-th vertex of the random walk
$(X_n)_{n\in\NN}$.

Let us define a measurable map $\Psi: \Omega\ra \Ga\bs\gengeod \XX$,
with image contained in the image of $\normalout\XX'$ by the canonical
projection $\gengeod \XX\ra \Ga\bs\gengeod \XX$, such that
$\Psi_*\PP$ is the normalised skinning measure
$\frac{\sigma^+_\D}{\|\sigma^+_\D\|}$ and that the following diagram
commutes for all $n\in\NN$~:
\begin{equation}\label{eq:commutrandomwalk}
\xymatrix{
\Omega\ar[rr]^{\Psi}\ar[dr]_{y_n}& &\;\Ga\bs\gengeod\XX\ar[dl]^{\pi\circ\flow{n}}
\\ & \;\;\YY\; &.
         }
\end{equation}
Assuming that we have such a map, we have  
$$
(y_n)_*\PP \;=\; (\pi_*\circ (\flow{n})_*\circ \Psi_*)\PP\;=\;
\pi_* (\flow{n})_* \;\frac{\sigma^+_\D}{\|\sigma^+_\D\|}\;=\;
\frac{1}{\|(\flow{n})_*\sigma^+_\D\|}\;\pi_*(\flow{n})_* \sigma^+_\D
$$ 
so that the convergence of the law of $y_n$ to
$\frac{1}{\Vol(\YY,G_*)}\; \vol_{\YY,G_*}$ follows from Corollary
\ref{coro:equidsimplicialregular} applied to $\Omega=(\normalout
D_i)_{i\in I}$.

\begin{center}
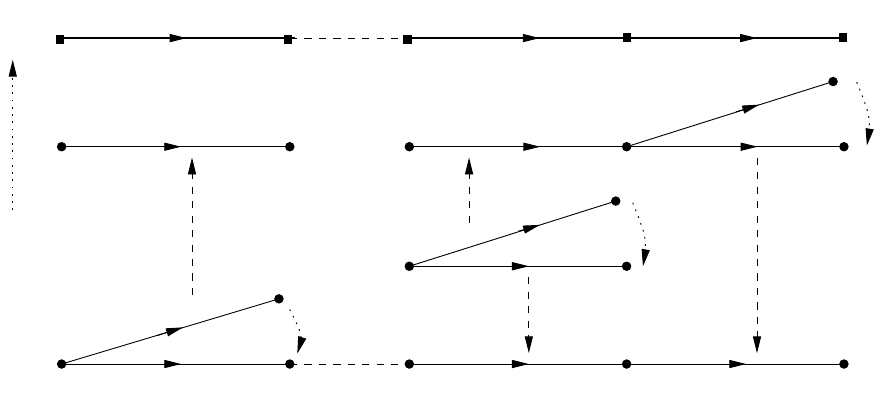
\end{center}

Let $(X_n= (f_n,\ga_n))_{n\in\NN}$ be a random path with origin
$y_0\in\YY'$, corresponding to $\omega\in\Omega$. Fix a representative
of $\ga_n$ in its right class for every $n\geq 1$, and a
representative of $\ga_0$ in its double class, that we still denote by
$\ga_n$ and $\ga_0$ respectively. Using ideas introduced for the
coding in Section \ref{subsec:codagesimplicial}, let us construct by
induction an infinite geodesic edge path $(e_n)_{n\in\NN}$ with origin
$o(e_0)=\wt{y_0}$ and a sequence $(\alpha_n)_{n\in\NN}$ in $\Ga$ such
that
\begin{equation}\label{eq:inducrandwalk}
e_n=\alpha_n\ga_n g_{\,\overline{f_n}}^{\;\;\;\;-1}\,\wt{f_n} \;.
\end{equation}

Let $\alpha_0=\id$ and $e_0=\ga_0\, g_{\,\overline{f_0}}^{\;\;\;\;-1}\; 
\wt{f_0}$. Since $o(\wt{f_0})=\, g_{\,\overline{f_0}} \;\wt{o(f_0)} =
\, g_{\,\overline{f_0}} \;\wt{y_0}$ by the construction of the lifts
and since $\ga_0\in G_{y_0}= \Ga_{\wt{y_0}}$, we have $o(e_0)=
\wt{y_0}$. Since the stabiliser of $g_{\,\overline{f_{0}}}^{\;\;\;\;-1} 
\;\wt{f_{0}}$ is $\rho_{\,\overline{f_{0}}}(G_{f_{0}})$, the edge $e_0$ 
does not depend on the choice of the representative $\ga_{0}$ modulo
$\rho_{\,\overline{f_{0}}}(G_{f_{0}})$ on the right, but depends on
the choice of the representative $\ga_{0}$ modulo $G'_{y_0}=
\Ga'_{\wt{y_0}}$ on the left.

The hypothesis that if $f_0\in E\YY'$ then $\ga_0\notin G'_{y_0}\,
\rho_{\,\overline{f_0}}(G_{f_0})$ ensures that the edge $e_0$ does not
belong to $E\XX'$. Indeed, assume otherwise that $e_0$ belongs to
$E\XX'$.  Then $f_0=p(e_0)\in E\YY'$, and by the assumptions on the
choice of lifts, the edges $g_{\,\overline{f_0}}^{\;\;\;\;-1}
\,\wt{f_0}$ and $e_0$ both belong to $E\XX'$. Since they are both
mapped to $f_0$ by the map $\XX'\ra\YY'=\Ga'\bs\XX'$, and by Equation
\eqref{eq:simplesubtree}, they are mapped one to the other by an
element of $\Ga'_{\,\wt{y_0}}= G'_{y_0}$. Let $\ga'_0\in \Ga'_{\,\wt{y_0}}$ 
be such that $\ga'_0\,e_0=g_{\,\overline{f_0}}^{\;\;\;\;-1}\,\wt{f_0}$.  
Then ${\ga'_0}^{-1} \ga_0$ belongs to the stabiliser in $\Ga$ of the
edge $g_{\,\overline{f_0}}^{\;\;\;\;-1} \,\wt{f_0}$, which is equal to
$g_{\,\overline{f_0}}^{\;\;\;\;-1}\,\Ga_{\wt{f_0}} g_{\,\overline{f_0}} = 
\rho_{\,\overline{f_0}}(G_{f_0})$. Therefore $\ga_0 \in G'_{y_0}\,
\rho_{\,\overline{f_0}}(G_{f_0})$, a contradiction.

Assume by induction that $e_n$ and $\alpha_n$ are constructed.
Define
$$ 
\alpha_{n+1}=\alpha_n\,\ga_n\,g_{\,\overline{f_n}}^{\;\;\;\;-1}\,g_{f_n}
$$
and
$$
e_{n+1}=\alpha_{n+1}\,\ga_{n+1}\;
g_{\,\overline{f_{n+1}}}^{\;\;\;\;\;\;-1}\;\wt{f_{n+1}}\;,
$$ 
so that the induction formula \eqref{eq:inducrandwalk} at rank
$n+1$ is satisfied.  By the construction of the lifts, since $y_{n+1}
=t(f_n)=o(f_{n+1})$, we have
$$
\wt{y_{n+1}}=g_{f_n}^{\;\;-1} \;t(\wt{f_n})=
g_{\,\overline{f_{n+1}}}^{\;\;\;\;\;\;-1} \;o(\wt{f_{n+1}})\;.
$$ 
Hence, since $\ga_{n+1}\in G_{y_{n+1}}=\Ga_{\wt{y_{n+1}}}$ fixes
$\wt{y_{n+1}}$, using the induction formula \eqref{eq:inducrandwalk}
at rank $n$ for the last equality,
\begin{align*}
o(e_{n+1})&=\alpha_{n+1}\,\ga_{n+1}\;
g_{\,\overline{f_{n+1}}}^{\;\;\;\;\;\;-1}\;o(\wt{f_{n+1}})=
\alpha_{n+1}\,\ga_{n+1}\;\wt{y_{n+1}}= \alpha_{n+1}\;\wt{y_{n+1}}
\\ & = \alpha_{n+1}\;g_{f_n}^{\;-1}\; t(\wt{f_n})
=\alpha_n\,\ga_n\,g_{\,\overline{f_n}}^{\;\;\;\;-1}\; t(\wt{f_n})
=t(e_n)\;.
\end{align*}
In particular, the sequence $(e_n)_{n\in\NN}$ is an edge path in
$\XX$.

Since the stabiliser of $g_{\,\overline{f_{n+1}}}^{\;\;\;\;\;\;-1} \;
\wt{f_{n+1}}$ is $$g_{\,\overline{f_{n+1}}}^{\;\;\;\;\;\;-1} \;
\Ga_{\wt{f_{n+1}}} g_{\,\overline{f_{n+1}}}= \rho_{\,\overline{f_{n+1}}} 
(G_{f_{n+1}})\;,$$ the edge $\ga_{n+1} \;
g_{\,\overline{f_{n+1}}}^{\;\;\;\;\;\;-1} \; \wt{f_{n+1}}$ does not
depend on the choice of the representative of the right coset
$\ga_{n+1}$. Let us prove that the length $2$ edge path
$(g_{f_n}^{\;-1}\;\wt{f_n},\; \ga_{n+1}\; g_{\,\overline{f_{n+1}}}
^{\;\;\;\;\;\;-1}\; \wt{f_{n+1}})$ is geodesic. Otherwise, the two
edges of this path are opposite one to another, hence $f_{n+1}=
\overline{f_n}$ by using the projection $p:\XX\ra\YY$, therefore
$g_{\,\overline{f_{n+1}}}= g_{f_n}$. Thus $\ga_{n+1}$ maps
$g_{f_n}^{\;-1}\; \wt{f_n}$ to $g_{f_n}^{\;-1}\;\overline{\wt{f_n}}$,
hence belongs to $\rho_{f_{n}}(G_{f_{n}})= \rho_{\,\overline{f_{n+1}}}
(G_{f_{n+1}})$, a contradiction by the assumptions on the random walk.

By construction, the element $\alpha_{n+1}$ of $\Ga$ sends the above
length $2$ geodesic edge path $(g_{f_n}^{\;-1}\;\wt{f_n},\; \ga_{n+1}\;
g_{\,\overline{f_{n+1}}}^{\;\;\;\;\;\;-1}\; \wt{f_{n+1}})$ to
$(e_n,e_{n+1})$. This implies on the one hand that the edge path
$(e_n,e_{n+1})$ is geodesic, and on the other hand that $\alpha_{n+1}$
is uniquely defined, since the action of $\Ga$ on $\XX$ is
$2$-acylindrical.

In particular, $(e_n)_{n\in\NN}$ is the sequence of edges followed by
a (discrete) geodesic ray in $\XX$, starting from a point of $\XX'$
but not by an edge of $\XX'$, that is, an element of $\normalout\XX'$.
Furthermore, this ray is well defined up to the action of
$\Ga'_{\wt{y_0}}$, hence its image, that we denote by
$\Psi(\omega)$, is well defined in $\Ga\bs \gengeod\XX$. Since
$p(o(e_n))=p(\wt{y_n})=y_n$ for all $n\in\NN$, the commutativity of
the diagram \eqref{eq:commutrandomwalk} is immediate.

For every $x\in V\XX'$, let $\normalout\XX'(x)$ be the subset of
$\normalout\XX'$ consisting of the elements $w$ with $w(0)=x$. By
construction, the above map from the subset of random paths in
$\Omega$ starting from $y_0$ to $\Ga'_{\wt{y_0}}\bs
\normalout\XX'(\wt{y_0})$, which associates to $(X_n)_{n\in\NN}$ the
$\Ga'_{\wt{y_0}}$-orbit of the geodesic ray with consecutive edges
$(e_n)_{n\in\NN}$, is clearly a bijection. This bijection maps the
measure $\PP$ to the normalised skinning measure
$\frac{\sigma^+_\D}{\|\sigma^+_\D\|}$, since by homogeneity, the
restriction to $\normalout\XX'(\wt{y_0})$ of $\wt\sigma^+_{\XX'}$ to
$\normalout\XX'(\wt{y_0})$, normalised to be a probability measure, is
the restriction to $\normalout\XX'(\wt{y_0})$ of the
$\Aut(\XX)_x$-homogeneous probability measure on the space of geodesic
rays with origin $x$ in the regular tree $\XX$. This proves the
result.  \cqfd

\medskip When $\YY$ is finite, all the groups $G_y$ for $y\in V\YY$
are trivial and $\YY'$ is reduced to a vertex,\footnote{The result for
general $\YY'$ follows by averaging.} the above random walk is the
non-backtracking simple random walk on the nonbipartite regular finite
graph $\YY$, and $\frac{1}{\Vol(\YY,G_*)}\; \vol_{\YY,G_*}$ is the
uniform distribution on $V\YY$. Hence this result (stated as Corollary
\ref{coro:introrandomwalk} in the Introduction) is classical.  See for
instance \cite[Theo.~1.2]{OrtWoe07} and \cite{AloBenLubSod07}, which
under further assumptions on the spectral properties of $\YY$ gives
precise rates of convergence, and also the book \cite{LyoPerBook},
including its Section 6.3 and its references.

\section{Rate of equidistribution for metric and simplicial 
trees}
\label{subsec:rateequidtrees}

In this Section, we give error terms for the equidistribution results
stated in Theorem \ref{theo:equid} for metric trees, and in Theorem
\ref{theo:equidsimplicial} for simplicial trees, under additional
assumptions required in order to get the  error terms for the mixing
property discussed in Chapter \ref{sec:mixingrate}.

We first consider the simplicial case, for the discrete time geodesic
flow. Let $\XX$, $X$, $\Ga$, $\wt c$, $c$, $\wt F_c$, $F_c$,
$\delta_c<+\infty$, $(\mu^\pm_x)_{x\in V\XX}$, $\wt m_c$, $m_c$ be as
in the beginning of Section \ref{subsec:mixingratesimpgraphs}.

\btheo\label{theo:equidsimplicialerrorterm} Assume that $\delta_c$ is
finite and that the Gibbs measure $m_{c}$ on $\Ga\bs\G \XX$ is
finite. Assume furthermore that
\begin{enumerate}
\item the families $(\Lambda\Ga, \mu^-_x, d_x)_{x\in V\C\Lambda\Ga}$
  and $(\Lambda\Ga, \mu^+_x, d_x)_{x\in V\C\Lambda\Ga}$ of metric
  measure spaces are uniformly doubling,\footnote{See Section
    \ref{subsec:Pattersondens} for definitions.}
\item there exists $\alpha\in\;]0,1]$ such that the discrete time geodesic
  flow on $(\Ga\bs\G \XX,m_c)$ is exponentially mixing for the
  $\alpha$-H\"older regularity.
\end{enumerate}
Let $\D=(D_i)_{i\in I}$ be a locally finite $\Ga$-equivariant family
of nonempty proper simplicial subtrees of $\XX$ with $\Ga\bs I$
finite. Let $\Omega= (\Omega_i)_{i\in I}$ be a locally finite
$\Ga$-equivariant family of measurable subsets of $\gengeod \XX$, with
$\Omega_i\subset \normalpm D_i$ for all $i\in I$ and
$\Omega_i=\Omega_j$ if $i\sim_\D j$.  Assume that
$\sigma^\pm_{\Omega}$ is finite and positive. Then there exists
$\kappa'>0$ such that for all $\psi\in \C_{\rm c}^\alpha
(\Ga\bs\gengeod\XX)$, we have, as $n\ra+\infty$,
$$
\frac{1}{\|(\flow{\pm n})_*\sigma^\pm_{\Omega}\|}\,
\int \psi\;d\,\pi_*(\flow{\pm n})_*\sigma^\pm_{\Omega}\;=\;
\frac{1}{\|m_{c}\|}\,\int \psi\;d\,m_{c} +
\bigO(\|\psi\|_\alpha\,e^{-\kappa'\,n})\;.
$$
\etheo

\noindent {\bf Remarks. } (1) If $\wt c=0$, if the simplicial subtree
$\XX'$ of $\XX$ satisfying $|\XX'|_1=\C\Lambda\Ga$ is uniform, if
$L_\Ga=\ZZ$ and if $\Ga$ is a lattice of $\XX'$, then we claim that
$\delta_c=\delta_\Ga$ is finite, $m_c=m_{\rm BM}$ is finite and
mixing, and $(\Lambda\Ga,\mu_x= \mu_x^\pm, d_x)_{x\in V\C\Lambda\Ga}$
is uniformly doubling.

Indeed, the above finiteness and mixing properties follow from the
results of Section \ref{subsec:ergodictrees}. Since $\XX'$ is uniform,
it has a cocompact discrete group of isometries $\Ga'$ whose Patterson
density (for the vanishing potential) is uniformly doubling on
$\Lambda\Ga'= \Lambda\Ga$, by Lemma \ref{lem:shadowlemma} (4). Since
$\wt c=0$ and $\Ga$ is a lattice, the Patterson densities of $\Ga$ and
of $\Ga'$ coincide (up a scalar multiple) by Proposition
\ref{prop:uniflatmBMfinie} (2).

\medskip (2) Assume that $\wt c=0$, that the simplicial subtree $\XX'$
of $\XX$ satisfying $|\XX'|_1=\C\Lambda\Ga$ is uniform without
vertices of degree $2$, that $L_\Ga=\ZZ$ and that $\Ga$ is a
geometrically finite lattice of $\XX$.  Then all assumptions of
Theorem \ref{theo:equidsimplicialerrorterm} are satisfied by the first
remark and by Corollary \ref{coro:expdecaygeomfinisimpl}. Therefore we
have an exponentially small error term in the equidistribution of the
equidistant levels sets.

\bigskip
\dem We only give the proof when $\pm=+$, the other case is treated
similarly. We follow the proof of Theorem \ref{theo:equid},
concentrating on the new features. We now have $\Ga\bs I$ finite by
assumption. Let $\eta\in\;]0,1]$ and $\psi\in\C^\alpha_{\rm c}(\Ga\bs
\gengeod X)$. We consider the constant $R>0$, the test function
$\Phi_\eta$, the fundamental domain $\Delta_\Ga$ and the lift
$\wt\psi$ of $\psi$ as in the proof of Theorem \ref{theo:equid}.

For all $n\in\NN$, all $w\in\G_+ X$ isometric on $[-n,+\infty[$ and
all $\ell\in V^+_{w,\,\eta,\,e^{-n}R}= B^+(w,e^{-n}R)$,\footnote{As
  said in Section \ref{subsec:trees}, the subsets $V^+_{w,\,\eta,\,s}$
  and $B^+(w,s)$ of the space of discrete geodesic lines $\G\XX$ are
  equal for every $s>0$ since $\XX$ is simplicial and $\eta<1$.} by
Lemma \ref{lem:estimdellw} where we can take $\eta=0$, we have
$d(\ell,w)= \bigO(e^{-n})$. Since $p$ is Lipschitz, the map $\wt\psi$
is $\alpha$-H\"older-continuous with $\alpha$-H\"older norm at most
$\|\psi\|_\alpha$. Hence for all $n\in\NN$, all $w\in\G_+ X$ isometric
on $[-n,+\infty[$ and all $\ell\in V^+_{w,\,\eta,\,e^{-n}R}$, we have
\begin{equation}\label{eq:controlholderequidequid}
\wt\psi(\ell)=\wt\psi(w) +\bigO(e^{-n\,\alpha}\;\|\psi\|_\alpha)\;.
\end{equation}

As in the proof of Theorem \ref{theo:equid} with $t$ replaced by $n$,
using Equation \eqref{eq:controlholderequidequid} instead of
Equation \eqref{eq:unifcont} in the series of equations
\eqref{ninelineeq}, since the symbols $w$ that appear in them are indeed
generalised geodesic lines isometric on $[-n,+\infty[\,$, we have
\begin{equation}\label{eq:prepaapplimixingexpo}
\frac{\int_{\Ga\bs\gengeod \XX} \;\psi\;d(\flow n)_*\sigma^+_{\Omega}}
{\|(\flow{n})_*\sigma^+_{\Omega}\|} =
\frac{\int_{\Ga\bs\G \XX} \;\psi\;\;\Phi_{\eta}\circ\flow{-n}\;dm_{c}}
{\int_{\Ga\bs\G \XX} \;\Phi_{\eta}\;dm_{c}}+\bigO(e^{-n\,\alpha}\;\|\psi\|_\alpha)\;.
\end{equation}

Let us now apply the assumption on the decay of correlations. In order
to do that, we need to regularise our test functions $\Phi_{\eta}$\,. 

\medskip 
By the definition of the Gibbs measures,\footnote{See Equation
  \eqref{eq:defigibbs}, using $x_0=\ell(0)$ as the basepoint for the
  Hopf parametrisation, and the fact that if $\epsilon>0$ is small
  enough and $\ell'\in B_d(\ell,\epsilon)$, then $\ell'(0)=\ell(0)$ as
  seen in the proof of Lemma \ref{lem:troisboulessimp}.} Lemma
\ref{lem:troisboulessimp} implies that, for all $\epsilon>0$ small
enough and $\ell\in\G\XX$,
\begin{align*}
\mu^-_{\ell(0)}
\big(B_{d_{\ell(0)}}&(\ell(-\infty), \frac{1}{c_0}\; \sqrt{\epsilon}
\,)\big)\; \mu^+_{\ell(0)}\big(
B_{d_{\ell(0)}}(\ell(+ \infty), \frac{1}{c_0}\; \sqrt{\epsilon}\,)\big)\\
\leq\;&
\wt m_c(B_d(\ell,\epsilon))
 \\ \leq \;& \mu^-_{\ell(0)}
\big(B_{d_{\ell(0)}}(\ell(-\infty), c_0\; \sqrt{\epsilon}\,)\big)\;
\mu^+_{\ell(0)}
\big(B_{d_{\ell(0)}}(\ell(+ \infty), c_0\; \sqrt{\epsilon}\,)\big)\;.
\end{align*}
Since the Patterson densities are uniformly doubling for basepoints in
$\C\Lambda\Ga$, since the footpoints of the geodesic lines in the
support of $\wt m_c$ belong to $\C\Lambda\Ga$, the Gibbs measure $m_c$
is hence doubling on its support. Let $\ov{m}_{c}=
\frac{m_{c}}{\|m_{c}\|}$.  As in the proof of Theorem
\ref{theo:expratesequid}, using discrete convolution approximation,
there exists $\kappa''>0$ and, for every $\eta>0$, a nonnegative
function $\reg\Phi_\eta \in \operatorname{C}_{\rm
  c}^\alpha(\Ga\bs\G\XX)$ such that
\begin{enumerate}
\item $\int_{\Ga\bs\G\XX}\;\reg\Phi_\eta\;d\,\overline{m}_{c}
  =\int_{\Ga\bs\G\XX} \;\Phi_\eta \;d\,\overline{m}_{c}$,
\item
$\int_{\Ga\bs\G\XX}\;|\reg\Phi_\eta-\Phi_\eta| \;d\,\overline{m}_{c} = 
\bigO\big(\eta\;\int_{\Ga\bs\G\XX}\Phi_\eta \;d\,\overline{m}_{c}\big)$,
\item  
$\|\reg\Phi_\eta\|_\alpha = \operatorname{O}\big(\eta^{-\kappa''}
\int_{\Ga\bs\G\XX}\;\Phi_\eta \;d\,\overline{m}_{c}\big)$.
\end{enumerate}

By Equation \eqref{eq:inttestskinmass}, the integral $\int_{\Ga\bs\G
  \XX} \Phi_{\eta}\;dm_{c}=\|\sigma^+_{\Omega}\|$ is constant (in
particular independent of $\eta$). All integrals below besides the
first one being over $\Ga\bs\G \XX$, and using
\begin{enumerate}
\item[$\bullet$] Equation \eqref{eq:prepaapplimixingexpo} and the
  above property (2) of the regularised map $\reg\Phi_\eta$ for the first
  equality,
\item[$\bullet$] the assumption of exponential decay of correlations
  for the second one, involving some constant $\kappa>0$, for the
  second equality,
\item[$\bullet$] the above properties (1) and (3) of the regularised
  map $\reg\Phi_\eta$ for the last equality,
\end{enumerate}
we hence have
\begin{align*}
\frac{\int_{\Ga\bs\gengeod \XX} \;\psi\;d(\flow n)_*\sigma^+_{\Omega}}
{\|(\flow{n})_*\sigma^+_{\Omega}\|} & =
\frac{\int \; \psi\;\reg\Phi_{\eta}\circ\flow{-n}\;d\,\ov{m}_{c}}
{\int \Phi_{\eta}\;d\,\ov{m}_{c}}+
\bigO(e^{-n\,\alpha}\;\|\psi\|_\alpha+\eta\;\|\psi\|_\infty)\\ & = 
\frac{\int \reg\Phi_{\eta}\;d\,\ov{m}_{c} \; \int  \psi\;d\,\ov{m}_{c}}
{\int \Phi_{\eta}\;d\,\ov{m}_{c}} +
\bigO(e^{-n\,\alpha}\;\|\psi\|_\alpha+\eta\;\|\psi\|_\infty +\\ &\;\;\;\;\;\;\;\;\;
\frac{1}{\int \Phi_{\eta}\;d\,\ov{m}_{c}} 
\;e^{-\kappa \,n}\|\reg\Phi_{\eta}\|_\alpha\|\psi\|_\alpha)\\
 & = \int  \psi\;d\,\ov{m}_{c} + \bigO\big((e^{-n\,\alpha}+\eta +
e^{-\kappa \,n}\;\eta^{-\kappa''})\|\psi\|_\alpha\big)\;.
\end{align*}
Taking $\eta=e^{-\lambda \,n}$ with $\lambda=
\frac{\kappa}{2\kappa''}$, the result follows with $\kappa'=
\min\{\alpha, \frac{\kappa}{2\kappa''}, \frac{\kappa}{2}\}$.
\cqfd

\bigskip
Let us now consider the metric tree case, for the continuous time
geodesic flow, where the main change is to assume a superpolynomial
decay of correlations and hence get a superpolynomial error term, for
observables which are smooth enough along the flow lines. Let
$(\XX,\lambda)$, $X$, $\Ga$, $\wt F$, $\delta=
\delta_{\Ga,\,F^\pm}<\infty$, $(\mu^\pm_x)_{x\in X}$, $\wt m_F$ and
$m_F$ be as in the beginning of Section
\ref{subsec:mixingratemetgraphs}.

\btheo\label{theo:equidmetricerrorterm} Assume that the Gibbs measure
$m_{F}$ on $\Ga\bs\G X$ is finite. Assume furthermore that
\begin{enumerate}
\item the families $(\Lambda\Ga,\mu^-_x, d_x)_{x\in \C\Lambda\Ga}$ and
  $(\Lambda\Ga, \mu^+_x, d_x)_{x\in \C\Lambda\Ga}$ of metric measure
  spaces are uniformly doubling,\footnote{See Section
    \ref{subsec:Pattersondens} for definitions.} and $\wt F$
  is bounded on $T^1\C\Lambda\Ga$,
\item there exists $\alpha\in\;]0,1]$ such that the (continuous time)
    geodesic flow on $(\Ga\bs\G X,m_F)$ has superpolynomial decay of
    $\alpha$-H\"older correlations.
\end{enumerate}
Let $\D=(D_i)_{i\in I}$ be a locally finite $\Ga$-equivariant family
of nonempty proper closed convex subsets of $X$ with $\Ga\bs I$
finite. Let $\Omega= (\Omega_i)_{i\in I}$ be a locally finite
$\Ga$-equivariant family of measurable subsets of $\gengeod X$, with
$\Omega_i\subset \normalpm D_i$ for all $i\in I$ and
$\Omega_i=\Omega_j$ if $i\sim_\D j$.  Assume that $\sigma^\pm_{\Omega}$
is finite and nonzero. Then for every $n\in\NN$, there exists
$k\in\NN$ such that for all $\psi\in \C^{k,\,\alpha}_{\rm c}
(\Ga\bs\gengeod X)$, we have, as $t\ra+\infty$,
$$
\frac{1}{\|(\flow{\pm t})_*\sigma^\pm_{\Omega}\|}\,
\int \psi\;d\,\pi_*(\flow{\pm t})_*\sigma^\pm_{\Omega}\;=\;
\frac{1}{\|m_{F}\|}\,\int \psi\;d\,m_{F} +
\bigO(\|\psi\|_{k,\,\alpha}\,t^{-n})\;.
$$
\etheo

\noindent {\bf Remarks. } (1) If $F=0$, if the metric subtree
$X'=\C\Lambda\Ga$ of $X$ is uniform, if the length spectrum of $\Ga$
on $X$ is not contained in a discrete subgroup of $\RR$ and if $\Ga$
is a lattice of $X'$, then we claim that $\delta=\delta_\Ga$ is
finite, $m_F= m_{\rm BM}$ is finite and mixing, and $(\Lambda\Ga,
\mu_x= \mu_x^\pm,d_x)_{x\in X'}$ is uniformly doubling.

Indeed, the above finiteness and mixing properties follow from
Proposition \ref{prop:uniflatmBMfinie} and Theorem \ref{theo:mixing}.
Since $X'$ is uniform, it has a cocompact discrete sugbroup of
isometries $\Ga'$ whose Patterson density (for the vanishing
potential) is uniformly doubling on $\Lambda\Ga'= \Lambda\Ga$, by
Lemma \ref{lem:shadowlemma} (4). Since $F=0$ and $\Ga$ is a lattice,
the Patterson densities of $\Ga$ and of $\Ga'$ coincide (up to a scalar
multiple) by Proposition \ref{prop:uniflatmBMfinie}.

\medskip\noindent (2) Assume that $F=0$, that the metric subtree
$X'=\C\Lambda\Ga$ of $X$ is uniform, that the length spectrum of $\Ga$
on $X$ is $4$-Diophantine and that $\Ga$ is a geometrically finite
lattice of $X'$.  Then all assumptions of Theorem
\ref{theo:equidmetricerrorterm} are satisfied by the first remark and
by Corollary \ref{coro:expdecaygeomfinimet}.  Therefore we have a
superpolynomially small error term in the equidistribution of the
equidistant levels sets.

\bigskip
\dem The proof is similar to the one of Theorem
\ref{theo:equidsimplicialerrorterm}, except for the doubling property
of the Gibbs measure on its support and the conclusion of the
proof. Let $X'=\C\Lambda\Ga$. The modification of Lemma
\ref{lem:troisboulessimp} used in the previous proof is now the third
assertion of Lemma \ref{lem:troisboulesmet}.

If the footpoints of $\ell,\ell'\in\G X'$ are at distance bounded by
$c_0\,\epsilon_0$, then by Proposition \ref{prop:continuGibbscocycle}
(2), since $|\wt F|$ is bounded on $T^1X'$ by assumption, the
quantities $|C^\pm_\xi(\ell(0), \ell'(0))|$ for $\xi\in \Lambda\Ga$
are bounded by the constant $c'_0=c_0\,\epsilon_0\,(\max_{T^1X'}|\wt
F-\delta|)$.  By the definition of the Gibbs measures (see Equation
\eqref{eq:defigibbs}), Assertion (3) of Lemma \ref{lem:troisboulesmet}
hence implies that if $\epsilon\leq \epsilon_0$ then for every $\ell\in\G\XX$,
\begin{align*}
e^{-2c'_0}\;& \epsilon\;\mu^-_{\ell(0)}
\big(B_{d_{\ell(0)}}(\ell(-\infty), \frac{1}{c_0}\; \sqrt{\epsilon}
\,)\big)\;\mu^+_{\ell(0)}\big(
B_{d_{\ell(0)}}(\ell(+ \infty), \frac{1}{c_0}\; \sqrt{\epsilon}\,)\big)
\\ \leq\; & \wt m_c(B_d(\ell,\epsilon))\\ 
\leq \; & e^{2c'_0}\;\epsilon\;\mu^-_{\ell(0)}
\big(B_{d_{\ell(0)}}(\ell(-\infty), c_0\; \sqrt{\epsilon}\,)\big)\;
\mu^+_{\ell(0)}
\big(B_{d_{\ell(0)}}(\ell(+ \infty), c_0\; \sqrt{\epsilon}\,)\big)\;.
\end{align*}
As in the simplicial case, since the Patterson densities are uniformly
doubling for basepoints in $X'$, the Gibbs measure $m_c$ is hence
doubling on its support.

\medskip 
Fix $n\in\NN$. As in the end of the proof of Theorem
\ref{theo:equidsimplicialerrorterm}, using the assumption of
superpolynomial decay of correlations, involving some degree of
regularity $k$ in order to have polynomial decay in $t^{-Nn}$ where
$N=\lceil\kappa''\rceil +1$, instead of the exponential one, we have,
for all $t\geq 1$ and $\psi\in\C_c^{k,\alpha}(\Ga\bs\gengeod X)$,
$$
\frac{\int_{\Ga\bs\gengeod \XX} \;\psi\;d(\flow t)_*\sigma^+_{\Omega}}
{\|(\flow{t})_*\sigma^+_{\Omega}\|}  = 
\int  \psi\;d\,\ov{m}_{c} + \bigO\big((e^{-t\,\alpha}+\eta +
t^{-Nn}\;\eta^{-\kappa''})\|\psi\|_{k,\,\alpha}\big)\;.
$$ 
Taking $\eta=t^{-n}$, by the definition of $N$, we hence have
$$
\frac{\int_{\Ga\bs\gengeod \XX} \;\psi\;d(\flow t)_*\sigma^+_{\Omega}}
{\|(\flow{t})_*\sigma^+_{\Omega}\|}  = 
\int  \psi\;d\,\ov{m}_{c} + \bigO\big(t^{-n}\;\|\psi\|_{k,\,\alpha}\big)\;.
$$
This proves Theorem \ref{theo:equidmetricerrorterm}. \cqfd

\chapter{Equidistribution of common perpendicular arcs}
\label{sec:equidarcs}

In this Chapter, we prove the equidistribution of the initial and
terminal vectors of common perpendiculars of convex subsets, at the
universal covering space level, for Riemannian manifolds and for
metric and simplicial trees. The results generalise
\cite[Theo.~8]{ParPau16ETDS}.

From now untill the end of Section \ref{subsec:equidcommperpcontIII}, we
consider the continuous time situation where $X$ is a proper
$\CAT(-1)$-space which is either an $\RR$-tree without terminal point
or a complete Riemannian manifold with pinched negative curvature at
most $-1$.  In Section \ref{subsec:equidcommperpdiscrtime}, $X$ will
be the geometric realisation of a simplicial tree $\XX$, and we will
consider the discrete time geodesic flow.

Let $\Ga$ be a nonelementary discrete group of isometries of $X$. Let
$x_0$ be any basepoint in $X$. Let $\wt F$ be a continuous
$\Ga$-invariant potential on $T^1X$, which is H\"older-continuous if
$X$ is a manifold.  Assume that $\delta= \delta_{\Ga,\,F^\pm}$ is
finite and positive and let $(\mu^\pm_x)_{x\in X}$ be (normalised)
Patterson densities for the pairs $(\Ga,F^\pm)$, with associated Gibbs
measure $m_{F}$. Let $\D^-=(D^-_i)_{i\in I^-}$ and $\D^+=(D^+_j)_{j\in
  I^+}$ be locally finite $\Ga$-equivariant families of nonempty
proper closed convex subsets of $X$.

For every $(i,j)$ in $I^-\!\times I^+$ such that the closures
$\overline{D^-_i}$ and $\overline{D^+_j}$ of $D^-_i$ and $D^+_j$ in
$X\cup\partial_\infty X$ have empty intersection, let $\lambda_{i,j}
=d(D^-_i,D^+_j)$ be the length of the common perpendicular from
$D^-_i$ to $D^+_j$, and let $\alpha^-_{i,\,j} \in\gengeod X$ be its
parametrisation: it is the unique map from $\RR$ to $X$ such that

$\bullet$~ $\alpha^-_{i,\,j}(t)= \alpha^-_{i,\,j}(0)\in D_i^-$ if $t\leq 0$,

$\bullet$~ $\alpha^-_{i,\,j}(t)= \alpha^-_{i,\,j} (\lambda_{i,\,j})
\in D_j^+$ if $t\geq \lambda_{i,\,j}$, and 

$\bullet$~ ${\alpha^-_{i,j}}|_{ [0,\,\lambda_{i,j}]} =\alpha_{i,\,j}$
is the shortest geodesic arc starting from a point of $D^-_i$ and
ending at a point of $D^+_j$.  

\noindent Let $\alpha^+_{i,j}= \flow{\lambda_{i,j}}
\alpha^-_{i,j}$. In particular, we have
$\flow{{\frac{\lambda_{i,j}}2}} \alpha^-_{i,j}=
\flow{{\frac{-\lambda_{i,j}}2 }}\alpha^+_{i,j}$.

\medskip 
We now state our main equidistribution result of common perpendiculars
between convex subsets in the continuous time and upstairs
settings. We will give the discrete time version in Section
\ref{subsec:equidcommperpdiscrtime}, and the downstairs version in
Chapter \ref{sec:equidcountdownstairs}.

\btheo\label{theo:mainequidup} Let $X$ be either a proper $\RR$-tree
without terminal points or a complete simply connected Riemannian
manifold with pinched negative curvature at most $-1$. Let $\Ga$ be a
nonelementary discrete group of isometries of $X$ and let $\wt F$ be a
bounded $\Ga$-invariant potential on $X$ which is H\"older-continuous
if $X$ is a manifold. Let $\D^\pm=(D^\pm_k)_{k\in I^\pm}$ be locally
finite $\Ga$-equivariant families of nonempty proper closed convex
subsets of $X$. Assume that the critical exponent $\delta$ is
positive,\footnote{It is finite since $\wt F$ is bounded, see Lemma
  \ref{lem:proprielemcritexpo} (6).} and that the Gibbs measure
$m_{F}$ is finite and mixing for the geodesic flow on $\Ga\bs\G X$.
Then
$$
\lim_{t\ra+\infty} \;\delta\;\|m_{F}\|\;e^{-\delta\, t}
\sum_{\substack{i\in I^-/_\sim,\; j\in I^+/_\sim, \;\ga\in\Ga\\ 
\ov{D^-_i}\cap \ov{D^+_{\ga j}}=\emptyset,\; 
\lambda_{i,\,\ga j}\leq t}} \;e^{\int_{\alpha_{i,\ga j}}\wt F}\; 
\Dirac_{\alpha^-_{i,\,\ga j}} \otimes\Dirac_{\alpha^+_{\ga^{-1}i,\,j}}\;=\; 
\wt\sigma^+_{\D^-}\otimes \wt\sigma^-_{\D^+}\;
$$
for the weak-star convergence of measures on the locally
compact space $\gengeod X\times \gengeod X$.  
\etheo

The proof of Theorem \ref{theo:mainequidup} follows that of
\cite[Theo.~8]{ParPau16ETDS}, which proves this result when $X$ is a
manifold and $F=0$. The first two and a half steps work for both trees
and manifolds and are given in Section
\ref{subsec:equidcommperpcontI}.  The differences begin in Step 3T.
After this, the steps for trees are called $3T$ and $4T$ and are given
in Section \ref{subsec:equidcommperpcontII} and the corresponding
steps for manifolds are $3M$ and $4M$, given in Section
\ref{subsec:equidcommperpcontIII}.

In the special case of $\D^-=(\ga x)_{\ga\in\Ga}$ and $\D^+=(\ga
y)_{\ga\in\Ga}$ for some $x,y\in X$, this statement,\footnote{or
  rather the following Equation \eqref{eq:reduconeij}} with Equation
\eqref{eq:skinnpoint}, gives the following version with potentials of
Roblin's double equidistribution theorem \cite[Theo.~4.1.1]{Roblin03}
when $F=0$, see \cite[Theo.~9.1]{PauPolSha15} for general $F$ when $X$
is a Riemannian manifold with pinched sectional curvature at most
$-1$.

\bcoro\label{coro:caspointempcont} Let $X$, $\Ga$, $\wt F$, $\delta$,
$m_F$ be as in Theorem \ref{theo:mainequidup}, and let $x,y\in X$. We
have
$$
\lim_{t\ra+\infty} \;\delta\;\|m_{F}\|\;e^{-\delta\, t}
\sum_{\ga\in\Ga\,:\,d(x,\ga y)\leq t} \;e^{\int_x^{\ga y}\wt F}\; 
\Dirac_{\ga y} \otimes\Dirac_{\ga^{-1}x}\;=\; 
\mu^+_{x}\otimes \mu^-_{y}\;
$$
for the weak-star convergence of measures on the compact space
$(X\cup\partial_\infty X)\times (X\cup\partial_\infty X)$. \cqfd
\ecoro

\medskip 
Let us give a version of Theorem \ref{theo:mainequidup} without
the assumption $\delta>0$.

\btheo\label{theo:mainequidupbis} Let $X$, $\Ga$, $\wt F$, $\delta$,
$m_F$ be as in Theorem \ref{theo:mainequidup}, except that the
critical exponent $\delta$ is not assumed to be positive.  Then for
every $\tau>0$, we have
$$
\lim_{t\ra+\infty} \;
\frac{\delta\;\|m_{F}\|}{1-e^{-\tau\,\delta}}\;e^{-\delta\, t}
\sum_{\substack{i\in I^-/_\sim,\; j\in I^+/_\sim, \;\ga\in\Ga\\ 
\ov{D^-_i}\cap \ov{D^+_{\ga j}}=\emptyset,\; 
t-\tau<\lambda_{i,\,\ga j}\leq t}} \;e^{\int_{\alpha_{i,\ga j}}\wt F}\; 
\Dirac_{\alpha^-_{i,\,\ga j}} \otimes\Dirac_{\alpha^+_{\ga^{-1}i,\,j}}
\;=\; \wt\sigma^+_{\D^-}\otimes \wt\sigma^-_{\D^+}\;
$$ 
for the weak-star convergence of measures on the locally compact
space $\gengeod X\times \gengeod X$.  
\etheo

\dem The key ingredient in order to deduce Theorem
\ref{theo:mainequidupbis} from Theorem \ref{theo:mainequidup} is the
following classical lemma (see \cite[Lem.~9.5]{PauPolSha15} for a
proof).

\blemm\label{lem:paups95} Let $I$ be a discrete set and let
$f,g:I\to[0,+\infty[$ be maps with $f$ proper. If $\delta+\kappa>0$
and, as $t\ra+\infty$,
$$
\sum_{i\in I,\; f(i)\le t}e^{\kappa\, f(i)}g(i)\;\sim\;
\frac{e^{(\delta+\kappa)t}}{\delta+\kappa}\,,
$$
then for every $c>0$, as $t\ra+\infty$,
$$
\sum_{i\in I,\; t-c<f(i)\le t} g(i)\;\sim\;
\frac{1-e^{-c\,\delta}}{\delta}\,e^{\delta t}\;.\;\;\;\Box
$$ 
\elemm

Let $\kappa>0$ be such that $\delta_{\Ga,\,F+\kappa}
=\delta_{\Ga,\,F}+\kappa>0$.  As the definition of the Gibbs measure
only involves the normalised potential, we have
$\|m_{F+\kappa}\|=\|m_F\|$. Thus, the statement of Theorem
\ref{theo:mainequidup} for the potential $F+\kappa$ is equivalent to
the claim that, as $t\ra+\infty$,
$$
\sum_{\substack{i\in I^-/_\sim,\; j\in I^+/_\sim, \;\ga\in\Ga\\ 
\ov{D^-_i}\cap \ov{D^+_{\ga j}}=\emptyset,\; 
\lambda_{i,\,\ga j}\leq t}} 
\;e^{\kappa\lambda_{i,\,\ga j}}e^{\int_{\alpha_{i,\ga j}}\wt F}\; 
\psi(\alpha^-_{i,\,\ga j},\alpha^+_{\ga^{-1}i,\,j})\;\sim\;
\frac{e^{(\delta+\kappa)\, t}}{(\delta+\kappa)\;\|m_{F}\|} 
\;\int \psi\;d(\wt\sigma^+_{\D^-}\otimes \wt\sigma^-_{\D^+})\;
$$ 
for all positive functions $\psi\in\C_{\rm c}(\gengeod X\times\gengeod
X)$.  Since the function $(i,j,\ga)\mapsto\lambda_{i,\ga j}$ is proper
by the local finiteness of the families $\D^-$ and $\D^+$, Lemma
\ref{lem:paups95} implies that for every $\tau>0$, as $t\ra+\infty$,
$$
\sum_{\substack{i\in I^-/_\sim,\; j\in I^+/_\sim, \;\ga\in\Ga\\ 
\ov{D^-_i}\cap \ov{D^+_{\ga j}}=\emptyset,\; 
t-\tau<\lambda_{i,\,\ga j}\leq t}} 
\;e^{\int_{\alpha_{i,\ga j}}\wt F}\; 
\psi(\alpha^-_{i,\,\ga j},\alpha^+_{\ga^{-1}i,\,j})
\;\sim\;
e^{\delta\, t}\frac{1-e^{\tau\,\delta}}{\delta\;\|m_{F}\|} 
\;\int \psi\;d(\wt\sigma^+_{\D^-}\otimes \wt\sigma^-_{\D^+})\;,
$$
which yields Theorem \ref{theo:mainequidupbis}.
\cqfd

\section{Part I of the proof of Theorem 
\ref{theo:mainequidup}: the common part}
\label{subsec:equidcommperpcontI}

\subsection*{Step 1: Reduction.}  
By additivity, by the local finiteness of the families $\D^\pm$, and
by the definition of $\wt \sigma^\pm_{\D^\mp} =\sum_{k\in I^\mp/_\sim}
\wt\sigma^\pm_{D^\mp_k}$, we only have to prove, for all fixed $i\in
I^-$ and $j\in I^+$, that
\begin{equation}\label{eq:reduconeij}
  \lim_{t\ra+\infty} \;\delta\;\|m_{F}\|\;e^{-\delta\, t}
  \sum_{\ga\in\Ga\,:\;0<\lambda_{i,\,\ga j}\leq t}
  \;\;e^{\int_{\alpha_{i,\ga j}}\wt F}\; 
  \Dirac_{\alpha^-_{i,\,\ga j}} \otimes\Dirac_{\alpha^+_{\ga^{-1}i,\,j}}\;=\; 
  \wt\sigma^+_{D^-_i}\otimes \wt\sigma^-_{D^+_j}\;
\end{equation}
for the weak-star convergence of measures on $\gengeod X\times
\gengeod X$.

Let $\Omega^-$ be a Borel subset of $\normalout D^-_i$ and let
$\Omega^+$ be a Borel subset of $\normalin D^+_j$.  In order to
simplify the notation, let
\begin{align}
  D^-=D^-_i, \;\;D^+=D^+_j, \;\;& \alpha_\ga=\alpha_{i,\ga j},\;\;
\alpha^-_\ga=\alpha^-_{i,\,\ga j}, \;\;\alpha^+_\ga=\alpha^+_{\ga^{-1}i,\, j}, 
\;\;\nonumber\\ &\lambda_\ga=  \lambda_{i,\,\ga j},\;\;
\wt\sigma^\pm= \wt\sigma^\pm_{D^\mp}\,.\label{eq:notationstep1}
\end{align}
Assume that $\Omega^-$ and $\Omega^+$ have positive finite skinning
measures and that their boundaries in $\normalout D^-$ and $\normalin
D^+$ have zero skinning measures (for $\wt\sigma^+$ and $\wt\sigma^-$
respectively). Let
\begin{equation}\label{eq:defiIOmegapm}
I_{\Omega^-,\,\Omega^+}(t)=\delta\;\|m_{F}\|\;e^{-\delta \,t}
\sum_{\substack{\ga\in\Ga\,:\;0<\lambda_{\ga }\leq t\\
\alpha^-_\ga|_{]0,\lambda_\ga]}\in\Omega^-|_{]0,\lambda_\ga]},\;\;
\alpha^+_\ga|_{]-\lambda_\ga,0]}\in\Omega^+|_{]-\lambda_\ga,0]}}}
\;e^{\int_{\alpha_{\ga}}\wt F}\;.
\end{equation}
We will prove the stronger statement, implying Equation
\eqref{eq:reduconeij} by restricting to $\Omega^\pm$ compact, and
useful in this generality for Chapter \ref{sec:equidcountdownstairs},
that, for every such $\Omega^\pm$, we have
\begin{equation}\label{eq:reducnarrowup}
\lim_{t\ra+\infty} \;I_{\Omega^-,\,\Omega^+}(t)\;=\; \wt\sigma^+(\Omega^-)\;
\wt\sigma^-(\Omega^+)\;.
\end{equation}

\subsection*{Step 2: First upper and lower bounds. } 
Using Lemma \ref{lem:defiR} (1), we may fix $R>e^2$ such that $\nu^\pm_w
(V^\mp_{w,\, \eta,\,R}) >0$ for all $\eta\in\;]0,1]$ and $w\in\normalmp
D^\pm$. Let $\phi^\mp_\eta=\phi^\mp_{\eta,\,R,\,\Omega^\pm}$ be the test
functions defined in  Equation \eqref{eq:defiphi}.

For all $t\geq 0$, let
\begin{equation}\label{eq:defiaetat}
a_\eta(t)=\sum_{\ga\in\Ga}\;
\int_{\ell\in\G X}\phi^-_\eta(\flow{-t/2}\ell)\;
\phi^+_\eta(\flow{t/2}\ga^{-1}\ell)\;d\wt m_{F}(\ell)\,.
\end{equation}
As in \cite{ParPau16ETDS}, the heart of the proof is to give two pairs
of upper and lower bounds, as $T\geq 0$ is large enough and $\eta\in
\;]0,1]$ is small enough, of the (Ces\`aro-type) quantity
\begin{equation}\label{eq:defiIetapmT}
i_{\eta}(T)=\int_0^{T}e^{\delta\,t}\; a_\eta(t)\;dt\,.
\end{equation}

By passing to the universal cover, the mixing property of the geodesic
flow on $\Ga\backslash\G X$ for the Gibbs measure $m_{F}$ gives that,
for every $\epsilon>0$, there exists $T_\epsilon=T_{\epsilon,\eta}
\geq 0$ such that for all $t\geq T_\epsilon $, we have
$$
\frac{e^{-\epsilon}}{\|m_{F}\|}\;
\int_{\G X} \phi^-_\eta\;d\wt m_{F}\int_{\G X} \phi^+_\eta\;d\wt m_{F}
\leq a_\eta(t)\leq \frac{e^\epsilon}{\|m_{F}\|}\;
\int_{\G X} \phi^-_\eta\;d\wt m_{F}\int_{\G X} \phi^+_\eta\;d\wt m_{F}\,.
$$
Hence by Lemma \ref{lem:integrable} (1), for all $\epsilon>0$ and
$\eta\in\;]0,1]$, there exists $c_\epsilon= c_{\epsilon,\eta}>0$ such
that for every $T\geq 0$, we have
\begin{equation}\label{eq:firstmajoI?????}
e^{-\epsilon}\; \frac{e^{\delta\,T}}{\delta\,\|m_{F}\|}\;
\wt \sigma^+(\Omega^-)\;\wt \sigma^-(\Omega^+)-c_\epsilon 
\leq i_{\eta}(T)\leq 
e^{\epsilon}\; \frac{e^{\delta\,T}}{\delta\,\|m_{F}\|}\;
\wt \sigma^+(\Omega^-)\;\wt \sigma^-(\Omega^+)+c_\epsilon\,.
\end{equation}

\subsection*{Step 3: Second upper and lower bounds. }  Let $T\geq 0$
and $\eta\in\;]0,1]$.  By Fubini's theorem for nonnegative measurable
maps, the definition\footnote{See Equation \eqref{eq:defiphi}.} of the
test functions $\phi^\pm_\eta$ and the flow-invariance\footnote{See
  Equation \eqref{eq:equivfibrationf}.}  of the fibrations
$f^\pm_{D^\mp}$, we have
\begin{align}
i_{\eta}(T)=\sum_{\ga\in\Ga}\;\int_0^{T}e^{\delta\,t}\;
\int_{\ell\in\G X}\;\; &h^-_{\eta,\,R}\circ f^+_{D^-}(\ell)\;
h^+_{\eta,\,R}\circ f^-_{D^+}(\ga^{-1}\ell)\nonumber\\ &
\mathbbm{1}_{\flow{t/2}\V^+_{\eta,\,R}(\Omega^-)}(\ell)\;
\mathbbm{1}_{\flow{-t/2}\V^-_{\eta,\,R}(\ga\Omega^+)}(\ell)
\;d\wt m_{F}\;dt\,.\label{eq:rewriteIsetp4}
\end{align}
We start the computations by rewriting the product term involving the
functions $h^\pm_{\eta,\,R}$. For all $\ga\in\Ga$ and $\ell\in
\U^+_{D^-}\cap \U^-_{\ga D^+}$, define (using Equation
\eqref{eq:equivfibrationf})
\begin{equation}\label{eq:defwpm}
w^-=f^+_{D^-}(\ell)\in\G_{+,\,0}X\;\;\;{\rm  and}\;\;\; 
w^+= f^-_{\ga D^+}(\ell)=\ga f^-_{ D^+}(\ga^{-1}\ell)\in\G_{-,\,0}X \;.
\end{equation}
This notation is ambiguous ($w^-$ depends on $\ell$, and $w^+$ depends
on $\ell$ and $\ga$), but it makes the computations less heavy. By
Equations \eqref{eq:commutflowh} and \eqref{eq:changemoinsplus}, we
have, for every $t\geq 0$,
\begin{align*}
h^-_{\eta,\,R}(w^-)&=
h^-_{\eta,\,R}\circ \flow{-t/2}(\flow{t/2}w^-)=
e^{\int_{w^-(0)}^{w^-(t/2)} (\wt F-\delta)}\; 
h^-_{\eta,\,e^{-t/2}R}(\flow{t/2}w^-)\;.
\end{align*}
Similarly, 
$$
h^+_{\eta,\,R}(\ga^{-1}w^+)=
e^{\int_{w^+(-t/2)}^{w^+(0)} (\wt F-\delta)}\; 
h^+_{\eta,\,e^{-t/2}R}(\flow{-t/2}w^+)\,.
$$
Hence, 
\begin{align}
h^-_{\eta,\,R}\circ f^+_{D^-}(\ell)& \;
h^+_{\eta,\,R}\circ f^-_{D^+}(\ga^{-1}\ell)\nonumber
\\ =\;& 
e^{-\delta\,t}\;
e^{\int_{w^-(0)}^{w^-(t/2)} \wt F+
\int_{w^+(-t/2)}^{w^+(0)} \wt F}\;
h^-_{\eta,\,e^{-t/2}R}(\flow{t/2}w^-)\,
h^+_{\eta,\,e^{-t/2}R}(\flow{-t/2}w^+)\,.\label{eq:commonuptohere}
\end{align}

\section{Part II of the proof of Theorem 
\ref{theo:mainequidup}: the metric tree case}
\label{subsec:equidcommperpcontII}

In this Section, we assume that $X$ is an $\RR$-tree and we will
consider the manifold case separately in Section
\ref{subsec:equidcommperpcontIII}.

\subsection*{Step 3T.} 
Consider the product term in Equation \eqref{eq:rewriteIsetp4}
involving the characteristic functions. By Lemma
\ref{lem:creationperp} (applied by replacing $D^+$ by $\ga D^+$),
there exists $t_0\geq 2\ln R+4$ such that for all $\eta\in \;]0,1]$,
$t\geq t_0$, $\ga\in\Ga$ and $\ell\in\G X$, if
$\mathbbm{1}_{\flow{t/2}\V^+_{\eta, \,R} (\Omega^-)} (\ell)\;
\mathbbm{1}_{\flow{-t/2}\V^-_{\eta,\,R}(\ga\Omega^+)} (\ell)\neq 0$,
or equivalently by Equation
\eqref{eq:flowbehavdynaneigh} if
$$
\ell\in \V^+_{\eta,\,e^{-t/2}R}(\flow{t/2}\Omega^-)\cap
\V^-_{\eta,\,e^{-t/2}R}(\ga\flow{-t/2}\Omega^+)\;,
$$
then the following facts hold.

\begin{enumerate}
\item[(i)] By the convexity of $D^\pm$, we have $\ell\in \U^+_{D^-}\cap
  \U^-_{\ga D^+}$.
\item[(ii)] By the definition\footnote{See Equation
  \eqref{eq:defwpm}.}  of $w^{\pm}$, we have $w^{-}\in \Omega^{-}$ and
  $w^{+}\in \ga \Omega^{+}$. The notation $(w^-,w^+)$ here coincides
  with the notation $(w^-, w^+)$ in Lemma \ref{lem:creationperp}.
\item[(iii)] There exists a common perpendicular $\alpha_\ga$ from
  $D^-$ to $\ga D^+$, whose length $\lambda_\ga$ satisfies
\begin{equation}\label{eq:tpreslambdaga}
|\;\lambda_\ga - t\;|\leq 2\eta\,,
\end{equation}
whose origin is $\alpha^-_\ga(0)=w^-(0)$, whose endpoint is
$\ga\,\alpha^+_\ga(0) =w^+(0)$, such that the points $w^- (\frac t2)$
and $w^+(-\frac t2)$ are at distance at most $\eta$ from
$\ell(0)\in\alpha_\ga$.  
\end{enumerate}
See the picture following Lemma \ref{lem:creationperp} (replacing
$D^+$ by $\ga D^+$). Hence, by Lemma \ref{lem:potentialcontinuity} and
since $\wt F$ is bounded,
\begin{align}
e^{-2\eta\,\|F\|_\infty}\;e^{\int_{\alpha_\ga}\wt F}\leq 
e^{\int_{w^-(0)}^{w^-(t/2)} \wt F+ \int_{w^+(-t/2)}^{w^+(0)} \wt F}\leq
e^{2\eta\,\|F\|_\infty}\;e^{\int_{\alpha_\ga}\wt F}\;.
\label{eq:controlpieceFlengthT}
\end{align}

For all $\eta\in \;]0,1]$, $\ga\in\Ga$ and $T\geq t_0$, let
$$
\A_{\eta,\ga}(T)=\big\{(t,\ell)\in [t_0,T] \times\G X:
\ell\in \V^+_{\eta,\,e^{-t/2}R}(\flow{t/2}\Omega^-)\cap
\V^-_{\eta,\,e^{-t/2}R}(\ga\flow{-t/2}\Omega^+)\big\}
$$
and
$$
  j_{\eta,\,\ga}(T)=\iint_{(t,\,\ell)\in \A_{\eta,\ga}(T)}
  h^-_{\eta,\,e^{-t/2}R}(\flow{t/2}w^-)\;
  h^+_{\eta,\,e^{-t/2}R}(\flow{-t/2}w^+)\;dt\;d\wt m_{F}(\ell)\,.
$$
By the above, since the integral of a function is equal to the
integral on any Borel set containing its support, and since the
integral of a nonnegative function is nondecreasing in the integration
domain, there hence exists $c_4>0$ such that for all $T\geq 0$ and
$\eta\in\;]0,1]$, we have
$$
i_{\eta}(T) \geq -\,c_4+
e^{-2\eta\|F\|_\infty}\sum_{\substack{\ga\in\Ga\,:\;
t_0+2\leq\lambda_{\ga}\leq T-2\eta\\
\alpha^-_\ga|_{[0,\,\lambda_\ga]}\in\Omega^-|_{[0,\,\lambda_\ga]}, \;\;
\alpha^+_\ga|_{[-\lambda_\ga,\,0]}\in\Omega^+|_{[-\lambda_\ga,\,0]}}}
\;e^{\int_{\alpha_\ga}\wt F}\;j_{\eta,\,\ga}(T)\;\;,
$$
and similarly, for every $T'\geq T$ (later on, we will take $T'$ to be
$T+ 4\eta$),
$$
i_{\eta}(T) \leq \,c_4+
e^{2\eta\,\|F\|_\infty}\sum_{\substack{\ga\in\Ga\,:\;
t_0+2\leq\lambda_{\ga}\leq 
T+2\eta\\
\alpha^-_\ga|_{[0,\,\lambda_\ga]}\in\Omega^-|_{[0,\,\lambda_\ga]}, 
\;\;\alpha^+_\ga|_{[-\lambda_\ga,\,0]}\in\Omega^+|_{[-\lambda_\ga,\,0]}}}\;
e^{\int_{\alpha_\ga}\wt F}\;j_{\eta,\,\ga}(T')\;.
$$

\subsection*{Step 4T: Conclusion. } Let $\epsilon>0$. Let
$\ga\in\Ga$ be such that $D^-$ and $\ga D^+$ do not intersect and the
length of their common perpendicular satisfies $\lambda_\ga\geq
t_0+2$. Let us prove that if $\eta$ is small enough and $\lambda_\ga$
is large enough,\footnote{with the enough's independent of $\ga$} then for
every $T\geq \lambda_\ga+ 2\eta$, we have
\begin{equation}\label{eq:step5T}
 1-\epsilon\leq j_{\eta,\,\ga}(T)\leq 1+\epsilon\,.
\end{equation}
This estimate proves the claim \eqref{eq:reducnarrowup}, as
follows. For every $\epsilon>0$, if $\eta>0$ is small enough, we have
$$
i_\eta(T+2\eta)\ge-c_4+e^{-2\eta\|F\|_\infty}(1-\epsilon)
\Big(\frac{I_{\Omega^-,\Omega^+}(T)}{\delta\,\|m_{F}\|
e^{-\delta\, T}}-\frac{I_{\Omega^-,\Omega^+}(t_0+2)}{\delta\,\|m_{F}\|
e^{-\delta (t_0+2)}}\Big)
$$
and by Equation \eqref{eq:firstmajoI?????}
$$
i_\eta(T+2\eta)\le c_{\epsilon}+
\frac{e^\epsilon\,\wt\sigma^+(\Omega^-)\,\wt\sigma^-(\Omega^+)}
{\delta\, \|m_{F}\|\,e^{-\delta(T+2\eta)}}\,.
$$
Thus, for $\eta$ small enough, 
$$
\wt\sigma^+(\Omega^-)\,\wt\sigma^-(\Omega^+)\ge
\frac{1-\epsilon}{e^{\epsilon}}\,I_{\Omega^-,\Omega^+}(T)+\smallo(1)
$$
as $T\to+\infty$, which gives 
$$
\limsup_{T\to+\infty}I_{\Omega^-,\Omega^+}(T)\le
\wt\sigma^+(\Omega^-)\,\wt\sigma^-(\Omega^+).
$$
The similar estimate for the lower limit proves the claim
\eqref{eq:reducnarrowup}.
  
In order to prove the claim \eqref{eq:step5T}, let $\eta\in\;]0,1]$
and $T\geq \lambda_\ga+2\eta$.  In order to simplify the notation, let
$$
r_t=e^{-t/2}R, \;\; w^-_t=\flow{t/2}w^-\;\;\;{\rm and}\;\;\;
w^+_t=\flow{-t/2}w^+\,.
$$ 
By the definition of $j_{\eta,\,\ga}$, using the inequalities
\eqref{eq:estimhsmaletaprimp} where the constant $c_1$ is uniform
since $\wt F$ is bounded (with the comment following them) and
the fact that $r_t=\bigO(e^{-\lambda_\ga/2})$ by Equation
\eqref{eq:tpreslambdaga}, we hence have
\begin{align}
j_{\eta,\,\ga}(T)&=
\iint_{(t,\ell)\in\A_{\eta,\ga}(T)} 
h^-_{\eta,\,r_t}(w^-_t)\;
h^+_{\eta,\,r_t}(w^+_t)\;dt\;d\wt m_{F}(\ell)\nonumber\\ & 
=\frac{e^{\bigO(e^{-\lambda_\ga/2})}}{(2\eta)^2}\!
\iint_{(t,\,\ell)\in\A_{\eta,\ga}(T)} 
\frac{dt\;d\wt m_{F}(\ell)}
{\mu_{\wss(w^-_t)}(B^+(w^-_t,\,r_t))\;
\mu_{\wsu(w^+_t)}(B^-(w^+_t,\,r_t))}\,.\label{eq:adaptatJetagaT}
\end{align}

Let $x_\ga$ be the midpoint of the common perpendicular $\alpha_\ga$,
so that $d(x_\ga,\ell(0))=\bigO(\eta)$ for every $(t,\ell)\in
\A_{\eta,\,\ga}(T)$ by the above Claim (iii).  Let us use the Hopf
parametrisation of $\G X$ with basepoint $x_\ga$, denoting by $s$ its
time parameter.  When $(t,\ell)\in \A_{\eta,\ga}(T)$, by Definition
\eqref{eq:defigibbs} of the Gibbs measure $\wt m_F$, by the $\RR$-tree
case of Proposition \ref{prop:continuGibbscocycle} (2), and since $\wt
F$ is bounded, we have
\begin{align}\label{eq:gibbseta}
d\wt m_{F}(\ell) & =
e^{C_{\ell_-}^-(x_\ga,\,\ell(0))\,+\,C^+_{\ell_+}(x_\ga,\,\ell(0))}\;
d\mu_{x_\ga}^-(\ell_-)\,d\mu^+_{x_\ga}(\ell_+)\,ds\nonumber\\
&= e^{\bigO(\eta)}d\mu_{x_\ga}^-(\ell_-)\,d\mu^+_{x_\ga}(\ell_+)\,ds\,.
\end{align}
  
Let $P_\ga$ be the plane domain of the $(t,s)\in\RR^2$ such that
$|\lambda_\ga-t|\leq 2\eta$ and there exist $s^\pm\in\;]-\eta,\eta[$
with $s^\mp= \frac{\lambda_\ga-t}{2} \pm s$. It is easy to see that
$P_\ga$ is a rhombus centred at $(\lambda_\ga,0)$ whose area is
$(2\eta)^2$.

\begin{center}
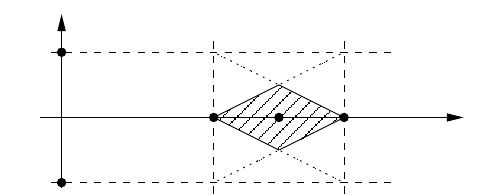
\end{center}

Let $\xi_\ga^\pm$ be the point at infinity of any fixed geodesic ray
from $x_\ga$ through $\alpha_\ga^\pm(0)$. If $A$ is a subset of $\G
X$, we denote by $A_\pm$ the subset $\{\ell_\pm\;:\;\ell\in A\}$ of
$\partial_\infty X$.

\blemm\label{lem:passeralinfini} For every $t\geq t_0$ such that
$|\lambda_\ga-t|\leq 2\eta$, we have
$$
(B^\pm(w^\mp_t,r_t))_\mp =
B_{d_{x_\ga}}(\xi^\mp_\ga, R\,e^{-\frac{\lambda_\ga}2})\,.
$$
\elemm

\dem We prove the statement for the negative endpoints, the proof of
the claim for positive endpoints is similar. Since $R\geq 1$ and
$d(x_\ga, \alpha_\ga^\mp(0))=\frac{\lambda_\ga}2$, the term on the right
hand side does not not depend on the choice of $\xi_\ga^\mp$.

Let us first prove the inclusion on the set on left hand side into the
one on the right hand side. Let $\ell'\in B^+(w^-_t,r_t)$, so
that\footnote{See the definition of the strong stable ball $B^+(w^-_t,r_t)$
  in Section \ref{subsec:unitbundle}.} there exists a geodesic line
$\wh w^-_t\in\G X$, extension of the geodesic ray $w^-_t: \mathopen{[}
-\frac{t}{2}, +\infty\mathclose{[}\ra X$, with
    $d_{W^+(w^-_t)}(\wh w^-_t,\ell')<r_t$. We may assume that $(\wh
    w^-_t)_-= \xi_\ga^-$ and $\ell'_-\neq \xi^-_\ga$. Let $p\in X$ be
    such that $[\ell'(0), \xi^-_\ga[
        \;\cap\,[\ell'(0),\ell'_-[\;=\;[\ell'(0),p]$.

\begin{center}
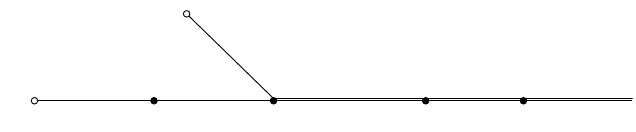
\end{center}

\noindent 
Since $t\geq t_0> 2\ln R$, we have $r_t<1$, hence\footnote{See the
  definition of the Hamenst\"adt distance $d_{W^+(w^-_t)}$ in Section
  \ref{subsec:unitbundle}.} $\ell'(0)=w^-_t(0)= w^-(t/2)$ and $p\in\;
\mathopen{]}\ell'(0), \xi^-_\ga\mathclose{[}\,$. Since 
$$
d(p, w^-(t/2))=-\ln d_{W^-+w^-_t)}(\wh w^-_t,\ell')>-\ln r_t=
\frac{t}{2}-\ln R\geq \frac{t_0}{2}-\ln R\geq 1\geq \eta
$$ 
and $d(x_\ga, w^-(t/2))=\big|\frac{\lambda_\ga}{2}-\frac{t}{2}\big|
\leq \eta$, we have $p\in [x_\ga, \xi^-_\ga[\,$. Hence
$$
d_{x_\ga}(\ell'_-,\xi^-_\ga)=e^{-d(p,\,x_\ga)}\leq\Big\{\begin{array}{ll}
e^{-d(p,\,w^-(\frac t 2)) -d(w^-(t/2),\,x_\ga)}& 
{\rm if}\;\frac{t}{2}\leq \lambda_\ga\\
e^{-d(p,\,w^-(\frac t 2)) +d(w^-(t/2),\,x_\ga)}& 
{\rm otherwise.}\end{array}
$$
In both cases,
$$
d_{x_\ga}(\ell'_-,\xi^-_\ga)=
d_{W^+(w^-_t)}(\wh w^-_t,\ell')\,e^{-\frac{\lambda_\ga}{2}+\frac{t}{2}}
< r_t\,e^{-\frac{\lambda_\ga}{2}+\frac{t}{2}}=  R\,e^{-\frac{\lambda_\ga}2}\;.
$$ 

Conversely, if $\xi\in B_{d_{x_\ga}}(\xi^-_\ga, R\,e^{-\frac{\lambda_\ga}2})$, 
let $\ell' \in \G X$ be such that $\ell'(0) = w^-(t/2)$ and $\ell'_-=
\xi$. We may assume that $\xi\neq \xi^-_\ga$. Let $\wh w^-_t$ be the
extension of $w^-_t$ such that $(\wh w^-_t)_-= \xi^-_\ga$.  Let $p\in
X$ be such that $[\ell'(0),\xi^-_\ga[ \;\cap\,[\ell'(0),\ell'_-[ \;
=\;[\ell'(0),p]$. Then as above, we have $R\, e^{-\frac{\lambda_\ga}2}<1$, 
hence
$$
d_{W^+(w^-_t)}(\wh w^-_t,\ell')=e^{-d(p,\,\ell'(0))}= 
d_{x_\ga}(\ell'_-,\xi^-_\ga)\;e^{\frac{\lambda_\ga}{2}-\frac{t}{2}}<  r_t\;,
$$ 
thus $\ell'\in B^+(w^-_t,r_t)$. \cqfd

\medskip It follows from this lemma that, for all $t\geq t_0$,
$s^\pm\in\;]-\eta,\eta[$ and $\ell\in\G X$, we have $\flow{\mp s^\mp}
\ell\in B^\pm(w^\mp_t,r_t)$ if and only if $d(\ell(0), \alpha^\pm_\ga
(0)) =s^\pm+\frac t2$ (or equivalently, by the definition of the time
parameter $s$ of $\ell$ in Hopf's parametrisation with basepoint
$x_\ga$, when $s^\pm+\frac t2=\frac {\lambda_\ga}{2}\pm s$), and
$\ell_\pm \in B_{d_{x_\ga}}(\xi_\ga^\pm,
R\,e^{\frac{-\lambda_\ga}2})$.  Thus,
$$
\A_{\eta,\ga}(T)=P_\ga\times 
B_{d_{x_\ga}}(\xi_\ga^-,R\,e^{-\frac{\lambda_\ga}2})
\times 
B_{d_{x_\ga}}(\xi_\ga^+,R\,e^{-\frac{\lambda_\ga}2})\,.
$$

\medskip To finish Step 4T and the proof of the theorem for
$\RR$-trees, note that by the definition of the skinning measure
(using again the Hopf parametrisation with basepoint $x_\ga$), by the
above Lemma \ref{lem:passeralinfini}, by the claim for $\RR$-trees of
Assertion (2) of Proposition \ref{prop:continuGibbscocycle} and the
boundedness of $\wt F$, we have
\begin{equation}\label{eq:penultimatePart4T}
\mu_{W^\pm(w^\mp_t)}(B^\pm(w^\mp_t,r_t))= e^{\bigO(\eta)}
\mu^\mp_{x_\ga}(B_{d_{x_\ga}}(\xi^\mp_\ga,R\,e^{-\frac{\lambda_\ga}2}))\,.
\end{equation}
Thus, by the above and by Equations \eqref{eq:adaptatJetagaT} and
\eqref{eq:gibbseta}, (and noting that
$\bigO(\eta)\pm\bigO(\eta)=\bigO(\eta)$)
\begin{align}
j_{\eta,\,\ga}(T) & =
\frac{e^{\bigO(e^{-\frac{\lambda_\ga}2})}}{4\eta^2} e^{\bigO(\eta)}(2\eta)^2
\frac{\mu^-_{x_\ga}(B_{d_{x_\ga}}(\xi^-_\ga,Re^{-\frac{\lambda_\ga}2}))
  \mu^+_{x_\ga}(B_{d_{x_\ga}}(\xi^+_\ga,Re^{-\frac{\lambda_\ga}2}))}
{\mu^-_{x_\ga}(B_{d_{x_\ga}}(\xi^-_\ga,Re^{-\frac{\lambda_\ga}2}))
  \mu^+_{x_\ga}(B_{d_{x_\ga}}(\xi^+_\ga,Re^{-\frac{\lambda_\ga}2}))}
\nonumber \\ & =e^{\bigO(\eta+e^{-\frac{\lambda_\ga}2})}\,,
\label{eq:effcontrojtree}
\end{align}
which gives the inequalities \eqref{eq:step5T}.
\cqfd

\bigskip 
The effective control on $j_{\eta,\,\ga}(T)$ given by Equation
\eqref{eq:effcontrojtree} is stronger than what is needed in order to
prove Equation \eqref{eq:step5T} in Step 4T. We will use it in Section
\ref{subsect:errormetricgraphgroup} in order to obtain error terms.

\section{Part III of the proof of Theorem 
\ref{theo:mainequidup}: the manifold case}
\label{subsec:equidcommperpcontIII}

The proof of Theorem \ref{theo:mainequidup} for manifolds is the same
one as for trees until Equation \eqref{eq:commonuptohere}.  The
remaining part of the proof that we give below is more technical than
for trees but the structure of the proof is similar. In this Section,
$X=\wt M$ is a Riemannian manifold, and we identify $\G X$ and $T^1 X$
with the standard unit tangent bundle of $\wt M$, as explained in Section
\ref{subsec:unitbundle}.

\subsection*{Step 3M.}  
Consider the product term in Equation \eqref{eq:rewriteIsetp4}
involving the characteristic functions. The quantity
$\mathbbm{1}_{\V^+_{\eta, \,R} (\Omega^-)} (\flow{-t/2}v)\;
\mathbbm{1}_{\V^-_{\eta,\,R}(\Omega^+)} (\ga^{-1} \flow{t/2}v)$ is
different from $0$ (hence equal to $1$) if and only if
$$
v\in \flow{t/2}\V^+_{\eta, \,R}
(\Omega^-)\cap \ga\flow{-t/2} \V^-_{\eta,\,R} (\Omega^+)=
\V^+_{\eta,\,e^{-t/2}R}(\flow{t/2}\Omega^-)\cap
\V^-_{\eta,\,e^{-t/2}R}(\ga\flow{-t/2}\Omega^+)\,,
$$ 
see Section \ref{subsect:nbhd} and in particular Equation
\eqref{eq:flowbehavdynaneigh}.  By Lemma \ref{lem:creationperpRiem} 
(applied by replacing $D^+$ by $\ga D^+$ and $w$ by $v$), there exist
$t_0,c_0>0$ such that for all $\eta\in \;]0,1]$ and $t\geq t_0$, for
all $v\in T^1\wt M$, if $\mathbbm{1}_{\V^+_{\eta, \,R} (\Omega^-)}
(\flow{-t/2}v)\; \mathbbm{1}_{\V^-_{\eta,\,R}(\Omega^+)} (\ga^{-1}
\flow{t/2}v)\neq 0$, then the following facts hold:

\begin{enumerate}
\item[(i)] by the convexity of $D^\pm$, we have $v\in \U^+_{D^-}\cap
  \U^-_{\ga D^+}$,
\item[(ii)] by the definition of $w^{\pm}$ (see Equation
  \eqref{eq:defwpm}), we have $w^{-}\in \Omega^{-}$ and $w^{+}\in \ga
  \Omega^{+}$ (the notation $(w^-,w^+)$ here coincides with the
  notation $(w^-, w^+)$ in Lemma \ref{lem:creationperpRiem}),
\item[(iii)] there exists a common perpendicular $\alpha_\ga$ from
  $D^-$ to $\ga D^+$,\footnote{and we then denote as previously by
  $v^-_\ga$ its tangent vector at its origin, by $v^+_\ga$ its tangent
  vector at its terminal point, and by $v^0_\ga$ its tangent vector at
  its midpoint} whose length $\lambda_\ga$ satisfies
\begin{equation}\label{eq:tpreslambdagamani}
|\;\lambda_\ga - t\;|\leq 2\eta + c_0\,e^{-t/2}\;,
\end{equation} 
whose origin $\pi(v^-_\ga)$ is at distance at most $c_0\,e^{-t/2}$
from $\pi(w^-)$, whose endpoint $\pi(v^+_\ga)$ is at distance at most
$c_0\,e^{-t/2}$ from $\pi(w^+)$, such that both points $\pi(\flow{t/2}
w^-)$ and $\pi(\flow{-t/2}w^+)$ are at distance at most $\eta+
c_0\,e^{-t/2}$ from $\pi(v)$, which is at distance at most
$c_0\,e^{-t/2}$ from some point $p_v$ of $\alpha_\ga$.
\end{enumerate}

\begin{center}
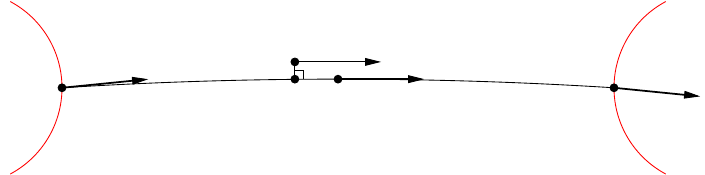
\end{center}

Using (iii) and the \eqref{eq:HC}-property\footnote{See Definition
  \ref{defi:HCproperty}.} which introduces a constant $\kappa_2\in
  \;]0,1]$, and since $\wt F$ is bounded, for all $\eta\in \;]0,1]$,
  $t\geq t_0$ and $v\in T^1\wt M$ for which
  $\mathbbm{1}_{\V^+_{\eta,\,R}(\Omega^-)} (\flow{-t/2}v)\;
  \mathbbm{1}_{\V^-_{\eta,\,R}(\Omega^+)} (\ga^{-1} \flow{t/2}v)\ne
  0$, we have
\begin{align}
e^{\int_{\pi(w^-)}^{\pi(\flow{t/2}w^-)} \wt F+
\int_{\pi(\flow{-t/2}w^+)}^{\pi(w^+)} \wt F}&=
e^{\int_{\pi(v^-_\ga)}^{p_v} \wt F+
\int_{p_v}^{\pi(v^+_\ga)} \wt F+
\operatorname{O}((\eta+e^{-t/2})^{\kappa_2})}\nonumber\\ & =
e^{\int_{\alpha_\ga}\wt F}\;e^{\operatorname{O}((\eta+e^{-\lambda_\ga/2})^{\kappa_2})}\;.
\label{eq:controlpieceFlength}
\end{align}

For all $\eta\in \;]0,1]$, $\ga\in\Ga$ and $T\geq t_0$, define
$$
\A_{\eta,\ga}(T)=\big\{(t,v)\in [t_0,T] \times T^1\wt M \;: \;
v\in \V^+_{\eta,\,e^{-t/2}R}(\flow{t/2}\Omega^-)\cap
\V^-_{\eta,\,e^{-t/2}R}(\ga\flow{-t/2}\Omega^+)\big\}\,,
$$ and
$$
  j_{\eta,\,\ga}(T)=\iint_{(t,\,v)\in \A_{\eta,\ga}(T)}
  h^-_{\eta,\,e^{-t/2}R}(\flow{t/2}w^-)\;
  h^+_{\eta,\,e^{-t/2}R}(\flow{-t/2}w^+)\;dt\;d\wt m_F(v)\,.
$$ 
By the above, since the integral of a function is equal to the
integral on any Borel set containing its support, and since the
integral of a nonnegative function is nondecreasing in the integration
domain, there hence exists $c_4>0$ such that for all $T\geq 0$ and
$\eta\in\;]0,1]$, we have
$$
i_{\eta}(T) \geq -\,c_4+\!\!\!\!\!\!\!\!\!
\sum_{\substack{\ga\in\Ga\,:\;
t_0+2+c_0\leq\lambda_{\ga}\leq 
T-\operatorname{O}(\eta+e^{-\lambda_\ga/2})\\
v^-_{\ga}\in\,\N_{-\operatorname{O}(\eta+e^{-\lambda_\ga/2})}\Omega^-,\;
v^+_{\ga}\in\,\ga\N_{-\operatorname{O}(\eta+e^{-\lambda_\ga/2})}\Omega^+}}
\!\!\!\!\!\!\!\!e^{\int_{\alpha_\ga}\wt F}\;j_{\eta,\,\ga}(T)\;
e^{-\operatorname{O}((\eta+e^{-\lambda_\ga/2})^{\kappa_2})}\;,
$$
and similarly, for every $T'\geq T$,
$$
i_{\eta}(T)\leq c_4+
\!\!\!\!\!\!\!\!\!
\sum_{\substack{\ga\in\Ga\,:\;
t_0+2+c_0\leq\lambda_{\ga}\leq 
T+\operatorname{O}(\eta+e^{-\lambda_\ga/2}) \\
v^-_{\ga}\in\,\N_{\operatorname{O}(\eta+e^{-\lambda_\ga/2})}\Omega^-,\;
v^+_{\ga}\in\,\ga\N_{\operatorname{O}(\eta+e^{-\lambda_\ga/2})}\Omega^+}}
\!\!\!\!\!\!\!\!\!
e^{\int_{\alpha_\ga}\wt F}\;j_{\eta,\,\ga}(T')\;
e^{\operatorname{O}((\eta+e^{-\lambda_\ga/2})^{\kappa_2})}\;.
$$
We will take $T'$ to be of the form $T+ \operatorname{O}(\eta
+e^{-\lambda_\ga/2})$, for a bigger $\operatorname{O}(\cdot)$ than the
one appearing in the index of the above summation.

\subsection*{Step 4M: Conclusion. } Let $\ga\in\Ga$ be such
that $D^-$ and $\ga D^+$ have a common perpendicular with length
$\lambda_\ga\geq t_0+2+c_0$. Let us prove that for all $\epsilon>0$, if
$\eta$ is small enough and $\lambda_\ga$ is large enough, then for every
$T\geq \lambda_\ga+ \operatorname{O}(\eta +e^{-\lambda_\ga/2})$ (with the
enough's and $\operatorname{O}(\cdot)$ independent of $\ga$), we have
\begin{equation}\label{eq:step5}
 1-\epsilon\leq j_{\eta,\,\ga}(T)\leq 1+\epsilon\,.
\end{equation}
Note that $\wt\sigma^\pm(\N_\varepsilon(\Omega^\mp))$ and
$\wt\sigma^\pm (\N_{-\varepsilon}(\Omega^\mp))$ tend to $\wt\sigma^\pm
(\Omega^\mp)$ as $\varepsilon\ra 0$ (since $\wt\sigma^\pm
(\partial\Omega^\mp)=0$ as required in Step 1). Using Steps 2, 3M and
4M, this will prove Equation \eqref{eq:reducnarrowup}, hence will
complete the proof of Theorem \ref{theo:mainequidup}.

We say that $(\wt M,\Ga,\wt F)$ {\em has radius-continuous strong
  stable/unstable ball masses}\index{radius-continuous ball masses} if
for every $\epsilon>0$, if $r>1$ is close enough to $1$, then for
every $v\in T^1\wt M$, if $B^-(v,1)$ meets the support of
$\mu^+_{\wsu(v)}$, then $$\mu^+_{\wsu (v)} (B^-(v,r))\leq e^\epsilon
\mu^+_{\wsu(v)} (B^-(v,1))$$ and if $B^+(v,1)$ meets the support of
$\mu^-_{\wss(v)}$, then
$$
\mu^-_{\wss(v)} (B^+(v,r))\leq e^\epsilon \mu^-_{\wss
  (v)}(B^+(v,1))\,.
$$ 
We say that $(\wt M,\Ga,\wt F)$ {\em has radius-H\"older-continuous
  strong stable/unstable ball masses}\index{radius-H\"older-continuous
  ball masses} if there exist $c\in\;]0,1]$ and $c'>0$ such that for
every $\epsilon\in\;]0,1]$, if $r>1$ is close enough to $1$, then for
every $v\in T^1\wt M$, if $B^-(v,1)$ meets the support of
$\mu^+_{\wsu(v)}$, then
$$
\mu^+_{\wsu (v)} (B^-(v,r))\leq e^{c'\epsilon^c} \mu^+_{\wsu(v)} (B^-(v,1))
$$ 
and if $B^+(v,1)$ meets the support of $\mu^-_{\wss(v)}$, then
$$
\mu^-_{\wss (v)}(B^+(v,r))\leq e^{c'\epsilon^c} \mu^-_{\wss(v)}(B^+(v,1))\,.
$$ 
Note that when $F=0$ and $M$ is locally symmetric with finite
volume, the conditional measures on the strong stable/unstable leaves
are homogeneous.  Hence $(\wt M,\Ga,\wt F)$ has
radius-H\"older-continuous strong stable/unstable ball masses.

When the sectional curvature of $\wt M$ has bounded derivatives and
when $(\wt M,\Ga,\wt F)$ has radius-H\"older-continuous strong
stable/unstable ball masses, we will prove a stronger statement: With
a constant $c_7>0$ and functions $\operatorname{O}(\cdot)$ independent
of $\ga$, for all $\eta\in\;]0,1]$ and $T\geq \lambda_\ga+
\operatorname{O}(\eta + e^{-\lambda_\ga/2})$, we have
\begin{equation}\label{eq:step5bis}
j_{\eta,\,\ga}(T)=
\Big(1+\operatorname{O}\Big(\frac{e^{-\lambda_\ga/2}}{2\eta}\Big)\Big)^2
e^{\operatorname{O}((\eta+e^{-\lambda_\ga/2})^{c_7})}\,.
\end{equation}
This stronger version will be needed for the error term estimate in
Section \ref{subsect:erroterms}. In order to obtain Theorem
\ref{theo:mainequidup}, only the fact that $j_{\eta,\,\ga}(T)$ tends
to $1$ as firstly $\lambda_\ga$ tends to $+\infty$, secondly $\eta$ tends
to $0$ is needed. A reader not interested in the error term
may skip many technical details below.

\medskip Given $a,b>0$ and a point $x$ in a metric space $X$ (with
$a,b,x$ depending on parameters), we will denote by
$B(x,a\,e^{\operatorname{O}(b)})$ any subset $Y$ of $X$ such that there
exists a constant $c>0$ (independent of the parameters) with
\begin{equation}\label{eq:defballbigO}
B(x,a\,e^{-c\,b})\subset  Y\subset B(x,a\,e^{c\,b})\,.
\end{equation}

Let $\eta\in\;]0,1]$ and $T\geq \lambda_\ga+\operatorname{O}(\eta+
e^{-\lambda_\ga/2})$.  In order to simplify the notation, let
$$
r_t=e^{-t/2}R, \;\; w^-_t=\flow{t/2}w^-\;\;\;{\rm and}\;\;\;
w^+_t=\flow{-t/2}w^+\,.
$$ 
By the definition of $j_{\eta,\,\ga}$, using the inequalities
\eqref{eq:estimhsmaletaprimp} where the constant $c_1$ is uniform
since $\wt F$ is bounded, and the fact that $r_t=
\bigO(e^{-\lambda_\ga/2})$ by Equation \eqref{eq:tpreslambdagamani},
we hence have
\begin{align}
j_{\eta,\,\ga}(T)&=
\iint_{(t,v)\in\A_{\eta,\ga}(T)} 
h^-_{\eta,\,r_t}(w^-_t)\;
h^+_{\eta,\,r_t}(w^+_t)\;dt\;d\wt m_F(v)\nonumber\\ & 
=\frac{e^{\operatorname{O}(e^{-\kappa_2\lambda_\ga/2})}}{(2\eta)^2}
\iint_{(t,v)\in\A_{\eta,\ga}(T)} 
\frac{dt\;d\wt m_F(v)}{{\scriptstyle
\mu_{\wss(w^-_t)}(B^+(w^-_t,\,r_t))\;
\mu_{\wsu(w^+_t)}(B^-(w^+_t,\,r_t))}}\,.\label{eq:adaptatJetaga}
\end{align}

We start the proof of Equation \eqref{eq:step5} by defining parameters
$s^+,s^-,s,v',v''$ associated with $(t,v)\in\A_{\eta,\ga}(T)$.  
\begin{center}
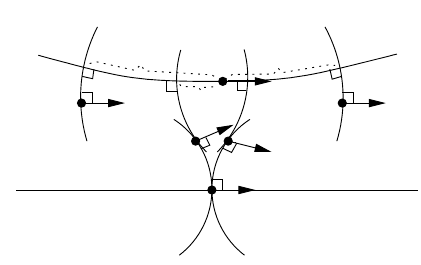
\end{center}

We have $(t,v)\in\A_{\eta,\ga}(T)$ if and only if there exist $s^\pm
\in\;]-\eta,\eta[\,$ such that
$$
\flow{\mp s^\mp}v\in B^\pm(\flow{\pm t/2} w^\mp,
e^{-t/2}R)\,.
$$ 

In order to define the parameters $s, v',v''$, we use the well known
local product structure of the unit tangent bundle in negative
curvature. If $v\in T^1M$ is close enough to $v^0_\ga$ (in particular,
$v_-\neq (v^0_\ga)_+$ and $v_+\neq (v^0_\ga)_-$), then let
$v'=f^+_{H\!B_-(v^0_\ga)}(v)$ be the unique element of $\wsu(v^0_\ga)$
such that $v'_+=v_+$, let $v''=f^-_{H\!B_+(v^0_\ga)}(v)$ be the unique
element of $\wss(v^0_\ga)$ such that $v''_-=v_-$, and let $s$ be the
unique element of $\RR$ such that $\flow{-s}v\in \wss(v')$. The map
$v\mapsto (s,v',v'')$ is a homeomorphism from a neighbourhood of
$v^0_\ga$ in $T^1\wt M$ to a neighbourhood of $(0,v^0_\ga,v^0_\ga)$ in
$\RR\times \wsu(v^0_\ga)\times \wss(v^0_\ga)$. Note that if
$v=\flow{r}v^0_\ga$ for some $r\in\RR$ close to $0$, then
$$
w^-=v^-_\ga,  \;w^+=v^+_\ga, \; s=r, \; v'=v''=v^0_\ga\,,
$$
$$
s^-=\frac{\lambda_\ga-t}{2}+s,\;s^+=\frac{\lambda_\ga-t}{2}-s\,.
$$
Up to increasing $t_0$ (which does not change Step 3M, up to increasing
$c_4$), we may assume that for every $(t,v)\in \A_{\eta,\ga}(T)$, the
vector $v$ belongs to the domain of this local product structure of
$T^1\wt M$ at $v^0_\ga$. 

The vectors $v,v',v''$ are close to $v^0_\ga$ if $t$ is large and $\eta$
small, as the following result shows. We denote (also) by $d$ the
Riemannian distance induced by Sasaki's metric on $T^1\wt M$.

\blemm \label{lem:vectproche} 
For every $(t,v)\in \A_{\eta,\ga}(T)$, we have $d(v,v^0_\ga),
d(v',v^0_\ga), d(v'',v^0_\ga)= \operatorname{O}(\eta+ e^{-t/2})$.
\elemm

\dem 
Consider the distance $d'$ on $T^1\wt M$, defined by
$$
\forall\;v_1,v_2\in T^1\wt M,\;\;\; d'(v_1,v_2) =
\max_{r\in[-1,0]} d\big(\pi(\flow{r}v_1),\pi(\flow{r}v_2)\big)\;.
$$ 
As seen in Claim (iii) of Step 3M, we have $d(\pi(w^\pm),\pi(v^\pm_\ga))$,
$d(\pi(v),\alpha_\ga)= \operatorname{O}(e^{-t/2})$, and furthermore,
$d(\pi(\flow{\pm t/2}w^\mp),\pi(v)), \frac{\lambda_\ga}{2}-\frac{t}{2}
=\operatorname{O}(\eta+e^{-t/2})$. Hence $d(\pi(v),\pi(v^0_\ga))=
\operatorname{O}(\eta+e^{-t/2})$. By Lemma \ref{lem:comparddHam}, 
we have
$$
d(\pi(\flow{-\frac{t}{2}-s^-}v),\pi(v^-_\ga))\leq 
d(\pi(\flow{-\frac{t}{2}-s^-}v),\pi(w^-))+
d(\pi(w^-),\pi(v^-_\ga))\leq 
R+c_0\,e^{-t/2}\;.
$$ 
By an exponential pinching argument, we hence have $d'(v,v^0_\ga)=
\operatorname{O}(\eta+ e^{-\lambda_\ga/2})$.  Since $d$ and $d'$ are
equivalent by Proposition \ref{prop:maniouholder},\footnote{In fact,
  Proposition \ref{prop:maniouholder} considers the distance
  $\delta_2(v,v')=\sup_{r\in[0,1]}d(\pi(\flow{r}v),\pi(\flow{r}v'))$
  instead of $d'$, but the argument is similar.} we therefore have
$d(v,v^0_\ga)= \operatorname{O}(\eta+ e^{-\lambda_\ga/2})$.

For all $w\in T^1\wt M$ and $V\in T_wT^1\wt M$, we may uniquely write
$V=V^{-}+V^0+V^{+}$ with $V^{-}\in T_w\wsu(w)$, $V^0\in
\RR\,\frac{d}{dt}\big|_{ t_0}\flow t w$ and $V^{+}\in T_w\wss(w)$.  By
\cite[Lem.~7.4]{PauPolSha15},\footnote{building on \cite{Brin95} whose
  compactness assumption on $M$ and torsion free assumption on $\Ga$
  are not necessary for this, the pinched negative curvature
  assumption is sufficient} Sasaki's metric (with norm $\|\cdot\|$) is
equivalent to the Riemannian metric with (product) norm
$$
\|V\|'=\sqrt{\|\,V^{-}\,\|^2+\|\,V^0\,\|^2+\|\,V^{+}\,\|^2}\,.
$$
By the dynamical local product structure of $T^1\wt M$ in the
neighbourhood of $v^0_\ga$ and by the definition of $v',v''$, the
result follows, since the exponential map of $T^1\wt M$ at $v^0_\ga$
is almost isometric close to $0$ and the projection to a factor of a
product norm is Lipschitz.  
\cqfd

\medskip
We now use the local product structure of the Gibbs measure to prove
the following result.

\blemm \label{lem:decompgibbs}
For every $(t,v)\in \A_{\eta,\ga}(T)$, we have
$$
 dt\,d\wt m_F(v) = 
e^{\operatorname{O}((\eta+e^{-\lambda_\ga/2})^{\kappa_2})}\;
dt\; ds\; d\musu{v^0_\ga}(v')\;d\muss{v^0_\ga}(v'')\,.
$$
\elemm

\dem By the definition of the measures (see Equations
\eqref{eq:defigibbs} and \eqref{eq:condtionelless}), since the above
parameter $s$ differs, when $v_-,v_+$ are fixed, only up to a constant
from the time parameter in Hopf's parametrisation with respect to the
basepoint $x_\ga=\pi(v^0_\ga)$, we have
\begin{align*}
d\wt m_{F}(v) & =
e^{C_{v_-}^-(x_\ga,\,\pi(v))\,+\,C^+_{v_+}(x_\ga,\,\pi(v))}\;
d\mu_{x_\ga}^-(v_-)\,d\mu^+_{x_\ga}(v_+)\,ds\\
d\musu{v^0_\ga}(v') & =
e^{C^+_{v'_{+}}(x_\ga,\,\pi(v'))}\,d\mu^+_{x_{\ga}}(v'_{+})\,,\\
d\muss{v^0_\ga}(v'') & =
e^{C^-_{v''_{-}}(x_\ga,\,\pi(v''))}\,d\mu^-_{x_{\ga}}(v''_{-})\;.
\end{align*}
By Proposition \ref{prop:continuGibbscocycle} \eqref{eq:cocycleLip}
since $F$ is bounded, we have $|\,C^{\pm}_\xi(z,z')\,|= \bigO
(d(z,z')^{\kappa_2})$ for all $\xi\in\partial_\infty\wt M$ and $z,z'\in\wt
M$ with $d(z,z')$ bounded.  Since the map $\pi:T^1\wt M\ra \wt M$ is
Lipschitz, and since $v_+=v'_+$ and $v_-=v''_-$, the result
follows from Lemma \ref{lem:vectproche}.  
\cqfd

\medskip When $\lambda_\ga$ is large, the submanifold
$\flow{\lambda_\ga/2} \Omega^-$ has a second order contact at
$v^0_\ga$ with $\wsu(v^0_\ga)$ and similarly, $\flow{-\lambda_\ga/2}
\Omega^+$ has a second order contact at $v^0_\ga$ with $\wss
(v^0_\ga)$. Let $P_\ga$ be the plane domain of $(t,s)\in\RR^2$ such
that $|\lambda_\ga-t|\leq 2\eta +c_0e^{-t/2}$ and there exist
$s^\pm\in\;]-\eta,\eta[$ with $s^\mp= \frac{\lambda_\ga-t}{2} \pm s +
\operatorname{O}(e^{-\lambda_\ga/2})$.  Note that its area is $(2\eta
+ \operatorname{O}(e^{-\lambda_\ga/2}))^2$.  By the above, we
have\footnote{with the obvious meaning of a double inclusion by
  Equation \eqref{eq:defballbigO}}
$$
\A_{\eta,\ga}(T)=P_\ga\times 
B^-(v^0_\ga,r_{\lambda_\ga}\,e^{\operatorname{O}(\eta+e^{-\lambda_\ga/2})})
\times 
B^+(v^0_\ga,r_{\lambda_\ga}\,e^{\operatorname{O}(\eta+e^{-\lambda_\ga/2})})\,.
$$
By Lemma \ref{lem:decompgibbs}, we hence have
\begin{align}
\int_{\A_{\eta,\ga}(T)} dt\,d\wt m_F(v) &= 
e^{\operatorname{O}((\eta+e^{-\lambda_\ga/2})^{\kappa_2})}\;
(2\eta + \operatorname{O}(e^{-\lambda_\ga/2}))^2
\;\times \nonumber \\  &
\musu{v^0_\ga}(B^-(v^0_\ga,r_{\lambda_\ga}\,e^{\operatorname{O}(\eta+e^{-\lambda_\ga/2})}))\;
\muss{v^0_\ga}(B^+(v^0_\ga,r_{\lambda_\ga}\,e^{\operatorname{O}(\eta+e^{-\lambda_\ga/2})}))\,.
\label{eq:volAetaga}
\end{align}

The last ingredient of the proof of Step 4M is the following
continuity property of the masses of balls in the strong stable and
strong unstable manifolds as their centre varies.  This result
generalises \cite[Lem.~11]{ParPau16ETDS}.  The precise control for the
error term is used in Section \ref{subsect:erroterms}.

\blemm\label{lem:variaboul} Assume that $(\wt M,\Ga,\wt F)$ has
radius-continuous strong stable/unstable ball mas\-ses. There exists
$c_5>0$ such that for every $\epsilon>0$, if $\eta$ is small enough
and $\lambda_\ga$ large enough, then for every $(t,v)\in \A_{\eta,\ga}(T)$,
we have
$$
\musu{w^+_t}(B^-(w^+_t,r_t))=e^{\operatorname{O}(\epsilon^{c_5})}\;
\musu{v^0_\ga}(B^-(v^0_\ga,r_{\lambda_\ga}))
$$
and 
$$
\muss{w^-_t}(B^+(w^-_t,r_t))=e^{\operatorname{O}(\epsilon^{c_5})}\;
\muss{v^0_\ga}(B^+(v^0_\ga,r_{\lambda_\ga}))\,.
$$
If we furthermore assume that the sectional curvature of $\wt M$ has
bounded derivatives and that $(\wt M,\Ga,\wt F)$ has
radius-H\"older-continuous strong stable/unstable ball masses, then we
may replace $\epsilon$ by $(\eta+ e^{-\lambda_\ga/2})^{c_6}$ for some
constant $c_6>0$.  
\elemm

\dem We prove the (second) claim for $W^{+}$, the (first) one
for $W^{-}$ follows similarly. The final statement is only used
for the error estimates in Section \ref{subsect:erroterms}.

\medskip
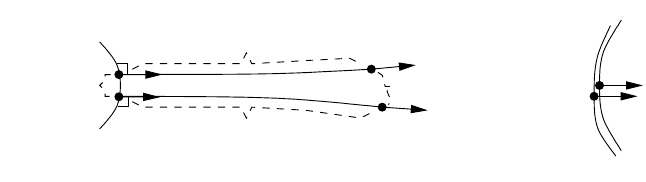

\medskip
Using respectively Equation \eqref{eq:flotVetaeta} since
$w^-_t=\flow{t/2}w^-$ and $r_t=e^{-t/2}R$, Equation
\eqref{eq:mussscaling} where $(\ell,t,w)$ is replaced by
$(v,t/2,w^-)$, and Equation \eqref{eq:changemoinsplus}, we have
\begin{align}
\muss{w^-_t}(B^+(w^-_t,r_t))&
=\int_{v\in B^+(w^-,\,R)} d\muss{\flow{t/2}w^-}(\flow{t/2}v)\nonumber
\\ & =\int_{v\in B^+(w^-,\,R)}e^{C^-_{v_-}(\pi(v),\;\pi(\flow{t/2}v))}
\;d\muss{w^-}(v)\nonumber
\\ & =\int_{v\in B^+(w^-,\,R)}e^{\int_{\pi(v)}^{\pi(\flow{t/2}v)}(\wt F-\delta)}
\;d\muss{w^-}(v)\,.\label{eq:boulboul}
\end{align}
Similarly, for every $a>0$, we have
\begin{equation}\label{eq:boulboulbis}
\muss{v^0_\ga}(B^+(v^0_\ga,ar_t))= 
\int_{v\in B^+(v^-_\ga,\,aR)}e^{\int_{\pi(v)}^{\pi(\flow{t/2}v)}(\wt F-\delta)}
\;d\muss{v^-_\ga}(v)\;.
\end{equation}

Let $h^-:B^+(w^-,R)\ra \wss(v^-_\ga)$ be the map such that
$(h^-(v))_-=v_-$, which is well defined and a homeomorphism onto its
image if $\lambda_\ga$ is large enough (since $R$ is fixed). By
Proposition \ref{prop:abscontskinmeas} applied with $D=H\!B_+(w^-)$ and
$D'=H\!B_+(v^-_\ga)$, we have, for every $v\in B^+(w^-,R)$,
$$
d\muss{w^-}(v)= e^{-C^-_{v_-}(\pi(v),\;\pi(h^-(v)))}\;d\muss{v^-_\ga}(h^-(v))\,.
$$

Let us fix $\epsilon>0$. The strong stable balls of radius $R$ centred
at $w^-$ and $v^-_\ga$ are very close (see the above picture). More
precisely, recall that $R$ is fixed, and that 
$$
d(\pi(w^-),\pi(v^-_\ga)) = \operatorname{O}(e^{-\lambda_\ga/2})
\;\;\;{\rm and}\;\;\;
d(\pi(\flow{t/2}w^-), \pi(\flow{\lambda_\ga/2}v^-_\ga))=
\operatorname{O}(\eta +e^{-\lambda_\ga/2})\;.
$$ 
Therefore we have $d(\pi(v),\pi(h^-(v)))\leq \epsilon$ for every
$v\in B^+(w^-,R)$ if $\eta$ is small enough and $\lambda_\ga$ is large
enough. If we furthermore assume that the sectional curvature has
bounded derivatives, then by Anosov's arguments, the strong stable
foliation is H\"older-continuous, see for instance
\cite[Theo.~7.3]{PauPolSha15}.  Hence we have $d(\pi(v),\pi(h^-(v)))
=\operatorname{O}((\eta+e^{-\lambda_\ga/2})^{c_6})$ for every $v\in
B^+(w^-,R)$, for some constant $c_6>0$, under the additional
regularity assumption on the curvature.  We also have
$h^-(B^+(w^-,R))= B^+(v^-_\ga, R\,e^{\operatorname{O}(\epsilon)})$
and, under the additional hypothesis on the curvature,
$h^-(B^+(w^-,R))= B^+(v^-_\ga,
R\,e^{\operatorname{O}((\eta+e^{-\lambda_\ga/2})^{c_6})})$.

In what follows, we assume that $\epsilon=
(\eta+e^{-\lambda_\ga/2})^{c_6}$ under the additional assumption on
the curvature. By Proposition \ref{prop:continuGibbscocycle}
\eqref{eq:cocycleLip} since $\wt F$ is bounded, we hence have, for
every $v\in B^+(w^-,R)$,
$$
d\muss{w^-}(v)= e^{\operatorname{O}(\epsilon^{\kappa_2})}\;
d\muss{v^-_\ga}(h^-(v))
$$
and, using the \eqref{eq:HC}-property and the boundedness of $\wt F$,
$$
\int_{\pi(v)}^{\pi(\flow{t/2}v)}(\wt F-\delta)-
\int_{\pi(h^-(v))}^{\pi(\flow{t/2}h^-(v))}(\wt F-\delta)= 
\operatorname{O}(\epsilon^{\kappa_2})\;.
$$
The result follows, by Equation \eqref{eq:boulboul} and
\eqref{eq:boulboulbis} and the continuity properties in the
radius of the strong stable/unstable ball masses. 
\cqfd

\medskip Now Lemma \ref{lem:variaboul} (with $\epsilon$ as in its
statement, and when its hypotheses are satisfied) implies that
\begin{align*}
\iint_{(t,v)\in\A_{\eta,\ga}(T)} \;
& \frac{dt\;d\wt m_F(v)}{\mu^-_{\wss(w^-_t)}(B^+(w^-_t,\,r_t))\;
 \mu^+_{\wsu(w^+_t)}(B^-(w^+_t,\,r_t))}\\
=\; & \frac{e^{\operatorname{O}(\epsilon^{c_5})}\iint_{(t,v)\in\A_{\eta,\ga}(T)} 
dt\;d\wt m_F(v)}{\mu^-_{\wss(v^0_\ga)}(B^+(v^0_\ga,\,r_t))\;
\mu^+_{\wsu(v^0_\ga)}(B^-(v^0_\ga,\,r_t))}\,.
\end{align*}
By Equation \eqref{eq:adaptatJetaga} and Equation
\eqref{eq:volAetaga}, we hence have
$$
j_{\eta,\,\ga}(T)=e^{\operatorname{O}((\eta+e^{-\lambda_\ga/2})^{\kappa_2})}
\;e^{\operatorname{O}(\epsilon^{c_5})}\;
\frac{(2\eta + \operatorname{O}(e^{-\lambda_\ga/2}))^2}{(2\eta)^2} 
$$ 
under the technical assumptions of Lemma \ref{lem:variaboul}. The
assumption on radius-continuity of strong stable/unstable ball
mas\-ses can be bypassed using bump functions, as explained in
\cite[page 81]{Roblin03}. This completes the proof of Equation
\eqref{eq:step5}, hence the proof of Theorem \ref{theo:mainequidup}.
\cqfd

\section{Equidistribution of common perpendiculars 
in simplicial trees}
\label{subsec:equidcommperpdiscrtime}

In this Section, we prove a version of Theorem \ref{theo:mainequidup}
for the discrete time geodesic flow on simplicial trees
(and we leave to the reader the version without the assumption that
the critical exponent of the system of conductances is positive).

Let $\XX$, $X$, $x_0$, $\Ga$, $\wt c$, $c$, $\wt F_c$, $F_c$,
$\delta_c<+\infty$, $(\mu^\pm_x)_{x\in V\XX}$, $\wt m_c$, $m_c$ be as
in the beginning of Section \ref{subsec:mixingratesimpgraphs}.  Let
$\D^-= (\DD^-_i)_{i\in I^-}$ and $\D^+=(\DD^+_j)_{j\in I^+}$ be
locally finite $\Ga$-equivariant families of nonempty proper
simplicial subtrees of $\XX$. We denote by $D^\pm_k=|\DD^\pm_k|_1$ the
geometric realisation of $\DD^\pm_k$ for $k\in I^\pm$.

For every edge path $\alpha=(e_1,\dots, e_n)$ in $\XX$, we set
$$
\wt c(\alpha)= \sum_{i=1}^n\wt c(e_i)\;.
$$

\btheo\label{theo:discretemainequidup} Assume that the critical 
exponent $\delta_c$ of $\wt c$ is positive and that the Gibbs measure
$m_{c}$ is finite and mixing for the discrete time geodesic flow on
$\Ga\bs\G\XX$.  Then
$$
\lim_{t\ra+\infty} \;\frac{e^{\delta_c}-1}{e^{\delta_c}}\;
\|m_{c}\|\;e^{-\delta_c\, t}
\sum_{\substack{i\in I^-/_\sim,\; j\in I^+/_\sim, \;\ga\in\Ga\\ 
{D^-_i}\cap {D^+_{\ga j}}=\emptyset,\; 
\lambda_{i,\,\ga j}\leq t}} \;e^{\wt c(\alpha_{i,\ga j})}\; 
\Dirac_{\alpha^-_{i,\,\ga j}} \otimes\Dirac_{\alpha^+_{\ga^{-1}i,\,j}}
\;=\; \wt\sigma^+_{\D^-}\otimes \wt\sigma^-_{\D^+}\;
$$
for the weak-star  convergence of measures on the locally
compact space $\gengeod \XX\times \gengeod \XX$.  
\etheo

\dem The proof is a modification of the continuous time proof for
metric trees in Sections \ref{subsec:equidcommperpcontI} and
\ref{subsec:equidcommperpcontII}. Here, we indicate the changes to
adapt the proof to the discrete time.  We use the conventions for the
discrete time geodesic flow described in Section \ref{subsec:trees}.

Note that for all $i\in I^-$, $j\in I^+$, $\ga\in\Ga$, the common
perpendicular $\alpha_{i,\ga j}$ is now an edge path from $D^-_i$ to
$D^+_{\ga j}$, and that by Proposition \ref{prop:integpotconduct}, we
have 
$$
\int_{\alpha_{i,\ga j}} \wt F_c=\wt c(\alpha_{i,\ga j})\;.
$$

In the definition of the bump functions in Section
\ref{subsec:skinningwithpot}, we assume (as we may) that $\eta <1$, so
that for all $\eta'\in\;]0,1[$ and $w\in\normalmp D^\pm$ such that
$w_\mp\in \Lambda\Ga$, we have 
$$
V^\pm_{w,\eta,\eta'}=B^\pm(w,\eta')\,,
$$ 
see Equation \eqref{eq:defidynamicalneighbourhood} and recall that we
are only considering discrete geodesic lines.  As $\ell(0)=w(0)$ for
every $\ell\in B^\pm(w,\eta')$ since $\eta'<1$, and as the time is now
discrete, Equations \eqref{eq:defihpm} and \eqref{eq:definumoins} give
\begin{equation}\label{eq:newhpmvalue}
h^\mp_{\eta,\,\eta'}(w)=\frac{1}{\mu_{W^\pm(w)}(B^\pm(w,\eta'))}\;.
\end{equation}
This is a considerable simplification compared with the inequalities
of Equation \eqref{eq:estimhsmaletaprimp}.

In the whole proof, we restrict to $t=n\in\NN$, $T=N\in\NN$. We keep
the notation of Equation \eqref{eq:notationstep1}, as well as the only
assumptions on the Borel sets $\Omega^\pm\subset \normalmp D^\pm$ to
have finite positive skinning measure, with boundary of zero skinning
measure. In Steps 1 and 2, we define instead of Equation
\eqref{eq:defiIOmegapm}
$$
I_{\Omega^-,\,\Omega^+}(N)=\;(e^{\delta_c}-1)\;\|m_{c}\|\;e^{-\delta_c(N+1)}
\sum_{\substack{\ga\in\Ga\,:\;0<\lambda_{\ga }\leq N\\
\alpha^-_\ga|_{]0,\lambda_\ga]}\in\Omega^-|_{]0,\lambda_\ga]},\;\;
\alpha^+_\ga|_{]-\lambda_\ga,0]}\in\Omega^+|_{]-\lambda_\ga,0]}}}
\;e^{\int_{\alpha_{\ga}}\wt F_c}\;,
$$
and instead of Equation \eqref{eq:defiaetat}
$$
a_\eta(n)=\sum_{\ga\in\Ga}\;
\int_{\ell\in\G \XX}\phi^-_\eta(\flow{-\lfloor n/2\rfloor}\ell)\;
\phi^+_\eta(\flow{\lceil n/2\rceil}\ga^{-1}\ell)\;d\wt m_{c}(\ell)\,.
$$
Equation \eqref{eq:defiIetapmT} is replaced by 
$$
i_{\eta}(N)= \sum_{n=0}^{N} e^{\delta_c\,n} \; a_\eta(n) \,,
$$
so that by a geometric sum argument, the pair of inequalities 
\eqref{eq:firstmajoI?????} becomes
$$
e^{-\epsilon}\; \frac{e^{\delta_c\,(N+1)}\,\wt \sigma^+(\Omega^-)\;
\wt \sigma^-(\Omega^+)}{(e^{\delta_c}-1)\,\|m_{c}\|}\;
-c_\epsilon 
\leq\; i_{\eta}(N)\leq\; 
e^{\epsilon}\; \frac{e^{\delta_c\,(N+1)}\,\wt \sigma^+(\Omega^-)\;
\wt \sigma^-(\Omega^+)}{(e^{\delta_c}-1)\,\|m_{c}\|}\;
+c_\epsilon\,.
$$

Step 3 is unchanged up to replacing $\int_0^T$ by $\sum_{n=0}^N$, $\wt
F$ by $\wt F_c$, $\delta$ by $\delta_c$ and $t/2$ by either $\lfloor
n/2\rfloor$ or $\lceil n/2 \rceil$, so that Equation
\eqref{eq:commonuptohere} becomes, since $\lfloor n/2\rfloor+\lceil
n/2 \rceil=n$,
\begin{align*}
h^-_{\eta,\,R}\circ & f^+_{D^-}(\ell)\;
h^+_{\eta,\,R}\circ f^-_{D^+}(\ga^{-1}\ell)\nonumber\\ =\;& 
e^{-\delta_c\,n}\;
e^{\int_{w^-(0)}^{w^-(\lfloor n/2\rfloor)} \wt F_c+
\int_{w^+(-\lceil n/2 \rceil)}^{w^+(0)} \wt F_c}\;
h^-_{\eta,\,e^{-\lfloor n/2\rfloor}R}(\flow{\lfloor n/2\rfloor}w^-)\,
h^+_{\eta,\,e^{-\lceil n/2 \rceil}R}(\flow{-\lceil n/2 \rceil}w^+)\,.
\end{align*}
The proof then follows in the same way as in Section
\ref{subsec:equidcommperpcontII}, with the simplifications in the
point (iii) that, taking $\eta<1/2$, we have $\lambda_\ga$ equal to
$t=n$, and the points $w^- (\lfloor \frac n2\rfloor)$, $w^+(-\lceil
\frac n2\rceil)$ and $\ell(0)$ are equal. In particular, Equation
\eqref{eq:controlpieceFlengthT} simplifies as
$$
e^{\int_{w^-(0)}^{w^-(\lfloor \frac t2\rfloor)} \wt F_c+ \int_{w^+(-\lceil
\frac t2\rceil)}^{w^+(0)} \wt F_c}=e^{\int_{\alpha_\ga}\wt F_c}\;,
$$
thus avoiding the assumption that $F_c$ (or equivalently $c$, see
Section \ref{subsec:potentials}) is bounded. We now define
$$
\A_{\eta,\ga}(N)=\big\{(n,\ell)\in [t_0,N] \times\G \XX:
\ell\in \V^+_{\eta,\,e^{-\lfloor \frac n2\rfloor}R}
(\flow{\lfloor \frac n2\rfloor}\Omega^-)\cap
\V^-_{\eta,\,e^{-\lceil\frac n2\rceil}R}
(\ga\flow{-\lceil\frac n2\rceil}\Omega^+)\big\}\;.
$$
The end of Step 3T simplifies as
$$
-\,c_4\leq\; i_{\eta}(N) -\sum_{\substack{\ga\in\Ga\,:\;
t_0+2\leq\lambda_{\ga}\leq N\\
\alpha^-_\ga|_{[0,\,\lambda_\ga]}\in\Omega^-|_{[0,\,\lambda_\ga]}, \;\;
\alpha^+_\ga|_{[-\lambda_\ga,\,0]}\in\Omega^+|_{[-\lambda_\ga,\,0]}}}
\;e^{\int_{\alpha_\ga}\wt F_c}\;j_{\eta,\,\ga}(N)\;\leq \,c_4\;.
$$

The statement of Step 4T now simplifies as
$$
j_{\eta,\ga}(N)=1\;,
$$ 
if $\eta<\frac12$, and if $\ga\in\Ga$ is such that $D^-$ and $\ga D^+$ do
not intersect and $\lambda_\ga$ is large enough. We introduce in its
proof the slightly modified notation
$$
r_n^-=e^{-\lfloor \frac n2\rfloor}R, \;\; r_n^+=e^{-\lceil \frac n2\rceil}R, \;\; 
w^-_n=\flow{\lfloor \frac n2\rfloor}w^-\;\;\;{\rm and}\;\;\;
w^+_n=\flow{-\lceil \frac n2\rceil}w^+\,.
$$ 
and we now take as $x_\ga$ the point at distance $\lfloor \frac
n2\rfloor$ from its origin on the common perpendicular $\alpha_\ga$.
Equation \eqref{eq:adaptatJetagaT} becomes (using Equation
\eqref{eq:newhpmvalue} instead of Equation
\eqref{eq:estimhsmaletaprimp})
$$
j_{\eta,\,\ga}(N)=\iint_{(n,\,\ell)\in\A_{\eta,\ga}(N)} 
\frac{dn\;d\wt m_{c}(\ell)}{{
\mu_{\wss(w^-_n)}(B^+(w^-_n,\,r_n^-))\;
\mu_{\wsu(w^+_n)}(B^-(w^+_n,\,r_n^+))}}\,.
$$
Since $\ell(0)=x_\ga$ if $(n,\ell)\in \A_{\eta,\,\ga}(N)$, Equation
\eqref{eq:gibbseta} simplifies as
$$
d\wt m_{c}(\ell)=d\mu_{x_\ga}^-(\ell_-)\,d\mu^+_{x_\ga}(\ell_+)\,ds\;,
$$ 
with $ds$ the counting measure on the Hopf parameter $s\in\ZZ$ of
$\ell$ (with basepoint $x_\ga$). If $\eta<\frac{1}{2}$, replacing
$P_\ga$ with its intersection with $\ZZ^2$ reduces it to one point
$(\lambda_\ga,0)$, and now $s=s^\pm=0$. Lemma \ref{lem:passeralinfini}
becomes
$$
(B^+(w^-_n,r_n^-))_- = 
B_{d_{x_\ga}}(\xi^-_\ga, R\,e^{-\lfloor\frac{\lambda_\ga}2\rfloor}),\;\;\;
(B^-(w^+_n,r_n^+))_+ = 
B_{d_{x_\ga}}(\xi^+_\ga, R\,e^{-\lceil\frac{\lambda_\ga}2\rceil})\,,
$$
so that
$$
\A_{\eta,\ga}(N)=\{(\lambda_\ga,0)\}\times 
B_{d_{x_\ga}}(\xi_\ga^-,R\,e^{-\lfloor\frac{\lambda_\ga}2\rfloor})
\times 
B_{d_{x_\ga}}(\xi_\ga^+,R\,e^{-\lceil\frac{\lambda_\ga}2\rceil})\,.
$$
Finally, since $\ell(0)= x_\ga$ if $(n,\ell)\in \A_{\eta,\,\ga}(N)$,
Equation \eqref{eq:penultimatePart4T} becomes
\begin{align*}
\mu_{W^+(w^-_n)}(B^+(w^-_n,r_n^-))&= 
\mu^-_{x_\ga}(B_{d_{x_\ga}}(\xi^-_\ga,R\,e^{-\lfloor\frac{\lambda_\ga}2\rfloor})),\\
\mu_{W^-(w^+_t)}(B^-(w^+_n,r_n^+))&= 
\mu^+_{x_\ga}(B_{d_{x_\ga}}(\xi^+_\ga,R\,e^{-\lceil\frac{\lambda_\ga}2\rceil}))\,.
\end{align*}
The last centred equation in Step 4T now reduces to $j_{\eta,\,\ga}(T)=1$.
\cqfd

\medskip 
For lattices in regular trees, we get more explicit expressions.

\bcoro \label{coro:twosubtrees} Let $\XX$ be a $(q+1)$-regular
simplicial tree (with $q\geq 2$) and let $\Ga$ be a lattice of $\XX$
such that $\Ga\bs\XX$ is not bipartite. Assume that the Patterson
density is normalised to be a family of probability measures. Let
$\DD^\pm$ be nonempty proper simplicial subtrees of $\XX$ with
stabilisers $\Ga_{\DD\pm}$ in $\Ga$, such that $\D^\pm=
(\ga\DD^\pm)_{\ga\in\Ga/\Ga_{\DD^\pm}}$ is locally finite. Let
$\sigma^\mp_{\D^\pm}$ be their skinning measures for the sero system
of conductances. Then
$$
\lim_{t\ra+\infty} \;\frac{q-1}{q+1}\; \Vol(\Ga\dbs\XX) \; q^{-t}
\sum_{\substack{(\alpha,\,\beta,\,\ga)\in
\Ga/\Ga_{\DD^-}\times\Ga/\Ga_{\DD^+}\times\Ga\\
 0< d(\alpha\DD^-, \,\ga\beta\DD^+)\leq t}}
\Dirac_{\alpha^-_{\alpha,\,\ga\beta}} \otimes
\Dirac_{\alpha^+_{\ga^{-1}\alpha,\,\beta}}\;=\;
\wt\sigma^+_{\D^-}\otimes \wt\sigma^-_{\D^+}\;,
$$
for the weak-star convergence of measures on the locally
compact space $\gengeod \XX\times \gengeod \XX$.

If the measure $\sigma^-_{\D^+}$ is nonzero and finite, then
$$
\lim_{t\ra+\infty} \;\frac{q-1}{q+1}\; 
\frac{\Vol(\Ga\dbs\XX)}{\|\sigma^-_{\D^+}\|} 
\; q^{-t} \sum_{\ga\in \Ga/\Ga_{\DD^+},\; 0< d(\DD^-, \, \ga \DD^+)\leq t}
\Dirac_{\alpha^-_{e,\,\ga}} \;=\;
\wt\sigma^+_{\DD^-}\;,
$$
for the weak-star convergence of measures on the locally
compact space $\gengeod \XX$. 
\ecoro

\dem In order to prove the first claim, we apply Theorem
\ref{theo:discretemainequidup} with $\wt c\equiv 0$, so that by
Proposition \ref{prop:uniflatmBMfinie}, Theorem
\ref{theo:uniflatmBMmixing}, and Proposition \ref{prop:computBM} (3),
we have $\delta_c=\ln q>0$, $m_{c}=m_{\rm BM}$ is finite and mixing,
and $\|m_{\rm BM}\|= \frac{q}{q+1} \Vol(\Ga\dbs\XX)$.

The second claim follows by restricting to $\alpha=\beta=e$ and
integrating on an appropriate fundamental domain (note that Equation
\eqref{eq:reducnarrowup} does not require $\Omega^+$ to be relatively
compact, just to have finite measure for $\wt\sigma^-$).  
\cqfd

\medskip 
The mixing assumption in Theorem \ref{theo:discretemainequidup}
implies that the length spectrum $L_\Ga$ of $\Ga$ is equal to $\ZZ$.
The next result considers the other case, when only the square of the
geodesic flow is mixing, while appropriately restricted. Note that the
smallest nonempty $\Ga$-invariant simplicial subtree of $\XX$ is
uniform, without vertices of degree $2$, for instance in the case when
$\XX$ is $(p+1,q+1)$-biregular with $p,q\geq 2$ and $\Ga$ is a lattice
of $\XX$.

\btheo\label{theo:discretemainequidupeven} Assume that the smallest
nonempty $\Ga$-invariant simplicial subtree of $\XX$ is uniform,
without vertices of degree $2$, and that the length spectrum $L_\Ga$
of $\Ga$ is $2\ZZ$.  Assume that the critical exponent $\delta_c$ of
$\wt c$ is positive, that the Gibbs measure $m_{c}$ is finite and
that its restriction to $\Ga\bs\Geven\XX$ is mixing for the square of
the discrete time geodesic flow on $\Ga\bs\Geven\XX$. Then
\begin{align*}
\lim_{t\ra+\infty} \;\frac{e^{2\,\delta_c}-1}{2\;e^{2\,\delta_c}}\;
\|m_{c}\|\;e^{-\delta_c\, t}
 \sum_{\substack{i\in I^-/_\sim,\; j\in I^+/_\sim, \;\ga\in\Ga\\ 
{D^-_i}\cap {D^+_{\ga j}}=\emptyset,\; 
\lambda_{i,\,\ga j}\leq t}} \;e^{\wt c(\alpha_{i,\ga j})}\; 
& \Dirac_{\alpha^-_{i,\,\ga j}} \otimes\Dirac_{\alpha^+_{\ga^{-1}i,\,j}}\\
\;=\;\;& \wt\sigma^+_{\D^-}\otimes \wt\sigma^-_{\D^+}\;
\end{align*}
for the weak-star convergence of measures on the locally
compact space $\gengeod \XX\times \gengeod \XX$.  
\etheo

\dem We denote by $\wt\sigma^\pm_{\D^\mp,\,{\rm even}}$ the
restriction of $\wt\sigma^\pm_{\D^\mp}$ to $\gengeodeven\XX$, and by
$\wt\sigma^\pm_{\D^\mp,\,{\rm odd}}$ the restriction of
$\wt\sigma^\pm_{\D^\mp}$ to $\gengeododd\XX= \gengeod\XX- \gengeodeven
\XX$.  We denote by $V_{\rm even}\XX$ the subset of $V\XX$ consisting
of the vertices at even distance from $x_0$, and by $V_{\rm odd} \XX=
V\XX-V_{\rm even}\XX$ its complement. The subsets $V_{\rm even}\XX$
and $V_{\rm odd}\XX$ are $\Ga$-invariant if $L_\Ga=2\ZZ$ by Equation
\eqref{eq:equivLamGaevendisteven}.

Let us first prove that
\begin{align}
\lim_{t\ra+\infty} \;\frac{e^{2\,\delta_c}-1}{2\;e^{2\,\delta_c}}\;
\|m_{c}\|\;e^{-\delta_c\, t}
 \sum_{\substack{i\in I^-/_\sim,\; j\in I^+/_\sim, \;\ga\in\Ga\\ 
\pi(\alpha^-_{i,\ga j}),\,\pi(\alpha^+_{\ga^{-1}i,j})\,\in\, V_{\rm even}\XX\\
{D^-_i}\cap {D^+_{\ga j}}=\emptyset,\; 
\lambda_{i,\,\ga j}\leq t}} & \;e^{\wt c(\alpha_{i,\ga j})}\; 
\Dirac_{\alpha^-_{i,\,\ga j}} \otimes\Dirac_{\alpha^+_{\ga^{-1}i,\,j}}\nonumber\\
\;=\;\;& \wt\sigma^+_{\D^-,\,{\rm even}}\otimes \wt\sigma^-_{\D^+,\,{\rm even}}
\;\label{eq:simplicialevenevencase}
\end{align}
for the weak-star convergence of measures on the locally
compact space $\gengeodeven \XX\times \gengeodeven \XX$.

\medskip The proof of this Equation \eqref{eq:simplicialevenevencase}
is a modification of the proof of the previous Theorem
\ref{theo:discretemainequidup}.  We now restrict to $t=2n\in\NN$,
$T=2N\in\NN$, and we replace $\wt m_{c}$ by $(\wt m_{c})\big|_{
  \Geven \XX}$ and $(\flow{t})_{t\in\ZZ}$ by $(\flow{2t})_{t\in\ZZ}$. Note
that since $\wt m_{c}$ is invariant under the time $1$ of the
geodesic flow, which maps $\Ga\bs \Geven \XX$ to $\Ga\bs \G\XX-\Ga\bs
\Geven \XX$, we have
\begin{equation}\label{eq:halfGibbsmass}
\big\|(m_{c})\big|_{\Ga\bs\Geven \XX}\big\|=\;\frac12\;\|m_{c}\|\;.
\end{equation}

Note that for all $i\in I^-$, $j\in I^+$ and $\ga\in\Ga$, if
$\pi(\alpha^-_{i,\ga j})$ and $\pi(\alpha^+_{\ga^{-1}i,j})$ belong to
$V_{\rm even}\XX$, then the distance between $D_i^-$ and $\ga D_j^+$
is even.\footnote{Indeed, for all $x,y,z$ in a simplicial tree, if $p$
  is the closest point to $x$ on $[y,z]$, then $d(y,z)= d(y,x)+d(x,z)-
  2\;d(x,p)$.}

In Steps 1 and 2, we now consider $\Omega^\pm$ two Borel subsets of
$\normalmp D^\pm\cap\gengeodeven\XX$, and we define instead of Equation
\eqref{eq:defiIOmegapm}
\begin{align*}
I_{\Omega^-,\,\Omega^+}(2N)=&
\;(e^{2\,\delta_c}-1)\;\frac{\|m_{c}\|}{2}\;e^{-2\,\delta_c(N+1)}\;\times\\ &
\sum_{\substack{\ga\in\Ga\,:\;0<\lambda_{\ga }\leq 2N,\;\;
\pi(\alpha^-_\ga),\,\pi(\alpha^+_\ga)\,\in\, V_{\rm even}\XX\\
\alpha^-_\ga|_{]0,\lambda_\ga]}\in\Omega^-|_{]0,\lambda_\ga]},\;\;
\alpha^+_\ga|_{]-\lambda_\ga,0]}\in\Omega^+|_{]-\lambda_\ga,0]}}}
\;e^{\int_{\alpha_{\ga}}\wt F_c}\;,
\end{align*}
and instead of Equation \eqref{eq:defiaetat}
$$
a_\eta(2n)=\sum_{\ga\in\Ga}\;
\int_{\ell\in\Geven \XX}\phi^-_\eta(\flow{-2\lfloor n/2\rfloor}\ell)\;
\phi^+_\eta(\flow{2\lceil n/2\rceil}\ga^{-1}\ell)\;d\wt m_{c}(\ell)\,.
$$
Equation \eqref{eq:defiIetapmT} is replaced by 
$$
i_{\eta}(2N)= \sum_{n=0}^{N} e^{\delta_c\,2n} \; a_\eta(2n) \,.
$$
The mixing property of the square of the geodesic flow on
$\Ga\backslash\Geven \XX$ for the restriction of the Gibbs measure
$m_{c}$ gives that, for every $\epsilon>0$, there exists
$T_\epsilon=T_{\epsilon,\eta} \geq 0$ such that for all $n\geq
T_\epsilon $, we have
\begin{align*}
&\frac{e^{-\epsilon}\,\int_{\Geven \XX} \phi^-_\eta\;
d\wt m_{c}\int_{\Geven \XX} \phi^+_\eta\;d\wt m_{c}}
{\|(m_{c})\big|_{\Ga\bs\Geven \XX}\|}
\\&\;\;\;\;\;\;\;\;\;\;\;\;\;\;\leq\;a_\eta(2n) 
\leq\;\frac{e^\epsilon\,\int_{\Geven \XX} \phi^-_\eta\;
d\wt m_{c}\int_{\Geven \XX} \phi^+_\eta\;d\wt m_{c}}
{\|(m_{c})\big|_{\Ga\bs\Geven \XX}\|}\,.
\end{align*}
Note that $\Geven \XX$ is saturated by the strong stable and strong
unstable leaves, since two points $x,y$ on a given horosphere of
centre $\xi\in \partial_\infty X$ are at even distance one from
another (equal to $2d(x,p)$ where $[x,\xi[\;\cap\,[y,\xi[\;=
[p,\xi[\,$).  By the disintegration statement in Proposition
\ref{prop:disintegration}, when $\ell$ ranges over $\U^\pm_D\cap
\Geven\XX$, we have
$$
d\wt m_{c}|_{\U^\pm_D\cap\, \Geven\XX}(\ell)=
\int_{\rho\in \partial^1_{\pm} D\cap\, \gengeodeven\XX}
d\nu^\mp_\rho(\ell)\, d\wt\sigma^\pm_{D}(\rho)\;.\
$$
Hence the proof of Lemma \ref{lem:integrable} extends to give
\begin{equation}\label{eq:lemmintegrableven}
\int_{\Geven \XX} \phi^\mp_\eta\;d\wt m_{c}=
\wt \sigma^\pm_{\rm even}(\Omega^\mp)\;,
\end{equation}
where in order to simplify notation $\wt \sigma^\pm_{\rm even}=\wt
\sigma^\pm_{D^\mp,\,{\rm even}}$.

Therefore, by Equations \eqref{eq:halfGibbsmass} and
\eqref{eq:lemmintegrableven}, and by a geometric sum argument,
the pair of inequalities \eqref{eq:firstmajoI?????} becomes
\begin{align*}
&\frac{2\;e^{-\epsilon}e^{2\delta_c\,(N+1)}\,
\wt \sigma^+_{\rm even}(\Omega^-)\,\wt \sigma^-_{\rm even}(\Omega^+)}
{(e^{2\delta_c}-1)\,\|m_{c}\|} -c_\epsilon 
\\&\;\;\;\;\;\;\;\;\;\;\;\;\;\;\leq i_{\eta}(2N)\leq
\frac{2\;e^{\epsilon}e^{2\delta_c\,(N+1)}\,
\wt \sigma^+_{\rm even}(\Omega^-)\,\wt \sigma^-_{\rm even}(\Omega^+)}
{(e^{2\delta_c}-1)\,\|m_{c}\|} +c_\epsilon\;.
\end{align*}
Up to replacing the summations from $n=0$ to $N$ to summations on even
numbers between $0$ to $2N$, and replacing $\lfloor n/2\rfloor$ by
$2\lfloor n/2\rfloor$ as well as $\lceil n/2 \rceil$ by $2\lceil
n/2\rceil$, the remaining part of the proof applies and gives the
result, noting that in Claim (iii) of Step 3T, we furthermore have
that the origin and endpoint of the constructed common perpendicular
$\alpha_\ga$ are in $V_{\rm even}\XX$.  This concludes the proof of
Equation \eqref{eq:simplicialevenevencase}.

\bigskip 
The remainder of the proof of Theorem
\ref{theo:discretemainequidupeven} consists in proving versions of the
equidistribution result Equation \eqref{eq:simplicialevenevencase} in
$\gengeododd\XX\times\gengeododd\XX$, $\gengeodeven\XX\times
\gengeododd\XX$, $\gengeododd\XX\times \gengeodeven\XX$ respectively ,
and in summing these four contributions.

\medskip 
By applying Equation \eqref{eq:simplicialevenevencase} by replacing
$x_0$ by a vertex $x'_0$ in $V_{\rm odd}\XX$, which exchanges $V_{\rm
  even} \XX$ and $V_{\rm odd}\XX$, $\;\gengeodeven\XX$ and
$\gengeododd\XX$, as well as $\wt\sigma^\pm_{\D^\mp,\,{\rm even}}$ and
$\wt\sigma^\pm_{\D^\mp,\,{\rm odd}}$, we have
\begin{align}
\lim_{t\ra+\infty} \;\frac{e^{2\,\delta_c}-1}{2\;e^{2\,\delta_c}}\;
\|m_{c}\|\;e^{-\delta_c\, t}
\sum_{\substack{i\in I^-/_\sim,\; j\in I^+/_\sim, \;\ga\in\Ga\\ 
\pi(\alpha^-_{i,\ga j}),\,\pi(\alpha^+_{\ga^{-1}i,j})\,\in\, V_{\rm odd}\XX\\
{D^-_i}\cap {D^+_{\ga j}}=\emptyset,\; 
\lambda_{i,\,\ga j}\leq t}} \;e^{\wt c\,(\alpha_{i,\ga j})}\; 
& \Dirac_{\alpha^-_{i,\,\ga j}} \otimes\Dirac_{\alpha^+_{\ga^{-1}i,\,j}}\nonumber\\
\;=\;\;& \wt\sigma^+_{\D^-,\,{\rm odd}}\otimes \wt\sigma^-_{\D^+,\,{\rm odd}}
\;\label{eq:simplicialoddoddcase}
\end{align}
for the weak-star convergence of measures on the locally
compact space $\gengeododd \XX\times \gengeododd \XX$.

\medskip 
Let us now apply Equation \eqref{eq:simplicialevenevencase} by
replacing $\D^-= (D^-_i)_{i\in I^-}$ by $\N_1\D^-= (\N_1 D^-_i)_{i\in
  I^-}$.  Let us consider the map $\varphi_+:\gengeod\XX\ra
\gengeod\XX$, which maps a generalised geodesic line $\ell$ to the
generalised geodesic line which coincides with $\flow{+1}\ell$ on
$[0,+\infty[$ and is constant (with value $\ell(1)$) on $]-\infty,0[$.
Note that this map is continuous and $\Ga$-equivariant, and that it
maps $\gengeodeven \XX$ in $\gengeododd \XX$ and $\gengeododd \XX$ in
$\gengeodeven\XX$.

Furthermore, by convexity, $\varphi_+$ induces for every $i\in I^-$ a
homeomorphism from $\normalout D^-_i$ to $\normalout \N_1D^-_i$, which
sends $\normalout D^-_i\cap \gengeododd\XX$ to $\normalout
\N_1D^-_i\cap \gengeodeven \XX$, such that, by Equation
\eqref{eq:pousstempunconduct}, for all $w\in \normalout D^-_i\cap
\gengeododd \XX$, if $e_w$ is the first  edge
followed by $w$
$$
d\wt\sigma^+_{D^-_i,\,{\rm odd}}(w) = e^{\wt c\,(e_w)-\delta_c} \; 
d\wt\sigma^+_{\N_1D^-_i,\,{\rm even}}(\varphi_+(w))\;.
$$

Note that for all $\ell>0$, there is a one-to-one correspondence
between the set of common perpendiculars of length $\ell$, with origin
and endpoint both in $V_{\rm even}$, between $\N_1D^-_i$ and $\ga
D^+_j$ for all $i\in I^-$, $j\in I^+$ and $\ga\in\Ga$, and the set of
common perpendiculars of length $\ell+1$, with origin in $V_{\rm odd}$
and endpoint in $V_{\rm even}$, between $D^-_i$ and $\ga D^+_j$ for
all $i\in I^-$, $j\in I^+$ and $\ga\in\Ga$. In particular,
$\varphi_+(\alpha^-_{i,\ga j})$ is the common perpendicular between
$\N_1D^-_i$ and $\ga D^+_j$, starting at time $t=0$ from $\N_1D^-_i$.

Therefore Equation \eqref{eq:simplicialevenevencase} applied
by replacing $\D^-= (D^-_i)_{i\in I^-}$ by $\N_1\D^-= (\N_1 D^-_i)_{i\in I^-}$
gives
\begin{align*}
  \lim_{t\ra+\infty} \;&\frac{e^{2\,\delta_c}-1}{2\;e^{2\,\delta_c}}\;
  \|m_{c}\|\;e^{-\delta_c\, t}\\ &
  \sum_{\substack{i\in I^-/_\sim,\; j\in I^+/_\sim, \;\ga\in\Ga\\
      \pi(\alpha^-_{i,\ga j})\,\in\, V_{\rm odd}\XX,\,
      \pi(\alpha^+_{\ga^{-1}i,j})\,\in\, V_{\rm even}\XX\\
      {D^-_i}\cap {D^+_{\ga j}}=\emptyset,\; \lambda_{i,\,\ga j}\leq
      t+1}} \;  e^{\wt c\,(e_{\alpha^-_{i,\,\ga j}})+\wt c\,(\varphi_+(\alpha^-_{i,\ga j}))}\;
   \Dirac_{\alpha^-_{i,\,\ga j}} \otimes\Dirac_{\alpha^+_{\ga^{-1}i,\,j}}\nonumber\\
  & \;\;=\;\;e^{\delta_c}\;\wt\sigma^+_{\D^-,\,{\rm odd}}\otimes
  \wt\sigma^-_{\D^+,\,{\rm even}} \;
\end{align*}
for the weak-star convergence of measures on the locally compact space
$\gengeododd \XX\times \gengeodeven \XX$. Since
$\wt c(e_{\alpha^-_{i,\,\ga j}})+\wt c(\varphi_+(\alpha^-_{i,\ga
  j}))=\wt c(\alpha^-_{i,\ga j})$, replacing $t$ by $t-1$ and simplifying
by $e^{\delta_c}$, we get
\begin{align}
  \lim_{t\ra+\infty} \;&\frac{e^{2\,\delta_c}-1}{2\;e^{2\,\delta_c}}\;
  \|m_{c}\|\;e^{-\delta_c\, t}\nonumber\\ &
  \sum_{\substack{i\in I^-/_\sim,\; j\in I^+/_\sim, \;\ga\in\Ga\\
\pi(\alpha^-_{i,\ga j})\,\in\, V_{\rm odd}\XX,\,
\pi(\alpha^+_{\ga^{-1}i,j})\,\in\, V_{\rm even}\XX\\
      {D^-_i}\cap {D^+_{\ga j}}=\emptyset,\; \lambda_{i,\,\ga j}\leq
      t}} \;e^{\wt c(\alpha^-_{i,\ga j})}\;
   \Dirac_{\alpha^-_{i,\,\ga j}} \otimes\Dirac_{\alpha^+_{\ga^{-1}i,\,j}}\nonumber\\
 & \;\;=\;\; \wt\sigma^+_{\D^-,\,{\rm odd}}\otimes
  \wt\sigma^-_{\D^+,\,{\rm even}} \;\label{eq:simplicialoddevencase}
\end{align}
for the weak-star convergence of measures on the locally compact space
$\gengeododd \XX\times \gengeodeven \XX$.

\medskip 
Now Theorem \ref{theo:discretemainequidupeven} follows by
summing Equation \eqref{eq:simplicialevenevencase}, Equation
\eqref{eq:simplicialoddoddcase}, Equation
\eqref{eq:simplicialoddevencase} and the formula, proven similarly,
obtained from Equation \eqref{eq:simplicialevenevencase} by replacing
$\D^+=(D^+_j)_{j\in I^+}$ by $(\N_1D^+_j)_{j\in I^+}$.
\cqfd

\medskip The following result for bipartite graphs (of groups) is used
in the arithmetic applications in Part \ref{sect:arithappli} (see
Section \ref{subsec:locconst}).

\bcoro\label{coro:twosubtreeseven} Let $\XX$ be a
$(p+1,q+1)$-biregular simplicial tree (with $p,q\geq 2$, possibly
with $p=q$), with corresponding partition $V\XX=V_p\XX\sqcup
V_q\XX$. Let $\Ga$ be a lattice of $\XX$ such that this partition is
$\Ga$-invariant. Assume that the Patterson density is normalised so
that $\|\mu_x\|= \frac{p+1}{\sqrt{p}}$ for every $x\in V_p\XX$. Let
$\DD^\pm$ be nonempty proper simplicial subtrees of $\XX$ with
stabilisers $\Ga_{\DD\pm}$ in $\Ga$, such that the families $\D^\pm=
(\ga\DD^\pm)_{\ga\in\Ga/\Ga_{\DD^\pm}}$ are locally finite. Let
$\sigma^\mp_{\D^\pm}$ be their skinning measures for the sero system
of conductances. Then
\begin{align*}
\lim_{t\ra+\infty} \;\frac{pq-1}{2}\;\TVol(\Ga\dbs\XX)\;\sqrt{pq}^{\,-t-2}
\sum_{\substack{(\alpha,\,\beta,\,\ga)\in 
\Ga/\Ga_{\DD^-}\times\Ga/\Ga_{\DD^+}\times\Ga\\
 0< d(\alpha\DD^-, \,\ga\beta\DD^+)\leq t}}
&  \Dirac_{\alpha^-_{\alpha,\,\ga\beta}} \otimes
  \Dirac_{\alpha^+_{\ga^{-1}\alpha,\,\beta}}\\ \;=\;\;&
\wt\sigma^+_{\D^-}\otimes \wt\sigma^-_{\D^+}\;
\end{align*}
for the weak-star convergence of measures on the locally
compact space $\gengeod \XX\times \gengeod \XX$.  

If the measure $\sigma^-_{\D^+}$ is nonzero and finite, then
$$
\lim_{t\ra+\infty} \;\frac{pq-1}{2}\;
\frac{\TVol(\Ga\dbs\XX)}{\|\sigma^-_{\D^+}\|}\;\sqrt{pq}^{\,-t-2}
\sum_{\substack{\ga\in \Ga/\Ga_{\DD^+}\\
0< d(\DD^-, \, \ga \DD^+)\leq t}}
\Dirac_{\alpha^-_{e,\,\ga}} \;=\;
\wt\sigma^+_{\DD^-}\;,
$$
for the weak-star convergence of measures on the locally
compact space $\gengeod \XX$. 
\ecoro

\dem In order to prove the first result, we apply Theorem
\ref{theo:discretemainequidupeven} with $\wt c\equiv 0$, so that by
Equation \eqref{eq:lyonshdim}, Proposition \ref{prop:uniflatmBMfinie}, 
Theorem \ref{theo:uniflatmBMmixing}, and Proposition
\ref{prop:computBM} (2), we have $\delta_c=\frac{1}{2}\ln (pq)>0$,
$m_{c}=m_{\rm BM}$ is finite and its restriction to $\Ga\bs \Geven\XX$
is mixing under the square of the geodesic flow, and $\|m_{\rm BM}\|=
\TVol(\Ga\dbs\XX)$.

The second claim follows as in the proof of Corollary
\ref{coro:twosubtrees}.
\cqfd

\medskip \rem In some special occasions, the measures involved in the
statements of Theorem \ref{theo:discretemainequidupeven} and Corollary
\ref{coro:twosubtreeseven} (whether skinning measures or Dirac masses)
are actually all supported on $\gengeodeven \XX$ (up to choosing
appropriately $x_0$). This is in particular the case if $\D^\pm=
(\ga\DD^\pm)_{\ga\in\Ga/\Ga_{\DD^\pm}}$ with $\DD^-, \DD^+$ horoballs at even
signed distance (see below), as the following lemma shows.

\begin{center}
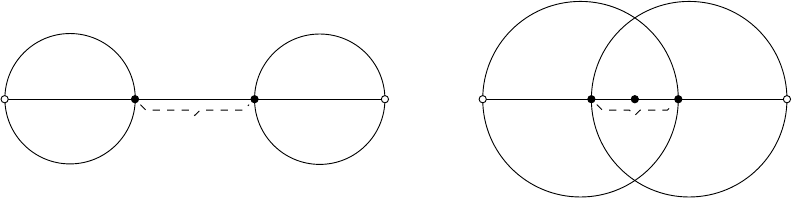
\end{center}

The {\em signed distance}\index{signed
  distance}\index{distance!signed} between horoballs $\H$ and $\H'$ in
an $\RR$-tree that are not centred at the same point at infinity
is the distance between them (that is, the length of their common
perpendicular) if they are disjoint, or the opposite of the diameter
of their intersection otherwise. Note that if nonempty, the
intersection of $\H$ and $\H'$ is a ball centred at the midpoint of
the segment contained in the geodesic line between the two points at
infinity of the horoballs, which lies in both horoballs.

\blemm \label{lem:horodisjoint}
Let $\XX$ be a simplicial tree, $\Ga$ a subgroup of $\Aut(\XX)$
and $\H,\H'$ two horoballs in $\XX$ (whose boundaries are contained in
$V\XX$), which either are equal or have distinct points at infinity.
If $\Lambda\Ga\subset 2\ZZ$ and $\H,\H'$ are at even signed distance,
then the signed distance between $\H$ and $\ga\H'$ is even for every
$\ga\in\Ga$ such that $\H$ and $\ga\H'$ do not have the same point at
infinity.  
\elemm

\dem Fix such a $\ga$. For every horoball $\H''$ and for all
$s\in\NN$, let $\H''[s]$ be the horoball contained in $\H$, whose
boundary is at distance $s$ from the boundary of $\H$. Shrinking the
horoballs $\H$ and $\H'$, by replacing them by the horoballs $\H[s]$
and $\H'[s]$ for any $s\in\NN$, only changes by $\pm 2s$ the
considered signed distances. Hence, taking $s$ large enough, we may
assume that $\H$ and $\ga\H'$ are disjoint, and that $\H$ and $\H'$
are disjoint or equal.  Let $[x,x']$ be the common perpendicular
between $\H$ and $\H'$ with $x\in\partial\H$, $x'\in\partial\H'$ if
$\H$ and $\H'$ are disjoint, and otherwise, let $x=x'$ be any point in
$\partial\H=\partial\H'$. Let $[y,y']$ be the common perpendicular
between $\H$ and $\ga\H'$, with $y\in\partial\H$, $y'\in
\partial(\ga\H')$. Note that $\ga x'\in \partial(\ga\H')$.

\begin{center}
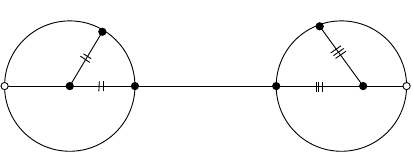
\end{center}

The distance between two points $x,y$ of a horosphere is always even
(equal to twice the distance from $x$ to the geodesic ray from $y$ to
the point at infinity of the horosphere). Since geodesic triangles in
trees are tripods, for all $a,b,c$ in a simplicial tree, since 
$$
d(a,c)=d(a,b)+d(b,c)-2d(b,[a,c])\;,
$$ 
if $d(a,b)$ and $d(b,c)$ are even, so is $d(a,c)$.

Since $\Lambda\Ga\subset 2\ZZ$, the distance between $x'$ and $\ga x'$
is even by Equation \eqref{eq:equivLamGaevendisteven}. Since $d(x,x')$
is even by assumption, we hence have that $d(x,\ga x')$ is even.
Therefore
$$
d(y,y')=d(x,\ga x') - d(x,y) - d(y',\ga x')
$$
is even.
\cqfd

\chapter{Equidistribution and counting of common 
perpendiculars in quotient spaces}
\label{sec:equidcountdownstairs}

In this Chapter, we use the results of Chapter \ref{sec:equidarcs} to
prove equidistribution and counting results in Riemannian manifolds
(or good orbifolds) and in metric and simplicial graphs (of groups).

Let $X$, $x_0$, $\Ga$ and $\wt F$ be as in the beginning of Chapter
\ref{sec:equidarcs}. We will need the following two notions in this
chapter.

Recall that the {\em narrow topology}\footnote{also called {\em weak
    topology}\index{convergence!narrow}%
\index{narrow topology}\index{topology!narrow}
\index{convergence!weak}%
\index{weak topology}\index{topology!weak} see for instance
\cite[p.~71-III]{DelMey78} or \cite{Billingsley68,Parthasarathy67}} on
the set $\M_{\rm f}(Y)$ of finite measures on a Polish space $Y$ is
the smallest topology such that the map from $\M_{\rm f}(Y)$ to $\RR$
defined by $\mu\mapsto \mu(g)$ is continuous for every bounded
continuous function $g:Y\ra\RR$. Since continous functions with
compact support are bounded, the narrow convergence of finite measures
implies their weak-star convergence.

 Recall that given a discrete group $G$ acting properly (but not
 necessarily freely) on a locally compact space $Z$, the {\it induced
   measure}\index{measure!induced}\footnote{See for instance \cite[\S
     2.6]{PauPolSha15} for details.} on $G\bs Z$ of a (positive,
 Radon) measure $\mu$ on $Z$ is a measure $\overline{\mu}$ which
 depends linearly and continuously for the weak-star topology on
 $\mu$, and satisfies $\overline{\Dirac_z}= \frac{1}{|G_z|}\;\Dirac_{G
   z}$ for every $z\in Z$.  The following observation on the behaviour
 of induced measures under quotients by properly discontinuous group
 actions will be used in the proofs of Corollary
 \ref{coro:mainequicountdown} and its analogues in Section
 \ref{subsect:equicountmetricgraphgroup}. Let $G$ be a discrete group
 that acts properly on a Polish space $\wt Y$ and let $Y=G\bs\wt Y$.
 Let $\wt \mu_k$ for $k\in\NN$ and $\wt \mu$ be $G$-invariant locally
 finite measures on $\wt Y$, with finite induced measures $\mu_k$ for
 $k\in\NN$ and $\mu$ on $Y$.  If for every Borel subset $B$ of $\wt Y$
 with $\wt \mu(B)$ finite and $\wt \mu(\partial B)=0$ we have
 $\lim_{k\ra\infty}\wt \mu_k(B)=\wt \mu(B)$, then the sequence
 $(\mu_k)_{k\in\NN}$ narrowly converges to $\mu$.

\section{Multiplicities and counting functions in Riemannian orbifolds}
\label{subsec:multandcount}

In this Section, we assume that $X= \wt M$ is a Riemannian
manifold. We denote its quotient Riemannian orbifold under $\Ga$ by
$M=\Ga\bs\wt M$, and the quotient Riemannian orbifold under $\Ga$ of
its unit tangent bundle by $T^1M=\Ga\bs T^1\wt M$. We use the
identifications $\G X=\G_{\pm,\,0} X= T^1X=T^1\wt M$ explained in
Chapter \ref{sec:negcurv}.

Let $\D=(D_i)_{i\in I}$ be a locally finite $\Ga$-equivariant family
of nonempty proper closed convex subsets of $\wt M$.  Let $\Omega=
(\Omega_i)_{i\in I}$ be a $\Ga$-equivariant family of subsets of
$T^1\wt M$, where $\Omega_i$ is a measurable subset of $\normalpm D_i$
for all $i\in I$ (the sign $\pm$ being constant) and
$\Omega_i=\Omega_j$ if $i\sim_\D j$.  The {\it
  multiplicity}\index{multiplicity} of an element $v\in T^1M$ with
respect to $\Omega$ is
$$
m_{\Omega}(v)= \frac{\card\,\{i\in I/_\sim\;:\;\wt v\in\Omega_i\}}
{\card(\stab_\Ga\wt v)}\;,
$$ 
for any preimage $\wt v$ of $v$ in $T^1\wt M$. The numerator and
the denominator are finite by the local finiteness of the family $\D$
and the discreteness of $\Ga$, and they depend only on the orbit of
$\wt v$ under $\Ga$.

The numerator takes into account the multiplicities of the images of
the elements of $\Omega$ in $T^1M$. The denominator of this
multiplicity is also natural, as any counting problem of objects
possibly having symmetries, the appropriate counting function consists
in taking as the multiplicity of an object the inverse of the
cardinality of its symmetry group.

\bexems The following examples illustrate the behaviour of the
multiplicity when $\Ga$ is torsion-free and $\Omega=\normalpm\D$.

\smallskip\noindent (1) If for every $i\in I$, the quotient
$\Ga_{D_i}\bs D_i$ of $D_i$ by its stabiliser $\Ga_{D_i}$ maps
injectively in $M$ by the map induced by the inclusion of $D_i$ in
$\wt M$, and if for every $i,j\in I$ such that $j\notin \Ga i$, the
intersection $D_i\cap D_j$ is empty, then the nonzero multiplicities
$m_{\Omega}(\ell)$ are all equal to $1$.

\medskip 
\noindent
\begin{minipage}{11.4cm} 
(2) Here is a simple example of a multiplicity different from $0$ or
  $1$.  Let $c$ be a closed geodesic in the Riemannian manifold $M$,
  let $\wt c$ be a geodesic line in $\wt M$ mapping to $c$ in $M$, let
  $\D=(\ga\,\wt c)_{\ga\in\Ga}$, let $x$ be a double point of $c$, let
  $v\in T^1_xM$ be orthogonal to the two tangent lines to $c$ at $x$
  (this requires the dimension of $\wt M$ to be at least $3$, if $x$
  is a transverse self-intersection point). Then
  $m_{\normalpm\D}(v)=2$.
\end{minipage}
\begin{minipage}{3.5cm}
\begin{center}
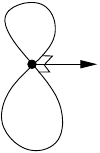
\end{center}
\end{minipage}
\eexems

\bigskip Given $t> 0$ and two unit tangent vectors $v,w\in T^1M$, we
define the number $n_t(v,w)$ of locally geodesic paths having $v$ and
$w$ as initial and terminal tangent vectors respectively, weighted by
the potential $F$, with length at most $t$, by
$$
n_t(v,w)=\sum_\alpha\;\;\card(\Ga_\alpha)\;e^{\int_\alpha F}\,,
$$ 
where the sum ranges over the locally geodesic paths
$\alpha:[0,s]\ra M$ in the Riemannian orbifold $M$ such that
$\dot\alpha(0)=v$, $\dot\alpha(s)=w$ and $s\in\;]0,t]$, and
$\Ga_\alpha$ is the stabiliser in $\Ga$ of any geodesic path $\wt
\alpha$ in $\wt M$ mapping to $\alpha$ by the quotient map $\wt M\ra
M$. If $F=0$ and $\Ga$ is torsion free, then $n_t(v,w)$ is precisely
the number of locally geodesic paths having $v$ and $w$ as initial and
terminal tangent vectors respectively, with length at most $t$.

\medskip
Let $\Omega^-=(\Omega^-_i)_{i\in I^-}$ and $\Omega^+= (\Omega^+_j)
_{j\in I^+}$ be $\Ga$-equivariant families of subsets of $T^1\wt M$,
where $\Omega^\mp_k$ is a measurable subset of $\normalpm D^\mp_k$ for
all $k\in I^\mp$ and $\Omega^\pm_k=\Omega^\pm_{k'}$ if $k\sim_{\D^\pm}
k'$. We will denote by $\N_{\Omega^-,\,\Omega^+,\,F}: \;]0,+\infty[\;
\ra\RR$ the following {\it counting function}%
\index{counting function}: for every $t>0$, let
$\N_{\Omega^-,\,\Omega^+,\,F}(t)$ be the number of common
perpendiculars whose initial vectors belong to the images in $T^1M$ of
the elements of $\Omega^-$ and terminal vectors to the images in
$T^1M$ of the elements of $\Omega^+$, counted with multiplicities and
weighted by the potential $F$, that is:
$$
\N_{\Omega^-,\,\Omega^+,\, F}(t)=
\sum_{v,\,w\,\in \,T^1M} m_{\Omega^-}(v)\;m_{\Omega^+}(w) \;n_t(v,w)\;.
$$
When $\Omega^\pm=\normalmp \D^\pm$, we denote
$\N_{\Omega^-,\,\Omega^+,\, F}$ by $\N_{\D^-,\,\D^+,\, F}$.

\brema\label{rem:pnpiclcs} Let $Y$ be a negatively curved complete
connected Riemannian manifold and let $\wt Y\to Y$ be its Riemannian
universal cover.  Let $D^\pm$ be a locally convex\footnote{not
  necessarily connected} geodesic metric space endowed with a
continuous map $f^\pm:D^\pm\ra Y$ such that if $\wt D^\pm\ra D^\pm$ is
a locally isometric universal cover and if $\wt f^\pm:\wt D^\pm\ra\wt
Y$ is a lift of $f^\pm$, then $\wt f^\pm$ is on each connected
component of $\wt D^\pm$ an isometric embedding whose image is a
proper nonempty closed locally convex subset of $\wt Y$, and the
family of images under the covering group of $\wt Y\ra Y$ of the
images by $\wt f^\pm$ of the connected components of $\wt D^\pm$ is
locally finite. Then $D^\pm$ (or the pair $(D^\pm,f^\pm)$) is a {\em
  proper nonempty properly immersed closed locally convex
  subset}\index{proper nonempty properly immersed closed locally
  convex subset} of $Y$.  

If $\Gamma$ is a discrete subgroup without torsion of isometries of a
CAT($-1$) Riemannian manifold $X$, if $\D^\pm=(\ga \wt D^\pm)_{\ga
  \in\Ga}$ where $\wt D^\pm$ is a nonempty proper closed convex subset
of $X$ such that the family $\D^\pm$ is locally finite, and if $D^\pm$
is the image of $\wt D^\pm$ by the covering map $X\ra\Ga\bs X$, then
$D^\pm$ is a proper nonempty properly immersed closed convex subset of
$\Ga\bs X$.  Under these assumptions, $\N_{\D^-,\,\D^+,\, F}$ is the
counting function $\N_{D^-,\,D^+,\, F}$ given in the introduction.
\erema

\medskip 
Let us continue fixing the notation used in Sections
\ref{subsec:downstairs} and \ref{subsect:erroterms}.  For every
$(i,j)$ in $I^-\!\times I^+$ such that $D^-_i$ and $D^+_j$ have a
common perpendicular\footnote{that is, whose closures
  $\overline{D^-_i}$ and $\overline{D^+_j}$ in $X\cup \partial_\infty
  X$ have empty intersection}, we denote by $\alpha_{i,\,j}$ this
common perpendicular, by $\lambda_{i,\,j}$ its length, by
$v^-_{i,\,j} \in \partial^1_+ D_i^-$ its initial tangent vector and by
$v^+_{i,\,j}\in \partial^1_- D_i^+$ its terminal tangent vector. Note
that if $i'\sim i$, $j'\sim j$ and $\ga\in\Ga$, then
\begin{equation}\label{eq:equivalphvpm}
\ga\,\alpha_{i',\,j'}=\alpha_{\ga i,\,\ga j},\;\;\;\lambda_{i',\,j'}=
\lambda_{\ga i,\,\ga j}\;\;\;{\rm and}\;\;\;
\ga\, v^\pm_{i',\,j'}=v^\pm_{\ga i,\,\ga j}\;.
\end{equation}
When $\Ga$ has no  torsion, we have, for the diagonal action of $\Ga$
on $I^-\times I^+$,
$$
\N_{\D^-,\,\D^+,\, F}(t)=
\sum_{
(i,\,j)\in \Ga\bs((I^-/_\sim)\times (I^+/_\sim))\,:\;
\overline{D^-_i}\cap \overline{D^+_j}=\emptyset,\; 
\lambda_{i,\, j}\leq t}
\;e^{\int_{\alpha_{i,\,j}} \wt F}\;.
$$

When the potential $F$ is zero and $\Ga$ has no torsion,
$\N_{\D^-,\,\D^+,\, F}(t)$ is the number of common perpendiculars of
length at most $t$, and the counting function $t\mapsto
\N_{\D^-,\,\D^+,\,0}(t)$ has been studied in various special cases of
negatively curved manifolds since the 1950's and in a number of recent
works, see the Introduction. The asymptotics of $\N_{\D^-,\,\D^+,\,0}
(t)$ as $t\ra+\infty$ in the case when $X$ is a Riemannian manifold
with pinched negative curvature are described in general in
\cite[Theo.~1]{ParPau16ETDS}, where it is shown that if the skinning
measures $\sigma^+_{\D^-}$ and $\sigma^{-}_{\D^+}$ are finite and
nonzero, then as $s\to+\infty$,
\begin{equation}\label{eq:countmanifoldF=0}
\N_{\D^-,\,\D^+,\,0}(s)\sim
\frac{\|\sigma^+_{\D^-}\|\,\|\sigma^{-}_{\D^+}\|}
{\|m_{\rm BM}\|}\,\frac{e^{\delta_\Ga \,s}}{\delta_\Ga}\,.
\end{equation}

\section{Common perpendiculars in Riemannian orbifolds}
\label{subsec:downstairs}

Corollary \ref{coro:mainequicountdown} below is the main result of
this text on the counting with weights of common perpendiculars and on
the equidistribution of their initial and terminal tangent vectors in
negatively curved Riemannian manifolds endowed with a H\"older-continuous
potential. We use the notation of Section \ref{subsec:multandcount}.

\bcoro\label{coro:mainequicountdown} Let $\wt M$ be a complete simply
connected Riemannian manifold with pinched negative sectional
curvature at most $-1$.  Let $\Ga$ be a nonelementary discrete group
of isometries of $\wt M$.  Let $\wt F:T^1\wt M\to\RR $ be a bounded
$\Ga$-invariant H\"older-continuous function with positive critical
exponent $\delta$. Let $\D^-=(D^-_i)_{i\in I^-}$ and
$\D^+=(D^+_j)_{j\in I^+}$ be locally finite $\Ga$-equivariant families
of nonempty proper closed convex subsets of $\wt M$. Assume that the
Gibbs measure $m_{F}$ is finite and mixing for the geodesic flow on
$T^1M$. Then,
\begin{align}
&\lim_{t\ra+\infty}\; \delta\;\|m_F\|\;e^{-\delta\, t}
\sum_{v,\,w\in T^1M} m_{\normalout\D^-}(v)\;m_{\normalin\D^+}(w)\; n_t(v,w) \;
\Dirac_{v} \otimes\Dirac_{w}\nonumber\\& =\;
\sigma^+_{\D^-}\otimes \sigma^-_{\D^+}\;\label{eq:equidistribdown}
\end{align}
for the weak-star convergence of measures on the locally compact space
$T^1M\times T^1M$. If $\sigma^+_{\D^-}$ and $\sigma^-_{\D^+}$ are
finite, the result also holds for the narrow convergence.

Furthermore, for all $\Ga$-equivariant families $\Omega^\pm=
(\Omega^\pm_k)_{k\in I^\pm}$ of subsets of $T^1\wt M$ with
$\Omega_k^\pm$ a Borel subset of $\partial^1_{\mp} D_k^\pm$ for all
$k\in I^\pm$ and $\Omega^\pm_k=\Omega^\pm_{k'}$ if $k\sim_{\D^\pm}
k'$, with nonzero finite skinning measure and with boundary
in $\partial^1_{\pm} D_k^\mp$ of zero skinning measure, we have, as
$t\ra+\infty$,
$$
\N_{\Omega^-,\,\Omega^+,\, F}(t)\;\sim\;
\frac{\|\sigma^+_{\Omega^-}\|\;\|\sigma^-_{\Omega^+}\|}{\delta\;\|m_F\|}\;
e^{\delta\, t}\;.
$$
\ecoro

\dem 
Note that the sum in Equation \eqref{eq:equidistribdown}
is locally finite, hence it defines a locally finite measure on
$T^1M\times T^1M$.  We are going to rewrite the sum in the statement
of Theorem \ref{theo:mainequidup} in a way which makes it easier to
push it down from $T^1\wt M\times T^1\wt M$ to $T^1M\times T^1M$.

For every $\wt v\in T^1\wt M$, let
$$
m^\mp(\wt v)=\card\,\{k\in I^\mp/_\sim\;:\;\wt v\in\normalpm D^\mp_k\}\;,
$$ 
so that for every $v\in T^1M$, the multiplicity of $v$ with respect
to the family $\normalpm \D^\mp$ is\footnote{See Section
  \ref{subsec:multandcount}.}
$$
m_{\normalpm \D^\mp}(v)=\frac{m^\mp(\wt v)}{\card(\stab_\Ga\wt v)}\;,
$$
for any preimage $\wt v$ of $v$ in $T^1\wt M$.

For all $\ga\in\Ga$ and $\wt v,\wt w\in T^1\wt M$, there exists
$(i,j)\in (I^-/_\sim)\times (I^+/_\sim)$ such that $\wt v=v^-_{i,\ga
  j}$ and $\wt w=v^+_{\ga^{-1}i,j}=\ga^{-1}v^+_{i,\ga j}$ if and only
if $\ga \wt w\in\flow{\RR}\,\wt v$, there exists $i'\in I^-/_\sim$
such that $\wt v\in \normalout D^-_{i'}$ and there exists $j'\in
I^+/_\sim$ such that $\ga\wt w\in \normalin D^+_{j'}$. Then the choice
of such elements $(i,j)$, as well as $i'$ and $j'$, is free.  We hence
have
\begin{align*}
&\sum_{\substack{i\in I^-/_\sim,\; j\in I^+/_\sim, \;\ga\in\Ga\\ 
0<\lambda_{i,\,\ga j}\leq t\,,\;
v^-_{i,\,\ga j}=\wt v\,,\;v^+_{\ga^{-1}i,\,j}=\wt
w}} \;e^{\int_{\alpha_{i,\,\ga j}}\wt F}\; 
\Dirac_{v^-_{i,\,\ga j}} \otimes\Dirac_{v^+_{\ga^{-1}i,\,j}}\\\;=&
\sum_{\substack{\ga\in\Ga,\;0<s\leq t\\ 
\ga\wt w=\flow{s}\wt v}}e^{\int_{\pi(\wt v)}^{\ga\pi(\wt w)}\wt F}\; 
\card\,\big\{(i,j)\in (I^-/_\sim)\times (I^+/_\sim) \;:\;
v^-_{i,\,\ga j}=\wt v\,,\;v^+_{\ga^{-1}i,\,j}=\wt w\big\}\;
\Dirac_{\wt v} \otimes\Dirac_{\wt w}\\\;=&
\sum_{\substack{\ga\in\Ga,\;0<s\leq t\\ 
\ga\wt w=\flow{s}\wt v}}
 \;e^{\int_{\pi(\wt v)}^{\ga\pi(\wt w)}\wt F}\; 
m^-(\wt v)\;m^+(\ga \wt w)\;
\Dirac_{\wt v} \otimes\Dirac_{\wt w}\;.
\end{align*}
Therefore
\begin{align*}
&\sum_{\substack{i\in I^-/_\sim,\; j\in I^+/_\sim, \;\ga\in\Ga\\ 
0<\lambda_{i,\,\ga j}\leq t}} \;e^{\int_{\alpha_{i,\,\ga j}}\wt F}\; 
\Dirac_{v^-_{i,\,\ga j}} \otimes\Dirac_{v^+_{\ga^{-1}i,\,j}}\\\;=&
\sum_{\wt v,\,\wt w\,\in\,T^1\wt M}\;
\Big(\sum_{\substack{\ga\in\Ga,\;0<s\leq t\\ 
\ga\wt w=\flow{s}\wt v}}
 \;e^{\int_{\pi(\wt v)}^{\ga\pi(\wt w)}\wt F}\;\Big)\; 
m^-(\wt v)\;m^+(\wt w)\; \Dirac_{\wt v} \otimes\Dirac{\wt w}\;.
\end{align*}

By definition, $\sigma^\pm_{\D^\mp}$ is the measure on $T^1M$ induced
by the $\Ga$-invariant measure $\wt \sigma^\pm_{\D^\mp}$. Thus
Corollary \ref{coro:mainequicountdown} follows from Theorem
\ref{theo:mainequidup} after a similar reduction\footnote{See Step 1
  of the proof of Theorem \ref{theo:mainequidup}.} as in Section
\ref{subsec:equidcommperpcontI}. The narrow convergence is obtained
when the skinning measures $\sigma^\pm_{\D^\mp}$ are finite, using the
continuity properties of the induced measures recalled at the beginning
of Chapter \ref{sec:equidcountdownstairs}, since no compactness
assumptions were made in Equation \eqref{eq:reducnarrowup} on
$\Omega^\pm$.

The counting statement follows from the equidistribution result by
integration.  
\cqfd

\bigskip 
In particular, if the skinning measures $\sigma^+_{\D^-}$ and
$\sigma^-_{\D^+}$ are positive and finite, Corollary
\ref{coro:mainequicountdown} gives, as $t\ra+\infty$, 
$$
\N_{\D^-,\,\D^+,\, F}(t)\;\sim\;
\frac{\|\sigma^+_{\D^-}\|\;\|\sigma^-_{\D^+}\|}{\delta\;\|m_F\|}\;
e^{\delta\, t}\;.
$$

\medskip 
\brema Under the assumptions of Corollary \ref{coro:mainequicountdown}
with the exception that $\delta$ may now be nonpositive, we have the
following asymptotic result as $t\ra+\infty$ for the growth of the
weighted number of common perpendiculars with lengths in $]t-\tau,t]$
  for every fixed $\tau>0$:
$$
\N_{\D^-,\,\D^+,\, F}(t)-\N_{\D^-,\,\D^+,\, F}(t-\tau)\;\sim\;
\frac{(1-e^{-\delta\, \tau})\,\|\sigma^+_{\D^-}\|\;\|\sigma^-_{\D^+}\|}
{\delta\;\|m_F\|}\; e^{\delta\, t}\;.
$$ 
This result follows by considering a large enough constant $\sigma$
such that $\delta_{\Ga,\,F+\sigma}=\delta+\sigma>0$, by applying
Corollary \ref{coro:mainequicountdown} with the potential $F+\sigma$
(see Remark \ref{rem:skinrem} (2)) as in the proof of Theorem
\ref{theo:mainequidupbis}.  
\erema

\medskip 
Using the continuity of the pushforwards of measures for the weak-star
and the narrow topologies, applied to the basepoint maps $\pi\times
\pi$ from $T^1\wt M\times T^1\wt M$ to $\wt M\times \wt M$, and from
$T^1M\times T^1M$ to $M\times M$, we have the following result of
equidistribution of the ordered pairs of endpoints of common
perpendiculars between two equivariant families of convex sets in $\wt
M$ or two families of locally convex sets in $M$.  When $M$ has
constant curvature and finite volume, $F=0$ and $\D^-$ is the
$\Ga$-orbit of a point and $\D^+$ is the $\Ga$-orbit of a totally
geodesic cocompact submanifold, this result is due to Herrmann
\cite{Herrmann62}.  When $\D^\pm$ are $\Ga$-orbits of points and $F$
is a H\"older-continuous potential, see \cite[Theo.~9.1,
  9.3]{PauPolSha15}, and we refer for instance to \cite{BoyMay16} for
an application of this particular case.

\bcoro\label{coro:mainequidbaspoint} Let $\wt M,\Ga,\wt F,\D^-,\D^+$
be as in Corollary \ref{coro:mainequicountdown}. Then
\begin{align*}
&
\lim_{t\ra+\infty} \;\delta\;\|m_F\|\;e^{-\delta\, t}\!\!
\sum_{\substack{i\in I^-/_\sim,\; j\in I^+/_\sim, \;\ga\in\Ga\\ 
0<\lambda_{i,\,\ga j}\leq t}} 
\!e^{\int_{\alpha_{i,\,\ga j}}\wt F}\; 
\Dirac_{\pi(v^-_{i,\,\ga j})} \otimes
\Dirac_{\pi(v^+_{\ga^{-1}i,\,j})}\\&=\; 
\pi_*\wt\sigma^+_{\D^-}\otimes \pi_*\wt\sigma^-_{\D^+}\;,
\end{align*}
for the weak-star convergence of measures on the locally compact space
$\wt M\times \wt M$, and
\begin{align*}
&\lim_{t\ra+\infty}\; \delta\;\|m_F\|\;e^{-\delta\, t}\!
\sum_{v,\,w\in T^1M} \!m_{\normalout\D^-}(v)\;m_{\normalin\D^+}(w)\; n_t(v,w) \;
\Dirac_{\pi(v)} \otimes\Dirac_{\pi(w)}\\&=\;\;
\pi_*\sigma^+_{\D^-}\otimes \pi_*\sigma^-_{\D^+}\;,
\end{align*}
for the weak-star convergence of measures on $M\times M$.  If the
measures $\sigma^\pm_{\D^\mp}$ are finite, then the above claim holds
for the narrow convergence of measures on $M\times M$. \cqfd \ecoro

\medskip 
We will now prove Theorems \ref{theo:mainintroequidis} and
\ref{theo:mainintrocount} (1) in the Introduction for
Riemannian manifolds.  Recall from Remark \ref{rem:pnpiclcs} the
definition of proper nonempty properly immersed closed locally convex
subsets $D^\pm$ in a pinched negatively curved complete connected
Riemannian manifold $Y$ and the associated maps $\wt f^\pm:\wt
D^\pm\ra\wt Y$.

\medskip\noindent{\bf Proof of Theorems \ref{theo:mainintroequidis}
  and \ref{theo:mainintrocount} (1) for Riemannian manifolds. } Let
$Y,F,D^\pm$ be as in these statements and assume that $Y$ is a
Riemannian manifold.  Let $\Ga$ be the covering group of a universal
Riemannian cover $\wt Y\ra Y$. Let $I^\pm=\Ga\times \pi_0(\wt D^\pm)$
with the action of $\Ga$ defined by $\ga\cdot (\alpha,c) = (\ga
\alpha,c)$ for all $\ga,\alpha\in \Ga$ and every connected component
$c$ of $\wt D^\pm$. Consider the families $\D^\pm=(D^\pm_k)_{k\in
  I^\pm}$ where $D^\pm_k= \alpha \;\wt f^\pm(c)$ if $k=(\alpha,c)$.
Then $\D^\pm$ are $\Ga$-equivariant families of nonempty proper closed
convex subsets of $\wt Y$, which are locally finite since $D^\pm$ are
properly immersed in $Y$. The conclusions in Theorems
\ref{theo:mainintroequidis} and \ref{theo:mainintrocount} (1) when $Y$
is a manifold then follow from Corollary \ref{coro:mainequicountdown},
applied with $\wt M=\wt Y$ and with $\wt F$ the lift of $F$ to $T^1\wt
M$.  \cqfd

\bcoro Let $\wt M,\Ga,\wt F,\D^-,\D^+$ be as in Corollary
\ref{coro:mainequicountdown}. Assume that $\sigma^\pm_{\D^\mp}$ are
finite and nonzero. Let
$$
n_{t,\D^+}(v)=\sum_{w\in T^1M}m_{\normalin\D^+}(w)\;n_t(v,w)
$$
be the number (counted with multiplicities) of locally geodesic paths
in $M$ of length at most $t$, with initial vector $v$, arriving
perpendicularly to $\D^+$. Then
$$
\lim_{s\to+\infty}\lim_{t\to+\infty}
\frac{\delta\,\|m_F\|^2\,e^{-\delta t}}
{\|\sigma^+_{\D^-}\|\|\sigma^-_{\D^+}\|} 
\sum_{v\in T^1M} m_{\normalout\D^-}(v)\;n_{t,\,\D^+}(v)\;
\Dirac_{\flow s v}=m_F\,,
$$ 
for the narrow convergence on $\Ga\bs\G \wt M$.
\ecoro

\dem For every $s\in\RR$, by Corollary \ref{coro:mainequicountdown},
using the continuity of the pushforwards of measures by the first
projection $(v,w)\mapsto v$ from $T^1M\times T^1M$ to $T^1M$, and by
the geodesic flow on $T^1M$ at time $s$, since
$(\flow{s})_*\Dirac_v=\Dirac_{\flow{s}v}$, we have
$$
\lim_{t\ra+\infty}\; \delta\;\|m_F\|\;e^{-\delta t}
\sum_{v\in T^1M} m_{\normalout\D^-}(v)\;n_{t,\D^+}(v) \;
\Dirac_{\flow{s}v} \;=\;
(\flow{s})_*\sigma^+_{\D^-}\|\sigma^-_{\D^+}\|\,.
$$
The result then follows from Theorem \ref{theo:equid} with
$\Omega=\normalout\D^-$.  \cqfd

\section{Error terms for  equidistribution and 
counting for Riemannian orbifolds}
\label{subsect:erroterms}

In Section \ref{subsec:mixingratemanifolds}, we discussed various
results on the rate of mixing of the geodesic flow for Riemannian
manifolds. In this Section, we apply these results to give error
bounds to the statements of equidistribution and counting of common
perpendicular arcs given in Section \ref{subsec:downstairs}. We use
again the notation of Section \ref{subsec:multandcount}.

\btheo\label{theo:expratecount} Let $\wt M$ be a complete simply
connected Riemannian manifold with pinched negative sectional
curvature at most $-1$.  Let $\Ga$ be a nonelementary discrete group
of isometries of $\wt M$.  Let $\wt F:T^1\wt M\to\RR $ be a bounded
$\Ga$-invariant H\"older-continuous function with positive critical
exponent $\delta$.  Assume that $(\wt M,\Ga,\wt F)$ has
radius-H\"older-continuous strong stable/unstable ball masses. Let
$\D^-=(D^-_i)_{i\in I^-}$ and $\D^+=(D^+_j)_{j\in I^+}$ be locally
finite $\Ga$-equivariant families of nonempty proper closed convex
subsets of $\wt M$ such that $\Ga\bs I^\pm$ are finite, with finite
nonzero skinning measure $\sigma_{\D^-}$ and $\sigma_{\D^+}$. Let
$M=\Ga\backslash\wt M$ and let $F:T^1M\ra \RR$ be the potential
induced by $\wt F$.

\smallskip\noindent (1) Assume that $M$ is compact and that the
geodesic flow on $T^1M$ is mixing with exponential speed for the
H\"older regularity for the potential $F$.  Then there exist
$\alpha\in\;]0,1]$ and $\kappa'>0$ such that for all nonnegative
$\psi^\pm\in \C_{\rm c}^\alpha(T^1M)$, we have, as $t\ra+\infty$,
\begin{multline*}
\frac{\delta\;\|m_F\|}{e^{\delta\, t}}
\sum_{v,\,w\in T^1M} m_{\normalout\D^-}(v)\;m_{\normalin\D^+}(w)\; n_t(v,w) \;
\psi^-(v)\,\psi^+(w)\\=
\int_{T^1M}\psi^-d\sigma^+_{\D^-} \int_{T^1M}\psi^+d\sigma^-_{\D^+}
+\operatorname{O}(e^{-\kappa't}\|\psi^-\|_\alpha\,\|\psi^+\|_\alpha)\,.
\end{multline*}

\smallskip\noindent (2) Assume that $\wt M$ is a symmetric space, that
$D^\pm_k$ has smooth boundary for every $k\in I^\pm$, that $m_F$ is
finite and smooth,\footnote{Recall that a measure on a smooth manifold
  $N$ is {\em smooth}\index{measure!smooth} if any local chart, it is
  absolutely continuous with respect to the Lebesgue measure, with
  smooth Radon-Nikodym derivative.} and that the geodesic flow on
$T^1M$ is mixing with exponential speed for the Sobolev regularity for
the potential $F$. Then there exist $\ell\in\NN$ and $\kappa'>0$ such
that for all nonnegative maps $\psi^\pm\in \C_{\rm c}^\ell(T^1M)$,
we have, as $t\ra+\infty$,
\begin{multline*}
\frac{\delta\;\|m_F\|}{e^{\delta\, t}}
\sum_{v,\,w\in T^1M} m_{\normalout\D^-}(v)\;m_{\normalin\D^+}(w)\; n_t(v,w) \;
\psi^-(v)\,\psi^+(w)\\=
\int_{T^1M}\psi^-d\sigma^+_{\D^-} \int_{T^1M}\psi^+d\sigma^-_{\D^+}
+\operatorname{O}(e^{-\kappa't}\|\psi^-\|_\ell\,\|\psi^+\|_\ell)\,.
\end{multline*}

Furthermore, if $\D^-$ and $\D^+$ respectively have nonzero finite
outer and inner skinning measures, and if $(\wt M, \Ga,\wt F)$
satisfies the conditions of (1) or of (2) above, then there exists
$\kappa''>0$ such that, as $t\ra+\infty$,
$$
\N_{\D^-,\,\D^+,\, F}(t)=
\frac{\|\sigma^+_{\D^-}\|\;\|\sigma^-_{\D^+}\|}{\delta\;\|m_F\|}\;
e^{\delta\,t}\big(1+\operatorname{O}(e^{-\kappa'' t})\big)\;.
$$
\etheo

The maps $\operatorname{O}(\cdot)$ depend on $\wt M,\Ga,F,\D,$ and
the speeds of mixing. 
The proof is a generalisation to nonzero potential of
\cite[Theo.~15]{ParPau16ETDS}.

\medskip
\dem We will follow the proofs of Theorem \ref{theo:mainequidup} and
Corollary \ref{coro:mainequicountdown} to prove generalisations of the
assertions (1) and (2) by adding to these proofs a regularisation
process of the test functions $\wt\phi^\pm_\eta$ as for the deduction
of Theorem \ref{theo:expratesequid} from Theorem \ref{theo:equid}.  We
will then deduce the last statement from these generalisations, again
using this regularisation process.

Let $\beta$ be either $\alpha\in\;]0,1]$ small enough in the H\"older
regularity case or $\ell\in\NN$ large enough in the Sobolev regularity
case.  We fix $i\in I^-$, $j\in I^+$, and we use the notation
$D^\pm,\alpha_\ga,\lambda_\ga$ and $\wt\sigma^\pm$ of Equation
\eqref{eq:notationstep1}. Let $v^\pm_\ga\in \normalmp D^\pm$ be the
initial and terminal tangent vectors to $\alpha_\ga$ and $\ga^{-1}
\alpha_\ga$ respectively.  Let $\wt \psi^\pm\in\C^\beta_{\rm c}
(\normalmp D^\pm)$. Under the assumptions of Assertion (1) or of
Assertion (2), we first prove the following avatar of Equation
\eqref{eq:reducnarrowup}, indicating only the required changes in its
proof: there exists $\kappa_0>0$ (independent of $\wt \psi^\pm$) such
that, as $T\ra+\infty$,
\begin{align}
&\delta\;\|m_F\|\;e^{-\delta\, T}
\sum_{\ga\in\Ga,\,0<\lambda_\ga\leq T} \;e^{\int_{\alpha_\ga}\wt F}\; 
\wt\psi^-(v^-_\ga)\,\wt\psi^+(v^+_\ga)
\nonumber\\=\; &
\int_{\normalout D^-} \wt\psi^-\, d\wt\sigma^+
\int_{\normalin D^+} \wt\psi^+\, d\wt\sigma^-
+\operatorname{O}(e^{-\kappa_0T}\|\wt\psi^-\|_\beta\,\|\wt\psi^+\|_\beta)
\label{eq:geneth13cor19errsimp}\,.
\end{align}
By  Lemma \ref{lem:fpmholder} and  the H\"older regularity of the
strong stable and unstable foliations under the assumptions of
Assertion (1), or by the smoothness of the boundary of $D^\pm$ under
the assumptions of Assertion (2), the maps $f^\pm_{D^\mp}:
\V_{\eta,\,R}^\pm(\normalpm D^\mp)\ra \normalpm D^\mp$ are
respectively H\"older-continuous or smooth fibrations, whose fiber over
$w\in \normalpm D^\mp$ is exactly $V_{w,\,\eta,\,R}^\pm$.  By applying
leafwise the regularisation process described in the proof of
Theorem \ref{theo:expratesequid} to characteristic functions, there
exist a constant $\kappa'_1>0$ and $\chi^\pm_{\eta,\,R}\in
\C^\beta(T^1\wt M)$ such that

$\bullet$~ $\|\chi^\pm_{\eta,\,R}\|_\beta=\bigO(\eta^{-\kappa'_1})$,

$\bullet$~ 
$\mathbbm{1}_{\V^\mp_{\eta \,e^{-\bigO(\eta)},
\,R\,e^{-\bigO(\eta)}}(\normalmp D^{\pm})}
\leq \chi^\pm_{\eta,\,R}\leq \mathbbm{1}_{\V^\mp_{\eta,\,R}(\normalmp D^{\pm})}$,

$\bullet$~ for every $w\in\normalmp D^{\pm}$, we have 
$$
\int_{V^\mp_{w,\,\eta,\,R}}\chi^\pm_{\eta,\,R}\,d\nu_w^\pm=
\nu_w^\pm(V^\mp_{w,\,\eta,\,R}) \,e^{-\bigO(\eta)}=
\nu_w^\pm(V^\mp_{w,\,\eta \,e^{-\bigO(\eta)},\,R \,e^{-\bigO(\eta)}}) 
\,e^{\bigO(\eta)}\,.
$$

We now define the new test functions (compare with Section
\ref{subsec:skinningwithpot}). For every $w\in \normalmp D^{\pm}$, let
$$
H^\pm_{\eta,\,R}(w)=
\frac{1}{\int_{V^\mp_{w,\,\eta,\,R}}\chi^\pm_{\eta,\,R}\,d\nu_w^\pm} \,.
$$
Let $\Phi^\pm_\eta:T^1\wt M\ra \RR$ be the map defined by
$$
\Phi^\pm_\eta=(H^\pm_{\eta,\,R}\; \wt\psi^\pm)\circ f^\mp_{D^\pm}
\;\; \chi^\pm_{\eta,\,R} \,.
$$
The support of this map is contained in $\V^\mp_{\eta,\,R}(\normalmp
D^{\pm})$. Since $M$ is compact in Assertion (1) and by homogeneity in
Assertion (2), if $R$ is large enough, by the definitions of the
measures $\nu_w^\pm$, the denominator of $H^\pm_{\eta,\,R}(w)$ is at
least $c \,\eta$ where $c>0$. The map $H^\pm_{\eta,\,R}$ is hence
H\"older-continuous under the assumptions of Assertion (1), and it  is
smooth under the assumptions of Assertion (2). Therefore
$\Phi^\pm_\eta\in \C^\beta(T^1\wt M)$ and there exists a constant
$\kappa'_2>0$ such that
$$
\|\Phi^\pm_\eta\|_\beta=
\bigO(\eta^{-\kappa'_2}\|\wt \psi^\pm\|_\beta)\,.
$$
As in Lemma \ref{lem:integrable}, the functions $\Phi^\pm_\eta$ are
measurable, nonnegative and satisfy
$$
\int_{T^1\wt M} \Phi^\pm_\eta\;d\wt m_F=
\int_{\normalmp D^{\pm}}\wt \psi^\pm\,d\wt\sigma^\mp\,.
$$  
As in  the proof of Theorem
\ref{theo:mainequidup}, we will estimate in two ways the quantity
\begin{equation}\label{eq:defiIetapmTerror}
I_{\eta}(T)=\int_0^{T}e^{\delta\,t}\;\sum_{\ga\in\Ga}\;
\int_{T^1\wt M}(\Phi^-_\eta\circ\flow{-t/2})\;
(\Phi^+_\eta\circ\flow{t/2}\circ\ga^{-1})\;d\wt m_F\;dt\,.
\end{equation}

We first apply the mixing property, now with exponential decay of
correlations, as in Step 2 of the proof of Theorem
\ref{theo:mainequidup}. For all $t\geq 0$, let
$$
A_\eta(t)=\sum_{\ga\in\Ga}\;
\int_{v\in T^1\wt M}\Phi^-_\eta(\flow{-t/2}v)\;
\Phi^+_\eta(\flow{t/2}\ga^{-1}v)\;d\wt m_F(v)\,.
$$
Then with $\kappa>0$ as in the definitions of the exponential mixing
for the H\"older or Sobolev regularity, we have
\begin{align*}
A_\eta(t)&=\frac{1}{\|m_F\|}\;
\int_{T^1\wt M}\Phi^-_\eta\,d\wt m_F\;
\int_{T^1\wt M}\Phi^+_\eta\,d\wt m_F\;+\;
\bigO\big(e^{-\kappa\, t}\|\Phi^-_\eta\|_\beta\|\Phi^+_\eta\|_\beta\big)
\\ & =\frac{1}{\|m_F\|}\;
\int_{\normalout D^{-}}\wt \psi^-\,d\wt\sigma^+
\int_{\normalin D^{+}}\wt \psi^+\,d\wt\sigma^-\;+\;
\bigO\big(e^{-\kappa \,t}\eta^{-2\kappa'_2}
\|\wt \psi^-\|_\beta\|\wt \psi^+\|_\beta\big)\,.
\end{align*}
Hence by integrating, 
\begin{equation}\label{eq:firstestimerror}
I_{\eta}(T)=
\frac{e^{\delta\,T}}{\delta\,\|m_F\|}\;\Big(
\int_{\normalout D^{-}}\wt \psi^-\,d\wt\sigma^+
\int_{\normalin D^{+}}\wt \psi^+\,d\wt\sigma^-\;+\;
\bigO\big(e^{-\kappa\, T}\eta^{-2\kappa'_2}
\|\wt \psi^-\|_\beta\|\wt \psi^+\|_\beta\big)\Big)\,.
\end{equation}

Now, as in Step 3 of the proof of Theorem \ref{theo:mainequidup}, we
exchange the integral over $t$ and the summation over $\ga$ in the
definition of $I_{\eta}(T)$, and we estimate the integral term
independently of $\ga$:
$$
I_{\eta}(T)=\sum_{\ga\in\Ga}\;\int_0^{T}e^{\delta\,t}\;
\int_{T^1\wt M}(\Phi^-_\eta\circ\flow{-t/2})\;
(\Phi^+_\eta\circ\flow{t/2}\circ\ga^{-1})\;d\wt m_F\;dt\,.
$$

Let $\wh\Phi^\pm_\eta=H^\pm_{\eta,\,R}\circ f^\mp_{D^\pm} \;
\chi^\pm_{\eta,\,R}$, so that $\Phi^\pm_\eta=\wt\psi^\pm\circ
f^\mp_{D^\pm}\;\wh\Phi^\pm_\eta$. By the last two properties of the
regularised maps $\chi^\pm_{\eta,\,R}$, we have, with
$\phi^\pm_{\eta',\,\eta'',\,\Omega^\pm}$ defined as in Equation
\eqref{eq:defiphi},
\begin{equation}\label{eq:controlwhPhi}
\phi^\pm_{\eta\,e^{-\bigO(\eta)},\,R \,e^{-\bigO(\eta)},\,\normalmp D^{\pm}}
\,e^{-\bigO(\eta)}\leq \wh\Phi^\pm_\eta\leq 
\phi^\pm_\eta\,e^{\bigO(\eta)}\;.
\end{equation}

If $v\in T^1\wt M$ belongs to the support of $(\Phi^-_\eta \circ
\flow{-t/2})\; (\Phi^+_\eta\circ \flow{t/2} \circ \ga^{-1})$, then we
have $v\in \flow{t/2}\V^+_{\eta,\,R}(\normalout D^{-})\cap \flow{-t/2}
\V^-_{\eta,\,R}(\ga\normalin D^{+})$. Hence the properties (i), (ii)
and (iii) of Step 3M of the proof of Theorem \ref{theo:mainequidup}
still hold (with $\Omega_-=\normalout D^{-}$ and $\Omega_+=\normalin
(\ga D^{+})$). In particular, if $w^-=f^+_{D^-}(v)$ and $w^+=f^-_{\ga
  D^+} (v)$, we have, by Assertion (iii) in Step 3M of the proof
of Theorem \ref{theo:mainequidup},\footnote{See also the picture at the
beginning of the proof of Lemma \ref{lem:variaboul}.} that
$$
d(w^\pm,v^\pm_\ga)=\bigO(\eta+e^{-\lambda_\ga/2})\,.
$$
Hence, with $\kappa'_3=\alpha$ in the H\"older case and $\kappa'_3=1$ in
the Sobolev case (we may assume that $\ell\geq 1$), we have
$$
|\,\wt\psi^\pm(w^\pm)-\wt\psi^\pm(v^\pm_\ga)\,|=
\bigO((\eta+e^{-\lambda_\ga/2})^{\kappa'_3}\|\wt\psi^\pm\|_\beta)\,.
$$
Therefore there exists a constant $\kappa'_4>0$ such that 
\begin{align*}
I_{\eta}(T)=\sum_{\ga\in\Ga}&\;\;\;\big(\wt\psi^-(v^-_\ga)\wt\psi^+(v^+_\ga)+
\bigO((\eta+e^{-\lambda_\ga/2})^{\kappa'_4}
\|\wt\psi^-\|_\beta\|\wt\psi^+\|_\beta)\big)\times
 \\ & \int_0^{T}e^{\delta\,t}\;\int_{v\in T^1\wt M} 
\;\wh\Phi^-_{\eta}(\flow{-t/2}v)\;
\wh\Phi^+_{\eta}(\ga^{-1}\flow{t/2}v)\;d\wt m_F(v)\;dt\,.
\end{align*}

Now, using the inequalities \eqref{eq:controlwhPhi}, Equation
\eqref{eq:geneth13cor19errsimp} follows as in Steps 3M and 4M of
the proof of Theorem \ref{theo:mainequidup}, by taking
$\eta=e^{-\kappa'_5 T}$ for some $\kappa'_5>0$ and using the effective
control given by Equation \eqref{eq:step5bis} in Step 4M.

\bigskip
In order to prove Assertions (1) and (2) of Theorem
\ref{theo:expratecount}, we may assume that the supports of $\psi^\pm$
are small enough, say contained in $B(x^\pm,\epsilon)$ for some
$x^\pm\in T^1M$ and $\epsilon$ small enough. Let $\wt x^\pm$ be lifts
of $x^\pm$ and let $\wt\psi^\pm\in \C^\beta_{\rm c}(T^1\wt M)$ with
support in $B(\,\wt x^\pm,\epsilon)$ be such that $\wt\psi^\pm =
\psi^\pm\circ Tp$ on $B(\,\wt x^\pm,\epsilon)$ where $p:\wt M\ra M$ is
the universal Riemanian orbifold cover. By a finite summation argument
since $\Ga\bs I^\pm$ are finite, and by Equation
\eqref{eq:geneth13cor19errsimp}, we have
\begin{align}
&\delta\;\|m_F\|\;e^{-\delta\, T}
\sum_{\substack{i\in I^-/\sim,\; j\in I^+/\sim,\;\ga\in \Ga\\
0<\lambda_{i,\,\ga j}\leq T}} \;e^{\int_{\alpha_\ga}\wt F}\; 
\wt\psi^-(v^-_\ga)\,\wt\psi^+(v^+_\ga)
\nonumber\\=\; &
\int_{\normalout D^-} \wt\psi^-\, d\wt\sigma^+
\int_{\normalin D^+} \wt\psi^+\, d\wt\sigma^-
+\bigO(e^{-\kappa_0T}\|\wt\psi^-\|_\beta\,\|\wt\psi^+\|_\beta)
\label{eq:geneth13cor19errsimpexpand}\,.
\end{align}

Assertions (1) and (2) are deduced from this equation in the same way
that Corollary \ref{coro:mainequicountdown} is deduced from Theorem
\ref{theo:mainequidup}. Taking the functions $\psi_k^\pm$ to be the
constant functions $1$ in Assertion (1) gives the last statement of
Theorem \ref{theo:expratecount} under the assumptions of Assertion
(1). An approximation argument gives the result under the assumptions
of Assertion (2).  
\cqfd

\section{Equidistribution and counting for quotient simplicial
and metric  trees}
\label{subsect:equicountmetricgraphgroup}

In this Section, we assume that $X$ is the geometric realisation of a
locally finite metric tree without terminal vertices $(\XX,\lambda)$,
and that $\Ga$ is a (nonelementary discrete) subgroup of $\Aut(\XX,
\lambda)$.  Let $\wt c:E\XX\ra\RR$ be a system of conductances for
$\Ga$, and let $c:\Ga\bs E\XX\ra\RR$ be its quotient function. We
assume in this Section that the potential $\wt F$ is the potential
$\wt F_c$ associated\footnote{See Section \ref{subsec:cond}.} with
$c$. Let $\delta_c=\delta_{\Ga,\,F_c}$ be the critical exponent of
$(\Ga,F_c)$ and let $\wt m_{c}=\wt m_{F_c}$ and $m_c=m_{F_c}$ be the
Gibbs measures of $F_c$ for the continuous time geodesic flow on
respectively $\G X$ and $\Ga\bs\G X$, as well as for the discrete time
geodesic flow on respectively $\G \XX$ and $\Ga\bs\G \XX$ when $(\XX,
\lambda)$ is simplicial, that is, if $\lambda$ is constant with value
$1$.

Let $\DD^\pm$ be simplicial subtrees of $\XX$, with the edge length
map induced by $\lambda$, and $D^\pm=|\DD^\pm|_\lambda$ its geometric
realisation, such that the $\Ga$-equivariant families $\D^\pm=(\ga
D^\pm)_{\ga\in \Ga/ \Ga_{\DD^\pm}}$ are locally finite in
$X$.\footnote{We leave to the reader the extension to more general
  locally finite families of subtrees, as for instance finite unions
  of those above.}

For all $\ga,\ga'$ in $\Ga$ such that $\ga D^-$ and $\ga'D^+$ are
disjoint, we denote by $\alpha_{\ga,\,\ga'}$ the common perpendicular
from $\ga D^-$ to $\ga' D^+$ (which is an edge path in $\XX$),
with length $\lambda_{\ga,\,\ga'}=d(\ga D^-,\ga' D^+)\in\NN$, and
by $\alpha^\pm_{\ga,\,\ga'} \in \gengeod\, \XX$ its parametrisations as
in the beginning of Chapter \ref{sec:equidarcs}: it is the unique map
from $\RR$ to $X$ such that $\alpha^-_{\ga,\,\ga'} (t)\in \ga V\DD^-$
is the origin $o(\alpha_{\ga,\,\ga'})$ of the edge path
$\alpha_{\ga,\,\ga'}$ if $t\leq 0$, $\alpha^-_{\ga,\,\ga'} (t)\in \ga'
V\DD^+$ is the endpoint $t(\alpha_{\ga,\,\ga'})$ of the edge path
$\alpha_{\ga,\,\ga'}$ if $t\geq \lambda_{\ga,\,\ga'}$, and
${\alpha^-_{\ga,\,\ga'}}|_{ [0,\,\lambda_{\ga,\,\ga'}]}$ is the
shortest geodesic arc starting from a point of $\ga D^-$ and ending
at a point of $\ga'D^+$.

For all $\ga,\ga'$ in $\Ga$ such that $\ga D^-$ and $\ga'D^+$ are
disjoint, we define the {\em multiplicity}\index{multiplicity}%
\index{common perpendicular!multiplicity} of the common perpendicular
$\alpha_{\ga,\,\ga'}$ from $\ga D^-$ to $\ga' D^+$ as
\begin{equation}\label{eq:multiplicitytree}
m_{\ga,\,\ga'}=
\frac 1{\card(\ga\Ga_{\DD^-}\ga^{-1}\cap\ga'\Ga_{\DD^+}{\ga'}^{-1})}\,.
\end{equation}
Note that $m_{\ga,\,\ga'}=1$ for all $\ga,\ga'\in\Ga$ when $\Ga$
acts freely on $E\XX$ (for instance when $\Ga$ is torsion-free).
Generalising the definition for simplicial trees in Section
\ref{subsec:equidcommperpdiscrtime}, we set
$$
\wt c(\alpha)= \sum_{i=1}^k \wt c\,(e_i)\,\lambda(e_i)\,,
$$ 
for any edge path $\alpha=(e_1,\dots, e_k)$ in $\XX$.

For $n\in\NN-\{0\}$, let
\begin{equation}\label{eq:NDDn}
\N_{\DD^-,\,\DD^+}(n)=
\sum_{\substack{[\ga]\in\, \Ga_{\DD^-}\bs\Ga/\Ga_{\DD^+}\\
0< d(D^-,\ga D^+)\leq n}} m_{e,\,\ga}\;e
^{\wt c\,(\alpha_{e,\,\ga})}\;,
\end{equation}
where $d$ is the the distance on $X= |\XX|_\lambda$.  When $\Ga$ is
torsion-free, $\N_{\DD^-,\,\DD^+} (n)$ is the number of edge paths in
the graph $\Ga\bs \XX$ of length at most $n$, starting by an outgoing
edge from the image of $\DD^-$ and ending by the opposite of an
outgoing edge from the image of $\DD^+$, with multiplicities coming
from the fact that $\Ga_{D^\pm}\bs \DD^\pm$ is not assumed to be
embedded in $\Ga\bs \XX$, and with weights coming from the
conductances.

In the next results, we distinguish the continuous time case (Theorem
\ref{theo:counting}) from the discrete time case (Theorem
\ref{theo:countingdiscrete}). We leave to the reader the versions
without the assumption $\delta_c>0$, giving for every
$\tau\in\NN-\{0\}$ an asymptotic on
$$
\N_{\DD^-,\,\DD^+,\,\tau}(n)=
\sum_{\substack{[\ga]\in\, \Ga_{\DD^-}\bs\Ga/\Ga_{\DD^+}\\
n-\tau< d(D^-,\ga D^+)\leq n}} m_{e,\,\ga}\;e^{\wt c\,(\alpha_{e,\,\ga})}
\;.
$$

When $\Ga\bs X$ is compact, $c=0$ and $\DD^\pm$ are reduced to points,
the counting results in Theorems \ref{theo:counting} and
\ref{theo:countingdiscrete} are proved in \cite{Guillope94}.  When
$\DD^\pm$ are singletons, Theorem \ref{theo:counting} is due to
\cite{Roblin03} if $c=0$. Otherwise,  the result seems to be new.

\btheo \label{theo:counting} Let $(\XX,\lambda)$, $\Ga$, $\DD^\pm$ and
$c$ be as in the beginning of this Section.  Assume that the critical
exponent $\delta_c$ is finite and positive, that the skinning measures
$\sigma^\pm_{\D^\mp}$ are finite and nonzero, and that the Gibbs
measure $m_{c}$ is finite and mixing for the continuous time geodesic
flow. Then as $t\ra+\infty$, the measures
$$
\delta_c\;\|m_{c}\|\;e^{-\delta_c\, t}
\sum_{\substack{[\ga]\,\in \Ga_{\DD^-}\bs\Ga/\Ga_{\DD^+}\\
0<d(\DD^-,\ga \,\DD^+)\leq t}} 
m_{e,\,\ga}\;e^{\wt c\,(\alpha_{e,\,\ga})}\;
\Dirac_{\Ga \alpha^-_{e,\ga}}
\otimes\Dirac_{\Ga\alpha^+_{\ga^{-1},e}}
$$ 
narrow converge to $\sigma^+_{\D^-}\otimes \sigma^-_{\D^+}$ in
$\Ga\bs \gengeod X\times \Ga\bs\gengeod X$, and
$$
\N_{\DD^-,\,\DD^+}(t)\sim 
\frac{\|\sigma^+_{\D^-}\|\;\|\sigma^-_{\D^+}\|}
{\delta_c\;\|m_{c}\|}\;\;e^{\delta_c\, t}\,.
$$
\etheo

\dem 
By Theorem \ref{theo:mainequidup}, we have
$$
\lim_{t\ra+\infty} \;\delta_c\;\|m_{c}\|\;e^{-\delta_c\, t}
\sum_{\substack{(a,b,\ga)\in 
\Ga/\Ga_{\DD^-}\times\Ga/\Ga_{\DD^+}\times\Ga\\ 
0<d(a \DD^-,\,\ga b \DD^+)\leq t}} \;e^{\int_{\alpha_{a,\ga b}}\wt F_c}\; 
\Dirac_{\alpha^-_{a,\ga b}} \otimes\Dirac_{\alpha^+_{\ga^{-1}a,b}}\;=\; 
\wt\sigma^+_{\D^-}\otimes \wt\sigma^-_{\D^+}\;,
$$ 
for the weak-star convergence on $\gengeod \,\XX\times
\gengeod\, \XX$.

The group $\Ga\times\Ga$ acts on $\Ga/\Ga_{\DD^-}\times
\Ga/\Ga_{\DD^+}\times\Ga$ by
$$
(a',b')\cdot (a,b,\ga)= (a'a,b'b,a'\ga(b')^{-1})\;.
$$ 
and the map from the discrete set $\Ga/\Ga_{\DD^-}\times \Ga/
\Ga_{\DD^+}\times\Ga$ to $\gengeod\,\XX\times \gengeod \,\XX$ which
sends $(a,b,\ga)$ to $(\alpha^-_{a,\ga b} \,, \alpha^+_{\ga^{-1}a,b})$
is $(\Ga\times\Ga)$-equivariant. In particular, the pushforward of
measures by this map sends the unit Dirac mass at $(a,b,\ga)$ to
$\Dirac_{\alpha^-_{a,\ga b}} \otimes\Dirac_{\alpha^+_{\ga^{-1}a,b}}$.

Every orbit of $\Ga\times\Ga$ on $\Ga/\Ga_{\DD^-}\times
\Ga/\Ga_{\DD^+} \times\Ga$ has a representative of the form
$(\Ga_{\DD^-}, \Ga_{\DD^+},\ga)$ for some $\ga\in\Ga$, since
$(a,b)\cdot(\Ga_{\DD^-}, \Ga_{\DD^+},a^{-1}\ga\,b)= (a\Ga_{\DD^-},
b\Ga_{\DD^-}, \ga)$. Furthermore the double class in
$\Ga_{\DD^-}\bs\Ga/\Ga_{\DD^+}$ of such a $\ga$ is uniquely defined,
and the stabiliser of $(\Ga_{\DD^-}, \Ga_{\DD^+},\ga)$ has cardinality
$|\Ga_{\DD^-}\cap \ga \Ga_{\DD^+} \ga^{-1}|$, since
$(a,b)\cdot(\Ga_{\DD^-}, \Ga_{\DD^+},\ga)= (\Ga_{\DD^-},
\Ga_{\DD^+},\ga')$ if and only if $a\in \Ga_{\DD^-}$, $b\in
\Ga_{\DD^+}$ and $a\ga \,b^{-1}=\ga'$. When $\ga'=\ga$, this
happens if and only if $b=\ga^{-1}a\ga$ and $a\in \Ga_{\DD^-}\cap
\ga \Ga_{\DD^+} \ga^{-1}$.

By using the properties recalled at the beginning of Chapter
\ref{sec:equidcountdownstairs} on the narrow convergence of induced
measures, and since no compactness assumptions were made in Equation
\eqref{eq:reducnarrowup} on $\Omega^\pm$, the measures
$$
\delta_c\;\|m_{c}\|\;e^{-\delta_c\, t}
\sum_{\substack{[\ga]\,\in \Ga_{\DD^-}\bs\Ga/\Ga_{\DD^+}\\
0<d(\DD^-,\ga \,\DD^+)\leq t}} 
\frac{1} {|\Ga_{\DD^-}\cap \ga \Ga_{\DD^+}\ga^{-1}|}
\;e^{\int_{\Ga\alpha_{e,\ga}}F_c}\; 
\Dirac_{\Ga \alpha^-_{e,\ga}}
\otimes\Dirac_{\Ga\alpha^+_{\ga^{-1},e}}
$$ 
hence narrow converge as $t\ra+\infty$ to $\sigma^+_{\D^-}\otimes
\sigma^-_{\D^+}$ in $\Ga\bs\gengeod \,\XX \;\times \Ga\bs \gengeod\,
\XX$.  By applying this convergence to the constant function $1$, and
by the finiteness and nonvanishing of $\sigma^+_{\D^-}$ and
$\sigma^-_{\D^+}$, the result follows using the defining property of
the potential $F_c$, see Proposition \ref{prop:integpotconduct}.
\cqfd

\medskip
In the remainder of this Section, we consider simplicial trees with
the discrete time geodesic flow. 

\btheo\label{theo:countingdiscrete} Let $(\XX,\lambda)$, $\Ga$, $\wt
c$ and $\DD^\pm$ be as in the beginning of this Section, with
$\lambda$ constant with value $1$.  Assume that the critical exponent
$\delta_c$ is finite and positive. If the Gibbs measure $m_{c}$ on the
space $\Ga\bs\G \XX$ of discrete geodesic lines modulo $\Ga$ is finite
and mixing for the discrete time geodesic flow, and if the skinning
measures $\sigma^\pm_{\D^\mp}$ are finite and nonzero, then as
$n\ra+\infty$, the measures
$$
\frac{e^{\delta_c}-1}{e^{\delta_c}}\;\|m_{c}\|\;e^{-\delta_c\, n}
\sum_{\substack{[\ga]\,\in \Ga_{\DD^-}\bs\Ga/\Ga_{\DD^+}\\
0<d(\DD^-,\ga \,\DD^+)\leq n}} 
m_{e,\,\ga}\;e^{\wt c\,(\alpha_{e,\,\ga})}\;
\Dirac_{\Ga \alpha^-_{e,\ga}}
\otimes\Dirac_{\Ga\alpha^+_{\ga^{-1},e}}
$$ 
narrow converge to $\sigma^+_{\D^-}\otimes
\sigma^-_{\D^+}$ in $\Ga\bs\gengeod\,\XX\times
\Ga\bs\gengeod\,\XX$
and
$$
\N_{\DD^-,\,\DD^+}(n)\sim 
\frac{e^{\delta_c}\;\|\sigma^+_{\D^-}\|\;\|\sigma^-_{\D^+}\|}
{(e^{\delta_c}-1)\;\|m_{c}\|}\;\;e^{\delta_c\, n}\,.
$$
\etheo

\dem 
The claims follow as in Theorem \ref{theo:counting}, replacing Theorem
\ref{theo:mainequidup} by Theorem \ref{theo:discretemainequidup}.
\cqfd

\brema A common perpendicular in a simplicial tree is, in the language
of graph theory, a non-backtracking walk. Among other
applications,\footnote{when restricting to groups $\Ga$ acting freely,
  which is never the case if $\Ga$ is a nonuniform lattice in the tree
  $\XX$, that is, when the quotient graph of groups $\Ga\bs\XX$ is
  infinite but has finite volume} Theorem \ref{theo:countingdiscrete}
gives a complete asymptotic solution to the problem of counting
non-backtracking walks from a given vertex to a given vertex of a
(finite) nonbipartite graph. See Theorem
\ref{theo:countingdiscreteeven} for the corresponding result in
bipartite graphs, and for example \cite[Th.~1.1]{AloBenLubSod07},
\cite[p.~4290,4302]{AngFriHoo15}, \cite[L.~2.3]{Friedman08},
\cite[Prop.~6.4]{Sodin07} for related results.  Anticipating on the
error terms that we will give in Section
\ref{subsect:errormetricgraphgroup}, note that the paper
\cite[Th.~1.1]{AloBenLubSod07} for instance gives a precise speed
using spectral properties, more precise than the ones we obtain.
\erema
 
\bexems\label{ex:cpregularsubgraphs} (1) Let $\XX,\Ga,\wt c$ be as in
Theorem \ref{theo:countingdiscrete}, and let $\DD^-=\{x\}$ and
$\DD^+=\{y\}$ for some $x,y\in V\XX$. If the Gibbs measure $m_{c}$ is
finite and mixing for the discrete time geodesic flow on $\Ga\bs\G
\XX$, then we have a discrete time version of Roblin's simultaneous
equidistribution theorem with potential,\footnote{See Corollary
  \ref{coro:caspointempcont} for the continuous time version.} and the
number $\N_{x,\,y}(n)$ of nonbacktracking edge paths of length at most
$n$ from the image of $x$ to the image of $y$ (counted with weights
and multiplicities) satisfies, by Equation \eqref{eq:skinning singleton},
$$
\N_{x,\,y}(n)\sim 
\frac{e^{\delta_c}\;\|\mu^+_{x}\|\;\|\mu^-_{y}\|}
{(e^{\delta_c}-1)\;\|m_{c}\|\,|\Ga_x|\,|\Ga_y|}\;\;e^{\delta_c\, n}\,.
$$

\smallskip\noindent (2) If $\YY$ is a finite connected nonbipartite
$(q+1)$-regular graph (with $q\geq 2$) and $\YY^\pm$ are points, then
the number of nonbacktracking edge paths from $\YY^-$ to $\YY^+$ of
length at most $n$ is equivalent as $n\ra+\infty$ to
\begin{equation}\label{eq:numberofwalksingraph}
\frac{q+1}{q-1}\,\frac{q^n}{|V\YY|} +\bigO(r^{n})\,
\end{equation}
for some $r<q$.  Indeed, by Theorem \ref{theo:countingdiscrete} with
$\XX$ the universal cover of $\YY$, $\Ga$ its covering group and
$c=0$, we have $\delta_c=\ln q$ and $m_{c}$ is the Bowen-Margulis
measure, so that normalising the Patterson measures to be probability
measures, we have $\|m_{c}\| =\frac{q}{q+1}\;|V\YY|$ by Equation
\eqref{eq:relatmBMvolhomog}. We refer to Section
\ref{subsect:errormetricgraphgroup} (see Remark (i) following the
proof of Theorem \ref{theo:equidsimplicialperperrorterm}) for the
error term.

Let $\YY$ be the figure $8$-graph with a single vertex and four directed
edges, and let $\YY^\pm$ be the singleton consisting of its
vertex. In this simple example, it is easy to count by hand that the
number of loops of length exactly $n$ without backtracking in $\YY$ is
$4\,3^{n-1}$. Thus the number $\N(n)$ of common perpendiculars of the
vertex to itself of length at most $n$ is by a simple geometric sum
$2(3^n-1)$. This agrees with Equation \eqref{eq:numberofwalksingraph}
that gives $\N(n)\sim 2\,3^n$ as $n\to+\infty$.


\smallskip \noindent (3) Let $\YY$ be a finite connected nonbipartite
$(q+1)$-regular graph (with $q\geq 2$). Let $\YY^\pm$ be regular
connected subgraphs of degrees $q^\pm\geq 0$, with $q^\pm<q+1$. Then
the number $\N(n)$ of nonbacktracking edge paths of length at most $n$ starting
transversally to $\YY^-$ and ending transversally to $\YY^+$ satisfies
$$
\N(n)=\frac{(q+1-q^-)(q+1-q^+)\,|V\YY^-|\;|V\YY^+|}{(q^2-1)\,|V\YY|}\;q^{n} 
+\bigO(r^{n})
$$ 
for some $r<q$.  This is a direct consequence of Theorem
\ref{theo:countingdiscrete}, using Proposition \ref{prop:computBM} (3)
and Proposition \ref{prop:computskin} (3), again refering to Section
\ref{subsect:errormetricgraphgroup}  for the error term.
\eexems

We refer for instance to Chapters \ref{sec:BruhatTitstrees} and
\ref{sec:rattionalequid} for examples of counting results in graphs of
groups where the underlying graph is infinite.

\medskip
In some applications (see the examples at the end of this Section), we
encounter bipartite simplicial graphs and, consequently, their
discrete time geodesic flow is not mixing. The following result
applies in this context.

Until the end of this section, we assume that the simplicial tree
$\XX$ has a $\Ga$-invariant structure of a bipartite graph, and we
denote by $V\XX = V_1\XX \sqcup V_2\XX$ the corresponding partition of
its set of vertices. For every $i\in\{1,2\}$, we denote by
$\gengeod_i\XX$ the space of generalised discrete geodesic lines
$\ell\in \gengeod\,\XX$ such that $\ell(0)\in V_i\XX$, so that we have
a partition $\gengeod\,\XX= \gengeod_1\XX\sqcup \gengeod_2\XX$. Note
that if the basepoint $x_0$ lies in $V_i\XX$, then $\G_{\rm even}\XX$
is equal to $\gengeod_i\XX\cap\G \XX$. Let $\sigma^\pm_{\D^\mp,\,i}=
\sigma^\pm_{\D^\mp}\mid_{\Ga\bs\gengeod_i\XX}\,$. For all
$i,j\in\{1,2\}$, we define
$$
\N_{\DD^-,\,\DD^+,\,i,\,j}(n)=
\sum_{\substack{[\ga]\in\, \Ga_{\DD^-}\bs\Ga/\Ga_{\DD^+}\\
0< d(D^-,\ga D^+)\leq n\\
o(\alpha_{e,\ga})\in V_i\XX,\;t(\alpha_{e,\ga})\in V_j\XX}} 
m_{e,\,\ga}\;e^{\wt c\,(\alpha_{e,\,\ga})}\;.
$$

\btheo\label{theo:countingdiscreteeven} Let $(\XX,\lambda)$, $\Ga$ and
$c$ be as in the beginning of this Section, with $\lambda$ constant
with value $1$.  Assume that the critical exponent $\delta_c$ is
finite and positive. Assume that $\XX$ has a $\Ga$-invariant structure
of a bipartite graph as above, and that the restriction to
$\Ga\bs\G_{\rm even}\XX$ of the Gibbs measure $m_{c}$ is finite and
mixing for the square of the discrete time geodesic flow. Then for all
$i,j\in\{1,2\}$ such that the measures $\sigma^-_{\D^-,\,i}$ and
$\sigma^+_{\D^+,\,j}$ are finite and nonzero, as $n$ tends to
$+\infty$ with $n\equiv i-j\mod 2$, the measures
$$
\frac{e^{2\,\delta_c}-1}{2\; e^{2\,\delta_c}}\;\|m_{c}\|\;e^{-\delta_c\, n}
\sum_{\substack{[\ga]\,\in \Ga_{\DD^-}\bs\Ga/\Ga_{\DD^+}\\
0<d(D^-,\ga \,D^+)\leq n\\
o(\alpha_{e,\ga})\in V_i\XX,\;t(\alpha_{e,\ga})\in V_j\XX}} 
m_{e,\,\ga}\;e^{\wt c\,(\alpha_{e,\,\ga})}\;
\Dirac_{\Ga \alpha^-_{e,\ga}}
\otimes\Dirac_{\Ga\alpha^+_{\ga^{-1},e}}
$$ 
narrow converge to $\sigma^+_{\D^-,\,i}\otimes \sigma^-_{\D^+,\,j}$
in $\Ga\bs\gengeod\,\XX\times \Ga\bs\gengeod\,\XX$ and
$$
\N_{\DD^-,\,\DD^+,\,i,\,j}(n)\sim 
\frac{2\;e^{2\,\delta_c}\;\|\sigma^+_{\D^-,\,i}\|\;\|\sigma^-_{\D^+,\,j}\|}
{(e^{2\,\delta_c}-1)\;\|m_{c}\|}\;\;e^{\delta_c\, n}\,.
$$
\etheo

\dem 
This Theorem is proved in the same way as the above Theorem
\ref{theo:countingdiscrete} using Theorem
\ref{theo:discretemainequidupeven}.  Note that we have a
$(\Ga\times\Ga)$-invariant partition
$$
\gengeod\,\XX\times \gengeod\,\XX=
\bigsqcup_{(i,\,j)\in\{1,\,2\}^2}\gengeod_i\XX\times \gengeod_j\XX\;,
$$ 
that $\alpha^-_{e,\,\ga}\in\gengeod_i\XX$ if and only if
$o(\alpha_{e,\,\ga}) \in V_i\XX$, and that $\alpha^+_{\ga^{-1},\,e} \in
\gengeod_j\XX$ if and only if $t(\alpha_{e,\,\ga})\in V_j\XX$, since
$\alpha^+_{\ga^{-1},\,e}(0)= \ga^{-1}\alpha^+_{e,\,\ga}(0)=
\ga^{-1}t(\alpha_{e,\,\ga})$.  
\cqfd

\bexems\label{ex:cpregularsubgraphsbi} (1) Let $\XX,\Ga,c$ be as in Theorem
\ref{theo:countingdiscreteeven}, and let $\DD^-=\{x\}$ and
$\DD^+=\{y\}$ for some vertices $x,y$ in the same $V_i\XX$ for
$i\in\{1,2\}$. If the restriction to $\Ga\bs\G_{\rm even}\XX$ of the
Gibbs measure $m_{c}$ is finite and mixing for the square of the
discrete time geodesic flow, then as $n$ is even and tends to
$+\infty$ ,
$$ 
\N_{\DD^-,\,\DD^+}(n)\sim\frac{2\;e^{2\,\delta_c}}{e^{2\,\delta_c}-1}
\;\frac{\|\mu^+_x\|\;\|\mu^-_y\|}{\|m_{c}\|\,|\Ga_x|\,|\Ga_y|}
\;\;e^{\delta_c\, n}\,.
$$ 
Indeed, we have $\N_{\DD^-,\,\DD^+}(n)=\N_{\DD^-,\,\DD^+,\,i,\,i}(n)$
and $\sigma^\pm_{\D^\pm,\,i}=\sigma^\pm_{\D^\pm}$, and we conclude as
in Example \ref{ex:cpregularsubgraphs} (1).

\medskip\noindent (2) Let $\YY$ be the complete biregular graph with
$q+1$ vertices of degree $p+1$ and $p+1$ vertices of degree $q+1$. Let
$\YY^\pm=\{y\}$ be a fixed vertex of degree $p+1$. Note that  $\YY$ being
bipartite, all common perpendiculars from $y$ to $y$ have even length
(the shortest one having length $4$). Then as $n$ is even
and tends to $+\infty$, we have
$$
\N_{\YY^-,\YY^+}(n) \sim\frac{q(p+1)}{(q+1)(pq-1)}(pq)^{n/2}\,.
$$ 
Indeed, the biregular tree $\XX_{p,q}$ of degrees $(p+1,q+1)$ is a
universal cover of $\YY$ with covering group $\Ga$ acting freely and
cocompactly, so that with $c=0$ we have $\delta_c=\ln\sqrt{pq}$ and
the Gibbs measure $m_{c}$ is the Bowen-Margulis measure $m_{\rm
  BM}$. If we normalise the Patterson density such that $\|\mu_y\|=
\frac{p+1}{\sqrt p}$, then by Proposition \ref{prop:computBM} (2), we
have $\|m_{\rm BM}\| = 2(p+1)(q+1)$. Thus the result follows from
the above Example (1). Note that if $p=q$, then
$$
\N_{\YY^-,\YY^+}(n)\sim \frac{q}{q^2-1}q^{n}\,,
$$ 
and the constant in front of $q^n$ is indeed different from that in
the nonbipartite case.  

\medskip\noindent (3) Let $\YY$ be a finite biregular graph with
vertices of degrees $p+1$ and $q+1$, where $p,q\geq 2$, and let
$V\YY=V_p\YY\sqcup V_q\YY$ be the corresponding partition.  If
$\YY^-=\{v\}$ where $v\in V_p\YY$ and $\YY^+$ is a cycle of length
$L\geq 2$, then as $N\to+\infty$, the number of common perpendiculars
of even length at most $2N$ from $\YY^-$ to $\YY^+$ is equivalent to
$$
\frac{L\;q\;(p-1)}{2\;(pq-1)\;|V_p\YY|}\;(pq)^N\;,
$$ 
and the number of common perpendiculars of odd length at most
$2N-1$ from $v$ to $\YY^+$ is equivalent to
$$
\frac{L\;(q-1)}{2\;(pq-1)\;|V_p\YY|}\;(pq)^{N}\;.
$$

\medskip
\dem The cycle $\YY^+$ has even length $L$ and has $\frac L2$ vertices
in both $V_p\YY$ and $V_q\YY$. A common perpendicular from $\YY^-$ to
$\YY^+$ has even length if and only if it ends at a vertex in
$V_p\YY$.

Let $\XX\ra \YY$ be a universal cover of $\YY$, whose covering group
$\Ga$ acts freely and cocompactly on $\XX$. Let $\DD^-=\{\wt v\}$
where $\wt v\in V_p\XX$ is a lift of $v$, and let $\DD^-$ be a
geodesic line in $\XX$ mapping to $\YY^+$. We use Theorem
\ref{theo:countingdiscreteeven} with $V_1\XX$ the (full) preimage of
$V_p\YY$ in $\XX$, with $V_2\XX$ the preimage of $V_q\YY$ in $\XX$ and
with $c=0$, so that $\delta_c=\ln\sqrt{pq}$ and $m_c=m_{\rm BM}$.  Let
us normalise the Patterson density of $\Ga$ as in Proposition
\ref{prop:computBM} (2), so that
$$
\|\sigma^+_{\D^-,\,1}\|=\|\mu_{\,\wt v}\|=\;\frac{p+1}{\sqrt p}\;.
$$ 
By the proof of Equation \eqref{eq:massskinlinebiregular}, the mass
for the skinning measure of the part of the inner unit normal bundle
of $\YY^+$ with basepoint in $V_p\YY$ is $\frac{p-1}{\sqrt p}
\frac{L}{2}$ and its complement has mass $\frac{L(q-1)}{2\sqrt q}$.
Recall also that, by Proposition \ref{prop:computBM} (2) and Remark
\ref{rem:TVoltorsionfree}, considering the graph $\YY$ as a graph of
groups with trivial groups,
$$
\|m_{\rm BM}\|=\TVol(\YY)=|E\YY|=2(p+1)|V_p\YY|=2(q+1)|V_q\YY|\,.
$$ 

The claim about the common perpendiculars of even length at most
$2N$ follows from Theorem \ref{theo:countingdiscreteeven} with $i=j=1$,
since
$$
\frac{2\;e^{2\,\delta_c}\;\|\sigma^+_{\D^-,\,i}\|\;\|\sigma^-_{\D^+,\,j}\|}
{(e^{2\,\delta_c}-1)\;\|m_{c}\|}
=\frac{2\;pq\;\frac{p+1}{\sqrt p}\;\frac{L(p-1)}{2\sqrt p}}
{(pq-1)\;\;2\;(p+1)\;|V_p\YY|}=
\frac{L\;q\;(p-1)}{2\;(pq-1)\;|V_p\YY|}\;.
$$

The claim about the common perpendiculars of odd length at most $2N-1$
follows similarly from Theorem \ref{theo:countingdiscreteeven} with
$i=1$ and $j=2$.  
\cqfd 

\medskip\noindent (4) Let $\YY$ be a finite biregular graph with
vertices of degrees $p+1$ and $q+1$, where $p,q\geq 2$, and let
$V\YY=V_p\YY\sqcup V_q\YY$ be the corresponding partition.  If $\YY^-$
and $\YY^+$ are cycles of length $L^-\geq 2$ and $L^+\geq 2$
respectively, then as $N\to+\infty$, the number of common
perpendiculars of even length at most $2N$ from $\YY^-$ to $\YY^+$ is
equal to
\begin{equation}\label{eq:bicyclette}
\frac{(p+q)\;L^-\;L^+}
{2\,(pq-1)\;|E\YY|}\;(pq)^{N+1} + \bigO(r^N)\;
\end{equation}
for some $r<\sqrt{pq}$. 

\medskip
\dem As in the above proof of Example (3), let $\XX\ra \YY$ be a
universal cover of $\YY$, with covering group $\Ga$ and let $\DD^\pm$
be a geodesic line in $\XX$ mapping to $\YY^\pm$. Let $V_1\XX$ be the
preimage of $V_p\YY$ and $V_2\XX$ be one of $V_q\YY$. We normalise the
Patterson density $(\mu_x)_{x\in V \XX}$ of $\Ga$ so that
$\|\mu_{x}\|= \frac{\deg_\XX(x)} {\sqrt{\deg_\XX(x)-1}}$. By
Proposition \ref{prop:computskin} (3) with $k=1$ and trivial vertex
stabilisers, and since a (simple) cycle of length $\len$ in a biregular
graph of different degrees $p+1$ and $q+1$ has exactly
$\frac{\len}{2}$ vertices of degree either $p+1$ or $q+1$, we have
$$
\|\sigma^\pm_{\D^\mp,1}\|=\sum_{\Ga x\in V_p\YY^\mp}\;
\frac{\|\mu_x\|\;(\deg_\XX(x)-k)}{\deg_\XX(x)}=
\sum_{y\in V_p\YY^\mp}\;\sqrt{p}
=\frac{L^\mp\sqrt{p}}{2}\;.
$$ Similarly, $ \|\sigma^\pm_{\D^\mp,2}\|==\frac{L^\mp\sqrt{q}}{2}
$. The result without the error term then follows from Theorem
\ref{theo:countingdiscreteeven}, using Proposition \ref{prop:computBM}
(2) and Remark \ref{rem:TVoltorsionfree}, since the number we are
looking for is $\N_{\DD^-,\DD^+,1,1}(2N)+\N_{\DD^-,\DD^+,2,2}(2N)$.
We refer to Section \ref{subsect:errormetricgraphgroup} (see Remark
(ii) following the proof of Theorem
\ref{theo:equidsimplicialperperrorterm}) for the error term.  \cqfd
\eexems
%

\section{Counting for simplicial graphs of groups}
\label{subsec:graphgroupsequidcount}

In this Section, we give an intrinsic translation ``a la Bass-Serre''
of the counting result in Theorem \ref{theo:countingdiscrete} using
graphs of groups (see \cite{Serre83} and Section \ref{subsec:trees}
for background information).  

Let $(\YY,G_*)$ be a locally finite, connected graph of finite groups,
and let $(\YY^\pm,G_*^\pm)$ be connected subgraphs of
subgroups.\footnote{See Section \ref{subsec:trees} for definitions and
  background.} Let $c:E\YY\ra \RR$ be a system of conductances on
$\YY$.

Let $\XX$ be the Bass-Serre tree of the graph of groups $(\YY,G_*)$
(with geometric realisation $X=|\XX|_1$) and $\Ga$ its fundamental
group (for an indifferent choice of basepoint). Assume that $\Ga$ is
nonelementary. We denote by $\G(\YY,G_*)= \Ga\bs\G\XX$ and
$\big(\flow{t}: \G(\YY,G_*)\ra \G(\YY,G_*) \big)_{t\in\ZZ}$ the
quotient of the (discrete time) geodesic flow on $\G\XX$, by $\wt
c:\XX\ra\RR$ the ($\Ga$-invariant) lift of $c$, with $\delta_c$ its
critical exponent and $\wt F_c: T^1X \ra \RR$ its associated
potential, by $m_c$ the Gibbs measure on $\G(\YY,G_*)$ associated with
a choice of Patterson densities $(\mu^\pm_x)_{x\in X}$ for the pairs
$(\Ga,F_c^\pm)$, by $\DD^\pm$ two subtrees in $\XX$ such that the
quotient graphs of groups $\Ga_{\DD^\pm} \dbs \DD^\pm$ identify with
$(\YY^\pm,G_*^\pm)$ (see below for precisions), and by
$\sigma^\pm_{(\YY^\mp,G^\mp_*)}$ the associated skinning measures.

\medskip
The {\it fundamental groupoid}\index{fundamental groupoid}
$\gls{fundamentalgroupoid}$ of $(\YY,G_*)$\footnote{denoted by
  $F(\YY,G_*)$ in \cite[\S 5.1]{Serre83}, called the {\it path group}
  in \cite[1.5]{Bass93}, see also \cite{Higgins76}} is the quotient of
the free product of the groups $G_v$ for $v\in V\YY$ and of the free
group on $E\YY$ by the normal subgroup generated by the elements $e\;
\overline{e}$ and $e\,\rho_e(g)\;\overline{e}\; \rho_{\,\overline{e}}
(g)^{-1}$ for all $e\in E\YY$ and $g\in G_e$. We identify each $G_x$
for $x\in V\YY$ with its image in $\pi(\YY,G_*)$.

Let $n\in\NN-\{0\}$. A {\em (locally) geodesic
  path}\index{geodesic!path} of length $n$ in the graph of groups
$(\YY,G_*)$ is the image $\alpha$ in $\pi(\YY,G_*)$ of a word, called
{\it reduced}\index{reduced} in \cite[1.7]{Bass93},
$$
h_0\,e_1\,h_1\, e_2\,\dots\, h_{n-1}\,e_n\,h_n
$$
with
\begin{enumerate}
\item[$\bullet$] $e_i\in E\YY$ and $t(e_i)= o(e_{i+1})$ for $1\leq
  i\leq n-1$ (so that $(e_1,\dots, e_n)$ is an edge path in the graph
  $\YY$);
\item[$\bullet$] $h_0\in G_{o(e_1)}$ and $h_i \in G_{t(e_i)}$ for
  $1\leq i\leq n$;
\item[$\bullet$] if $e_{i+1}=\overline{e_i}$ then $h_i$ does not
  belong to $\rho_{e_i}(G_{e_i})$, for $1\leq i\leq n-1$.
\end{enumerate}
Its {\em origin}\index{common perpendicular!origin}\index{origin} is
$o(\alpha)=o(e_1)$ and its {\it endpoint}%
\index{common perpendicular!endpoint}\index{endpoint} is $t(\alpha)=
t(e_n)$. They do not depend on the chosen words with image $\alpha$ in
$\pi(\YY,G_*)$.

A {\em common perpendicular}\index{common perpendicular} of length $n$
from $(\YY^-,G_*^-)$ to $(\YY^+,G_*^+)$ in the graph of groups
$(\YY,G_*)$ is the double coset
$$
[\alpha]=G^-_{o(\alpha)}\,\alpha\; G^+_{t(\alpha)}
$$ 
of a geodesic path $\alpha$ of length $n$ as above, such that:

\smallskip
$\bullet$~ $\alpha$ {\em starts transversally from $(\YY^-, G_*^-)$},%
\index{common perpendicular!starting transversally}%
\index{starting transversally} that is, its origin $o(\alpha)=o(e_1)$
belongs to $V\YY^-$ and $h_0\notin G^-_{o(e_1)}\,
\rho_{\,\overline{e_1}}(G_{e_1})$ if $e_1\in E\YY^-$,

\smallskip
$\bullet$~ $\alpha$ {\em ends transversally in
  $(\YY^+,G_*^+)$},\index{common perpendicular!ending
  transversally}\index{ending transversally} that is, its endpoint
$t(\alpha)=t(e_n)$ belongs to $\YY^+$ and $h_n\notin
\rho_{e_n}(G_{e_n})\,G^+_{t(e_n)}$ if $e_n\in E\YY^+$.  

\smallskip
Note that these two notions do not depend on the representative of the
double coset $G^-_{o(\alpha)}\,\alpha \;G^+_{t(\alpha)}$, and we also
say that the double coset $[\alpha]$ {\em starts transversally from
  $(\YY^-, G_*^-)$} or {\em ends transversally in $(\YY^-,G_*^-)$}.

\medskip
We denote by $\Perp((\YY^\pm,G_*^\pm),n)$ the set of common
perpendiculars in $(\YY,G_*)$ of length at most $n$ from
$(\YY^-,G_*^-)$ to $(\YY^+,G_*^+)$. We denote by
$$
c(\alpha)= \sum_{i=1}^n c(e_i)
$$ 
the {\em conductance}\index{conductance} of a geodesic path
$\alpha$ as above, which depends only on the double class $[\alpha]$.
We define the {\it multiplicity}\index{multiplicity} $m_\alpha$ of a
geodesic path $\alpha$ as above by
$$
m_\alpha=\frac{1}{\card(G_{o(\alpha)}^-\cap\alpha \,G_{t(\alpha)}^+ \,\alpha^{-1})}\;.
$$ 
It depends only on the double class $[\alpha]$ of $\alpha$. We
define the {\it counting function}\index{counting function} of the
common perpendiculars in $(\YY,G_*)$ of length at most $n$ from
$(\YY^-,G_*^-)$ to $(\YY^+,G_*^+)$ (counted with multiplicities and
with weights given by the system of conductances $c$) as
$$
\N_{(\YY^-,G_*^-),\,(\YY^+,G_*^+)}(n)=
\sum_{[\alpha]\in\Perp((\YY^\pm,G_*^\pm),n)} m_{\alpha}\;e^{c(\alpha)}\;.
$$

\btheo\label{theo:countinggraphgroup} Let $(\YY,G_*)$,
$(\YY^\pm,G_*^\pm)$ and $c$ be as in the beginning of this
Section. Assume that the critical exponent $\delta_c$ of $c$ is finite
and positive and that the Gibbs measure $m_{c}$ on $\G(\YY,G_*)$ is
finite and mixing for the discrete time geodesic flow. Assume that the
skinning measures $\sigma^\pm_{(\YY^\mp,G^\mp_*)}$ are finite and
nonzero.  Then as $n\in\NN$ tends to $\infty$
$$
\N_{(\YY^-,G_*^-),\,(\YY^+,G_*^+)}(n)\sim 
\frac{e^{\delta_c}\;\|\sigma^+_{(\YY^-,G^-_*)}\|\;\|\sigma^-_{(\YY^+,G^+_*)}\|}
{(e^{\delta_c}-1)\;\|m_{c}\|}\;\;e^{\delta_c\, n}\;.
$$
\etheo

\dem 
Let $\XX$ be the Bass-Serre tree of $(\YY,G_*)$ and $\Ga$ its
fundamental group (for an indifferent choice of basepoint). As seen in
Section \ref{subsec:trees}, the Bass-Serre trees $\DD^\pm$ of
$(\YY^\pm,G_*^\pm)$, with fundamental groups $\Ga^\pm$, identify with
simplicial subtrees $\DD^\pm$ of $\XX$, such that $\Ga^\pm$ are the
stabilisers $\Ga_{\DD^\pm}$ of $\DD^\pm$ in $\Ga$, and that the
maps $(\Ga_{\DD^\pm} \bs \DD^\pm) \ra (\Ga\bs\XX)$ induced by the
inclusion maps $\DD^\pm\ra\XX$ by taking quotient, are injective.

As in Definition \ref{exem:quotgraphgroup}, for all $z\in V\YY\cup
E\YY$ and $e\in E\YY$, we fix a lift $\wt z\in V\XX\cup E\XX$ of $z$
and $g_e\in\Ga$, such that $\overline{\wt e}= \wt{\overline{e}}$, $g_e
\,\wt{t(e)}= t(\wt e)$, $G_z=\Ga_{\wt z}$, and the monomorphism
$\rho_e:G_e\ra G_{t(e)}$ is $\ga\mapsto g_e^{-1}\ga g_e$. We assume,
as we may, that $\wt z\in V\DD^\pm\cup E\DD^\pm $ if $z\in
V\YY^\pm\cup E\YY^\pm$. We assume, as we may using Equation
\eqref{eq:simplesubtree}, that if $e\in E\YY^\pm$, then $g_e\in
\Ga_{\DD^\pm}$. We denote by $p:\XX\ra\YY=\Ga\bs\XX$ the canonical
projection.

\begin{center}
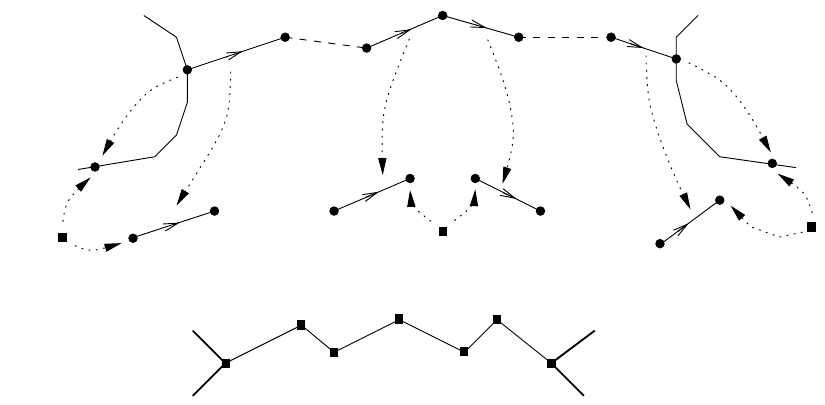
\end{center}

For all $\ga,\ga'\in\Ga$ such that $\ga\DD^-$ and $\ga' \DD^+$ are
disjoint, the common perpendicular $\alpha_{\ga\DD^-,\,\ga' \DD^+}$
from $\ga\DD^-$ to $\ga' \DD^+$ is an edge path $(f_1,f_2,\dots, f_k)$
with $o(f_1)\in \ga\DD^-$ and $t(f_k) \in \ga'\DD^+$. Note that
$\ga^{-1} \,o(f_1)$ and $\reallywidetilde{p(o(f_1))}$ are two vertices
of $\DD^-$ in the same $\Ga$-orbit, and that ${\ga'}^{-1} \,t(f_k)$
and $\reallywidetilde{p(t(f_k))}$ are two vertices of $\DD^+$ in the
same $\Ga$-orbit. Hence by Equation \eqref{eq:simplesubtree}, we may
choose $\ga_0\in \Ga_{\DD^-}$ such that $\ga_0\ga^{-1} \,o(f_1) =
\reallywidetilde{p(o(f_1))}$ and $\ga_{k+1}\in \Ga_{\DD^-}$ such that
$\ga_{k+1}{\ga'}^{-1} \,t(f_k)=\reallywidetilde{p(t(f_k))}$. For
$1\leq i\leq k$, choose $\ga_i\in \Ga$ such that $\ga_if_i=
\widetilde{p(f_i)}$. We define (see the above picture)

\medskip
$\bullet$~ $e_i=p(f_i)$ for $1\leq i\leq k$,

\smallskip 
$\bullet$~ $h_i=g_{e_i}^{-1} \ga_i{\ga_{i+1}}^{\!-1}g_{\,\overline{e_{i+1}}}$, 
which belongs to $\Ga_{\widetilde{t(e_i)}} =G_{t(e_i)}$ for
$1\leq i\leq k-1$,

\smallskip 
$\bullet$~ $h_0= \ga_0\ga^{-1}\ga_1^{-1}g_{\,\overline{e_1}} =
\ga^{-1} (\ga\ga_0\ga^{-1})\ga_1^{-1}g_{\,\overline{e_1}} $, which
belongs to $\Ga_{\widetilde{o(e_1)}} =G_{o(e_1)}$,

\smallskip $\bullet$~ $h_k= g_{e_k}^{-1} \ga_k\ga'{\ga_{k+1}}^{\!-1} =
g_{e_k}^{-1} \ga_k(\ga'\ga_{k+1}{\ga'}^{-1})^{-1}\ga'$, which belongs
to $\Ga_{\widetilde{t(e_k)}} = G_{t(e_k)}$.

\blemm \label{lem:bassserrerie}\mbox{}\begin{enumerate}
\item[(1)] The word $h_0e_1h_1\dots h_{k-1}e_kh_k$ is reduced. Its
  image $\alpha$ in the fundamental groupoid $\pi(\YY,G_*)$ does not
  depend on the choices of $\ga_1,\dots,\ga_k$, and starts
  transversally from $(\YY^-,G_*^-)$ and ends transversally in
  $(\YY^+,G_*^+)$. The double class $[\alpha]$ of $\alpha$ is
  independent of the choices of $\ga_0$ and $\ga_{k+1}$.
\item[(2)] The map $\wt\Theta$ from the set of common perpendiculars
  in $\XX$ between disjoint images of $\DD^-$ and $\DD^+$ under
  elements of $\Ga$, into the set of common perpendiculars in
  $(\YY,G_*)$ from $(\YY^-,G_*^-)$ to $(\YY^+,G_*^+)$, sending
  $\alpha_{\ga \DD^-,\,\ga' \DD^+}$ to $[\alpha]$, is invariant under
  the action of $\Ga$ at the source, and preserves the lengths and the
  multiplicities.
\item[(3)] The map $\Theta$ induced by $\wt\Theta$ from the set of
  $\Ga$-orbits of common perpendiculars in $\XX$ between disjoint
  images of $\DD^-$ and $\DD^+$ under elements of $\Ga$ into the set
  of common perpendiculars in $(\YY,G_*)$ from $(\YY^-,G_*^-)$ to
  $(\YY^+,G_*^+)$ is a bijection, preserving the lengths and the
  multiplicities.
\end{enumerate}
\elemm

\dem (1) If $e_{i+1}=\overline{e_i}$, then by the definition of $h_i$,
we have
\begin{align*}
h_i\in\rho_{e_i}(G_{e_i})=g_{e_i}^{-1}\Ga_{\wt{e_i}}g_{e_i} & 
\Longleftrightarrow g_{e_i}\,h_i\,g_{e_i}^{-1}\;\wt{e_i} = \wt{e_i} \\ &
\Longleftrightarrow g_{e_i}\,g_{e_i}^{-1} \ga_i{\ga_{i+1}}^{\!-1}
g_{\,\overline{e_{i+1}}}\,g_{e_i}^{-1}\;\wt{e_i} = \wt{e_i}\\ &
\Longleftrightarrow \ga_i{\ga_{i+1}}^{\!-1}\;\wt{\overline{e_{i+1}}} 
= \wt{e_i} \\ &
\Longleftrightarrow \overline{{\ga_{i+1}}^{\!-1}\;\wt{e_{i+1}}} = 
\ga_i^{-1}\wt{e_i} \;\Longleftrightarrow\; \overline{f_{i+1}} = f_i\;.
\end{align*}
Hence the word $h_0e_1h_1\dots h_{k-1}e_kh_k$ is reduced.

The element $\ga_i$ for $i\in \{1, \dots, k\}$ is uniquely determined
up to multiplication on the left by an element of
$\Ga_{\wt{e_i}}=G_{e_i}$. If we fix\footnote{We leave to the reader
  the verification that the changes induced by various $i$'s do not
  overlap.} $i\in \{1, \dots, k\}$ and if we replace $\ga_i$ by
$\ga'_i=\alpha\,\ga_i$ for some $\alpha\in G_{e_i}$, then only the
elements $h_{i-1}$ and $h_i$ change, replaced by elements that we
denote by $h'_{i-1}$ and $h'_i$ respectively. We have (if $2\leq i\leq
k-1$, but otherwise the argument is similar by the definitions of
$h_0$ and $h_k$)
\begin{align*}
h'_{i-1}\,e_i\,h'_i&=
g_{e_{i-1}}^{\!-1} \ga_{i-1}\,{\ga_{i}}^{-1}\alpha^{-1}g_{\,\overline{e_{i}}}
\;e_i\;g_{e_i}^{-1} \alpha\,\ga_i\,{\ga_{i+1}}^{\!-1}g_{\,\overline{e_{i+1}}}\\ &
=g_{e_{i-1}}^{\!-1} \ga_{i-1}\,{\ga_{i}}^{-1}g_{\,\overline{e_{i}}}
\,\rho_{\,\overline{e_{i}}}(\alpha)^{-1}\;e_i\;\rho_{e_i}(\alpha) \,
g_{e_i}^{-1} \ga_i\,{\ga_{i+1}}^{\!-1}g_{\,\overline{e_{i+1}}}\;.
\end{align*}
Since $\rho_{\,\overline{e_{i}}}(\alpha)^{-1}\;e_i\;\rho_{e_i}(\alpha)$
is equal to $\overline{e_i}^{\;-1}=e_i$ in the fundamental groupoid, the
words $h'_{i-1}\,e_i\,h'_i$ and $h_{i-1}\,e_i\,h_i$ have the same
image in $\pi(\YY,G_*)$. Therefore $\alpha$ is does not depend on the
choices of $\ga_1,\dots,\ga_k$.

\medskip
We have $o(\alpha)=o(e_1)\in V\YY^-$ and $t(\alpha)=t(e_k)\in V\YY^+$,
hence $\alpha$ starts from $\YY^-$ and ends in $\YY^+$.

Assume that $e_1\in E\YY^-$. Let us prove that $h_0\in G^-_{o(e_1)}\,
\rho_{\,\overline{e_1}}(G_{e_1})$ if and only if $f_1\in \ga\,E\DD^-$.
\begin{center}
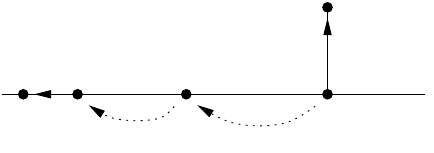
\end{center}

By the definition of $\rho_{\,\overline{e_1}}$, we have $h_0\in
G^-_{o(e_1)}\, \rho_{\,\overline{e_1}} (G_{e_1})$ if and only if there
exists $\alpha\in \Ga_{\,\wt{o(e_1)}} \cap \Ga_{\DD^-}$ such that
$\alpha^{-1}\,h_0\in g_{\,\overline{e_1}}^{\;-1} \,\Ga_{\,\wt{e_1}}\;
g_{\,\overline{e_1}}$. By the definition of $h_0$ and since $\ga_1$
maps $f_1$ to $\wt{e_1}$, we have
\begin{align*}
\alpha^{-1}\,h_0\in 
g_{\,\overline{e_1}}^{\;-1} \,\Ga_{\,\wt{e_1}} \;g_{\,\overline{e_1}}& 
\Longleftrightarrow
g_{\,\overline{e_1}}\,\alpha^{-1}\,
\big(\ga_0\ga^{-1}\ga_1^{-1}g_{\,\overline{e_1}}\big)\,
g_{\,\overline{e_1}}^{\;-1} \;\wt{e_1}\;=\;\wt{e_1}\\ & 
\Longleftrightarrow
f_1=\ga\,\ga_0^{-1}\,\alpha\,g_{\,\overline{e_1}}^{\;-1} \;\wt{e_1}\;.
\end{align*}
Since $\wt {e_1}\in E\DD^-$ and $\ga_0$, $\alpha$,
$g_{\,\overline{e_1}}$ all belong to $\Ga_{\DD^-}$, this last
condition implies that $f_1\in \ga\,E\DD^-$. Conversely (for future
use), if $f_1\in\ga\,E\DD^-$, then (see the above picture)
$\ga_0\,\ga^{-1} f_1$ is an edge of $\DD^-$ with origin
$\wt{o(e_1)}$, in the same $\Ga$-orbit than the edge
$g_{\,\overline{e_1}}^{\;-1} \;\wt{e_1}$ of $\DD^-$, which also has
origin $\wt{o(e_1)}$. By Equation \eqref{eq:simplesubtree}, this
implies that there exists $\alpha\in \Ga_{\,\wt{o(e_1)}} \cap
\Ga_{\DD^-}$ such that $f_1=\ga\,\ga_0^{-1}\,\alpha\,
g_{\,\overline{e_1}}^{\;-1} \;\wt{e_1}$. By the above equivalences, we
hence have that $h_0\in G^-_{o(e_1)}\, \rho_{\,\overline{e_1}}
(G_{e_1})$.

Similarly, one proves that if $e_k\in E\YY^+$, then $h_k\in
\rho_{e_k}(G_{e_k})\,G^+_{t(e_k)}$ if and only if $f_k\in \ga'E\DD^+$.
Since $(f_1,\dots, f_n)$ is the common perpendicular edge path from
$\ga\,\DD^-$ to $\ga'\DD^+$, this proves that $\alpha$ starts
transversally from $\YY^-$ and ends transversally in $\YY^-$.

\medskip
Note that the element $\ga_0\in \Ga_{\DD^-}$ is uniquely defined up to
multiplication on the left by an element of $\Ga_{\wt{o(e_1)}}\,\cap\,
\Ga_{\DD^-} =G^-_{o(e_1)}$, and appears only as the first letter in
the expression of $h_0$. Note that the element $\ga_{k+1}\in
\Ga_{\DD^+}$ is uniquely defined up to multiplication on the left by
an element of $\Ga_{\wt{t(e_k)}} \cap\,\Ga_{\DD^+}=G^+_{t(e_k)}$,
hence $\ga_{k+1}^{\!-1}$ is uniquely defined up to multiplication on
the right by an element of $G^+_{t(e_k)}$, and appears only as the
last letter in the expression of $h_k$. Therefore $\alpha$ is uniquely
defined in the fundamental groupoid $\pi(\YY,G_*)$ up to
multiplication on the left by an element of $G^-_{o(e_1)}$ and
multiplication on the right by an element of $G^+_{t(e_k)}$, that is,
the double class $[\alpha]\subset \pi(\YY,G_*)$ is uniquely
defined.

\medskip 
(2) Let $\beta$ be an element in $\Ga$ and let $x=\alpha_{\ga \DD^-,
  \,\ga' \DD^+}$ be a common perpendicular in $\XX$ between disjoint
images of $\DD^-$ and $\DD^+$ under elements of $\Ga$. Let us prove
that $\wt\Theta(\beta\,x) =\wt\Theta(x)$.

Since $\beta\,x=\alpha_{\beta\,\ga \DD^-,\,\beta\,\ga' \DD^+}$, in the
construction of $\wt\Theta(\beta\,x)$, we may take, instead of the
elements $\ga_0,\ga_1,\dots,\ga_k,\ga_{k+1}$ used to construct
$\wt\Theta(x)$, the elements
$$
\ga_0^\sharp=\ga_0,\;\;\ga^\sharp_1=\ga_1\,\beta^{-1},\dots,\;\;
\ga^\sharp_k=\ga_k\,\beta^{-1},\;\;\ga^\sharp_{k+1}=\ga_{k+1}.
$$ 
And instead of $\ga$ and $\ga'$, we now may use
$\ga^\sharp=\beta\,\ga$ and ${\ga'}^\sharp=\beta\,\ga'$.

The only terms involving $\ga,\ga',\ga_1,\dots,\ga_k$ in the
construction of $\wt\Theta(x)$ come under the form $\ga^{-1}
\ga_1^{-1}$ in $h_0$, $\ga_i\,{\ga_{i+1}}^{\!-1}$ in $h_i$ for $1\leq
i\leq k-1$, and $\ga_k\ga'$ in $h_k$. Since $(\ga^\sharp)^{-1}
(\ga_1^\sharp)^{-1} =\ga^{-1} \ga_1^{-1}$, $(\ga_i^\sharp)
(\ga_{i+1}^\sharp)^{-1}=\ga_i\,{\ga_{i+1}}^{\!-1}$ for $1\leq i\leq
k-1$, and $(\ga_k^\sharp) ({\ga'}^\sharp)=\ga_k\ga'$, this proves that
$\wt\Theta(\beta\,x) =\wt\Theta(x)$, as wanted.

It is immediate that if the length of $\alpha_{\ga \DD^-, \,\ga'
  \DD^+}$ is $k$, then the length of $[\alpha]$ is $k$.

Let us prove that the multiplicity, given in Equation
\eqref{eq:multiplicitytree},
$$
m_{\ga,\,\ga'}=
\frac 1{\card(\ga\Ga_{\DD^-}\ga^{-1}\cap\ga'\Ga_{\DD^+}{\ga'}^{-1})}
$$ 
of the common perpendicular $\alpha_{\ga \DD^-, \,\ga' \DD^+}$ in
$\XX$ between $\ga\,\DD^-$ and $\ga'\,\DD^+$ is equal to the
multiplicity
$$
m_\alpha=
\frac{1}{\card(G_{o(\alpha)}^-\cap \alpha \,G_{t(\alpha)}^+ \,\alpha^{-1})}
$$ 
of the common perpendicular $\alpha$ in $(\YY,G_*)$ from
$(\YY^-,G_*^-)$ to $(\YY^+,G_*^+)$.

Since the multiplicity $m_{\ga,\,\ga'}$ is invariant under the
diagonal action by left translations of $\ga_0^{-1}\ga^{-1}\in \Ga$ on
$(\ga,\ga')$, we may assume that $\ga=\ga_0=\id$. Since the
multiplicity $m_{\ga,\,\ga'}$ is invariant under right translation by
${\ga_{k+1}}^{\!-1}$, which stabilises $\DD^+$, on the element $\ga'$,
we may assume that $\ga_{k+1}=\id$. In particular, we have
$$
o(f_1)=\wt{o(e_1)}\;\;\;\;{\rm and}\;\;\;\;t(f_k)=\ga'\;\wt{t(e_k)}\;.
$$ 

We use the basepoint $x_0=o(e_1)$ in the construction of the
fundamental group and the Bass-Serre tree of $(\YY,G_*)$, so that
(see in particular \cite[Eq.~(1.3)]{Bass93})
$$
V\XX=\coprod_{\beta\in \pi(\YY,\,G_*)\;:\; o(\beta)=x_0} \beta \;G_{t(\beta)}
$$
and
$$
\Ga=\pi_1(\YY,G_*)=\{\beta\in\pi(\YY,G_*)\;:\;o(\beta)=t(\beta)=x_0\}\;.
$$

Since an element in $\Ga$ which preserves $\DD^-$ and $\ga'\,\DD^+$
fixes pointwise its (unique) common perpendicular in $\XX$, we have
\begin{align*}
\Ga_{\DD^-}\,\cap\,\ga'\Ga_{\DD^+}{\ga'}^{-1} & =
\Ga_{\DD^-}\,\cap\,\Ga_{\ga'\,\DD^+}=
(\Ga_{o(f_1)}\cap\Ga_{\DD^-}) \,\cap\,(\Ga_{t(f_k)}\cap\Ga_{\ga'\,\DD^+})
\\ & = (\Ga_{\wt{o(e_1)}}\cap\Ga_{\DD^-}) \,\cap\,
(\Ga_{\ga'\,\wt{t(e_k)}}\cap\Ga_{\ga'\,\DD^+})\;.
\end{align*}
Note that $\Ga_{\wt{o(e_1)}}\cap\Ga_{\DD^-}=G^-_{o(e_1)}$. By the
construction of the edges in the Bass-Serre tree of a graph of groups
(see \cite[page 11]{Bass93}), the vertex $\alpha \;G_{t(e_k)}$ is
exactly the vertex $t(f_k)=\ga'\,\wt{t(e_k)}$. By
\cite[Eq.~(1.4)]{Bass93}, we hence have
$$ 
\alpha \,G_{t(e_k)}\,\alpha^{-1}=\operatorname{Stab}_{\pi_1(\YY,G_*)}
(\alpha \,G_{t(e_k)})=\Ga_{\ga'\,\wt{t(e_k)}}\;.
$$
Therefore $m_{\ga,\,\ga'}=m_\alpha$.

\medskip
(3) Let $[\alpha]=G_{o(\alpha)}\,\alpha\,G_{t(\alpha)}$ be a common
perpendicular in $(\YY,G_*)$ from $(\YY^-,G_*^-)$ to $(\YY^+,G_*^+)$,
with representative $\alpha\in\pi(\YY,G_*)$, and let $h_0 \,e_1\,
h_1\, \dots \,e_k \,h_k$ be a reduced word whose image in
$\pi(\YY,G_*)$ is $\alpha$.

We define

\medskip
$\bullet$~ $\ga_1=g_{\,\overline{e_1}}\;h_0^{-1}$,

\smallskip 
$\bullet$~ $f_1=\ga_1^{-1}\,\wt{e_1}\,$,

\smallskip 
$\bullet$~ assuming that $\ga_i$ and $f_i$ for some $1\leq i\leq k-1$
are constructed, let
$$
\ga_{i+1}=g_{\,\overline{e_{i+1}}}\;h_i^{-1}\,g_{e_i}^{-1}\,\ga_i
\;\;\;\;{\rm and}\;\;\;\;
f_{i+1}={\ga_{i+1}}^{\!-1}\,\wt{e_{i+1}}\;,
$$

\smallskip 
$\bullet$~ with $\ga_k$ and $f_k$ constructed by induction, finally
let $\ga'=\ga_k^{-1}\,g_{e_k}\,h_k$.

\medskip
It is easy to check, using the equivalences in the proof of Lemma
\ref{lem:bassserrerie} (1) with $\ga=\ga_0=\ga_{k+1}=\id$, that the
sequence $(f_1,\dots, f_k)$ is the edge path of a common perpendicular
in $\XX$ from $\DD^-$ to $\ga'\DD^+$ with origin $\wt{o(e_1)}$ and
endpoint $\ga'\,\wt{t(e_k)}$.

If $h_0$ is replaced by $\alpha\, h_0$ with $\alpha\in G^-_{o(e_1)}$,
then by induction, $f_1,f_2\dots, f_k$ are replaced by $\alpha f_1,
\alpha f_2, \dots, \alpha f_k$ and $\ga'$ is replaced by $\alpha\ga'$.
Note that $(\alpha f_1,\alpha f_2, \dots, \alpha f_k)$ is then the
common perpendicular edge path from $\DD^-=\alpha\DD^-$ to $\alpha
\ga'\DD^+$. If $h_k$ is replaced by $h_k\,\alpha $ with $\alpha \in
G^+_{t(e_k)}$, then $f_1,f_2\dots, f_k$ are unchanged, and $\ga'$ is
replaced by $\ga'\,\alpha $. Note that $\ga'\,\alpha \,\DD^+=\ga'\,
\DD^+$.

Hence the map which associates to $[\alpha]$ the $\Ga$-orbit of the
common perpendicular in $\XX$ from $\DD^-$ to $\ga'\DD^+$ with edge
path $(f_1,\dots, f_k)$ is well defined. It is easy to see by
construction that this map is the inverse of $\Theta$.  
\cqfd

\medskip 
Theorem \ref{theo:countinggraphgroup} now follows from Theorem
\ref{theo:countingdiscrete}. 
\cqfd

\section{Error terms for equidistribution and counting for metric 
and simplicial graphs of groups}
\label{subsect:errormetricgraphgroup}

In this Section, we give error terms to the equidistribution and
counting results of Section \ref{subsect:equicountmetricgraphgroup},
given by Theorem \ref{theo:counting} for metric trees (and their
continuous time geodesic flows) and by Theorem
\ref{theo:countingdiscrete} for simplicial trees (and their discrete
time geodesic flows), under appropriate assumptions on bounded
geometry and the rate of mixing.

\medskip
Let $(\XX,\lambda)$, $X$, $\Ga$, $\wt c$, $c$, $\wt F_c$, $F_c$,
$\delta_c$, $\DD^\pm$, $D^\pm$, $\D^\pm$, $\lambda_{\ga,\ga'}$, $\alpha_{\ga,\ga'}$,
$\alpha^\pm_{\ga,\ga'}$, $m_{\ga,\ga'}$ be as in Section
\ref{subsect:equicountmetricgraphgroup}.  We first consider the
simplicial case (when $\lambda=1$), for the discrete time geodesic
flow.

\btheo\label{theo:equidsimplicialperperrorterm} Let $\XX$ be a locally
finite simplicial tree without terminal vertices, let $\Ga$ be a
nonelementary discrete subgroup of $\Aut(\XX)$, let $\wt c$ be a
system of conductances on $\XX$ for $\Ga$ and let $\DD^\pm$ be
nonempty proper simplicial subtrees of $\XX$.  Assume that the critical
exponent $\delta_c$ is finite and positive, that the Gibbs measure
$m_{c}$ (for the discrete time geodesic flow) is finite and that the
skinning measures $\sigma^\pm_{\D^\mp}$ are finite and nonzero. Assume
furthermore that
\begin{enumerate}
\item at least one of the following holds :

$\bullet$~ $\Ga_{\DD^\pm}\bs \partial D^\pm$ is compact
 
$\bullet$~ $\C\Lambda\Ga$ is uniform and $\Ga$ is a lattice of
  $\C\Lambda\Ga$,
\item there exists $\beta\in\;]0,1]$ such that the discrete time geodesic
      flow on $(\Ga\bs\G \XX,m_c)$ is exponentially mixing for the
      $\beta$-H\"older regularity.
\end{enumerate}
Then there exists $\kappa'>0$ such that for all $\psi^\pm\in \C_{\rm
  c}^\beta (\Ga\bs\gengeod\XX)$, we have, as $n\ra+\infty$,
\begin{align*}
 \frac{e^{\delta_c}-1}{e^{\delta_c}}&\;\|m_{c}\|\;e^{-\delta_c\, n}
\sum_{\substack{[\ga]\,\in \Ga_{\DD^-}\bs\Ga/\Ga_{\DD^+}\\
0<d(\DD^-,\ga \,\DD^+)\leq n}} 
m_{e,\,\ga}\;e^{\wt c\,(\alpha_{e,\,\ga})}\;
\psi^-(\Ga \alpha^-_{e,\ga})\;\psi^+(\Ga\alpha^+_{\ga^{-1},e}) \\= \;&
\int \psi^-\;d\sigma^+_{\D^-}\;\int \psi^+\;d\sigma^-_{\D^+} 
+ \bigO\big(e^{-\kappa' \,n}\;\|\psi^-\|_\beta\;\|\psi^+\|_\beta\big)
\end{align*}
and if $\Ga_{\DD^\pm}\bs \partial D^\pm$ is compact, then
$$
\N_{\DD^-,\,\DD^+}(n)=
\frac{e^{\delta_c}\;\|\sigma^+_{\D^-}\|\;\|\sigma^-_{\D^+}\|}
{(e^{\delta_c}-1)\;\|m_{c}\|}\;\;e^{\delta_c\, n} + 
\bigO\big(e^{(\delta_c-\kappa')n}\big)\,.
$$
\etheo

\medskip
\dem 
We follow the scheme of proof of Theorem \ref{theo:expratecount},
replacing aspects of Riemannian manifolds by aspects of simplicial
trees as in the proof of Theorem \ref{theo:discretemainequidup}.  Let
$\wt \psi^\pm\in \C^\beta_{\rm c}(\gengeod \XX)$.
In order to simplify the notation, let $\lambda_\ga=\lambda_{e,\ga}$,
$\alpha_\ga= \alpha_{e,\ga}$, $\alpha^-_\ga= \alpha^-_{e,\ga}$,
$\alpha^+_\ga=\alpha^+_{\ga^{-1},e}$ and $\wt \sigma^\pm=\wt
\sigma^\pm_{\DD^\mp}$.

Let us first prove the following avatar of Equation
\eqref{eq:geneth13cor19errsimp}, indicating only the required changes
in its proof: there exists $\kappa_0>0$ (independent of $\wt
\psi^\pm$) such that, as $n\ra+\infty$,
\begin{align}
\frac{e^{\delta_c}-1}{e^{\delta_c}}&\;\|m_c\|\;e^{-\delta_c\, n}
\sum_{\ga\in\Ga,\,0<\lambda_\ga\leq n} \;e^{\wt c\,(\alpha_\ga)}\; 
\wt\psi^-(\alpha^-_\ga)\,\wt\psi^+(\alpha^+_\ga)
\nonumber\\=\; &
\int_{\normalout D^-} \wt\psi^-\, d\wt\sigma^+
\int_{\normalin D^+}  \wt\psi^+\, d\wt\sigma^-
+\bigO(e^{-\kappa_0n} \;\|\wt\psi^-\|_\beta\,\|\wt\psi^+\|_\beta)
\label{eq:geneth13cor19errsimpavatar}\,.
\end{align}
Most of the new work to be done in order to prove this formula
concerns regularity properties of the test functions that will be
introduced later on.

We fix $R>e^2$ large enough. Let $0<\eta<1$. We introduce the
following modification of the test functions
$\phi^\pm_\eta$:\footnote{See Equation \eqref{eq:defiphi} for the
  definition of $\phi^\pm_\eta$ and Equation \eqref{eq:defihpm} for
  the definition of $h^\pm_{\eta,\,R}$, that simplifies as
  $h^\mp_{\eta,\,R}(w)=(\mu_{\wssu(w)}(B^\pm(w,R)))^{-1}$ since $\XX$
  is simplicial, as seen in Equation \eqref{eq:newhpmvalue}.}
$$
\Phi^\pm_\eta=(h^\pm_{\eta,\,R}\; \wt\psi^\pm)\circ f^\mp_{D^\pm}
\;\; \mathbbm{1}_{\V^\mp_{\eta,\,R}(\normalmp D^{\pm})} \,.
$$
As in Lemma \ref{lem:integrable}, the functions $\Phi^\pm_\eta$ are
measurable and satisfy
\begin{equation}\label{eq:controlL1testfun}
\int_{\G\XX} \Phi^\pm_\eta\;d\wt m_c=
\int_{\normalmp \DD^{\pm}}\wt \psi^\pm\,d\wt\sigma^\mp\,.
\end{equation}
\blemm The maps $\Phi^\pm_\eta$ are $\beta$-H\"older-continuous with
\begin{equation}\label{eq:conrolPhinormhold}
\|\Phi_\eta^\pm\|_\beta=\bigO(\|\wt \psi^\pm\|_\beta)\;.
\end{equation}
\elemm

\dem Since $\XX$ is a simplicial tree and $\eta<1$, we have
$V^\pm_{w,\,\eta,\,R}=B^\pm(w,R)$ for every $w\in \normalpm
D^{\mp}$. By the proof of Lemma
\ref{lem:troisboulessimp},\footnote{See also the proof of Lemmas
  \ref{lem:troisboulesmet} and \ref{lem:fibrationonelip}.} there
exists $c_R>0$ depending only on $R$ such that if $\ell'\in\G\XX$
satisfy $d(\ell,\ell')\leq c_R$, then $\ell'$ coincides with $\ell$ on
$\pm[0,\ln R+1]$. Therefore, if $\ell\in B^\pm(w,R)$ for some
$w\in\normalpm D^\mp$ and $d(\ell,\ell')\leq c_R$, then $\ell'$
coincides with $w$ on $\pm[\ln R,\ln R+1]$, thus $\ell'\in
B^\pm(w',R)$ where $w'\in \normalpm D^\mp$ is the geodesic ray with
$w'(0)=w(0)$ and $w'_\pm=\ell'_\pm$. Hence (see Section
\ref{subsec:holdercont}) the characteristic function
$\mathbbm{1}_{\V^\pm_{\eta,\,R}(\normalpm D^{\mp})}$ is
$c_R$-locally constant, thus $\beta$-H\"older-continuous by Remark
\ref{rem:locconstholder}.

By Assumption (1) in the statement of Theorem
\ref{theo:equidsimplicialperperrorterm}, the denominator of
$$
h^\mp_{\eta,\,R}(w) =\frac{1}{\mu_{\wssu(w)}(B^\pm(w,R))}
$$ 
is at least a positive constant depending only on $R$, hence
$h^\mp_{\eta,\,R}$ is bounded by a constant depending only on
$R$. Since the map $\mathbbm{1}_{B^\pm(w,R)}$ is $c_R$-locally
constant, so is the map $h^\mp_{\eta,\,R}$.  The result then follows
from Lemma \ref{lem:fibrationonelip} and Equation
\eqref{eq:prodholdernorm}. \cqfd

\medskip
In order to prove Equation \eqref{eq:geneth13cor19errsimpavatar}, as
in the proofs of Theorems \ref{theo:expratecount} and
\ref{theo:discretemainequidup}, for all $N\in\NN$, we estimate in two
ways the quantity
\begin{equation}\label{eq:defiIetapmTerrorsimplicial}
I_{\eta}(N)= \sum_{n=0}^{N} e^{\delta_c\,n}\sum_{\ga\in\Ga}\;
\int_{\ell\in\G \XX}\Phi^-_\eta(\flow{-\lfloor n/2\rfloor}\ell)\;
\Phi^+_\eta(\flow{\lceil n/2\rceil}\ga^{-1}\ell)\;d\wt m_{c}(\ell)\,.
\end{equation}

On the one hand, as in order to obtain Equation
\eqref{eq:firstestimerror}, using now Assumption (2) in the statement
of Theorem \ref{theo:equidsimplicialperperrorterm} on the exponential
mixing for the discrete time geodesic flow, a geometric sum argument
and Equations \eqref{eq:controlL1testfun} and
\eqref{eq:conrolPhinormhold}, we have
\begin{align}
I_{\eta}(N)=
\frac{e^{\delta_c(N+1)}}{(e^{\delta_c}-1)\,\|m_c\|}\; \Big(
\int_{\normalout \DD^{-}}&\wt \psi^-\,d\wt\sigma^+
\int_{\normalin \DD^{+}}\wt \psi^+\,d\wt\sigma^-\;+\;\nonumber\\ &
\bigO(e^{-\kappa\, N}
\|\wt \psi^-\|_\beta\|\wt \psi^+\|_\beta)\Big)\,.
\label{eq:firstestimerrorsimplicial}
\end{align}

On the other hand, exchanging the summations over $\ga$ and $n$ in the
definition of $I_{\eta}(N)$, we have
$$
I_{\eta}(N)=\sum_{\ga\in\Ga}\;\sum_{n=0}^N e^{\delta_c\,n}\;
\int_{\G\XX}\Phi^-_\eta(\flow{-\lfloor n/2\rfloor}\ell)\;
\Phi^+_\eta(\flow{\lceil n/2\rceil}\ga^{-1}\ell)\;d\wt m_{c}(\ell)\,.
$$
With the simplifications in Step 3T of the proof of Theorem
\ref{theo:mainequidup} given by the proof of Theorem
\ref{theo:discretemainequidup}, if $\eta< \frac{1}{2}$, if
$\ell\in\G\XX$ belongs to the support of
$\Phi^-_\eta\circ\flow{-\lfloor n/2\rfloor}\;
\Phi^+_\eta\circ\flow{\lceil n/2\rceil}\circ\ga^{-1}$, setting $w^-=
f^+_{D^-}(\ell)$ and $w^+=f^-_{\ga D^+} (\ell)$, we then have
$\lambda_\ga=n$, $w^-(0)=\alpha^-_\ga(0)$, $w^+(0)=\ga\alpha^+_\ga(0)$ and
$$
w^-(\lfloor n/2\rfloor)=w^+(-\lceil n/2\rceil)=\ell(0)=
\alpha^-_\ga(\lfloor n/2\rfloor)=\ga\alpha^+_\ga(-\lceil n/2\rceil)\;.
$$
Hence by the triangle inequality
\begin{align*}
d(w^-,\alpha^-_\ga)&=
\int_{\lfloor n/2\rfloor}^{+\infty} d(w^-(s),\alpha^-_\ga(s))\;e^{-2s}\;ds
\leq e^{-2\lfloor n/2\rfloor}\int_0^{+\infty} 2s'\,e^{-2s'}\;ds'
\\ & =\bigO(e^{-\lambda_\ga})\,.
\end{align*}
Similarly, $d(w^+,\ga\alpha^+_\ga)= \bigO(e^{-\lambda_\ga})$.
Therefore, since $\wt\psi^\pm$ is $\beta$-H\"older-continuous,
$$
|\,\wt\psi^-(w^-)-\wt\psi^-(v^-_\ga)\,|,\;\;
|\,\wt\psi^+(\ga^{-1}w^+)-\wt\psi^+(v^\pm_\ga)\,|=
\bigO(e^{-\beta\lambda_\ga}\;\|\wt\psi^\pm\|_\beta)\,.
$$ 
Note that now $\Phi^\pm_\eta=\wt\psi^\pm \circ f^{\mp}_{D^\pm} \;\phi_
\eta^\pm$, so that
\begin{align*}
I_{\eta}(N)=
\sum_{\ga\in\Ga}&\;\;\;\big(\wt\psi^-(\alpha^-_\ga)\wt\psi^+(\alpha^+_\ga)+
\bigO(e^{-2\beta\lambda_\ga}\;\|\wt\psi^-\|_\beta
\|\wt\psi^+\|_\beta)\big)\;\times \\ & \sum_{n=0}^N e^{\delta_c\,n}
\int_{\G\XX}\phi^-_\eta(\flow{-\lfloor n/2\rfloor}\ell)\;
\phi^+_\eta(\flow{\lceil n/2\rceil}\ga^{-1}\ell)\;d\wt m_{c}(\ell)\,.
\end{align*}

Now if $\eta<\frac{1}{2}$, Equation
\eqref{eq:geneth13cor19errsimpavatar} with $\kappa_0=
\min\{2\beta, \kappa\}$ follows as in Steps 3T and 4T of the
proof of Theorem \ref{theo:mainequidup} with the simplifications given
by the proof of Theorem \ref{theo:discretemainequidup}.

The end of the proof of the equidistribution claim of Theorem
\ref{theo:equidsimplicialperperrorterm} follows from Equation
\eqref{eq:geneth13cor19errsimpavatar} as the one of Theorem
\ref{theo:expratecount} from Equation \eqref{eq:geneth13cor19errsimp}.

The counting claim follows from the equidistribution one by taking
$\psi^\pm=\mathbbm{1}_{\Ga\V_{\eta,R}(\normalmp\DD^\pm)}$, which has
compact support since $\Ga_{\DD^\pm}\bs \partial\DD^\pm$ is assumed to
be compact, and is $\beta$-H\"older-continuous by previous arguments.
\cqfd

\medskip
\noindent {\bf Remarks. } (i) Assume that $\wt c=0$, that the
simplicial tree $\XX'$ with $|\XX'|_1=\C\Lambda\Ga$ is uniform without
vertices of degree $2$, that $L_\Ga=\ZZ$ and that $\Ga$ is a
geometrically finite lattice of $\XX'$.  Then all assumptions of
Theorem \ref{theo:equidsimplicialperperrorterm} are satisfied by the
results of Section \ref{subsec:ergodictrees} and by Corollary
\ref{coro:expdecaygeomfinisimpl}. Therefore we have an exponentially
small error term in the (joint) equidistribution of the common
perpendiculars, and in their counting if $\Ga_{\DD^\pm}\bs
\partial\DD^\pm$ is compact, see Examples \ref{ex:cpregularsubgraphs}
(2) and \ref{ex:cpregularsubgraphsbi} (4).

\medskip\noindent (ii) Assume in this remark that Assumption (2) of
the above theorem is replaced by the assumptions that $\C\Lambda\Ga$
is uniform without vertices of degree $2$, that $L_\Ga=2\ZZ$, and that
there exists $\beta\in\;]0,1]$ such that the square of the discrete
time geodesic flow on $(\Ga\bs\Geven \XX,m_c)$ is exponentially mixing
for the $\beta$-H\"older regularity, for instance by Corollary
\ref{coro:expdecaygeomfinisimpl} (2) if $\Ga$ is geometrically finite.
Then a similar proof (replacing the references to Theorem
\ref{theo:discretemainequidup} by references to Theorem
\ref{theo:discretemainequidupeven}) shows that there exists
$\kappa'>0$ such that for all $\psi^\pm\in \C_{\rm c}^\beta
(\Ga\bs\gengeod\XX)$, we have, as $n\ra+\infty$,
\begin{align*}
&\frac{e^{2\delta_c}-1}{2\,e^{2\delta_c}}\;\|m_{c}\|\;e^{-\delta_c\, n}
\sum_{\substack{[\ga]\,\in \Ga_{\DD^-}\bs\Ga/\Ga_{\DD^+}\\
0<d(\DD^-,\ga \,\DD^+)\leq n}} 
m_{e,\,\ga}\;e^{\wt c\,(\alpha_{e,\,\ga})}\;
\psi^-(\Ga \alpha^-_{e,\ga})\;\psi^+(\Ga\alpha^+_{\ga^{-1},e}) \\= \;&
\int \psi^-\;d\sigma^+_{\D^-}\;\int \psi^+\;d\sigma^-_{\D^+} 
+ \bigO\big(e^{-\kappa' \,n}\;\|\psi^-\|_\beta\;\|\psi^+\|_\beta\big)
\end{align*}
and if $\Ga_{\DD^\pm}\bs \partial\DD^\pm$ is compact, then
$$
\N_{\DD^-,\,\DD^+}(n)=
\frac{2\,e^{2\delta_c}\;\|\sigma^+_{\D^-}\|\;\|\sigma^-_{\D^+}\|}
{(e^{2\delta_c}-1)\;\|m_{c}\|}\;\;e^{\delta_c\, n} + 
\bigO\big(e^{(\delta_c-\kappa')n}\big)\,.
$$

\bigskip
Let us now consider the metric tree case, for the continuous time
geodesic flow, where the main change is to assume a superpolynomial
decay of correlations and hence get a superpolynomial error term. We
refer to the beginning of Section \ref{subsec:mixingratemetgraphs} for
the definitions of the function space $\C_{\rm c}^{k,\,\beta}
(\Ga\bs\gengeod X)$ and the superpolynomial mixing.

\btheo\label{theo:equidmetricperperrorterm} Let $(\XX,\lambda)$,
$\Ga$, $\wt c$ and $\DD^\pm$ be as in the beginning of this Section,
and let $D^\pm=|\DD^\pm|_\lambda$.  Assume that the critical exponent
$\delta_c$ is finite and positive, that the Gibbs measure $m_{c}$ (for
the continuous time geodesic flow) is finite and that the skinning
measures $\sigma^\pm_{\D^\mp}$ are finite and nonzero. Assume
furthermore that
\begin{enumerate}
\item at least one of the following holds :

$\bullet$~ $\Ga_{D^\pm}\bs \partial D^\pm$ is compact
 
$\bullet$~ the metric subtree $\C\Lambda\Ga$ is uniform and $\Ga$ is a
  lattice of $\C\Lambda\Ga$,
\item there exists $\beta\in\;]0,1]$ such that the continous time
geodesic flow on $(\Ga\bs\G X,m_c)$ has superpolynomial decay of
$\beta$-H\"older correlations.
\end{enumerate}
Then for every $n\in\NN$ there exists $k\in\NN$ such that for all
$\psi^\pm\in \C_{\rm c}^{k,\,\beta} (\Ga\bs\gengeod X)$, we have, as
$T\ra+\infty$,
\begin{align*}
&\delta_c\;\|m_{c}\|\;e^{-\delta_c\, T}
\sum_{\substack{[\ga]\,\in \Ga_{D^-}\bs\Ga/\Ga_{D^+}\\
0<d(D^-,\ga \,D^+)\leq T}} 
m_{e,\,\ga}\;e^{\wt c\,(\alpha_{e,\,\ga})}\;
\psi^-(\Ga \alpha^-_{e,\ga})\;\psi^+(\Ga\alpha^+_{\ga^{-1},e}) \\= \;&
\int_{\Ga\bs\gengeod X} \psi^-\;d\sigma^+_{\D^-}\;
\int_{\Ga\bs\gengeod X} \psi^+\;d\sigma^-_{\D^+} 
+ \bigO\big(T^{-n}\;\|\psi^-\|_{k,\,\beta}\;\|\psi^+\|_{k,\,\beta}\big)
\end{align*}
and if $\Ga_{D^\pm}\bs \partial D^\pm$ is compact, then for every $n\in\NN$
$$
\N_{D^-,\,D^+}(T)=
\frac{\|\sigma^+_{\D^-}\|\;\|\sigma^-_{\D^+}\|}
{\delta_c\;\|m_{c}\|}\;\;e^{\delta_c\, T} + 
\bigO\big(e^{\delta_c\, T}T^{-n}\big)\,.
$$
\etheo

\noindent {\bf Remark. } Assume that $\wt c=0$, that the metric tree
$\C\Lambda\Ga$ is uniform, either that $\Ga\bs\XX$ is finite and the
length spectrum $L_\Ga$ of $\Ga$ is $2$-Diophantine, or that $\Ga$ is
a geometrically finite lattice of $\C\Lambda\Ga$ and $L_\Ga$ is
$4$-Diophantine.  Then all assumptions of Theorem
\ref{theo:equidmetricperperrorterm} are satisfied by the results of
Section \ref{subsec:ergodictrees} and by Corollary
\ref{coro:expdecaygeomfinimet}. Therefore we have a superpolynomially
small error term in the (joint) equidistribution of the common
perpendiculars (and in their counting if $\Ga_{D^\pm}\bs \partial
D^\pm$ is compact).

\medskip
\dem The proof is similar to the one of Theorem
\ref{theo:equidsimplicialperperrorterm}, except that since the time is
now continuous, we need to regularise our test functions also in the
time direction in order to obtain the regularity required for the
application of the assumption on the mixing rate. We again use the
simplifying notation $\lambda_\ga=\lambda_{e,\ga}$, $\alpha_\ga=
\alpha_{e,\ga}$, $\alpha^-_\ga= \alpha^-_{e,\ga}$,
$\alpha^+_\ga=\alpha^+_{\ga^{-1},e}$ and $\wt \sigma^\pm=\wt
\sigma^\pm_{D^\mp}$.

We fix $n\in\NN-\{0\}$. Using the rapidly mixing property, there exists a
regularity $k$ such that for all $\psi,\psi'\in\C^{k,\,\beta}_{\rm
  b} (\Ga\bs\G X)$ we have as $t\ra+\infty$
\begin{equation}\label{eq:rapidmixingpourperp}
\operatorname{cov}_{\,\overline{m_c},\,t}\,(\psi,\psi')=
\bigO(t^{-N\,n}\;\|\psi\|_{k,\,\beta}\;\|\psi'\|_{k,\,\beta})\;,
\end{equation} 
where $N\in\NN-\{0\}$ is a constant which will be made precise later on.

Let us first prove that for all $\wt \psi^\pm\in \C^{k,\,\beta}_{\rm
  c} (\gengeod X)$, we have, as $T\ra+\infty$,
\begin{align}
&\delta_c\;\|m_c\|\;e^{-\delta_c\, T}
\sum_{\ga\in\Ga,\,0<\lambda_\ga\leq T} \;e^{\wt c\,(\alpha_\ga)}\; 
\wt\psi^-(\alpha^-_\ga)\,\wt\psi^+(\alpha^+_\ga)
\nonumber\\=\; &
\int_{\normalout D^-} \wt\psi^-\, d\wt\sigma^+
\int_{\normalin D^+}  \wt\psi^+\, d\wt\sigma^-
+\bigO(T^{-n} \;\|\wt\psi^-\|_{k,\,\beta}\,\|\wt\psi^+\|_{k,\,\beta})
\label{eq:geneth13cor19errmetavatar}\,.
\end{align}

In order to prove this formula, we introduce modified test functions
with bounded H\"older-continuous derivatives up to order $k$ (by a
standard construction) in the time direction (the stable leaf and
unstable leaf directions remain discrete). We fix $R>0$ large enough.

For every $\eta\in\;]0,1[\,$, there exists a map $\wh{\mathbbm{1}_\eta}:
\RR\ra[0,1]$ which has bounded $\beta$-H\"older-continuous derivatives
up to order $k$, which is equal to $0$ if $t\notin[-\eta,\eta]$ and to
$1$ if $t\in [-\eta\, e^{-\eta},\eta\, e^{-\eta}]$ (when $k=0$, just
take $\wh{\mathbbm{1}_\eta}$ to be continuous and affine on each
remaining segment $[-\eta,-\eta\, e^{-\eta}]$ and $[\eta\, e^{-\eta},
\eta]$), such that, for some constant $\kappa'_1>0$,
$$
\|\wh{\mathbbm{1}_\eta}\|_{k,\,\beta}=\bigO(\eta^{-\kappa'_1})\;.
$$ 
Using leafwise this regularisation process, there exists
$\chi^\pm_{\eta,\,R}\in \C^{k,\,\beta}_{\rm b}(\G X)$ such that

$\bullet$~ $\|\chi^\pm_{\eta,\,R}\|_{k,\,\beta}=\bigO(\eta^{-\kappa'_1})$,

$\bullet$~ 
$\mathbbm{1}_{\V^\mp_{\eta \,e^{-\eta},\,R}(\normalmp D^{\pm})}
\leq \chi^\pm_{\eta,\,R}\leq \mathbbm{1}_{\V^\mp_{\eta,\,R}(\normalmp D^{\pm})}$,

$\bullet$~ for every $w\in\normalmp D^{\pm}$, we have 
$$
\int_{V^\mp_{w,\,\eta,\,R}}\chi^\pm_{\eta,\,R}\,d\nu_w^\pm=
\nu_w^\pm(V^\mp_{w,\,\eta,\,R}) \,e^{-\bigO(\eta)}=
\nu_w^\pm(V^\mp_{w,\,\eta \,e^{-\eta},\,R}) \,e^{\bigO(\eta)}\,.
$$ 
As in the proof of Theorem \ref{theo:expratecount} in the manifold
case, the new test functions are defined, with
$$
H^\pm_{\eta,\,R}\;:\; w\in \normalmp D^{\pm}\mapsto
\frac{1}{\int_{V^\mp_{w,\,\eta,\,R}}\chi^\pm_{\eta,\,R}\,d\nu_w^\pm} \;,
$$
by
$$
\Phi^\pm_\eta=(H^\pm_{\eta,\,R}\; \wt\psi^\pm)\circ f^\mp_{D^\pm}
\;\; \chi^\pm_{\eta,\,R}\;\;:\;\G X \ra \RR\,.
$$ 
Let $\wh\Phi^\pm_\eta=H^\pm_{\eta,\,R}\circ f^\mp_{D^\pm} \;
\chi^\pm_{\eta,\,R}$, so that $\Phi^\pm_\eta=\wt\psi^\pm\circ
f^\mp_{D^\pm}\;\wh\Phi^\pm_\eta$. By the last two properties of the
regularised maps $\chi^\pm_{\eta,\,R}$, we have, with $\phi^\mp_\eta$
defined as in Equation \eqref{eq:defiphi},
\begin{equation}\label{eq:controlwhPhimettree}
\phi^\pm_{\eta\,e^{-\eta}}\,e^{-\bigO(\eta)}\leq \wh\Phi^\pm_\eta\leq 
\phi^\pm_\eta\,e^{\bigO(\eta)}\;.
\end{equation}

By Assumption (1), if $R$ is large enough, by the definitions of
the measures $\nu_w^\pm$, the denominator of $H^\pm_{\eta,\,R}(w)$ is
at least $c \,\eta$ where $c>0$. As in the proof of Theorem
\ref{theo:expratecount},
$$
\int_{\G X} \Phi^\pm_\eta\;d\wt m_c=
\int_{\normalmp D^{\pm}}\wt \psi^\pm\,d\wt\sigma^\mp\,
$$ 
and  there exists $\kappa''>0$ such that
$$
\|\Phi^\pm_\eta\|_{k,\,\beta}=
\bigO(\eta^{-\kappa''}\|\wt \psi^\pm\|_{k,\,\beta})\,.
$$ 
We again estimate in two ways as $T\ra+\infty$ the quantity
\begin{equation}\label{eq:defiIetapmTerrormetric}
I_{\eta}(T)= \int_{0}^{T} e^{\delta_c\,t}\sum_{\ga\in\Ga}\;
\int_{\ell\in\G X}\Phi^-_\eta(\flow{-t/2}\ell)\;
\Phi^+_\eta(\flow{t/2}\ga^{-1}\ell)\;d\wt m_{c}(\ell)\,dt\,.
\end{equation}
Note that as $T\ra+\infty$,
\begin{align*}
e^{-\delta_c\,T}\;\int_1^T e^{\delta_c\,t}\, t^{-N\,n}\,dt &=
e^{-\delta_c\,T}\;\int_1^{T/2} e^{\delta_c\,t}\, t^{-N\,n}\,dt+
e^{-\delta_c\,T}\;\int_{T/2}^T e^{\delta_c\,t}\, t^{-N\,n}\,dt\\ & 
=\bigO(e^{-\delta_c\,T/2})+\bigO(T^{-N\,n+1})=\bigO( T^{-(N-1)\,n})\;.
\end{align*}
Using Equation \eqref{eq:rapidmixingpourperp}, an integration argument
and the above two properties of the test functions, we hence have
\begin{align}
I_{\eta}(T)=
\frac{e^{\delta_c\,T}}{\delta_c\,\|m_c\|}\;\Big(
\int_{\normalout D^{-}}&\wt \psi^-\,d\wt\sigma^+
\int_{\normalin D^{+}}\wt \psi^+\,d\wt\sigma^-\;+\;\nonumber\\ &
\bigO(T^{-(N-1)\,n}\eta^{-2\kappa''}
\|\wt \psi^-\|_{k,\,\beta}\|\wt \psi^+\|_{k,\,\beta})\Big)\,.
\label{eq:firstestimerrormetric}
\end{align}

As in Step 3T of the proof of Theorem \ref{theo:mainequidup}, for all
$\ga\in\Ga$ and $t>0$ large enough, if $\ell\in\G X$ belongs to the
support of $\Phi^-_\eta\circ\flow{-t/2}\; \Phi^+_\eta\circ \flow{t/2}
\circ\ga^{-1}$ (which is contained in the support of
$\phi^-_\eta\circ\flow{-t/2}\; \phi^+_\eta\circ \flow{t/2}
\circ\ga^{-1}$~), then we may define $w^-= f^+_{D^-}(\ell)$ and $w^+=
f^-_{\ga D^+} (\ell)$.

By the property (iii) in Step 3T of the proof of Theorem
\ref{theo:mainequidup}, the generalised geodesic lines $w^-$ and
$\alpha_ \ga^-$ coincide, besides on $]-\infty, 0]$, at least on
$[0, \frac{t}{2}-\eta]$, and similarly, $w^+$ and $\ga\alpha_\ga^+$
coincide, besides on $[0,+\infty[\,$, at least on $[-\frac{t}{2}
+\eta, 0]$. Therefore, by an easy change of variable and since 
$|\frac{t}2-\frac{\lambda_\ga}2|\leq \eta$,
\begin{align*}
d(w^-,\alpha_ \ga^-) &\leq \int_{\frac{t}{2}-\eta}^{+\infty}
d(w^-(s),\alpha_ \ga^-(s))\,e^{-2s}\,ds \leq
e^{-2(\frac{t}{2}-\eta)}\;\int_0^{+\infty} 2s\,e^{-2s}\,ds\\ & = 
\bigO(e^{-t}) = \bigO(e^{-\lambda_\ga}) \;.
\end{align*}
Similarly, $d(w^+,\ga\alpha_\ga^+)= \bigO(e^{-\lambda_\ga})$. Hence
since $\wt\psi^\pm$ is $\beta$-H\"older-continuous, we have
$$
|\,\wt\psi^-(w^-)-\wt\psi^-(\alpha^-_\ga)\,|,\;\;
|\,\wt\psi^+(\ga^{-1}w^+)-\wt\psi^+(\alpha^+_\ga)\,|=
\bigO(e^{-\beta \lambda_\ga} \|\wt\psi^\pm\|_\beta)\,.
$$

Therefore, as in the proof of Theorem
\ref{theo:equidsimplicialperperrorterm}, we have
\begin{align*}
I_{\eta}(T)=\sum_{\ga\in\Ga}&\;\;\;\big(\wt\psi^-(\alpha^-_\ga)
\wt\psi^+(\alpha^+_\ga) + \bigO(e^{-2\beta \lambda_\ga}
\|\wt\psi^-\|_\beta\|\wt\psi^+\|_\beta)\big)\times
 \\ & \int_0^{T}e^{\delta\,t}\;\int_{\ell\in \G X} 
\;\wh\Phi^-_{\eta}(\flow{-t/2}\ell)\;
\wh\Phi^+_{\eta}(\ga^{-1}\flow{t/2}\ell)\;d\wt m_c(\ell)\;dt\,.
\end{align*}

Finally, Equation \eqref{eq:geneth13cor19errmetavatar} follows as in
the end of the proof of Equation \eqref{eq:geneth13cor19errsimp},
using Equations \eqref{eq:controlwhPhimettree} and
\eqref{eq:effcontrojtree} instead of Equations \eqref{eq:controlwhPhi}
and \eqref{eq:step5bis}, by taking $\eta=T^{-n}$ and $N=
2(\lceil\kappa''\rceil+1)$.

The end of the proof of the equidistribution claim of Theorem
\ref{theo:equidmetricperperrorterm} follows from Equation
\eqref{eq:geneth13cor19errmetavatar} as the one of Theorem
\ref{theo:expratecount} from Equation \eqref{eq:geneth13cor19errsimp}.

The counting claim follows from the equidistribution one by taking
$\psi^\pm$ to be $\beta$-H\"older-continuous plateau functions around
$\Ga\V_{\eta,R}(\normalmp\DD^\pm)$.
\cqfd

\bigskip
We are now in a position to prove one of the counting results in
the introduction.

\medskip
\noindent{\bf Proof of Theorem \ref{theo:countintrograph}.} %
Let $\XX$ be the universal cover of $\YY$, with fundamental group
$\Ga$ for an indifferent choice of basepoint, and let $\DD^\pm$ be
connected components of the preimages of $\YY^\pm$ in $\XX$.
Assertion (1) of Theorem \ref{theo:countintrograph} follows from
Theorem \ref{theo:equidmetricperperrorterm} and its subsequent
Remark. Assertion (2) of Theorem \ref{theo:countintrograph} follows
from Theorem \ref{theo:equidsimplicialperperrorterm} and its
Remarks (ii) and (i) following its proof, respectively, if $\YY$ is
bipartite or not.  
\cqfd

\chapter{Geometric applications}
\label{sect:geomappli}

In this final Chapter of Part \ref{part:equid}, we apply the
equidistribution and counting results obtained in the previous
Chapters in order to study geometric equidistribution and counting
problems for metric and simplicial trees concerning conjugacy classes
in discrete isometry groups and closed orbits of the geodesic flows.

\section{Orbit counting in conjugacy classes 
for groups acting on  trees}
\label{sect:orbcountconjugacyclass}

In this Section, we study the orbital counting problem for groups
acting on metric or simplicial trees when we consider only the images
by elements in a given conjugacy class. We refer to the Introduction
for motivations and previously known results for manifolds (see
\cite{Huber56} and \cite{ParPau15MZ}) and graphs (see \cite{Douma11}
and \cite{KenSha15}). The main tools we use are Theorem
\ref{theo:counting} for the metric tree case and Theorem
\ref{theo:equidsimplicialperperrorterm} for the simplicial tree case,
as well as their error terms. In particular, we obtain a much more
general version of Theorem \ref{theo:genekenisonsharpintro} in the
Introduction.

\medskip
Let $(\XX,\lambda)$ be a locally finite metric tree without terminal
vertices, let $X=|\XX|_\lambda$ be its geo\-metric realisation, let
$x_0\in V\XX$ and let $\Ga$ be a nonelementary discrete subgroup of
$\Aut(\XX,\lambda)$.\footnote{See Section \ref{subsec:trees} for
  definitions and notation.}  Let $\wt c:E\XX\ra\RR$ be a
$\Ga$-invariant system of conductances, let $\wt F_c$ and $F_c$ be its
associated potentials on $T^1X$ and $\Ga\bs T^1X$ respectively, and
let $\delta_c=\delta_{\Ga,\,F_c^\pm}$ be its critical
exponent.\footnote{See Section \ref{subsec:cond} for definitions and
  notation.}  Let $(\mu^\pm_x)_{x\in X}$ (respectively
$(\mu^\pm_x)_{x\in V\XX}$) be (normalised) Patterson densities for the
pairs $(\Ga,F_c^\pm)$, and let $\wt m_{c}=\wt m_{F_c}$ and
$m_c=m_{F_c}$ be the associated Gibbs measures on $\G X$ and $\Ga\bs
\G X$ (respectively $\G \XX$ and $\Ga\bs \G \XX$) for the continuous
time geodesic flow (respectively the discrete time geodesic flow, when
$\lambda\equiv 1$).\footnote{See Sections
  \ref{subsec:pattersongibbstrees} and \ref{subsec:ergodictrees} for
  definitions and notation.}

Recall that the {\em virtual centre}\index{virtual centre} $Z^{\rm
  virt}(\Ga)$ of $\Ga$ is the finite (normal) subgroup of $\Ga$
consisting of the elements $\ga\in\Ga$ acting by the identity on the
limit set $\Lambda\Ga$ of $\Ga$ in $\partial_\infty X$, see for
instance \cite[\S 5.1]{Champetier00}. If $\Lambda\Ga=\partial_\infty
X$ (for instance if $\Ga$ is a lattice), then $Z^{\rm virt} (\Ga)=
\{\id\}$.

For any nontrivial element $\ga$ in $\Ga$ with translation length
$\len(\ga)$ in $X$, let $C_\ga$ be

\smallskip
$\bullet$~ the translation axis of $\ga$ if $\ga$ is loxodromic on
$X$,

\smallskip
$\bullet$~ the fixed point set of $\ga$ if $\ga$ is elliptic on $X$,

\smallskip
\noindent and let $\Ga_{C_\ga}$ be the stabiliser of $C_\ga$ in
$\Ga$. In the simplicial case (that is, when $\lambda\equiv 1$),
$C_\ga$ is a simplicial subtree of $\XX$. Note that
$\lambda(\ga)=\lambda(\ga'\ga(\ga')^{-1})$ and $\ga' C_{\ga} =
C_{\ga'\ga(\ga')^{-1}}$ for all $\ga'\in \Ga$, and that for any
$x_0\in X$
\begin{equation}\label{eq:treeequitranslat}
d(x_0,C_\ga)=\frac{d(x_0,\ga x_0)-\len(\ga)}{2}\;.
\end{equation}

Let $\D=(\ga' C_\ga)_{\ga'\in \Ga/\Ga_{C_\ga}}$, which is a locally
finite $\Ga$-invariant family of nonempty proper (since $\ga\neq \id$)
closed convex subsets of $X$.\footnote{See Section
  \ref{subsec:equivfammult} for definitions and notations.}  By the
equivariance properties of the skinning measures, the total mass of
the skinning measure\footnote{See the previous footnote.}
$\sigma^-_{\D}$ depends only on the conjugacy class $\KKK$ of $\ga$ in
$\Ga$, and will be denoted by $\|\sigma^-_\KKK\|$. This quantity,
called the {\em skinning measure}\index{skinning measure} of $\KKK$,
is positive unless $\partial_\infty C_\ga=\Lambda\Ga$, which is
equivalent to $\ga\in Z^{\rm virt}(\Ga)$ (and implies in particular
that $\ga$ is elliptic).  Furthermore, $\|\sigma^-_\KKK\|$ is finite
if $\ga$ is loxodromic, and it is finite if $\ga$ is elliptic and
$\Ga_{C_\ga}\bs (C_\ga\cap \C\Lambda\Ga)$ is compact. This last
condition is in particular satisfied if $C_\ga\cap \C\Lambda\Ga$
itself is compact, and this is the case for instance if, for some
$k\geq 0$, the action of $\Ga$ on $X$ is {\em
  $k$-acylindrical}\index{acylindrical} (see for instance
\cite{Sela97,GuiLev11}), that is, if any element of $\Ga$ fixing a
segment of length $k$ in $\C\Lambda\Ga$ is the identity.

For every $\ga\in \Ga-\{e\}$, we define 
$$
m_\ga= \frac{1}{\card(\Ga_{x_0}\cap \Ga_{C_\ga})}\,,
$$ 
which is a natural multiplicity of $\ga$, and equals $1$ if the
stabiliser of $x_0$ in $\Ga$ is trivial (for instance if $\Ga$ is
torsion-free).  Note that for every $\beta\in\Ga$, the real number
$m_{\beta\ga\beta^{-1}}$ depends only on the double coset of
$\beta$ in $\Ga_{x_0}\bs\Ga/\Ga_{C_{\ga}}$.

The centraliser $Z_\Ga(\ga)$ of $\ga$ in $\Ga$ is contained in the
stabiliser of $C_\ga$ in $\Ga$. The index 
$$
i_\KKK=[\Ga_{C_\ga}:Z_\Ga(\ga)]
$$
depends only on the conjugacy class $\KKK$ of $\ga$; it will be
called the {\em index}\index{index} of $\KKK$. The index
$i_\KKK$ is finite if $\ga$ is loxodromic (the stabiliser of its
translation axis $C_\ga$ is then virtually cyclic), and also finite if
$C_\ga$ is compact (as for instance if the action of $\Ga$ on $X$ is
$k$-acylindrical for some $k\geq 0$).

We define 
$$
c_\ga= \sum_{i=1}^k \wt c\,(e_i)\,\lambda(e_i)\,,
$$ 
where $(e_1,\dots, e_k)$ is the shortest edge path from $x_0$ to
$C_\ga$.

\medskip We finally define the  orbital counting function in
conjugacy classes, counting with multiplicities and weights
  coming from the system of conductances, as
$$
N_{\KKK,\,x_0}(t)=
\sum_{\alpha\in\KKK,\; d(x_0,\,\alpha x_0)\leq t} m_\alpha\;e^{c_\alpha}\,.
$$ 
for $t\in [0,+\infty[$ (simply $t\in\NN$ in the simplicial case). When
the stabiliser of $x_0$ in $\Ga$ is trivial and when the system of
conductances $c$ vanishes, we recover the definition of the
Introduction (above Theorem \ref{theo:genekenisonsharpintro}).

\btheo\label{orbcountconjuggene} Let $\KKK$ be the conjugacy class of
a nontrivial element $\ga_0$ of $\Ga$, with finite index $i_\KKK$, and
with positive and finite skinning measure $\|\sigma^-_\KKK\|$. Assume
that $\delta_c$ is finite and positive.

\medskip\noindent (1) Assume that $m_{c}$ is finite and mixing for the
continuous time geodesic flow on $\Ga\bs X$.  Then, as $t\ra+\infty$,
$$
N_{\KKK,\,x_0}(t)\sim 
\frac{i_{\KKK}\,\|\mu^+_{x_0}\|\,\|\sigma^-_{\KKK}\|\;
e^{-\frac{\len(\ga_0)}{2}}}{\delta_c\,\|m_{c}\|}
\;e^{\frac{\delta_c}2\,t }\,.
$$ 
If $\Ga_{C_\ga}\bs (C_\ga\cap\C\Lambda\Ga)$ is compact when $\ga\in
\KKK$ is elliptic and if there exists $\beta\in\;]0,1]$ such that the
continous time geodesic flow on $(\Ga\bs\G X,m_c)$ has superpolynomial
decay of $\beta$-H\"older correlations, then the error term is
$\bigO\big(t^{-n}\,e^{\frac{\delta_c}2\,t }\big)$ for every $n\in\NN$.

\medskip\noindent(2) Assume that $\lambda\equiv 1$ and that $m_{c}$ is
finite and mixing for the discrete time geodesic flow on $\Ga\bs \G \XX$.
Then, as $n\ra+\infty$,
$$
N_{\KKK,\,x_0}(n)\sim
\frac{e^{\delta_c}\;i_{\KKK}\,\|\mu^+_{x_0}\|\,\|\sigma^-_{\KKK}\|\;
}{(e^{\delta_c}-1)\,\|m_{c}\|}
\;e^{\delta_c\lfloor\frac{n-\len(\ga_0)}2\rfloor}\,.
$$ 
If $\Ga_{C_\ga}\bs (C_\ga\cap\C\Lambda\Ga)$ is compact when $\ga\in
\KKK$ is elliptic and if there exists $\beta\in\;]0,1]$ such that the
discrete time geodesic flow on $(\Ga\bs\G \XX,m_c)$ is exponentially
mixing for the $\beta$-H\"older regularity, then the error term is
$\bigO\big(e^{(\delta_c-\kappa)n/2}\big)$ for some $\kappa>0$.  
\etheo

One can also formulate a version of the above result for groups acting
on bipartite simplicial trees based on Theorem
\ref{theo:countingdiscreteeven} and Remark (ii) following the proof of
Theorem \ref{theo:equidsimplicialperperrorterm}.

The error term in Assertion (1) holds for instance if $\wt c=0$, $X$
is uniform, and either $\Ga\bs X$ is compact and the length spectrum
$L_\Ga$ is $2$-Diophantine or $\Ga$ is a geometrically finite lattice
of $X$ whose length spectrum $L_\Ga$ is $4$-Diophantine, by the Remark
following Theorem \ref{theo:equidmetricperperrorterm}. When $\Ga\bs X$
is compact and $\Ga$ has no torsion (in particular, $\Ga$ has then a
very restricted group structure, as it is then a free group), we thus
recover a result of \cite{KenSha15}.

The error term in Assertion (2) holds for instance if $\wt c=0$, $\XX$
is uniform with vertices of degrees at least $3$, $\Ga$ is a
geometrically finite lattice of $\XX$ with length spectrum equal to
$\ZZ$, by Remark (i) following the proof of Theorem
\ref{theo:equidsimplicialperperrorterm}.

Theorem \ref{theo:genekenisonsharpintro} in the introduction follows
from this theorem, using Proposition \ref{prop:uniflatmBMfinie} (3)
and Theorem \ref{theo:uniflatmBMmixing}.

\medskip
\dem 
We only give a full proof of Assertion (1) of this theorem, Assertion
(2) follows similarly using Theorems \ref{theo:countingdiscrete} and
\ref{theo:equidsimplicialperperrorterm} instead of Theorems
\ref{theo:counting} and \ref{theo:equidmetricperperrorterm}.

\medskip 
The proof is similar to the proof of \cite[Theo.~8]{ParPau15MZ}.  Let
$D^-=\{x_0\}$ and $D^+=C_{\ga_0}$. Let $\D^-=(\ga D^-)_{\ga'\in
  \Ga/\Ga_{D^-}}$ and $\D^+=(\ga D^+)_{\ga\in \Ga/\Ga_{D^+}}$. By
Equation \eqref{eq:skinning singleton}, we have
$$
\|\sigma^+_{\D^-}\|= \frac{\|\mu^+_{x_0}\|}{|\Ga_{x_0}|}\,.
$$ 
By Equation \eqref{eq:treeequitranslat}, by the
definition\footnote{See Equation \eqref{eq:NDDn} in Section
  \ref{subsect:equicountmetricgraphgroup}.} of the counting function
$\N_{D^-,\,D^+}$ and by the last claim of Theorem \ref{theo:counting},
we have, as $t\ra+\infty$,
\begin{align*}
\sum_{\alpha\in\KKK,\; 0<d(x_0,\,\alpha x_0)\leq t} m_\alpha\;e^{c_\alpha}&=
\sum_{\alpha\in\KKK,\; 0<d(x_0,\,C_{\alpha})\leq \frac{t-\len(\ga_0)}{2}} m_\alpha\;e^{c_\alpha} 
\\ &= \sum_{\ga\in\Ga/Z_\Ga(\ga_0),\; 0<d(x_0,\,\ga C_{\ga_0})\leq \frac{t-\len(\ga_0)}{2}} 
m_{\ga\ga_0\ga^{-1}}\;e^{c_{\ga\ga_0\ga^{-1}}}\\ & = |\Ga_{x_0}|\,i_\KKK\;
\sum_{\ga\in\Ga_{x_0}\bs\Ga/\Ga_{C_{\ga_0}},\; 0<d(x_0,\,\ga C_{\ga_0})\leq \frac{t-\len(\ga_0)}{2}} 
m_{\ga\ga_0\ga^{-1}}\;e^{c_{\ga\ga_0\ga^{-1}}} \\  &
=|\Ga_{x_0}|\,i_\KKK\;\N_{D^-,\,D^+}\big(\frac{t-\len(\ga_0)}{2}\big)
\\  &\sim
|\Ga_{x_0}|\,i_\KKK\;\frac{\|\sigma^+_{\D^-}\|\;\|\sigma^-_{\D^+}\|}
{\delta_c\;\|m_{c}\|}\; e^{\delta_c\;\frac{t-\len(\ga_0)}{2}}\;.
\end{align*}
Assertion (1) without the error term follows, and
the error term statement follows similarly from Theorem
\ref{theo:equidmetricperperrorterm}.
\cqfd

\medskip
Theorem \ref{orbcountconjuggene} (1) without an explicit form of the
multiplicative constant in the asymptotic is due to \cite{KenSha15}
under the strong restriction that $\Ga$ is a free group acting freely
on $X$ and $\Ga\bs X$ is a finite graph.  The following result is due
to \cite[Theo.~1]{Douma11} in the special case when $\XX$ is a regular
tree and the group $\Ga$ has no torsion and finite quotient $\Ga\bs\XX$.

\bcoro Let $\XX$ be a regular simplicial tree with vertices of degree
$q+1\geq 3$, let $x_0\in V\XX$, let $\Ga$ be a lattice of $\XX$ such
that $\Ga\bs \XX$ is nonbipartite, and let $\KKK$ be the conjugacy
class of a loxodromic element  $\ga_0\in\Ga$.  Then, as $n\ra+\infty$,
$$
\sum_{\alpha\in\KKK,\; d(x_0,\,\alpha x_0)\leq n} m_\alpha \sim \frac{\len(\ga_0)}
{[Z_\Ga(\ga_0):\ga_0^\ZZ]\;\Vol(\Ga\dbs \XX)}
\;q^{\lfloor\frac{n-\len(\ga_0)}2\rfloor}\,.
$$ 
If we assume furthermore that $\Ga$ has no torsion and that $\ga_0$
is primitive, then
we have as $n\ra+\infty$,
$$
\card\{\alpha\in\KKK : d(x_0,\,\alpha x_0)\leq n\}\sim \frac{\len(\ga_0)}
{|\Ga\bs V\XX|}
\;q^{\lfloor\frac{n-\len(\ga_0)}2\rfloor}\,.
$$
\ecoro

\medskip \dem Under these assumptions, taking $c\equiv 0$ in Theorem
\ref{orbcountconjuggene} so that the Gibbs measure is the
Bowen-Margulis measure, the discrete time geodesic flow on $\Ga\bs\G \XX$
is finite and mixing by Proposition \ref{prop:uniflatmBMfinie} (3) and
Theorem \ref{theo:uniflatmBMmixing}. We also have $\delta_c=\ln q$. Using
the normalisation of the Patterson density $(\mu^\pm_x)_{x\in V \XX}$
to probability measures, Proposition \ref{prop:computBM} (3) and
Equation \eqref{eq:massskinline}, the result follows, since when $\ga$
is loxodromic,
$$
\Vol(\Ga_{C_\ga}\dbs C_\ga)=
\frac{\Vol(\ga^\ZZ\dbs C_\ga)}{[\Ga_{C_\ga}:\ga^\ZZ]}=
\frac{\len(\ga)}{[\Ga_{C_\ga}:Z_\Ga(\ga)]\;[Z_\Ga(\ga):\ga^\ZZ]}\,.\;.\;\;\;\Box
$$ 

\bigskip
The value of $C'$ given below Theorem \ref{theo:genekenisonsharpintro}
in the Introduction follows from this corollary.

\bigskip
We leave to the reader an extension with nonzero potential $F$ of the
results for manifolds in \cite{ParPau15MZ}, along the lines of the
above proofs.

\section{Equidistribution and counting of closed 
orbits on metric and simplicial graphs (of groups)}
\label{sec:equidclosedorbits}

Classically, an important characterisation of the Bowen-Margulis
measure on compact negatively curved Riemannian manifolds is that it
coincides with the weak-star limit of properly normalised sums of
Lebesgue measures supported on periodic orbits, see
\cite{Bowen72a}. Under much weaker assumptions than compactness, this
result was extended to $\CAT(-1)$ spaces with zero potential in
\cite{Roblin03} and to Gibbs measures in the manifold case in
\cite[Theo.~9.11]{PauPolSha15}. As a corollary of the simultaneous
equidistribution results from Chapter \ref{sec:equidarcs}, we prove in
this Section the equidistribution towards the Gibbs measure of
weighted closed orbits in quotients of metric and simplicial graphs of
groups and as a corollary of this result, in the standard manner, we
obtain asymptotic counting results for weighted (primitive) closed
orbits.

\medskip
Let $(\XX,\lambda)$ be a locally finite metric tree without terminal
vertices, and let $X=|\XX|_\lambda$ be its geometric realisation.  Let
$\Ga$ be a nonelementary discrete subgroup of $\Aut(\XX,\lambda)$.
Let $\wt c:E\XX\ra\RR$ be a $\Ga$-invariant system of conductances,
and $c:\Ga\bs E\XX\ra\RR$ its induced function.

Given a periodic orbit $g$ of the geodesic flow on $\Ga\bs \G X$, if
$(e_1,\dots, e_k)$ is the sequence of edges followed by $g$, we denote
by $\gls{lebmeasperio}$ the Lebesgue measure along $g$, by
$\gls{periodlength}$ the length of $g$ and by $\gls{periodconduct}$
its {\it period}\index{period} for the system of conductances $c$:
$$ 
\lambda(g)= \sum_{i=1}^k\lambda(e_i)\;\;\;{\rm and}\;\;\;
c(g)=\sum_{i=1}^k\lambda(e_i)\,c(e_i)\;.
$$ 
Let $\periodic(t)$ be the set of periodic orbits of the continuous
time geodesic flow on $\Ga\bs\G X$ with length at most $t$ and let
$\periodic'(t)$ be the subset of primitive ones.

\btheo\label{equidclosedorbitscont} Assume that the critical exponent
$\delta_c$ of $c$ is finite and positive and that the Gibbs measure
$m_{c}$ of $c$ is finite and mixing for the continuous time geodesic
flow.  As $t\to+\infty$, the measures
$$
\delta_c\,e^{-\delta_c\,t}\sum_{g\in\periodic'(t)}e^{c(g)}\,\Lebmeas_g
$$
and 
$$
\delta_c\,t\,e^{-\delta_c\,t}
\sum_{g\in\periodic'(t)}e^{c(g)}\,\frac{\Lebmeas_g}{\lambda(g)}
$$ 
converge to $\frac{m_{c}}{\|m_{c}\|}$ for the weak-star convergence
of measures.  If $\Ga$ is geometrically finite, the convergence holds
for the narrow convergence.  
\etheo

We conjecture that if $\Ga$ is geometrically finite and if its length
spectrum is $4$-Diophantine\footnote{See the definition in Section
  \ref{subsec:mixingratemetgraphs}.}, then for all $n\in\NN$ and
$\beta\in\;]0,1]$, there exist $k\in\NN$ and an error term of the
form $\bigO(t^n\,\|\psi\|_{k,\,\beta})$ for these equidistribution
claims evaluated on any $\psi\in\C^{k,\,\beta}_{\rm c}(\Ga\bs\G X)$. But since
we will not need this result and since the proof is likely to be very
long, we do not address the problem here.

\medskip
\dem Let $\wt F_c$ and $F_c$ be the potentials on $T^1X$ and $\Ga\bs
T^1X$ respectively associated\footnote{See Section \ref{subsec:cond}.}
with $c$, and note that the period of a periodic orbit $g$ for the
geodesic flow on $\Ga\bs\G X$ satisfies\footnote{See Proposition
  \ref{prop:integpotconduct} and Section \ref{subsec:potentials}.}
$$
c(g)= \Lebmeas_g(F_c^\sharp)=\per_{F_c}(\ga)\;,
$$ 
where $F_c^\sharp$ is the composition of the canonical map
$\Ga\bs\G X\ra \Ga\bs T^1X$ with $F_c:\Ga\bs T^1X\ra\RR$, and
$\ga\in\Ga$ is the loxodromic element of $\Ga$ whose conjugacy class
corresponds to $g$.

Let $\H_{\Ga,t}$ be the subset of $\Ga$ that consists in the
loxodromic elements whose translation length is at most $t$, and let
$\H'_{\Ga,t}$ be the subset of $\H_{\Ga,t}$ that consists in the
primitive ones. For every $\ga\in\H_{\Ga,t}$, we denote by $g_\ga$ its
corresponding periodic orbit in $\Ga\bs \G X$. The first claim is
equivalent to the following assertion: we have
\begin{equation}\label{eq:primitiveperiods}
\delta_c\,e^{-\delta_c\,t}\sum_{\ga\in\H'_{\Ga,t}}e^{\per_{F_c}(\ga)}\,\Lebmeas_{g_\ga}
\weakstar \frac{m_{c}}{\|m_{c}\|}
\end{equation}
as $t\to+\infty$.  We proceed with the proof of the convergence
claimed in Equation \eqref{eq:primitiveperiods} as in
\cite[Theo.~9.11]{PauPolSha15}. We first prove that
\begin{equation}\label{eq:allperiods}
\nu''_t=\delta_c\,\|m_{c}\|\,e^{-\delta_c\,t}
\sum_{\ga\in\H_{\Ga,t}}e^{\per_{F_c}(\ga)}\,\Lebmeas_{g_\ga}
\weakstar m_{c}\,.
\end{equation}
We then refer to Step 2 of the proof of \cite[Theo.~9.11]{PauPolSha15}
for the fact that the contribution of the non primitive elements is
negligible, so that Equation \eqref{eq:allperiods} implies Equation
\eqref{eq:primitiveperiods}.  Although the proof of this deduction in
loc.~cit.~is written for manifolds, the arguments are directly
applicable for any $\CAT(-1)$ space $X$ and potential $F$ satisfying
the \ref{eq:HC}-property.\footnote{See Definition
  \ref{defi:HCproperty}.}  In particular, the use of Proposition 5.13
(i) and (ii) of op.~cit.~in the proof of Step 2 in loc.~cit.~is
replaced now by the use of Theorem \ref{theo:HTSR} (1) and (4)
respectively.

Let us fix $x\in X$. Let 
$$
V(x)=\big\{(\xi,\eta)\in (X\cup\partial_\infty X)^2: 
\xi\ne\eta, \ x\in\;]\xi,\eta[\big\}\,,
$$ 
which is an open subset of $X\cup\partial_\infty X$. Note that the
family $(V(y))_{y\in X}$ covers the set of pairs of distinct points of
$\partial_\infty X$. For every $t >0$, let $\nu_t$ be the measure
on $(X\cup\partial_\infty X)^2$ defined by
$$
\nu_t=\delta_c\,\|m_{c}\|\,e^{-\delta_ct}
\sum_{\ga\in\Ga\;:\;d(x,\ga x)\le t}e^{\int_x^{\ga x}\wt F_c}
\Dirac_{\ga^{-1}x}\otimes\Dirac_{\ga x}\,.
$$ 
The measures $\nu_t$ weak-star converge to $\mu^-_x\otimes\mu^+_x$
as $t\to+\infty$ by Corollary \ref{coro:caspointempcont} (taking in
its statement $y=x$).

Let $\ga_\pm$ be the attracting and repelling fixed points of any
loxodromic element $\ga\in\Ga$. Let
$$
\nu'''_t=\delta_c\,\|m_{c}\|\,e^{-\delta_c\,t}
\sum_{\ga\in\H_{\Ga,\,t}}e^{\per_{F_c}(\ga)}
\Dirac_{\ga_-}\otimes\Dirac_{\ga_+}\,.
$$ 
Since $\XX$ is an $\RR$-tree, every element $\ga\in\Ga$ such that
$x\in\;]\ga^{-1}x,\ga x[$ is loxodromic, and for such $x$ and $\ga$ we
have
$$
d(x,\ga x)= \len(\ga)\;\;\;{\rm and}\;\;\;\int_x^{\ga x}\wt F_c=
\per_{F_c}(\ga)\;.
$$  
If furthermore $d(x,\ga x)$ is large, then $\ga^{-1}x$ and $\ga x$ are
respectively close to $\ga_-$ and $\ga_+$ in $X\cup\partial_\infty X$.

Hence, for every continuous map $\psi:(X\cup\partial_\infty X)^2 \ra
[0,+\infty[$ with compact support contained in $V(x)$, and for every
$\epsilon>0$, if $t$ is large enough, we have
$$
e^{-\epsilon}\nu_t(\psi)\leq \nu'''_t(\psi)\leq e^{\epsilon}\nu_t(\psi)\;.
$$
Using Hopf's parametrisation with basepoint $x$, and by Equation
\eqref{eq:changemoinsplus}, the measures $\wt m_c$ and
$\mu^-_x\otimes\mu^+_x\otimes ds$ are equal on $V(x)\times \RR$. Hence
the measures $\nu_t\otimes\, ds$ weak-star converge to $\wt m_c$ on
$V(x)\times \RR$ as $t\to+\infty$.  For every continuous $\psi':\G
X\to\RR$ which is a product of continuous functions in each variable
with compact support in $V(x)\times\RR$ and for every $\epsilon>0$, if
$t$ is large enough, we hence have
$$
e^{-\epsilon}\wt m_c(\psi')\leq \nu'''_t\otimes ds(\psi')\leq 
e^{\epsilon}\wt m_c(\psi')\;.
$$ 
Note that the support of any continuous function with compact
support on $\G X$ may be covered by finitely many open sets
$V(x)\times\RR$ where $x\in X$.  The induced measure\footnote{See the
  beginning of Chapter \ref{sec:equidcountdownstairs}.} of the
$\Ga$-invariant measure $\nu'''_t\otimes\, ds$ on $\G X$ is the
measure $\nu''_t$ on $\Ga\bs\G X$. Since $m_c$ is the induced mesasure
of $\wt m_c$, this proves Equation \eqref{eq:allperiods}, hence gives the
first claim of Theorem \ref{equidclosedorbitscont}.

\medskip
The second claim follows from the first one in the same
way as in \cite[Theo.~9.11]{PauPolSha15}: Consider the measures
$$
m'_t=\delta_c e^{-\delta_c t}
\sum_{g\in\periodic'(t)} \;e^{\L_g(F)} \;\L_{g}\;\;\;{\rm and}\;\;\; 
m''_t=\delta_c t e^{-\delta_c t} \sum_{g\in\periodic'(t)} \;e^{\L_g(F)}\;
\frac{\L_{g}}{\lambda(g)}
$$ 
on $\Ga\bs\G X$.  Fix a continuous map $\psi:\Ga\bs \G X\ra
\mathopen{[}0, +\infty \mathclose{[}$ with compact support. For every
$\epsilon>0$, for every $t>0$, we have, since $\lambda(g)\geq
e^{-\epsilon}t$ for all $g$ in the second sum below,
\begin{align*}
  m''_t(\psi) & \geq m'_t(\psi)\geq \delta_c e^{-\delta_c t}
  \sum_{g\in\periodic'(t)-\periodic'(e^{-\epsilon}t)}
  \;e^{\L_g(F)} \;\L_{g}(\psi)\\
  & \geq e^{-\epsilon}\delta_c t e^{-\delta_c t}
  \sum_{g\in\periodic'(t)-\periodic'(e^{-\epsilon}t)} \;e^{\L_g(F)}
  \;\frac{\L_{g}(\psi)}{\lambda(g)}\\
&= e^{-\epsilon} m''_t(\psi) - e^{-\epsilon}\delta_c t \;e^{-\delta_c t}
  \sum_{g\in\periodic'(e^{-\epsilon}t)} \;e^{\L_g(F)}
  \;\frac{\L_{g}(\psi)}{\lambda(g)}\;.
\end{align*}
By the local finiteness of $X$, the closed orbits meeting the support
of $\psi$ have a positive lower bound on their lengths. Thus, by the
first claim of Theorem \ref{equidclosedorbitscont}, there exists a
constant $C>0$ such that the second term of the above difference is at
most $C\,t\,e^{-\delta_c t} e^{\delta_c \,e^{-\epsilon}t}$, which tends to
$0$ as $t$ tends to $+\infty$.  Hence by applying twice the first claim
of Theorem \ref{equidclosedorbitscont}, we have
$$
\frac{m_c(\psi)}{\|m_c\|}= \lim_{t\ra+\infty}m'_t(\psi)\leq
\liminf_{t\ra+\infty}m''_t(\psi) \leq \limsup_{t\ra+\infty}
m''_t(\psi)\leq \lim_{t\ra+\infty} e^{\epsilon} m'_t(\psi)=
e^{\epsilon}\frac{m_c(\psi)}{\|m_c\|}\;,
$$ 
and the result follows by letting $\epsilon$ go to $0$.  (and
writing any continuous map $\psi:T^1M\ra \RR$ with compact support
into the sum of its positive and negative parts).

\medskip
In order to prove the last claim of Theorem
\ref{equidclosedorbitscont}, assume that $\Ga$ is geometrically
finite.  The narrow convergence follows as in
\cite[Theo.~9.16]{PauPolSha15}, using the fact that there exists a
compact subset of $\Ga\bs \G X$ meeting every periodic orbit of the
geodesic flow, and replacing Lemma 3.10, Lemma 3.2, Equation (112),
Theorem 8.3 and Corollary 5.15 of op.~cit.~by respectively Lemma
\ref{lem:shadowlemma} (1), the \eqref{eq:HC}-property, Proposition
\ref{prop:continuGibbscocycle} \eqref{eq:cocycleombre}, Theorem
\ref{theo:DOPB} and Corollary \ref{coro:finitudeGibbsdivuniq} (1).
\cqfd

\medskip
In a similar way, replacing in the above proof Corollary
\ref{coro:caspointempcont} of Theorem \ref{theo:mainequidup} by the
similar corollary of Theorem \ref{theo:discretemainequidup} with
$\D^-=(\ga x)_{\ga\in\Ga}$ and $\D^-=(\ga y)_{\ga\in\Ga}$ for any
$x,y\in V\XX$, we get the following analogous result for simplicial
trees. For every $n\in\NN$, let now $\periodic(n)$ be the set of
periodic orbits of the discrete time geodesic flow on $\Ga\bs\G \XX$
with length at most $n$ and let $\periodic'(n)$ be the subset of
primitive ones.

\btheo\label{equidclosedorbitsdisc} Let $\XX$ be a locally finite
simplicial tree without terminal vertices, let $\Ga$ be a
nonelementary discrete subgroup of $\Aut(\XX)$ and let $\wt
c:E\XX\ra\RR$ be a $\Ga$-invariant system of conductances.  Assume
that the critical exponent $\delta_c$ of $c$ is finite and positive
and that the Gibbs measure $m_{c}$ is finite and mixing for the
discrete time geodesic flow.  As $n\to+\infty$, the measures
$$
\frac{e^{\delta_c}-1}{e^{\delta_c}}\,e^{-\delta_cn}
\sum_{g\in\periodic'(n)}e^{c(g)}\,\Lebmeas_g
$$
and 
$$
\frac{e^{\delta_c}-1}{e^{\delta_c}}\,n\,e^{-\delta_cn}
\sum_{g\in\periodic'(n)}e^{c(g)}\,\frac{\Lebmeas_g}{\lambda(g)}
$$ 
converge to $\frac{m_c}{\|m_c\|}$ for the weak-star convergence of
measures.  If $\Ga$ is geometrically finite, the convergence holds for
the narrow convergence.  
\cqfd 
\etheo

\medskip
In the special case when $\Ga\bs X$ is a compact graph and $F=0$, the
following immediate corollary of Theorem \ref{equidclosedorbitscont}
is proved in \cite{Guillope94}, and it follows from the results of
\cite{ParPol90}.\footnote{See the introduction of \cite{Sharp10} for
  comments.}  There are also some works on non-backtracking random
walks with related results. For example, for regular finite graphs,
\cite{LubPhiSar86} and \cite{Friedman91} (see
\cite[Lem.~2.3]{Friedman08}) give an expression of the irreducible
trace which is the number of closed walks of a given length.

\bcoro Let $(\XX,\lambda)$ be a locally finite metric tree without
terminal vertices.  Let $\Ga$ be a geometrically finite discrete
subgroup of $\Aut(\XX,\lambda)$.  Let $c:E\XX\ra\RR$ be a
$\Ga$-invariant system of conductances, with finite and positive
critical exponent $\delta_c$.

\smallskip
\noindent (1) If the Gibbs measure $m_c$ is finite and
mixing for the continuous time geodesic flow, then
$$
\sum_{g\in\periodic'(t)}e^{c(g)}\sim \frac{e^{\delta_c\, t}}{\delta_c\, t}
$$
as $t\to+\infty$.

\smallskip
\noindent (2) If $\lambda= 1$ and if the Gibbs measure $m_c$ is
finite and mixing for the discrete time geodesic flow, then
$$
\sum_{g\in\periodic'(n)}e^{c(g)}\sim \frac{e^{\delta_c}}{e^{\delta_c}-1}\,
\frac{e^{\delta_c\, n}}{n}
$$ 
as $n\to+\infty$. 
\cqfd 
\ecoro

%% file: meridonIII.tex

\part{Arithmetic applications}
\label{sect:arithappli}

\chapter{Fields with discrete valuations}
\label{sec:fieldsandvaluations}

Let $\wh K$ be a non-Archimedean local field.  Basic examples of such
fields are the field of formal Laurent series over a finite field, and
the field of $p$-adic numbers (see Examples \ref{ex:laurent} and
\ref{ex:padic}). In Part \ref{sect:arithappli} of this book, we apply
the geometric equidistribution and counting results for simplicial
trees given in Part \ref{part:equid}, in order to prove arithmetic
equidistribution and counting results in such fields $\wh K$.  The
link between the geometry and the algebra is provided by the
Bruhat-Tits tree of $(\PGL_2, \wh K)$, the construction of which is
recalled in Section \ref{subsec:BruhatTitstrees}.  We will only use
the system of conductances equal to $0$ in this Part
\ref{sect:arithappli}.
%
%

In the present Chapter, before embarking on our arithmetic
applications, we recall basic facts on local fields for the
convenience of the geometer reader. For more details, we refer for
instance to \cite{Serre62, Goss98}.  We refer
to \cite{BroParPau16CRAS} for an announcement of the results of
Part \ref{sect:arithappli}, with a presentation different from the one
in the Introduction.

We will only give results for the algebraic group $\underline G=
\PGL_2$ over $\wh K$ and special discrete subgroups $\Ga$ of
$\PGL_2(\wh K)$, even though the same methods give equidistribution
and counting results when $\underline G$ is any semisimple connected
linear algebraic group over $\wh K$ of $\wh K$-rank $1$ and $\Ga$ any
lattice in $G=\underline G(\wh K)$.

\section{Local fields and valuations}
\label{subsec:valuatedlocalfields}

Let $F$ be a field and let $F^\times=(F-\{0\},\times)$ be its
multiplicative group.  A surjective group morphism $v:F^\times\to\ZZ$
to the additive group $\ZZ$, that satisfies
$$
v(a+b)\ge\min\{v(a),v(b)\}\;
$$ 
for all $a,b\in F^\times$, is a (normalised discrete) {\em
  valuation}\index{valuation} $v$ on $F$.  We make the usual
convention and extend the definition of $v$ to $F$ by setting
$v(0)=+\infty$. Note that $v(a+b)=\min\{v(a),v(b)\}$ if $v(a)\ne
v(b)$.  When $F$ is an extension of a finite field $k$, the valuation
$v$ vanishes on $k^\times$.

The subring
$$
\OOO_v=\{x\in F\;:\;v(x)\ge 0\}
$$
is the {\em valuation ring}\index{valuation!ring} (or local ring) of
$v$ (or of $F$ is $v$ is implicit). 

The maximal ideal
$$
\mathfrak m_v=\{x\in F\;:\;v(x)> 0\}
$$
of $\OOO_v$ is principal and it is generated as an ideal of $\OOO_v$
by any element $\gls{uniformizergene}\in F$ with 
$$
v(\pi_v)=1
$$
which is called a {\em uniformiser}\index{uniformiser} of $F$.

The {\em residual field}\index{residual field} of the valuation $v$ is
$$
k_v=\OOO_v/\mathfrak m_v\;.
$$ 
When $k_v$ is finite, the valuation $v$ defines a (normalised,
non-Archimedean) {\em absolute value}\index{absolute!value}
$|\cdot|_v$ on $F$ by
$$
|x|_v=|k_v|^{-v(x)}\;,
$$ 
with the convention that $|k_v|^{-\infty}=0$.  This absolute value
induces an ultrametric distance on $F$ by
$$
(x,y)\mapsto |x-y|_v\;.
$$ 
Let $F_v$ be the completion of $F$ with respect to this distance.  The
valuation $v$ of $F$ uniquely extends to a (normalised discrete)
valuation on $F_v$, again denoted by $v$.

\bexem\label{ex:laurent} Let $K=\FF_q(Y)$ be the field of rational
functions in one variable $Y$ with coefficients in a finite field
$\FF_q$ of order a positive power $q$ of a positive prime $p$ in
$\ZZ$, let $\FF_q[Y]$ be the ring of polynomials in one
variable $Y$ with coefficients in $\FF_q$, and let
$\gls{valuationatinfinity}: K^\times \ra \ZZ$ be the {\em valuation at
infinity}\index{valuation!at infinity} of $K$, defined on every
$P/Q\in K$ with $P\in\FF_q[Y]$ and $Q\in\FF_q[Y]-\{0\}$  by
$$
v_\infty(P/Q)=\deg Q-\deg P\;.
$$ 
The absolute value associated with $v_\infty$ is
$$
|P/Q|_\infty=q^{\deg P-\deg Q}\;.
$$

The completion of $K$ for $v_\infty$ is the field $K_{v_\infty}=
\FF_q((Y^{-1}))$ of formal Laurent series in one variable $Y^{-1}$
with coefficients in $\FF_q$. The elements $x$ in $\FF_q((Y^{-1}))$
are of the form
$$
x= \sum_{i\in\ZZ} x_i\,Y^{-i}
$$ 
where $x_i\in \FF_q$ for every $i\in\ZZ$, and $x_i=0$ for $i$ small
enough. The valuation at infinity of $\FF_q((Y^{-1}))$ extending the
valuation at infinity of $\FF_q(Y)$ is
$$
v_\infty(x)=\sup\{i\in\ZZ\;:\;\forall\;j<i,\;\;\;x_j=0\}\,,
$$ 
that is, 
$$
v_\infty\big(\sum_{i=i_0}^\infty x_iY^{-i}\big)=i_0
$$ 
if $x_{i_0}\ne 0$.  The valuation ring of $v_\infty$ is the ring
$\OOO_{v_\infty}= \FF_q[[Y^{-1}]]$ of formal power series in one
variable $Y^{-1}$ with coefficients in $\FF_q$. The element
$$
\pi_{v_\infty}=Y^{-1}
$$ 
is a uniformiser of $v_ \infty$, the residual field
$\OOO_{v_\infty}/ \pi_{v_\infty}\OOO_{v_\infty}$ of $v_\infty$ is
$k_{v_\infty}=\FF_q$.  
\eexem

\bexem \label{ex:padic} Given a positive prime $p\in\ZZ$, the field of
$p$-adic numbers $\QQ_p$ is the completion of $\QQ$ with respect to
the absolute value $|\cdot|_p$ of the {\it $p$-adic
  valuation}\index{valuation!adic@$p$-adic} $v_p$ defined by setting
$$
v_p\big(p^n\frac ab\big)=n\,,
$$ 
when $n\in\ZZ$, $a,b\in\ZZ-\{0\}$ are not divisible by $p$.  Then
the valuation ring $\OOO_{v_p}$, denoted by $\ZZ_p$, of $\QQ_p$ is the
closure of $\ZZ$ for the absolute value $|\cdot|_p$, the element
$\pi_{v_p}=p$ is a uniformiser, and the residual field is
$k_{v_p}=\ZZ_p/p\ZZ_p=\FF_p$, a finite field of order $p$.  \eexem

A field endowed with a valuation is a {\em non-Archimedean local
  field}\index{non-Archimedean local field}%
\index{local field!non-Archimedean} if it is complete with respect to
its absolute value and if its residual field is finite.\footnote{There
  are also two Archimedean local fields $\CC$ and $\RR$, see for
  example \cite{Cassels86}.} Its valuation ring is then a compact open
additive subgroup.  Any non-Archimedean local field is isomorphic to a
finite extension of the $p$-adic field $\QQ_p$ for some prime $p$, or
to the field $\FF_q((Y^{-1}))$ of formal Laurent series in one
variable $Y^{-1}$ over $\FF_q$ for some positive power $q$ of a prime
$p$.

These formal Laurent series fields may occur as completions of
numerous (global) fonctions fields, that we now define.  The basic
case is described in Example \ref{ex:laurent} above, and the general
case is detailed in Section \ref{subsec:valuedfields} below. The
geometer reader may skip Section \ref{subsec:valuedfields} and use
only Example \ref{ex:laurent} in the remainder of Part
\ref{sect:arithappli} (using $g=0$ when the constant $g$ occurs).

\section{Global function fields}
\label{subsec:valuedfields}

In this Section, we fix a finite field $\gls{finitefield}$ with $q$
elements, where $q$ is a positive power of a positive prime $p\in\ZZ$,
and we recall the definitions and basic properties of a function field
$K$ over $\FF_q$, its genus $g$, its valuations $v$, its completion
$K_v$ for the associated absolute value $|\cdot|_v$ and the associated
affine function ring $R_v$. See for instance \cite{Goss98,Rosen02} for
the content of this Section.

\medskip 
Let $\gls{functionfield}$ be a (global) {\it function
  field}\index{function field} over $\FF_q$, which can be defined in
two equivalent ways as
\begin{enumerate}
\item the field of rational functions on a geometrically irreducible
  smooth projective curve $\gls{smoothprojcurve}$ over $\FF_q$, or
\item an extension of $\FF_q$ of transcendence degree $1$, in
which $\FF_q$ is algebraically closed. 
\end{enumerate} 

\smallskip\noindent 
There is a bijection between the set of closed points of ${\bf C}$ and
the set of (normalised discrete) valuations of its function field $K$,
the valuation of a given element $f\in K$ being the order of the zero
or the opposite of the order of the pole of $f$ at the given closed
point. We fix such an element $v$ from now on. We denote by
$\gls{genus}$ the genus of the curve ${\bf C}$.

In the basic Example \ref{ex:laurent}, ${\bf C}$ is the projective
line $\PP^1$ over $\FF_q$, which is a curve of genus $g=0$, and the
closed point associated with the valuation at infinity $v_\infty$ is
the point at infinity $[1:0]$.

\medskip
We denote by $\gls{completfunctionfield}$ the completion of $K$ for
$v$, and by 
$$
\gls{valuationring}= \{x\in K_v\;:\;v(x)\geq 0\}
$$ 
the valuation ring of (the unique extension to $K_v$) of $v$.  We
choose a uniformiser $\gls{uniformizer}\in K$ of $v$.
We denote by $\gls{residualfield}=\OOO_v/\pi_v\OOO_v$ the residual
field of $v$, which is a finite field of order 
$$
\gls{orderresidualfield}=|k_v|\;.
$$ 
The field $k_v$ is from now on identified with a fixed lift in
$\OOO_v$ (see for instance \cite[Théo.~1.3]{Colmez05}), and is an
extension of the field of constants $\FF_q$. The degree of this
extension is denoted by $\gls{degreevaluation}$, so that
$$
q_v=q^{\Deg v}\;.
$$
We denote by $\gls{absolutevalue}$ the (normalised) absolute value
associated with $v\,$: for every $x\in K_v$, we have
$$
|x|_v=(q_v)^{-v(x)}=q^{-v(x)\Deg v}\;.
$$
Every element $x\in K_v$ is\footnote{See for instance
\cite[Coro.~1.6]{Colmez05}.} a (converging) Laurent series $x=
\sum_{i\in\ZZ} x_i\, (\pi_v)^{i}$ in the variable $\pi_v$ over $k_v$,
where $x_i\in k_v$ is zero for $i\in\ZZ$ small enough. We then have
\begin{equation}\label{eq:defiabvalv}
|x|_v=(q_v)^{-\sup\{j\in\ZZ\;:\;\forall\;i<j,\;x_i=0\}}\;,
\end{equation}
and $\OOO_v$ consists of the (converging) power series $x=
\sum_{i\in\NN} x_i\, (\pi_v)^{i}$ (where $x_i\in k_v$) in the variable
$\pi_v$ over $k_v$.

\medskip 
We denote by $\gls{affinering}$ the affine algebra of the affine curve
${\bf C}-\{v\}$, consisting of the elements of $K$ whose only poles
are at the closed point $v$ of ${\bf C}$. Its field of fractions is
equal to $K$, hence we will often write elements of $K$ as $x/y$
with $x,y\in R_v$ and $y \neq 0$. In the basic Example 
\ref{ex:laurent}, we have $R_{v_\infty}=\FF_q[Y]$. Note that
\begin{equation}\label{eq:interRvOv}
R_v\cap \OOO_v=\FF_q\;,
\end{equation}
since the only rational functions on ${\bf C}$ whose only poles are at
$v$ and whose valuation at $v$ is nonnegative are the constant ones.
We have (see for instance \cite[II.2 Notation]{Serre83}, \cite[page
63]{Goss98})
\begin{equation}\label{eq:inversRv}
(R_v)^\times=(\FF_q)^\times \;.
\end{equation}

\medskip
The following result is immediate when ${\bf C}=\PP^1$, since then 
$R_v+ \OOO_v=K_v$.

\blemm \label{lem:bost}
The dimension of the quotient vector space $K_v/(R_v+ \OOO_v)$
over $\FF_q$ is equal to the genus $g$ of ${\bf C}$.  
\elemm

\dem (indicated by J.-B.~Bost) We refer for instance to \cite{Serre55} for 
background on sheaf cohomology. We denote in the same way the
valuation $v$ and the corresponding closed point on ${\bf C}$. 

Let $\OOO=K\cap\OOO_v$ be the discrete valuation ring of $v$
restricted to $K$.  Since $K$ is dense in $K_v$ and $\OOO_v$ is open
and contains $0$, we have $K_v=K+\OOO_v$. Therefore the canonical map
$$
K/(R_v+\OOO)\;\ra K_v/(R_v+\OOO_v)
$$
is a linear isomorphism over $\FF_q$. Let us hence prove that
$\dim_{\FF_q} K/(R_v+\OOO) =g$.

In what follows, $\V$ ranges over the affine Zariski-open
neighbourhoods of $v$ in ${\bf C}$, ordered by inclusion. Let
$\OOO_{\bf C}$ be the structural sheaf of ${\bf C}$. Note that by the
definition of $R_v$, since the zeros of elements of $K^\times$ are
isolated and by the relation between valuations of $K$ and closed
points of ${\bf C}$,
$$
R_v=H^0({\bf C}-\{v\},\OOO_{\bf C}),\;\;\;\;
K=\varinjlim_\V \;H^0(\V-\{v\},\OOO_{\bf C}) \;\;\;\;{\rm and}\;\;\;\; 
\OOO=\varinjlim_\V \;H^0(\V,\OOO_{\bf C})\;.
$$ 
Since $\V$ and ${\bf C}-\{v\}$ are affine curves, we have $H^1({\bf C}
-\{v\},\OOO_{\bf C})=H^1(\V,\OOO_{\bf C})=0$.  By the Mayer-Vietoris
exact sequence for the covering $\{{\bf C}-\{v\},\V\}$ of ${\bf C}$,
we hence have an exact sequence
$$
H^0({\bf C},\OOO_{\bf C})\ra H^0({\bf C}-\{v\},\OOO_{\bf C})\times
H^0(\V,\OOO_{\bf C})\ra H^0(\V-\{v\},\OOO_{\bf C}) \ra H^1({\bf C},
\OOO_{\bf C})\;.
$$
Therefore
\begin{align*}
K/(R_v+\OOO)&=\varinjlim_\V \;\;H^0(\V-\{v\},\OOO_{\bf C})/
\big(H^0({\bf C}-\{v\},\OOO_{\bf C})+H^0(\V,\OOO_{\bf C})\big)
\\ & \simeq H^1({\bf C},\OOO_{\bf C})\;.
\end{align*}

Since $\dim_{\FF_q}H^1({\bf C},\OOO_{\bf C})=g$ by one definition of
the genus of ${\bf C}$, the result follows.  \cqfd

\bigskip 
Recall that $R_v$ is a Dedekind ring.\footnote{See for instance
  \cite[II.2 Notation]{Serre83}. We refer for instance to \cite[\S
    1.1]{Narkiewicz04} for background on Dedekind rings.} In
particular, every nonzero ideal (respectively fractional ideal) $I$ of
$R_v$ may be written uniquely as $I=\prod_\ppp \ppp^{v_\ppp(I)}$ where
$\ppp$ ranges over the prime ideals in $R_v$ and $v_\ppp(I)\in\NN$
(respectively $v_\ppp(I)\in\ZZ$), with only finitely many of them
nonzero. By convention $I=R_v$ if $v_\ppp(I)=0$ for all $\ppp$.  For
all $x,y\in R_v$ (respectively $x,y\in K$), we denote by
$$
\langle x,\,y\rangle=x\,R_v+y\,R_v
$$
the ideal (respectively fractional ideal) of $R_v$ generated by $x,y$.
If $I,J$ are nonzero fractional ideals of $R_v$, we have
\begin{equation}\label{eq:interidealdedeking}
I\cap J= \prod_{\ppp} \ppp^{\max\{v_\ppp(I), \,v_\ppp(J)\}}
\;\;\;{\rm and}\;\;\;
I+J= 
\prod_{\ppp} \ppp^{\min\{v_\ppp(I),\, v_\ppp(J)\}}\;.
\end{equation}

The (absolute) {\it norm}\index{norm} of a nonzero ideal
$I=\prod_\ppp \ppp^{v_\ppp(I)}$ of $R_v$ is
$$
\gls{normideal}=[R_v:I]=|R_v/I|=\prod_\ppp q^{v_\ppp(I)\deg\ppp }\;,
$$
where $\deg\ppp$ is the degree of the field $R_v/\ppp$ over $\FF_q$,
so that $N(R_v)=1$. By convention $N(0)=0$. This norm is
multiplicative:
$$
N(IJ)=N(I)N(J)\;,
$$ and the norm of a nonzero fractional ideal $I=\prod_\ppp
\ppp^{v_\ppp(I)}$ of $R_v$ is defined by the same formula
$N(I)=\prod_\ppp q^{v_\ppp(I)\deg\ppp }$. Note that if $(a)$ is the
principal fractional ideal in $R_v$ generated by a nonzero element
$a\in K$, we define $N(a)=N\big((a)\big)$. We have (see for instance
\cite[page 63]{Goss98})
\begin{equation}\label{eq:normprincipideal}
N(a)=|a|_v\;.
\end{equation}

{\em Dedekind's zeta function}\index{Dedekind's zeta
  function}\index{zeta function} of $K$ is (see for instance \cite[\S
  7.8]{Goss98} or \cite[\S 5]{Rosen02})
$$
\gls{zetafunctionfield}=\sum_I \frac{1}{N(I)^s}
$$
if $\Re \;s>1$, where the summation is over the nonzero ideals $I$ of
$R_v$. By for instance \cite[\S 5]{Rosen02}, it has an analytic
continuation on $\CC-\{0,1\}$ with simple poles at $s=0,s=1$. It is
actually a rational function of $q^{-s}$. In particular, if $K=\FF_q(Y)$,
then (see \cite[Theo.~5.9]{Rosen02})
\begin{equation}\label{zetamoinsun}
\zeta_{\FF_q(Y)}(-1)=\frac{1}{(q-1)(q^2-1)}\;.
\end{equation}


\medskip We denote by $\gls{haarKv}$ the Haar measure of the (abelian)
locally compact topological group $(K_v,+)$, normalised so that
$\haar_{K_v}(\OOO_\nu) =1$.\footnote{Other normalisations are useful
when considering Fourier transforms, see for instance Tate's
thesis \cite{Tate67}.}  The Haar measure scales as follows under
multiplication: for all $\lambda,x\in K_v$, we have
\begin{equation}\label{eq:homothetyhaar}
d\haar_{K_v}(\lambda x)=|\lambda|_v\;d\haar_{K_v}(x)\;.
\end{equation}

Note that any nonzero fractional ideal $I$ of $R_v$ is a discrete
subgroup of $(K_v,+)$, and we will again denote by $\haar_{K_v}$ the
Haar measure on the compact group $K_v/I$ which is induced by the
above normalised Haar measure of $K_v$.

\blemm \label{lem:covolideal} For every nonzero fractional ideal $I$ of
$R_v$, we have
$$
\haar_{K_v} (K_v/I)= q^{\,g-1}\;N(I) \;.
$$
\elemm

\dem By the scaling properties of the Haar measure, we may assume that
$I$ is an ideal in $R_v$.  By Lemma \ref{lem:bost}, we have $\card\;
K_v/(R_v+ \OOO_v)=q^g$. By Equation \eqref{eq:interRvOv} and by the
normalisation of the Haar measure, we have 
$$
\haar_{K_v} (R_v+\OOO_v)/R_v=\haar_{K_v} \OOO_v/(R_v\cap \OOO_v)
=\haar_{K_v} \OOO_v/\FF_q=\frac{1}{q}\;.
$$
Hence
$$
\haar_{K_v} (K_v/R_v)= 
\card \big(K_v/(R_v+ \OOO_v)\big)\;
\haar_{K_v} (R_v+ \OOO_v)/R_v = q^{g-1}\;.
$$
Since $\haar_{K_v} (K_v/I)= N(I)\;\haar_{K_v} (K_v/R_v)$, the result
follows. 
\cqfd

\chapter{Bruhat-Tits trees and modular groups}
\label{sec:BruhatTitstrees}

In this Chapter, we give background information and preliminary
results on the main link between the geometry and the algebra used for
our arithmetic applications: the (discrete time) geodesic flow on
quotients of Bruhat-Tits trees by arithmetic lattices. We refer to Section \ref{subsec:trees} 
for the basic notation and terminology for trees and graphs (of groups).

We denote the image in $\PGL_2$ of an element 
$\begin{pmatrix} a&b\\c&d\end{pmatrix}\in\GL_2$ by
$\begin{bmatrix}a&b\\c&d \end{bmatrix}\in\PGL_2$.

\section{Bruhat-Tits trees}
\label{subsec:BruhatTitstrees}

Let $K_v$ be a non-Archimedean local field, with valuation $v$,
valuation ring $\OOO_v$, choice of uniformiser $\pi_v$, and residual
field $k_v$ of order $q_v$ (see Section
\ref{subsec:valuatedlocalfields} for definitions).

In this Section, we recall the construction and basic properties of
the Bruhat-Tits tree $\XX_v$ of $(\PGL_2,K_v)$, see for instance
\cite{Tits79}. We use its description given in \cite{Serre83}, to
which we refer for proofs and further information.

\medskip 
An {\it $\OOO_{v}$-lattice}\index{lattice@$\OOO_{v}$-lattice}
$\Lambda$ in the $K_v$-vector space $K_v\times K_v$ is a rank $2$ free
$\OOO_{v}$-submodule of $K_v\times K_v$, generating $K_v\times K_v$ as
a vector space. The Bruhat-Tits tree $\gls{BruhatTitstree}$ of
$(\PGL_2,K_v)$ is the graph whose set of vertices $V\XX_v$ is the set
of homothety classes (under $(K_v)^\times$) $[\Lambda]$ of
$\OOO_{v}$-lattices $\Lambda$ in $K_v\times K_v$, and whose
nonoriented edges are the pairs $\{x,x'\}$ of vertices such that there
exist representatives $\Lambda$ of $x$ and $\Lambda'$ of $x'$ for
which $\Lambda \subset \Lambda'$ and $\Lambda'/\Lambda$ is isomorphic
to $\OOO_{v}/{\pi_{v}} \OOO_{v}$. We again denote by $\XX_v$ the
geometric realisation of $\XX_v$.\footnote{giving length $1$ to (the
  geometric realisation of) each nonoriented edge, see Section
  \ref{subsec:trees}.} Two (oriented) edges are naturally associated
with each nonoriented edge.  If $K$ is any field endowed with a
valuation $v$ whose completion is $K_v$, then the similarly defined
Bruhat-Tits tree of $(\PGL_2,K)$ coincides with $\XX_v$, see
\cite[p.~71]{Serre83}.

The graph $\XX_v$ is a regular tree of degree $|\PP_1(k_{v})|= q_v+1$.
In particular, the Bruhat-Tits tree of $(\PGL_2,\QQ_p)$ is regular of
degree $p+1$, and if $K_v=\FF_q((Y^{-1}))$ and $v=v_\infty$, then the
Bruhat-Tits tree $\XX_v$ of $(\PGL_2,K_v)$ is regular of degree
$q+1$. More generally, if $K_v$ is the completion of a function field
over $\FF_q$ endowed with a valuation $v$ as in Section
\ref{subsec:valuedfields}, then the Bruhat-Tits tree of $(\PGL_2,K_v)$
is regular of degree $q_v+1=q^{\deg v}+1$.

The {\it standard base point}\index{standard base point}
$\gls{basepointBTt}$ of $\XX$ is the homothety class
$[\OOO_v\times\OOO_v]$ of the $\OOO_{v}$-lattice $\OOO_{v} \times
\OOO_{v}$, generated by the canonical basis of $K_v \times K_v$.  In
particular, we have
\begin{equation}\label{eq:distooov}
d(*_v,[\OOO_v\times x\OOO_v]) = |v(x)|
\end{equation}
for every $x\in (K_v)^\times$. The {\it
  link}\index{link}
$$
\gls{link}(*_v)=\{y\in V\XX_v\;:\;d(y,*_v)=1\}
$$ 
of $*_v$ in $\XX_v$ identifies with the projective line $\PP_1(k_v)$.

The left linear action of $\GL_2(K_v)$ on $K_v\times K_v$ induces a
faithful, vertex-transitive left action by automorphisms of $\PGL_2
(K_v)$ on the Bruhat-Tits tree $\XX_v$. The stabiliser in
$\PGL_2(K_v)$ of $*_v$ is $\PGL_2 (\OOO_{v})$, acting projectively on
$\lk(*_v)=\PP_1(k_v)$ by reduction modulo $\pi_v\OOO_v$ of the
coef\-ficients. We will hence identify $\PGL_2(K_v)/ \PGL_2(\OOO_{v})$
with $V\XX_v$ by the map $g\,\PGL_2(\OOO_{v})\mapsto g\,*_v$.

We identify the projective line $\PP_1(K_v)$ with $K_v \cup
\{\infty\}$ using the map $K_v(x,y) \mapsto \frac{x}{y}$, so that
$$
\gls{pointprojectifinfty}=[1:0]\;.
$$ 
The projective action of $\GL_2(K_v)$ or $\PGL_2(K_v)$ on
$\PP^1(K_v)$ is the action by {\it homographies}%
\index{homography}\footnote{or linear fractional transformations} on
$K_v\cup\{\infty\}$, given by $(g,z)\mapsto g\cdot z
=\frac{a\,z+b}{c\,z+d}$ if $g=\begin{pmatrix} a & b \\ c &
d\end{pmatrix}\in \GL_2(K_v)$, or $g=\begin{bmatrix} a & b \\ c &
d\end{bmatrix}\in\PGL_2(K_v)$. As usual we define
$\infty\mapsto\frac ac$ and $-\frac dc\mapsto\infty$.

There exists a unique homeomorphism between the boundary at infinity
$\partial_\infty \XX_v$ of $\XX_v$ and $\PP_1(K_v)$ such that the
(continuous) extension to $\partial_\infty \XX_v$ of the isometric
action of $\PGL_2(K_v)$ on $\XX_v$ corresponds to the projective
action of $\PGL_2(K_v)$ on $\PP_1(K_v)$. From now on, we identify
$\partial_\infty \XX_v$ and $\PP_1(K_v)$ by this homeomorphism.  Under
this identification, $\OOO_v$ consists of the positive endpoints
$\ell_+$ of the geodesic lines $\ell$ of $\XX_v$ with negative
endpoint $\ell_-=\infty$ that pass through the vertex $*_v$ (see the
picture below).

 \begin{center}
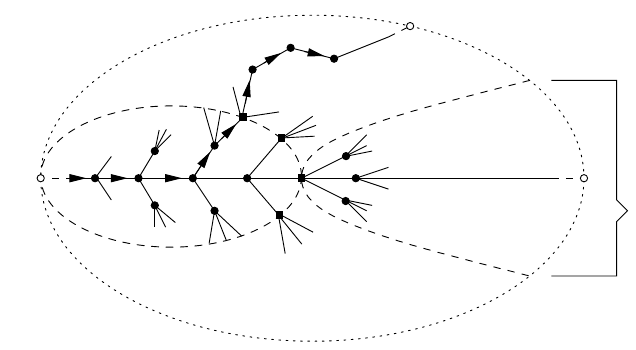
\end{center}

Let $\H_\infty$ be the horoball centred at $\infty\in\partial_\infty
\XX_v$ whose associated horosphere passes through $*_v$.  
There is a unique labeling of the edges of $\XX_v$ by elements of
$\PP_1(k_v)=k_v\cup\{\infty\}$ such that (see the above picture)
\begin{itemize}
\item the label of any edge of $\XX_v$ pointing towards
  $\infty\in\partial_\infty \XX_v$ is $\infty$,
\item for any $x=\sum_{i\in\ZZ} x_i\,(\pi_v)^{i} \in K_v$, the
  sequence $(x_i)_{i\in\ZZ}$ is the sequence of the labels of the
  (directed) edges that make up the geodesic line $]\infty,x[$
      oriented from $\infty$ towards $x$,
\item $x_0$ is the label of the edge of $]\infty,x[$ exiting the
      horoball $\H_\infty$.
\end{itemize}
We refer to \cite[Sect.~5]{Paulin02} for a detailed treatment of the
case $K_v=\FF_q((Y^{-1}))$ and $v=v_\infty$.

\bigskip 
For all $\eta,\eta'\in K_v=\partial_\infty \XX_v -\{\infty\}$, we have
\begin{equation}\label{eq:valuationham}
 |\eta-\eta'|_v=d_{\H_\infty}(\eta,\eta')^{\ln q_v}
\end{equation}
by the definitions of the absolute value $|\cdot|_v$ and of
Hamenst\"adt's distance, see Equation \eqref{eq:defiabvalv}, the above
geometric interpretation, and Equation \eqref{eq:defidisthamenbord}.
Note that in \cite{Paulin02}, Hamenst\"adt's distance in a regular
tree is defined in a different way: In that reference, the distance
$|\eta-\eta'|_v$ equals Hamenst\"adt's distance between $\eta$ and
$\eta'$.

In particular, the Hölder norms\footnote{See the definition in Section
  \ref{subsec:holdercont}.} $\|\psi\|_{\beta,\,|\,\cdot\,-\,\cdot\,|_v}$ and
$\|\psi\|_{\beta',\,d_{\H_\infty}}$ of a function $\psi:K_v\ra\RR$,
respectively for the distances $(x,y)\mapsto |x-y|_v$ and
$d_{\H_\infty}$ on $K_v$, are related by the following formula:
\begin{equation}\label{eq:relatholdernormhamabsval}
\forall\;\beta\in\;\big]0,\frac{1}{\ln q_v}\big],\;\;\;\;
 \|\psi\|_{\beta,\,|\,\cdot\,-\,\cdot\,|_v}=\|\psi\|_{\beta\ln q_v,\,d_{\H_\infty}}\;.
\end{equation}

\bigskip 
The group $\PGL_2(K_v)$ acts simply transitively on the set of ordered
triples of distinct points in $\partial_\infty \XX_v= \PP_1(K_v)$. In
particular, it acts transitively on the space $\G\XX_v$ of (discrete)
geodesic lines in $\XX_v$. The stabiliser under this action of the
geodesic line (from $\infty=[1:0]$ to $0=[0:1]$)
$$
\gls{ellstar}:n\mapsto [\OOO_v\times(\pi_v)^{-n}\OOO_v]
$$ 
is the maximal compact-open subgroup 
$$
\gls{maxcompactopendiag}
=\Bigg\{\begin{bmatrix}a&0\\0&d\end{bmatrix}\;:\;
a,d\in (\OOO_v)^\times\Bigg\}
$$ 
of the diagonal group
$$
\uA(K_v)=\Bigg\{\begin{bmatrix}a&0\\0&d\end{bmatrix}\;:\;
a,d\in (K_v)^\times\Bigg\}\,.
$$
We will hence identify $\PGL_2(K_v)/ \uA(\OOO_v)$ with $\G\XX_v$ by
the mapping $\wt\Xi: g\uA(\OOO_v)\mapsto g\,\ell^*$. Define
$$
a_v=\begin{bmatrix}1&0\\0&\pi_v^{-1}\end{bmatrix},
$$ 
which belongs to $\uA(K_v)$ and centralises $\uA(\OOO_v)$. The
homeomorphism $\wt\Xi$ is equivariant for the actions on the left of
$\PGL_2(K_v)$ on $\PGL_2(K_v)/ \uA(\OOO_v)$ and $\G\XX_v$. It is also
equivariant for the actions on $\PGL_2(K_v)/ \uA(\OOO_v)$ under
translations on the right by $(a_v)^\ZZ$ and on $\G\XX_v$ under the
discrete geodesic flow $(\flow{n})_{n\in\ZZ}$: for all $n\in\ZZ$ and
$x\in \PGL_2(K_v)/ \uA(\OOO_v)$, we have
\begin{equation}\label{eq:commutgeodfowrightasubv}
\wt\Xi\,(x\;{a_v}^n)=\flow{n}\,\wt\Xi\,(x)\;.
\end{equation}

Furthermore, the stabiliser in $\PGL_2(K_v)$ of the ordered pair of
endpoints $(\ell^*_-=\infty, \ell^*_+=0)$ of $\ell^*$ in
$\partial_\infty \XX_v =\PP_1(K_v)$ is $\uA(K_v)$.  Therefore any
element $\ga\in \PGL_2(K_v)$ which is loxodromic on $\XX_v$ is
diagonalisable over $K_v$. By \cite[page 108]{Serre83}, the
translation length on $\XX_v$ of $\ga_0=\begin{bmatrix}a & 0 \\ 0 &
  d \end{bmatrix}$ is
\begin{equation}
\lambda(\ga_0)= |v(a)-v(d)|\;.
\end{equation}
Note that if $\wt{\ga_0}=\begin{pmatrix}a & 0 \\ 0 & d \end{pmatrix}
\in\GL_2(K_v)$ is a representative of $\ga_0$ such that $\det
\wt{\ga_0}\in(\OOO_v)^\times$, then $0=v(\det\wt{\ga_0})=v(ad)
=v(a)+v(d)$, so that $v(d)=-v(a)$ and $\lambda(\ga_0)= 2|v(a)|$. Since
$v(a)\neq v(d)$ if $\lambda(\ga_0)\neq 0$, we have
$v(\operatorname{tr}\wt{\ga_0})=v(a+d)=\min\{v(a),v(d)\}=-|v(a)|$. Thus,
\begin{equation}\label{eq:translatlength}
\lambda(\ga_0)= 2|v(\operatorname{tr}\wt{\ga_0})|\;.
\end{equation}
By conjugation, this formula is valid if $\ga_0\in\PGL_2(K_v)$ is
loxodromic on $\XX_v$ and represented by $\wt{\ga_0}\in\GL_2(K_v)$
such that $\det\wt{\ga_0}\in(\OOO_v)^\times$.

\bigskip
Let $\H$ be a horoball in $\XX_v$ whose boundary is contained in
$V\XX_v$ and whose point at infinity $\xi$ is different from $\infty$.
With $\beta_\infty:V\XX_v\times V\XX_v\ra\ZZ$ the Busemann function at
$\infty$ (see Equation \eqref{eq:buscocycastree}), the {\it
  height}\index{heigth} of $\H$ is
$$
\gls{height}(\H)=
\max\{\beta_\infty(x,*_v)\;:\;x\in\partial\H\} \in\ZZ\;,
$$ 
which is the signed distance between $\H_\infty$ and
$\H$.\footnote{See the definition of the signed distance just above
  Lemma \ref{lem:horodisjoint}.} It is attained at the intersection
point with $\partial\H$ of the geodesic line from $\infty$ to $\xi$,
which is then called the {\it highest point}\index{highest point} of
$\H$. Note that the height of $\H$ is invariant under the action of
the stabiliser of $\H_\infty$ in $\PGL_2(K_v)$ on the set of such
horoballs $\H$.

The following lemma is a generalisation of \cite[Prop.~6.1]{Paulin02}
that only covered the particular case of $K=\FF_q(Y)$ and $v=v_\infty$.

\blemm \label{lem:imagehoroinftyparpglderv} 
Assume that $K_v$ is the completion of a function field $K$ over
$\FF_q$ endowed with a valuation $v$, with associated affine function
ring $R_v$.  For every $\ga=
\begin{bmatrix}a&b\\c&d\end{bmatrix}\in\PGL_2(K)$ with $a,b,c,d\in
K$ such that $ad-bc\in (\OOO_v)^\times$ and $c\neq 0$, the image of
$\H_\infty$ by $\ga$ is the horoball centred at $\frac ac \in
K \subset K_v = \partial_\infty \XX_v - \{\infty\}$ with height
$$
\height_\infty(\ga\H_\infty)= -2\;v(c)\;.
$$
\elemm

\dem It is immediate that $\ga\infty=\frac ac$ under the projective
action.  Up to multiplying $\ga$ on the left by $\begin{bmatrix} 1& -
\frac{a}{c} \\0&1\end{bmatrix}\in\PGL_2(K)$, which does not change $c$
nor the height of $\ga\H_\infty$, we may assume that $a=0$ and that
$b$ has the form $c^{-1}u$ with $u=bc-ad\in(\OOO_v)^\times$.
Multiplying $\ga$ on the right by $\begin{bmatrix} 1& - \frac{d}{c}
\\0&1\end{bmatrix}\in\PGL_2(K)$ preserves $\ga\H_\infty$ and does not
change $a=0$, $b=c^{-1}u$ or $c$. Hence we may assume that
$d=0$. Since $\ga$ then exchanges the points $\infty$ and $0$ in
$\partial_\infty \XX_v$, the highest point of $\ga\H_\infty$ is $\ga
*_v$.  Assuming first that $0,\ga*_v,*_v,\infty$ are in this order on
the geodesic line from $0$ to $\infty$, we have by
Equation \eqref{eq:distooov}
\begin{align*}
\height_\infty(\ga\H_\infty)&= d(*_v,\ga *_v)=
d([\OOO_v\times\OOO_v],[c^{-1}u\OOO_v\times c\OOO_v])\\& =
d([\OOO_v\times\OOO_v],[\OOO_v\times c^2\OOO_v]) =-v(c^2)=-2\;v(c)\;.
\end{align*}
If $0,\ga*_v,*_v,\infty$ are in the opposite order, then the same
computation holds, up to replacing the distance $d$ by its opposite
$-d$.  
\cqfd

\section{Modular graphs of groups}
\label{subsec:modulargroup}

Let $K$ be a function field over $\FF_q$, let $v$ be a (normalised
discrete) valuation of $K$, let $K_v$ be the completion of $K$
associated with $v$, and let $R_v$ be the affine function ring
associated with $v$ (see Section \ref{subsec:valuedfields} for
definitions).

The group $\gls{gammav}=\PGL_2(R_v)$ is a lattice in the locally compact
group $\PGL_2(K_v)$, and a lattice\footnote{See
Section \ref{subsec:trees} for a definition.} of the Bruhat-Tits tree
$\XX_v$ of $(\PGL_2, K_v)$, called the {\it modular
group}\index{modular!group} at $v$ of $K$. The quotient graph
$\Ga_v\bs\XX_v$ is called the {\it modular graph}\index{modular!graph}
at $v$ of $K$, and the quotient graph of groups\footnote{See
Section \ref{subsec:trees} for a definition.} $\Ga_v\dbs\XX_v$ is
called the {\it modular graph of groups}\index{modular!graph of
groups} at $v$ of $K$.  We refer to \cite{Serre83} for background
information on these objects, and for instance to \cite{Paulin02} for
a geometric treatment when $K=\FF_q(Y)$ and $v= v_\infty$.


By for instance \cite{Serre83}, the {\it set of cusps}\index{cusp}
$\Ga_v\bs\PP_1(K)$ is finite, and $\Ga_v\bs \XX_v$ is the disjoint
union of a finite connected subgraph containing $\Ga_v\,*_v$ and of
maximal open geodesic rays $h_z (\,]0,+\infty[)$, for $z=\Ga_v\wt z\in
\Ga_v\bs \PP_1(K)$, where $h_z$ (called a cuspidal ray, see Section
\ref{subsec:trees}) is the injective mage by the canonical projection
$\XX_v\ra \Ga_v \bs \XX_v$ of a geodesic ray whose point at infinity
in $\PP_1(K) \subset \partial_\infty \XX_v$ is equal to $\wt z$.
Conversely, any geodesic ray whose point at infinity lies in $\PP_1(K)
\subset \partial_\infty \XX_v$ contains a subray that maps injectively
by the canonical projection $\XX_v\ra \Ga_v \bs \XX_v$.

The group $\Ga_v=\PGL_2(R_v)$ is a geometrically finite
lattice\footnote{See Section \ref{subsec:trees} for a definition and
  for instance \cite{BasLub01} for a profusion of geometrically
  infinite lattices in simplicial trees.} by for instance
\cite{Paulin04b}.  The set of bounded parabolic fixed points of
$\Ga_v$ is exactly $\PP_1(K)\subset \partial_\infty \XX_v$, and the
set of conical limit points of $\Ga_v$ is $\PP_1(K_v)-\PP_1(K)$.

Let us denote by $\widehat{\Ga_v\bs\XX_v}= (\Ga_v\bs \XX_v)\sqcup
\E_v$ Freudenthal's compactification of $\Ga_v \bs\XX_v$ by its finite
set of ends $\E_v$, see \cite{Freudenthal31}. This set of ends is
indeed finite, in bijection with $\Ga_v\bs\PP_1(K)$ by the map which
associates to $z\in \Ga_v\bs\PP_1(K)$ the end towards which the
cuspidal ray $h_z$ converges. See for instance \cite{Serre83} for a
geometric interpretation of $\E_v$ in terms of the curve ${\bf C}$.

Let $\gls{fracidealclassset}$ be the set of classes of fractional
ideals of $R_v$.  The map which associates to an element
$[x:y]\in\PP_1(K)$ the class of the fractional ideal $xR_v+yR_v$
generated by $x,y$ induces a bijection from the set of cusps
$\Ga\bs\PP_1(K)$ to $\I_v$.

The volume\footnote{See Section \ref{subsec:trees} for a definition.}
of the modular graph of groups $\Ga_v\dbs\XX_v$ can be computed using
Equation \eqref{eq:inversRv} and Exercice 2 b) in
\cite[II.2.3]{Serre83}:
\begin{equation}\label{eq:covolPGLdeR}
\Vol(\PGL_2(R_v)\dbs\XX_v)=(q-1)\Vol(\GL_2(R_v)\dbs\XX_v)=
2\;\zeta_K(-1)\;.
\end{equation}

\medskip If $K=\FF_q(Y)$ is the rational function field over $\FF_q$
and if we consider the valuation at infinity $v=v_\infty$ of $K$, then
the {\em Nagao lattice}\index{Nagao lattice}\footnote{This lattice was
studied by Nagao in \cite{Nagao59}, see also \cite{Mozes99,BasLub01}.
It is called the modular group in \cite{Weil70}.}
$\Ga_v= \PGL_2(\FF_q[Y])$ acts transitively on $\PP^1(K)$. Its
quotient graph of groups $\Ga_v\dbs\XX_v$ is the following {\it
modular ray}\index{modular!ray} (with associated edge-indexed graph)

\begin{center}
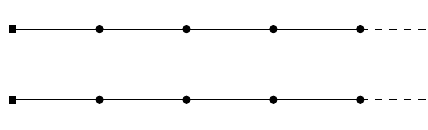
\end{center}

\noindent 
where $\Ga_{-1}=\PGL_2(\FF_q)$, $\Ga'_0=\Ga_0\cap\Ga_{-1}$ and, for
every $n\in\NN$,
$$
\Ga_n=\bigg\{\begin{bmatrix}a&b\\0&d\end{bmatrix}
\in\PGL_2(\FF_q[Y])
\;:\; a,d\in \FF_q^\times, b\in \FF_q[Y], \deg b\leq n+1\bigg\}\;.
$$

Note that even though $\PGL_2(K_v)$ has inversions on $\XX_v$, its
subgroup $\Ga_v=\PGL_2(R_v)$ acts without inversion on $\XX_v$ (see
for instance \cite[II.1.3]{Serre83}). In particular, the quotient
graph $\Ga_v\bs\XX_v$ is then well defined.

\section{Computations of measures for Bruhat-Tits trees}
\label{subsec:volcomputBruhatTitstrees}

In this Section, we compute explicit expressions for the skinning
measures\footnote{See Section \ref{subsec:skinningmeasures}.} of
horoballs and geodesic lines, and for the Bowen-Margulis
measures,\footnote{See Section \ref{subsec:gibbsmeasure}.}  when
considering lattices of Bruhat-Tits trees.  See \cite[Section
  7]{ParPau16LMS} and \cite[Section 4]{ParPau16MA} for analogous
computations in the real and complex hyperbolic spaces respectively,
and \cite{BroPau07JLMS} for related computations in the tree case.

Let $(K_v,v)$ be as in the beginning of Section
\ref{subsec:BruhatTitstrees}.  Let $\Ga$ be a lattice of the
Bruhat-Tits tree $\XX_v$ of $(\PGL_2,K_v)$.  Since $\XX_v$ is regular
of degree $q_v+1$, the critical exponent of $\Ga$ is
\begin{equation}\label{eq:expocritlatticeBT}
  \delta_\Ga=\ln q_v 
\end{equation}
by Proposition \ref{prop:uniflatmBMfinie} and Equation 
\eqref{eq:lyonshdim}.

We normalise the Patterson density $(\mu_x)_{x\in V\XX_v}$ of $\Ga$ as
follows.\footnote{Since the system of conductances is $0$, note that
  $\mu^+_x=\mu^-_x=\mu_x$ is the standard Patterson measure of $\Ga$.}
Let $\H_\infty$ be the horoball in $\XX_v$ centred at $\infty$ whose
associated horosphere passes through $*_v$. Let $t\mapsto x_t$ be the
geodesic ray in $\XX_v$ such that $x_0=*_v$ and which converges to
$\infty$.

\begin{center}
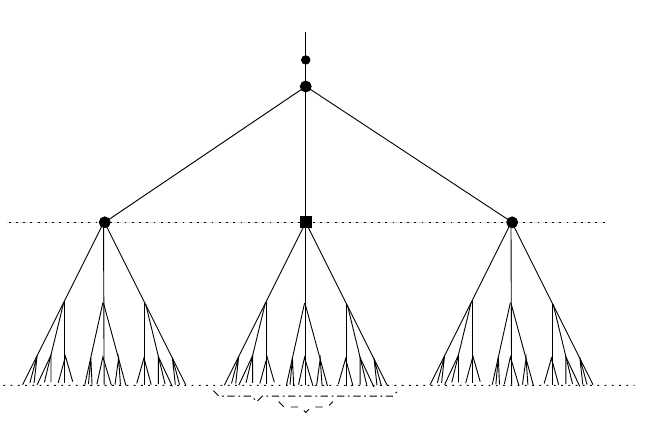
\end{center}

Hamenst\"adt's measure\footnote{See
Equation \eqref{eq:defiHamensmeasure}.} associated with $\H_\infty$
$$
\mu_{\H_\infty}=\lim_{t\to+\infty}e^{\delta_\Ga t}\mu_{x_t}=
\lim_{t\to+\infty} {q_v}^t\mu_{x_t}\;
$$ 
is a Radon measure on $\partial_\infty \XX_v-\{\infty\}=K_v$, invariant
under all isometries of $\XX_v$ preserving $\H_\infty$, since $\Ga$ is
a lattice. Hence it is invariant under the translations by the
elements of $K_v$. By the uniqueness property of Haar measures,
$\mu_{\H_\infty}$ is a constant multiple of the chosen Haar
measure\footnote{Recall that we normalise the Haar measure of
$(K_v,+)$ such that $\haar_{K_v}(\OOO_v) =1$.} of $K_v$, and we
normalise the Patterson density $(\mu_x)_{x\in V\XX_v}$ so that
\begin{equation}\label{eq:hamenhaar}
\mu_{\H_\infty}=\haar_{K_v}\;.
\end{equation}

We summarise the various measure computations in the following result.

\bprop\label{prop:mescomputBT} Let $\Ga$ be a lattice of the
Bruhat-Tits tree $\XX_v$ of $(\PGL_2,K_v)$, with Patterson density
normalised as above.
\begin{enumerate}
\item The outer/inner skinning measures of the singleton $\{*_v\}$ are
  given by
  $$
  d\,\wt\sigma^\pm_{\{*_v\}}(\rho)= d\mu_{*_v}(\rho_\pm)=
  \big(\max\{1,|\rho_\pm|_v\}\big)^{-2}\;d\haar_{K_v}(\rho_\pm)\;
  $$
  on the set of $\rho\in\normalpm\{*_v\}$ such that
  $\rho_\pm\neq \infty$.
\item The total mass of the Patterson density is
  $$
  \|\mu_{x}\|=\frac{q_v+1}{q_v}
  $$ 
  for all $x\in V\XX_v$.
\item The skinning measure of the horoball $\H_\infty$ is the
  projection of the Haar measure of $K_v$: For all
  $\rho\in\normalpm\H_\infty$, we have
$$
d\wt\sigma^\pm_{\H_\infty}(\rho)= d\mu_{\H_\infty}(\rho_\pm)=
d\haar_{K_v}(\rho_\pm)\;.
$$
\item If $\infty$ is a bounded parabolic fixed point of $\Ga$, with
  $\Ga_\infty$ its stabiliser in $\Ga$, if $\D=
  (\ga\H_\infty)_{\ga\in\Ga/\Ga_\infty}$, we have
$$
\|\sigma^\pm_{\D}\|=\haar_{K_v}(\Ga_\infty\bs K_v)=
\Vol(\Ga_\infty\dbs\partial\H_\infty)\;.
$$
\item 
  Let $L$ be a geodesic line in $\XX_v$ with endpoints $L_\pm\in
  K_v= \partial_\infty \XX_v-\{\infty\}$.  Then on the set of
  $\rho\in \normalout L$ such that $\rho_+\in
  K_v= \partial_\infty \XX_v - \{\infty\}$ and $\rho_+\neq L_\pm$, the
  outer skinning measure of $L$ is
$$
d\,\wt\sigma^+_L(\rho)=
\frac{|L_+-L_-|_v}{|\rho_+-L_-|_v\,|\rho_+-L_+|_v}
\;d\haar_{K_v}(\rho_+)\,.
$$
\item Let $L$ be a geodesic line in $\XX_v$, let $\Ga_L$ be the
  stabiliser in $\Ga$ of $L$, and assume that $\Ga_L\bs L$ has finite
  length. Then with $\D=(\ga L)_{\ga \in\Ga/\Ga_L}$, we have
$$
\|\sigma^\pm_\D\|=\frac{q_v-1}{q_v}\;\Vol(\Ga_L\dbs L)\,.
$$
\end{enumerate}
\eprop

\dem (1) For every $\xi\in K_v$, by the description of the geodesic
lines in the Bruhat-Tits tree $\XX_v$ starting from $\infty$ given in
Section \ref{subsec:BruhatTitstrees}, we have $\xi\in \OOO_v$ if and
only if $P_{\H_\infty}(\xi)=*_v$.\footnote{Recall that
  $P_{\H_\infty}:\partial_v\XX_v-\{\infty\}\ra \partial\H_\infty$ is
  the closest point map to the horoball $\H_\infty$, see Section
  \ref{subsect:nbhd}.}

\begin{center}
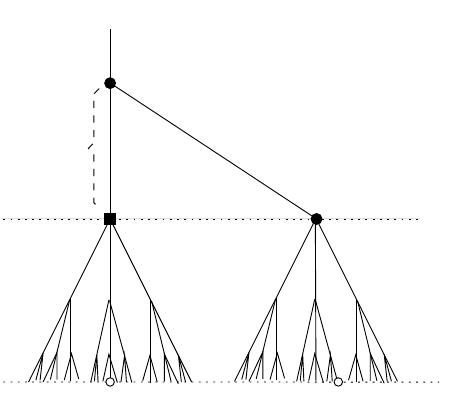
\end{center}
For every $\xi\in K_v-\OOO_v$, by Equations
\eqref{eq:valuationham} and \eqref{eq:defidisthamenbord}, we have
\begin{equation}\label{eq:calcvaldishamen}
|\xi|_v = d_{\H_\infty}(0,\xi)^{\ln q_v}= q_v^{\;\,\frac{1}{2}
  \,d(*_v,\,P_{\H_\infty}(\xi))}\;.
\end{equation}

On the set of geodesic rays $\rho\in\normalpm\{*_v\}$ such that
$\rho_\pm\neq \infty$, by Equation \eqref{eq:skinnpoint}, by the last
claim of Proposition \ref{prop:HPun}, by
Equation \eqref{eq:expocritlatticeBT},\footnote{Since the potential is
zero, the Gibbs cocycle is the critical exponent times the Busemann
cocycle.} since $P_{\H_\infty}(\rho_{\pm})$ belongs to the geodesic
ray $[*_v,\rho_\pm[$ (even when $\rho_\pm\in\OOO_v$), and by Equation
\eqref{eq:hamenhaar}, we have
\begin{align*}
  d\wt\sigma^\pm_{\{*_v\}}(\rho)&= \;d\mu_{*_v}(\rho_{\pm})
  =e^{C_{\rho_{\pm}}(P_{\H_\infty}(\rho_{\pm}),\,*_v)}\;\;
  d\mu_{\H_\infty}(\rho_{\pm})\\ &
  =e^{\delta_\Ga\,\beta_{\rho_{\pm}}(P_{\H_\infty}(\rho_{\pm}),\,*_v)}\;\;
  d\mu_{\H_\infty}(\rho_{\pm})\\ & =
q_v^{-d(P_{\H_\infty}(\rho_{\pm}),\,*_v)}\;\;
  d\haar_{K_v}(\rho_{\pm}) \;.
\end{align*}
Therefore, if $\rho\in\normalpm\{*_v\}$ is such that
$$
\rho_\pm\in\OOO_v=\{\xi\in K_v\;:\;|\xi|_v\leq 1\}
=\{\xi\in K_v\;:\;P_{\H_\infty}(\xi)=*_v\}\;,
$$ 
then $d\,\wt\sigma^\pm_{\{*_v\}}(\rho)= d\haar_{K_v}(\rho_\pm)\;$.  If
$\rho_\pm\in K_v-\OOO_v$, Equation \eqref{eq:calcvaldishamen} gives
the claim.

\medskip\noindent(2) This Assertion follows from Assertion (1) by a
geometric series argument, but we give a direct proof.

As $\Ga$ is a lattice, the family $(\mu_x)_{x\in V\XX_v}$ is actually
equivariant under $\Aut(\XX_v)$,\footnote{See Proposition 
\ref{prop:uniflatmBMfinie} (2).} which acts transitively on the vertices
of $\XX_v$, and the stabiliser in $\Aut(\XX_v)$ of the standard base 
point $*_v$ acts transitively on the edges starting from $*_v$.

Since $\XX_v$ is $(q_v+1)$-regular, since the set of points at
infinity of the geodesic rays starting from $*_v$, whose initial edge
has endpoint $0\in\lk(*_v)=\PP_1(k_v)$, is equal to $\pi_v\OOO_v$,
since all geodesic lines from $\infty\in\partial_\infty \XX_v$ to
points of $\pi_v\OOO_v\subset\partial_\infty \XX_v$ pass through $*_v$,
and by the normalisation of the Patterson density and of the Haar
measure, we have
\begin{align*}
  \|\mu_{*_v}\|&=(q_v+1)\;\mu_{*_v}(\pi_v\OOO_v)=
  (q_v+1)\;\mu_{\H_\infty}(\pi_v\OOO_v)=(q_v+1)\haar_{K_v}(\pi_v\OOO_v)\\
  &= \frac{q_v+1}{q_v} \haar_{K_v}(\OOO_v)=
  \frac{q_v+1}{q_v}\;.
\end{align*}

\noindent 
(3) This follows from Equation \eqref{eq:relatskinhoroballHerPau}, and
from the normalisation $\mu_{\H_\infty}=\haar_{K_v}$ of the Patterson
density.

\medskip
\noindent (4) This follows from Assertion (3) and from
Equation \eqref{eq:massskinhorosph} (where the normalisation of the
Patterson density was different than the one given by Assertion (2)).

\medskip
\noindent(5) Let $L,L_+,L_-$ be as in the statement of Assertion (5),
see the picture below. The result follows from Lemma \ref{lem:hamline}
applied with $\H=\H_\infty$, from Equations \eqref{eq:hamenhaar},
\eqref{eq:expocritlatticeBT} and \eqref{eq:valuationham}.

\begin{center}
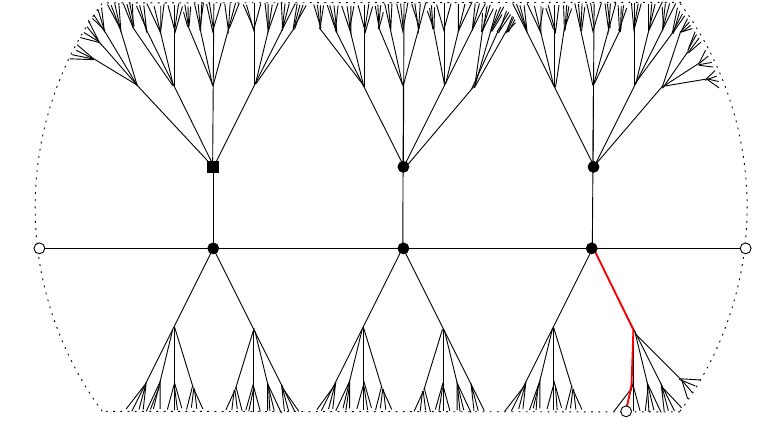
\end{center}

\noindent (6) This follows from Equation \eqref{eq:massskinline}
(where the normalisation of the Patterson density was different),
since $\XX_v$ is $(q_v+1)$-regular, and from Assertion (2):
$$
\|\sigma^\pm_\D\|=
\|\mu_{*_v}\|\;\frac{q_v-1}{q_v+1}\;\Vol(\Ga_L\dbs L)
=\frac{q_v-1}{q_v}\;\Vol(\Ga_L\dbs L)\;.\;\;\;
\Box
$$

\medskip We now turn to measure computations for arithmetic lattices
$\Ga$ in $\XX_v$ in the function field case. We still assume that the
Patterson density of $\Ga$ is normalised so that $\mu_{\H_\infty}
(\OOO_v) =1$, and we denote by $m_{\rm BM}$ the Bowen-Margulis measure
of $\Ga$ associated with this choice of Patterson density.

\bprop\label{prop:q-skinmasses} 
Let $K$ be a function field over $\FF_q$ of genus $g$ and let $v$ be a
valuation of $K$. Let $\Ga$ be a finite index subgroup of
$\Ga_v=\PGL_2(R_v)$, with Patterson density normalised such that
$\mu_{\H_\infty}=\haar_{K_v}$.
\begin{enumerate}
\item We have 
$$
\|m_{\rm BM}\|=\frac{(q_v+1)\;[\Ga_v:\Ga]}{q_v}\;
\Vol(\Ga_v\dbs\XX_v)= 
\frac{2\;(q_v+1)\;\zeta_K(-1)\;[\Ga_v:\Ga]}{q_v}\,,
$$
and if $K=\FF_q(Y)$ and $v=v_\infty$ is the valuation at infinity of
${\bf C}=\PP_1$, then
$$
\|m_{\rm BM}\| =\frac{2\;[\PGL_2(\FF_q[Y]):\Ga]}{q\,(q-1)^2}\,.
$$
\item Let $\Ga_\infty$ be the stabiliser in $\Ga$ of $\infty
  \in \partial_\infty \XX_v$, and let $\D=(\ga\H_\infty)_{\ga\in
    \Ga/\Ga_\infty}$. We have
$$
\|\sigma^\pm_{\D}\|=\frac{q^{\,g-1}\;[(\Ga_v)_\infty:\Ga_\infty]}{q-1}\,.
$$
\end{enumerate}
\eprop

\dem (1) Recall that $\Ga_v=\PGL_2(R_v)$ acts without inversion on
$\XX_v$. By Equation \eqref{eq:relatmBMvolhomog} (which uses a
different normalisation of the Patterson density of $\Ga$), and by
Proposition \ref{prop:mescomputBT} (2), we have
$$
\|m_{\rm BM}\|=\|\mu_{*_v}\|^2\frac {q_v}{q_v+1}\Vol(\Ga\dbs\XX_v)
=\frac {(q_v+1)\;[\Ga_v:\Ga]}{q_v}\;
\Vol(\Ga_v\dbs\XX_v)\;.
$$
The first claim of Assertion (1) hence follows from Equation
\eqref{eq:covolPGLdeR}.  

If $K=\FF_q(Y)$ and $v=v_\infty$, then the second claim of Assertion
(1) follows either from the first claim where the value of
$\zeta_K(-1)$ is given by Equation \eqref{zetamoinsun}, or from the
fact that $q_v=q$ and that the covolume $\Vol(\PGL_2(\FF_q[Y])
\dbs\XX_{v_\infty})$ of the Nagao lattice $\PGL_2(\FF_q[Y])$ is
\begin{equation}\label{eq:q-vol}
\Vol(\PGL_2(\FF_q[Y])\dbs\XX_{v_\infty})= \frac{2}{(q-1)(q^2-1)}\,,
\end{equation}
 as an
easy geometric series computation shows using the description of the
modular ray in Section \ref{subsec:modulargroup} (see also
\cite[Sect.~10.2]{BasLub01}).

\medskip
\noindent(2) Let us prove that
\begin{equation}\label{eq:haarcovolinfty}
\haar_{K_v}((\Ga_v)_\infty\bs K_v)=\frac{q^{\,g-1}}{q-1}\,. 
\end{equation}
The result then follows by Proposition \ref{prop:mescomputBT} (4)
since
$$
\|\sigma^\pm_{\D}\|=\haar_{K_v}(\Ga_\infty\bs K_v)=
[(\Ga_v)_\infty: \Ga_\infty]\;
\haar_{K_v}((\Ga_v)_\infty\bs K_v)\;.
$$

The stabiliser of $\infty=[1:0]$ in $\Ga_v$ acts on $K_v$ exactly by
the set of transformations $z\mapsto az+b$ with $a\in (R_v)^\times$
and $b\in R_v$. Since\footnote{See Equation \eqref{eq:inversRv}.}
$(R_v)^\times=(\FF_q)^\times$ acts freely by multiplication on the
left on $(K_v-R_v)/R_v$, and by Lemma \ref{lem:covolideal}, we have
$$
\haar_{K_v}((\Ga_v)_\infty\bs K_v)=\frac{1}{q-1}\haar_{K_v}(K_v/R_v)=
\frac{q^{\,g-1}}{q-1}\;. 
$$
This proves Equation \eqref{eq:haarcovolinfty}.
\cqfd

\section{Exponential decay of correlation and error terms for 
arithmetic quotients of Bruhat-Tits trees}
\label{subsec:locconst}

As in the beginning of Section \ref{subsec:BruhatTitstrees}, let $K_v$
be a non-Archimedean local field, with valuation $v$, valuation ring
$\OOO_v$, choice of uniformiser $\pi_v$, and residual field $k_v$ of
order $q_v$. Let $\Ga$ be a lattice of the Bruhat-Tits tree $\XX_v$ of
$(\PGL_2,K_v)$. In this Section, we discuss the error terms estimates
that we will use in Part \ref{sect:arithappli}.

\medskip
Partly in order to simplify the references, we start by summarizing
in the next statement the only results from the geometric Part
\ref{part:equid} of this book, on geometric equidistribution and
counting problems, that we will use in this algebraic Part
\ref{sect:arithappli}. We state it with the normalisation introduced
in Section \ref{subsec:volcomputBruhatTitstrees} which will be useful
in what follows.

\btheo\label{theo:algebrogeometricequid} Let $\Ga$ be a lattice of
$\XX_v$ whose length spectrum $L_\Ga$ is equal to $2\ZZ$.  Assume that
the Patterson density of $\Ga$ is normalised so that $\|\mu_x\|=
\frac{q_v+1}{q_v}$ for every $x\in V\XX_v$. Let $\DD^\pm$ be nonempty
proper simplicial subtrees of $\XX_v$ with stabilisers $\Ga_{\DD^\pm}$
in $\Ga$, such that the families $\D^\pm=
(\ga\DD^\pm)_{\ga\in\Ga/\Ga_{\DD^\pm}}$ are locally finite in $\XX_v$.
For every $\ga\in\Ga$ such that $\DD^-$ and $\ga \DD^+$ are disjoint,
let $\alpha^-_{e,\,\ga}$ be the generalised geodesic line, isometric
exactly on $[0,d(\DD^-, \, \ga \DD^+)]$, whose image is the common
perpendicular between $\DD^-$ and $\ga \DD^+$. If the measure
$\sigma^-_{\D^+}$ is nonzero and finite, then
$$
\lim_{n\ra+\infty} \;\frac{({q_v}^2-1)(q_v+1)}{2\;{q_v}^3}\;
\frac{\Vol(\Ga\dbs\XX_v)}{\|\sigma^-_{\D^+}\|}\;q_v^{\,-n}
\sum_{\substack{\ga\in \Ga/\Ga_{\DD^+}\\
0< d(\DD^-, \, \ga \DD^+)\leq n}}
\Dirac_{\alpha^-_{e,\,\ga}} \;=\;
\wt\sigma^+_{\DD^-}\;,
$$
for the weak-star convergence of measures on the locally
compact space $\gengeod \XX_v$.

Furthermore, if $\Ga$ is geometrically finite, then for every $\beta
\in \;]0,1]$, there exists an error term for this equidistribution
claim when evaluated on $\wt\psi \in\C_c^\beta(\gengeod\XX)$ of the
form $\bigO(\,\|\wt\psi\,\|_\beta\;e^{-\kappa\,  n})$ for some
$\kappa>0$.  
\etheo

As recalled at the end of Section \ref{subsec:trees}, lattices in
$\PGL_2(\KK_v)$ are geometrically finite, see \cite{Lubotzky91}. We
will hence be able to use the error term in Theorem
\ref{theo:algebrogeometricequid} in particular when

$\bullet$~ $K_v$ is the completion of a function field $K$ over
$\FF_q$ with respect to a (normalised discrete) valuation $v$ of $K$
and $\Ga$ is a finite index subgroup of $\PGL_2(R_v)$ with $R_v$ the
affine function ring associated with $v$,\footnote{See Section
\ref{subsec:valuedfields} for definitions.} as in Chapters
\ref{sec:rattionalequid} and \ref{subsec:integralrepresentation}, and
in Sections \ref{subsec:quadirratposchar} and
\ref{subsec:countcrossratquadirrat};

$\bullet$~ when $K_v=\QQ_p$ and $\Ga$ is an arithmetic lattice in
$\PGL_2(K_v)$ derived from a quaternion algebra, see Sections
\ref{subsec:quadirratzerochar} and \ref{subsec:countcrossratquadirrat}.

\medskip
\dem In order to prove the first claim, we apply Corollary
\ref{coro:twosubtreeseven} with $\XX=\XX_v$ and $p=q=q_v$. Since
$L_\Ga=2\ZZ$, the lattice $\Ga$ leaves invariant the partition of
$V\XX_v$ into vertices at even distance from a basepoint $x_0\in
V\XX_v$ and vertices at odd distance from $x_0$. Since the Patterson
density is now normalised so that $\|\mu_{x_0}\|= \frac{q_v+1}{q_v}$
(instead of $\|\mu_{x_0}\|= \frac{q_v+1}{\sqrt{q_v}}$ in Corollary
\ref{coro:twosubtreeseven}), the skinning measures
$\wt\sigma^\mp_{\D^\pm}$ are now $\frac{1}{\sqrt{q_v}}$ times the ones
in the statement of Corollary \ref{coro:twosubtreeseven}. Hence the
second assertion of Corollary \ref{coro:twosubtreeseven} gives
$$
\lim_{n\ra+\infty} \;\frac{{q_v}^2-1}{2\;{q_v}^2}\;
\frac{\TVol(\Ga\dbs\XX_v)}{\sqrt{q_v}\;\|\sigma^-_{\D^+}\|}\;{q_v}^{-n}
\sum_{\substack{\ga\in \Ga/\Ga_{\DD^+}\\ 0< d(\DD^-, \, \ga \DD^+)\leq n}}
\Dirac_{\alpha^-_{e,\,\ga}} \;=\;
\sqrt{q_v}\;\,\wt\sigma^+_{\DD^-}\;.
$$
By Equation \eqref{eq:volTvol}, we have
$$
\TVol(\Ga\dbs\XX_v)=(q_v+1)\,\Vol(\Ga\dbs\XX_v)\;.
$$
The first claim follows.

\medskip
The last claim concerning error terms follows from Remark (ii)
following the proof of Theorem
\ref{theo:equidsimplicialperperrorterm}.  
\cqfd

\medskip
In the last four Chapters \ref{sec:rattionalequid},
\ref{sec:nonarchquadratequid}, \ref{sec:crossrat} and
\ref{subsec:integralrepresentation} of this book, we will need to
push to infinity the measures appearing in the statement of
Theorem \ref{theo:algebrogeometricequid}. We regroup in the following
two remarks the necessary control tools  for such a pushing.

\brema{\rm
With the notation of Theorem \ref{theo:algebrogeometricequid}, we will
use Lemma \ref{lem:holderextensiontoray} when $X=|\XX_v|_1$ is the
geometric realisation of the simplicial tree $\XX_v$, $\alpha=
\alpha^-_{e,\,\ga}$ is\footnote{the generalised geodesic line
  isometric exactly on $[0,d(\DD^-,\ga\DD^+)]$ parametrising} the
common perpendicular between $\DD^-$ and $\ga\DD^+$ for $\ga\in\Ga$
(when it exists), and $\rho= \rho_\ga$ is any extension of $\alpha$ to
a geodesic ray, or rather to a generalised geodesic line isometric
exactly on $[0,+\infty[\,$. Under the assumptions of Theorem
    \ref{theo:algebrogeometricequid}, we have by Lemma 
\ref{lem:holderextensiontoray}
\begin{equation}\label{eq:extensgeodray}
\frac{({q_v}^2-1)(q_v+1)}{2\;{q_v}^3}\;
\frac{\Vol(\Ga\dbs\XX_v)}{\|\sigma^-_{\D^+}\|}\;q_v^{\,-n}
\sum_{\substack{\ga\in \Ga/\Ga_{\DD^+}\\
0< d(\DD^-, \, \ga \DD^+)\leq n}}
\Dirac_{\rho_\ga} \;\;\stackrel{*}{\rightharpoonup}\;\;
\wt\sigma^+_{\DD^-}\;,
\end{equation}
with, if $\Ga$ is geometrically finite, an error term when
evaluated on $\wt\psi \in\C_c^\beta(\gengeod\XX_v)$ of the form
$\bigO(e^{-\kappa \,n}\;\|\wt\psi\,\|_\beta)$ for some $\kappa>0$
small enough (depending in particular on $\beta\in\;]0,1]$).
}
\erema

\medskip
From now on in this book, for every subtree $\DD$ of $\XX_v$ with
geometric realisation $D=|\DD|_1$, we endow $\partial_\infty \XX_v-
\partial_\infty D$ with the distance-like map $d_D$ defined in
Equation \eqref{eq:distancelike}. We use this map $d_D$ in order to
define both the $\beta$-H\"older-continuity of maps with values in
$(\partial_\infty \XX_v- \partial_\infty D,d_D)$ and the
$\beta$-H\"older-norm of a function defined on $(\partial_\infty
\XX_v- \partial_\infty D,d_D)$.

\brema\label{rem:utidD} {\rm With the notation of Theorem
  \ref{theo:algebrogeometricequid}, we will use Proposition
  \ref{prop:partialplusholder} when $\XX= \XX_v$ and $\DD=\DD^-$.
  Under the assumptions of Theorem \ref{theo:algebrogeometricequid},
  with $\rho_\ga$ any extension to a geodesic ray of the common
  perpendicular $\alpha^-_{e,\,\ga}$ between $\DD^-$ and $\ga\DD^+$
  for $\ga\in\Ga$, since pushing forward measures on $\normalout
  \DD^-$ by the homeomorphism $\partial^+:\rho\mapsto \rho_+$
  introduced in Proposition \ref{prop:partialplusholder} is
  continuous, we have by Equation \eqref{eq:extensgeodray}
\begin{equation}\label{eq:extensgeodraytoinfty}
\frac{({q_v}^2-1)(q_v+1)}{2\;{q_v}^3}\;
\frac{\Vol(\Ga\dbs\XX_v)}{\|\sigma^-_{\D^+}\|}\;q_v^{\,-n}
\sum_{\substack{\ga\in \Ga/\Ga_{\DD^+}\\
0< d(\DD^-, \, \ga \DD^+)\leq n}}
\Dirac_{(\rho_\ga)_+} \;\;\stackrel{*}{\rightharpoonup}\;\;
(\partial^+)_*\wt\sigma^+_{\DD^-}\;.
\end{equation}
If $\Ga$ is geometrically finite, for all $\beta\in\;]0,1]$ and $\psi
\in\C_c^\beta(\partial_\infty \XX_v-\partial_\infty \DD^-)$, using the
error term in Equation \eqref{eq:extensgeodray} with regularity
$\frac{\beta}{2}$ when evaluated on $\wt\psi =\psi\circ\partial^+$,
which belongs to $\C_c^{\frac{\beta}{2}}(\gengeod\XX_v)$ by the first
claim of Proposition \ref{prop:partialplusholder}, we have by the last
claim of Proposition \ref{prop:partialplusholder} an error term in
Equation \eqref {eq:extensgeodraytoinfty} evaluated on $\psi$ of the
form $\bigO(e^{-\kappa \,n}\;\|\psi\,\|_\beta)$ for some $\kappa>0$
small enough.  } \erema

A stronger assumption than the H\"older regularity is the locally
constant regularity, that has been defined at the end of Section
\ref{subsec:holdercont}, and is applicable here since the involved
spaces are totally disconnected.  Several error terms estimates in the
literature use this stronger regularity (see for instance
\cite{AthGhoPra12, KemPauSch17}). The following result of decay of
correlations under locally constant regularity follows from Corollary
\ref{coro:expdecaygeomfinisimpl} by Remark \ref{rem:locconstholder}
and since\footnote{See the end of Section \ref{subsec:trees}.} any
lattice of $\PGL_2(K_v)$ is geometrically finite.\footnote{Note that
  being a lattice of $\PGL_2(K_v)$ is much stronger than been a
  lattice in $\Aut(\XX_v)$.}

\bprop\label{prop:expodecaylc} Assume that $\Ga$ is a lattice of
$\PGL_2(K_v)$, and let $\beta\in\mathopen{]}0,1\mathclose{]}$.

(1) If $L_\Ga=\ZZ$, there exist $C,\kappa>0$ such that for every
$\epsilon\in\;]0,1]$, for all $\epsilon$-locally constant maps
$\phi,\psi : \Ga\bs\G\XX_v\ra\RR$ and $n\in\ZZ$, we have
\begin{align*}
\Big|\int_{\Ga\bs\G\XX_v}\phi\circ\flow{-n}\;\psi\; d\, m_{\rm BM}
- &\frac{1}{\|m_{\rm BM}\|}
\int_{\Ga\bs\G\XX_v} \phi\; d\, m_{\rm BM}
\int_{\Ga\bs\G\XX_v} \psi\;d\, m_{\rm BM}\;\Big| \\ &\le
C\;e^{-\kappa|n|}\;\|\phi\|_{\epsilon\,{\rm lc},\,\beta}\;
\|\psi\|_{\epsilon\,{\rm lc},\,\beta}\;. 
\end{align*}

(2) If $L_\Ga=2\ZZ$, then there exist $C,\kappa >0$ such that for
every $\epsilon\in\;]0,1]$, for all $\epsilon$-locally constant maps
$\phi,\psi : \Ga\bs\Geven\XX_v\ra\RR$ and $n\in\ZZ$, we have
\begin{align*}
\Big|\int_{\Ga\backslash \Geven \XX}\phi\circ\flow{-2n}\;\psi\; d\, m_{\rm BM}
- &\frac{1}{ m_{\rm BM}(\Ga\backslash \Geven \XX)}
\int_{\Ga\backslash \Geven \XX} \phi\; d\, m_{\rm BM}
\int_{\Ga\backslash \Geven \XX} \psi\;d\, m_{\rm BM}\;\Big| \\ &\le\; 
C\;e^{-\kappa|n|}\;\|\phi\|_{\epsilon\,{\rm lc},\,\beta}\;
\|\psi\|_{\epsilon\,{\rm lc},\,\beta}\;.\;\;\;
\Box
\end{align*}
\eprop

We will not use the following result in this book, but its Assertion
(2) is used in the announcement \cite{BroParPau16CRAS} which only
considers the locally constant regularity.

\btheo\label{theo:algebroequid} Assume that $\Ga$ is a lattice of
$\PGL_2(K_v)$, and let $\beta\in\mathopen{]}0,1\mathclose{]}$.

(1) If $L_\Ga=\ZZ$, there exists $\kappa>0$ such that for every
$\epsilon\in\;]0,1]$ and every $\epsilon$-locally constant map
$\wt\psi : \gengeod\XX_v\ra\RR$, we have, as $n\ra+\infty$,
\begin{align*}
\frac{({q_v}-1)}{(q_v+1)\;{q_v}}\; 
\frac{\Vol(\Ga\dbs\XX_v)}{\|\sigma^-_{\D^+}\|}\;q_v^{\,-n}
& \sum_{\substack{\ga\in \Ga/\Ga_{\DD^+}\\
0< d(\DD^-, \, \ga \DD^+)\leq n}} 
\wt\psi (\alpha^-_{e,\,\ga})\\=\;&
\int_{\gengeod\XX_v}\;\wt\psi\;d\wt\sigma^+_{\DD^-}+
\bigO\big(e^{-\kappa\, n}\;\|\wt\psi\,\|_{\epsilon\,{\rm lc},\,\beta}\big) \;.
\end{align*}

(2) If $L_\Ga=2\ZZ$, there exists $\kappa>0$ such that for every
$\epsilon\in\;]0,1]$ and every $\epsilon$-locally constant map
$\wt\psi : \gengeod\XX_v\ra\RR$, we have, as $n\ra+\infty$,
\begin{align*}
\frac{({q_v}^2-1)(q_v+1)}{2\;{q_v}^3}\;
\frac{\Vol(\Ga\dbs\XX_v)}{\|\sigma^-_{\D^+}\|}\;q_v^{\,-n}
& \sum_{\substack{\ga\in \Ga/\Ga_{\DD^+}\\
0< d(\DD^-, \, \ga \DD^+)\leq n}}
\wt\psi(\alpha^-_{e,\,\ga})\\=\;&
\int_{\gengeod\XX_v}\;\wt\psi\;d\wt\sigma^+_{\DD^-}+
\bigO\big(e^{-\kappa\, n}\;\|\wt\psi\,\|_{\epsilon\,{\rm lc},\,\beta}\big)\;.
\end{align*}
\etheo

\dem This result follows by replacing in the proof of Theorem
\ref{theo:equidsimplicialperperrorterm} (or rather Remark (ii)
following its proof) the use of the exponential decay of
$\beta$-Hölder correlations given by Corollary
\ref{coro:expdecaygeomfinisimpl} by the use of Proposition
\ref{prop:expodecaylc}.  
\qed

\medskip
We conclude this Section by giving an algebraic second proof of
Theorem \ref{theo:algebroequid}, in the following situation. Let $K_v$
be the completion of a function field $K$ over $\FF_q$ with respect to
a valuation $v$ of $K$ and let $\Ga$ be a nonuniform\footnote{This
  assumption is introduced in order to apply the following Theorem
  \ref{theo:edmc}. Note that the existence of a nonuniform lattice in
  $G_v=\PGL_2(K_v)$ forces the characteristic to be positive, see for
  instance \cite{Lubotzky91}.}  lattice of $G_v=\PGL_2(K_v)$. We only
obtain a version using an exponent $\beta\geq\ln q_v$ for the locally
constant norm. For simplicity, we consider only Assertion (1) of
Theorem \ref{theo:algebroequid}, though Assertion (2) xould be treated
similarly.

\medskip  
The group $G_v$ acts (on the left) on the complex vector space of maps
$\psi$ from $\Ga\bs G_v$ to $\RR$, by right translation on the source:
For every $g\in G_v$, we have $g\psi:x\mapsto \psi (xg)$.  A function
$\psi:\Ga\bs G_v\ra\RR$ is {\it algebraically locally
constant}\index{algebraically locally constant} if there exists a
compact-open subgroup $U$ of $\PGL_2(\OOO_v)$ which leaves $\psi$
invariant:
$$
\forall\;g\in U,\;\;\; g\psi=\psi\;,
$$
or equivalently, if $\psi$ is constant on each orbit of $U$ under the
right action of $G_v$ on $\Ga\bs G_v$. Note that $\psi$ is then
continuous, since the orbits of $U$ are compact-open subsets. We
define
$$
d_\psi=\dim \big(\operatorname{Vect}_\RR (\PGL_2(\OOO_v) \psi)\big)
$$
as the dimension of the complex vector space generated by the images
of $\psi$ under the elements of $\PGL_2(\OOO_v)$, which is finite,
and even satisfies 
$$
d_\psi\leq [\PGL_2(\OOO_v):U]\;.
$$ 

We define the {\it ${\rm alc}$-norm}\index{alc@${\rm alc}$-norm}%
\index{norm!${\rm alc}$} of every bounded algebraically
locally constant map $\psi:\Ga\bs G_v\ra\RR$ by
$$
\|\psi\|_{{\rm alc}}=\sqrt{d_\psi }\,\|\psi\|_\infty\;.
$$
Though the ${\rm alc}$-norm does not satisfy the triangle inequality, we
have $\|\lambda\,\psi\|_{{\rm alc}}=|\lambda|\,\|\psi\|_{{\rm alc}}$ for every
$\lambda\in\RR$. We denote by ${\rm alc}(\Ga\bs G_v)$ the vector space of
bounded algebraically locally constant maps $\psi$ from $\Ga\bs G_v$
to $\RR$.

For every $n\in\NN$, let $U_n$ be the compact-open subgroup of
$\PGL_2(\OOO_v)$ which is the kernel of the morphism
$\PGL_2(\OOO_v)\ra \PGL_2(\OOO_v/ {\pi_v}^n\OOO_v)$ of reduction
modulo ${\pi_v}^n$. Note that any compact-open subgroup $U$ of
$\PGL_2(\OOO_v)$ contains $U_n$ for some $n\in\NN$. Hence $\psi:\Ga\bs
G_v\ra\RR$ is algebraically locally constant if and only if there
exists $n\in\NN$ such that $\psi$ is constant on each right orbit of
$U_n$. For every $n\in\NN$, since the order of $\PGL_2(\OOO_v/
{\pi_v}^n\OOO_v)$ is at most the order of $(\OOO_v/
{\pi_v}^n\OOO_v)^4$, which is ${q_v}^{4n}$, if $\psi:\Ga\bs G_v\ra\RR$
is constant on each right orbit of $U_n$, then
\begin{equation}\label{eq:relatalcepsilonlc}
\|\psi\|_{{\rm alc}}\leq {q_v}^{2n}\,\|\psi\|_\infty\;.
\end{equation}

Recall\footnote{See Section \ref{subsec:BruhatTitstrees}.} that we
have a natural homeomorphism $\,\Xi: \Ga g\uA(\OOO_v) \mapsto \Ga g
\,\ell^*$ between the double coset space $\Ga\bs G_v/\uA(\OOO_v)$ and
$\Ga\bs\G\XX_v$. We denote by $p_{\G}:\Ga\bs G_v\ra \Ga\bs\G\XX_v$ the
composition map of the canonical projection $\big(\Ga\bs G_v\big)\ra
\big(\Ga\bs G_v/\uA(\OOO_v)\big)$ and of $\Xi$. By Equation
\eqref{eq:commutgeodfowrightasubv}, for all $x\in \Ga\bs G_v$ and
$n\in\NN$, we have
\begin{equation}\label{eq:commutgeodfowrightasubvd}
p_{\G}\,(x\;{a_v}^n)=\flow{n}\,p_{\G}(x)\;.
\end{equation}

\blemm 
For every $\epsilon\in\;]0,1]$, for every $\epsilon$-locally
constant function $\psi : \Ga\bs\G\XX_v\ra\RR$, if $n=\lceil -
\frac{1}{2}\ln \epsilon \rceil$, then the map $\psi \circ p_{\G} : 
\Ga\bs G_v \ra \RR$ is $U_n$-invariant and
\begin{equation}\label{eq:relatalcepsionlc}
\|\psi\circ p_{\G}\|_{{\rm alc}} \leq 
{q_v}^2\;\|\psi\|_{\epsilon\,{\rm lc},\,\ln q_v}\;.
\end{equation}
\elemm

\dem 
Let $\epsilon,\psi,n$ be as in the statement. Let us first prove that
if $\ell,\ell'\in\G \XX_v$ satisfy $\ell_{[-n,+n]}=\ell'_{[-n,+n]}$, then
$d(\ell,\ell')\leq \epsilon$.

If $\ell_{[-n,+n]}=\ell'_{[-n,+n]}$, then $d(\ell(t),\ell'(t))=0$ for
$t\in[-n,n]$ and by the triangle inequality $d(\ell(t),\ell'(t))\leq
2(|t|-n)$ if $|t|\geq n$, hence
\begin{align*}
d(\ell,\ell') &\leq 2\int_{n}^{+\infty} \;2\,(t-n)\;e^{-2\,t}\;dt=
2\;e^{-2\,n}\int_{0}^{+\infty} \;u\;e^{-u}\;\frac{du}{2}=
e^{-2\,n}\\ &\leq e^{-2(-\frac{1}{2}\ln \epsilon)}=
\epsilon\;,
\end{align*}
as wanted.

\medskip
In order to prove that $\psi \circ p_{\G} : \Ga\bs G_v \ra \RR$ is
$U_n$-invariant, let $x,x'\in \Ga\bs G_v$ be such that $x'\in x\,U_n$.
Since $U_n$ acts by the identity map on the ball of radius $n$ in the
Bruhat-Tits tree $\XX_v$, the geodesic lines $p_{\G}(x)$ and
$p_{\G}(x')$ in $\Ga\bs\G\XX_v$ coincide (at least) on $[-n,n]$.
Hence, as we saw in the beginning of the proof, we have $d(p_{\G}(x),
p_{\G}(x'))\leq \epsilon$.  Therefore $\psi (p_{\G}(x)) =\psi
(p_{\G}(x'))$ since $\psi $ is $\epsilon$-locally constant.

Now, using Equation \eqref{eq:relatalcepsilonlc}, we have 
\begin{align*}
\|\psi \circ p_{\G}\|_{{\rm alc}}& \leq 
{q_v}^{2n}\,\|\psi \circ p_{\G}\|_\infty\\ &\leq 
{q_v}^{2(1-\frac{1}{2}\ln\epsilon)}\,\|\psi\|_\infty=
{q_v}^2\;\epsilon^{-\ln q_v}\,\|\psi\|_\infty 
={q_v}^2\,\|\psi\|_{\epsilon\,lc,\,\ln q_v}\;. \;\;\;\;\Box
\end{align*}

\medskip 
Now, we will use an algebraic result of exponential decay of
correlations, Theorem \ref{theo:edmc} (see for instance
\cite{AthGhoPra12}). We first recall some definitions and notation,
useful for its statement.

\medskip 
Recall that the left action of the locally compact unimodular group
$G_v$ on the locally compact space $\G\XX_v$ is continuous and
transitive, and that its stabilisers are compact hence unimodular.
Since $\Ga$ is a lattice, the (Borel, positive, regular)
Bowen-Margulis measure $\wt m_{\rm BM}$ on $\G\XX_v$ is
$G_v$-invariant (see Proposition \ref{prop:uniflatmBMfinie} (2)).
Hence by \cite{Weil65} (see also \cite[Lem.~5]{GorPau14}), there
exists a unique Haar measure on $G_v$, which disintegrates by the
evaluation map $\wt{p_\G}:G_v\ra \G\XX_v$ defined by $g\mapsto
g\ell_*$, with conditional measure on the fiber over $\ell=
g\ell_*\in\G\XX_v$ the probability Haar measure on the stabiliser
$g\uA(\OOO_v)g^{-1}$ of $\ell$ under $G_v$. Hence, taking the quotient
under $\Ga$ and normalising in order to have probability measures, if
$\mu_v$ is the right $G_v$-invariant probability measure on $\Ga\bs
G_v$, we have
\begin{equation}\label{eq:pushforwardpcalG}
(p_\G)_*\mu_v=\frac{m_{\rm BM}}{\|m_{\rm BM}\|}\;.
\end{equation}

For every $g\in G_v$, we denote by $R_g:\Ga\bs G_v\ra \Ga\bs G_v$ the
right action of $g$, and for all bounded continuous functions $\wt
\psi,\wt \psi'$ on $\Ga\bs G_v$, we define
$$
\gls{correlationcoeffalgeb}(\wt \psi,\wt \psi')=
\int_{\Ga\bs G_v}\wt \psi\circ R_g\;\wt \psi'\;d\mu_v
\;-\;\int_{\Ga\bs G_v}\wt \psi\;d\mu_v\;
\int_{\Ga\bs G_v}\wt \psi'\;d\mu_v\;.
$$
Note that by Equations \eqref{eq:pushforwardpcalG} and
\eqref{eq:commutgeodfowrightasubvd}, for all bounded continuous
functions $\psi,\psi':\Ga\bs \G\XX_v\ra\RR$ and $n\in\ZZ$, we
have\footnote{See Section \ref{subsec:mixingratesimpgraphs} for a
definition of $\operatorname{cov}_{\mu,\,n}$.}
\begin{equation}\label{eq:covalgcovgeo}
\operatorname{cov}_{\frac{m_{\rm BM}}{\|m_{\rm BM}\|},\,n}(\psi,\psi')
=\operatorname{cov}_{\mu_v,\,{a_v}^n}(\psi\circ p_\G,\psi'\circ p_\G)\;.
\end{equation}

Recall that the adjoint representation of $G_v=\PGL_2(K_v)$ is the
continuous morphism $G_v\ra\GL(\M_2(K_v))$ defined by $[h]\mapsto
\{x\mapsto hxh^{-1}\}$, which is independent of the choice of the
representative $h\in \GL_2(K_v)$ of $[h]\in\PGL_2(K_v)$. For every
$g\in G_v$, we denote by $|g|_v$ the operator norm of the adjoint
representation of $g$.  For instance, recalling that
$a_v=\begin{bmatrix} 1 & 0 \\ 0 & \pi_v^{-1}\end{bmatrix}$, we have,
for all $n\in\ZZ$,
\begin{equation}\label{eq:controlnormav}
|{a_v}^n|_v={q_v}^{|n|}\;.
\end{equation}

We refer for instance to \cite{AthGhoPra12} for the following result
of exponential decay of correlations.

\btheo\label{theo:edmc} Let $\Ga$ be a nonuniform lattice of
$G_v$. There exist $C',\kappa'>0$ such that, for all bounded locally
constant functions $\wt \psi,\wt \psi':\Ga\bs G_v\ra\RR$ and $g\in
G_v$,
\begin{equation}\label{eq:1102}
\big|\operatorname{cov}_{\mu_v,\,g}(\wt \psi,\wt \psi')\,\big|
\leq \;C'\;\|\wt \psi\|_{{\rm alc}}\;\|\wt \psi'\|_{{\rm alc}}\; 
{|g|_v}^{-\kappa'}\;.\;\;\;\Box
\end{equation}
\etheo

Proposition \ref{prop:expodecaylc} (1) and therefore Theorem
\ref{theo:algebroequid} (1) with $\beta\geq\ln q_v$ follows from this
result applied to $\wt \psi=\psi\circ p_\G$, $\wt \psi'=\psi'\circ
p_\G$ and $g={a_v}^n$ by using Equations \eqref{eq:covalgcovgeo},
\eqref{eq:relatalcepsionlc} and \eqref{eq:controlnormav} and by taking
$C=C'q_v^4$ and $\kappa=\kappa'\ln q_v$.  \cqfd

\bigskip
\rem 
There is a similar relationship between locally constant functions on
$K_v$ in an algebraic sense and the ones in the metric sense.

The additive group $(K_v,+)$ acts on the complex vector space of
functions from $K_v$ to $\RR$, by translations on the source: for all
$y\in K_v$ and $\psi:K_v\ra\RR$, the function $y\cdot \psi$ is equal
to $x\mapsto \psi(x+y)$. A function $\psi:K_v\ra\RR$ is {\it
  algebraically locally constant}%
\index{algebraically locally constant} if there exists $k\in\NN$ such
that $\psi$ is invariant under the action of the compact-open subgroup
$(\pi_v)^k\OOO_v$ of $K_v$, that is, if for all $x\in K_v$ and $y\in
(\pi_v)^k\OOO_v$, we have $\psi(x+y)=\psi (x)$. Note that a locally
constant function from $K_v$ to $\RR$ is continuous.

For any locally constant function $\psi:K_v\to\RR$, the complex vector
space $\operatorname{Vect}_\RR (\OOO_v\cdot\psi)$ generated by the
images of $\psi$ under the elements of $\OOO_v$ is finite
dimensional. Its dimension $d_\psi$ satisfies, with $k$ as above,
$$
d_\psi\leq [\OOO_v:(\pi_v)^k\OOO_v]={q_v}^k\;.
$$
We define the {\it ${\rm alc}$-norm}\index{alc@${\rm alc}$-norm}%
\index{norm!${\rm alc}$} of every bounded algebraically locally constant 
function $\psi:K_v\ra\RR$
by
$$
\|\psi\|_{{\rm alc}}= \sqrt{d_\psi}\;\|\psi\|_\infty\;.
$$
Though the ${\rm alc}$-norm does not satisfy the triangle inequality, we
have $\|\lambda \psi\|_{{\rm alc}}=|\lambda |\,\|\psi\|_{{\rm alc}}$ for every
$\lambda\in\RR$, and the set of bounded algebraically locally constant
maps from $K_v$ to $\RR$ is a real vector space.

Actually, a function $\psi:K_v\ra\RR$ is algebraically locally constant
if and only if it is locally constant. More precisely, for every
$\epsilon\in\;]0,1]$, since the closed balls of radius $q_v^{-k}$ in
$K_v$ are the orbits by translations under $(\pi_v)^k\OOO_v$, every
$\epsilon$-locally constant function $\psi:K_v\ra\RR$ is constant
under the additive action of $(\pi_v)^k\OOO_v$ for
$k=\lceil\frac{-\ln \epsilon}{\ln q_v}\rceil$, hence
$$
\|\psi\|_{{\rm alc}}\leq \|\psi\|_{\epsilon\,{\rm lc},\,\frac{1}{2}}\;.
$$

\section{Geometrically finite lattices with 
infinite Bowen-Margulis mea\-sure}
\label{subsec:exampleinfiniteBM}

This Section is a digression from the theme of arithmetic
applications, in which we use the Nagao lattice defined in Section
\ref{subsec:modulargroup} in order to construct a geometrically finite
discrete group of automorphisms of a simplicial tree which has
infinite Bowen-Margulis measure. This example was promised after
Proposition \ref{prop:uniflatmBMfinie}, where we saw that such
examples do not exist for uniform trees.

We will equivariantly change the lengths of the edges of a simplicial
tree $\XX$ endowed with a geometrically finite (nonuniform) lattice
$\Ga$ in order to turn $\XX$ into a metric tree in which the group
$\Ga$ remains a geometrically finite lattice, but now has a
geometrically finite subgroup with infinite Bowen-Margulis
measure. This example is an adaptation of the negatively curved
manifold example of \cite[\S 4]{DalOtaPei00}. The simplicial example
is obtained as a modification of the metric tree example.

\btheo \label{theo:geominfinite} 
There exists a geometrically finite discrete group of automorphisms of
a metric tree with constant degrees, whose Bowen-Margulis measure is
infinite.

There exists a geometrically finite discrete group of automorphisms of
a simplicial tree with uniformly bounded degrees whose Bowen-Margulis
measure is infinite.  
\etheo

\dem 
Let $v=v_\infty$ be the valuation at infinity of $K=\FF_q(Y)$, $K_v=
\FF_q((Y^{-1}))$, $\OOO_v=\FF_q[[Y^{-1}]]$ and $R_v=\FF_q[Y]$.  Let
$\XX_v$ be the Bruhat-Tits tree of $(\PGL_2,K_v)$ with base point
$*_v=[\OOO_v\times\OOO_v]$. Let $\Ga_v=\PGL_2(R_v)$, which is a
lattice of $\XX_v$, with quotient graph of groups the modular ray
$\Ga_v\dbs\XX_v$ described in Section \ref{subsec:modulargroup}.  We
denote by $(y_i)_{i=-1,0,\dots}$ the ordered vertices along $\Ga_v\dbs
\XX_v$ with vertex stabilisers $(\Ga_i)_{i=-1,0,\dots}$, and by
$(e_i)_{i\in\NN}$ the ordered edges along $\Ga_v\dbs\XX_v$ (pointing
away from the origin of the modular ray).

The subgroup 
$$
P=\bigcup_{i\ge 0}\Ga_i=
\bigg\{\begin{bmatrix}a&Q\\0&d\end{bmatrix}:Q\in\FF_q[Y], a,d\in
\FF_q^\times\bigg\}
$$ 
is the stabiliser in $\Ga$ of $\infty\in\partial_\infty \XX_v$.  Let
$P_0$ be the finite index subgroup of $P$ consisting of the elements
$\begin{bmatrix}1&Q\\0&1\end{bmatrix}$ with $Q(0)=0$.  Observing that
$d(\ga *_v,*_v)=2(i+1)$ for any $\ga\in \Ga_i - \Ga_{i-1}$ and that the
cardinality of $(\Ga_i- \Ga_{i-1})\cap P_0$ is $(q-1)q^{i+1}$, it is
easy to see that the Poincar\'e series
$$
\Q_{P_0,\,0,\,*_v,\,*_v}(s)=\sum_{\ga\in P_0} e^{-s\,d(*_v, \ga *_v)}
$$ 
of the discrete (though elementary) subgroup $P_0$ of $\Isom(\XX_v)$
is (up to a multiplicative constant) equal to $\sum_{i=0}^\infty q^i
e^{-2si}$, which gives $\delta_{P_0}=\frac{\ln q}2$ for the critical
exponent of $P_0$ on $\XX_v$.

Let $h$ be an element of $\Ga_v$ which is loxodromic on $\XX_v$ and
whose fixed points belong to the open subset $Y^{-1}\OOO_v$ of
$K_v=\partial_\infty \XX_v-\{\infty\}$.  Hence the horoball
$\H_\infty$ centred at $\infty\in\partial_\infty \XX_v$, whose
horosphere contains $*_v$, is disjoint from the translation axis
$\Ax_h$ of $h$.  Note that the stabiliser of $\H_\infty$ in $\Ga_v$ is
$P$ and that $P_0$ acts freely on the edges exiting $\H_\infty$.  Let
$x_0\in V\XX_v$ be the closest point on $\Ax_h$ to $\H_\infty$, let
$e_*$ be the edge with origin $*_v$ pointing towards $x_0$, and let
$e_-, e_+$ be the two edges with origin $x_0$ on $\Ax_h$.

\begin{center}
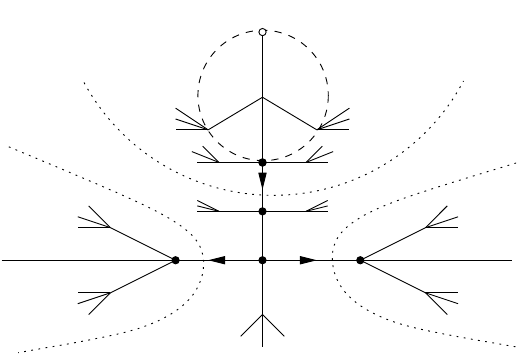
\end{center}

Let $U_h$ be the set of points $x$ in $V\XX_v-\{x_0\}$ such that the
geodesic segment from $x_0$ to $x$ starts either by the edge $e_-$ or
by $e_+$.  Let $U_{P_0}$ be the set of points $y$ in $V\XX_v-\{t(e_*)\}$
such that the geodesic segment from $t(e_*)$ to $y$ starts by the edge
$\overline{e_*}$. We have
\begin{enumerate}
\item $U_h\cap U_{P_0}=\emptyset$ and $x_0\notin U_h\cup U_{P_0}$,
\item $h^k(V\XX_v-U_h)\subset U_h$ for every $k\in\ZZ-\{0\}$ and
  $w(V\XX_v-U_{P_0})\subset U_{P_0}$ for every $w\in P_0-\{\id\}$,
\item $d(x,y)=d(x,x_0)+d(x_0,y)$ for all $x\in U_h$ and $y\in U_{P_0}$.
\end{enumerate}

Let $\Ga'$ be the subgroup of $\Ga_v$ generated by $P_0$ and $h$. By a
ping-pong argument, $\Ga'$ is a free product of $P_0$ and of the
infinite cyclic group generated by $h$, and $\Ga'$ is geometrically
finite (see for instance \cite[Theorem C.2 (xi)]{Maskit88} for
Kleinian groups). Hence every element $\ga$ in $\Ga'-\{e\}$ may be
written uniquely as a word $w_0h^{n_0}w_1h^{n_1}\dots w_kh^{n_k}$ with
$k\in\NN$, $w_i\in P_0, n_i\in\ZZ$ with $w_i\neq e$ if $i\neq 0$ and
$n_i\neq 0$ if $i\neq k$. Using the above properties, we have by
induction
\begin{equation}\label{eq:calcdop}
d(x_0,w_0h^{n_0}w_1h^{n_1}\dots
w_kh^{n_k}x_0)=\sum_{0\leq i\leq k} d(x_0,h^{n_i}x_0)+
\sum_{0\leq i\leq k}d(x_0,w_i x_0)\;.
\end{equation}

Let $\lambda:E\XX_v\to\RR_+$ be the $\Ga'$-invariant length map on the
set of edges of $\XX_v$ such that for every $i\in\NN$, the length of
$e\in E\XX_v$ is $1$ if $e$ is not contained in $\bigcup_{\ga\in
  \Ga'}\ga\H_\infty$, and otherwise, if $e$ maps to $e_i$ or to
$\overline{e_i}$ under the canonical map $\XX_v \ra \Ga_v\bs \XX_v$,
then $\lambda(e)=1+\ln\frac{i+1}{i}$ if $i\geq 1$ and $\lambda(e)=1$
if $i=0$. Note that the distance in the metric graph $|\XX_v|_\lambda$
from $*_v$ to the vertex on the geodesic ray from $*_v$ to $\infty$
originally at distance $i$ from $*_v$ is now $i+\ln\, i$. The
distances along the translation axis of $h$ have not changed. Equation
\eqref{eq:calcdop} remains valid with the new distance. Let us denote
by $d'$ this (new) distance on $|\XX_v|_\lambda$, and by $\Q'_{\Ga'}(s) 
=\Q'_{\Ga',\,0,\,x_0,\,x_0}(s)$ and $\Q'_{P_0}(s)$ the Poincaré series
for the actions of $\Ga'$ and $P_0$ on $(|\XX_v|_\lambda,d')$ (taking
$(x_0,x_0)$ as pair of basepoints and the zero system of
conductances).

\medskip
Let us now prove that the discrete subgroup $\Ga'$ of automorphisms of
the metric tree $(\XX_v, \lambda)$ (with degrees all equal to $q+1$)
satisfies the first claim of Theorem \ref{theo:geominfinite}.

\medskip 
By $\Ga'$-invariance of $\lambda$, the group $\Ga'$ remains a discrete
subgroup of $\Aut(\XX_v, \lambda)$.  The elements of $\Ga'\infty$, or
equivalently, the points at infinity of the horoballs in the
$\Ga'$-equivariant family of horoballs $(\ga\H_\infty)_{\ga\in
  \Ga'/\Ga'_\infty}$ with pairwise disjoint interiors in $\XX_v$,
remain bounded parabolic fixed points of $\Ga'$, and the other limit
points remain conical limit points of $\Ga'$. Hence $\Ga'$ remains a
geometrically finite discrete subgroup of $\Aut(\XX_v,\lambda)$.

The Poincar\'e series $\Q'_{P_0}$ is (up to a multiplicative and
additive constant) $\sum_{i=1}^\infty q^i \, e^{-2si} \,i^{-2s}$,
which has the same critical exponent $\delta_{P_0}= \frac{\ln q}2$ as
previously, but it is easy to see that the discrete group $P_0$ is now
of convergence type if $q\geq 3$.

Let $\lambda(h)$ be the (old and new) translation length of $h$.
Using Equation \eqref{eq:calcdop} and partitionning the above words by
the number of nonzero powers of $h$ they contain, we have, for
every $s\geq 0$,
\begin{align*}
& \Q'_{\Ga'}(s) =\sum_{\ga\in \Ga'} e^{-s\,d'(x_0, \ga x_0)}\\
&\leq \Big(1+
\sum_{k=1}^{+\infty}\Big(\sum_{i\in\ZZ-\{0\}}e^{-s\, |i|\,\lambda(h)}  
\sum_{w\in P_0}e^{-s\,d'(x_0,w x_0)}\Big)^k\Big)
\sum_{w\in P_0} e^{-s\, d'(x_0,wx_0)}\,.
\end{align*}
At $s=\delta_{P_0}$, the above sums over $i\in\ZZ-\{0\}$ and over $w
\in P_0$ converge, and by replacing $h$ by a high enough power if
necessary, we may assume that $\sum_{i\in\ZZ-\{0\}} e^{-s\,|i|\lambda(h)} 
\sum_{w\in P_0} e^{-s\,d'(x_0,w x_0)} <1$ if $s= \delta_{P_0}$, which
makes the Poincar\'e series of $\Ga'$ converge at $s=\delta_{P_0}$.
Since $P_0$ is a subgroup of $\Ga'$ with critical exponent $\delta_P$,
we have that $\Q'_{\Ga'}(s)\geq \Q'_{P_0}(s) = +\infty$ if
$s<\delta_{P_0}$. We conclude that the critical exponent of $\Ga'$ is
equal to $\delta_{P_0}$ and that $\Ga'$ is of convergence type.  By
Corollary \ref{coro:finitudeGibbsdivuniq} (1), the Bowen-Margulis measure
of $\Ga'$ is infinite.

\medskip In order to prove the second claim of Theorem
\ref{theo:geominfinite}, we first define a new length map
$\lambda:E\XX_v\to\RR_+$ which coincides with the previous one on
every edge $e$ of $\XX_v$, unless $e$ maps to $e_i$ or to
$\overline{e_i}$ for any $i\in\NN$ under the canonical map $\XX_v \ra
\Ga\bs \XX_v$, in which case we set
$$
\lambda(e)=1+\lfloor\ln(i+1)\rfloor-\lfloor\ln i\rfloor
$$ 
if $i\geq 1$ and $\lambda(e)=1$ if $i=0$ (where $\lfloor \cdot\rfloor$
is the largest previous integer map). This map $\lambda$ now has
values in $\{1,2\}$, and we subdivide each edge of length $2$ into two
edges of length $1$. The tree $\YY$ thus obtained has uniformly
bounded degrees (although it is no longer a uniform tree), and the
group $\Ga'$ defines a geometrically finite discrete subgroup of
$\Aut(\YY)$ with infinite Bowen-Margulis measure.  To see that $P_0$
is again of convergence type, observe that each distance in the
simplicial case appearing in the Poincar\'e series differs by at most
$1$ from the corresponding distance in the metric tree case.  The
remainder of the argument is the same as in the metric tree case.
\cqfd

\chapter{Equidistribution and  counting of rational points in 
completed function fields}
\label{sec:rattionalequid}

Let $K$ be a (global) function field over $\FF_q$ of genus $g$, let
$v$ be a (normalised discrete) valuation of $K$, let $K_v$ be the
associated completion of $K$ and let $R_v$ be the affine function ring
associated with $v$.\footnote{See Section \ref{subsec:valuedfields}
  for definitions and background.}  In this Chapter, we prove
analogues of the classical results on the counting and
equidistribution towards the Lebesgue measure on $\RR$ of the Farey
fractions $\frac pq$ with $(p,q)\in\ZZ\times(\ZZ-\{0\})$ relatively
prime.\footnote{See for instance \cite{Neville49}, as well as
  \cite{ParPau14AFST} for an approach using methods similar to the
  ones in this text.}  In particular, we prove various
equidistribution results of locally finite families of elements of $K$
towards the Haar measure on $K_v$, using the geometrical work on
equidistribution of common perpendiculars done in Section
\ref{subsec:equidcommperpdiscrtime} and recalled in Section
\ref{subsec:locconst}.

\section{Counting and equidistribution of non-Archimedian Farey fractions}
\label{subsec:mertens}

The first result of this Section is an analog in function fields of
the equidistribution of Farey fractions to the Lebesgue measure in
$\RR$, see the Introduction, and for example
\cite[p.~978]{ParPau14AFST} for the precise statement and a
geometric proof.  For every $(x_0,y_0)\in R_v\times R_v- \{(0,0)\}$,
let
$$
m_{v,\,x_0,\,y_0}=\card\{a\in (R_v)^\times\,:\;
\exists\;b\in x_0 R_v\cap y_0R_v,\,(a-1)x_0y_0-bx_0\in y_0^2R_v\}\,.
$$
For future use, note that by Equation \eqref{eq:inversRv}
\begin{equation}\label{eq:mvxoyocasspecial}
m_{v,\,1,\,0}=q-1\;.
\end{equation}
For every $(a,b)\in R_v\times R_v$ and every subgroup $H$ of
$\GL_2(R_v)$, let $H_{(a,b)}$ be the stabiliser of $(a,b)$ for the
linear action of $H$ on $R_v\times R_v$.

\btheo\label{theo:Mertensfunctionfieldfissgp} 
Let $G$ be a finite index subgroup of $\GL_2(R_v)$, and let
$(x_0,y_0)\in R_v\times R_v- \{(0,0)\}$. Let
$$
c=\frac{({q_v}^2-1)\;(q_v+1)\;\zeta_K(-1)\;m_{v,\,x_0,\,y_0}\; (N\langle
  x_0,y_0\rangle)^2\;[\GL_2(R_v):G]}
{(q-1)\;q^{\,g-1}\;q_v^3\;
  [\GL_2(R_v)_{(x_0,y_0)}:G_{(x_0,y_0)}]}\;.
$$
Then, as $s\to+\infty$, 
$$
c\;s^{-2} \sum_{(x,\,y)\in G(x_0,\,y_0),\; |y|_v\leq s}
\Delta_{\frac xy}\;\;\weakstar\;\; \haar_{K_v}\;.
$$
\etheo

For every $\beta\in\;]0,\frac{1}{\ln q_v}]$, there exists $\kappa>0$
such that for every $\beta$-Hölder-continuous function $\psi:
K_v \ra \RR$ with compact support,\footnote{where $K_v$ is endowed
with the distance $(x,y)\mapsto |x-y|_v$} as for instance if
$\psi:K_v\ra\RR$ is locally constant with compact support (see
Remark \ref{rem:locconstholder}), there is an error term in the
equidistribution claim of Theorem
\ref{theo:Mertensfunctionfieldfissgp} evaluated on $\psi$, of the form
$\bigO(s^{-\kappa}\|\psi\|_{\beta})$.

It is remarkable that due to the general nature of our geometrical
tools, we are able to work with any finite index subgroup $G$ of
$\GL_2(R_v)$, and not only with its congruence subgroups. In this
generality, the usual techniques (for instance involving analysis of
Eisenstein series) are not likely to apply. Also note that the
H\"older regularity for the error term is a much weaker assumption
than the locally constant one that is usually obtained by analytic
number theory methods.

Theorem \ref{theo:mertensfunctionfieldintro} in the Introduction
follows from this result, by taking $K=\FF_q(Y)$ (so that $g=0$),
$v=v_\infty$, $q_v=q$, $(x_0,y_0)=(1,0)$ and $s=q^t$, and by using
Equations \eqref{zetamoinsun} and \eqref{eq:mvxoyocasspecial} in order
to simplify the constant.

\medskip
Before proving Theorem \ref{theo:Mertensfunctionfieldfissgp}, let us
give a counting result which follows from this equidistribution result.

The additive group $R_v$ acts on $R_v\times R_v$ by the horizontal
shears (transvections): 
$$
\forall\;z\in R_v,\;\forall\; (x,y)\in R_v\times R_v,\;\;\;\; 
z\cdot (x,y)=(x+zy,y)\;,
$$ 
and this action preserves the absolute value $|y|_v$ of the
vertical coordinate $y$. Let $G$ be a finite index subgroup of
$\GL_2(R_v)$, and let $x_0,y_0\in R_v$. Let $R_{v,\,G}$ be the finite
index additive subgroup of $R_v$ consisting of the elements $x\in R_v$
such that $\begin{pmatrix} 1 & x \\ 0 & 1\end{pmatrix}\in G$, acting
by horizontal shears on $G(x_0,\,y_0)$. Note that $R_{v,\,G}=R_v$ if
$G= \GL_2(R_v)$. We may then define a {\it counting
  function}\index{counting function} $\Psi_{G,\,x_0,y_0}$ of the elements
of $K$ in an orbit by homographies under $G$, as
$$ 
\Psi_{G,\,x_0,\,y_0}(s)=\card\;\;R_{v,\,G}\,\bs 
\big\{(x,\,y)\in G(x_0,\,y_0),\;\; |y|_v\leq s\}\;.
$$

\bcoro \label{coro:mertenscount} Let $G$ be a finite index subgroup of
$\GL_2(R_v)$, and let $(x_0,y_0)\in R_v\times R_v-\{(0,0)\}$. Then
there exists $\kappa>0$ such that, as $s\to+\infty$,
\begin{align*}
& \Psi_{G,\,x_0,y_0}(s)\\=\; & \frac{(q-1)\;q^{\,2g-2}\;q_v^3\; 
[\GL_2(R_v)_{(x_0,y_0)}:G_{(x_0,y_0)}]\;[R_v:R_{v,\,G}]}
{({q_v}^2-1)\;(q_v+1)\;\zeta_K(-1)\;m_{v,\;x_0,\,y_0}\;
(N\langle x_0,y_0\rangle)^2\;[\GL_2(R_v):G]}\;s^{2} 
+\bigO(s^{2-\kappa})\;.
\end{align*}
\ecoro

\dem This follows from Theorem \ref{theo:Mertensfunctionfieldfissgp},
by considering the locally constant characteristic function of a
closed and open fundamental domain of $K_v$ modulo the action by
translations of $R_{v,\,G}$, and by using Lemma \ref{lem:covolideal}
with $I=R_v$. 
\cqfd

\medskip 
Let us fix some notation for this Section.  For every subgroup $H$ of
$\GL_2(R_v)$, we denote by $\overline{H}$ its image in $\Ga_v=
\PGL_2(R_v)$. Let $\XX_v$ be the Bruhat-Tits tree\footnote{See Section
  \ref{subsec:BruhatTitstrees}.} of $(\PGL_2,K_v)$, which is regular of
degree $q_v+1$.  Let
$$
r=\frac{x_0}{y_0}\in K\cup\{\infty\}\,.
$$ 
If $y_0= 0$, let $g_r=\id\in\GL_2(K)$, and if $y_0\neq 0$, let
$$
g_r=\begin{pmatrix} r& 1\\1& 0\end{pmatrix}\in\GL_2(K)\;.
$$

\medskip
\noindent
{\bf Proof of Theorem \ref{theo:Mertensfunctionfieldfissgp}. }  We
apply Theorem \ref{theo:algebrogeometricequid} with
$\Ga= \overline{G}$, $\DD^-= \H_\infty$ and
$\DD^+=g_r\H_\infty$. Recall that $\H_\infty$ is the horoball in
$\XX_v$ centred at $\infty$ whose boundary contains $*_v$ (see
Section \ref{subsec:BruhatTitstrees}).

Note that $\Ga$ has finite index in $\Ga_v$ and, in particular, it is
a lattice of $\XX_v$.  By \cite[II.1.2, Coro.]{Serre83}, for all $x\in
V\XX_v$ and $\ga\in \GL_2(R_v)$, the distance $d(x,\ga x)$ is even
since $v(\det \ga)=0$. Hence by the equivalence in Equation
\eqref{eq:equivLamGaevendisteven}, the length spectrum $L_{\Ga_v}$ of
$\Ga_v$ is $2\ZZ$. The length spectrum of $\Ga$ is also $2\ZZ$, since
it is contained in $L_{\Ga_v}$.

Note that $\DD^+$ is a horoball in $\XX_v$ centred at
$r=\frac{x_0}{y_0} \in \partial_\infty \XX_v$, by Lemma
\ref{lem:imagehoroinftyparpglderv}.  The stabiliser $\Ga_{\DD^-}$ of
$\DD^-$ (respectively $\Ga_{\DD^+}$ of $\DD^+$) coincides with the
stabiliser $\Ga_\infty$ of $\infty \in \partial_\infty \XX_v$
(respectively the stabiliser $\Ga_r$ of $r$) in $\Ga$.  Note that the
families $\D^\pm= (\ga\DD^\pm)_{\ga\in\Ga/\Ga_{\DD^\pm}}$ are locally
finite, since $\Ga_v$, and hence its finite index subgroup $\Ga$, is
geometrically finite,\footnote{See Section \ref{subsec:modulargroup}.}
and since $\infty$ and $r\in K$ are bounded parabolic limit points of
$\Ga_v$, hence of its finite index subgroup $\Ga$.

For every $\ga\in\Ga/\Ga_r$ such that $\DD^-$ and $\ga\DD^+$ are
disjoint, let $\alpha^-_{e,\,\ga}$ be the generalised geodesic line,
isometric exactly on $[0,d(\DD^-, \, \ga \DD^+)]$, whose image is the
common perpendicular between $\DD^-$ and $\ga \DD^+$, and let
$\rho_\ga$ be the geodesic ray starting (at time $t=0$) from
$\alpha^-_{e,\,\ga} (0)$ and ending at the point at infinity $\ga\cdot
r$ of $\ga\DD^+$.  Note that $\rho_\ga$ and $\alpha^-_{e,\,\ga}$
coincide on $[0, d(\DD^-,\ga\DD^+)]$.

Since the Patterson densities of lattices of $\XX_v$ have total mass
$\frac{q_v+1}{q_v}$ by Proposition \ref{prop:mescomputBT} (2), they
are normalised as in Theorem \ref{theo:algebrogeometricequid}. Then by
Equation \eqref{eq:extensgeodraytoinfty}, we have
\begin{equation}\label{eq:applitwosubtreevenMert1}
\lim_{n\ra+\infty} \;\frac{({q_v}^2-1)(q_v+1)}{2\;q_v^3}\;
\frac{\Vol(\Ga\dbs\XX_v)}{\|\sigma^-_{\D^+}\|}\;{q_v}^{-n}
\sum_{\substack{\ga\in \Ga/\Ga_{r}\\
0< d(\DD^-, \, \ga \DD^+)\leq n}}
\Dirac_{(\rho_\ga)_+} \;=\;
(\partial^+)_*\wt\sigma^+_{\DD^-}\;.
\end{equation}
Furthermore, for every $\beta\in\;]0,\frac{1}{\ln q_v}]$, by the
comment following Equation \eqref{eq:extensgeodraytoinfty}, we have an
error term of the form $\bigO(e^{-\kappa\,n}\,\|\psi\|_{\beta\ln q_v})$
for some $\kappa>0$ in the above formula when evaluated on
$\psi\in\C^{\beta\ln q_v}_c(\partial_\infty \XX_v-\{\infty\})$, where
$\partial_\infty \XX_v-\{\infty\}$ is endowed with Hamenstädt's
distance\footnote{See Equation \eqref{eq:defidisthamenbord}, and note
  that with $D=\H_\infty$, Hamenstädt's distance $d_{\H_\infty}$
  coincides with the distance-like function $d_D$ introduced in
  Equation \eqref{eq:distancelike} and used in Remark \ref{rem:utidD},
  since $D$ is a horoball.}  $d_{\H_\infty}$.  Hence we have an error
term $\bigO(e^{-\kappa\,n}\,\|\psi\|_{\beta} )$ for some $\kappa>0$ in
the above formula when evaluated on $\psi\in\C^{\beta}_c(K_v)$, where
$K_v=\partial_\infty \XX_v -\{\infty\}$ is endowed with the distance
$(x,y)\mapsto |x-y|_v$, see Equation
\eqref{eq:relatholdernormhamabsval}.

By Proposition \ref{prop:mescomputBT} (3), we have
$$
(\partial^+)_*\wt\sigma^+_{\DD^-} = \haar_{K_v}\;.
$$
Hence Equation \eqref{eq:applitwosubtreevenMert1}
gives, with the appropriate error term,
\begin{equation}\label{eq:applitwosubtreevenMert2}
\lim_{n\ra+\infty} \;\frac{({q_v}^2-1)(q_v+1)}{2\,q_v^3}\;
\frac{\Vol(\Ga\dbs\XX_v)}{\|\sigma^-_{\D^+}\|}\;{q_v}^{-n}
\sum_{\substack{\ga\in \Ga/\Ga_{r}\\
0< d(\DD^-, \, \ga \DD^+)\leq n}}
\Dirac_{\ga \cdot r} \;=\;
\haar_{K_v}\;.
\end{equation}

Let $g\in \GL_2(K)$ be such that $g\infty\neq \infty$. This condition
is equivalent to asking that the $(2,1)$-entry $c=c(g)$ of $g$ is
nonzero. By Lemma \ref{lem:imagehoroinftyparpglderv}, the signed
distance between the horospheres $\H_\infty$ and $g\H_\infty$ is
\begin{equation}\label{eq:signdisthoroinfty}
d(\H_\infty, g\H_\infty) = -2\;v(c) = 2\;\frac{\ln |c|_v}{\ln q_v}\;.
\end{equation}
If $y_0\ne 0$, then $(x,y) = g(x_0,y_0)$ if and only if $(\frac
x{y_0}, \frac y{y_0})= gg_r (1,0)$, and the $(2,1)$-entry of
$gg_r$ is $\frac y{y_0}$. If $y_0= 0$ (which implies that
$g_r=\id$ and $x_0 \neq 0$), then $(x,y) = g(x_0,y_0)$ if and only
if $(\frac x{x_0}, \frac y{x_0})= g (1,0)$, and the $(2,1)$-entry of
$g=gg_r$ is $\frac y{x_0}$.  Let 
$$
z_0=\left\{\begin{array}{ll} y_0 & {\rm if~} y_0\neq 0\\x_0& 
{\rm otherwise.}\end{array}\right.
$$
By Equation \eqref{eq:signdisthoroinfty}, the signed distance
between $\DD^-= \H_\infty$ and $g\DD^+= g\,g_r \H_\infty$ is
$$
d(\DD^-,g\DD^+)=
\frac{2}{\ln q_v}\; \ln\Big|\frac{y}{z_0}\Big|_v\;.
$$

By discreteness, there are only finitely many double classes $[g]\in
G_{(1,0)}\bs G/ G_{(x_0,y_0)}$ such that $\DD^-= \H_\infty$ and
$g\DD^+= g\,g_r\H_\infty$ are not disjoint. Let $Z(G)$ be the
centre of $G$, which is finite. Since $Z(G)$ acts trivially on
$\PP_1(K_v)$, the map $G/G_{(x_0,y_0)}\ra \Ga/\Ga_{r}$ induced by the
canonical map $\GL_2(R_v)\ra \PGL_2(R_v)$ is $|Z(G)|$-to-$1$.
Using the change of variable
$$
s =|z_0|_v\;{q_v}^{\frac{n}{2}}\;,
$$
so that ${q_v}^{-n}={|z_0|_v}^2\;s^{-2}$, Equation
\eqref{eq:applitwosubtreevenMert2} gives, with the
appropriate error term,
\begin{align}
\lim_{s\ra+\infty} 
\;\frac{({q_v}^2-1)\;(q_v+1)\;{|z_0|_v}^2}{2\;{q_v}^{3}\;|Z(G)|}\;
\frac{\Vol(\Ga\dbs\XX_v)}{\|\sigma^-_{\D^+}\|} &\;s^{-2}
\sum_{(x,\,y)\in G(x_0,\,y_0),\;\; |y|_v\leq s}
\Dirac_{\frac xy} \nonumber\\ & =\;
\haar_{K_v}\;.\label{eq:applitwosubtreevenMert3}
\end{align}

\medskip
The order of the centre $Z(GL_2(R_v))=(R_v)^\times\id$ is $q-1$ by
Equation \eqref{eq:inversRv}. The map $\GL_2(R_v)/G\ra \Ga_v/\Ga$
induced by the canonical map $\GL_2(R_v)\ra \PGL_2(R_v)$ is hence
$\frac{q-1}{|Z(G)|}$-to-$1$. By Equation \eqref{eq:covolPGLdeR}, we
hence have
\begin{align}
\Vol(\Ga\dbs\XX_v)&=[\Ga_v:\Ga]
\Vol(\Ga_v\dbs\XX_v)=2\;\zeta_K(-1)\,[\Ga_v:\Ga]\nonumber\\ 
&= \frac{2}{q-1}\;\zeta_K(-1)\;|Z(G)|\;[\GL_2(R_v):G]\;.
\label{eq:covoltotMertens}
\end{align}

\medskip Theorem \ref{theo:Mertensfunctionfieldfissgp} follows from
Equations \eqref{eq:applitwosubtreevenMert3} and
\eqref{eq:covoltotMertens} and from Lemma \ref{lem:totmassskinhoro}
below. \cqfd

\blemm \label{lem:totmassskinhoro} We have
$$
\|\sigma^-_{\D^+}\|=\frac{q^{\,g-1}\;{|z_0|_v}^2\;
[\GL_2(R_v)_{(x_0,y_0)}:G_{(x_0,y_0)}]}{m_{v,\,x_0,\,y_0}\;
(N\langle x_0,y_0\rangle)^2}\;.
$$
\elemm

\dem 
Let $\ga_r$ be the image of $g_r$ in $\PGL_2(K)$.  Let us define
$\Ga'={\ga_r}^{-1}\Ga\ga_r$, which is a finite index subgroup in
$\Ga'_v={\ga_r}^{-1}\Ga_v\ga_r$ and a lattice of $\XX_v$. Since
$\ga_r$ maps $\infty$ to $r$, the point $\infty$ is a bounded
parabolic limit point of $\Ga'$, and we have $(\Ga')_\infty=
{\ga_r}^{-1}\Ga_r\ga_r$. Since the canonical map $\GL_2(R_v)\ra
\PGL_2(R_v)$ is injective on the stabiliser $\GL_2(R_v)_{(x_0,y_0)}$,
we have
$$
[(\Ga'_v)_\infty:(\Ga')_\infty]=[(\Ga_v)_r:\Ga_r]=
[\GL_2(R_v)_{(x_0,y_0)}:G_{(x_0,y_0)}] \;.
$$

Since the Patterson density of a lattice does not depend on the
lattice (see Proposition \ref{prop:uniflatmBMfinie} (1)), the skinning
measures $\wt \sigma^\pm_\H$ of a given horoball $\H$ do not depend on
the lattice. Thus
$$
\ga_*\;\wt \sigma^\pm_\H=\wt\sigma^\pm_{\ga\H}
$$
for every $\ga\in\Aut(\XX_v)$. Let $\D_1^+=
(\ga'\H_\infty)_{\ga'\in\Ga'/\Ga'_\infty}$, which is a locally finite
$\Ga'$-equivariant family of horoballs.  We hence have, using
Proposition \ref{prop:mescomputBT} (4) for the third equality,
\begin{align}
\|\sigma^-_{\D^+}\|&=\|\sigma^-_{\ga_r\D_1^+}\|=\|\sigma^-_{\D_1^+}\|
=\haar_{K_v}((\Ga')_\infty\bs K_v)\nonumber\\ &
=[(\Ga'_v)_\infty:(\Ga')_\infty]\;\haar_{K_v}((\Ga'_v)_\infty\bs K_v)
\nonumber\\ &=[\GL_2(R_v)_{(x_0,y_0)}:G_{(x_0,y_0)}]\;
\haar_{K_v}((\Ga'_v)_\infty\bs K_v)\;.\label{eq:lemmtotmassskinhoro2}
\end{align}

Every element in the stabiliser of $\infty$ in $\PGL_2(K_v)$ can be
uniquely written in the form 
$\alpha=\begin{bmatrix} a & b\\0 & 1\end{bmatrix}$
with $(a,b)\in (K_v)^\times\times K_v$. Note that
$$
\begin{pmatrix} r & 1\\1 & 0\end{pmatrix}
\begin{pmatrix} a & b\\0 & 1\end{pmatrix}
\begin{pmatrix} 0 & 1\\1 & -r\end{pmatrix}
=\begin{pmatrix} br+1 & ar-br^2-r\\b & a-br\end{pmatrix}\;.
$$
When $x_0= 0$ or $y_0= 0$, we have $\alpha\in \Ga'_v$  if and only if 
$$
b\in R_v\;\;\;{\rm and}\;\;\; a\in (R_v)^\times\;.
$$ 
When $x_0,y_0\neq 0$, we have $\alpha\in \Ga'_v$ if and only if
$\ga_r\alpha{\ga_r}^{-1}\in \Ga_v$, hence if and only if
$$
b\in R_v\cap \frac{1}{r} R_v,\;\;\;\; a\in (R_v)^\times,\;\;\;\; 
ar-br^2-r\in R_v\;.
$$

Let $U'_\infty$ be the kernel of the group morphism from $(\Ga'_v)_\infty$ to
$(K_v)^\times$ sending $\begin{bmatrix} a & b\\0 & 1\end{bmatrix}$ to
$a$, and let $m_v$ be its index in $(\Ga'_v)_\infty$. If $x_0= 0$ or
$y_0= 0$, then $m_v$ is equal to $|(R_v)^\times|$, so that, by
Equation \eqref{eq:inversRv},
$$
m_v=|(R_v)^\times| =|(\FF_q)^\times| =q-1\;.
$$
If $x_0,y_0\neq 0$, we have
$$
m_v=\card\;\{a\in(R_v)^\times\;:\;\exists\;
b\in R_v\cap \frac{1}{r} R_v,\;\; ar-br^2-r\in R_v\}\;.
$$
Note that the notation $m_v$ coincides with the constant
$m_{v,\,x_0,\,y_0}$ defined before the statement of Theorem
\ref{theo:Mertensfunctionfieldfissgp} in both cases.

If $I_{(x_0,y_0)}$ is the nonzero fractional ideal
$$
I_{(x_0,y_0)}=\begin{cases}
R_v & {\rm if}\; x_0= 0\;{\rm or}\; y_0= 0\,,\\
R_v\cap \frac{1}{r} R_v\cap \frac{1}{r^2} R_v & {\rm otherwise,}
\end{cases}
$$
then 
$$
U'_\infty=\bigg\{\begin{bmatrix} 1 & b\\0 &
  1\end{bmatrix}\;:\;b\in I_{(x_0,y_0)}\bigg\}\,.
$$ 
Therefore by Lemma \ref{lem:covolideal},
\begin{equation}\label{eq:lemmtotmassskinhoro3}
\haar_{K_v}((\Ga'_v)_\infty\bs K_v)=
\frac{\haar_{K_v}(I_{(x_0,y_0)}\bs K_v)}{[(\Ga'_v)_\infty:U'_\infty]}
=\frac{q^{\,g-1}\;N(I_{(x_0,y_0)})}{m_v}\;.
\end{equation}
Let $(x_0)=\prod_{\ppp} \ppp^{\nu_\ppp(x_0)}$ and $(y_0)= \prod_{\ppp}
\ppp^{\nu_\ppp(y_0)}$ be the prime decompositions of the principal
ideals $(x_0)$ and $(y_0)$. By the formulas of the prime
decompositions of intersections, sums and products of ideals in
Dedekind rings (see for instance \cite[\S 1.1]{Narkiewicz04} and
Equation \eqref{eq:interidealdedeking}), we have
$$
(x_0^2)\cap (x_0y_0) \cap (y_0^2)=(x_0^2)\cap(y_0^2)= 
\prod_{\ppp} \ppp^{2\max\{\nu_\ppp(x_0),\, \nu_\ppp(y_0)\}}
$$ 
and 
$$\langle x_0,y_0 \rangle = \prod_{\ppp} \ppp^{\min\{\nu_\ppp(x_0),
  \, \nu_\ppp(y_0)\}}\,.
$$
By the definition of the ideal $I_{(x_0,y_0)}$, by the
multiplicativity of the norm, and by Equation
\eqref{eq:normprincipideal}, we hence have if $x_0\neq 0$ and $y_0\neq
0$
\begin{equation}\label{eq:lemmtotmassskinhoro4}
\frac{N(I_{(x_0,y_0)})\; (N\langle x_0,y_0\rangle)^2}{{|y_0|_v}^2}=
N\Big(\big((x_0^2)\cap(x_0y_0)\cap (y_0^2)\big)\langle
x_0,y_0\rangle^2(x_0)^{-2}(y_0)^{-2}\Big)=1\;.
\end{equation}
If $x_0=0$ or $y_0= 0$, then
\begin{equation}\label{eq:lemmtotmassskinhoro5}
N(I_{(x_0,y_0)})= N(R_v)=1\;.
\end{equation}

Lemma \ref{lem:totmassskinhoro} follows from Equations
\eqref{eq:lemmtotmassskinhoro2}, \eqref{eq:lemmtotmassskinhoro3} and
\eqref{eq:lemmtotmassskinhoro4} if $x_0\neq 0$ and $y_0\neq 0$ or
\eqref{eq:lemmtotmassskinhoro5} if $x_0=0$ or $y_0=0$.  
\cqfd

\bigskip 
Let us state one particular case of Theorem
\ref{theo:Mertensfunctionfieldfissgp} in an arithmetic setting, using
a congruence sugbroup.

\btheo\label{theo:Mertensfunctionfield} Let $I$ be a nonzero ideal of
$R_v$. Then as $t\to+\infty$, we have
\begin{align*}
\frac{(q_v^{\;2}-1)\;(q_v+1)\;\zeta_K(-1)\;N(I)
\prod_{\ppp|I}(1+\frac{1}{N(\ppp)})}{q^{\,g-1}\;q_v^{\;3}}&
\;(q_v)^{-2t} \sum_{\substack{(x,y)\in R_v\times I\\ \langle
    x,\,y\rangle=R_v,\;v(y)\geq -t}} \Delta_{\frac xy} 
 \\ &\weakstar\;\; \haar_{K_v}\;,
\end{align*}
where the product ranges over the prime factors $\ppp$ of the ideal $I$.
Furthermore, if
$$
\Psi(t)=\card\;\;R_v\,\bs\big\{(x,y)\in R_v\times
 I\;:\; \langle x,y\rangle = R_v,\;v(y)\geq -t\big\}\,,
$$ 
then there exists $\kappa>0$ such that, as $t\ra +\infty$,
$$
\Psi(t)=
\frac{q^{\,2g-2}\;q_v^{\;3}}{({q_v}^2-1)\;(q_v+1)\;\zeta_K(-1)\;N(I)
\prod_{\ppp|I}(1+\frac{1}{N(\ppp)})}
\;{q_v}^{2t}+\bigO({q_v}^{(2-\kappa)t})\;.
$$
\etheo

For every $\beta\in\;]0,\frac{1}{\ln q_v}]$, there exists $\kappa>0$
such that for every $\psi\in\C^\beta_c(K_v)$ there is an error term in
the above equidistribution claim evaluated on $\psi$, of the form
$\bigO(q_v^{\;-\kappa\, t}\,\|\psi\|_{\beta})$.

\medskip \dem The counting claim is deduced from the equidistribution
claim in the same way that Corollary \ref{coro:mertenscount} is
deduced from Theorem \ref{theo:Mertensfunctionfieldfissgp}, noting
that the action of $R_v$ by horizontal shears preserves $R_v\times I$.

\medskip 
In order to prove the equidistribution claim, we apply
Theorem \ref{theo:Mertensfunctionfieldfissgp} with $(x_0,y_0)=(1,0)$
and with $G$ the Hecke congruence subgroup
\begin{equation}\label{eq:heckecongruencesubgroup}
G_I=\bigg\{\begin{pmatrix} a&b\\c&d\end{pmatrix}\in\GL_2(R_v)\;:\;
c\in I\bigg\}\,,
\end{equation}
which is the preimage of the upper triangular subgroup of
$\GL_2(R_v/I)$ by reduction modulo $I$.  In this case, the constant
$m_{v,x_0,y_0}$ appearing in the statement of Theorem
\ref{theo:Mertensfunctionfieldfissgp} is equal to $q-1$ by Equation
\eqref{eq:mvxoyocasspecial}. The group $G_I$ has finite index in
$\GL_2(R_v)$. The following result is well-known to arithmetic readers
(see for instance \cite[page 24]{Shimura71} when $R_v$ is replaced by
$\ZZ$), we only give a sketch of proof (indicated to us by J.-B.~Bost)
for the sake of the geometer readers.

\blemm \label{lem:indexHeckesubgrou} We have
$$
[\GL_2(R_v):G_I]=N(I)\prod_{\ppp|I}(1+\frac{1}{N(\ppp)})\;,
$$ 
where the product ranges over the prime factors $\ppp$ of the ideal
$I$.  
\elemm

\dem 
Recall that we denote by $|E|$ the cardinality of a finite set
$E$. For every commutative ring $A$ with finite group of invertible
elements $A^\times$, we have a disjoint union
$$
\GL_2(A)=\bigcup_{a\in A^\times}
\begin{pmatrix} a & 0 \\ 0 & 1\end{pmatrix} \SL_2(A)\;.
$$
Hence $[\GL_2(A):\SL_2(A)]=|A^\times|$. Since $\begin{pmatrix} a & 0
  \\ 0 & 1\end{pmatrix}$ belongs to $G_I$ for all $a\in(R_v)^\times$,
we have 
$$
G_I=\bigcup_{a\in (R_v)^\times}
\begin{pmatrix} a & 0 \\ 0 & 1\end{pmatrix} G_I\cap \SL_2(R_v)\;,
$$
so that $[\GL_2(R_v):G_I]=[\SL_2(R_v):G_I\cap \SL_2(R_v)]$.

The group morphism of reduction modulo $I$ from $\SL_2(R_v)$ to
$\SL_2(R_v/I)$ is onto, by an argument of further reduction to the
various prime power factors of $I$ and of lifting elementary matrices.
The order of the upper triangular subgroup of $\SL_2(R_v/I)$ is
$|(R_v/I)^\times|\;|R/I|$, where $(R_v/I)^\times$ is the group of
invertible elements of the ring $R_v/I$ (that we will see again
below). Hence
\begin{align}
[\GL_2(R_v):G_I]&=[\SL_2(R_v):G_I\cap \SL_2(R_v)]=
\frac{|\SL_2(R_v/I)|}{|(R_v/I)^\times|\;|R/I|}\nonumber\\ &
=\frac{|\GL_2(R_v/I)|}{|(R_v/I)^\times|^2|R/I|}\;.
\label{eq:premcalcindHeckegroup}
\end{align}
By the multiplicativity of the norm and by the Chinese remainder
theorem,\footnote{saying that the rings $R_v/I$ and $\prod_{\ppp}
  R_v/\ppp^{v_\ppp(I)}$ are isomorphic, see for instance \cite[page
  11]{Narkiewicz04}} one reduces the result to the case when
$I=\ppp^n$ is the $n$-th power of a fixed prime ideal $\ppp$ with norm
$N(\ppp)=N$, where $n\in\NN$. Note that since $R_v/\ppp$ is a field,
we have
$$
|\GL_1(R_v/\ppp)|=|(R_v/\ppp)^\times|=|R_v/\ppp|-1=N-1
$$
and
$$
|\GL_2(R_v/\ppp)|= (|R_v/\ppp|^2-1)(|R_v/\ppp|^2-|R_v/\ppp|)=
N^2(N-1)^2\big(1+\frac{1}{N}\big)\;.
$$
For $k=1$ or $k=2$, the kernel of the morphism of reduction modulo
$\ppp$ from $\GL_k(R_v/I)=\GL_k(R_v/\ppp^n)$ to $\GL_k(R_v/\ppp)$ has
order $N^{k^2(n-1)}$. Hence
$$
|\GL_2(R_v/I)|= N^{4(n-1)}N^2(N-1)^2\big(1+\frac{1}{N}\big)\;,
$$
and
$$
|(R_v/I)^\times|= N^{n-1}(N-1)\;.
$$
Therefore, by Equation \eqref{eq:premcalcindHeckegroup}, we have
$$
[\GL_2(R_v):G_I]=\frac{N^{4(n-1)}N^2(N-1)^2\big(1+\frac{1}{N}\big)}
{N^{2(n-1)}(N-1)^2N^n}= N^n\big(1+\frac{1}{N}\big)\;.
$$
This proves the result.   \cqfd

\bigskip 
We can now conclude the proof of Theorem
\ref{theo:Mertensfunctionfield}. Note that $\GL_2(R_v)_{(1,0)}=
(G_I)_{(1,0)}$. The result then follows from Theorem
\ref{theo:Mertensfunctionfieldfissgp} and its Corollary
\ref{coro:mertenscount}, using the change of variables $s=(q_v)^t$,
since
$$
G_I(1,0)=\{(x,y)\in R_v\times  I\;:\; 
\langle x,y\rangle = R_v\}\;.\;\;\;\Box
$$

\medskip 
The following result is a particular case of Theorem
\ref{theo:Mertensfunctionfield}.

\bcoro\label{coro:Mertensfractions} Let $P_0$ be a nonzero element of
the polynomial ring $R=\FF_q[Y]$ over $\FF_q$, and let $P_0=
a_0\prod_{i=1}^k (P_i)^{n_i}$ be the prime decomposition of $P_0$.
Then as $t\to+\infty$,
$$
\frac{(q+1)\prod_{i=1}^kq^{n_i\deg P_i}(1+q^{-\deg P_i})}{(q-1)\,q^2}
\;q^{-2t} \sum_{\substack{(P,Q)\in R\times (P_0R)\\PR+QR=R,\;\deg Q\le  t}}
\Delta_{\frac PQ}\;\;\weakstar\;\;\haar_{\FF_q((Y^{-1}))}\;.
$$
\ecoro

For every $\beta\in\;]0,\frac{1}{\ln q_v}]$, there exists $\kappa>0$
such that for every $\psi\in\C^\beta_c\big(\FF_q((Y^{-1}))\big)$ there is an
error term in the above equidistribution claim evaluated on $\psi$, of
the form $\bigO(q^{\;-\kappa\,t}\,\|\psi\|_{\beta})$.

\medskip
\dem In this statement, we use the standard convention that $k=0$ if
$P_0$ is constant, $a_0\in(\FF_q)^\times$ and $P_i\in R$ is monic.

We apply the first claim of Theorem \ref{theo:Mertensfunctionfield}
with $K=\FF_q(T)$ and $v=v_\infty$ so that $g=0$, $q_v=q$ and $R_v=R$,
and with $I= P_0R$, so that $N(I)=\prod_{i=1}^kq^{n_i\deg P_i}$. The
result follows from Equation \eqref{zetamoinsun}.  
\cqfd

\section{Mertens's formula in function fields}
\label{subsec:mertensinfunctionfields}

In this Section, we recover the function field analogue of Mertens's
classical formula on the average order of the Euler function. We begin
with a more general counting and equidistribution result.

Let $\mmm$ be a (nonzero) fractional ideal of $R_v$, with norm
$N(\mmm)$. Note that the action of the additive group $R_v$ on
$K_v\times K_v$ by the horizontal shears $z\cdot(x,y)=(x+zy,y)$
preserves $\mmm\times\mmm$.
We consider the counting function $\psi_{\mmm}:[0,+\infty[\;\ra \NN$
defined by
$$
\psi_{\mmm}(s)=\card\big(R_v\,\bs
\big\{(x,y)\in\mmm\times\mmm\;:\;0<N(\mmm)^{-1}N(y)\leq s,
\; \langle x,y\rangle=\mmm\big\}\big)\,.
$$ 
Note that $\psi_\mmm$ depends only on the ideal class of $\mmm$ and
thus we can assume in the computations that $\mmm$ is integral, that
is, contained in $R_v$.

\bcoro\label{coro:appliEulerFunctField} There exists $\kappa>0$ such
that, as $s\to+\infty$, 
$$
\psi_{\mmm}(s)= \frac{(q-1)\;q^{\,2g-2}\;{q_v}^3}
{({q_v}^2-1)\;(q_v+1)\;\zeta_K(-1)\;m_{v,\,x_0,\,y_0}\;}\;\;s^2+
\bigO(s^{2-\kappa})\,,
$$
where $\mmm=\langle x_0,y_0\rangle$. Furthermore, as $s\to+\infty$,
$$
\frac{({q_v}^2-1)\;(q_v+1)\;\zeta_K(-1)\;m_{v,\,x_0,\,y_0}}
{(q-1)\;q^{\,g-1}\;{q_v}^3}\;s^{-2} 
\sum_{\substack{(x,y)\in\mmm\times\mmm\;\\
    \;N(\mmm)^{-1}N(y)\leq s, \;  \langle x,y\rangle=\mmm}} 
\Delta_{\frac xy}\;\weakstar\;\haar_{K_v}\;.
$$ 
\ecoro

For every $\beta\in\;]0,\frac{1}{\ln q_v}]$, there exists $\kappa>0$
such that for every $\psi\in\C^\beta_c\big(K_v)\big)$ there is an
error term in the above equidistribution claim evaluated on $\psi$, of
the form $\bigO(s^{-\kappa}\,\|\psi\|_{\beta})$.

\medskip
Theorem \ref{theo:fractidealintro} in the Introduction follows from
this result, by taking $K=\FF_q(Y)$ (so that $g=0$) and $v=v_\infty$
(so that $q_v=q$).  In order to simplify the constant, we use Equation
\eqref{zetamoinsun} and the fact that the ideal class number of $K$,
that equals the number of orbits of $\PGL_2(\FF_q[Y])$ on
$\PP^1(\FF_q(Y))$, is $1$.  Thus, if $\mmm=\langle x_0,y_0\rangle$
then the constant $m_{v,\,x_0,\,y_0}$ is equal to $m_{v,\,1,\,0}$,
which is $q-1$ by Equation \eqref{eq:mvxoyocasspecial}.

\medskip
\dem 
Every nonzero ideal $I$ in $R_v$ is of the form $I=xR_v+yR_v$ for
some $(x,y)\in R_v\times R_v-\{(0,0)\}$, see for instance
\cite[Coro.~5, page 11]{Narkiewicz04}. For all $(x,y)$ and $(z,w)$ in
$R_v\times R_v$, we have $x \,R_v+y\,R_v=z\,R_v+w\,R_v$ if and only if
$(z,w)\in \GL_2(R_v)(x,y)$. The ideal class group of $K$ corresponds
bijectively to the set $\PGL_2(R_v)\bs\PP^1(K)$ of cusps of the
quotient graph of groups $\PGL_2(R_v)\dbs\XX_v$,\footnote{where
  $\XX_v$ is the Bruhat-Tits tree of $(\PGL_2,K_v)\,$} by the map
induced by $I=x\,R_v+y\,R_v \mapsto [x:y]\in \PP^1(K)$.

Given a fixed ideal $\mmm$ in $R_v$, we apply Theorem
\ref{theo:Mertensfunctionfieldfissgp} with $G=\GL_2(R_v)$ and
$(x_0,y_0) \in R_v\times R_v-\{(0,0)\}$ a fixed pair such that
$x_0\,R_v+y_0\,R_v=\mmm$, so that $G(x_0,y_0)=\mmm\times\mmm$. Using
therein the change of variable $s \mapsto N(\mmm) s$ and Equation
\eqref{eq:normprincipideal}, the result follows from Theorem
\ref{theo:Mertensfunctionfieldfissgp} and its Corollary
\ref{coro:mertenscount}.  
\cqfd

\medskip 
As already encountered in the proof of Lemma
\ref{lem:indexHeckesubgrou}, the {\em Euler function}\index{Euler
  function} $\gls{Eulerfunct}$ of $R_v$ is defined on the set of
(nonzero, integral) ideals $I$ of $R_v$ by setting\footnote{See for
  example \cite[\S 1]{Rosen02}.}
$$
\varphi_{R_v}(I)= \card ((R_v/I)^\times)\;,
$$ 
and we denote $\varphi_{R_v}(y)=\varphi_{R_v}(y\,R_v)$ for every
$y\in R_v$.  Thus, by the definition of the action of $R_v$ on
$R_v\times R_v$ by horizontal shears, we have
\begin{align}
\psi_{R_v}(s)&=\sum_{y\in R_v,\;0<N(y)\le s} 
\card\{x\in R_v/yR_v\;:\;\langle x,y\rangle=R_v\}\nonumber\\ &= 
\sum_{y\in R_v,\;0<N(y)\le s} \varphi_{R_v}(y)
\label{eq:relataverEulercounting}\,.
\end{align}

As a particular application of Corollary
\ref{coro:appliEulerFunctField}, we get a well-known asymptotic result
on the number of relatively prime polynomials in $\FF_q[Y]$. The Euler
function of the ring of polynomials $R=\FF_q[Y]$ is then the map
$\phi_q:R-\{0\}\to\NN$ defined by
$$
\phi_q(Q)=\big|\big(R/QR\big)^\times\big|=
\card\big\{P\in R: \langle P,Q\rangle =R,\; \deg P<\deg Q\big\}\,.
$$
Note that $\phi_q(\lambda Q)=\phi_q(Q)$ for every
$\lambda\in(\FF_q)^\times$.

\bcoro [Mertens's formula for polynomials]\label{coro:mertenspolyn} We
have
$$
\lim_{n\to+\infty}\frac {1}{q^{2n}}\sum_{Q\in\FF_q[X],
  \, \deg Q\le n}\phi_q(Q)=\frac{q\,(q-1)}{q+1}\,.
$$
\ecoro

\medskip \dem We apply the first claim of Corollary
\ref{coro:appliEulerFunctField}, in the special case when $K=\FF_q(T)$
and $v=v_\infty$ so that $g=0$, $q_v=q$ and $R_v=R$, and with
$\mmm=R_v$, so that $m_{v,\,x_0,\,y_0}=q-1$, in order to obtain the
asymptotic value of $\psi_{R_v}(s)$ with the change of variable
$s=q^n$. The result follows from Equations 
\eqref{eq:relataverEulercounting} and \eqref{zetamoinsun}.
\cqfd

\medskip 
The above result is an analog of Mertens's formula when $K$ is
replaced by $\QQ$ and $R_v$ by $\ZZ$, see \cite[Theo.~330]{HarWri08}.
See also \cite[Satz 2]{Grotz79}, \cite[\S 4.3]{Cosentino99}, as well
as \cite{ParPau14AFST} and \cite[\S 5]{ParPau16MA} for further
developments.

A much more precise result than Corollary \ref{coro:mertenspolyn} can
be obtained by purely number theo\-retical means as follows.  The
average value of $\phi_q$ is computed in \cite[Prop.~2.7]{Rosen02}:
For $n\ge 1$,
$$
\sum_{\deg f=n, \;f \textrm{ monic}}\phi_q(f)=q^{2n}(1-\frac 1q)\,.
$$
This gives $\sum_{\deg f=n}\phi_q(f)=q^{2n}\frac{(q-1)^2}{q}$, so that
$$
\sum_{0<\deg f\le n}\phi_q(f)=
\sum_{k=1}^{n} q^{2k}\frac{(q-1)^2}{q}=q(q-1)^2\frac{q^{2n}-1}{q^2-1}
=\frac{q(q-1)(q^{2n}-1)}{(q+1)}\,,
$$ 
from which Corollary \ref{coro:mertenspolyn} easily follows.


%% file: meridonIV.tex

\chapter{Equidistribution and counting of quadratic irrational points
in \\ non-Archimedean local fields}
\label{sec:nonarchquadratequid}

Let $K_v$ be a non-Archimedean local field, with valuation $v$,
valuation ring $\OOO_v$, choice of uniformiser $\pi_v$, and residual
field $k_v$ of order $q_v$; let $\XX_v$ be the Bruhat-Tits tree of
$(\PGL_2,K_v)$.\footnote{See Sections \ref{subsec:valuatedlocalfields}
  and \ref{subsec:BruhatTitstrees}.} In this Chapter, we give counting
and equidistribution results in $K_v= \partial_\infty
\XX_v-\{\infty\}$ of an orbit under a lattice of $\PGL_2(K_v)$ of a
fixed point of a loxodromic element of this lattice. We use these
results to deduce equidistribution and counting results of quadratic
irrational elements in non-Archimedean local fields.

When $\XX_v$ is replaced by a real hyperbolic space, or by a more
general simply connected complete Riemannian manifold with negative
sectional curvature, there are numerous quantitative results on the
density of such an orbit, see the works of Patterson, Sullivan, Hill,
Velani, Stratmann, Hersonsky-Paulin, Parkkonen-Paulin. See for
instance \cite{ParPau16LMS} for references.  The arithmetic
applications when $\XX_v$ is replaced by the upper halfspace model of
the real hyperbolic space of dimension $2$, $3$ or $5$ are counting
and equidistribution results of quadratic irrational elements in
$\RR$, $\CC$ and the Hamiltonian quaternions. See for instance
\cite[Coro.~3.10]{ParPau12JMD} and \cite{ParPau14AFST}.

\section{Counting and equidistribution of loxodromic fixed points}
\label{subsec:loxofix}

An element $\ga\in\PGL_2(K_v)$ is said to be {\em
  loxodromic}\index{loxodromic} if it is loxodromic\footnote{See
  Section \ref{subsec:catmoinsun}.} on the (geometric realisation of
the) simplicial tree $\XX_v$.  Its translation length is
$$
\lambda(\ga)= \min_{x\in V\XX_v} d(x,\ga x)>0\,,
$$ 
and the subset
$$ 
\Ax_{\ga}=\{x\in V\XX_v\;:\; d(x,\ga x)= \lambda(\ga)\}
$$ 
is the image of a (discrete) geodesic line in $\XX_v$, which we call
the {\em (discrete) translation axis}\index{translation!axis} of
$\ga$. The points at infinity of $\Ax_{\ga}$ are denoted by $\ga^-$
and $\ga^+$, chosen so that $\ga$ translates away from $\ga^-$ and
towards $\ga^+$ on $\Ax_{\ga}$. Note that for every
$\ga'\in\PGL_2(K_v)$, we have
$$
\ga'\Ax_{\ga}=\Ax_{\ga'\ga\,(\ga')^{-1}}\;\;\;\;{\rm and }\;\;\;\; 
\ga'\,\ga^\pm=(\ga'\ga\,(\ga')^{-1})^\pm\;.
$$

If $\Ga$ is a discrete subgroup of $\PGL_2(K_v)$ and if $\alpha$ is
one of the two fixed points of a loxodromic element of $\Ga$, we
denote the other fixed point of this element by $\gls{galconjloxo}$. Since
$\Ga$ is discrete, the translation axes of two loxodromic elements of
$\Ga$ coincide if they have a common point at infinity. Hence
$\alpha^\sigma$ is uniquely defined. For every $\ga\in\Ga$, we hence have
\begin{equation}\label{eq:gacommutgal}
(\ga\cdot\alpha)^\sigma=\ga\cdot(\alpha^\sigma)\;.
\end{equation}
We define the {\it complexity}\index{complexity} $h(\alpha)$ of the
loxodromic fixed point $\alpha$ by
\begin{equation}\label{eq:complexityloxofix}
\gls{complexityloxo}=\frac{1}{|\alpha-\alpha^\sigma|_v}
\end{equation}
if $\alpha,\alpha^\sigma\neq \infty$, and by $h(\alpha)=0$ if $\alpha$
or $\alpha^\sigma$ is equal to $\infty$.  We define
$\gls{reciprocityindexloxo}\in\{1,2\}$ by $\iota_\alpha=2$ if there
exists an element $\ga\in\Ga$ such that $\ga\cdot\alpha=
\alpha^\sigma$, and $\iota_\alpha=1$ otherwise.\footnote{Recall that
  the groups $\GL_2(K_v)$ and $\PGL_2(K_v)$ act on $\PP^1(K_v)
  =K_v\cup \{\infty\}$ by homographies, and that these actions are
  denoted by $\cdot\;$, see Section \ref{subsec:BruhatTitstrees}.}

Following \cite[II.1.2]{Serre83}, we denote by $\gls{PGLplus}$ the
kernel of the group morphism $\PGL_2(K_v)\ra\ZZ/2\ZZ$ defined by
$\ga=[g]\mapsto v(\det g) \mod 2$. The definition does not depend on
the choice of a representative $g\in \GL_2(K_v)$ of an element $\ga\in
\PGL_2(K_v)$, since 
$$
v\big(\det\begin{pmatrix} \lambda & 0 \\ 0 &
  \lambda \end{pmatrix}\big)=2\, v(\lambda)
$$ 
is even for every $\lambda\in (K_v)^\times$.  Note that when $K_v$ is
the completion of a function field over $\FF_q$ endowed with a
valuation $v$, with associated affine function ring $R_v$, the group
$\Ga_v= \PGL_2(R_v)$ is contained in $\PGL_2(K_v)^+$: For every $g\in
\GL_2(R_v)$, since $\det g\in(R_v)^\times=(\FF_q)^\times$,\footnote{See 
Equation \eqref{eq:inversRv}.} we have $v(\det g) =0$.

\medskip 
The following result proves the equidistribution in $K_v$ of the
loxodromic fixed points with complexity at most $s$ in a given orbit
by homographies under a lattice in $\PGL_2(K_v)$ as $s\ra+\infty$, and
its associated counting result. If $\xi\in \partial_\infty \XX_v =
\PP_1(K_v)$ and $\Ga$ is a subgroup of $\PGL_2(K_v)$, we denote by
$\Ga_\xi$ the stabiliser in $\Ga$ of $\xi$.

\btheo \label{theo:equidloxofix} Let $\Ga$ be a lattice in
$\PGL_2(K_v)^+$, and let $\ga_0\in\Ga$ be a loxodromic element of
$\Ga$. Then as $s\ra+\infty$,
$$
\frac{(q_v+1)^2\;\Vol(\Ga\dbs \XX_v)}
{2\;q_v^2\;\Vol(\Ga_{\ga_0^-}\dbs \Ax_{\ga_0})}\; s^{-1}
\sum_{\alpha\,\in\,\Ga\cdot \ga_0^-,\;h(\alpha)\leq s}\Delta_\alpha
\;\;\weakstar\;\; \haar_{K_v}
$$
and there exists $\kappa=\kappa_\Ga>0$ such that
$$
\card\{\alpha\in(\Ga\cdot \ga_0^-)\cap \OOO_v\;:\;h(\alpha)\leq s\}
= \frac{2\;q_v^2\;\Vol(\Ga_{\ga_0^-}\dbs \Ax_{\ga_0})}
{(q_v+1)^2\;\Vol(\Ga\dbs \XX_v)}\; s+\bigO(s^{1-\kappa})\;.
$$
\etheo

For every $\beta\in\;]0,\frac{1}{\ln q_v}]$, there is an error term of
the form $\bigO(s^{-\kappa}\|\psi\|_{\beta})$ for some $\kappa>0$ in
the equidistribution claim evaluated on any $\beta$-Hölder-continuous
function $\psi:\OOO_v\ra\CC$.

\medskip 
\dem The second result follows from the first one by integrating on the
characteristic function of the compact-open subset $\OOO_v$, whose
Haar measure is $1$.

\medskip In order to prove the equidistribution result, we apply
Theorem \ref{theo:algebrogeometricequid} with $\DD^-= \{*_v\}$ and
$\DD^+=\Ax_{\ga_0}$.  The families $\D^\pm= (\ga\DD^\pm)_{\ga\in\Ga/
  \Ga_{\DD^\pm}}$ are locally finite, since $\Ga$ is discrete and the
stabiliser $\Ga_{\DD^+}$ of $\DD^+ $ acts cocompactly on $\DD^+$.
Furthermore, $\|\sigma^-_{\D^+}\|$ is finite and nonzero by Equation
\eqref{eq:massskinline}.  Since $\Ga$ is contained in $\PGL_2(K_v)^+$,
the length spectrum $L_{\Ga}$ of $\Ga$ is contained in $2\ZZ$ by
\cite[II.1.2, Coro.]{Serre83}. Hence, it is equal to $2\ZZ$ by the
equivalence given by Equation \eqref{eq:equivLamGaevendisteven}.

For every $\ga\in \Ga$ such that $d(\DD^-,\ga\DD^+)>0$,\footnote{that
  is, such that $*_v\notin \ga\DD^+$} let $\alpha^-_{e,\,\ga}$ be the
generalised geodesic line, isometric exactly on $[0,d(\DD^-, \, \ga
  \DD^+)]$, whose image is the common perpendicular between $\DD^-$
and $\ga \DD^+$, and let $\rho_\ga$ be the geodesic ray starting at
time $0$ from the origin of $\alpha^-_{e,\,\ga}$ (which is $*_v$) with
point at infinity $\ga \cdot \ga_0^-$. Since $\XX_v$ is a tree and
$\ga \cdot \ga_0^-$ is one of the two endpoints of $\ga\DD^+$, the
geodesic segment $\alpha^-_{e,\,\ga} |_{[0,d(\DD^-, \ga\DD^+)]}$ is an
initial subsegment of $\rho_\ga$.\footnote{It connects $*_v$ to its
  closest point $P_{\ga\DD^+}(*_v)$ on $\ga\DD^+$, with
  $P_{\;\cdot\;}(\cdot)$ defined in Section \ref{subsect:nbhd}.}
Therefore, by Equation \eqref{eq:extensgeodraytoinfty}, for the
weak-star convergence of measures on $\normalout\DD^-$, we have
\begin{equation}\label{eq:demoequidloxofixThet}
\lim_{n\ra+\infty} \;\frac{({q_v}^2-1)(q_v+1)}{2\,q_v^3}\;
\frac{\Vol(\Ga\dbs\XX_v)}{\|\sigma^-_{\D^+}\|}\;{q_v}^{-n}
\sum_{\substack{\ga\in \Ga/\Ga_{\DD^+}\\
0< d(\DD^-, \, \ga \DD^+)\leq n}}
\Dirac_{\ga \cdot \ga_0^-} \;=\;
(\partial^+)_*\wt\sigma^+_{\DD^-}\;.
\end{equation}
Furthermore, for every $\beta\in\;]0,\frac{1}{\ln q_v}]$, by the
comment following Equation \eqref{eq:extensgeodraytoinfty} and since
$\Ga$, being a lattice in $\PGL_2(K_v)$, is geometrically
finite,\footnote{See the end of Section \ref{subsec:trees}.} we have
an error term of the form $\bigO(e^{-\kappa\,n}\,\|\psi\|_{\beta\ln
  q_v})$ for some $\kappa>0$ in the above formula when evaluated on
$\psi\in \C^{\beta\ln q_v}_c(\partial_\infty \XX_v)$, where
$\partial_\infty \XX_v$ is endowed with the visual distance
$d_{*_v}$.\footnote{Note that when $D=\{x\}$ is a singleton, the
  distance-like map $d_D$ used in Remark \ref{rem:utidD} coincides
  with the visual distance $d_x$, as said after Equation
  \eqref{eq:distancelike}.}  Note that on $\OOO_v$, the visual
distance $d_{*_v}$ and the distance $(x,y)\mapsto |x-y|_v$ are related
by
$$
|x-y|_v= d_{\H_\infty}(x,y)^{\ln q_v}=d_{*_v}(x,y)^{\ln q_v}\;,
$$ 
using Equation \eqref{eq:valuationham} for the first equality. Hence
we have an error term $\bigO(e^{-\kappa\,n}\,\|\psi\|_{\beta})$ for
some $\kappa>0$ in the above formula when evaluated on $\psi\in
\C^{\beta}_c(\OOO_v)$, where $\OOO_v$ is endowed with the distance
$(x,y)\mapsto |x-y|_v$.

\begin{center}
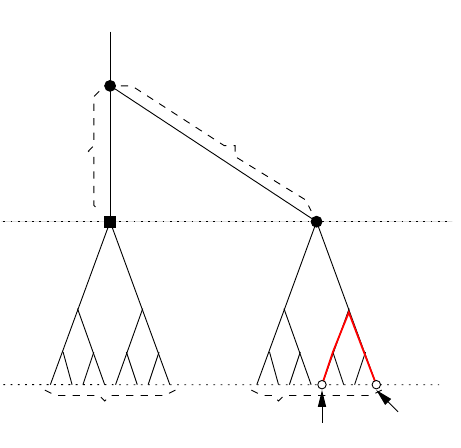
\end{center}

Let us fix for the moment $k\in \NN$. For every $\xi\in\pi_v^{-k}+
\OOO_v$, we have $|\xi|_v=q_v^{-v(\xi)}=q_v^k$ if $k\geq 1$ and
$|\xi|_v\leq 1$ if $k=0$. By restricting the measures to the
compact-open subset $\pi_v^{-k}+\OOO_v$ and by Proposition
\ref{prop:mescomputBT} (1), we have, with the appropriate error term
when $k=0$,
\begin{align}
\lim_{n\ra+\infty} \;\frac{({q_v}^2-1)(q_v+1)}{2\,q_v^3}\;
\frac{\Vol(\Ga\dbs\XX_v)}{\|\sigma^-_{\D^+}\|}\;&{q_v}^{-n}
\sum_{\substack{\ga\in \Ga/\Ga_{\DD^+}\\
\ga \cdot \ga_0^-\in \pi_v^{-k}+\OOO_v\\
0< d(\DD^-, \, \ga \DD^+)\leq n}}\Dirac_{\ga \cdot \ga_0^-} 
\nonumber\\ =& \;q_v^{-2\,k}\;{\haar_{K_v}}\big|_{(\pi_v)^{-k}+\OOO_v}\;.
\label{eq:applitwosubtreevenloxofix}
\end{align}

If $\beta\in\Ga$ is loxodromic and satisfies $\beta^-\in \pi_v^{-k} +
\OOO_v$ and $\beta^+\notin \pi_v^{-k}+\OOO_v$, then the translation
axis of $\beta$ passes at distance at most $2k$ from $*_v$, since it
passes through $P_{\H_\infty}(\pi_v^{-k})$ which is the closest point
on $\H_\infty$ to $\pi_v^{-k}$. If $\beta\in\Ga$ is loxodromic and
satisfies $\beta^-,\beta^+\in \pi_v^{-k} + \OOO_v$, then
$$
d(*_v, \Ax_\beta)=2k+d(\H_\infty,\Ax_\beta)\;.
$$
Furthermore, we  have, by Equations
\eqref{eq:valuationham} and \eqref{eq:defidisthamenbord}
$$
|\beta^- - \beta^+|_v=d_{\H_\infty}(\beta^-, \beta^+)^{\ln q_v}
=q_v^{-d(\H_\infty,\,\Ax_\beta)}\;.
$$
Therefore by the definition of the complexity in Equation
\eqref{eq:complexityloxofix}, we have for these elements 
\begin{equation}\label{eq:distbetapm}
h(\beta^-)= \frac{1}{|\beta^--\beta^+|_v}=q_v^{d(\H_\infty,\,\Ax_\beta)}
=q_v^{\,d(*_v, \,\Ax_\beta)-2k}\;.
\end{equation}

Since the family $\D^+= (\ga\DD^+)_{\ga\in\Ga/\Ga_{\DD^+}}$ is locally
finite, there are only finitely many elements $\ga\in\Ga/\Ga_{\DD^+}$
such that $\ga\DD^+=\Ax_{\ga\ga_0\ga^{-1}}$ is at distance at most
$2k$ from $*_v$. Hence for all but finitely many $\ga\in
\Ga/\Ga_{\DD^+}$ such that $\ga \cdot \ga_0^-=(\ga\ga_0\ga^{-1})^-\in
\pi_v^{-k}+\OOO_v$, we have $\ga\cdot \ga_0^+= (\ga\ga_0\ga^{-1})^+
\in \pi_v^{-k}+\OOO_v$ and, using Equation \eqref{eq:distbetapm} with
$\beta=\ga\ga_0\ga^{-1}$,
$$
h(\ga \cdot \ga_0^-)= q_v^{\,d(\DD^-, \ga\DD^+)-2k}\;.
$$
Therefore, using the change of variable $s=q_v^{\;n-2k}$, Equation
\eqref{eq:applitwosubtreevenloxofix} becomes
\begin{equation}\label{eq:applitwosubtreevenloxofix2}
\lim_{s\ra+\infty} \;\frac{({q_v}^2-1)(q_v+1)}{2\,q_v^3}\;
\frac{\Vol(\Ga\dbs\XX_v)}{\|\sigma^-_{\D^+}\|}\;s^{-1}
\sum_{\substack{\ga\in \Ga/\Ga_{\DD^+}\\
\ga \cdot \ga_0^-\in \pi_v^{-k}+\OOO_v\\
0< h(\ga \cdot \ga_0^-)\leq s}}
\Dirac_{\ga \cdot \ga_0^-} \;=\; {\haar_{K_v}}_{\;\mid (\pi_v)^{-k}+\OOO_v}\;.
\end{equation}

Note that the stabiliser $\Ga_{\ga_0^-}$ of $\ga_0^-$ in $\Ga$ has
index $\iota_{\ga_0^-}$ in $\Ga_{\DD^+}$ by the definition of
$\iota_{\ga_0^-}$ and that $\Ga/\Ga_{\ga_0^-}$ identifies with
$\Ga\cdot \ga_0^-$ by the map $\ga\Ga_{\ga_0^-}\mapsto \ga\cdot
\ga_0^-$. Since $\big((\pi_v)^{-k}+\OOO_v\big)_{k\in\NN}$ is a
countable family of pairwise disjoint compact-open subsets covering
$K_v$, and since the support of any continuous function with compact
support is contained in finitely many elements of this family, we have
\begin{equation}\label{eq:applitwosubtreevenloxofix3}
\lim_{s\ra+\infty} \;\frac{({q_v}^2-1)(q_v+1)}{2\,q_v^3\;\iota_{\ga_0^-}}\;
\frac{\Vol(\Ga\dbs\XX_v)}{\|\sigma^-_{\D^+}\|}\;s^{-1}
\sum_{\substack{\alpha\in \Ga\cdot \ga_0^-\\
0< h(\alpha)\leq s}} \Dirac_{\alpha} \;=\; \haar_{K_v}\;,
\end{equation}
with the appropriate error term.

Recall that by Equation \eqref{eq:massskinline}, if the Patterson
measures are normalised to be probability measures, then
$$
\|\sigma^-_{\D^+}\| = \frac{q_v-1}{q_v+1}\;\Vol(\Ga_{\DD^+}\dbs\DD^+)\;.
$$
Hence if instead the Patterson densities are normalised to have total
mass $\frac{q_v+1}{q_v}$ as in Proposition \ref{prop:mescomputBT} (2),
then
$$
\|\sigma^-_{\D^+}\| = \frac{q_v-1}{q_v}\;\Vol(\Ga_{\DD^+}\dbs\DD^+)\;.
$$
Note that, since $\iota_{\ga_0^-}=[\Ga_{\Ax_{\ga_0}}:\Ga_{\ga_0^-}]$,
$$
\Vol(\Ga_{\ga_0^-}\dbs \Ax_{\ga_0})=
\iota_{\ga_0^-}\;\Vol(\Ga_{\Ax_{\ga_0}}\dbs \Ax_{\ga_0})\;.
$$
Equation \eqref{eq:applitwosubtreevenloxofix3} thus gives the
equidistribution result in Theorem \ref{theo:equidloxofix}.  
\cqfd

\bigskip
In the following two Sections, we use Theorem \ref{theo:equidloxofix}
to deduce counting and equidistribution results of elements of
non-Archimedean local fields that are quadratic irrational over
appropriate subfields, when an appropriate algebraic complexity tends
to infinity.

\section{Counting and equidistribution of quadratic irrationals 
in positive characteristic}
\label{subsec:quadirratposchar}

Let $K$ be a (global) function field over $\FF_q$, let $v$ be a
(normalised discrete) valuation of $K$, let $K_v$ be the associated
completion of $K$ and let $R_v$ be the affine function ring associated
with $v$.\footnote{See Section \ref{subsec:valuedfields}.}

\medskip An element $\beta\in K_v$ is {\it quadratic
  irrational}\index{quadratic irrational} over $K$ if $\beta\notin K$
and $\beta$ is a root of a quadratic polynomial $a\beta^2+b\beta+c$
for some $a,b,c\in K$ with $a\ne 0$.  The {\em Galois
  conjugate}\index{Galois conjugate} $\gls{galconj}$ of $\beta$ is the
other root of the same polynomial. Let 
$$
\gls{trace}
=\beta+\beta^\sigma\;\;\;{\rm  and}\;\;\;\gls{norm} =\beta\beta^\sigma
$$ 
be the relative {\it trace}\index{trace!of quadratic irrational} and
relative {\it norm}\index{norm!of quadratic irrational} of $\beta$. It
is easy to check that $\beta^\sigma\neq\beta$, as the following lemma
shows.

\medskip 
Recall that for every field $F$, a polynomial $P\in F[Y]$ is {\it
  separable}\index{separable} if its roots in any algebraic closure of
$F$ are simple, and {\it inseparable}\index{inseparable} otherwise. It
is well known (see for instance \cite[\S V.6]{Lang02}) that, with $p$
the characteristic of $F$, any irreducible quadratic polynomial $P$
over $F$ is separable when $p\neq 2$, and is inseparable when $p=2$ if
and only if $P=a(Y^2-b)$ with $a\in F^\times$ and $b\in F$ which is
not a square in $F$.

\blemm 
An irreducible quadratic polynomial $P$ over $K$ which splits over
$K_v$ is separable.  
\elemm

\dem 
The result is immediate if $q$ is odd. Otherwise, assume for a
contradiction that $P$ is inseparable, so that $P=a(X^2-b)$ with $a\in
F^\times$ and $b\in K$ not a square in $K$. Since $P$ splits over
$K_v$, the element $b$ is a square in $K_v$. Since $K_v$ is isomorphic
to the field $\FF_{q_v}((\pi_v))$ of formal Laurent series over
$\FF_{q_v}$ with variable the uniformiser $\pi_v$, which may be
assumed to belong to $K$, there exist $m\in\ZZ$ and a sequence
$(a_n)_{n\in\NN}$ in $\FF_{q_v}$ such that, by the properties of the
Frobenius automorphism $x\mapsto x^2$,
$$
b=\big(\sum_{n\in\NN} a_n\pi_v^{m+k}\big)^2= 
\pi_v^{2m}\big(\sum_{n\in\NN} a_n^2\pi_v^{2k}\big)\;.
$$ 
Since $b\in K$, this implies that $a_n^2=0$ for $n$ large enough,
hence that $a_n=0$ for $n$ large enough, so that $b$ is a square in
$K$, a contradiction.  \cqfd

\medskip
The next proposition gives a characterisation
of quadratic irrationals over $K$.

\bprop \label{prop:caracquadrairrat}
Let $\beta\in K_v$. The following assertions are equivalent:
\begin{enumerate}
\item $\beta$ is quadratic irrational over $K$,
\item $\beta$ is a fixed point of a loxodromic element of
  $\PGL_2(R_v)$.
\end{enumerate}
\eprop

\dem The fact that (2) implies (1) is immediate since $\PGL_2(R_v)$
acts by homographies. The converse is classical once we know that
$\beta\neq \beta^\sigma$, see for instance \cite[Lem.~6.2]{ParPau11MZ}
in the Archimedean case and \cite{BerNak00} above its Section 5 when
$K=\FF_q(Y)$ and $v=v_\infty$.  \cqfd

\medskip
If $\beta\in K_v$ is quadratic irrational over $K$, its Galois
conjugate $\beta^\sigma$ is the other fixed point of any loxodromic
element of $\PGL_2(R_v)$ fixing $\beta$, hence the notations
$\beta^\sigma$ in this Section and in Section \ref{subsec:loxofix}
coincide.

The actions by homographies of the groups $\GL_2(R_v)$ and
$\PGL_2(R_v)$ on $K_v\cup\{\infty\}$ preserve the set of quadratic
irrationals over $K$. Contrary to the case of rational points, both
groups act with infinitely many orbits.

The {\em complexity}\index{quadratic irrational!complexity}%
\index{complexity} of a quadratic irrational $\alpha\in K_v$ over $K$ is
$$
\gls{complexityquadirrat}=\frac{1}{|\alpha-\alpha^\sigma|_v}\;,
$$ 
see for instance \cite[\S 6]{HerPau10} for motivations and results
when $K=\FF_q(Y)$ and $v=v_\infty$. Note that this complexity is
invariant under the action of the stabiliser $\GL_2(R_v)_\infty$ of
$\infty$ in $\GL_2(R_v)$, which is its upper triangular subgroup.  In
particular, it is invariant under the action of $R_v$ by
translations.\footnote{This is a particular case of Proposition
  \ref{prop:propricomplexity} (2) below.}  In \cite{ParPau12JMD},
where $K$ and $|\cdot|_v$ are replaced by $\QQ$ and its Archimedean
absolute value, there was, for convenience, an extra factor $2$ in the
numerator of the complexity, which is not needed here. We refer for
instance to \cite[Rem.~3, p.~136]{ParPau11MZ} for the connection of
this complexity to the standard height, and to \cite[\S 4.2,
  4.4]{ParPau12JMD} and \cite[\S 6.1]{ParPau11MZ} for studies using
this complexity.

The complexity $h(\cdot)$ satisfies the following elementary
properties, giving in particular its behaviour under the action of
$\PGL_2(R_v)$ by homographies on the quadratic irrationals in $K_v$
over $K$. We also give the well-known computation of the Jacobian of
the Haar measure for the change of variables given by homographies,
and prove the invariance of a measure which will be useful in Section
\ref{subsec:crossratioloxofix}.

For all $g=\begin{pmatrix} a & b\\ c & d \end{pmatrix}\in\GL_2(K_v) $
and $z\in K_v$ such that $g\cdot z\neq \infty$, let
$$
jg(z)=\frac{|\det g|_v}{|c\,z+d|_v^{\;2}}\,.
$$

\bprop\label{prop:propricomplexity} Let $\alpha\in K_v$ be a quadratic
irrational over $K$.
\begin{enumerate}
\item We have $\displaystyle h(\alpha)=
\frac{1}{\sqrt{|\tr(\alpha)^2-4\n(\alpha)|_v}}$~.
\item For every $g=\begin{pmatrix} a & b\\ c & d \end{pmatrix}\in
  \GL_2(K)$ with $|\det g|_v=1$, we have
$$
h(g\cdot\alpha)=|\n(d+c\,\alpha)|_v\;h(\alpha)\;.
$$
\item If $Q_\alpha:R_v\times R_v\ra [0,+\infty[$ is the map
  $(x,y)\mapsto |\n(x-y\,\alpha)|_v$, then for every $g\in
  \GL_2(R_v)$, we have
$$
Q_{g\cdot\alpha}=\frac{h(\alpha)}{h(g\cdot\alpha)}\;Q_\alpha\circ g^{-1}\;.
$$
In particular, if $g\in \GL_2(R_v)$ fixes $\alpha$, then
$$
Q_\alpha\circ g=Q_{\alpha}\;.
$$
\item For all $x,y,z\in K_v$ and $g\in \GL_2(K_v)$ such that $g\cdot
  x,\,g\cdot y,\,g\cdot z\neq \infty$, we have
$$
|g\cdot x-g\cdot y|_v^{\;2}\;= \;|x-y|_v^{\;2} \;\;jg(x)\;\;jg(y)
$$
and
$$
jg(z)=
\frac{d(g^{-1})_*\haar_{K_v}}{d\haar_{K_v}}(z)\;.
$$
\item The measure 
$$
d\mu(z)=\frac{d\haar_{K_v}(z)}{|z-\alpha|_v\;|z-\alpha^\sigma|_v}
$$
on $K_v-\{\alpha,\alpha^\sigma\}$ is invariant under the stabiliser of 
$\alpha$ in $\PGL_2(R_v)$ 
\end{enumerate}
\eprop

\dem (1) This follows from the formula $(\alpha-\alpha^\sigma)^2=
(\alpha+\alpha^\sigma)^2-4\,\alpha\,\alpha^\sigma$.  

\medskip\noindent  
(2) Since $g$ has rational coefficients (that is,
coefficients in $K$), we have 
\begin{align*}
g\cdot\alpha -(g\cdot \alpha)^\sigma&=g\cdot\alpha -g\cdot \alpha^\sigma
=\frac{a\alpha+b}{c\alpha+d}-\frac{a\alpha^\sigma+b}{c\alpha^\sigma+d}\\ & 
= \frac{(ad-bc)(\alpha-\alpha^\sigma)}{(c\alpha+d)(c\alpha^\sigma+d)}
= \frac{(\det g)(\alpha-\alpha^\sigma)}{\n(d+c\alpha)}\;.
\end{align*}
Taking absolute values and inverses, this gives Assertion (2).

\medskip\noindent  
(3) Let $g=\begin{pmatrix} a & b\\ c & d \end{pmatrix}\in
\GL_2(R_v)$. Note that $g^{-1}\cdot\alpha=
\frac{d\alpha-b}{a-c\alpha}$. For all $x,y\in R_v$, we hence have
\begin{align*}
\n\big((ax+by)-(cx+dy)\alpha\big)& =
\n\big(x(a-c\alpha)-y(d\alpha-b)\big)\\ &=
\n\big(x(a-c\alpha)-y(a-c\alpha)\,g^{-1}\cdot\alpha\big)\\ &=
\n(a-c\alpha)\;\n(x-y\;g^{-1}\cdot\alpha)\;.
\end{align*}
Taking absolute values and using Assertion (2), we have
$$
Q_\alpha\circ g=
\frac{h(g^{-1}\cdot\alpha)}{h(\alpha)}\;Q_{g^{-1}\cdot\alpha}\;.
$$
Assertion (3) follows by replacing $g$ by its inverse.

\medskip\noindent  
(4) Let $g=\begin{pmatrix} a & b\\ c & d \end{pmatrix}\in
\GL_2(K_v)$. As seen in the proof of Assertion (2), we have
$$
g\cdot x -g\cdot y= \frac{(\det g)(x-y)}{(c\,x+d)(c\,y+d)}\;.
$$
Taking absolute values and squares, this gives the first claim of
Assertion (4).

Recall that a homography $z\mapsto \frac{a\,z+b}{c\,z+d}$ is
holomorphic\footnote{We refer for instance to \cite{Serre92} for
  background on holomorphic functions on non-Archimedean local
  fields.} on $K_v-\{-\frac{d}{c}\}$, with derivative $z\mapsto
\frac{ad-bc}{(cz+d)^2}$. Hence infinitesimally close to $z$, the
homography acts (up to translations which leave the Haar measure
invariant) by a homothety of ratio $\frac{ad-bc}{(cz+d)^2}$. By
Equation \eqref{eq:homothetyhaar}, this proves that
$$
d\haar_{K_v}(g\cdot z)= 
\frac{|\det g|_v}{|c\,z+d|_v^{\;2}}\;d\haar_{K_v}(z)\;,
$$
as wanted.

\medskip\noindent 
(5) Let $g=\begin{pmatrix} a & b\\ c & d \end{pmatrix}\in \GL_2(R_v)$
fixing $\alpha$. Note that an element of $\GL_2(R_v)$ which fixes
$\alpha$ also fixes $\alpha^\sigma$.  By Assertion (4), we have
\begin{align*}
d\mu(g\cdot z)& = \frac{d\haar_{K_v}(g\cdot z)}
{|g\cdot z-\alpha|_v\;|g\cdot z-\alpha^\sigma|_v}
= \frac{d\haar_{K_v}(g\cdot z)}
{|g\cdot z-g\cdot \alpha|_v\;|g\cdot z-g\cdot \alpha^\sigma|_v}\\ &
=\frac{jg(z)\;d\haar_{K_v}(z)}
{|z-\alpha|_v\;\sqrt{jg(z)\;jg(\alpha)}\;
|z-\alpha^\sigma|_v\;\sqrt{jg(z)\;jg(\alpha^\sigma)}}
\\ &= \frac{1}{\sqrt{jg(\alpha)\;jg(\alpha^\sigma)}}\;d\mu(z)\;.
\end{align*}
By Assertion (4) applied with $x=\alpha$ and $y=\alpha^\sigma$,
we have
$$
\sqrt{jg(\alpha)\;jg(\alpha^\sigma)}=
\frac{|g\cdot\alpha - g\cdot \alpha^\sigma|_v}{|\alpha - \alpha^\sigma|_v}
= 1\;.
$$
The result follows. 
\cqfd

\medskip 
Let $G$ be a finite index subgroup of $\GL_2(R_v)$. We say that a
quadratic irrational $\beta\in K_v$ over $K$ is {\it
  $G$-reciprocal}\index{quadratic irrational!reciprocal}%
\index{reciprocal} (simply {\it reciprocal}\index{reciprocal} if
$G=\GL_2(R_v)$) if some element of $G$ maps $\beta$ to $\beta^\sigma$.
We define the {\it $G$-reciprocity index}\index{reciprocity index}
$\gls{reciprocityindex}$ as $2$ if $\beta$ is $G$-reciprocal and $1$
otherwise.  Similarly, we say that a loxodromic element $\ga$ of $G$
is {\it $G$-reciprocal}\index{loxodromic!reciprocal}\index{reciprocal}
(simply {\it reciprocal}\index{reciprocal} if $G=\GL_2(R_v)$) if there
exists an element in $G$ that switches the two fixed points of $\ga$.

\bprop \label{prop:caracreciprocal}
Let $G$ be a finite index subgroup of $\GL_2(R_v)$, and let $\ga$ be a 
loxodromic element of $G$. The following assertions are equivalent:
\begin{enumerate}
\item $\ga$ is conjugate in $G$ to $\ga'\ga^{-1}$ for some $\ga'\in
  G$ pointwise fixing $\Ax_{\ga}$,
\item the loxodromic element $\ga$ is $G$-reciprocal,
\item the quadratic irrational $\ga^-$ is $G$-reciprocal.
\end{enumerate}
When $G=\GL_2(R_v)$, Assertions (1), (2) and (3)  are also equivalent to 
\begin{enumerate}
\item[(4)] the image of $\ga''\ga$ in $\PGL_2(R_v)$, for some $\ga''\in
  G$ pointwise fixing $\Ax_{\ga}$, is conjugate to the
  image  in $\PGL_2(R_v)$ of $\;^t\ga$.
\end{enumerate}
\eprop

\dem Most of the proofs are similar to the ones when $R_v$, $K$ and
$|\cdot|_v$ are replaced by $\ZZ$, $\QQ$ and its Archimedean absolute
value, see for instance \cite{ParPau12JMD}.  We only give hints for
the sake of completeness.  Let $\alpha=\ga^-$.

If $\alpha$ is $G$-reciprocal, then let $\beta\in G$ be such that
$\beta\cdot\alpha=\alpha^\sigma$. Since $R_v\subset K$, we have
$\beta\cdot\alpha^\sigma=\alpha$. Hence $\beta\ga\beta^{-1}$ is a
loxodromic element of $G$ fixing $\alpha$ and $\alpha^\sigma$, having the 
same translation length as $\ga$, but translating in the opposite
direction on $\Ax_{\ga}$. Hence $\ga'=\beta\ga\beta^{-1}\ga$
fixes pointwise $\Ax_{\ga}$. Therefore (3) implies (1).

If $\beta\in G$ conjugates $\ga$ to $\ga'\ga^{-1}$ for some $\ga'\in
G$ pointwise fixing $\Ax_{\ga}$, then $\beta$ preserves the set
$\{\alpha,\alpha^\sigma\}$. Hence, it preserves the translation axis
of $\ga$ but it switches $\alpha$ and $\alpha^\sigma$ since $\ga$ and
$\ga'\ga^{-1}$ translate in opposite directions on
$\Ax_{\ga}$. Therefore (1) implies (2).

The fact that (2) implies (3) is immediate, since $\alpha^\sigma =\ga^+$.

\medskip The equivalence between (1) and (4) when $G=\GL_2(R_v)$
follows from the fact that the stabiliser of $\Ax_{\ga}$ normalises
the pointwise stabiliser of $\Ax_{\ga}$, and from the formula $$
^t\begin{pmatrix}a & b \\ c & d\end{pmatrix}^{-1}=\frac{1}{ad-bc}\;
\begin{pmatrix}0 & 1 \\ -1 & 0\end{pmatrix}
\begin{pmatrix}a & b \\ c & d\end{pmatrix}
\begin{pmatrix}0 & 1 \\ -1 & 0\end{pmatrix}^{-1}\;
$$
which is valid over any field.
\cqfd

\bigskip The following result says that any orbit of a given quadratic
irrational in $K_v$ over $K$, by homographies under a given finite
index subgroup of the modular group $\PGL_2(R_v)$, equidistributes to
the Haar measure on $K_v$. Again,  note that we are not assuming
the finite index subgroup to be a congruence subgroup.   
  
\btheo \label{theo:caracposequidquadirr}
Let $G$ be a finite index subgroup of $\GL_2(R_v)$.  Let
$\alpha_0\in K_v$ be a quadratic irrational over $K$. Then, as
$s\to+\infty$,
$$
\frac{(q_v+1)^2\;\zeta_K(-1)\;m_0\;[\GL_2(R_v):G]}
{2\;q_v^2\;(q-1)\;|v(\operatorname{tr} g_0)|}\;s^{-1}
\sum_{\alpha\in G\cdot\alpha_0\;:\;h(\alpha)\le s}\Delta_{\alpha}\;\;
\weakstar\;\;\haar_{K_v}\;,
$$
where $g_0\in G$ fixes $\alpha_0$ with $v(\operatorname{tr} g_0)\neq
0$, and where $m_0$ is the index of $g_0^\ZZ$ in the stabiliser of
$\alpha_0$ in $G$.  Furthermore, there exists $\kappa>0$ such that, as
$s\to+\infty$,
$$
\card\{\alpha\in (G\cdot\alpha_0)\cap\OOO_v\;:\;h(\alpha)\le s\}
= \frac{2\;q_v^2\;(q-1)\;|v(\operatorname{tr} g_0)|}{(q_v+1)^2\;
\zeta_K(-1)\;m_0\;[\GL_2(R_v):G]}\;s +\bigO(s^{1-\kappa})\;.
$$
\etheo

For every $\beta\in\;]0,\frac{1}{\ln q_v}]$, there exists $\kappa>0$
such that for every $\psi\in\C^\beta_c(K_v)$ there is an
error term in the above equidistribution claim evaluated on $\psi$, of
the form $\bigO(s^{-\kappa}\|\psi\|_{\beta})$.

\medskip 
\dem We apply Theorem \ref{theo:equidloxofix} with $\Ga$ the
image of $G$ in $\Ga_v=\PGL_2(R_v)$ and with $\ga_0$ the image in
$\Ga_v$ of the element $g_0$ introduced in the statement.  Note
that $\Ga$, which is contained in $\Ga_v$, is indeed contained in
$\PGL_2(K_v)^+$.

By Equation \eqref{eq:translatlength}, for every $g\in\GL_2(R_v)$, the
translation length of $g$ in $\XX_v$ is $2\;|v(\operatorname{tr} g)|$,
and $g$ is loxodromic if and only if $v(\operatorname{tr} g)\neq
0$. This implies that $g_0$ exists, since $G$ has finite index in
$\GL_2(R_v)$, and such an element exists in $\GL_2(R_v)$ by
Proposition \ref{prop:caracquadrairrat}. Up to replacing $g_0$ by its
inverse, which changes neither $|v(\operatorname{tr} g_0)|$ nor $m_0$,
we assume that $\ga_0^-=\alpha_0$. Furthermore
$$
\lambda(\ga_0)= 2\;|v(\operatorname{tr} g_0)|\;.
$$

Since the centre of $\GL_2(K_v)$ acts trivially by homographies, we
have
$$
G\cdot \alpha_0= \Ga\cdot \alpha_0\;.
$$
For every $\alpha\in G\cdot \alpha_0$, the complexities $h(\alpha)$,
when $\alpha$ is considered as a quadratic irrational or when
$\alpha$ is considered as a loxodromic fixed point, coincide.

Since the centre $Z(G)$ of $G$ acts trivially by homographies, by the
definition of $m_0$ in the statement, we have
$$
[\Ga_{\ga_0^-}: \ga_0^\ZZ] = \frac{[G_{\alpha_0}:g_0^\ZZ]}{|Z(G)|}= 
\frac{m_0}{|Z(G)|}\;.
$$
Therefore, 
\begin{align}
\Vol(\Ga_{\ga_0^-}\dbs\Ax_{\ga_0})& = \frac{1}{[\Ga_{\ga_0^-}:\ga_0^\ZZ]}
\;\Vol(\ga_0^\ZZ\dbs\Ax_{\ga_0}) =
\frac{\lambda(\ga_0)}{[\Ga_{\ga_0^-}: \ga_0^\ZZ]}\nonumber\\ & 
=\frac{2\;|v(\operatorname{tr} g_0)|\;|Z(G)|}{m_0}\;.
\label{eq:covolcycle}
\end{align}

Theorem \ref{theo:caracposequidquadirr} now follows from Theorem
\ref{theo:equidloxofix} using Equations \eqref{eq:covoltotMertens} and
\eqref{eq:covolcycle}.  
\cqfd

\bexem 
(1) Theorem \ref{theo:appliquadirratintro} in the Introduction
follows from this result, by taking $K=\FF_q(Y)$ and $v=v_\infty$, and
by using Equation \eqref{zetamoinsun} in order to simplify the
constant.

\medskip
\noindent (2) Let $G_I$ be the Hecke congruence subgroup associated with a
nonzero ideal $I$ of $R_v$, see Equation
\eqref{eq:heckecongruencesubgroup}. By Lemma
\ref{lem:indexHeckesubgrou}, we have, as $s\to+\infty$,
$$
\frac{(q_v+1)^2\;\zeta_K(-1)\;m_0\;N(I)\prod_{\ppp|I}(1+\frac{1}{N(\ppp)})}
{2\,q_v^2\;(q-1)\;|v(\operatorname{tr} g_0)|}\;s^{-1}
\sum_{\alpha\in G_I\cdot\alpha_0\;:\;h(\alpha)\le s} \Delta_{\alpha}\;\;
\weakstar\;\;\haar_{K_v}\;.
$$
\eexem

\medskip
We conclude this Section by a characterisation of quadratic
irrationals and reciprocal quadratic irrationals in the field of
formal Laurent series $\FF_q((Y^{-1}))$ in terms of continued
fractions. When $\FF_q[Y]$, $\FF_q(Y)$ and $v_\infty$ are replaced by
$\ZZ$, $\QQ$ and its Archimedean absolute value, we refer for instance
to \cite{Sarnak07} and \cite[Prop.~4.3]{ParPau12JMD} for
characterisations of reciprocal quadratic irrationals.

\medskip
Recall that Artin's {\it continued fraction 
expansion}\index{continued fraction} of $f\in \FF_q((Y^{-1}))-
\FF_q(Y)$ is the sequence $(a_i=a_i(f))_{i\in\NN}$ in $\FF_q[Y]$ with
$\deg a_i>0$ if $i>0$ such that
$$
f=a_0+\frac{1}{\displaystyle a_1+\frac{1}{\displaystyle 
a_2+\frac{1}{\displaystyle a_3+\frac{1}{\ddots}}}}\;.
$$ 
See for instance the surveys \cite{Lasjaunias00,Schmidt00}, as well as
\cite{Paulin02} for a geometric interpretation.  We say that the
continued fraction expansion of $f$ is {\it eventually
  periodic}\index{continued fraction!eventually periodic} if there
exist $n\in\NN$ and $N\in\NN-\{0\}$ such that $a_{n+i}=a_{n+N+i}$ for
every $i\in\NN$, and we write
$$
f=[a_0,\dots,a_{n-1},\overline{a_n,\dots,a_{n+N-1}}]\;.
$$
Such a sequence $a_n,\dots,a_{n+N-1}$ is called a {\it
  period}\index{period} of $f$, and if of minimal length, it is well
defined up to cyclic permutation.

Two elements $\beta, \beta' \in \FF_q((Y^{-1}))$ are in the same
$\PGL_2(\FF_q[Y])$-orbit if and only if their continued fraction
expansions have equal tails up to an invertible element of $\FF_q[Y]$
by \cite[Theo.~1]{Schmidt00} or \cite[Theo.~1]{BerNak00} (and even
before that by \cite[Sect.~IV.3]{deMathan70}).  More precisely, $\beta,
\beta' \in K_v$ are in the same $\PGL_2(\FF_q[Y])$-orbit if and only
if there exist $m,n\in\NN$ and $x\in\FF_q^\times$ such that for every
$k\in\NN$, we have $a_{n+k}(\beta')= x^{(-1)^k}\,a_{m+k}(\beta)$.

\bprop Assume that $K=\FF_q(Y)$ and $v=v_\infty$.
\begin{enumerate}
\item An element $\alpha\in K_v-K$ is quadratic irrational over
  $K$ if and only if its continued fraction expansion of $\beta$ is
  eventually periodic, and  if and only if it is a fixed point of a
  loxodromic element of $\PGL_2(\FF_q[Y])$.
\item A quadratic irrational $\alpha\in K_v$ is reciprocal if and
  only if the period $a_0,\dots, a_{N-1}$ of the continued fraction
  expansion of $\alpha$ is palindromic up to cyclic permutation and
  invertible elements, in the sense that there exist $x\in
  \FF_q^\times$ and $p\in\NN$ such that for $k=0,\dots,N-1$, we have
  $a_{k+p}=x^{(-1)^k}a_{N-k-1}$ (with indices modulo $N$).
\end{enumerate}
\eprop
 
\dem (1) The equivalence of being quadratic irrational and having an
eventually periodic continued fraction expansion is well-known, see
for instance the survey \cite[Theo.~3.1]{Lasjaunias00}. The second
part of the claim follows from Proposition
\ref{prop:caracquadrairrat}.

\medskip \noindent
(2) The proof is similar to the Archimedean case in \cite[\S
23]{Perron13}.\footnote{See also \cite[Coro.~1]{BerNak00} by relating,
  using twice the period, what the authors call the $-$ continued
  fraction expansion to the standard expansion.} For every quadratic
irrational $f\in \FF_q((Y^{-1}))$, up to the action of
$\GL_2(\FF_q[Y])$, we may assume that $f, (f^\sigma)^{-1} \in
Y^{-1}\FF_q[[Y^{-1}]]$ and $f= [0,\; \overline{a_1,a_2,\dots,
  a_n}]$. Then we may define by induction quadratic irrationals
$f_2,\dots,f_n\in \FF_q((Y^{-1}))$ over $\FF_q(Y)$ such that
$$
\frac{1}{f}= a_1+ f_2,\;\;\frac{1}{f_2}= a_2+ f_3, \;\dots, \;
\frac{1}{f_{n+1}}= a_{n+1}+ f_n, \;\;
\frac{1}{f_n}= a_n+ f\;.
$$
Passing to the Galois conjugates, we have
$$
\frac{1}{f^\sigma}= a_1+ f_2^\sigma, \;
\;\frac{1}{f_2^\sigma}= a_2+ f_3^\sigma, \;\dots, \;
\frac{1}{f_n^\sigma}= a_n+ f^\sigma\;.
$$
Taking these equations in the reverse order, we have
$$
\frac{1}{-\frac{1}{f^\sigma}}= a_n - \frac{1}{f_n^\sigma},\;\;
\frac{1}{-\frac{1}{f_n^\sigma}}= a_{n-1} - \frac{1}{f_{n-1}^\sigma},
\;\dots, \;
\frac{1}{-\frac{1}{f_2^\sigma}}= a_1- \frac{1}{f^\sigma}\;,
$$
so that, since $-\frac{1}{f^\sigma}\in Y^{-1}\FF_q[[Y^{-1}]]$, we have
$$
-\frac{1}{f^\sigma}=[0,\;\overline{a_n,\dots, a_2, a_1}]\;.
$$ 
Therefore $f^\sigma=[\overline{-a_n,\dots, -a_2, -a_1}]$.  Thus, if
$f$ and $f^\sigma$ are in the same orbit, the periods are palindromic
up to cyclic permutation and invertible elements by
\cite[Theo.~1]{Schmidt00}, \cite[Theo.~1]{BerNak00}.  
\cqfd

%
%

\section{Counting and equidistribution of quadratic 
irrationals in $\QQ_p$}
\label{subsec:quadirratzerochar}

There are interesting arithmetic (uniform) lattices of $\PGL_2(\QQ_p)$
constructed using quaternion algebras.  In this Section, we study
equidistribution properties of loxodromic fixed points elements of
these lattices. See for instance \cite{LedPol05} for an
equidistribution result of the eigenvalues of the loxodromic elements.  We
use \cite{Vigneras80} as our standard reference on quaternion
algebras.

\medskip
Let $F$ be a field and let $a,b\in F^\times$. Let $D= \big(\frac{a,b}F\big)$
be the quaternion algebra over $F$ with basis $1,i,j,k$ as a
$F$-vector space such that $i^2=a$, $j^2=b$ and $ij=ji=-k$.  If
$x=x_0+x_1i+x_2j+x_3k\in D$, then its {\em conjugate}\index{conjugate}
is
$$
\gls{conjugate}= x_0-x_1i-x_2j-x_3k\,,
$$
its (reduced) {\em norm}\index{norm!of a quaternion} is 
$$
\gls{normquaternion}= x\,\overline{x}\,=\,\overline{x}\, x= x_0^2-a\,x_1^2
-b\,x_2^2+ab\,x_3^2
$$ 
and its (reduced) {\em trace}\index{trace!of a quaternion} is 
$$
\gls{tracequaternion}= x+\overline{x} =2\,x_0\,.
$$
   
Let us fix two negative rational integers $a,b$ and let
$D=\big(\frac{a,b}\QQ\big)$ .  For every field extension $E$ of $\QQ$, we
denote by $D_E$ the quaternion algebra $D\otimes_\QQ E$ over $E$, and
we say that $D$ {\it splits over}\index{splitting over} $E$ if the
$E$-algebra $D\otimes_\QQ E$ is isomorphic to $M_2(E)$.  The
assumption that $a,b$ are negative implies that $D$ does not split
over $\RR$.  Furthermore, when $p\in\NN $ is an odd prime, $D$ splits
over $\QQ_p$ if and only if the equation $a\,x^2+b\,y^2=1$ has a
solution in $\QQ_p$, see \cite[page 32]{Vigneras80}.

The {\em reduced discriminant}\index{discriminant!reduced} of $D$ is 
$$
\gls{discred}=\prod_{q\in{\rm Ram}(D)}q\,.
$$ 
where ${\rm Ram}(D)$ is the finite set of primes $p$ such that $D$
does not split over $\QQ_p$.

For instance, the quaternion algebra $D=\big(\frac{-1,-1}\QQ\big)$ splits over
$\QQ_p$ if and only if $p\neq 2$, hence it has reduced discriminant $2$.

Assume from now on that $p\in\NN$ is a positive rational prime such
that $D$ splits over $\QQ_p$ and, for simplicity, that $\QQ_p$
contains square roots $\sqrt{a}$ and $\sqrt{b}$ of $a$ and $b$. For
example, if $a=b=-1$, this is satisfied if $p\equiv 1\mod 4$.  We
then have an isomorphism of $\QQ_p$-algebras $\theta
=\theta_{a,\,b}:D_{\QQ_p}\to M_2(\QQ_p)$ defined by
\begin{equation}\label{eq:defisigma}
\theta(x_0+x_1i+x_2j+x_3k)=
\begin{pmatrix}
x_0+x_1\sqrt{a} & \sqrt{b}\,(x_2+\sqrt{a}\,x_3)\vspace{.2cm}\\
\sqrt{b}\,(x_2-\sqrt{a}\,x_3) & x_0-x_1\sqrt{a}
\end{pmatrix}\,, 
\end{equation}
so that 
$$
\det(\theta(x))=\redn(x)\;\;\;{\rm and}\;\;\;
\operatorname{tr} (\theta(x)) =\redtr(x)\;.
$$ 
If the assumption on the existence of the square roots in $\QQ_p$
is not satisfied, we can replace $\QQ_p$ by an appropriate finite
extension, and prove equidistribution results in this extension.

Let $\OOO$ be a {\it $\ZZ\big[\frac{1}{p}\big]$-order}\index{order} in
$D_{\QQ_p}$, that is, a finitely generated
$\ZZ\big[\frac{1}{p}\big]$-submodule of $D_{\QQ_p}$ generating
$D_{\QQ_p}$ as a $\QQ_p$-vector space, which is a subring of
$D_{\QQ_p}$. Let $\OOO^1$ be the group of elements of norm $1$ in
$\OOO$. Then the image $\Ga^1_\OOO$ of $\theta(\OOO^1)$ in
$\PGL_2(\QQ_p)$ is a cocompact lattice, see for instance
\cite[Sect.~IV.1]{Vigneras80}. In fact, this lattice is contained in
$\PSL_2(\QQ_p)$, hence in $\PGL_2(\QQ_p)^+$.  In this Section
\ref{subsec:quadirratzerochar}, we denote by $\XX_{p}$ the Bruhat-Tits
tree of $(\PSL_2,\QQ_p)$, which is $(p+1)$-regular.

The next result computes the covolume of this lattice.\footnote{The
  index $q$ ranges over the primes dividing $\discr_D$, that is, over
  the elements of ${\rm Ram}(D)\,$.}

\bprop\label{pro:padictvol} Let $D$ be a quaternion algebra over $\QQ$
which splits over $\QQ_p$ and does not split over $\RR$, and let
$\OOO$ be a $\ZZ\big[\frac{1}{p}\big]$-order in $D_{\QQ_p}$. If
$\OOO_{\rm max}$ is a maximal $\ZZ\big[\frac{1}{p}\big]$-order in
$D_{\QQ_p}$ containing $\OOO$, then
$$
\Vol(\Ga^1_{\OOO}\dbs\XX_p)=[\OOO^1_{\rm max}:\OOO^1]\;
\frac{p}{12}\;\prod_{q\divides\discr_D}(q-1)\,.
$$
\eprop

\dem 
We refer to \cite[page 53]{Vigneras80} for the (common)
definition of the {\em discriminant}\index{discriminant}
$\discr(\QQ_p)$ of the local field $\QQ_p$ and $\discr(D_{\QQ_p})$ of
the quaternion algebra $D_{\QQ_p}$ over the local field $\QQ_p$. We
will only use the facts that $\discr(\QQ_p)=1$ as it easily follows
from the definition, and that
\begin{equation}\label{eq:calcdiscquatalg}
\discr(D_{\QQ_p})=\discr(\QQ_p)^4\big(N(p\ZZ_p)\big)^2=p^2
\end{equation} 
which follows by \cite[Lem.~4.7, page 53]{Vigneras80} and
\cite[Cor.~1.7, page 35]{Vigneras80} for the first equality and
$N(p\ZZ_p)=\card(\ZZ_p/p\ZZ_p)=\card(\ZZ/p\ZZ)=p$ for the second one.

We refer to \cite[Sect.~II.4]{Vigneras80} for the definition of the
{\em Tamagawa measure}\index{Tamagawa measure}\index{measure!Tamagawa}
$\mu_{\Tamagawa}$ on $X^\times$ when $X=D_{\QQ_p}$ or $X=\QQ_p$. It is
a Haar measure of the multiplicative locally compact group $X^\times$,
and understanding its explicit normalisation is the main point of this
proposition. By \cite[Lem.~4.6, page 52]{Vigneras80},\footnote{See
  more precisely the top of page 55 in op.~cit.} with $dx$ the Haar
measure on the additive group $X$,\footnote{with a normalisation that
  does not need to be made precise} with $\|x\|$ the {\it module} of
the left multiplication by $x\in X^\times$ on the additive group
$X$,\footnote{so that $(M_x)_*dx=\|x\|\,dx$ where $M_x:y\mapsto xy$ is
  the left multiplication by $x$ on $X$} we have
$$
d\mu_{\Tamagawa}(x)=\frac{1}{\sqrt{\discr(X)}\;\|x\|}\;dx\;.
$$ 
By \cite[proof of Lem.~4.3, page 50]{Vigneras80}, identifying
$D_{\QQ_p}$ to $M_2(\QQ_p)$ by $\theta$, the measure of $\GL_2(\ZZ_p)$
for the measure $\frac{1}{(1-p^{-1})\;\|x\|}\;dx$ is
$1-p^{-2}$. Hence, by scaling and by Equation
\eqref{eq:calcdiscquatalg}, we have
$$
\mu_{\Tamagawa}(\GL_2(\ZZ_p))=
\frac{(1-p^{-2})(1-p^{-1})}{\sqrt{\discr(D_{\QQ_p})}}
=\frac{(p^2-1)(p-1)}{p^4}\,.
$$ 
By \cite[Lem.~4.3, page 49]{Vigneras80}, the mass of $\ZZ_p^\times$
for the measure $\frac{1}{(1-p^{-1})\;\|x\|} \;dx$ on $\QQ_p^\times$
is $1$, hence by scaling
$$
\mu_{\Tamagawa}(\ZZ_p^\times)=\frac{1-p^{-1}}{\sqrt{\discr(\QQ_p)}}
=\frac{p-1}{p}\,.
$$
By \cite[pages 53--54]{Vigneras80}, since we have an exact sequence
$$
1\longrightarrow \SL_2(\QQ_p)\longrightarrow \GL_2(\QQ_p)
\stackrel{\det}{\longrightarrow}\;\QQ_p^\times\longrightarrow 1\;,
$$ 
the Tamagawa measure of $\GL_2(\QQ_p)$ disintegrates by the
determinant over the Tamagawa measure of $\QQ_p^\times$ with
conditional measures the translates of a measure on $\SL_2(\QQ_p)$,
called the {\it Tamagawa measure}\index{Tamagawa measure}%
\index{measure!Tamagawa} of $\SL_2(\QQ_p)$ and again denoted
by $\mu_{\Tamagawa}$.  Thus,
$$
\mu_{\Tamagawa}(\SL_2(\ZZ_p))=
\frac {\mu_{\Tamagawa}(\GL_2(\ZZ_p))}{\mu_{\Tamagawa}(\ZZ_p^\times)}=\frac{p^2-1}
{p^3}
$$

By Example 3 on page 108 of \cite{Vigneras80}, since the
$\ZZ\big[\frac{1}{p}\big]$-order $\OOO_{\rm max}$ is maximal, we have,
with $G=\theta(\OOO^1_{\rm max})$,
$$
\mu_{\Tamagawa} (G\bs\SL_2(\QQ_p)) =
\frac{1}{24}\;(1-p^{-2})\;\prod_{q\divides\discr_D}(q-1)\;.
$$

Since $\GL_2(\QQ_p)$ acts transitively on $V\XX_p$ with stabiliser of
the base point $*=[\ZZ_p\times\ZZ_p]$ the maximal compact subgroup
$\GL_2(\ZZ_p)$,\footnote{See Section \ref{subsec:BruhatTitstrees}.}
and by the centred equation mid-page 116 of \cite{Serre83}, we have
\begin{align*}
\Vol(G\dbs\XX_p) &=\sum_{[x]\in G\bs V\XX_p}\frac{1}{|G_x|}=
\frac{\mu_{\Tamagawa}(G\bs\GL_2(\QQ_p))}{\mu_{\Tamagawa}(\GL_2(\ZZ_p))}
=\frac{\mu_{\Tamagawa} (G\bs\SL_2(\QQ_p) )}
{\mu_{\Tamagawa}(\SL_2(\ZZ_p))}\\ &=
\frac{p}{24}\;\prod_{q|\discr_D}(q-1)\,.
\end{align*} 
The natural homomorphism $G=\theta(\OOO^1_{\rm max})\ra
\Ga^1_{\OOO_{\rm max}}$ is $2$-to-$1$, so that
$$
\Vol(\Ga^1_{\OOO_{\rm max}}\dbs\XX_p)=2\,\Vol(G\dbs\XX_p)\,.
$$ 
Since $[\Ga^1_{\OOO_{\rm max}}: \Ga^1_{\OOO}] = [\OOO^1_{\rm max} :
  \OOO^1]$, Proposition \ref{pro:padictvol} follows.  
\cqfd

\bigskip 
Note that the fixed points $z$ for the action on $\PP^1(\QQ_p)=
\QQ_p\cup\{\infty\}$ by homographies of the elements in the image of
$\theta(D)$ are quadratic over $\QQ(\sqrt{a},\sqrt{b})$.  More
precisely, $\frac{z}{\sqrt{b}}$ is quadratic over $\QQ(\sqrt{a})$.  An
immediate application of Theorem \ref{theo:equidloxofix}, using
Proposition \ref{pro:padictvol}, gives the following result of
equidistribution of quadratic elements in $\QQ_p$ over
$\QQ(\sqrt{a},\sqrt{b})$.

\btheo\label{theo:padicequidistrib} Let $\Ga$ be a finite index
subgroup of $\Ga^1_\OOO$, and let $\ga_0\in\Ga$ be a loxodromic
element of $\Ga$. Then as $s\ra+\infty$,
\begin{align*}
\frac{(p+1)^2\;\prod_{q|\discr_D}(q-1)\;[\OOO^1_{\rm max}:\OOO^1]\;
[\Ga^1_\OOO:\Ga]}
{24\;p\;\Vol(\Ga_{\ga_0^-}\dbs \Ax_{\ga_0})}\; & s^{-1}
\sum_{\alpha\,\in\,\Ga\cdot \ga_0^-,\;h(\alpha)\leq s}\Delta_\alpha
\\ & \weakstar\;\; \haar_{\QQ_p}\;,
\end{align*} 
where $\OOO_{\rm max}$ is a maximal $\ZZ\big[\frac{1}{p}\big]$-order
in $D_{\QQ_p}$ containing $\OOO$, and there exists $\kappa>0$ such
that as $s\ra+\infty$
\begin{align*}
\card& \{\alpha\in(\Ga\cdot \ga_0^-)  \cap  \ZZ_p\;:\;h(\alpha)\leq s\}
\\& =\; \frac{24\;p\;\Vol(\Ga_{\ga_0^-}\dbs \Ax_{\ga_0})}
{(p+1)^2\;\prod_{q|\discr_D}(q-1)\;[\OOO^1_{\rm max}:\OOO^1]\;
[\Ga^1_\OOO:\Ga]}\; s + \bigO(s^{1-\kappa})\;. \; \; \; \;\Box
\end{align*} 
\etheo

Assume furthermore that the positive rational prime $p\in\NN$ is such
that $p\equiv 1 \mod 4$ and that the integer $\frac{p^2-1}{4}$ is not
of the form $4^a(8b+7)$ for $a,b\in \NN$ (for instance $p=5$).  By
Legendre's three squares theorem (see for instance \cite{Grosswald85}),
there exist $x'_1,x'_2, x'_3 \in\ZZ$ such that $\frac{p^2-1}{4}=
{x'_1}^2 + {x'_2}^2 + {x'_3}^2$.  Hence there are $x_1,x_2,x_3\in
2\ZZ$ such that $p^2-1= {x_1}^2 + {x_2}^2 + {x_3}^2$.

A standard consequence of Hensel's theorem says that when $p$ is odd,
a number $n\in\ZZ$ has a square root in $\ZZ_p$ if $n$ is relatively
prime to $p$ and has a square root modulo $p$, see for instance
\cite[page 351]{Knapp07}. Thus, $1-p^2$ has a square root in $\ZZ_p$,
that we denote by $\sqrt{1-p^2}$.  As noticed above, since $p\equiv 1
\mod 4$, the element $-1$ has a square root in $\QQ_p$, that we denote
by $\varepsilon$.  The element
$$
\alpha_0=
\frac{\varepsilon\, x_1+\sqrt{1-\,p^2}}{x_3+\varepsilon \,x_2}
$$ 
is a quadratic irrational in $\QQ_p$ over $\QQ(\varepsilon)$.

\medskip
The following result is a counting and equidistribution result of 
quadratic irrationals over $\QQ(\varepsilon)$ in $\QQ_p$. We
denote by $\gls{galconjqp}$ the Galois conjugate of a quadratic
irrational $\alpha$ in $\QQ_p$ over $\QQ(\varepsilon)$, and by
$$
\gls{complexityquadirratqp}=\frac{1}{|\alpha-\alpha^\sigma|_p}
$$
the {\it complexity}\index{complexity} of $\alpha$.

\btheo\label{theo:padicequidistribexemplx} Let $D=\big(\frac{-1,-1}\QQ\big)$
be Hamilton's quaternion algebra over $\QQ$.  Let $p\in\NN$ be a
positive rational prime with $p\equiv 1 \mod 4$ such that
  $p^2-1= {x_1}^2 + {x_2}^2 + {x_3}^2$ for some $x_1,x_2,x_3\in
2\ZZ$, and 
let $\OOO$ be the $\ZZ\big[\frac{1}{p}\big]$-order\footnote{This order
  plays an important role in the construction of Ramanujan graphs by
  Lubotzky, Phillips and Sarnak \cite{LubPhiSar86,LubPhiSar88} (see
  also \cite[\S 7.4]{Lubotzky94}), and in the explicit construction of
  free subgroups of $\SO(3)$ in order to construct
  Hausdorff-Banach-Tarsky paradoxical decompositions of the
  $2$-sphere, see for instance \cite[page 11]{Lubotzky94}.}
$$
\OOO=\big\{y\in\ZZ\big[{\textstyle \frac{1}{p}}\big]+
\ZZ\big[{\textstyle\frac{1}{p}}\big]\;i+
\ZZ\big[{\textstyle\frac{1}{p}}\big]\;j+
\ZZ\big[{\textstyle\frac{1}{p}}\big]\;k\;:\; 
y\equiv 1\mod 2\big\}\;
$$ 
in $D_{\QQ_p}$.  Let $\Ga$ be a finite index subgroup of
$\Ga^1_\OOO$. Then as $s\ra+\infty$,
$$
\frac{(p+1)^2\;[\Ga^1_\OOO:\Ga]}{2\; p^2\;k_\Ga}\; s^{-1}
\sum_{\alpha\,\in\,\Ga \cdot \alpha_0,\;h(\alpha)\leq s}\Delta_\alpha
\;\;\weakstar\;\; \haar_{\QQ_p}\;,
$$ 
where $k_\Ga$ is the smallest positive integer such that
$\begin{bmatrix} 1+\varepsilon\,x_1 &
  -x_3+\varepsilon\,x_2\\ x_3+\varepsilon\,x_2 &
  1-\varepsilon\,x_1 \end{bmatrix}^{k_\Ga}\in \Ga$. Furthermore, there
exists $\kappa>0$ such that as $s\ra+\infty$
$$
\card\{\alpha\in(\Ga\cdot \alpha_0)\cap \ZZ_p\;:\;h(\alpha)\leq s\}
= \frac{2\;p^2\;k_\Ga}{(p+1)^2\;[\Ga^1_\OOO:\Ga]}\; s +\bigO(s^{1-\kappa})\;.
$$
\etheo

\dem The group $\OOO^\times$ of invertible elements of $\OOO$ is
$$
\OOO^\times=\big\{x\in\OOO\;:\; \redn(x)\in p^\ZZ\big\}\;.
$$ 
The centre of $\OOO^\times$ is $Z(\OOO^\times)=\{\pm p^n\;:\;n\in
\ZZ\}$ and the centre  of $\OOO^1$ is $Z(\OOO^1)=\{\pm 1\}$. We
identify $\OOO^1/Z(\OOO^1)$ with its image in $\OOO^\times/
Z(\OOO^\times)$.  The quotient group $\OOO^\times/Z(\OOO^\times)$ is a
free group on $s=\frac{p+1}{2}$ generators $\ga_1,\ga_2,\dots, \ga_s$,
which are the images modulo $Z(\OOO^\times)$ of some elements of
$\OOO$ of norm $p$, see for instance
\cite[Coro.~2.1.11]{Lubotzky94}.\footnote{The group
  $\OOO^\times/Z(\OOO^\times)$ is denoted by $\Lambda(2)$ in
  \cite[page 11]{Lubotzky94}.}

Since $\redn(p)=p^2$, any reduced word of even length in $S=
\{\ga_1^\pm, \ga_2^\pm, \dots ,\ga_s^\pm\}$ belongs to
$\OOO^1/Z(\OOO^1)$. Two distinct elements in $S$ differ by a reduced
word of length $2$, and $\ga_1$ does not belong to $\OOO^1/Z(\OOO^1)$.
Hence $\{1,\ga_1\}$ is a system of left coset representatives of
$\OOO^1/Z(\OOO^1)$ in $\OOO^\times/Z(\OOO^\times)$, and the index of
$\OOO^1/Z(\OOO^1)$ in $\OOO^\times/Z(\OOO^\times)$ is 
\begin{equation}\label{eq:indexOOOuntimes}
[\OOO^\times/Z(\OOO^\times):\OOO^1/Z(\OOO^1)] = 2\;.
\end{equation}

Let
$$
g_0=
\begin{pmatrix} 
\frac{1+\varepsilon\,x_1}{p} & \frac{-x_3+\varepsilon\,x_2}{p}\\ 
\frac{x_3+\varepsilon\,x_2}{p} & \frac{1-\varepsilon\,x_1}{p} 
\end{pmatrix}\;.
$$ By the definition of the isomorphism $\theta$ in Equation
\eqref{eq:defisigma} (with $\sqrt{a}=\sqrt{b}=\varepsilon$) and of the
integers $x_1,x_2,x_3$, the element $g_0$ belongs to $\theta(\OOO)$
since $x_1,x_2,x_3$ are even (and $p$ is odd), and $\det g_0=1$. Hence
$g_0\in \theta(\OOO^1)$. Its fixed points for its action by homography
on $\PP^1(\QQ_p)$ are, by an easy computation,
$$
\frac{\varepsilon\, x_1\pm\sqrt{1-\,p^2}}{x_3+\varepsilon \,x_2}
\;.
$$
In particular, $\alpha_0$ is one of these two fixed points. Note that
$\operatorname{tr} g_0=\frac{2}{p}$, hence $|v_p(\operatorname{tr}
g_0)|=1$, and  the image $[g_0]$ of $g_0$ in $\PGL_2(\QQ_p)$ is a
primitive loxodromic element of $\Ga^1_\OOO$.

Let us define 
$$
\ga_0=[g_0]^{u\;k_\Ga}
$$
where $u\in\{\pm 1\}$ is chosen so that $\ga_0^-=\alpha_0$ and
where $k_\Ga$ is defined in the statement of Theorem
\ref{theo:padicequidistribexemplx}.  Since $\Ga$ has finite index in
$\Ga^1_\OOO$, some power of $[g_0]$ does belong to $\Ga$, hence
$k_\Ga$ exists (and note that $k_\Ga=1$ if $\Ga=\Ga^1_\OOO$). By
the minimality of $k_\Ga$, the element $\ga_0$ is a primitive
loxodromic element of $\Ga$.  We will apply Theorem
\ref{theo:equidloxofix} to this $\ga_0$.

The algebra isomorphism $\theta$ induces a group isomorphism from
$\OOO^\times/ Z(\OOO^\times)$ onto its image in $\PGL_2(\QQ_p)$, that
we denote by $\Ga^\times_\OOO$.\footnote{This group is denoted by
  $\Ga(2)$ in \cite[page 95]{Lubotzky94}.} By
\cite[Lem.~7.4.1]{Lubotzky94}, the group $\Ga^\times_\OOO$ acts simply
transitively on the vertices of the Bruhat-Tits tree $\XX_p$. 

In particular, $\Ga^1_\OOO$ acts freely on $\XX_p$, and by Equation
\eqref{eq:indexOOOuntimes}, we have
\begin{align} 
\Vol(\Ga^1_\OOO\dbs\XX_p)& =
[\Ga^\times_\OOO:\Ga^1_\OOO]\;\Vol(\Ga^\times_\OOO\dbs\XX_p)
\nonumber\\ & =[\OOO^\times/Z(\OOO^\times):\OOO^1/Z(\OOO^1)]\;
\card\big(\Ga^\times_\OOO\bs V\XX_p\big)= 2\;.
\label{eq:calTVolpadic}
\end{align}

Again since $\Ga^1_\OOO$ (hence $\Ga$) acts freely on $\XX_p$, since
$\ga_0$ is primitive loxodromic in $\Ga$, and by Equation
\eqref{eq:translatlength}, we have
\begin{align} 
\Vol(\Ga_{\ga_0^-}\dbs \Ax_{\ga_0})& = \card
\big(\Ga_{\ga_0^-}\bs V\!\Ax_{\ga_0}\big)=\lambda(\ga_0)\nonumber
\\ & =k_\Ga\;\lambda([g_0])=2\;k_\Ga\;|v_p(\operatorname{tr}g_0)|
=2\;k_\Ga\;.\label{eq:calVolAxpadic}
\end{align}
Using Equations \eqref{eq:calTVolpadic} and \eqref{eq:calVolAxpadic},
the result now follows from Theorem \ref{theo:equidloxofix}.
\cqfd

\chapter{Equidistribution and counting of crossratios}
\label{sec:crossrat}

We use the same notation as in Chapter \ref{sec:nonarchquadratequid}:
$K_v$ is a non-Archimedean local field, with valuation $v$, valuation
ring $\OOO_v$, choice of uniformiser $\pi_v$, residual field $k_v$ of
order $q_v$, and $\XX_v$ is the Bruhat-Tits tree of $(\PGL_2,K_v)$.
Let $\Ga$ be a lattice in $\PGL_2(K_v)$.

In this Chapter, we give counting and equidistribution results in
$K_v=\partial_\infty\XX_v-\{\infty\}$ of orbit points under $\Ga$,
using a complexity defined using crossratios, which is different from
the one in Chapter \ref{sec:nonarchquadratequid}. We refer to
\cite{ParPau14AFST} for the development when $K_v$ is $\RR$ or $\CC$
with its standard absolute value.

Recall that the {\em crossratio}\index{crossratio} of four pairwise
distinct points $a,b,c,d$ in $\PP_1(K_v)=K_v\cup\{\infty\}$ is
$$
\gls{crossratio}=\frac{(c-a)\,(d-b)}{(c-b)\,(d-a)}\in (K_v)^\times\;,
$$
with the standard conventions when one of the points is $\infty$.
Adopting Ahlfors's terminology in the complex case, the {\em absolute
  crossratio}\index{absolute!crossratio}\index{crossratio!absolute} of
four pairwise distinct points $a,b,c,d\in\PP_1(K_v)=K_v\cup\{\infty\}$
is
$$
\gls{absolcrossratio}=|[a,b,c,d]|_v=
\frac{|c-a|_v\,|d-b|_v}{|c-b|_v\,|d-a|_v}\,,
$$ 
with conventions analogous to the previous ones when one of the
points is $\infty$. As in the classical case, the crossratio and the
absolute crossratio are invariant under the diagonal projective action
of $\GL_2(K_v)$ on the set of quadruples of pairwise distinct points
in $\PP_1(K_v)$.\footnote{The logarithm in base $q_v$ of this absolute
  crossratio is up to the order equal to the (logarithmic) crossratio
  $\ldbrack  \cdot,\cdot,\cdot,\cdot\rdbrack$ introduced in
  Section \ref{subsec:trees}: More precisely, if
  $\xi_1,\xi_2,\xi_3,\xi_4$ are pairwise distinct points in the
  boundary of $\XX_v$, then
$$
\ldbrack\xi_1,\xi_2,\xi_3,\xi_4\rdbrack =
\log_{q_v}|\xi_1,\xi_4,\xi_3,\xi_2|_v\;.
$$
We will not use this relationship in this book.}

\section{Counting and equidistribution of 
crossratios of loxodromic fixed points}
\label{subsec:crossratioloxofix}

Let $\alpha,\beta\in K_v$ be loxodromic fixed points of $\Ga$. Recall
that $\alpha^\sigma, \beta^\sigma$ is the other fixed point of a
loxodromic element of $\Ga$ fixing $\alpha,\beta$, respectively.  The
{\em relative height}\index{relative height} of $\beta$ with respect
to $\alpha$ is\footnote{The factor $|\alpha-\alpha^\sigma|_v$ in the
  denominator, that did not appear in \cite{ParPau14AFST} in the
  analogous definition for the case when $K_v$ is $\RR$ or $\CC$, is
  there in order to simplify the statements below.}
$$
\gls{relheightloxo} =
\frac{1}{|\alpha-\alpha^\sigma|_v\,|\beta-\beta^\sigma|_v}\;
\max\big\{|\beta-\alpha|_v\,|\beta^\sigma-\alpha^\sigma|_v,\;
|\beta-\alpha^\sigma|_v\,|\beta^\sigma-\alpha|_v\big\}\;.
$$
When $\beta\notin\{\alpha,\alpha^\sigma\}$, we have
\begin{equation}\label{eq:relhei}
h_{\alpha}(\beta)=\max\{|\alpha,\beta,\beta^\sigma,\alpha^\sigma |_v,\;
|\alpha,\beta^\sigma,\beta,\alpha^\sigma |_v \}=\frac{1}
{\min\{|\alpha,\beta,\alpha^\sigma,\beta^\sigma |_v,\;
|\alpha,\beta^\sigma,\alpha^\sigma,\beta |_v \}}\,.
\end{equation}
We will use the relative height as a complexity when $\beta$ varies in
a given orbit of $\Ga$ (and $\alpha$ is fixed).

The following properties of relative heights are easy to check using
the definitions, the invariance properties of the crossratio and
Equation \eqref{eq:gacommutgal}.

\blemm\label{lem:proprirelatheight} Let $\alpha,\beta\in K_v$ be
loxodromic fixed points of $\Ga$. Then

\smallskip\noindent(1) $h_{\alpha^\rho}(\beta^\tau)=
h_{\alpha}(\beta)$ for all $\rho,\tau\in\{\id,\sigma\}$.

\smallskip\noindent(2) If $\beta \in\{\alpha,\alpha^\sigma\}$, then
$h_{\alpha}(\beta)=1$.

\smallskip\noindent(3) $h_{\ga\cdot\alpha}(\ga \cdot\beta)= h_{\alpha}(\beta)$ for
every $\ga\in\Ga$.  

\smallskip\noindent(4) $h_{\alpha}(\ga\cdot \beta)= h_{\alpha}(\beta)$ for
every $\ga\in\stab_{\Ga}(\{\alpha,\alpha^\sigma\})$.  \cqfd
\elemm 

The following result relates the relative height of two loxodromic
fixed points with the distance between the two translation axes.

\bprop \label{prop:relatheighlongcomperp}
Let $\alpha,\beta\in K_v$ be loxodromic fixed points of
$\Ga$ such that $\beta\notin\{\alpha,\alpha^\sigma\}$. Then 
$$
h_\alpha(\beta)=q_v^{\;d(]\alpha,\alpha^\sigma[\,,\,]\beta,\beta^\sigma[)}\,.
$$
\eprop

In particular, we have $h_\alpha(\beta)> 1$ if and only if the
geodesic lines $]\alpha,\alpha^\sigma[$ and $]\beta,\beta^\sigma[$ in
$\XX_v$ are disjoint, and $h_\alpha(\beta)=1$ otherwise (using Lemma
\ref{lem:proprirelatheight} (2) when
$\beta\in\{\alpha,\alpha^\sigma\}$).

\medskip 
\dem 
Up to replacing $\alpha,\beta,\alpha^\sigma, \beta^\sigma$ by their
images under a large enough power $\ga$ of a loxodromic element in $\Ga$
with attracting fixed point in $\OOO_v$, we may assume that these four
points belong to $\OOO_v$. Note that $\ga$ exists since
$\Lambda\Ga=\partial_\infty \XX_v$, and it preserves the relative
height by Lemma \ref{lem:proprirelatheight} (3) as well as the
distances between translation axes.

Let $A=\;]\alpha,\alpha^\sigma[$ and $B=\;]\beta,\beta^\sigma[$. Let
$u$ be the closest point to $*_v$ on $A$, so that
$$
v(\alpha-\alpha^\sigma)= d(u,*_v)\;.
$$

We will consider five configurations.

\begin{center}
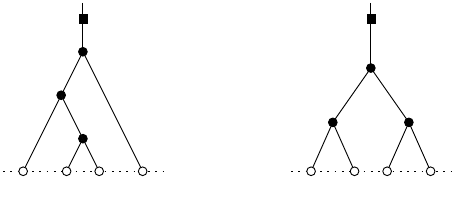
\end{center}

\noindent{\bf Case 1. } First assume that $A$ and $B$ are
disjoint. Let $[a,b]$ be the common perpendicular from $A$ to $B$,
with $a\in A$, so that, by the geometric interpretation of elements in
$\OOO_v$ given in Section \ref{subsec:BruhatTitstrees}, we have
$$
d(A,B)=d(a,b)\;.
$$

First assume that $u\neq a$. Up to exchanging $\alpha,\alpha^\sigma$
(which does not change $d(A,B)$ or $h_\alpha(\beta)$ by Lemma
\ref{lem:proprirelatheight} (1)), we may assume that $a\in
[u,\alpha[\,$. Then (see the picture on the left above),
$$
v(\beta-\beta^\sigma)=d(b,*_v),\;\;\;\;
v(\alpha-\beta)=v(\alpha-\beta^\sigma)=d(a,*_v)
$$
and
$$
v(\alpha^\sigma-\beta)=v(\alpha^\sigma-\beta^\sigma)=d(u,*_v)\;.
$$
Therefore
\begin{align*}
|\alpha,\beta,\alpha^\sigma,\beta^\sigma |_v&=
|\alpha,\beta^\sigma,\alpha^\sigma,\beta |_v=
q_v^{\;v(\alpha-\beta)+v(\alpha^\sigma-\beta^\sigma)-
v(\alpha-\alpha^\sigma)-v(\beta-\beta^\sigma)}
\\ &=q_v^{d(a,\,*_v)-d(b,\,*_v)}=q_v^{-d(a,b)}=q_v^{-d(A,B)}\;,
\end{align*}
which proves the result by Equation \eqref{eq:relhei}.

\medskip 
Assume on the contrary that $u=a$. Let $u'\in V\XX_v$ be such that
$[a,*_v]\cap [a,b]=[a,u']$. Note that $u'\in [*_v,b]$ since
$\beta,\beta^\sigma\in\OOO_v$.  Then (see the picture on the right
above),
$$
v(\beta-\beta^\sigma)=d(b,*_v),\;\;\;\;
v(\alpha-\beta)=v(\alpha-\beta^\sigma)=
v(\alpha^\sigma-\beta)=v(\alpha^\sigma-\beta^\sigma)=d(u',*_v)\;.
$$
Therefore
\begin{align*}
|\alpha,\beta,\alpha^\sigma,\beta^\sigma |_v&=
|\alpha,\beta^\sigma,\alpha^\sigma,\beta |_v=
q_v^{\;v(\alpha-\beta)+v(\alpha^\sigma-\beta^\sigma)-
v(\alpha-\alpha^\sigma)-v(\beta-\beta^\sigma)}
\\ &=q_v^{2d(u',\,*_v)-d(a,\,*_v)-d(b,\,*_v)}=q_v^{-d(a,\,b)}=q_v^{-d(A,\,B)}\;,
\end{align*}
which proves the result by Equation \eqref{eq:relhei}.

\begin{center}
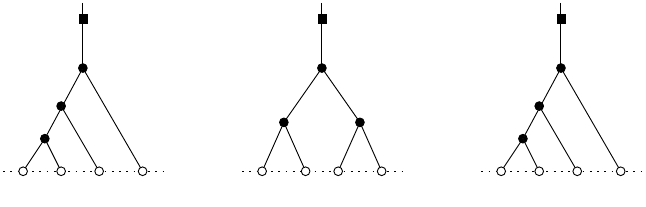
\end{center}

\noindent{\bf Case 2. } Now assume that $A$ and $B$ are
not disjoint, so that
$$
d(A,B)=0\;.
$$ 
Since $\beta\notin\{\alpha,\alpha^\sigma\}$, the intersection
$A\cap B$ is a compact segment $[a,b]$ (possibly with $a=b$) in
$\XX_v$. Up to exchanging $a$ and $b$, $\alpha$ and $\alpha^\sigma$,
as well as $\beta$ and $\beta^\sigma$ (which does not change $d(A,B)$
or $h_\alpha(\beta)$ by Lemma \ref{lem:proprirelatheight} (1)), we may
assume that $\alpha,a,b,\alpha^\sigma$ and $\beta,a,b,\beta^\sigma$
are in this order on $A$ and $B$ respectively, and that $a\in
[u,\alpha[\,$.

\medskip Assume first that $b\in\;]u,\alpha[\,$.  Then (see the
picture on the left above),
$$
v(\alpha-\beta)=d(a,*_v),\;\;\;\;
v(\alpha-\beta^\sigma)=v(\beta-\beta^\sigma)=d(b,*_v)
$$
and
$$
v(\beta-\alpha^\sigma)=v(\beta^\sigma-\alpha^\sigma)=d(u,*_v)\;.
$$
Therefore
$$
|\alpha,\beta^\sigma,\alpha^\sigma,\beta |_v=
q_v^{\;v(\alpha-\beta)+v(\alpha^\sigma-\beta^\sigma)-
v(\alpha-\alpha^\sigma)-v(\beta-\beta^\sigma)}
=q_v^{d(a,\,*_v)-d(b,\,*_v)}=q_v^{d(a,b)}\geq 1\;,
$$
and 
$$
|\alpha,\beta,\alpha^\sigma,\beta^\sigma |_v=
q_v^{\;v(\alpha-\beta^\sigma)+v(\alpha^\sigma-\beta)-
v(\alpha-\alpha^\sigma)-v(\beta-\beta^\sigma)}
=q_v^{0}=1=q_v^{-d(A,B)}\;,
$$
which proves the result by Equation \eqref{eq:relhei}.

\medskip Assume  that $b\in\;]u,\alpha^\sigma[\,$.  Then (see
the picture in the middle above),
$$
v(\alpha-\beta)=d(a,*_v),\;\;\;\;
v(\alpha^\sigma-\beta^\sigma)=d(b,*_v)
$$
and
$$
v(\alpha-\beta^\sigma)=v(\beta-\alpha^\sigma)=
v(\beta-\beta^\sigma)=d(u,*_v)\;.
$$
Therefore
\begin{align*}
|\alpha,\beta^\sigma,\alpha^\sigma,\beta |_v&=
q_v^{\;v(\alpha-\beta)+v(\alpha^\sigma-\beta^\sigma)-
v(\alpha-\alpha^\sigma)-v(\beta-\beta^\sigma)}\\ &
=q_v^{d(a,\,*_v)+d(b,\,*_v)-2d(u,*_v)}=q_v^{d(a,b)}\geq 1\;,
\end{align*}
and 
$$
|\alpha,\beta,\alpha^\sigma,\beta^\sigma |_v=
q_v^{\;v(\alpha-\beta^\sigma)+v(\alpha^\sigma-\beta)-
v(\alpha-\alpha^\sigma)-v(\beta-\beta^\sigma)}
=q_v^{0}=1=q_v^{-d(A,B)}\;,
$$
which proves the result by Equation \eqref{eq:relhei}.

\medskip Assume at last that $b=u$. Let $u'\in V\XX_v$ be such that
$[b,*_v]\cap[b,\beta^\sigma[\;=[b,u']$. Then (see the picture
on the right above),
$$
v(\alpha-\beta)=d(a,*_v),\;\;\;\;
v(\alpha^\sigma-\beta)=d(u,*_v)
$$
and
$$
v(\alpha-\beta^\sigma)=v(\beta-\beta^\sigma)=
v(\alpha^\sigma-\beta^\sigma)=d(u',*_v)\;.
$$
Therefore
\begin{align*}
|\alpha,\beta^\sigma,\alpha^\sigma,\beta |_v&=
q_v^{\;v(\alpha-\beta)+v(\alpha^\sigma-\beta^\sigma)-
v(\alpha-\alpha^\sigma)-v(\beta-\beta^\sigma)}\\ &
=q_v^{d(a,\,*_v)-d(u,*_v)}=q_v^{d(a,b)}\geq 1\;,
\end{align*}
and 
$$
|\alpha,\beta,\alpha^\sigma,\beta^\sigma |_v=
q_v^{\;v(\alpha-\beta^\sigma)+v(\alpha^\sigma-\beta)-
v(\alpha-\alpha^\sigma)-v(\beta-\beta^\sigma)}
=q_v^{0}=1=q_v^{-d(A,B)}\;,
$$
which proves the result by Equation \eqref{eq:relhei}.
\cqfd

\medskip 
The next result says that the relative height is an appropriate
complexity on a given orbit under $\Ga$ of a loxodromic fixed point,
and that the counting function we will study is well defined. We
denote by $\Ga_\xi$ the stabiliser in $\Ga$ of a point $\xi
\in \partial_\infty\XX_v=\PP_1(K_v)$.

\blemm Let $\alpha,\beta\in K_v$ be loxodromic fixed points of
$\Ga$. Then for every $s>1$, the set
$$
E_s=\{\beta'\in \Ga_\alpha \bs \Ga\cdot\beta\;:\; 
1 < h_\alpha(\beta')\leq s\}
$$
is finite.
\elemm

\dem 
The set $E_s$ is well defined by Lemma \ref{lem:proprirelatheight}
(4).  Recall that a loxodromic fixed point is one of the two points at
infinity of a unique translation axis.  By local finiteness, there
are, up to the action of the stabiliser of a fixed translation axis
$A$, only finitely many images under $\Ga$ of another translation axis
$B$ at distance at most $\frac{\ln s}{\ln q_v}$ from $A$. Since the
stabiliser of $A$ contains the stabiliser of either of its points at
infinity with index at most $2$, the result then follows from
Proposition \ref{prop:relatheighlongcomperp}.  
\cqfd

\medskip 
We now state our main counting and equidistribution result of orbits
of loxodromic fixed points, when the complexity is the relative height
with respect to a fixed loxodromic fixed point.

\btheo \label{theo:mainrelatheight} Let $\Ga$ be a lattice in
$\PGL_2(K_v)^+$. Let $\alpha_0,\beta_0 \in K_v$ be loxodromic fixed
points of $\Ga$. Then for the weak-star convergence of measures on
$K_v-\{\alpha_0,\alpha_0^\sigma\}$, as $s\ra+\infty$,
$$
\frac{(q_v+1)^2\;\Vol(\Ga\dbs \XX_v)}
{2\;q_v^2\;|\alpha_0-\alpha_0^\sigma|_v\;
\Vol(\Ga_{\beta_0}\dbs \;]\beta_0,\beta_0^\sigma[)}\;s^{-1}
\sum_{\beta\in \Ga\cdot\beta_0\;:\;h_{\alpha_0}(\beta)\le s}\Delta_{\beta}
\;\;\weakstar\;\;
\frac{d\haar_{K_v}(z)}{|z-\alpha_0|_v\,|z-\alpha_0^\sigma|_v}\;.
$$
Furthermore, there exists $\kappa>0$ such that, as $s\ra+\infty$,
\begin{align*}
\card&\;\;\Ga_{\alpha_0}\bs\{\beta\in \Ga\cdot\beta_0\;:\;
h_{\alpha_0}(\beta)\le s\} \\ & =\; 
\frac{2\;q_v\;(q_v-1)\;\Vol(\Ga_{\alpha_0}\dbs\;]\alpha_0,\alpha_0^\sigma[)\;
\Vol(\Ga_{\beta_0}\dbs \;]\beta_0,\beta_0^\sigma[)}
{(q_v+1)^2\;\Vol(\Ga\dbs \XX_v)}\;s+\bigO(s^{1-\kappa})\;.
\end{align*}
\etheo

For every $\beta'\in\; ]0,1]$, there
exists $\kappa>0$ such that for every $\psi\in\C^{\beta'}_c(K_v-
\{\alpha_0,\alpha_0^\sigma\})$, where $K_v- \{\alpha_0,
\alpha_0^\sigma\}$ is endowed with the distance-like map
$d_{]\alpha_0, \,\alpha_0^\sigma[}$,\footnote{See Equation
  \eqref{eq:distancelike}.}  there is an error term
in the equidistribution claim of Theorem \ref{theo:mainrelatheight}
when evaluated on $\psi$, of the form $\bigO(s^{-\kappa}
\|\psi\|_{\beta'})$. This result applies for instance if $\psi:K_v-
\{\alpha_0,\alpha_0^\sigma\}\ra\RR$ is locally constant with compact
support, see Remark \ref{rem:locconstholder}.

\medskip 
\dem The proof of the equidistribution claim is similar to the one of
Theorem \ref{theo:equidloxofix}.  We now apply Theorem
\ref{theo:algebrogeometricequid} with $\DD^-:= \;]\alpha_0,
\alpha_0^\sigma[$ and $\DD^+:= \; ]\beta_0,\beta_0^\sigma[\,$.
Since $\Ga$ is contained in $\PGL_2(K_v)^+$, the length spectrum
$L_{\Ga}$ of $\Ga$ is equal to $2\ZZ$. The families $\D^\pm=
(\ga\DD^\pm)_{\ga\in\Ga/\Ga_{\DD^\pm}}$ are locally finite, and
$\|\sigma^-_{\D^+}\|$ is finite and nonzero.

Arguing as in the proof of Theorem
\ref{theo:equidloxofix},\footnote{See Equation
  \eqref{eq:demoequidloxofixThet}.} we have
\begin{equation}\label{eq:dixseptquatre}
\lim_{n\ra+\infty} \;\frac{({q_v}^2-1)(q_v+1)}{2\,{q_v}^3}\;
\frac{\Vol(\Ga\dbs\XX_v)}{\|\sigma^-_{\D^+}\|}\;{q_v}^{-n}
\sum_{\substack{\ga\in \Ga/\Ga_{\DD^+}\\
0< d(\DD^-, \, \ga \DD^+)\leq n}}
\Dirac_{\ga \cdot \beta_0} \;=\;
(\partial^+)_*\wt\sigma^+_{\DD^-}\;
\end{equation}
for the weak-star convergence of measures on $\partial_\infty \XX_v-
\partial_\infty \DD_-$. Futhermore, for every $\beta'\in\; ]0,1]$,
there exists $\kappa>0$ such that for every $\beta$-Hölder-continuous
function $\psi\in\C^{\beta'}_c (\partial_\infty \XX_v-\partial_\infty
\DD_-)$, where $\partial_\infty \XX_v-\partial_\infty \DD_-$ is
endowed with the distance-like map $d_{\DD^-}$, there is an error term
in the equidistribution statement of Equation \eqref{eq:dixseptquatre}
when evaluated on $\psi$, of the form
$\bigO(s^{-\kappa}\|\psi\|_{\beta'})$.

By Proposition \ref{prop:relatheighlongcomperp},
we have
$$
h_{\alpha_0}(\ga\cdot\beta_0)= q_v^{\;d(\DD^-,\,\ga\DD^+)}\,.
$$
By Proposition \ref{prop:mescomputBT} (5), we have 
$$
(\partial^+)_*\wt\sigma^+_{\DD^-}(z)=
\frac{|\alpha_0-\alpha_0^\sigma|_v}{|z-\alpha_0|_v\,|z-\alpha_0^\sigma|_v}
\;d\haar_{K_v}(z)\,
$$ 
for $z$ in the full measure subset $K_v-\{\alpha_0,\alpha_0^\sigma\}$
of $\partial_\infty \XX_v$.  Hence, using the change of variable
$s=q_v^{\;n}$, we have, with the appropriate error term,
$$
\lim_{s\ra+\infty} 
\;\frac{({q_v}^2-1)(q_v+1)}{2\,q_v^3}\;\frac{\Vol(\Ga\dbs\XX_v)}
{|\alpha_0-\alpha_0^\sigma|_v\;\|\sigma^-_{\D^+}\|}\;s^{-1}
\sum_{\substack{\ga\in \Ga/\Ga_{\DD^+}\\
1< h_{\alpha_0}(\ga \cdot \beta_0)\leq s}}
\Dirac_{\ga \cdot \beta_0} \;=\; 
\frac{d\haar_{K_v}(z)}{|z-\alpha_0|_v\,|z-\alpha_0^\sigma|_v}\;\;.
$$

We again denote by $\iota_{\alpha_0}$ the index
$$
\iota_{\alpha_0}=[\Ga_{\{\alpha_0,\,\alpha_0^\sigma\}}:\Ga_{\alpha_0}]\;,
$$ and similarly for $\beta_0$. Since the stabiliser $\Ga_{\beta_0}$
of $\beta_0$ in $\Ga$ has index $\iota_{\beta_0}$ in $\Ga_{\DD^+}$ and
$\Ga/\Ga_{\beta_0}$ identifies with $\Ga\cdot \beta_0$ by the map
$\ga\mapsto \ga\cdot \beta_0$, we have, with the appropriate error
term,
$$
\lim_{s\ra+\infty} \;
\frac{({q_v}^2-1)(q_v+1)}{2\,q_v^3}\;\frac{\Vol(\Ga\dbs\XX_v)}
{|\alpha_0-\alpha_0^\sigma|_v\;\|\sigma^-_{\D^+}\|\;\iota_{\beta_0}}
\;s^{-1}\sum_{\substack{\beta\in  \Ga\cdot \beta_0\\
1< h_{\alpha_0}(\beta)\leq s}}
\Dirac_{\beta} \;=\; 
\frac{d\haar_{K_v}(z)}{|z-\alpha_0|_v\,|z-\alpha_0^\sigma|_v}\;\;.
$$
As in the end of the proof of Theorem \ref{theo:equidloxofix}, we have
$$
\|\sigma^-_{\D^+}\| = \frac{q_v-1}{q_v\;\iota_{\beta_0}}\;
\Vol(\Ga_{\beta_0}\dbs \;]\beta_0,\beta_0^\sigma[)\;.
$$ 
This proves the equidistribution claim, and its error term.

In order to obtain the counting claim, we note that since $\wt
\sigma^+_{\DD^-}$ is invariant under the stabiliser in $\Ga$ of
$\DD^-$, hence under $\Ga_{\alpha_0}$, the measures on both sides of
the equidistribution claim in Theorem \ref{theo:mainrelatheight} are
invariant under $\Ga_{\alpha_0}$, see Proposition
\ref{prop:propricomplexity} (5) for the invariance of the right hand
side.  By Proposition \ref{prop:mescomputBT} (5) and (6), and by the
definition of $\iota_{\alpha_0}$, we have
\begin{align}
\int_{\Ga_{\alpha_0}\bs (K_v-\{\alpha_0,\,\alpha_0^\sigma\})}\;
\frac{d\haar_{K_v}(z)}{|z-\alpha_0|_v\,|z-\alpha_0^\sigma|_v}
&=\frac{\iota_{\alpha_0}}{|\alpha_0-\alpha_0^\sigma|_v}\;
\int_{\Ga_{\DD^-}\bs\normalout\DD^-} d\wt\sigma^+_{\DD^-}\nonumber\\ & =
\frac{(q_v-1)\;\iota_{\alpha_0}\;\Vol(\Ga_{\DD^-}\dbs\DD^-)}
{q_v\;|\alpha_0-\alpha_0^\sigma|_v}\nonumber\\ & 
=\frac{(q_v-1)\;\Vol(\Ga_{\alpha_0}\dbs\;]\alpha_0,\alpha_0^\sigma[)}
{q_v\;|\alpha_0-\alpha_0^\sigma|_v}\;.\label{eq:calcdomfondmodGalph}
\end{align}
The counting claim follows by evaluating the equidistribution claim on
the characteristic function $\psi$ of a compact-open fundamental
domain for the action of $\Ga_{\alpha_0}$ on $K_v-\{\alpha_0,
\alpha_0^\sigma\}$.  This characteristic function is locally constant,
hence $\beta'$-Hölder-continuous for the distance-like function
$d_{\DD^-}$, as seen end of Section \ref{subsec:holdercont}. 
\cqfd

\section{Counting and equidistribution of crossratios of 
quadratic irrationals}
\label{subsec:countcrossratquadirrat}

In this Section, we give two arithmetic applications of Theorem
\ref{theo:mainrelatheight}.

\medskip 
Let us first consider an application in positive characteristic.  Let
$K$ be a (global) function field over $\FF_q$, let $v$ be a
(normalised discrete) valuation of $K$, let $K_v$ be the associated
completion of $K$ and let $R_v$ be the affine function ring associated
with $v$.\footnote{See Section \ref{subsec:valuedfields}.}  Given two
quadratic irrationals $\alpha,\beta \in K_v$ over $K$, with Galois
conjugates $\alpha^\sigma, \beta^\sigma$ respectively, such that
$\beta\notin \{\alpha, \alpha^\sigma\}$, we define the {\em relative
  height}\index{relative height} of $\beta$ with respect to $\alpha$
by
\begin{equation}\label{eq:relheiarith}
\gls{relheightquadirrat}=\frac{1}
{\min\{|\alpha,\beta,\alpha^\sigma,\beta^\sigma |_v,\;
|\alpha,\beta^\sigma,\alpha^\sigma,\beta |_v \}}\,.
\end{equation}

\medskip 
The following result says that the orbit of any quadratic
irrational in $K_v$ over $K$, by homographies under a given finite
index subgroup of the modular group $\PGL_2(R_v)$, equidistributes,
when its complexity is given by the relative height with respect to
another fixed quadratic irrational $\alpha_0$.  The limit measure is
absolutely continuous with respect to the Haar measure on $K_v$ and it
is invariant under the stabiliser of $\alpha_0$ in $\PGL_2(R_v)$ by
Proposition \ref{prop:propricomplexity} (5).

\btheo \label{theo:caracposequidquadirrrelathei} Let $G$ be a finite
index subgroup of $\GL_2(R_v)$.  Let $\alpha_0,\beta_0\in K_v$ be
quadratic irrationals over $K$. Let $g_0,h_0$ be elements in $G$
fixing $\alpha_0,\beta_0$ with $v(\operatorname{tr} g_0),
v(\operatorname{tr} h_0)\neq 0$, and let $m_0,n_0$ be the index of
$g_0^\ZZ,h_0^\ZZ$ in the stabiliser of $\alpha_0,\beta_0$ in $G$
respectively. Then, as $s\to+\infty$,
\begin{align*}
\frac{(q_v+1)^2\;\zeta_K(-1)\;n_0\;[\GL_2(R_v):G]}
{2\;q_v^2\;(q-1)\;|\alpha_0-\alpha_0^\sigma|_v\;
|v(\operatorname{tr} h_0)|}\;&s^{-1}
\sum_{\beta\in G\cdot\beta_0\;:\;h_{\alpha_0}(\beta)\le s}\Delta_{\beta}
\\ &\weakstar\;\;
\frac{d\haar_{K_v}(z)}{|z-\alpha_0|_v\,|z-\alpha_0^\sigma|_v}\;,
\end{align*}
and there exists $\kappa>0$ such that, as $s\to+\infty$,
\begin{align*}
\card&\;\;\Ga_{\alpha_0}\bs\{\beta\in G\cdot\beta_0\;:\;
h_{\alpha_0}(\beta)\le s\}\\&=\;
\frac{4\;q_v\;(q_v-1)\;(q-1)\;
|v(\operatorname{tr} g_0)|\;|v(\operatorname{tr} h_0)|\;|Z(G)|}
{(q_v+1)^2\;\zeta_K(-1)\;m_0\;n_0\;[\GL_2(R_v):G]}
\;s +\bigO(s^{1-\kappa})\;.
\end{align*}
\etheo

\medskip
\dem This follows, as in the proof of Theorem
\ref{theo:caracposequidquadirr}, from Theorem
\ref{theo:mainrelatheight} using Equations \eqref{eq:covoltotMertens}
and \eqref{eq:covolcycle}, as well as Equation
\eqref{eq:calcdomfondmodGalph} for the counting claim.  
\cqfd

\smallskip
\bexem (1) Theorem \ref{theo:applitwoquadirratdintro} in the
Introduction follows from this result, by taking $K=\FF_q(Y)$ and
$v=v_\infty$ (so that $q_v=q$), and by using Equation
\eqref{zetamoinsun} in order to simplify the constant.

\medskip 
(2) If $G_I$ is the Hecke congruence subgroup associated with
a nonzero ideal $I$ of $R_v$ (see Equation
\eqref{eq:heckecongruencesubgroup}), using Lemma
\ref{lem:indexHeckesubgrou}, we have, as $s\to+\infty$,
$$
\frac{(q_v+1)^2\;\zeta_K(-1)\;n_0\;N(I)\prod_{\ppp|I}(1+\frac{1}{N(\ppp)})}
{2\;q_v^2\;(q-1)\;|\alpha_0-\alpha_0^\sigma|_v\;
|v(\operatorname{tr} h_0)|}\;s^{-1}
\sum_{\beta\in G_I\cdot\beta_0\;:\;h_{\alpha_0}(\beta)\le s} \Delta_{\beta}\;\;
\weakstar\;\;
\frac{d\haar_{K_v}(z)}{|z-\alpha_0|_v\,|z-\alpha_0^\sigma|_v}\;.
$$
\eexem

\bigskip 
The second arithmetic application of Theorem
\ref{theo:mainrelatheight} is in $\QQ_p$. We use the notation of
Section \ref{subsec:quadirratzerochar}.

Let $p\in\NN$ be a positive rational prime with $p\equiv 1 \mod 4$
such that $\frac{p^2-1}4$ is not of the form $4^a(8b+7)$ for $a,b\in
\NN$ (for instance $p=5$). Let $\varepsilon$ be a square root of $-1$
in $\QQ_p$. Let $x_1,x_2,x_3\in 2\ZZ$ be such that $p^2-1= {x_1}^2 +
{x_2}^2 + {x_3}^2$. We again consider
$$ 
\alpha_0= 
\frac{\varepsilon\,x_1+\sqrt{1-\,p^2}}{x_3+\varepsilon \,x_2} \;,
$$
which is a quadratic irrational in $\QQ_p$ over $\QQ(\varepsilon)$. We
denote by $\alpha^\sigma$ the Galois conjugate of a quadratic
irrational $\alpha$ in $\QQ_p$ over $\QQ(\varepsilon)$, and by
\begin{equation}\label{eq:relheiarithpadic}
h_{\alpha}(\beta)=\frac{1}
{\min\{|\alpha,\beta,\alpha^\sigma,\beta^\sigma |_p,\;
|\alpha,\beta^\sigma,\alpha^\sigma,\beta |_p\}}
\end{equation}
the {\em relative height}\index{relative height} of a quadratic
irrational $\beta$ in $\QQ_p$ over $\QQ(\varepsilon)$ with respect to
$\alpha$, assuming that $\beta\notin\{\alpha,\alpha^\sigma\}$. We again
consider Hamilton's quaternion algebra $D=(\frac{-1,-1}\QQ)$ over
$\QQ$ and its $\ZZ\big[\frac{1}{p}\big]$-order
$$
\OOO=\big\{y\in\ZZ\big[{\textstyle \frac{1}{p}}\big]+
\ZZ\big[{\textstyle\frac{1}{p}}\big]\;i+
\ZZ\big[{\textstyle\frac{1}{p}}\big]\;j+
\ZZ\big[{\textstyle\frac{1}{p}}\big]\;k\;:\; 
y\equiv 1\mod 2\big\}\;.
$$

\medskip 
The following result says that the orbit of $\alpha_0$ in $\QQ_p$ by
homographies under a given finite index subgroup of the arithmetic
group $\Ga^1_\OOO$ (defined in Section \ref{subsec:quadirratzerochar})
equidistributes, when its complexity is given by the relative height
with respect to $\alpha_0$, to a measure absolutely continuous with
respect to the Haar measure on $\QQ_p$.

\btheo \label{theo:padicquidquadirrrelathei} With the above notation,
let $\Ga$ be a finite index subgroup of $\Ga^1_\OOO$.  Then, as
$s\to+\infty$,
$$
\frac{(p+1)^2\;[\Ga^1_\OOO:\Ga]}
{2\;p^2\;k_\Ga\;|\alpha_0-\alpha_0^\sigma|_p}\;s^{-1}
\sum_{\alpha\in \Ga\cdot\alpha_0\;:\;h_{\alpha_0}(\alpha)\le s}\Delta_{\alpha}\;\;
\weakstar\;\;
\frac{d\haar_{\QQ_p}(z)}{|z-\alpha_0|_p\,|z-\alpha_0^\sigma|_p}\;,
$$
where $k_\Ga$ is the smallest positive integer such that
$\begin{bmatrix} 1+\varepsilon\,x_1 & -x_3+\varepsilon\,x_2\\
x_3+\varepsilon\,x_2 & 1-\varepsilon\,x_1 \end{bmatrix}^{k_\Ga}\in
\Ga$. Furthermore, there exists $\kappa>0$ such that, as $s\to+\infty$,
$$
\card\;\;\Ga_{\alpha_0}\bs\{\alpha\in \Ga\cdot\alpha_0\;:\;
h_{\alpha_0}(\alpha)\le s\}=\frac{4\;p\;(p-1)\;(k_\Ga)^2}
{(p+1)^2\;[\Ga^1_\OOO:\Ga]} \;s +\bigO(s^{1-\kappa})\;.
$$
\etheo

\dem This follows, as in the proof of Theorem
\ref{theo:padicequidistribexemplx}, from Theorem
\ref{theo:mainrelatheight} using Equations \eqref{eq:calTVolpadic} and
\eqref{eq:calVolAxpadic}, as well as Equation
\eqref{eq:calcdomfondmodGalph} and again Equation
\eqref{eq:calVolAxpadic} for the counting claim.  
\cqfd

\chapter{Equidistribution and counting of 
integral representations by quadratic norm forms}
\label{subsec:integralrepresentation}

In the final Chapter of this book, we give another arithmetic
equidistribution and counting result of rational elements in
non-Archimedean local fields of positive characteristic, again using
our geometric equidistribution and counting results of common
perpendiculars in trees summarised in Section \ref{subsec:locconst}.
We use here a complexity defined using the norm forms associated with
fixed quadratic irrationals.  In particular, the complexity in this
Chapter is different from that used in the Mertens type of results in
Section \ref{subsec:mertens}. We refer for instance to \cite[\S
  5.3]{ParPau14AFST} for motivations and results in the Archimedean
case, and also to \cite{GorPau14} for higher dimensional norm forms.

\medskip 
Let $K$ be a (global) function field over $\FF_q$ of genus $g$, let
$v$ be a (normalised discrete) valuation of $K$, let $K_v$ be the
associated completion of $K$ and let $R_v$ be the affine function ring
associated with $v$.\footnote{See Section \ref{subsec:valuedfields}.}
Let $\alpha\in K_v$ be a quadratic irrational over $K$. The {\it norm
  form}\index{norm!form} $\n_\alpha$ associated with $\alpha$ is the
quadratic form $K\times K\ra K$ defined by
$$
(x,y)\mapsto \n(x-y\alpha)=(x-y\alpha)(x-y\alpha^\sigma)=
x^2-xy\tr(\alpha)+y^2\n(\alpha)\;.
$$ 
See Proposition \ref{prop:propricomplexity} (3) for elementary
transformation properties under elements of $\GL_2(R_v)$ of this norm
form.

\medskip A pair $(x,y)\in R_v\times R_v$ is an {\em integral
  representation} of an element $z\in K$ by the quadratic norm form
$\n_\alpha$ if $\n_\alpha(x,y)=z$.  The following result describes the
projective equidistribution as $s\ra+\infty$ of the integral
representations by $\n_\alpha$ of elements with absolute value at most
$s$. For every $(x_0,y_0)\in R_v\times R_v$, let $H_{(x_0,y_0)}$ be
the stabiliser of $(x_0,y_0)$ for the linear action of any subgroup
$H$ of $\GL_2(R_v)$ on $R_v\times R_v$. We use the notation $N\langle
x_0,y_0\rangle$ for the norm of the ideal $\langle x_0,y_0\rangle$
generated by $x_0,y_0$ (see Section \ref{subsec:valuedfields}) and the
notation $m_{v,\,x_0,\,y_0}$ introduced above Theorem
\ref{theo:Mertensfunctionfieldfissgp}.

\btheo \label{theo:mainnormform} Let $G$ be a finite index
subgroup of $\GL_2(R_v)$, let $\alpha\in K_v$ be a quadratic
irrational over $K$, and let $(x_0,y_0)\in R_v\times R_v-
\{(0,0)\}$. Let 
$$
c'=\frac{(q_v-1)\;(q_v+1)^2\;\zeta_K(-1)\;m_{v,\,x_0,\,y_0}\;
(N\langle x_0,\,y_0\rangle)^2\;
[\GL_2(R_v):G]}{q_v^3\;(q-1)\;q^{g-1}\;
[\GL_2(R_v)_{(x_0,\,y_0)}:G_{(x_0,\,y_0)}]}\;.
$$
Then for the weak-star convergence of measures on
$K_v-\{\alpha,\alpha^\sigma\}$, we have
$$
\lim_{s\ra+\infty}\; c' \;s^{-1} \sum_{(x,\,y)\in G(x_0,\,y_0),\;
  |\n(x-y\alpha)|_v \leq s} \Delta_{\frac xy}\;=\;
\frac{d\haar_{K_v}(z)}{|z-\alpha|_v\,|z-\alpha^\sigma|_v}\;.
$$
\etheo

For every $\beta\in\; ]0,1]$, there exists $\kappa>0$ such that for
every $\psi\in\C^{\beta}_c(K_v-\{\alpha,\alpha^\sigma\})$, where
$K_v-\{\alpha,\alpha^\sigma\}$ is endowed with the distance-like
map $d_{]\alpha,\,\alpha^\sigma[}$,\footnote{See Equation
  \eqref{eq:distancelike}.}  there is an error term in the
equidistribution claim of Theorem \ref{theo:mainnormform} when
evaluated on $\psi$, of the form $\bigO(s^{-\kappa}\|\psi\|_{\beta})$.
This holds  for instance if $\psi:K_v- \{\alpha,
\alpha^\sigma\}\ra\RR$ is locally constant with compact support (see
Remark \ref{rem:locconstholder}).

\bexems (1) Let $(x_0,y_0)=(1,0)$, $K=\FF_q(Y)$ and $v=v_\infty$ (so
that $g=0$ and $q_v=q$). Theorem \ref{theo:normformintro} in the Introduction
follows from Theorem \ref{theo:mainnormform}, using Equations
\eqref{zetamoinsun} and \eqref{eq:mvxoyocasspecial} to simplify the
constant $c'$.

\smallskip\noindent
(2) Let $(x_0,y_0)=(1,0)$ and let $G=G_I$ be the
Hecke congruence subgroup of $\GL_2(R_v)$ defined in Equation
\eqref{eq:heckecongruencesubgroup}. The index in $[\GL_2(R_v):G_I]$ is
given by Lemma \ref{lem:indexHeckesubgrou} and $G_I$ satisfies
$(G_I)_{(1,0)}= \GL_2(R_v)_{(1,0)}$. For every nonzero ideal $I$ of
$R_v$, for the weak-star convergence of measures on
$K_v-\{\alpha,\alpha^\sigma\}$, we have
$$
\lim_{s\ra+\infty}\; c_I \;s^{-1} 
\sum_{(x,\,y)\in R_v\times I,\; \langle x,\,y\rangle=R_v,\;
  |\n(x-y\alpha)|_v \leq s} \Delta_{\frac xy}\;=\;
\frac{d\haar_{K_v}(z)}{|z-\alpha|_v\,|z-\alpha^\sigma|_v}\;,
$$
where 
$$
c_I=\frac{(q_v-1)\;(q_v+1)^2\;\zeta_K(-1)\;
N(I)\prod_{\ppp|I}(1+\frac{1}{N(\ppp)})}
{q_v^3\;q^{g-1}}\;.
$$

\smallskip\noindent (3) This third example is only interesting when the
ideal class number is larger than $1$.  Given any fractional ideal
$\mmm$ of $R_v$, taking $(x_0,y_0)\in R_v\times R_v$ such that the
fractional ideals $\langle x_0,y_0\rangle$ and $\mmm$ have the same
ideal class and $G=\GL_2(R_v)$, using the change of variables
$s\mapsto s N(\mmm)^2$ in the statement of Theorem
\ref{theo:mainnormform}, for the weak-star convergence of
measures on $K_v-\{\alpha,\alpha^\sigma\}$, with the same error term
as for Theorem \ref{theo:mainnormform}, we have
$$
\lim_{s\ra+\infty}\; c_\mmm \;s^{-1} 
\sum_{(x,\,y)\in \mmm \times \mmm,\; \langle x,\,y\rangle=\mmm,\;
  N(\mmm)^{-2}|\n(x-y\alpha)|_v \leq s} \Delta_{\frac xy}\;=\;
\frac{d\haar_{K_v}(z)}{|z-\alpha|_v\,|z-\alpha^\sigma|_v}\;,
$$
where 
$$
c_\mmm=\frac{(q_v-1)\;(q_v+1)^2\;\zeta_K(-1)\;m_{v,\,x_0,\,y_0}}
{q_v^3\;(q-1)\;q^{g-1}}\;.
$$
\eexems

\medskip 
Before proving Theorem \ref{theo:mainnormform}, let us give a counting
result which follows from it.  Any subgroup of $G$ acts on the left on
any orbit of $G$. Furthermore, the stabiliser $G_\alpha$ of $\alpha$
in $G$ preserves the map $(x,y) \mapsto |\n(x-y\alpha)|_v$, by
Proposition \ref{prop:propricomplexity} (3). We may then define a {\it
  counting function}\index{counting function} $\Psi'(s) =
\Psi'_{G,\,\alpha,\,x_0,y_0}(s)$ of elements of $R_v\times R_v$ in a
linear orbit under a finite index subgroup $G$ of $\GL_2(R_v)$ on
which the absolute value of the norm form associated with $\alpha$ is
at most $s$, as
$$ 
\Psi'(s) = \card\;\;G_{\alpha}\,\bs 
\big\{(x,\,y)\in G(x_0,\,y_0),\;\; |\n(x-y\alpha)|_v \leq s\}\;.
$$

\bcoro \label{coro:represintegralcount} Let $G$ be a finite index
subgroup of $\GL_2(R_v)$, let $\alpha\in K_v$ be a quadratic
irrational over $K$,  and let $(x_0,y_0)\in R_v\times R_v-\{(0,0)\}$.
Let $g_0\in G_{\alpha}$ with $v(\operatorname{tr} g_0)\neq
0$ and let $m_0$ be the index of $g_0^\ZZ$ in $G_{\alpha}$. Let
$$
c''=\frac{2\;q_v^2\;(q-1)\;q^{g-1}\;|Z(G)|\;|v(\operatorname{tr} g_0)|
\;[\GL_2(R_v)_{(x_0,\,y_0)}:G_{(x_0,\,y_0)}]}
{(q_v+1)^2\;\zeta_K(-1)\;|\alpha-\alpha^\sigma|_v\;m_0\;m_{v,\,x_0,\,y_0}
\;(N\langle x_0,\,y_0\rangle )^2\;[\GL_2(R_v):G]}\;.
$$
Then there exists $\kappa>0$ such that, as $s\to+\infty$,
$$
\Psi'(s)=c''\;s +\bigO(s^{1-\kappa})\;.
$$
\ecoro

\dem Using Equation \eqref{eq:calcdomfondmodGalph} (with $\Ga$ the
image of $G$ in $\PGL_2(R_v)$) and Equation \eqref{eq:covolcycle}, we
have
$$
\int_{G_{\alpha}\bs (K_v-\{\alpha,\,\alpha^\sigma\})}\;
\frac{d\haar_{K_v}(z)}{|z-\alpha|_v\,|z-\alpha^\sigma|_v}=
\frac{2\; (q_v-1)\;|Z(G)|\;|v(\operatorname{tr} g_0)|}
{q_v\;|\alpha-\alpha^\sigma|_v\;m_0}\;.
$$
The corollary then follows by applying the equidistribution claim in
Theorem \ref{theo:mainnormform} to the characteristic function of a
compact-open fundamental domain of $K_v-\{\alpha,\alpha^\sigma\}$
modulo the action by homographies of $G_{\alpha}$. \cqfd

\bexem Let $(x_0,y_0)=(1,0)$, $K=\FF_q(Y)$, $v=v_\infty$ (so that
$g=0$ and $q_v=q$) and $G=\GL_2(\FF_q[Y])$.  Using Equations
\eqref{zetamoinsun} and \eqref{eq:mvxoyocasspecial}, Proposition
\ref{prop:propricomplexity} (1), the change of variable $s=q^t$ and
the fact that $|Z(G)|=q-1$ in order to simplify the constant $c''$ of
Corollary \ref{coro:represintegralcount}, and recalling the expression
of the absolute value at $\infty$ in terms of the degree from Section
\ref{subsec:valuedfields}, we get the following counting result: For
every integral quadratic irrational $\alpha\in \FF_q((Y^{-1}))$ over
$\FF_q(Y)$, there exists $\kappa>0$ such that, as $t\to+\infty$,
\begin{align*}
\card& \;\;_{{\textstyle\GL_2(\FF_q[Y])_{\alpha}}}\!\Big\bs 
  \bigg\{(x,\,y)\in \FF_q[Y]\times\FF_q[Y]\;:\; 
\begin{array}{l}\langle x,\,y\rangle= \FF_q[Y],\\[1mm]
\deg(x^2-xy\tr (\alpha) +y^2\n(\alpha))\leq t\end{array}\bigg\}
\\&
=\; \frac{2\;q\;(q-1)^3}{m_0\;(q+1)}\;\deg(\operatorname{tr} g_0)\;
q^{-\frac{1}{2}\deg(\tr(\alpha)^2-4\n(\alpha))}\;\;q^t 
+\bigO(q^{t-\kappa})\;,
\end{align*}
where $g_0\in \GL_2(\FF_q[Y])$ fixes $\alpha$ with
$\deg(\operatorname{tr} g_0)\neq 0$ and $m_0$ is the index of
$g_0^\ZZ$ in the stabiliser $\GL_2(\FF_q[Y])_{\alpha}$ of $\alpha$ in
$ \GL_2(\FF_q[Y])$.  
\eexem

\bigskip \noindent {\bf Proof of Theorem \ref{theo:mainnormform}. }
The proof is similar to that of Theorem
\ref{theo:Mertensfunctionfieldfissgp}. Let $r=\frac{x_0}{y_0}\in
K\cup\{\infty\}$. If $y_0= 0$, let $g_r=\id\in\GL_2(K)$, and if
$y_0\neq 0$, let
$$
g_r=\begin{pmatrix} r& 1\\1& 0\end{pmatrix}\in\GL_2(K)\;.
$$ 
We apply Theorem \ref{theo:algebrogeometricequid} with
$\Ga=\overline{G}$ the image of $G$ in $\PGL_2(R_v)$, $\DD^-=\;
  ]\alpha,\alpha^\sigma[$ the (image of any) geodesic line in $\XX_v$
with points at infinity $\alpha$ and $\alpha^\sigma$, and
$\DD^+=\ga_r\H_\infty$, where $\ga_r$ is the image of $g_r$ in
$\PGL_2(R_v)$.

We have $L_{\Ga_v}=2\ZZ$ and the family $\D^+= (\ga\DD^+)_{\ga \in
  \Ga/\Ga_{\DD^+}}$ is locally finite, as seen in the beginning of the
proof of Theorem \ref{theo:Mertensfunctionfieldfissgp}. The family
$\D^-= (\ga\DD^-)_{\ga\in\Ga/\Ga_{\DD^-}}$ is locally finite as seen
in the beginning of the proof of Theorem \ref{theo:equidloxofix}.

By Proposition \ref{prop:mescomputBT} (5), we have (on the full
measure subset $K_v-\{\alpha,\alpha^\sigma\}$ of $\partial_\infty \XX_v$)
$$
(\partial^+)_*\wt\sigma^+_{\DD^-}=
\frac{|\alpha-\alpha^\sigma|_v}{|z-\alpha|_v\,|z-\alpha^\sigma|_v}
\;d\haar_{K_v}(z)\,.
$$

As in order to obtain Equation \eqref{eq:applitwosubtreevenMert2},
since the point at infinity of $\ga\DD^+$ is $\ga\cdot r$, we have,
with an error term for every $\beta\in\; ]0,1]$ of the form
    $\bigO(s^{-\kappa}\|\psi\|_{\beta})$ for some $\kappa>0$ when
    evaluated on $\psi\in\C^{\beta}_c (\partial_\infty
    \XX_v-\partial_\infty\DD^-)$,
\begin{multline}
\lim_{n\ra+\infty} \;\frac{({q_v}^2-1)(q_v+1)}{2\,q_v^3}\;
\frac{\Vol(\Ga\dbs\XX_v)}{\|\sigma^-_{\D^+}\|}\;{q_v}^{-n}
\sum_{\substack{\ga\in \Ga/\Ga_{r}\\
0< d(\DD^-, \, \ga \DD^+)\leq n}}
\Dirac_{\ga\cdot r} \\=\;
\frac{|\alpha-\alpha^\sigma|_v}{|z-\alpha|_v\,|z-\alpha^\sigma|_v}
\;d\haar_{K_v}(z)\label{eq:applitwosubtreevennormform1}\;.
\end{multline}

We use the following result in order to switch from counting over
elements $\ga\in\Ga/\Ga_{r}$ for which $0< d(\DD^-, \, \ga \DD^+)\leq
t$ to counting over integral representations with bounded value of the
norm form. See \cite[page 1054]{ParPau11BLMS} for the analogous result
for the real hyperbolic $3$-space and indefinite binary Hermitian
forms.

\blemm \label{lem:cacldDpDmrepint}
Let $g\in \GL_2(R_v)$ and let $\ga$ be the image of $g$ in
$\PGL_2(K)$. Let $z_0=y_0$ if $y_0\neq 0$ and $z_0=x_0$
otherwise. Let $(x,y)=g(x_0,y_0)$. If $d(\DD^-,\ga\DD^+)>0$, then
$$
d(\DD^-,\ga\DD^+)= \frac{1}{\ln q_v}\;
\ln \Big(|\n(x-y\alpha)|_v\;\frac{h(\alpha)}{|z_0|_v^{\;2}}\Big)
\;.
$$
\elemm

\dem 
We start by showing that
$$
g\,g_r(1,0)=\Big(\frac{x}{z_0},\frac{y}{z_0}\Big)\;.
$$ 
Indeed, if $y_0\neq 0$, we have
$$
g\,g_r(1,0)=g(r,1)=\frac{1}{y_0}\;g(x_0,y_0)
$$  
and  otherwise
$$
g\,g_r(1,0)=g(1,0)=\frac{1}{x_0}\;g(x_0,0)= \frac{1}{x_0} \;
g(x_0,y_0)\,.
$$
In particular, 
$$
(g\, g_r)^{-1}=\begin{pmatrix} \;\;*&*\\-\frac{y}{z_0} &\frac{x}{z_0}
\end{pmatrix}\,.
$$ 
Note that $g\, g_r\in\GL_2(K)$ and $|\det (g\,g_r)|_v=|\det g|_v
\;|\det g_r|_v =1$ since $g\in \GL_2(R_v)$.  By Proposition
\ref{prop:propricomplexity} (2), we hence have
\begin{equation}\label{eq:cacldDpDmrepint1}
h((g\, g_r)^{-1}\cdot \alpha)=
\Big|\n\Big(\frac{x}{z_0}-\frac{y}{z_0}\;\alpha\Big)\Big|_v\;h(\alpha)
=|\n(x-y\;\alpha)|_v\;\frac{h(\alpha)}{|z_0|_v^{\;2}}\;.
\end{equation}

With $\beta_{\;\cdot\:}(\cdot,\cdot)$ the Busemann function defined in
Equation \eqref{eq:buscocycastree}, we use the signed distance $d(L,
H)=\min_{x\in L} \beta_{\xi}(x, x_H)$ between a geodesic line $L$ and
a horoball $H$ centred at $\xi\neq L^\pm$, where $x_H$ is any point of
the boundary of $H$. Now, by Equations \eqref{eq:valuationham} and
\eqref{eq:defidisthamenbord}, we have
\begin{align}
d(\DD^-,\ga\DD^+) & = d(]\alpha,\alpha^\sigma[\,,\ga\ga_r\H_\infty)=
d\big(](g\, g_r)^{-1}\cdot \alpha,(g\, g_r)^{-1}\cdot \alpha^\sigma[\,,
\H_\infty\big)
\nonumber\\ & =
v\big((g\, g_r)^{-1}\cdot \alpha-(g\, g_r)^{-1}\cdot \alpha^\sigma\big)
\nonumber\\ &
= \frac{-\ln \big|(g\, g_r)^{-1}\cdot \alpha-
(g\, g_r)^{-1}\cdot \alpha^\sigma\big|_v}{\ln q_v}
 = \frac{\ln h((g\, g_r)^{-1}\cdot \alpha)}{\ln q_v}\;.
\label{eq:cacldDpDmrepint2}
\end{align}
Combining Equations \eqref{eq:cacldDpDmrepint1} and
\eqref{eq:cacldDpDmrepint2} gives the result.  
\cqfd

\medskip By discreteness, there are only finitely many double classes
$[g]\in G_{\alpha}\bs G/ G_{(x_0,y_0)}$ such that $\DD^-=
\;]\alpha,\alpha^\sigma[$ and $g\DD^+= g\,g_r\H_\infty$ are not
disjoint. Let $Z(G)$ be the centre of $G$, which is finite. Since
$Z(G)$ acts trivially on $\PP_1(K_v)$, the map $G/G_{(x_0,y_0)}\ra
\Ga/\Ga_{r}$ induced by the canonical map $\GL_2(R_v)\ra
\PGL_2(R_v)$ is $|Z(G)|$-to-$1$.  Using the change of variable
$$
s =\frac{{|z_0|_v}^2}{h(\alpha)}\;{q_v}^{n}\;,
$$ 
and Lemma \ref{lem:cacldDpDmrepint} since $\ga\cdot r=\frac{x}{y}$
with the notation of this lemma, Equation
\eqref{eq:applitwosubtreevennormform1} gives
\begin{align*}
\lim_{s\ra+\infty} 
\;\frac{({q_v}^2-1)\;(q_v+1)\;{|z_0|_v}^2}{2\;{q_v}^{3}\;|Z(G)|
}\;
\frac{\Vol(\Ga\dbs\XX_v)}{\|\sigma^-_{\D^+}\|}\;s^{-1}&
\sum_{(x,\,y)\in G(x_0,\,y_0),\;\; |\n(x-y\alpha)|_v\leq s}
\Dirac_{\frac xy} \\=\;\; &
\frac{d\haar_{K_v}(z)}{|z-\alpha|_v\,|z-\alpha^\sigma|_v}\;,
\end{align*}
with the appropriate error term. Replacing $\Vol(\Ga\dbs\XX_v)$ and
$\|\sigma^-_{\D^+}\|$ by their values respectively given by Equation
\eqref{eq:covoltotMertens} and Lemma \ref{lem:totmassskinhoro}, the
claim of Theorem \ref{theo:mainnormform} follows. \cqfd

%% file: appendixBuzzi.tex
\addcontentsline{toc}{chapter}{{\large ~~~~Appendix}}

\chapter{A weak Gibbs measure is the unique equilibrium, by J. Buzzi}
\label{appendixBuzzi}


\newcommand{\barsig}{\ov\sigma}
\newcommand{\barmu}{{\ov\mu}}
\newcommand{\barm}{{\overline{m}}}
\newcommand{\barphi}{\ov\phi}
\newcommand{\Fix}{{\operatorname{Fix}}}
\newcommand\Prob{\PP}
\newcommand\Proberg{\PP_{\rm erg}}
\newcommand{\var}{\operatorname{var}}

\newtheorem{step}{\bf Step}


In this Appendix, for a transitive topological Markov shift endowed with a
Hölder-con\-ti\-nuous potential, we prove that a weak Gibbs measure is
the unique equilibrium measure.

\section{Introduction}

Let $\sigma:\Sigma\to\Sigma$ be a topological Markov shift (possibly
one- or two-sided), see for instance Section \ref{subsec:TMS}. More
precisely, we consider the one-sided and two-sided vertex-shifts
defined by a countable oriented graph $G$ with set of vertices $V_G$
and set of arrows $A_G\subset V_G\times V_G$. We assume that $\Sigma$
is transitive, that is, that $G$ is connected (as an oriented graph).

We denote by $\Prob(\Sigma)$ the set of $\sigma$-invariant probability
measures on $\Sigma$ and by $\Proberg(\Sigma)$ the subset of ergodic
ones. Recall that, for all $n\in\NN$, the {\it $n$-cylinders} are the
following subsets of $\Sigma$, where $x$ varies in $\Sigma$:
$$
C_n(x)=[x_0,\dots , x_{n-1}]=
\{y\in\Sigma\;:\;\forall\;k\in\{0,\dots, n-1\},\; y_k=x_k\}\;,
$$
so that the $1$-cylinders are $[v]=\{y\in\Sigma\;:\; y_0=v\}$ for all
$v\in V_G$.  The points of $\Sigma$ admitting $n\in\NN$ as period
under the shift $\sigma$ form the set
$$
\Fix_n(\Sigma)=\{x\in\Sigma\;:\;\sigma^nx=x\}\;.
$$
We fix a {\it potential} on $\Sigma$, that is, a continuous function
$\phi:\Sigma\to\RR$. We do not assume that $\phi$ is bounded. We
define $\phi^-=\max\{-\phi,0\}$ and, for all $n\in\NN-\{0\}$,
$$
\var_n(\phi)=\sup_{x,y\in\Sigma\,:\; \forall\,k \in\{0,\dots, n-1\},\;x_k=y_k} 
\;|\phi(y)-\phi(x)|
$$
if $(\Sigma,\sigma)$ is one-sided and otherwise
$$
\var_n(\phi)=\sup_{x,y\in\Sigma\,:\; \forall\,k \in\{-n+1,\dots, n-1\},\;x_k=y_k} 
\;|\phi(y)-\phi(x)|\;.
$$
We say that $\phi$ has {\it summable variations} if $\sum_{n\geq1}
\var_n(\phi)<\infty$. This is in particular the case if $\phi$ is
H\"older-continuous. Let $S_n\phi= 
\sum_{i=0}^{n-1}\phi\circ\sigma^i$ for all $n\in\NN$.

\bdefi
A \emph{weak Gibbs measure}\index{Gibbs!measure!weak} for the
potential $\phi$ is a $\sigma$-invariant Borel probability measure $m$
on $\Sigma$ such that there exists a number $c(m)\in\RR$ such that for
every $v\in V_G$, there exists $C\geq 1$ with
\begin{equation}\label{eq.wGibbs} \forall\;
n\geq1,\;\forall x\in\Fix_n(\Sigma)\cap[v],\quad
C^{-1}\leq \frac{m(C_n(x))}{\exp\left(S_n\phi(x)-c(m)n\right)} \leq
C\;.  
\end{equation}
\edefi

Note that $c(m)$ is then unique; it is called the {\it Gibbs
constant}\index{Gibbs!constant} of $m$. Let us stress that we do not
assume the so-called \emph{Big Image Property} \cite{Sarig03} and
hence using the above weakened Gibbs property (that is, allowing $C$ to
depend on $v$) is necessary.

Note that if $\Sigma$ is locally compact, that is, if every vertex of
$G$ has finite degree (the number of arrows arriving or leaving from
the given vertex), then the above condition is equivalent to the fact
that for any nonempty compact subset $K$ in $\Sigma$, there exists
$C\geq 1$ with
$$
\forall\; n\geq 1,\;\forall\; x\in\Fix_n(\Sigma)\cap K,\quad 
C^{-1}\leq\frac{m(C_n(x))}{\exp\left(S_n\phi(x)-c(m)n\right)} \leq C\;.
$$

The \emph{pressure}\index{pressure} $P(\phi,\nu)$ of an element
$\nu\in \Prob(\Sigma)$ such that $\int\phi^-\,d\nu<+\infty$ is
$$
P(\phi,\nu)=h_\nu(\sigma)+\int\phi\,d\nu\;.
$$
An \emph{equilibrium measure}\index{equilibrium measure} $\mu_{eq}$
for $(\Sigma,\phi)$ is an element $\mu_{eq}\in\Prob(\Sigma)$ such that
$\int\phi^-\,d\mu_{eq}<+\infty$ and
$$
P(\phi,\mu_{eq}) = 
\sup \{P(\phi,\nu): \nu\in\Prob(\Sigma) \text{ and } 
\int\phi^-\,d\nu<+\infty\}\;.
$$
The \emph{Gurevi\v{c} pressure}\index{pressure!Gurevi\v{c}}%
\index{Gurevi\v{c} pressure} is 
$$
P_G(\phi)=\limsup_{n\to\infty}\;\frac1n
\log\sum_{x\in\Fix_n(\Sigma)\cap[v]} e^{S_n\phi(x)}
$$
for any vertex $v\in V_G$. Note that the Gurevi\v{c} pressure does not
depend on $v$. Let us recall a few results on the above notions.

\btheo[Iommi-Jordan \mbox{\cite[Theorem 2.2]{IomJor13}}]\label{t-vp}
If $\phi$ has summable variations, the following variational principle holds:
$$
P_G(\phi)= \sup \{P(\phi,\nu): \nu\in\Prob(\Sigma) \text{ and } 
\int\phi^-\,d\nu<+\infty\}\;.\;\;\;\Box
 $$
\etheo

\btheo[Buzzi-Sarig \mbox{\cite[Theorem 1.1]{BuzSar03}}]\label{t-unique}
Assume that $\phi$ has summable variations.

If $P_G(\phi)<\infty$, then there exists at most one equilibrium
measure.

If there exists an equilibrium measure $\mu$, then $d\mu=h\,d\nu$ where
$h:\Sigma\to\RR$ is a continuous, positive function and $\nu$ is a
positive measure with full support on $\Sigma$ such that 
\begin{enumerate} 
\item[$\bullet$] $L_\phi\, h=e^{P_G(\phi)} h$, and $L_\phi^*\,\nu=
e^{P_G(\phi)}\nu$ where $L_\phi$ is the transfer operator defined by
$L_\phi \,u\,(x) = \sum_{y\in\sigma^{-1}x} \,e^{\phi(y)}\,u(y)$.
\item[$\bullet$] $\nu$ is finite on each cylinder. \cqfd
\end{enumerate} 
\etheo

We note that \cite{BuzSar03} assumed $\sup\phi<\infty$, but this was
only used to justify the variational principle and so this condition
can be removed by using Theorem \ref{t-vp}.

We now state the main result of this appendix.

\btheo\label{theo:mainappendix_onesided}
Let $(\Sigma,\sigma)$ be a one-sided transitive topological Markov
shift and let $\phi:\Sigma\to\RR$ be a potential with summable
variations. Let $m$ be a $\sigma$-invariant probability measure on
$\Sigma$ such that $\int \phi^-dm<+\infty$.

Then $m$ is a weak Gibbs measure if and only if it is an equilibrium
measure. In this case, the Gibbs constant $c(m)$ is equal to the
Gurevi\v{c} pressure and the equilibrium measure is unique.
\etheo

By a classical argument that follows, this result extends to two-sided
topological Markov shifts (up to a slight strengthening of the
regularity assumption on $\phi$, still satisfied if $\phi$ is
H\"older-continuous).

\bcoro\label{coro:mainappendix_twosided}
Let $(\Sigma,\sigma)$ be a two-sided transitive topological Markov
shift and let $\phi:\Sigma\to\RR$ be a potential with $\sum_{n\geq1}
n\;\var_n(\phi)<\infty$. Let $m$ be a $\sigma$-invariant probability
measure on $\Sigma$ such that $\int \phi^-dm<+\infty$.

Then $m$ is a weak Gibbs measure if and only if it is an equilibrium
measure. In this case, the Gibbs constant $c(m)$ is equal to the
Gurevi\v{c} pressure and the equilibrium measure is unique.
\ecoro

\rem
The case of the full shift $\NN^\ZZ$ has been treated
in \cite[Sec. 3]{PesSenZha16}. More generally, assuming the Big Image
Property, the above result follows from \cite{Sarig03}
and \cite{BuzSar03} along the lines of \cite{PesSenZha16}.

\bigskip
\noindent{\bf Proof of Corollary \ref{coro:mainappendix_twosided}. } 
Let $(\Sigma,\sigma)$, $\phi$, and $m$ be as in the statement of this
Corollary.  Let $\pi:\Sigma\to\Sigma_+$ with $(x_n)_{n\in\ZZ}\mapsto
(x_n)_{n\in\NN}$ be the obvious factor map onto the one-sided
topological Markov shift $(\Sigma_+,\sigma_+)$ defined by the same
graph $G$ as for $(\Sigma,\sigma)$, called the {\it natural
extension}.\index{natural extension}

First, we replace $\phi$ by a potential $\ov\phi$ depending only on
future coordinates. The proof of \cite[Lemma 1.6]{Bowen75} applies to
our non-compact setting without changes. To be more precise, for each
vertex $a\in V_G$, choose $z^a\in\Sigma$ with $z^a_0=a$. Define
$r:\Sigma\to\Sigma$ by $r(x)=y$ with $y_n=x_n$ for $n\geq0$ and
$y_n=z^{x_0}_n$ for $n\leq 0$. For every $x\in\Sigma$,  let
$$
u(x) = \sum_{k\geq0} (\phi\circ\sigma^k-\phi\circ\sigma^k\circ r)(x)\;.
$$
This defines a bounded real function on $\Sigma$ since
$|\phi \circ \sigma^k -\phi\circ\sigma^k\circ r|\leq\var_{k+1}(\phi)$
and $\phi$ has summable variations. Moreover, $u$ itself has summable
variations since, given $x,y\in\Sigma$ with $x_k=y_k$ for $|k|<n$, we
have
$$\begin{aligned}
|u(x)-u(y)| &\leq \sum_{0\leq k< \lfloor n/2\rfloor} 
\left(|\phi(\sigma^kx)-\phi(\sigma^ky)|
+ |\phi(\sigma^k(rx))-\phi(\sigma^k(ry))|\right) \\
&\quad + \sum_{k\geq\lfloor n/2\rfloor} 
\left(|\phi\circ\sigma^k(x)-\phi\circ\sigma^k(r x)|+
|\phi\circ\sigma^k(y)-\phi\circ\sigma^k(r y)|\right)\\
&\leq 4\sum_{k\geq \lfloor n/2\rfloor} \var_{k+1}(\phi)\;,
\end{aligned}$$
so that
$$
\sum_{n\geq1}\var_n(u)\leq 8\sum_{k\geq1} k\var_k(\phi)<\infty\;.
$$ 
Now define $\ov\phi:\Sigma\to\RR$ by
$$
\ov\phi = \phi + u\circ\sigma - u\;.
$$
The function $\ov\phi$ is continuous with summable variations.
Following \cite{Bowen75}, let us prove that $\ov\phi=\ov\phi\circ r$. We
have
$$\begin{aligned}
    \ov\phi &= \phi
     + \sum_{k\geq0}(\phi\circ\sigma^{k+1}-\phi\circ\sigma^k\circ r\circ \sigma) 
     - \sum_{k\geq0} (\phi\circ\sigma^k-\phi\circ\sigma^{k}\circ r)\\
     & = \phi-\phi -\sum_{k\geq0}
\left(\phi\circ\sigma^k\circ r\circ\sigma-\phi\circ\sigma^k\circ r\right)\\
     &=\sum_{k\geq0}\left(\phi\circ\sigma^k\circ r-
\phi\circ\sigma^k\circ r\circ\sigma\right)\;.
\end{aligned}$$
Now, $r^2=r$ and $r\circ\sigma\circ r=r\circ\sigma$. Hence
$\ov\phi \circ r=\ov\phi$ as claimed.  Thus, $\ov\phi$ induces on the
one-sided shift a function $\wt\phi:\Sigma_+\ra\RR$ defined by
$$
\tilde\phi: (x_n)_{n\in\NN}\mapsto 
\ov\phi(\dots z^{x_0}_{-2} z^{x_0}_{-1}x_0x_1\dots) \;,
$$
satisfying $\ov\phi=\wt\phi\circ\pi$.

To conclude, observe that $S_n\ov\phi(x)-S_n\phi(x)=
S_n(u\circ\sigma-u)(x) =0$ if $x\in \Fix_n(\sigma)$, and that
cylinders defined by the same finite words have the same measure for
an invariant probability measure $m$ on the two-sided shift
$(\Sigma,\sigma)$ and for its image $\pi_*m$ on the one-sided shift
$(\Sigma_+,\sigma_+)$. Therefore $m$ is a weak Gibbs measure for
$\phi$ if and only if $\pi_*m$ is a weak Gibbs measure for $\wt\phi$,
and their Gibbs constant are then equal.

By construction $\pi_*m(\wt\phi)=m(\ov\phi)=m(\phi)$ since $m$ is
invariant. As it is well-known, the natural extension $\pi$ preserves
the entropy. Thus, the measure $m$ is an equilibrium measure with
respect to $\phi$ if and only if $\pi_*m$ is an equilibrium measure
with respect to $\wt\phi$.

The reduction to one-sided topological Markov shifts is thus complete.
\cqfd

\section{Proof of the main result Theorem \ref{theo:mainappendix_onesided}}

The uniqueness of the equilibrium state is given by
Theorem \ref{t-unique}.  We need to prove that weak Gibbs measures and
equilibrium measures coincide under the integrability assumption on
$\phi^-$ and that the number $c(m)$ is equal to the Gurevi\v{c} pressure.

\begin{step}\label{s.converse}
If $m$ is an equilibrium measure, then it is a weak Gibbs measure.
\end{step}

This is a routine consequence of Theorem \ref{t-unique}.  Our definition
of an equilibrium measure $m$ enforces $\int\phi^-\,dm<+\infty$
(hences excludes the concomitance of $h_m(\sigma)=+\infty$ and
$\int \phi\,dm=-\infty$).

Recall from Theorem \ref{t-unique} that $dm=h\,d\nu$ with $h$ and
$\nu$ as mentionned. For $v\in V_G$ and $x\in\Fix_n(\Sigma)\cap[v]$,
we have
\begin{align*}
m(C_n(x)) &= \int h\,\mathbbm{1}_{C_n(x)}\, d\nu
= e^{-n\, P_G(\phi)} \int h\,\mathbbm{1}_{C_n(x)} \, d((L_\phi^*)^n\nu)\\ &
= e^{-n\, P_G(\phi)} \int L_\phi^n(h\,\mathbbm{1}_{C_n(x)}) \, d\nu\;.
\end{align*}
By definition,
$$
L_\phi^n(h\,\mathbbm{1}_{C_n(x)})(z)=\exp (S_n\phi(x_0\dots
x_{n-1}z))\;h(x_0\dots x_{n-1}z)
$$
for all $z\in \sigma^n(C_n(x))=\sigma([v])$ (and
$L_\phi^n(h\,\mathbbm{1}_{C_n(x)})(z)=0$ otherwise). Hence\footnote{We
use, for all $u,v,d\geq 0$ and $c> 0$, the notation $u=c^{\pm d} v$ if
$\frac{1}{c^d}v\leq u \leq c^dv$.}
$$
m(C_n(x))= e^{-n \,P_G(\phi)}
\exp\Big( S_n\phi(x)\pm\sum_{k=1}^{n}\var_k(\phi)\Big) 
\int_{\sigma([v])} h\, d\nu\;.
$$
As $0< \int_{\sigma([v])}h\,d\nu<+\infty$ and $\sum_{k=1}^{+\infty}\var_k(\phi)
<+\infty$, the measure $m$ is a weak Gibbs measure for $\phi$ with
Gibbs constant $c(m)=P_G(\phi)$.

\medskip
We now turn to the converse implication. Let $m$ be a weak Gibbs
measure for $\phi$ such that $\int\phi^-\,dm<+\infty$.

\medskip
The weak Gibbs condition only controls the cylinders that start and
end with the same symbol. Passing to an induced system (that is,
considering a first return map on a $1$-cylinder) will remove this
restriction. More precisely, let $a\in V_G$ be a vertex of $G$ and let
$\mu$ be an invariant probability measure on $(\Sigma,\sigma)$ with
$\mu([a])>0$.  The \emph{induced system}\index{induced!system} on the
$1$-cylinder $[a]=\{x\in\Sigma\;:\;x_0=a\}$ is the map $\barsig:[a]\to[a]$
(almost everywhere) defined as follows:
\begin{enumerate}
\item[$\bullet$]
let $\tau(x)=\inf\{n\geq1\;:\;\sigma^nx\in[a]\}$ be the {\it
first-return time} in $[a]$, that we also denote by $\tau_{[a]}(x)$ when
we want to emphasize $[a]$;
\item[$\bullet$]
let $\barsig(x)=\sigma^{\tau(x)}(x)$ if $\tau(x)<\infty$; 
\item[$\bullet$]
let $\barmu(B)=\mu(B\cap[a])/\mu([a])$ for every Borel subset $B$ of
$\Sigma$ be the restriction of $\mu$ to $[a]$ normalized to be a
probability measure.
\end{enumerate}
We also define $\tau^0(x)=0$ and by induction $\tau^{n+1}(x)
= \tau(x)+\tau^n(\barsig x)$ for every $n\in\NN$. Note that $\barsig$
can only be iterated on the subset
$$
\{x\in[a]\;:\;\forall\; n\geq1,\;\tau^n(x)<\infty\}\;.
$$
By Poincar\'e's recurrence theorem, this is a full measure subset of
$[a]$, hence the distinction will be irrelevant for our purposes.

The \emph{induced partition}\index{induced!partition} is 
$$
\beta= \{[a, \xi_1, \dots, \xi_{n-1}, a]\ne\emptyset:n\geq1,\xi_i\ne a\}\;.
$$
We note that $\barsig:[a]\to[a]$ is topologically Bernoulli with
respect to the partition $\beta$ (that is, $\barsig:b\to[a]$ is a
homeomorphism for each $b\in\beta$). For every integer $N\geq1$, we
define the $N$-th {\it iterated partition}\index{iterated partition} 
$\beta^N$ of $\beta$ by
$$
\beta^N=\{b_0\cap\barsig^{-1}b_1\cap\dots\cap\barsig^{-N+1}b_{N-1}
\ne\emptyset\;:\;  b_0,\dots,b_{N-1}\in\beta\}
$$
and we write $\beta^N(x)$ for the element of the partition $\beta^N$
that contains $x$.

\begin{step}\label{s.mixing}
The topological Markov shift may be assumed to be topologically
mixing.
\end{step} 

This follows from the spectral decomposition for topological Markov
shifts, see for instance \cite[Lem.~2.2]{BuzSar03}.

\medbreak
\begin{step}\label{s.ergodicity}
The Gibbs property implies full support and ergodicity.
\end{step}

Let $A$ be a $\sigma$-invariant ($\sigma^{-1}(A)=A$) measurable subset of
$\Sigma$ with $m(A)>0$ and let us prove that $m(A)=1$.

Observe that the Gibbs property, together with the transitivity of
$\Sigma$, implies that any cylinder has positive measure for $m$,
hence that $m$ has full support. Let $a\in V_G$ be such that
$m(A\cap[a]) >0$.

As $m([a])>0$, we may consider the induced system on $[a]$. Let $N\geq
1$. When $f$ is a homeomorphism between topological spaces, let $f^*$
denote the pushforwards of measures by $f^{-1}$. First note that, for
almost every $x\in [a]$ and every $N\in\NN-\{0\}$, since $\barsig^N$
is an homeomorphism from $\beta^N(x)$ onto $[a]$, we have
$$
\frac{m(A\cap[a])}{m([a])}
= \frac{m(\barsig^N(A\cap\beta^N(x)))}{m(\barsig^N(\beta^N(x)))}
= \frac{\int_{\beta^N(x)\cap A} \frac{d(\barsig^N)^*m}{dm} \, dm}
{\int_{\beta^N(x)} \frac{d(\barsig^N)^*m}{dm}\, dm}\;.
$$
Now, observe that for
$m$-almost every $y\in\beta^N(x)$:
$$
\frac{d(\barsig^N)^*m}{dm}(y)= \lim_{n\to\infty} 
\frac{m(\sigma^{\tau^N(y)}[y_0,\dots ,y_n])}{m([y_0,\dots ,y_n])}=
C^{\pm2} e^{-S_{\tau^N(y)}\phi(y)+\tau^N(y)\,c(m)}\;.
$$
Hence, since $\tau^N$ is constant on $\beta^N(x)$,
$$
\frac{m(A\cap[a])}{m([a])}
= C^{\pm4}\; \;e^{\pm \sum_{k\geq 1} \var_k(\phi)}\;
\;\frac{m(A\cap\beta^N(x))}{m(\beta^N(x))}\;.
$$
By Doob's increasing martingale convergence theorem (see for
instance \cite{Petersen89}), for $m$-almost every $x\in[a]-A$,
the ratio on the right hand side converges to $0$ as
$N\to\infty$. Thus $[a]$ is contained in $A$ modulo $m$. Therefore
$A=\bigcup_{a\in W} [a]$ modulo $m$ for some subset $W$ of $V_G$.

Since $\Sigma$ is topologically mixing, for any vertex $b$, the
intersection $[a]\cap\sigma^{-i}[b]\cap\sigma^{-j}[a]$ is not empty
for some integers $0<i<j$. Pick some point $x$ in that set. By
invariance, $m([b])\geq m(\sigma^i(C_j(x)))\geq m(C_j(x))$. But this
last number is positive by the weak Gibbs property.  Thus $[b]$ is
contained in $A$ modulo $m$. Hence $m(A)=1$, proving the ergodicity of
$m$.

\begin{step}\label{s.pressure}
The Gurevi\v{c} pressure $P_G(\phi)$ is equal to $c(m)$, hence is
finite. Furthermore $h_m(\sigma)<\infty$ and $\phi\in \LL^1(m)$.
\end{step}

Fix $v\in V_G$ and let $K=[v]$. Note that $m(K)>0$.  The ergodicity of
$m$ gives a Cesaro convergence: as $n\to\infty$, we have
 $$\begin{aligned}
   \frac1n\sum_{k=0}^{n-1}m(K\cap\sigma^{-k}K) \longrightarrow m(K)^2>0\;.
 \end{aligned}$$
The Gibbs property implies that, for all $n\geq1$,
\begin{align}
m(K\cap \sigma^{-n}K) & = 
C^{\pm 1} \sum_{x\in\Fix_n(\Sigma)\cap K} e^{S_n\phi(x)-c(m)n} \nonumber\\ &  = 
C^{\pm1} \Big(\sum_{x\in\Fix_n(\Sigma)\cap K} e^{S_n\phi(x)}\Big)\, e^{-c(m)n}
\label{eq.mKK}\;.
\end{align}
If we write $Z_n$ for the term between the parenthesis, we have by the
definition of the Gurevi\v{c} pressure:
 $$
   P_G(\phi) = \limsup_{n\to\infty}\frac1n\log Z_n\;.
 $$ 
As the value of the left hand side of Equation \eqref{eq.mKK} is less
than one, we see that $c(m)\geq P_G(\phi)$. If this was a strict
inequality, then the left hand side of Equation \eqref{eq.mKK} would
converge to zero, contradicting its Cesaro convergence to $m(K)^2>0$.
Therefore $P_G(\phi)= c(m)$.

Since $c(m)$ is finite, so is $P_G(\phi)$. Hence Theorem \ref{t-vp}
implies that, for any $\nu\in\Prob(\Sigma)$ with
$\int \phi^-\,d\nu<+\infty$, we have $h_\nu(\sigma)<\infty$ and $\phi$
is $\nu$-integrable. In particular, this holds for $\nu=m$, which
finishes the proof of Step \ref{s.pressure}.

\begin{step}\label{s.rokhlin}
If the mean entropy $H_{\barmu}(\beta)=-\sum_{b\in\beta} \barmu(b) 
\log\barmu(b)$ is finite, then
$$
h_\barmu(\barsig) = -\int\log\frac{d\barmu}{d\barsig^*\barmu}\, d\barmu\;,
$$
where $\barsig^*\barmu$ is the measure on $\Sigma$ defined by
$B\mapsto \sum_{b\in\beta}\barmu(\barsig(B\cap b))$, with respect to
which $\mu$ is absolutely continuous: $\mu\llcurly\barsig^*\barmu$.
\end{step}

This is a classical formula, sometimes called the {\it Rokhlin
formula}\index{Rokhlin formula} (see for instance \cite{BuzSar03}),
which follows from the computation of the entropy in terms of the
information function when the mean entropy is finite
$$\begin{aligned}
h_\barmu(\barsig)&=-\int \sum_{b\in\beta} \mathbbm{1}_b(x) 
\lim_{n\to\infty} \log \mathbb E_\barmu(\mathbbm{1}_b\;|\;
\barsig^{-1}\beta\vee\dots\vee\barsig^{-n}\beta)(x)\, d\barmu(x)
 \end{aligned}$$
and from the identity, for $x\in b$,
$$
    \mathbb E_\barmu(\mathbbm{1}_b\;|\;
    \barsig^{-1}\beta\vee\dots\vee\barsig^{-n}\beta)(x)
    = \frac{\barmu(\beta^{n+1}(x))}{\barmu(\beta^n(\barsig x))}
    =  \frac{\barmu(\beta^{n+1}(x))}{\barsig^*\barmu(\beta^{n+1}(x))}\;.
$$
The absolute continuity follows from a direct computation and ensures
that the integral above is well-defined.
\bigbreak

\begin{step}\label{s.Bernoulli}
For all $a\in V_G$, $N\geq1$ and $\mu\in\Proberg(\Sigma)$ with
$\int \phi^-\,d\mu<+\infty$, we have
\begin{equation}\label{eq.Rokhlin}
 h_\mu(\sigma)=-\mu([a])\int_{[a]}\frac1N
\log\frac{d\,\barmu}{d((\barsig^{N})^*\,\barmu)}\,
 d\,\barmu 
\end{equation}
\end{step}

We use arguments from the proof of \cite[Theorem 1.1]{BuzSar03}: the
key is to see that the induced partition $\beta$ of $[a]$ has finite
mean entropy for the induced measure $\bar\mu$ using a Bernoulli
approximation.

Let us consider the Bernoulli measure $\barmu_B$ for $([a],\barsig)$
defined by
$$
\barmu_B\big(\bigcap_{i=0}^{n-1}\barsig^{-i}B_i\big)=\prod_{i=0}^{n-1}\barmu(B_i)
$$ 
for all $B_i\in\beta$. We construct from it an invariant measure
$\mu_B$ on $(\Sigma,\sigma)$: For every Borel subset $A$, let
$$
\mu_B(A) = \mu([a]) \int_{[a]} 
\sum_{i=0}^{\tau_{[a]}-1} \mathbbm{1}_A\circ\sigma^i \, d\,\barmu_B\;.
$$
Note that $\mu_B$ is ergodic, since $\barmu_B$ is ergodic and
$\bigcup_{i\geq 1}\barsig^{-i}([a])$ has full measure. Define
$\barphi=\sum_{i=0}^{\tau_{[a]}-1} \phi\circ\sigma^i$. Let
$C'=\sup_{k\geq 1}\var_k(\phi)$, which is finite since $\phi$ has
summable variations. Since every $b\in\beta$ is a cylinder of length
$\tau_{[a]}(x)+1$ for every $x\in b$, the conditional expectation
$$
\mathbb{E}_{\barmu}(\barphi\,|\, \beta)= \sum_{b\in\beta} \mathbbm{1}_b\;
\frac{1}{\barmu(b)}\;\int_b\barphi\,d\bar\mu
$$
satisfies $\| \barphi-\mathbb{E}_{\barmu}(\barphi \,|\,
\beta)\|_\infty\leq C'$. Hence
\begin{align*}
\int \phi \,d\mu_B &
=\mu([a]) \int_{[a]} \barphi\; d\,\barmu_B \geq \mu([a])\int_{[a]}
\big(\mathbb{E}_{\barmu}(\barphi \,|\, \beta)-C'\big)\; d\,\barmu_B
\\ &
\geq
\mu([a])\int_{[a]}  \barphi\; d\,\barmu-C'
=\int \phi\, d\mu -C' >-\infty\;.
\end{align*}
Therefore, the last paragraph of the proof of Step \ref{s.pressure}
applies to $\nu=\mu_B$ and $h_{\mu_B}(\sigma)<+\infty$. Since $\mu_B$ is
ergodic, Abramov's formula yields
$$
h_{\mu_B}(\sigma)=\mu([a])\, h_{\barmu_B}(\barsig)\;.
$$ 
Since $\barmu_B$ is Bernoulli, the right hand side of this equality is
equal to
$$
\mu([a])\,H_{\barmu_B}(\beta)=\mu([a])H_{\barmu}(\beta)\;,
$$ 
so that $H_{\barmu}(\beta)$ is proven to be finite. Thus,
Step \ref{s.rokhlin} applies:
$$
h_\barmu(\barsig)= 
-\int_{[a]} \log\frac{d\barmu}{d(\barsig^*\barmu)}\, d\barmu\;.
$$ 
This formula extends to $h_\barmu(\barsig^N)$ for all integers $N\geq
1$.  Using Abramov's formula this time for $\mu$ and $\barmu$ (since
$\mu$ is ergodic), we have
$$
h_\mu(\sigma)=\mu([a]) \,h_\barmu(\barsig)=
\frac{\mu([a])}{N} \,h_\barmu(\barsig^N) =
-\mu([a])\int_{[a]} \frac{1}{N}\,
\log\frac{d\barmu}{d((\barsig^N)^*\barmu)} \, d\barmu \;,
$$
as claimed.

\begin{step}\label{s.entropy}
The entropy of $m$ is equal to $c(m)-\int\phi\,dm$.
\end{step}

In order to prove this, we apply Step \ref{s.Bernoulli} with $\mu=m$
(which is possible, since $m$ has been proven to be ergodic in
Step \ref{s.ergodicity}). As in the proof of Step \ref{s.ergodicity},
the Radon-Nikodym derivative is almost everywhere
$$
\frac{d\barm}{d((\barsig^N)^*\barm)}(x) = \lim_{n\to\infty} 
\frac{\barm(\beta^n(x))}{\barm(\barsig^N(\beta^n(x)))}
= C^{\pm2} \exp \big(S_{\tau^N(x)}\phi(x)-c(m)\tau^N(x)\big)\;.
$$
Therefore, using Step \ref{s.Bernoulli} and the fact that
$m_{\mid[a]}=m([a])\,\barm$, we have
\begin{equation}\label{eq.hint}
h_m(\sigma) = \lim_{N\to\infty}\frac1N \Big(\int_{[a]} 
\big(c(m)\,\tau^N(x)-S_{\tau^N(x)}\phi(x)\big)\, dm\pm2\log C\Big)\;.
\end{equation}
Note that $\tau^N(x)$ can be seen as a Birkhoff sum for the induced
system on $[a]$ and the function $\tau$ and that, by Kac's theorem
(see for instance \cite[Sect.~2.4]{Petersen89}),
$$
\int_{[a]} \tau\;d\,\barm = m([a])^{-1}\;.
$$ 
Therefore, Birkhoff's ergodic theorem yields, with convergence in
$\LL^1(\barm)$,
$$
\lim_{N\to\infty}\frac{c(m)\,\tau^N(x)}N = 
\frac{c(m)}{m([a])}\;.
  $$
In order to analyze the second term in Equation \eqref{eq.hint}, let
$\hat\phi(x)=\sum_{k=0}^{\tau(x)-1}\phi(\sigma^kx)$ and observe that,
by a variation of the proof of Kac's theorem, $\hat\phi \in 
\LL^1(\barm)$ with $\barm(\hat\phi) = m([a])^{-1}m(\phi)$.  
Indeed, passing to the natural extension, one can assume the system to
be invertible and use the partition modulo $m$ given by
$$
\bigcup_{n\geq1,\,0\leq k<n} \sigma^k(\{x\in[a]\;:\;\tau(x)=n\})\;.
$$
Since $S_{\tau^N(x)}\phi(x)$ coincides with the Birkhoff sum $\bar
S_N\hat\phi(x)$ for the induced system, Birkhoff's ergodic theorem
yields, with convergence in $\LL^1(\barm)$,
$$
\lim_{N\to\infty} \frac1N S_{\tau^N(x)}\phi(x) = \lim_{N\to\infty} 
\frac1N\bar S_N\hat\phi(x)=m([a])^{-1}m(\phi)\;.
$$
The claim follows.

\begin{step}
Conclusion: any weak Gibbs measure is an equilibrium measure and
$c(m)=P_G(\phi)$.
\end{step}

Steps \ref{s.pressure} and \ref{s.entropy} prove that $h_m(\sigma)
+ \int\phi\, dm$ is well-defined and equal to $c(m)$, which by
Step \ref{s.pressure} is equal to $P_G(\phi)$, which is equal to
$\sup \{P(\phi,\nu): \nu\in\Prob(\Sigma) \text{ and }
\int \phi^-\,d\nu<+\infty\}$ by Theorem \ref{t-vp}, so that $m$ is an
equilibrium measure. This completes the proof of
Theorem \ref{theo:mainappendix_onesided}. \cqfd